\documentclass[11pt,letterpaper]{article}
\usepackage[T1]{fontenc}
\usepackage[utf8]{inputenc}
\usepackage{lmodern}
\usepackage[margin=1in]{geometry}
\usepackage{amsmath,amssymb,amsthm,mathtools}
\usepackage{aliascnt,needspace}
\usepackage{array,booktabs,longtable}
\usepackage{enumitem,microtype,xcolor}
\usepackage{xurl}
\usepackage[colorlinks=true,linkcolor=blue!45!black,citecolor=blue!45!black,urlcolor=blue!45!black]{hyperref}
\usepackage[nameinlink,noabbrev]{cleveref}
\usepackage{listings}
\usepackage{fancyhdr}
\setlist{topsep=4pt,itemsep=3pt,parsep=1pt}

\newtheorem{theorem}{Theorem}[section]
\newaliascnt{lemma}{theorem}
\newtheorem{lemma}[lemma]{Lemma}
\aliascntresetthe{lemma}
\newaliascnt{proposition}{theorem}
\newtheorem{proposition}[proposition]{Proposition}
\aliascntresetthe{proposition}
\newaliascnt{corollary}{theorem}
\newtheorem{corollary}[corollary]{Corollary}
\aliascntresetthe{corollary}
\theoremstyle{definition}
\newaliascnt{definition}{theorem}

\aliascntresetthe{definition}
\newaliascnt{remark}{theorem}

\aliascntresetthe{remark}
\newtheorem*{example}{Example}

\crefname{theorem}{Theorem}{Theorems}
\Crefname{theorem}{Theorem}{Theorems}

\crefname{lemma}{Lemma}{Lemmas}
\Crefname{lemma}{Lemma}{Lemmas}

\crefname{proposition}{Proposition}{Propositions}
\Crefname{proposition}{Proposition}{Propositions}

\crefname{corollary}{Corollary}{Corollaries}
\Crefname{corollary}{Corollary}{Corollaries}

\crefname{definition}{Definition}{Definitions}
\Crefname{definition}{Definition}{Definitions}

\crefname{remark}{Remark}{Remarks}
\Crefname{remark}{Remark}{Remarks}

\crefname{equation}{equation}{equations}
\Crefname{equation}{Equation}{Equations}

\crefname{section}{Section}{Sections}
\Crefname{section}{Section}{Sections}
\crefname{subsection}{Section}{Sections}
\Crefname{subsection}{Section}{Sections}

\crefname{table}{Table}{Tables}
\Crefname{table}{Table}{Tables}

\newcommand{\Q}{\mathbb Q}
\newcommand{\Z}{\mathbb Z}

\numberwithin{equation}{section}

\title{\textbf{Efficient differencing and linear bounds\\ in the Waring--Siegel problem}}
\author{Morley Davidson}
\date{27 November 2026\\{\normalsize Revised 1 December 2026}}
\hypersetup{pdfauthor={Morley Davidson},pdftitle={Efficient differencing and linear bounds in the Waring--Siegel problem},pdfsubject={Additive number theory}}
\begin{document}
\maketitle
\begin{abstract}
Let $G_{\mathbb K}(k)$ be the Waring number for totally positive
summands in a number field, and let $H_{\mathbb K}(k)$ be the
corresponding primitive local width. We prove that
$G_{\mathbb K}(k)\le\max\{H_{\mathbb K}(k),\lfloor6.873k\rfloor\}$
for every field and every $k\ge3$. For large $k$, with a cutoff
independent of the field, the analytic coefficient is $6.216216302$.
We give effective coefficients beginning at $k=120$, including
$6.258$ for $k\ge10^{12}$. For each fixed field the limiting upper
ratio to $k$ is at most $6.216216302$. Uniformly over all fields
one has $H_{\mathbb K}(k)\ll k\log(2k)$. On families with bounded
ramification at each fixed residue cardinality the width is $o(k\log k)$.
We determine how often the
local width exceeds a prescribed fraction of $k$. An extension of
efficient differencing to fractional ideal lattices gives leading
coefficient one in the historical $k\log k$ estimate for every
class number. Seven cubes suffice in every number field. Among
eligible targets, four cubes have relative density $(23/27)^m$,
where $m$ is the number of triadic places of cubic width five.
The bounds for powers of degrees four through ten are respectively
$16,17,24,31,38,44,50$. Twenty-five tenth powers suffice for
almost all eligible targets. For arithmetic progressions of mixed
exponents we give an explicit number of summands and a necessary and
sufficient finite-place condition.
\end{abstract}

\begingroup
\setlength{\parskip}{0pt}
\setcounter{tocdepth}{1}
\tableofcontents
\endgroup

\section{Introduction}\label{sec:intro}
A congruence obstruction appears even in the formulation of Waring's
problem over algebraic number fields. If $\mathcal O$ is the ring of
integers of $\mathbb K$, write $J_k$ for the additive group generated
by the $k$th powers in $\mathcal O$. Siegel asked for representations
of totally positive integers $\nu\in J_k$ of sufficiently large
absolute norm by totally positive $k$th powers
\cite{hsieg44,hsieg45}. We write $G_{\mathbb K}(k)$ for the least
number of summands, with the size convention \eqref{eq:hnorm}.
Over $\mathbb Q$ this gives the usual $G(k)$. In the field problem
the exponent, the ramification of the field at primes dividing it,
and the choice of archimedean region all enter the circle method.
Only the first two can affect the number of summands in the
main bound below.

Let $H_{\mathbb K}(k)$ be the least number of powers needed at every
finite completion to represent every locally eligible integer with at
least one unit among the roots. For example, a fourth power is $0$ or
$1$ modulo $16$ according as its root is even or odd. In a primitive
representation of $0$ modulo $16$ at least one root is odd, so the
number of odd roots is a positive multiple of $16$. This already
forces $H_{\mathbb Q}(4)\ge16$.

For each fixed proper extension $\mathbb K$ of $\mathbb Q$, the strict
inequality $H_{\mathbb K}(k)<H_{\mathbb Q}(k)$ occurs at infinitely
many exponents, by Corollary~\ref{cor:comparison-fixed-field}. At a
fixed exponent the comparison depends on the completions. For example,
the cubic widths of $\mathbb Q(\sqrt5)$, $\mathbb Q$, and
$\mathbb Q(\sqrt{-3})$ are $3$, $4$, and $5$, respectively, by
Corollary~\ref{nf3:all-local}. Proposition~\ref{prop:comparison-fixed-degree}
gives the completion criterion at a fixed exponent.

The height threshold for global
representability is allowed to depend on $\mathbb K$ and $k$.
Put \(C_1=6.216216302\).
\begin{theorem}\label{thm:intro-main}
There is a cutoff independent of \(\mathbb K\) such that, for \(k\)
above it,
\[
 G_{\mathbb K}(k)\le
       \max\{H_{\mathbb K}(k),\lfloor C_1k-\sqrt{k}\rfloor\}.
\]
For every \(k\ge3\) one may use \(\lfloor6.873k\rfloor\) as the
analytic entry. The effective bounds of
Theorem~\ref{thm:nf-two-coefficient-finite} lower this entry to
\(\lfloor6.830k\rfloor\) for \(k\ge120\) and
\(\lfloor6.258k\rfloor\) for \(k\ge10^{12}\), with intermediate
bounds in \cref{tab:nf-two-coefficient-finite}.
For each fixed \(\mathbb K\),
\[
           \limsup_{k\to\infty}G_{\mathbb K}(k)/k\le C_1.
\]
Over all number fields and all \(k\ge2\) one has
\(G_{\mathbb K}(k)\ll k\log(2k)\) with an absolute implied constant.
\end{theorem}

The analytic assertions follow from
Theorems~\ref{thm:nf-optimized-moment}, \ref{thm:nf-all-linear}
and~\ref{thm:nf-two-coefficient-finite}. For each fixed field,
Theorem~\ref{thm:local-fixed-ramification} gives
\(H_{\mathbb K}(k)\le4k+O_{\mathbb K}(1)\). Over \(\Q\) this
coefficient is attained at every \(k=2^a\), \(a\ge2\): a primitive sum
representing zero modulo \(4k\) requires at least \(4k\) odd roots.
When the field varies with \(k\), Theorem~\ref{thm:hhybrid} gives
\(H_{\mathbb K}(k)\ll k(1+\max_{p\mid k}v_p(k))\). Whether the
logarithm in the resulting field-uniform Waring bound has the best
possible order is open. By Corollary~\ref{cor:nf-log-sharpness}, this
is equivalent to the existence of fields \(\mathbb K_j\) and exponents
\(k_j\to\infty\) with \(H_{\mathbb K_j}(k_j)\gg k_j\log k_j\).

The same limiting coefficient \(C_1\) holds uniformly at sufficiently
large prime-power exponents \(p^a\) with \(p\ge5\), and whenever the
largest ramification index above a prime dividing \(k\) is at most
\(k\exp(-c\sqrt{\log k})\), for fixed \(c>\sqrt{\log2}\).
With bounded ramification indices, every local width above
\((1+\eta)k\) for fixed \(\eta>0\) is bounded by the exact rational
width once \(k\) is sufficiently large (Theorem~\ref{thm:rational-envelope}).
More generally, suppose that in a family of fields the ramification
indices are bounded at each fixed residue cardinality. Then
\(H_{\mathbb K}(k)=o(k\log k)\) uniformly throughout the family.
If the ramification index is bounded by a fixed power of the residue
cardinality times a constant, the local upper bound is
\(4k+1+\log_2 k\) for large \(k\) (Theorem~\ref{thm:local-family-residue}).
Several of our analytic bounds depend on a choice of moment orders
and differencing parameters, and we have not sought the best choice
in every case.

\begin{theorem}\label{thm:intro-small}
For every number field the following summand counts suffice.
\[
\begin{array}{c|rrrrrrrr}
 k&3&4&5&6&7&8&9&10\\\hline
 G_{\mathbb K}(k)\le&7&16&17&24&31&38&44&50
\end{array}
\]
For a totally real field the quartic bound is attained precisely when
a dyadic completion is $\mathbb Q_2$.
Almost all eligible targets are sums of at most $25$ tenth powers.
\end{theorem}
The cubic proof has three parts. Over totally real fields three pairs
of cubes lead to a positive ternary quadratic form. A valuation prescribed
at a fixed auxiliary prime excludes the spinor exceptions. At mixed
signature, reduction by integral rotations controls the complex
coordinates of a representing vector. At totally imaginary fields
three successive approximations to complex cube roots leave a
multiple of six small enough for a four-cube identity.
Four cubes have exact relative density $(23/27)^{m_5(\mathbb K)}$ among eligible targets, where $m_5$ counts the triadic places with local width five (Theorem~\ref{thm:four-cube-exact-density}). The higher exponents use smooth mean values and primitive local
solutions. Theorems~\ref{nf3:realG7}, \ref{thm:nf-small}
and~\ref{thm:nf-quartic-lower} give the details.

The local part also affects density. Over a fixed field the ordinary
and primitive local widths agree for every sufficiently large $k$.
Outside the powers of two, their limiting upper ratio to $k$ is at most
$3/2$. Along powers of two, a residue-degree-one completion can
require $4k$ summands. At large dyadic exponents we determine the
relative density for every summand count between the analytic
density threshold and $4k$. One consequence is that, for $7k/2$
summands, the mean of these densities over real quadratic fields
in discriminant order tends to $57/64$. The order of the limits and
the local splitting calculation are given in
Corollary~\ref{cor:quadratic-density}. For composita of $r$ independent real quadratic fields sampled by their individual discriminants, the proportion with a surviving obstruction is asymptotic to $2^{-r}$ as $r$ increases (Corollary~\ref{cor:multiquadratic-density}). The exact three-cube local densities also show that their product among degree-$d$ fields is least when three splits completely, with value $(7/9)^d$ (Corollary~\ref{cor:three-local-density-degree}).

The history begins with the trace torus of Siegel, the local work of
Stemmler, Birch and Ramanujam \cite{hstemmler61,hbirch61,hram63},
and the estimates of K\"orner, Tatuzawa, Eda and Wang
\cite{hkorn61,htat73,heda67,hwang}.
Diminishing ranges had already given a degree-independent leading
term of order $k\log k$ \cite{hdav98}. The ideal-class transfer and
the minor-arc estimate for smooth algebraic integers in
\cite{hdav99,hdavpenn} give two further points of comparison.
In Section~\ref{sec:ed} we adapt the efficient differencing of
Wooley and of Vaughan and Wooley \cite{edwooley,edvw} to arbitrary
ideal classes. Dividing an element by an ideal moves the variable
to a fractional ideal lattice. Counting the generators that lie in
the original conjugate region costs a power of a logarithm.
This gives the analytic summand count
$k(\log k+\log\log k+2+o(1))$ using the methods available in 1995.

Sections~\ref{sec:nf-local} and~\ref{sec:arith-local} study the
local widths, first uniformly in the field and then at a fixed
completion. A concrete starting point is the calculation for
$\mathbb Z_2[\sqrt2]$ in Section~\ref{sec:sqrt-two-local}. Sections~\ref{sec:hbirch} and~\ref{sec:powerbasis}
relate Birch's formulation to the conjugate regions and give transfers
from rational and intermediate fields. Sections~\ref{sec:nf-small-powers}
and~\ref{nf3:cubes} treat small exponents. For cubes one can begin with
Theorem~\ref{nf3:global}: splitting a target into two sums of four gives
a lower bound for the number of eight-cube representations. The
seven-cube construction of Theorem~\ref{nf3:realG7} follows a different
route. The three-cube density
question occupies Section~\ref{sec:field-cube-outlook}, including the
distinction between a lower bound of order $X^{1-\varepsilon}$ and
positive density. For the linear analytic bound one can turn from
Section~\ref{sec:ed} directly to Sections~\ref{sec:nf-global}--\ref{sec:nf-global-envelope}.
Sections~\ref{sec:nf-global}
and~\ref{sec:nf-fibres} combine the number-field Vinogradov mean
value theorem \cite{NFWooley}, the factors and fibres of Chow and
Dietmann \cite{NFChow}, and arithmetic exponent pairs
\cite{NFMili}. The marked primes are unramified of residue degree one.
Their local phases admit a rational $p$-adic transformation, while the
Fourier transform and all the moments still live on the trace torus
of $\mathbb K$. For mixed powers the moments are taken at separate root
scales and a common target scale. Corollary~\ref{cor:nf-mixed-consecutive}
shows that fewer than $499.805k$ consecutive exponents beginning at $k$
suffice for every sufficiently large totally positive target in every
number field, for all sufficiently large $k$. For a progression of arbitrary step $d$, Theorem~\ref{thm:uniform-mixed-criterion} gives a necessary and sufficient finite-place condition and an explicit number of powers. Theorem~\ref{thm:nf-local-global-envelope} gives the final comparison for each prescribed signature and degree.

Section~\ref{sec:nf-small-powers} proves the small-power estimates,
including the local bounds in \S\ref{subsec:aud-small-local}.
The finite calculations for the local rings and the smooth moments
are specified there and in Appendix~\ref{app:certificates}.
The mean-value and circle-method arguments occur in
Sections~\ref{sec:ed}, \ref{sec:nf-small-powers},
\ref{sec:nf-global} and~\ref{sec:nf-fibres}.

The two-coefficient bound of Wooley \cite{Wnew} extends to
principal ideal moduli by packing on the trace lattice. With prime
multipliers, joint completion forces the frequencies onto a line
over \(\mathbb K\). The estimates in
\cref{sec:joint-conductors,sec:joint-completion,sec:refined-conductors} use that line to
improve the conductor profile from \cite{HexWaring} and to
obtain the same analytic coefficient over the rationals and
over number fields. Further small-power calculations are given in
the companion supplement \cite{NFSmallSupplement}.

\paragraph{Formal verification.}
A full Lean formalization is underway. Lean has verified the uniform
$8k^5$ bound for $k\ge2$, the uniform eight-cube bound with positivity
and norm control, and an earlier degree-uniform analytic theorem under
its primitive local hypothesis, including its unconditional consequence
for odd prime exponents. The earlier classical bounds $G(8)\le38$ and
$G(9)\le46$ have also been verified in Lean.

\paragraph{Acknowledgement.}
Some thirty years ago, as a doctoral student at the University of Michigan,
I had the good fortune to learn about Waring's problem from Robert Vaughan
and my advisor, Trevor Wooley. I remain grateful to them both for their
teaching and guidance, and for the clarity and beauty of their writing
on the subject. In revising this preprint, I hope to bring its arguments
into clearer order and to give readers something of the pleasure I found
in their papers.

\paragraph{Use of generative AI.}
OpenAI's ChatGPT (GPT-6) was used extensively in the
development of this work, including mathematical exploration,
formulation and refinement of conjectures and proofs,
computational experimentation, preparation of numerical
certificates, assistance with Lean formalization, literature
research, and drafting and revision of the manuscript.
The author
assumes full responsibility for the mathematical correctness,
originality, references, and final contents of this article.

\section{The formulation of Siegel and the early bounds}\label{sec:history}
Let $\mathbb K$ be a number field of degree $d=r_1+2r_2$, with ring of
integers $\mathcal O$ and different $\mathfrak d$. We write $N$ for the
absolute norm. The order
\[
 J_k=\mathbb Z[\{x^k:x\in\mathcal O\}]
\]
is also the additive group generated by the $k$th powers, since a product
of powers is a power. It contains $k!\mathcal O$. Indeed, polarizing
the monomial $x^k$ gives
\[
 k!x_1\cdots x_k
 =\sum_{S\subseteq\{1,\dots,k\}}(-1)^{k-|S|}
       \Bigl(\sum_{j\in S}x_j\Bigr)^k,
\]
and we may take $x_2=\cdots=x_k=1$. For example,
$J_2(\mathbb Z[i])=\mathbb Z+2i\mathbb Z$, as one sees by squaring
$a+bi$ and by subtracting $1^2$ and $i^2$ from $(1+i)^2$.

Following \cite{hdav99}, $G_{\mathbb K}(k)$ is the least $s$ for which
every totally positive $\nu\in J_k$ of sufficiently large norm is a sum
of at most $s$ $k$th powers of nonzero totally positive algebraic integers
$x_i$ with
\begin{equation}\label{eq:hnorm}
                   N(x_i)\ll_{\mathbb K,k,s}N(\nu)^{1/k}.
\end{equation}
Total positivity refers to the real embeddings. Deleting zero summands
is equivalent to the usual convention of allowing totally nonnegative
summands. At a complex embedding the squared modulus contributes to the
absolute norm in \eqref{eq:hnorm}. The dependence permitted in the cutoff
$N(\nu)$ and in \eqref{eq:hnorm} will be understood throughout.

\subsection{The trace torus and the Archimedean density}
Put $V=\mathbb K\otimes_{\mathbb Q}\mathbb R=
\mathbb R^{r_1}\times\mathbb C^{r_2}$.
With $e(t)=\exp(2\pi i t)$, the character
$e(\operatorname{Tr}_{V/\mathbb R}(\alpha x))$ is periodic in $\alpha$
modulo $\mathfrak d^{-1}$ for $x\in\mathcal O$. For a bounded region
$B\subset V$ with positive real coordinates, write
\[
 f_B(\alpha;T)=\sum_{x\in\mathcal O\cap TB}
               e(\operatorname{Tr}_{V/\mathbb R}(\alpha x^k)).
\]
Orthogonality on $V/\mathfrak d^{-1}$, with Haar probability measure,
gives the exact number of representations in $(TB)^s$ as
\begin{equation}\label{eq:horth}
 \int_{V/\mathfrak d^{-1}} f_B(\alpha;T)^s
                 e(-\operatorname{Tr}_{V/\mathbb R}(\alpha\nu))\,d\alpha .
\end{equation}
Siegel's dissection \cite{hsieg44,hsieg45} approximates $\alpha$ by
elements $\gamma\in\mathbb K/\mathfrak d^{-1}$. The conductor ideal
$\{a\in\mathcal O:a\gamma\in\mathfrak d^{-1}\}$ records the finite
modulus. Completion modulo this ideal and rescaling $x=Ty$ separate
finite exponential sums from the Fourier transform of the pushforward
of volume by $y\mapsto y^k$. The scale of the main term in
\eqref{eq:horth} is $T^{d(s-k)}$, with factors given by the finite singular
series at $\nu$ and the Archimedean density at $T^{-k}\nu$.
For a finite prime $\mathfrak p$ with residue cardinality $q$, the
local factor in probability normalization is
\[
 \lim_{j\to\infty} q^{-j(s-1)}
 \#\left\{\mathbf x\in(\mathcal O/\mathfrak p^j)^s:
        \sum_i x_i^k\equiv\nu\pmod{\mathfrak p^j}\right\}
\]
whenever the limit exists. The Archimedean factor can be identified
without a multidimensional Dirichlet limit.

\begin{proposition}\label{prop:harch}
Let $k\ge2$, let $B=(0,1)^{r_1}\times\{z:|z|<1\}^{r_2}$, and use
ordinary Lebesgue measure in real and complex coordinates. For $s>k$
the $s$-fold convolution of the power-map pushforward of this measure
has a continuous density on $V$, given by an absolutely convergent
Fourier integral. It is bounded below by a positive constant on each
compact subset of
\[
       (0,s)^{r_1}\times\{z:|z|<s\}^{r_2}.
\]
\end{proposition}
\begin{proof}
At a real place the density on $(0,1)$ is $k^{-1}t^{1/k-1}$.
At a complex place it is $k^{-1}|w|^{2/k-2}$ on the unit disk.
In polar coordinates the substitution $r\mapsto r^k$ contributes
$k^{-1}|w|^{2/k-1}\,d|w|$, and the angle traverses the circle
$k$ times. The product $g$ of these densities belongs to
$L^1(V)\cap L^p(V)$ for $1\le p<k/(k-1)$. Set $p=s/(s-1)$.
Hausdorff--Young gives $\widehat g\in L^s$, so $(\widehat g)^s\in L^1$.
Fourier inversion gives the asserted continuous representative of
$g^{*s}$. The trace pairing changes the ordinary dual coordinates
by an invertible real linear map, with the corresponding constant
Jacobian.

For a target $u$ in the stated open set, choose all $s$ output
variables near $u/s$. This point lies in the interior of the product
of intervals and disks. One may vary the first $s-1$ outputs in fixed
small neighborhoods so that their sum subtracted from $u$ supplies
the last output in the same interior. On these neighborhoods $g$
has a strictly positive essential lower bound, including when a
complex coordinate equals zero. Integrating there gives a lower
bound valid on a neighborhood of $u$. A finite covering gives the
assertion on compact sets.
\end{proof}

Tatuzawa \cite[pp.~263--264]{htat75} withdrew the multidimensional
mean-value step used in \cite{htat73}. The cut point selected in a
one-dimensional second mean-value theorem can depend on the outer
variables. Wang refers to the original argument for the multidimensional
Dirichlet integral in \cite[Lemma~15, p.~138]{audWang1987}. The Fourier integral in
Proposition~\ref{prop:harch} gives the required density directly.
There is a separate normalization at a complex place in
\cite[Theorem~6.2, pp.~72--73]{hwang}, since the map $z\mapsto z^k$ has $k$ inverse branches off the
origin. Thus the factor from all substitutions, in the real and complex
definitions on his p.~73, is $k^{-s(r_1+r_2)}$. The same correction
applies to the explicit Archimedean product in his \S7.1.
It preserves the positivity used in representation arguments.

\subsection{Published bounds and the role of the field degree}
Here $d=[\mathbb K:\mathbb Q]$. The following bounds illustrate the
different ways $d$ has entered the argument. The stated ranges for $k$
are part of the comparisons.
\begin{itemize}
\item Siegel \cite{hsieg44,hsieg45} obtained
$G_{\mathbb K}(k)\le dk(2^{k-1}+d)+1$, as recorded also in
\cite[\S16]{hvw}. Tatuzawa's first treatment \cite{htat58}
gave the sufficient condition $s>8dk(d+k)$.
\item K\"orner \cite{hkorn61,hkorn62} introduced mean values and
Vinogradov's approach into this setting. For large $k$ his bound has
leading term $3dk\log k$. The comparison in \cite[p.~31]{heda75}
gives $dk\{3\log k+3\log((d^2+1)/3)+1\}$.
\item Birch \cite{hbirch61} and Ramanujam \cite{hram63} gave the
degree-independent bound
$G_{\mathbb K}(k)\le\max\{2^k+1,8k^5\}$.
Ramanujam proved the underlying $8k^5$ theorem for the order
generated by powers in a complete discrete valuation ring of
characteristic zero with finite residue field.
\item Tatuzawa \cite[Theorem~5]{htat73} gave the sufficient range
$s>4dk+2dG(k)$. Wang \cite[Theorem~7.3]{hwang} gives
$G_{\mathbb K}(k)\le 4dk+2dG(k)$.
\item Eda \cite{heda67,heda75} adapted Vinogradov's mean values to
number fields. His 1975 main theorem states,
for sufficiently large $k$ and $d>1$, the sufficient range
\begin{align*}
s>{}& 2dk\log k+6dk\log\log k\\
   &+2(2d\log d+d\log\log d+2\log d+14)k+1 .
\end{align*}
We record this as Eda's stated bound. There are independent difficulties
in its derivation. The moment hypothesis of his Lemma~4 (p.~39) contains
$d/(d-1)$, which is absent from Lemma~6 (p.~43) when Lemma~4 is
applied on p.~45. Also the ratio $\theta=(1-s_0/2)/(k-1/2)$
on pp.~43--44 is negative for $s_0\ge3$. The choice $Y_0=X_0^\theta$
on p.~54 then gives $Y_0\to0$ as $X_0\to\infty$, so the conjugate box
eventually contains only zero.
\item The diminishing-ranges argument of \cite{hdav98} gives a leading
coefficient independent of $d$. In the explicit form recorded in
\cite[\S16]{hvw}, it gives
\begin{align}
 G_{\mathbb K}(k)&\le (3+o(1))k\log k+c_d k,\label{eq:hdav98}\\
 c_d&=4d+3\log((d^2+1)/2)+7,\nonumber
\end{align}
where the $o(1)$ is independent of $d$. The rational-to-field argument
in \cite{hdav99} gives
\[
 G_{\mathbb K}(k)\le 2d\bigl(\widehat G(k)+2k\bigr).
\]
Here $\widehat G(k)$ denotes the threshold for positive density among
rational integers, with the convention of \cite{hdav99,hvw}.
\end{itemize}

In his Penn State proceedings article~\cite{hdavpenn}, the author
formulates the summation region by bounding the individual conjugates.
His mean-value estimate for smooth algebraic integers gives a minor-arc
bound with improved dependence on both $k$ and $d$. In the author's
earlier work on efficient differencing, the aim was a coefficient
$1$ multiplying $k\log k$ in class number one. A Vinogradov
mean-value approach suggested a coefficient $2$ for general ideal
classes. One difficulty in the differencing argument already occurs
in class number one. Generators of the same ideal in a bounded
conjugate region can differ by a unit. At height $T$ the logarithmic
unit lattice leaves $O_{\mathbb K}((1+\log T)^{r_1+r_2-1})$
possibilities, a loss absorbed by $T^\varepsilon$. The
diminishing-ranges result \eqref{eq:hdav98} gives an absolute
leading coefficient and an additive multiple of $k$ that grows
with the degree. We next examine how the ideal classes enter the
differencing argument.

\section{Efficient differencing in arbitrary ideal classes}\label{sec:ed}
Let $n=[\mathbb K:\mathbb Q]$, and put $H(x)=\max_{v\mid\infty}|x|_v$,
using the ordinary complex modulus. It is convenient to measure the size
of a box by $X$, so that its sides have length of order $X^{1/n}$.
We adapt the repeated differences of \cite{edwooley,edvw}. Two features
of the argument over $\mathbb K$ are the choice of the lattice after division
by an ideal and the treatment of primes at which a polynomial has a
singular fibre. We give the details since they determine whether the
exponent may be taken independently of $n$.

Fix $k\ge2$. For $0<\eta<1$, let $\mathcal A(X,\eta)$ be the set of
nonzero totally positive $x\in\mathcal O$ with $H(x)\le X^{1/n}$ and
\[
 X^{\eta/2}<N\mathfrak p\le X^\eta
           \qquad(\mathfrak p\mid(x)).
\]
Write $S_s(\mathcal A)$ for the number of solutions in $\mathcal A^{2s}$
of $\sum x_i^k=\sum y_i^k$. Define $\delta_1=k-1$ and
\begin{equation}\label{eq:eddefect}
 \delta_{s+1}=\delta_s-
  \frac{2\delta_s}{k+\delta_s}
  \left(1-\left(\frac{k-\delta_s}{2k}\right)^k\right),
 \qquad \lambda_s=2s-k+\delta_s .
\end{equation}

\begin{theorem}\label{thm:edmoment}
For integers $k\ge2$, $s\ge1$ and $\varepsilon>0$, every sufficiently
small fixed $\eta>0$ gives
\[
 S_s(\mathcal A(X,\eta))\ll X^{\lambda_s+\varepsilon}.
\]
For all sufficiently small fixed $\eta$ there is $L=L(\mathbb K,\eta)$
such that $|\mathcal A(X,\eta)|\gg X(\log X)^{-L}$. In particular, the
$k$th powers of elements of $\mathcal A(X,\eta)$ represent at least
$X^{k-\delta_s-\varepsilon}$ distinct sums of $s$ terms for all sufficiently
large $X$, after decreasing $\eta$ if necessary.
\end{theorem}

\subsection{Division by an ideal}
Choose a finite set $\mathcal L$ of integral ideals representing the narrow
ideal classes, including $\mathcal O$ for the trivial class. For
$0\ne x\in L\in\mathcal L$ put $I_L(x)=(x)L^{-1}$. Given an integral
ideal $D\mid I_L(x)$, choose $L'\in\mathcal L$ and totally positive $a$
with $LD=aL'$. Then
\begin{equation}\label{eq:edclass}
 y=x/a\in L',\qquad I_{L'}(y)=I_L(x)/D,\qquad
 |a|_v\asymp_{\mathbb K,\mathcal L}(ND)^{1/n}\quad(v\mid\infty).
\end{equation}
For the last assertion, multiply any such $a$ by a totally positive
unit. The logarithms of these units form a lattice in
$\sum_{v\mid\infty}n_v t_v=0$, with $n_v=1$ or $2$. Reducing the vector
$\bigl(\log|a|_v-n^{-1}\log |N_{\mathbb K/\mathbb Q}a|\bigr)_v$
in a bounded fundamental domain proves the assertion. Fix one choice
for each $L,D$. For the moment inequality alone one could equally use
ordinary ideal classes, since the sets in the divided lattice may have
any archimedean signs. The narrow choice preserves signs throughout.
This also proves, for every nonzero integral ideal $J$,
\begin{equation}\label{eq:edlattice}
 \#\{z\in J:H(z)\le T^{1/n}\}\ll 1+T/NJ .
\end{equation}
Indeed $J$ is a balanced scalar multiple of a lattice in a fixed finite
list, and in any one such lattice disjoint balls of a fixed radius can be
placed about its points.

There is an alternative description in terms of irreducible factors.
The prime-ideal classes of a narrow principal ideal form a sequence with
sum zero in the narrow class group. Each minimal zero-sum subsequence
has length at most the group order $h^+$, by taking successive partial
sums. Thus an ideal with all prime factors of norm at most $R$ splits
into narrow principal blocks of norm at most $R^{h^+}$. This is the
zero-sum description of atoms, applied to the monoid of totally positive
elements modulo totally positive units; see \cite{edhk84,edghk}. Multiplying
blocks until a given size $M$ is reached costs an overshoot at most
$R^{h^+}$. Formula \eqref{eq:edclass} lets us multiply the prime ideals
themselves, reducing this overshoot to $R$.

Fix a large parameter $W$ and write $Y=W^{\eta/2}$ and $R=W^\eta$.
In a lattice $L\in\mathcal L$, an allowed set of size parameter $U$
will mean any subset of
\[
 \{x\in L\setminus\{0\}:H(x)\le U^{1/n},
       \quad Y<N\mathfrak p\le R\text{ for every }\mathfrak p\mid I_L(x)\}.
\]
Allowing subsets makes positivity, residue conditions and restrictions
on the archimedean shape available at each division. We need the estimate
\begin{equation}\label{eq:edsmall}
 \#\{x\in L\setminus\{0\}:H(x)\ll W^C,\ NI_L(x)\le M\}
             \ll_{C,\mathbb K} M(1+\log W)^{C_1}
                \qquad (1\le M\ll W^{C_2}),
\end{equation}
with an exponent depending only on the field. There are $O_{\mathbb K}(M)$
integral ideals of norm at most $M$, by the classical ideal-counting
estimate \cite{ednark}. For a fixed principal ideal, the logarithms of
its generators in the indicated box differ by unit logarithms in a box
of side $O(1+\log W)$: each logarithmic coordinate has a lower bound of
that order in absolute value by the norm formula and the upper bounds
at the other places. Counting points of the unit lattice and its finite
torsion gives \eqref{eq:edsmall}.

\subsection{An estimate with a polynomial variable}
For a finite frequency set in $\mathbb K$ we use character orthogonality
on the dual of any fractional ideal containing the frequencies. Integrals
below have Haar mass one. We write $e_\alpha(t)$ for the character
evaluated at $t$. Let $\mathcal A$ be allowed at $U$, and let $Z$ be a
subset of $\{z\in\Lambda:H(z)\le T^{1/n}\}$ for some
$\Lambda\in\mathcal L$. Write
\[
 f=\sum_{x\in\mathcal A}e_\alpha(x^k),\qquad
 F=\sum_{z\in Z}e_\alpha(\Phi(z)),\qquad E=\int |F|^2|f|^{2s}.
\]
Here $\Phi\in\mathbb K[t]$. Polynomial coefficients are allowed
numerators and denominators of height $O(W^C)$ for a fixed $C$.
This means that one rational integer of size $O(W^C)$ clears
denominators and that the resulting algebraic integers have all
conjugates $O(W^C)$. All polynomial substitutions below preserve such
a bound with a possibly larger fixed $C$.

In the estimates that follow a factor $W^{O_{k,s}(\eta)+o(1)}$ is
suppressed, the constant in $O_{k,s}(\eta)$ being independent of
$\eta$. Constants and logarithmic factors in the $o(1)$ may depend on
the field, $C$ and fixed $\eta$. The scales satisfy
$1\le U\ll T\ll W^{C_0}$, with fixed constants. When a divided scale
is bounded it can be replaced by $1$ at the expense of a fixed factor.

\begin{lemma}\label{lem:edpoly}
Suppose $s\le\lambda\le2s$ and $S_s(\mathcal A)\ll U^\lambda$, with the
suppressed factor, for every allowed set. Put
\[
 A_0=1,\qquad A_j=\tfrac12+\frac{2s-\lambda}{2k}A_{j-1}.
\]
If $2s-\lambda\le k$, then for $\deg\Phi=j$, $1\le j\le k$, one has
\begin{equation}\label{eq:edpoly}
 E\ll T^{2A_j}U^\lambda .
\end{equation}
For $0\le j\le k$ the number of solutions to
$\Phi(z)=\sum_{i=1}^s(y_i^k-x_i^k)$ is
$\ll T^{A_j}U^\lambda$.
\end{lemma}

\begin{proof}
For constant $\Phi$ the autocorrelation at any fixed shift is at most
its value at zero by Cauchy's inequality, and $|Z|\ll T$. For $j\ge1$
the final assertion follows from \eqref{eq:edpoly} by Cauchy's inequality
and the hypothesis on $S_s$. We prove \eqref{eq:edpoly} by induction
on $j$. Put $M=T^{A_{j-1}/k}$; then $M^k\le T$ since $A_i\le1$.
At most $j$ elements have any given value of $\Phi(z)$, so that
$\int |F|^2\ll T$. If $U\le M$, the trivial bound $|f|\ll U$ gives
$E\ll TU^{2s}\le T M^{2s-\lambda}U^\lambda$. This is
\eqref{eq:edpoly}. For $U\ge M$ the part of $f$ on which
$NI_L(x)\le M$ can be discarded at the same cost. Indeed its supremum
is $\ll M(1+\log W)^{C_1}$ by \eqref{eq:edsmall}, and
$TM^{2s}\le T M^{2s-\lambda}U^\lambda$.
We may also omit the at most $j-1$ critical points $\Phi'(z)=0$, at a
cost $O_j(S_s(\mathcal A))$. Both splittings use
$|a+b|^m\le2^{m-1}(|a|^m+|b|^m)$.

Call a prime of norm in $(Y,R]$ exceptional if it divides one of the
fixed class representatives, a common denominator of $\Phi$, or the
numerator of its leading derivative coefficient. There are $O(1)$ such
primes, since the relevant nonzero norms are $W^{O(1)}$ and $Y$ is a
fixed power of $W$. All other primes in the interval are called ordinary.
The derivative reduces to a nonzero polynomial at each ordinary prime.
We first estimate the part of the solution count $E$ in which no
exceptional prime occurs in the relative ideals of the smooth variables
and no ordinary prime occurring there divides either $\Phi'(z)$ or
$\Phi'(w)$, with $z,w$ the two polynomial variables.

Independently color each ordinary prime in the interval red or blue,
with equal probabilities. For any tuple in this part, the primes from
the smooth variables are disjoint from those dividing the two derivative
values. The union has bounded cardinality, since its members have norm
exceeding $Y$ and all values, after clearing denominators, have norm
$W^{O(1)}$. Thus with probability bounded below all smooth primes in
the tuple are blue and all derivative primes red. By summation over
tuples it suffices to bound, uniformly in a coloring, the full moment
for the corresponding blue smooth set and red polynomial set. This
averaging concerns the nonnegative solution counts.

For each remaining $x$, select a divisor $D$ of $I_L(x)$ by multiplying
its prime ideals in a fixed order until $M\le ND\le MR$. If $M=1$
one may take $D=\mathcal O$. There are $O(MR)$ possible ideals. For
each $D$ write $x=a y$ as in \eqref{eq:edclass}; the resulting set of
$y$ is allowed at $U'\ll U/M$. H\"older's inequality costs
$O((MR)^{2s-1})$ to replace $|f|^{2s}$ by the sum of the individual
$2s$th powers. In each resulting moment fix the value
\[
 b=\Phi(z)+a^k\sum_{i=1}^s y_i^k .
\]
Locally at each $\mathfrak p^v\Vert D$ one has $v_{\mathfrak p}(a)=v$;
hence $\Phi(z)\equiv b\pmod{\mathfrak p^{kv}}$. The permitted roots
modulo $\mathfrak p$ have nonvanishing derivative, so at most $j$ of
them lift, each uniquely, modulo $\mathfrak p^{kv}$. The Chinese
remainder theorem gives at most $j^{\omega(D)}=O(1)$ classes modulo
$D^k$, since all its prime divisors have norm $>Y$. Apply
$|\sum_r c_r|^2\le j^{\omega(D)}\sum_r |c_r|^2$ to the multiplicities
for this fixed $b$. The resulting pairs have $w-z=h\in D^k$.

For $h=0$ their count for each $D$ is $\ll T(U')^\lambda$. For $h\ne0$
the polynomial
\[
               \frac{\Phi(z+h)-\Phi(z)}{a^k}
\]
has degree $j-1$. Its coefficients satisfy the required height bounds by
\eqref{eq:edclass}; denominators of $a$ and $a^{-1}$ have polynomial
size by the ideal identity and the norm bound. The inductive assertion
for a single polynomial variable gives $\ll T^{A_{j-1}}(U')^\lambda$
solutions for each $h$. By \eqref{eq:edlattice} there are
$O(1+T/M^k)=O(T/M^k)$ possible $h$. Summing first over $h$ and then
over $D$ bounds the regular contribution, with the stated suppressed
factors, by
\[
 M^{2s}\left(T+\frac{T}{M^k}T^{A_{j-1}}\right)
                 (U/M)^\lambda
          \ll T^{2A_j}U^\lambda .
\]

It remains to include tuples containing a common prime of the two kinds
just separated, or an exceptional prime in a smooth variable.
For the latter, let $f_p$ be the sum restricted by
$\mathfrak p\mid I_L(x)$ and put $q=N\mathfrak p$. A specified slot
contributes at most
\[
  \int |F|^2 |f_p| |f|^{2s-1}
  \le E^{1-1/(2s)}
       \left(\int |F|^2 |f_p|^{2s}\right)^{1/(2s)} .
\]
If the exceptional tuples dominate, some such prime, of which there
are only $O(1)$, therefore gives
$E\ll E'(T,CU/q)$, using \eqref{eq:edclass} to divide the smooth
variables. The polynomial is divided by the corresponding $a^k$.
Here and below $E'$ denotes another moment in the stated family, with
the same polynomial degree.

For an ordinary prime shared by a specified smooth slot and polynomial
slot put $F_p=\sum_{z:\mathfrak p\mid\Phi'(z)} e_\alpha(\Phi(z))$,
where $z$ ranges in the remaining polynomial set. Define
\[
 B_p=\int |F_p|^2|f|^{2s},\qquad
 C_p=\int |F_p|^2|f_p|^{2s}.
\]
The contribution of these tuples is bounded by the integral of
$|F F_p f_p|\,|f|^{2s-1}$. Two applications of H\"older give the bound
$E^{1/2}B_p^{(s-1)/(2s)}C_p^{1/(2s)}$. By orthogonality
$\sum_p B_p\ll E$: any tuple counted on the left has $\mathfrak p$
dividing a specified nonzero derivative value, and this is possible for
only $O(1)$ primes in $(Y,R]$. If the ordinary collisions dominate,
choose a dyadic interval $Q<q\le2Q$ carrying at least a constant
multiple of $E/(1+\log W)$. There are $O(Q)$ ideals in that interval.
H\"older in the sum over $p$ now gives
\[
 E\ll (1+\log W)^{2s} Q^{s}\sum_{Q<q\le2Q} C_p
 \ll (1+\log W)^{2s} Q^{s+1}\max_{Q<q\le2Q} C_p .
\]
The same calculation applies when $s=1$, with the factor involving
$B_p$ omitted. The roots of $\Phi'$ modulo such a prime lie in at most
$j-1$ residue classes. By Cauchy's inequality we may take both
polynomial variables in one of them. Choose $z_0$ there from the
polynomial set; writing $\Lambda\mathfrak p=b\Lambda'$ as in
\eqref{eq:edclass} gives $z=z_0+bz'$ with
$H(z')\ll (T/q)^{1/n}$. Here $\mathfrak p$ is coprime to $\Lambda$.
Division in the smooth slots gives $x=ay$ with $H(y)\ll (U/q)^{1/n}$.
Consequently
\begin{equation}\label{eq:edcollision}
 E\ll (1+\log W)^{2s}q^{s+1} E'(C T/q,C U/q),\qquad
 \Phi_{\rm new}(z')=a^{-k}\Phi(z_0+bz').
\end{equation}

For clarity, the last two reductions yield a finite descent. At a node
either the regular estimate already bounds the moment, or an exceptional
or ordinary reduction can be chosen. In each reduction the smooth scale
drops by at least $Y/C$, and it was initially at most $W^{C_0}$.
Thus the number of reductions is bounded independently of $W$ for fixed
$\eta$. All transformed coefficients consequently retain polynomial
height in $W$. If the smooth scale is bounded the elementary estimate
$E\ll T U^{2s}$ suffices. The ratio $U/T$ remains bounded along this
descent, and the chosen prime has $q\ll U$ whenever $f_p$ is nonempty.
The factors $C$ in the new boxes therefore cause only fixed constants.
The quantity $T^{2A_j}U^\lambda$ pays the cost $q^{s+1}$ in
\eqref{eq:edcollision}, since $2A_j\ge1$ and $\lambda\ge s$.
It increases on reversing the exceptional reduction as well.
The logarithms accumulated in a bounded number of reductions are
$W^{o(1)}$. The factor $R^{2s}$ from the count of divisor ideals
occurs only at a regular step, where the polynomial degree decreases.
There are at most $k$ such steps in the induction. This proves the
claimed estimate with exponent loss $O_{k,s}(\eta)+o(1)$.
\end{proof}

\begin{proof}[Proof of Theorem~\ref{thm:edmoment}]
For a single smooth variable the moment is at most $k|\mathcal A|$.
Suppose $\lambda_s$ is available for all allowed sets at every scale.
Apply Lemma~\ref{lem:edpoly} with $T=U$, $\Phi(z)=z^k$, and
$Z=\mathcal A$. Then $\lambda_{s+1}=\lambda_s+2A_k$ is admissible,
with the indicated loss. If $a=(k-\delta_s)/(2k)$, the solution of the
recursion is
\[
 A_j=\frac{k}{k+\delta_s}
      +\frac{\delta_s}{k+\delta_s}a^j .
\]
This gives \eqref{eq:eddefect}. The recursion preserves
$0\le\delta_s\le k$, $s\le\lambda_s\le2s$, starting at $\lambda_1=1$:
indeed $2A_k\ge1$, and the multiplier of $\delta_s$ in
\eqref{eq:eddefect} lies in $[0,1]$ for $k\ge2$.
Taking $W=X$ and $\eta$ sufficiently small proves the upper estimate.

Let $e$ be the exponent of the narrow ideal class group. Choose a
multiple $L$ of $e$ with $1/\eta<L<2/\eta$, possible for small fixed
$\eta$. The prime ideal theorem gives $\gg X^{1/L}/\log X$ prime ideals
of norm in $(cX^{1/L},2cX^{1/L}]$ for fixed $c>0$. At least a fixed
proportion belong to one narrow class. Each product of $L$ distinct
primes from this class is narrow principal, and unique factorization of
ideals distinguishes their unordered products. Balance a totally
positive generator of each product as in \eqref{eq:edclass}. Choosing
$c$ small enough places every generator in the required box. All primes
lie in the prescribed interval for large $X$. Finally Cauchy's
inequality bounds the number of distinct sums below by
$|\mathcal A|^{2s}/S_s(\mathcal A)$, proving the theorem.
\end{proof}

\subsection{Size of the auxiliary sumset}
Since $0\le a\le1/2$, equation \eqref{eq:eddefect} gives
\[
 k\log(\delta_s/\delta_{s+1})+\delta_s-\delta_{s+1}
 \ge (k/\delta_s+1)(\delta_s-\delta_{s+1})
 \ge 2(1-2^{-k}).
\]
Summation yields
\[
 \delta_s\le k\exp\!\left(1-
               \frac{2(1-2^{-k})(s-1)}{k}\right).
\]
Thus $s=\tfrac12 k\log k+O(k)$ can give bounded defect, and
$s=\tfrac12 k(\log k+\log\log k)+O(k)$ can give defect
$O(1/\log k)$. The implied constants in these two variable counts
are absolute. Thresholds for the archimedean box may depend on the field.

Using an exponent this close to the natural size of the sumset calls
for a corresponding minor-arc estimate. We establish one on products
of a smooth element and a prime element. The prime element makes it
possible to count the singular residue classes without a class-number
restriction. The use of bilinear sums and the large sieve here is in
the tradition of \cite{edwooleysws}.

\subsection{A bilinear sum on the supplementary domain}
Let $X=P^n$ and let $\mathfrak D$ be the different. Write
$B(t)=\prod_{v\mid\infty}\max(1,|t_v|)^{n_v}$. For a character
$\gamma\in\mathbb K/\mathfrak D^{-1}$ put
$\mathfrak q_\gamma=\{a\in\mathcal O:a\gamma\in\mathfrak D^{-1}\}$.
Define $\mathfrak M$ on the trace torus to consist of the points having
a lift $\alpha=\gamma+\beta$ with
\begin{equation}\label{eq:edmajor}
       (N\mathfrak q_\gamma) B(P^k\beta)\le X^{3/4}.
\end{equation}
This includes all lifts of the centers; take $\mathfrak m$ to be its
complement. Fix $0<\theta<1/2$, put $M=X^\theta$, $U=X/M$, and take
one balanced totally positive generator from each narrow principal
prime ideal of norm in $(cM,2cM]$, where $c>0$ is fixed. There are
$\gg M/\log M$ such ideals by the prime ideal theorem in narrow classes.
Partition the choices into finitely many sectors of angular diameter
less than $\pi/(3k)$ at each complex place, and let $\mathcal P$ be
a cell of largest cardinality. Thus, with constants uniform over the
cells, every $v\in\mathcal P$ has $|v|_w\asymp M^{1/n}$ at each
infinite place $w$, and
\begin{equation}\label{eq:edlip}
 |v^k-v'^k|_w\gg M^{(k-1)/n}|v-v'|_w\qquad(v,v'\in\mathcal P).
\end{equation}
To see this, divide the factor $(v^k-v'^k)/(v-v')$ by $v'^{k-1}$;
the real parts of the terms of the resulting geometric sum are positive
and its constant term is $1$.

Let $\mathcal A_U$ be the totally positive allowed set in $\mathcal O$
at $U$, with $W=X$, and put
\[
               g(\alpha)=\sum_{v\in\mathcal P}
                         \sum_{u\in\mathcal A_U} e_\alpha((vu)^k).
\]
For fixed $\theta$, the counting proof of Theorem~\ref{thm:edmoment}
at $U$ applies with parameter $\eta/(1-\theta)$, giving
$g(0)\gg X(\log X)^{-C}$. The same conclusions hold on replacing the
size bound for $u$ by any fixed positive multiple of $U^{1/n}$.

\begin{proposition}\label{prop:edbilinear}
Given $r\ge1$ with $\theta>(1-\theta)\delta_r$, put
\[
             \rho=\frac{\theta-(1-\theta)\delta_r}{2r}.
\]
For every $\varepsilon>0$ the choice of a sufficiently small fixed
$\eta>0$ gives $\sup_{\mathfrak m}|g|\ll X^{1-\rho+\varepsilon}$.
\end{proposition}

\begin{proof}
We first bound the number of $v\in\mathcal P$ for which $\alpha v^k$
lies in any one torus cell of coordinate side $O(U^{-k/n})$.
Choose $v_0$ in this cell. There are $c_v\in\mathfrak D^{-1}$ such that
\[
 d_v=v^k-v_0^k,\qquad H(\alpha d_v-c_v)\ll U^{-k/n},
\]
where for an archimedean vector $H$ denotes its maximum modulus. All
determinants $d_vc_{v'}-d_{v'}c_v$ belong to $\mathfrak D^{-1}$ and
have absolute norm $O((M/U)^k)=o(1)$. They vanish. Distinct generators
in $\mathcal P$ have distinct $k$th powers, since their ideals are
distinct. Unless the cell has a single member, we obtain
$\gamma\in\mathbb K$ with $c_v=\gamma d_v$ for every $v$ in it.
In particular all the $v^k$ lie in one class modulo
$\mathfrak q=\mathfrak q_\gamma$, with $N\mathfrak q\ll M^k$.
Only $O_k(1)$ of the prime ideals $(v)$ can divide $\mathfrak q$, since
each has norm comparable to $M$.

Among the remaining elements the power congruence occupies
$O_{\mathbb K,k,\varepsilon}((N\mathfrak q)^\varepsilon)$ residue
classes for $v$. At primes away from $k$ this follows from the simple
root bound for $u^k=1$. At a prime above $k$, set
$d=v_{\mathfrak p}(k)$. For a unit root modulo $\mathfrak p^m$ with
$m>2d$, Hensel's lemma gives a root in the completion agreeing with
it modulo $\mathfrak p^{m-d}$. There are at most $k$ such roots in the
field, and at most $(N\mathfrak p)^d$ lifts of each residue to level
$m$. The remaining levels are bounded in terms of the fixed prime and
$k$. The Chinese remainder theorem proves the estimate.

Put $\beta=\alpha-\gamma$. Within any one of these residue classes,
\eqref{eq:edlip} bounds all pairwise differences by
\[
             |v-v'|_w\ll
             M^{1/n}\min\{1,(P^k|\beta_w|)^{-1}\}.
\]
Indeed $H(\beta(v^k-v'^k))\ll U^{-k/n}$ and $UM=X=P^n$.
An ideal coset of norm $Q$ has $O(1+V/Q)$ points in a product of
intervals and discs of volume $V$, uniformly in its aspect ratios.
After dividing each coordinate by its radius, a nonzero difference
has maximum modulus $\gg(Q/V)^{1/n}$ by the field norm, so disjoint
small balls give this bound. Here $V\ll M/B(P^k\beta)$ and
$\alpha\in\mathfrak m$ gives $Q B(P^k\beta)>X^{3/4}$.
This proves that the cell contains $O_\varepsilon(X^\varepsilon)$
prime generators.

Let $b_y$ count the sums of $r$ powers of elements of $\mathcal A_U$.
The trigonometric polynomial $H_0(\xi)=\sum_y b_y e_\xi(y)$ has
coordinate degrees $O(U^{k/n})$ in an integral trace basis.
For a polynomial of degrees $O(L)$ on an $n$-torus, the sum of
$\sup |H_0|^2$ over grid cells of side $\asymp 1/L$ is
$O(L^n\sum_y|b_y|^2)$. Indeed, on each interval the fundamental
theorem of calculus gives the supremum bound by a constant times
$L\int(|H_0|^2+L^{-2}|H'_0|^2)$. Iteration in the coordinates and
Parseval give the assertion, including the mixed derivatives.
H\"older's inequality and the cell count now give
\[
 |g(\alpha)|^{2r}
 \ll M^{2r-1}\sum_{v\in\mathcal P}|H_0(\alpha v^k)|^2
 \ll X^\varepsilon M^{2r-1} U^k S_r(\mathcal A_U)
 \ll X^{2r+\varepsilon}M^{-1}U^{\delta_r}.
\]
Replacing $\varepsilon$ by a smaller number proves the claim.
\end{proof}

\subsection{The resulting Waring bound}
\begin{theorem}\label{thm:edglobal}
Let $k\ge2$ and $n=[\mathbb K:\mathbb Q]$. With any $r,\theta,\rho$ as
in Proposition~\ref{prop:edbilinear} and any $s\ge1$,
\begin{equation}\label{eq:edfinite}
 G_{\mathbb K}(k)\le 4nk+1+2s+\lfloor\delta_s/\rho\rfloor+1 .
\end{equation}
In particular, as $k\longrightarrow\infty$,
\begin{equation}\label{eq:edglobal}
 G_{\mathbb K}(k)\le 4nk+
 k\left(\log k+\log\log k+2+\log2+
   O\!\left(\frac{\log\log k}{\log k}\right)\right)
\end{equation}
with an absolute implied constant. The norm threshold for
representation may depend on $\mathbb K$ and $k$.
\end{theorem}

\begin{proof}
We record the major-arc argument to specify how the restriction on the
prime factors of a summand enters. Write $b=4nk+1$. Use the full Weyl
sum $L_w$ at scale $P$ with a nonnegative smooth product weight $w$
supported in fixed annuli at the complex places and positive intervals
away from zero at the real places. A fixed finite set of such weights
will suffice. The major-arc kernel for $\ell\in J_k$ is
\begin{equation}\label{eq:edkernel}
 \int_{\mathfrak M}L_w(\alpha)^b e_\alpha(-\ell)\,d\alpha
 = c_{\mathbb K,b}X^{b-k}
       \mathfrak S_b(\ell)J_w(\ell/P^k)+O(X^{b-k-\tau})
\end{equation}
for some $\tau>0$, with the positive factor $c_{\mathbb K,b}$
determined by the measures. Here $J_w$ is the $b$-fold convolution of
the power-map pushforward of $w$ evaluated at the normalized target.

For details, distinct lifted major arcs in \eqref{eq:edmajor} are
disjoint for large $X$: their distinct centers have norm difference
$\gg 1/(Q_1Q_2)$, where $Q_i=N\mathfrak q_i$, whereas a point of
intersection would give
$X^k/(Q_1Q_2)\ll B(P^k\beta_1)B(P^k\beta_2)\le
 X^{3/2}/(Q_1Q_2)$. The lower bound follows from
$(\gamma_1-\gamma_2)\mathfrak q_1\mathfrak q_2\subset\mathfrak D^{-1}$.
At $\gamma+\beta$ apply Poisson summation on each coset of
$\mathfrak q=\mathfrak q_\gamma$. The zero dual mode is
$c_{\mathbb K}X\,C(\gamma)I_w(P^k\beta)$, with
\[
 C(\gamma)=(N\mathfrak q)^{-1}
              \sum_{y\bmod\mathfrak q}e_\gamma(y^k),
 \qquad I_w(t)=\int w(y)e_t(y^k)\,dy.
\]
The nonzero modes together contribute $O_A(X^{-A})$ for every fixed
$A$. In fact a nonzero $\xi\in\mathfrak q^{-1}\mathfrak D^{-1}$ has
absolute norm $\gg 1/N\mathfrak q$. Relative to the radii
$(1+P^k|\beta_v|)/P$ it therefore has maximum size
$L\gg (X/(N\mathfrak q B(P^k\beta)))^{1/n}\gg X^{1/(4n)}$.
Integration by parts at a place attaining this maximum bounds the
oscillatory integral, with linear twist $-P\xi$, by $O_N(L^{-N})$.
The product weight permits the other places to be bounded trivially.
The ideal-coset packing bound above counts
$O(1+L^n N\mathfrak q B(P^k\beta)/X)$ frequencies in the associated
box. Summing dyadic shells proves the error estimate.

Classically $|C(\gamma)|\ll_\varepsilon
(N\mathfrak q)^{-1/k+\varepsilon}$ \cite[pp.~38--39]{htat73}.
Also $|I_w(t)|\ll B(t)^{-1/2}$. At a real place this is the
second-derivative test on the support. At a complex place the Hessian
of $\operatorname{Re}(t z^k)$ is nonsingular there, giving the
two-dimensional estimate $O((1+|t|)^{-1})$; integration by parts
away from the origin also suffices for this untwisted integral.
There are at most $Q$ characters of any given conductor of norm $Q$.
Since $b>2k$, the product of the ideal sum and the archimedean
integral still converges absolutely with an extra factor
$(QB(t))^{\tau_0}$ for some $\tau_0>0$. Extending the ranges in
\eqref{eq:edmajor} therefore gives \eqref{eq:edkernel}, uniformly
in $\ell$, by Fourier inversion. For $k\ge3$, \cite[Theorem~2]{htat73} gives
$\mathfrak S_b(\ell)\gg_{\mathbb K,k}1$ for $\ell\in J_k$ in this
range of $b$. For $k=2$ the same assertion follows from the elementary
local argument below. At a completion over $p$, the $\mathbb Z_p$-module
spanned by squares has a basis of at most $n$ squares: choose a basis
from the generating set modulo $p$ and use Nakayama's lemma.
Every $c\in\mathbb Z_p$ is a sum of four $p$-adic integer squares by
Lagrange's theorem and compactness. Expanding $\ell-1$ in the chosen
basis and appending $1$ gives a primitive solution in at most $4n+1$
variables. The local-density argument following this proof applies.

Multiply the target by a $k$th power of a totally positive unit, as
in the balancing argument \eqref{eq:edclass}. With
$P=|N_{\mathbb K/\mathbb Q}\nu|^{1/(nk)}$, this puts $\nu/P^k$
in a fixed compact set bounded away from zero at every place, and
positive at the real places. For a given point of this compact set,
choose the weight $w$ positive near a root of $y_v^k=(\nu/P^k)_v/b$
at each place, taking the positive root at a real place. The power map has invertible derivative there. Thus the
convolution $J_w$ has a positive lower bound in a neighborhood of
this target, by integrating over small neighborhoods of $b-1$ roots.
A finite subcover yields weights and uniform neighborhoods for all
balanced targets.

Let $f$ be the smooth sum in Theorem~\ref{thm:edmoment}. Take
$t=\lfloor\delta_s/\rho\rfloor+1$ and integrate
$L_w^b f^{2s}g^t e_\alpha(-\nu)$. Use sufficiently small fixed
dilations of the smooth and bilinear root boxes, so the sum of all
their powers, when divided by $P^k$, lies in a ball small enough for
the preceding lower bound on the residual target. This changes only
fixed constants in the estimates. Expanding $f^{2s}g^t$ first, each
residual $\ell$ lies in $J_k$ and hence each major-arc kernel is
$\gg X^{b-k}$. Their sum is
$\gg X^{b+2s+t-k}(\log X)^{-C}$ by the lower bounds for the two
masses. On the supplementary domain, trivial estimation of $L_w$,
the unrestricted moment of $f$, and Proposition~\ref{prop:edbilinear}
give an upper bound
\[
                 \ll X^{b+2s+t-k+\delta_s-t\rho+\varepsilon}.
\]
The strict inequality defining $t$ proves positivity for large $X$.
Every summand has the required signs and height $O(P)$ before the unit
normalization is reversed. Its absolute norm is unchanged under
that last step. This proves \eqref{eq:edfinite}.

For large $k$ take $\theta=1/2-1/\log k$ and
$r=\lceil k(\log k+\log\log k+1)/2\rceil$. The defect estimate above
gives $\delta_r=O(1/\log k)$ and thus
$\rho^{-1}=2k(\log k+O(\log\log k))$. Taking
$s=\lceil k(\log k+\log\log k+1+\log2)/2\rceil$
gives $\delta_s\le (1+O(1/k))/(2\log k)$ and
$t\le k(1+O(\log\log k/\log k))$. Equation~\eqref{eq:edglobal}
follows.
\end{proof}

There are two useful variants of the last argument. First, historical
local inputs permit the replacement
\[
 4nk+1\quad\hbox{by}\quad\min\{4nk+1,8k^5+1\}
\]
in \eqref{eq:edfinite}. Indeed Ramanujam \cite{hram63} gives a local
representation of $\ell-1$ at every finite place with at most $8k^5$
summands, so that appending $1$ gives a solution with a unit coordinate.
More generally any uniform primitive local width $h$ permits
$b=\max(2k+1,h)$. At the finitely many small places Hensel's lemma in
the unit coordinate gives a target-independent lower density: if
$e=v_{\mathfrak p}(k)$, one can fix the other $b-1$ coordinates
modulo $\mathfrak p^{2e+1}$ and lift arbitrarily. At all larger places
the absolute complete-sum majorant gives an Euler factor
$1+O((N\mathfrak p)^{-1-\tau})$. These bounds furnish the uniform
positive singular series used in \eqref{eq:edkernel}.

Second, one can dispense with the bilinear sum if a quadratic estimate
in $k$ is sufficient. Choose $b=4nk+1$ and use Tatuzawa's ordinary
minor-arc estimate \cite[Lemma~2 and the proof of Theorem~5]{htat73},
$L\ll P^{n-1/2^k}$, with its associated basic domains. Let $V$ sum one
character for each distinct auxiliary $s$-term sum at a small fixed
dilation of the root box. The minor contribution of $L^bV^2$ is
$O(P^{nb-b/2^k}V(0))$ by orthogonality. On the basic domains the main
term is $\gg P^{n(b-k)}V(0)^2$, by expansion over both auxiliary sums
and the uniform shifted-target estimate of that theorem, using the
archimedean density of Proposition~\ref{prop:harch}. Consequently
$G_{\mathbb K}(k)\le b+2s$ for every $s$ with
$\delta_s<b/(n2^k)$. In particular one may take
\[
 s=2+\left\lceil\frac{k}{2(1-2^{-k})}
           \{1+\max(0,\log(kn2^k/b))\}\right\rceil .
\]
This gives $4nk+(\log2)k^2+O(k)$ and identifies the saving due to the
bilinear estimate in \eqref{eq:edglobal}.

\subsection{Separating the local and analytic widths}

A positive-density version of the preceding construction permits the
local width to enter through a maximum. Write $h$ for any integer such
that every element of the local additive span of $k$th powers has a
representation by $h$ such powers with a unit coordinate, at every
finite place. For the finitely many ramified residue characteristics
the unit coordinate is understood in the local ring. Extra zero
coordinates can always be appended.

\begin{theorem}\label{thm:edmaximum}
Fix $r,\theta,\rho$ as in Proposition~\ref{prop:edbilinear}. For
integers $d\ge2$, $a,c\ge1$ satisfying
\[
 S=2d(a+c)>2k,\qquad 2(d-1)c\rho>\delta_{(d-1)a},
\]
one has
\begin{equation}\label{eq:edmaxfinite}
                  G_{\mathbb K}(k)\le\max\{h,S\}.
\end{equation}
In particular the historical local estimates give
\begin{equation}\label{eq:edmaxhistorical}
 G_{\mathbb K}(k)\le\max\!\left\{
 \min(4nk+1,8k^5+1),\;
 k\left(\log k+\log\log k+2+\log2+
      O\!\left(\frac{\log\log k}{\log k}\right)\right)\right\}.
\end{equation}
The implied constant is absolute. The first entry in this maximum
can also be replaced by $h$. In particular if $h\le Ck$ on a family
of fields and exponents, the second entry alone suffices for all
large $k$, uniformly in that family for each fixed $C$.
\end{theorem}

\begin{proof}
We first construct a smooth weight on totally positive roots at
volume $Y$ with mass comparable to $Y$. For a small fixed $\eta>0$
choose an integer $m$ with $\eta/2<1/m<\eta$. Then choose
$\eta/2<\kappa_0<\kappa_1<1/m$ so close to $1/m$ that
$1-m\kappa_0<\eta$. Put $\kappa_* = 1-m\kappa_1>0$.
Choose one balanced positive generator for each narrow principal
prime ideal. Form products $D$ of $m$ such generators, with prime
norms in $(Y^{\kappa_0},Y^{\kappa_1}]$. The roots have the form
$Dz$, where $z\in\mathcal O$ carries weight
$\psi(z/Z^{1/n})$, $Z=Y/|ND|$. Here $N$ on an element denotes
its field norm, and $\psi$ is smooth, nonnegative, nonzero and
supported in a compact box away from zero at each place, positive
at the real places. In particular
\[
             Y^{\kappa_*}\le Z\le Y^\eta .
\]
The $m$-tuples may be ordered. Their contribution to any root has
bounded multiplicity, since it has norm $O(Y)$ and each of these
primes has norm greater than $Y^{\kappa_0}$. Divide all weights
by a fixed constant so that the coefficients on individual roots
are at most one. Denote the resulting sum by $f_Y$. Counting $z$
in its dilated box gives mass comparable to $Z$, uniformly in $D$,
and hence
\[
  f_Y(0)\asymp Y\!
  \left(\sum_{\substack{Y^{\kappa_0}<N\mathfrak p\le Y^{\kappa_1}\\
              [\mathfrak p]^+=1}}\frac1{N\mathfrak p}\right)^m
  \asymp Y.
\]
The last estimate follows by partial summation from the prime ideal
theorem in narrow classes. The roots have modulus comparable to
$Y^{1/n}$ at every place, since each generator used in $D$ is
balanced and their number is fixed.

For each fixed $j$ the defect for $f_Y$ is still $\delta_j$. Indeed,
there are $O(Y^\eta)$ possible outer summands when the roots are
grouped by $z$. For a fixed $z$, the inner sum over $D$ has bounded
coefficients, satisfies the prime restriction in
Theorem~\ref{thm:edmoment} at $Y$, and has height $O(Y^{1/n})$.
Multiplication by $z^k$ preserves normalized Haar measure on the
dual torus, since its action on the dual lattice is injective.
The triangle inequality in $L^{2j}$ and the moment estimate give
\begin{equation}\label{eq:eddensemoment}
           \int |f_Y|^{2j}\ll Y^{2j-k+\delta_j+\varepsilon}
\end{equation}
after choosing $\eta$ sufficiently small. The outer triangle inequality
costs $Y^{2j\eta}$ in the moment, so take $2j\eta<\varepsilon/2$
and apply Theorem~\ref{thm:edmoment} with loss $\varepsilon/2$.
Only the inner $D$-sum has the two-sided prime restriction.
All these assertions hold
under a fixed small dilation of the root boxes.

For the bilinear sum, fix $\theta<\theta_1<1/2$. Take dyadic
intervals $(M,2M]$ contained in $[X^\theta,X^{\theta_1}]$ and put
$U=X/M$ for each interval. Let
\[
  g=\sum_M\ \sum_v f_U(\alpha v^k),
\]
where $v$ runs through the balanced positive generators of the narrow
principal primes in this interval. The prime ideal theorem gives
$g(0)\asymp X$, by summing the reciprocal prime norms in this
exponent window. Each individual root contributes bounded weight.
In fact the marked large prime is uniquely determined by the root:
all prime factors of $D$, and all prime factors of $z$ by its norm
bound, have size $O(X^\eta)<X^\theta$. We decrease $\eta$ to ensure
this inequality. Divide $g$ by a fixed constant if necessary.
Partition each interval into the finitely many complex angular
sectors used in \eqref{eq:edlip}. The spacing proof of
Proposition~\ref{prop:edbilinear} applies to each interval and cell,
now using \eqref{eq:eddensemoment} for its inner sum. The constants
are uniform for $X^\theta\le M\le X^{\theta_1}$, since the determinant
bound uses $(M/U)^k=o(1)$ uniformly there. Summing $O(\log X)$ intervals
yields
\begin{equation}\label{eq:eddenseminor}
                    \sup_{\mathfrak m}|g|\ll X^{1-\rho+\varepsilon}.
\end{equation}

Set $f=f_X$, $s=da$, $u=dc$ and $b=S/d$. Both weighted root sets lie
in a fixed compact product of annuli at scale $P=X^{1/n}$.
Take a nonnegative smooth product weight $w$ on this region that
majorizes their coefficients, and let $L=L_w$. Its support can be
enlarged within a compact set separated from zero at every place.
Set
\[
                       A=\int |f|^{2s}|g|^{2u}.
\]
In the orthogonality expansion all coefficients are nonnegative.
Replacing $a$ roots from each side in the $f$ variables and $c$
from each side in the $g$ variables by the majorizing ordinary
weight gives
\[
       A\le\int |L|^b |f|^{2(s-a)} |g|^{2(u-c)} .
\]
On $\mathfrak m$ the last integral is
$O(X^{S-k+\delta_{s-a}-2(u-c)\rho+\varepsilon})$ by
\eqref{eq:eddensemoment} and \eqref{eq:eddenseminor}.
On $\mathfrak M$, H\"older bounds it by
$(\int_{\mathfrak M}|L|^S)^{1/d}A^{1-1/d}$.
For $S>2k$ the absolute complete sum and archimedean integral
estimates in the proof of \eqref{eq:edkernel} give
$\int_{\mathfrak M}|L|^S\ll X^{S-k}$. It follows that
\begin{equation}\label{eq:edlossless}
                            A\ll X^{S-k}.
\end{equation}
Write $\mathfrak M(D_0)$ for the restricted major arcs
$Q B(P^k\beta)\le D_0$. For some $\tau>0$,
\[
      \int_{\mathfrak M\setminus\mathfrak M(D_0)}|L|^S
                    \ll X^{S-k}D_0^{-\tau}
\]
for $1\le D_0\le X^{3/4}$. The same H\"older inequality and
\eqref{eq:edlossless} give a bound
$O(X^{S-k}D_0^{-\tau/d})$ for the mixed integrand on this shell.

Take $v_0=\max(h-S,0)$ and $s_0=S+v_0$. Integrate
$L^{b+v_0} f^{2(s-a)}g^{2(u-c)}e_\alpha(-\nu)$. The minor and
shell estimates just proved acquire only a factor $O(X^{v_0})$.
We give the details on $\mathfrak M(D_0)$ since the weights there
must now contribute to the local density. At a fixed lifted center
$\gamma$ with conductor norm $Q\le D_0$, group a restricted sum by
its marked primes and the product $D$. The remaining summation over
$z$ has scale $Z^{1/n}\gg X^{(1-\theta_1)\kappa_*/n}$.
The coefficient multiplying $z$ is a unit modulo the conductor for
large $X$. Lattice summation on its residue classes therefore supplies
the same factor $C(\gamma)$ as for $L$. More precisely each such sum,
after division by $Z$, is its archimedean integral times
$c_{\mathbb K}C(\gamma)$ with an error $o(1)$ uniform in the marked
primes and $P^k\beta$ for fixed $D_0$. To verify the uniformity,
scale $z$ by $Z^{1/n}$: the normalized coefficient of $z$ in the
root, divided by $P$, is in a compact set bounded away from zero
at each place, by balance. The smooth weights and phases have
uniformly bounded derivatives on their common compact support.
Riemann summation in each of the fixed finitely many residue lattices
gives the assertion. Summing the errors costs $o(X)$, since the sum
of their scale factors $Z$ is $O(X)$ for both $f$ and $g$.
Thus $X^{-1}$ times either weighted sum on this fixed core has
arithmetic factor $C(\gamma)$ and an archimedean factor which is
the Fourier transform of a positive measure of bounded mass.
The latter measure is obtained by pushing the boxes forward by
the normalized power map. Its mass is bounded below uniformly in
$X$, by the two estimates at zero.

Use small fixed dilations for all restricted roots. Enlarge $w$ so
that, at each place, it is positive in a neighborhood of the roots
of $(\nu/P^k)_v/(b+v_0)$ chosen as in the proof of
\eqref{eq:edkernel}. A finite collection suffices after unit
balancing. The archimedean convolution of the $b+v_0$ ordinary
weights has a positive lower bound on every target left after the
restricted roots have been subtracted. Convolution with their
positive measures preserves a uniform lower bound. The absolutely
convergent singular series has $s_0\ge h$ and $s_0>2k$ and hence
is bounded below as above. The full product of this series and
convolution differs from its truncation at $QB(t)\le D_0$ by
$O(D_0^{-\tau'})$, uniformly in the restricted measures, for some
$\tau'>0$: use $|C(\gamma)|^{s_0}$ arithmetically and the bounded
masses times $|I_w(t)|^{b+v_0}$ at infinity. Here $b=2(a+c)\ge4$.
First choose $D_0$ large to absorb both tails, then let $X$ tend
to infinity in the fixed-core estimate. The resulting representation
count is $\gg X^{s_0-k}$. This proves \eqref{eq:edmaxfinite}.

For \eqref{eq:edmaxhistorical}, use the preceding choices of
$\theta,r$ and take $d=\lceil\sqrt{k}\rceil$.
Round the previous $s$ up to a multiple $(d-1)a$ and its auxiliary
number of bilinear roots up to an even multiple $2(d-1)c$.
The integer rounding and the replacement of $d-1$ by $d$ raise the
total by $O(\sqrt{k}\log k)$, which is absorbed by the stated error.
The local choices $4nk+1$ and $8k^5+1$ were proved above. This
establishes the theorem.
\end{proof}

\subsection{An asymmetric pair of smooth factors}

Wooley's asymmetric use of two smooth mean values
\cite[Sections~3--4]{edwooleysws} improves the absolute term of order
$k$. We give the lattice argument for passing between the two
estimates in arbitrary signature. In this subsection replace the
exponent $3/4$ in \eqref{eq:edmajor} by $\omega=\omega(k)<1$.
For each fixed $k$ choose the three window endpoints
$\lambda_-<\lambda_+<\omega$ before letting $X$ tend to infinity.
Disjointness and the zero-mode approximation in the proof of
\eqref{eq:edkernel} still hold: $X^{2\omega}=o(X^k)$, and the
scaled shortest nonzero dual vector there has length
$\gg X^{(1-\omega)/n}$.

\begin{theorem}\label{thm:edasymmetric}
Let $h$ be a primitive local width as in Theorem~\ref{thm:edmaximum}.
As $k$ tends to infinity,
\begin{equation}\label{eq:edasymptotic}
  G_{\mathbb K}(k)\le\max\!\left\{h,\,
    k\left(\log k+\log\log k+2+
             O\!\left(\frac{\log\log k}{\log k}\right)\right)\right\},
\end{equation}
with an absolute implied constant. In particular
$h=\min(4nk+1,8k^5+1)$ is available from the historical local
estimates.
\end{theorem}

\begin{proof}
Take $1/2<\lambda_-<\lambda_+<\omega<1$, with
$\lambda_+/\lambda_-$ close to one. We construct a weight of mass
$\asymp X$ from products $vu$. Choose a large integer $m$ and form
$v$ as an ordered product of the balanced positive generators of
$m$ narrow principal primes, each of norm in
$[X^{\lambda_-/m},X^{\lambda_+/m}]$. Split the products into
dyadic ranges $M<|Nv|\le2M$ and complex angular cells as in
\eqref{eq:edlip}. For each range use the independent inner weight
$f_U$ of \eqref{eq:eddensemoment}, where $U=X/M$. Write $g$ for the
sum of $e_\alpha((vu)^k)$ with these weights. Any one root has bounded
multiplicity, since all outer prime divisors have norm at least
$X^{\lambda_-/m}$, and their generators have been fixed.
After a constant normalization the coefficients are at most one.
Summing $U\asymp X/|Nv|$ over the $m$ marked primes gives $g(0)\asymp X$
by partial summation, exactly as before.

Any subset of the outer products in one dyadic range has an
admissible $2w$-moment $\ll M^{2w-k+\delta_w+\varepsilon}$.
Indeed choose $1/m$ strictly between $\eta/2$ and $\eta$, and
choose the exponent window close enough that
\[
 \frac{\lambda_-}{m\lambda_+}>\eta/2,
 \qquad \frac{\lambda_+}{m\lambda_-}<\eta .
\]
Constant factors at the ends of a dyadic interval are harmless by
strictness. Apply Theorem~\ref{thm:edmoment} at $M$ and its
bounded-coefficient version, decreasing $\eta$ for the required
moments. We work with a single dyadic range and angular cell below.
Let $\lambda=\log M/\log X$. An error $O(1/\log X)$ at the two
endpoints can be absorbed in the constants.

For fixed sufficiently large $C$, consider the pairs
\begin{equation}\label{eq:edjointlattice}
 (d,c)\in\mathcal O\times\mathfrak D^{-1},\qquad
     H(d)\le CM^{k/n},\quad H(\alpha d-c)\le CU^{-k/n}.
\end{equation}
Two possibilities arise. If these pairs contain two linearly
independent vectors over $\mathbb K$, their number is
$O(M^k/U^k)$. To see this, send $(d,c)$ to
$(M^{-k/n}d,U^{k/n}(\alpha d-c))$. The resulting real lattice
has covolume comparable to $(U/M)^k$. Each of the two independent
bounded vectors, multiplied by the members of a fixed integral
basis, gives $n$ bounded vectors of the image lattice. The $2n$
resulting vectors are real-linearly independent: before applying
the displayed $\mathbb K_\infty$-linear map, their independence follows
by scalar extension from the invertible two-by-two matrix over
$\mathbb K$. They span a sublattice with bounded covering radius.
The full image lattice consequently has bounded covering radius as
well. Its Voronoi cells about the points in the fixed box have
disjoint interiors, are all contained in a bounded enlargement,
and have volume equal to its covolume. This proves the count.
The same count holds if \eqref{eq:edjointlattice} contains only the
zero vector.

Here is the consequence for the bilinear sum in this case. Write
$h(\xi)=f_U(\xi)$ and let $b_y$ be the coefficients of $h^t$ for an
integer $t\ge1$. Choose complex numbers $a_v$ of modulus at most one
to absorb the arguments of $h(\alpha v^k)^t$ and the bounded
multiplicities of the outer weight (enlarging a fixed constant if
necessary). With $V_y=\sum_v a_v e_\alpha(yv^k)$, two applications
of H\"older give
\[
 |g_M|^{2tw}\ll M^{2w(t-1)}
       \left(\sum b_y\right)^{2w-2}
       \left(\sum b_y^2\right)
       \sum_{H(y)\ll U^{k/n}} |V_y|^{2w}.
\]
The first application writes the $t$th power of the sum in terms of
$\sum b_y V_y$ after taking absolute values at the outer roots.
For the second, apply H\"older to
$b_y^{1-1/w}\,b_y^{1/w}\,|V_y|$.
To bound the last sum, majorize its box by a nonnegative Schwartz
function $\Psi(y/U^{k/n})$ with Fourier transform of compact support.
One may use a constant multiple of the squared modulus of the Fourier
transform of a smooth bump of sufficiently small support. The
coefficient at frequency $d$ in $|V_y|^{2w}$ has absolute value at
most the untwisted $2w$-moment of the outer roots, by Cauchy's
inequality on each translate of the $w$-fold representation function.
Poisson summation over $y\in\mathcal O$ bounds the last sum by this
moment times $O(U^k)$ times the number of pairs in
\eqref{eq:edjointlattice}, with $C$ fixed sufficiently large.
As $\sum b_y\ll U^t$ and $\sum b_y^2\ll U^{2t-k+\delta_t+\varepsilon}$,
the resulting estimate is
\begin{equation}\label{eq:edlargeasym}
 |g_M|\ll X^{1+\varepsilon}
     \left(M^{\delta_w}U^{-k+\delta_t}\right)^{1/(2tw)}.
\end{equation}

It remains to consider rank one. Since $U$ tends to infinity, a
nonzero pair has $d\ne0$, and all the pairs satisfy $c=\gamma d$
for one $\gamma\in\mathbb K$. Put $\mathfrak q=\mathfrak q_\gamma$
and $\beta=\alpha-\gamma$, so that $N\mathfrak q\ll M^k$. For $v$
with $\alpha v^k$ in one torus cell of side $O(U^{-k/n})$, any pair
of their powers has difference $d$ as in \eqref{eq:edjointlattice}.
Thus these powers agree modulo $\mathfrak q$, and
$H(\beta(v^k-v'^k))\ll U^{-k/n}$ for all pairs in the cell.
Separate the outer prime factors of $v$ which divide $\mathfrak q$.
Their product of chosen generators is $d_0$, say $D=|Nd_0|$.
There are $O_{k,m,\lambda_-}(1)$ possibilities for $d_0$: the number
of prime divisors of $\mathfrak q$ of norm at least
$X^{\lambda_-/m}$ is bounded, and only $m$ marked factors are used.
For each possibility write $v=d_0 v'$. There are $O(M/D)$ choices
for $v'$, by balance and lattice counting; $D\ll M$. Moreover $v'$
is coprime to $\mathfrak q$. Its $k$th power in any one torus cell
lies in a single class modulo
\[
    \mathfrak q'=\mathfrak q/\gcd(\mathfrak q,(d_0)^k),
              \qquad N\mathfrak q'\ge N\mathfrak q/D^k .
\]
The unit-root bound in Proposition~\ref{prop:edbilinear} gives
$O(X^\varepsilon)$ classes for $v'$ modulo $\mathfrak q'$.
The sector inequality and the bound for
$\beta(v^k-v'^k)$, applied first to two original outer generators,
place their quotients by $d_0$ in a translated box of volume
$O(M/(D B(P^k\beta)))$. Hence the occupancy of any frequency
cell in the given $d_0$ group is
\[
 \ll X^\varepsilon\!\left(1+
   \frac{M D^{k-1}}{N\mathfrak q\,B(P^k\beta)}\right)
 \ll X^\varepsilon D^{k-1}
         \qquad(\alpha\in\mathfrak m).
\]
In the last estimate $N\mathfrak q\, B(P^k\beta)>X^\omega$ and
$M\ll X^{\lambda_+}$. Apply the box large sieve in
Proposition~\ref{prop:edbilinear} with moment $r$.
H\"older now has prefactor $(M/D)^{2r-1}$, so the extra factor in
the $2r$th power is $D^{k-2r}\le1$ as soon as $2r\ge k$.
This proves
\begin{equation}\label{eq:edsmallasym}
             |g_M|\ll X^{1+\varepsilon}
                     (M^{-1}U^{\delta_r})^{1/(2r)} .
\end{equation}
Summing over the intervals and cells, one obtains
$\sup_{\mathfrak m}|g|\ll X^{1-\rho+\varepsilon}$ for every
$\rho$ at most the minimum, over the chosen exponent window, of
\[
 \frac{(1-\lambda)(k-\delta_t)-\lambda\delta_w}{2tw},
       \qquad \frac{\lambda-(1-\lambda)\delta_r}{2r}.
\]
At a fixed conductor all the $m$ outer marked primes are units;
$u$ has the same freely summed final variable as before. Thus the
fixed-core and shell arguments of Theorem~\ref{thm:edmaximum}
apply without change.

For completeness write $L=\log k$ and $l=\log L$, and take
$r=\lceil kL/2\rceil$, $t=\lceil kl/L\rceil$,
$w=\lceil kl\rceil$. Choose a sufficiently narrow exponent window
inside $[1-6l/L,1-5l/L]$. The defect estimate gives
$\delta_r=O(1)$ and $\delta_w=O(k/L^2)$. Also $k-\delta_t\ge t/2$
for large $k$: since $\lambda_j\ge j$, one has $\delta_j\ge k-j$
and in particular $\delta_j\ge k/2$ for $j\le t$; each decrement
in \eqref{eq:eddefect} is then at least $1/2$.
Consequently the first displayed saving is at least
$(5/4-o(1))/(kL)$, while the second is at least
$(1-O(l/L))/(kL)$. We may therefore take
$\rho^{-1}=k(\log k+O(\log\log k))$. In
\eqref{eq:edmaxfinite}, round a smooth half-moment of size
$k(\log k+\log\log k+1)/2$ and its required number of bilinear
variables up as before. Its defect is at most
$(1+O(1/k))/\log k$. The bilinear width is at most
$k(1+O(l/L))$ and rounding adds $O(\sqrt{k}\log k)$,
as required.
\end{proof}

\section{The local width at exponents of bounded prime multiplicity}
\label{sec:nf-local}

For the integer ring $A$ of a finite extension of $\mathbb Q_p$, write
$J_k(A)=\sum_{x\in A}\mathbb Z x^k$. This is a ring, since products of
powers are powers. Let $w_k(A)$ be the least number of powers needed for
every element of $J_k(A)$, and $h_k(A)$ the corresponding width when
at least one of the roots must be a unit, even at target zero. Zero
summands are allowed. Put
\[
 H_{\mathbb K}(k)=\sup_{v\nmid\infty}h_k(\mathcal O_{\mathbb K_v}),
 \qquad W_{\mathbb K}(k)=\sup_{v\nmid\infty}w_k(\mathcal O_{\mathbb K_v}).
\]
For each completion one has $w_k(A)\le h_k(A)\le w_k(A)+1$:
represent $\nu-1$ and append $1$. Moreover
$W_{\mathbb K}(k)\le G_{\mathbb K}(k)$. Indeed, modulo any power of
a finite prime each local generator lifts to a global generator, and
so does every element of their additive span. Adding large positive
rational integers divisible by the modulus gives totally positive
global targets of arbitrarily large norm. Reduction of their
representations gives solutions at each level. Compactness gives
a representation in the completion.

Put $A(k)=\max_{p\mid k}v_p(k)$ for $k\ge2$. The following asymptotic
is uniform both in the field and in its degree.
\begin{theorem}\label{thm:nf-bounded-multiplicity}
For every fixed positive integer $A_0$ one has
\[
 \limsup_{\substack{k\to\infty\\A(k)\le A_0}}
 \frac1k\sup_{\mathbb K}H_{\mathbb K}(k)\le1.
\]
For every $\varepsilon>0$, the positive integers $k$ for which
$\sup_{\mathbb K}H_{\mathbb K}(k)>(1+\varepsilon)k$ form a set of
natural density zero.
\end{theorem}

\subsection{Residue-field coefficients}
Suppose $k=p^a m$, $(m,p)=1$, and $A$ has residue field
$F=\mathbb F_q$. Set $H=(F^\times)^m$, $E=\mathbb F_p(H)$,
$t=[F:E]$, and $d=(|E|-1)/|H|$. Denote by $L$ the least number of
members of $H\cup\{0\}$ whose sums fill $E$. The additive span of $H$
is a field, since it contains $1$ and is closed under multiplication.
The successive sumsets are invariant under $H$ and grow by at least
one nonzero multiplicative coset until they fill $E$; equality of
consecutive sumsets implies invariance under the additive span of $H$.
Thus $L\le d$. Also $td\le m$, by $q-1\ge t(|E|-1)$ and
$(q-1)/|H|\le m$. A basis of $E$ chosen in $H$ gives
$L\le(p-1)[E:\mathbb F_p]$.

We will use both estimates
\begin{equation}\label{eq:nf-residue-width}
 tL\le m,\qquad tL\ll_p\log(2m).
\end{equation}
For $q\le 16(m+1)^4$, the basis bound proves the second. For larger
$q$, a proper subfield containing $H$ would give
$(q-1)/m\le\sqrt q-1$, a contradiction. For $b\in F^\times$ write
$d_0=(m,q-1)$. If $d_0=1$, a single term suffices. Otherwise the
smooth projective Fermat curve $X^{d_0}+Y^{d_0}=bZ^{d_0}$
has at least $q+1-(d_0-1)(d_0-2)\sqrt q$ rational points by
Hasse--Weil \cite{NFWeilCurves}. At most $d_0$ of them are at infinity.
There is therefore an affine solution. Since the images of the
$m$th and $d_0$th power maps agree, $L\le2$. This proves
\eqref{eq:nf-residue-width}. When $a=0$, reserving one unit before
representing the residue and applying simple-root Hensel lifting gives
$h_k(A)\le k+1$. Ramanujam \cite[Lemma~1]{hram63} traces the local
representation theorem for $p\nmid k$ to Theorem~7 of Stemmler's
report \cite{hstemmler61} and to Birch \cite[Lemma~5]{hbirch61}.

\subsection{A uniformizer calculation}
We give the local calculation in some detail since ramification can be
arbitrarily large. Suppose $a\ge1$, set $r=a+1$, let $\delta=1$ for
odd $p$ and $\delta=2$ for $p=2$, and put $B=W(F)$ for the unramified
coefficient ring. The $p$-typical Witt vectors used here are indexed
from $0$. We use their usual ghost polynomials, the Teichm\"uller vector
$[u]=(u,0,\ldots,0)$, and the Frobenius and Verschiebung identities
in characteristic $p$ \cite{NFHazewinkel}.

Let $R=(B/p^\delta)[[T]]$, $D=\sum_{j=0}^{a-1}(a-j)p^j$,
and take $N=2\delta D+1$. Write $I$ for the kernel of
$W_r(R)\longrightarrow W_r(R/T^N)$ and $G$ for the additive group
generated by $[y^m]$, $y\in R$. Evaluation at a uniformizer $\pi$
followed by the last ghost polynomial gives a ring map
\[
 W_r(R)\longrightarrow A/p^{a+\delta}A,\qquad
 (z_j)\longmapsto\sum_{j=0}^a p^j\widetilde z_j^{\,p^{a-j}} .
\]
Here $\widetilde z_j$ is any lift to $A$ of $z_j(\pi)$ modulo $p^\delta$.
Lifts give the same result by the binomial theorem:
$(u+p^\delta v)^{p^b}\equiv u^{p^b}\pmod {p^{b+\delta}}$.
The ghost identities prove multiplicativity and additivity.
Evaluation $R\to A/p^\delta A$ is onto, since $A=B[\pi]$.
The image of $G$ is consequently the reduction of $J_k(A)$.

Choose $b_0=1$, $b_i=1+T^i$ for $1\le i\le a$.
We claim that every element of $I$ is of the form
\begin{equation}\label{eq:nf-chart}
 \sum_{i=0}^a[b_i+z_i]-\sum_{i=0}^a[b_i],\qquad z_i\in TR .
\end{equation}
These expressions belong to $G$ since each element of $1+TR$ has
an $m$th root there. The centered Witt sum in \eqref{eq:nf-chart}
is an integral polynomial map in the $z_i$ with linear matrix
$J=(b_i^{p^j-1})_{0\le j,i\le a}$. To see the latter, differentiate
the ghost identities over the universal torsion-free coefficient ring
at $\mathbf z=0$; lower-coordinate powers in the $j$th output ghost
have zero linear term. Modulo $p$, multiply column $i$ by $b_i$,
then subtract column $0$ from every other column. The resulting
matrix has entries $1$ in its first column and $T^{ip^j}$ elsewhere.
Its Vandermonde determinant in the numbers $T^{p^j}$ has $T$-adic
order $D$. Thus $J$ has a right inverse over $R[T^{-1}]$ such that
$C=T^{\delta D}J^{-1}$ has entries in $R$. Indeed,
$\det J=T^D u+p v$ for a unit $u\in R$; expand its inverse in the
finite geometric series in $pv/(T^D u)$. Substitute
$\mathbf z=T^{\delta D+1}C\mathbf x$. The linear term equals
$T^N\mathbf x$ and all higher terms are multiples of $T^{N+1}$.
After division by $T^N$ the map is $\mathbf x+T\mathbf P(\mathbf x)$
with integral coefficients. Successive approximation in the complete
ring $R$ solves the equation for every right-hand side. This proves
\eqref{eq:nf-chart}.

Put $M=W_r(R/T^N)$. It is a module over $B$ via the constant unramified
ring. Explicitly, if $g$ generates $F^\times$, the double
Teichm\"uller lift $[\widehat g]$ satisfies $X^{q-1}-1=0$ and
has residue $g$. In the finite local Witt ring, Hensel factorization
of $X^{q-1}-1$ over $\mathbb Z_p$ makes all factors except the
unramified factor of $g$ units at $[\widehat g]$. This supplies the
map from $B$. The same construction for $H$ gives the action of
$W(E)$ on the image of $G$ in $M$. The case $F=\mathbb F_2$
uses the constant integer ring. Write $\ell$ for length over $B$.
For $p$ odd the formula
\[
 p^s(z_0,\ldots,z_a)
  =(0,\ldots,0,z_0^{p^s},\ldots,z_{a-s}^{p^s})
\]
in characteristic $p$ yields
$\ell(p^sM)=(a+1-s)\lceil N/p^s\rceil$ for $0\le s\le a$.
It follows by counting the free coefficients of the surviving
$T$-polynomials. Put
\[
 U(p,a)=rN+
 (p-1)\sum_{s=1}^{a}p^{s-1}(r-s)\left\lceil N/p^s\right\rceil
 \quad(p>2).
\]
For $p=2$ put $M_1=W_r(F[[T]]/(T^N))$ and define
\[
 U(2,a)=rN+2\left(rN+\sum_{s=1}^a
 2^{s-1}(r-s)\left\lceil N/2^s\right\rceil\right).
\]
Reduction $M\to M_1$ is surjective with kernel of $B$-length $rN$.
The kernel has coordinatewise addition since products of two coefficients
in $2B/4B$ vanish in all Witt carries. It is thus killed by $2$,
and multiplication by $2$ induces a surjection $M_1\to2M$.

We now show that
\begin{equation}\label{eq:nf-uniform-witt}
 h_k(A)\le tL\,U(p,a)+a+1.
\end{equation}
Build a filtration of the image of $G$ by adjoining actual generators
$[y^m]$, as a module over $W(E)$. Each cyclic successive quotient has
length $e_i$. The number of permitted constant powers needed to
represent a coefficient modulo $p^{e_i}$ is at most
$L(1+p+\cdots+p^{e_i-1})$: represent each residue digit as a sum in
$H$, lift by Teichm\"uller constants, and use $p^j$ repetitions at
digit $j$. A permitted constant times $[y^m]$ is again an allowed
power. Lift a target backwards through the filtration. For any $s$
the sum of $(e_i-s)_+$ is bounded by the length of $p^sM$ over
$W(E)$, since $p^sM_i/p^sM_{i-1}$ maps onto the $p^s$-multiple
of the $i$th quotient. Restricting scalars multiplies $B$-length by
$t$. The identity
\[
 1+p+\cdots+p^{e-1}
  = e+(p-1)\sum_{s\ge1}p^{s-1}(e-s)_+
\]
now gives the stated cost $tL\,U(p,a)$ in odd characteristic.
At $2$, the cost of $M$ is at most $rN$ plus twice the cost of $M_1$:
insert $\ell(M)=rN+\ell(M_1)$ and
$\ell(2^sM)\le\ell(2^{s-1}M_1)$ for $s\ge1$ in the same
cost identity. The characteristic-two coordinate formula for $M_1$
gives the displayed $U(2,a)$. Apply this construction after
subtracting the fixed baseline $\sum_i[b_i]$ from the target. The
remaining correction in $I$ is \eqref{eq:nf-chart}. Its $r$ new
terms have unit roots on evaluation. Finally multiplication by
$k$ on the logarithms of $1+p^\delta A$ has image
$p^{a+\delta}A$, so that $(1+p^\delta A)^k=1+p^{a+\delta}A$.
Moving one unit root absorbs the remaining error. For this lifting
step Ramanujam \cite[Lemma~2]{hram63} points to Theorem~12 of
Stemmler's 1959 report \cite{hstemmler61}. This proves
\eqref{eq:nf-uniform-witt}.

\begin{proof}[Proof of Theorem~\ref{thm:nf-bounded-multiplicity}]
Fix $A_0$. For odd $p$ tending to infinity and $1\le a\le A_0$,
we have $N=O_{A_0}(p^{a-1})$ and
\[
 U(p,a)=p^a+O_{A_0}(p^{a-1}).
\]
Indeed for all sufficiently large $p$ the term $s=a$ is
$(p-1)p^{a-1}$. For $s<a$, use $\lceil N/p^s\rceil\le1+N/p^s$
to bound the remaining terms by $O_{A_0}(p^{a-1})$, including
$rN$. The first inequality of \eqref{eq:nf-residue-width} in
\eqref{eq:nf-uniform-witt} therefore gives
$h_k(A)\le(1+O_{A_0}(p^{-1}))k+O_{A_0}(1)$. On the other hand,
for $p$ in any fixed finite set and $a\le A_0$, the second
inequality of \eqref{eq:nf-residue-width} gives
$h_k(A)\ll_{p,A_0}\log(2m)$. First choose the finite set of small
primes and then let $k$ grow. Combined with $h_k(A)\le k+1$
at good primes this proves the first assertion, with field-independent
threshold for the local inequality.

For a fixed $A_0$ and fixed $\varepsilon$, all sufficiently large
exceptional exponents must be divisible by $p^{A_0+1}$ for some
prime $p$. Their upper density is at most
$\sum_p p^{-A_0-1}$, by counting multiples directly up to $X$.
Let $A_0\to\infty$.
\end{proof}

The primitive condition has a visible cost at dyadic places. In
$\mathbb Z_2$ fourth powers modulo $16$ are $0$ and $1$, and a unit
contributes $1$. A primitive representation of zero needs at least
sixteen terms. Conversely reserve one odd fourth power, match the
residue of the remaining target modulo $16$ with at most fifteen
copies of $1$, and use $(1+4\mathbb Z_2)^4=1+16\mathbb Z_2$ to
lift. Thus $h_4(\mathbb Z_2)=16$. More generally the same argument
gives $h_{2^a}(\mathbb Z_2)=2^{a+2}$ for $a\ge2$. Every odd unit
is a sign times a member of $1+4\mathbb Z_2$, so its $2^a$th power
is $1$ modulo $2^{a+2}$, and even roots give zero there. This sequence
also shows why the restriction on prime multiplicities in
Theorem~\ref{thm:nf-bounded-multiplicity} matters.

\subsection{Exact dyadic widths in residue degree one}
\label{sec:nf-exact-dyadic}

Let $A$ be a fixed dyadic integer ring with residue field $\mathbb F_2$,
and let $F$ be its fraction field. Set
\[
        \Lambda=\log(A^\times)\subset F,\qquad
        \Lambda\cap\mathbb Q_2=2^b\mathbb Z_2.
\]
Here $b$ is an integer depending on $F$. The logarithm is defined
by the convergent series on $1+4A$ and then by
$\log u=2^{-j}\log(u^{2^j})$ for $j$ sufficiently large.
A common $j$ works for all units, since
$A^\times/(1+4A)$ is a finite $2$-group. This definition makes
$\Lambda$ a compact additive subgroup, hence a $\mathbb Z_2$-module.
It contains $4A$, since the logarithm maps $1+4A$ onto $4A$.
In particular $b$ is well-defined and $b\le2$.

\begin{theorem}\label{thm:nf-exact-dyadic}
For all sufficiently large $a$ depending on $A$, and every odd $m\ge1$,
\[
             w_{2^am}(A)=h_{2^am}(A)=2^{a+b}.
\]
When $[F:\mathbb Q_2]$ is odd, $b=2$. Thus any number field with
a dyadic completion of odd degree and residue field $\mathbb F_2$
has $G_{\mathbb K}(2^a)\ge4\cdot2^a$ for all sufficiently large
$a$, in every signature.
\end{theorem}
\begin{proof}
The odd power map is an automorphism of $A^\times$, as follows
by passage to its finite $2$-group quotients. For $a$ sufficiently
large one has the equality of sets
\begin{equation}\label{eq:nf-linear-log-image}
             (A^\times)^{2^a}=1+2^a\Lambda .
\end{equation}
Indeed, $u^{2^a}=\exp(2^a\log u)$. By compactness, fix $j_0$ with
$\Lambda\subseteq2^{-j_0}A$. For $\ell\in\Lambda$ the terms beyond
the first in $(\exp(2^a\ell)-1)/2^a$ belong to $4A$ uniformly when
$a$ is large. This follows term by term using $v_2(n!)\le n-1$.
Hence the quotient is in $\Lambda$. Conversely
$2^{-a}\log(1+2^a\ell)\equiv\ell\pmod{4A}$ uniformly, using the
same bound for $v_2(n)$. This element of $\Lambda$ is the logarithm
of a unit, which proves the reverse inclusion. Nonunit $2^am$th
powers belong to $2^a\Lambda$, since $\mathfrak p_A^{2^a}\subseteq
2^{a+2}A$ for large $a$.

It follows from \eqref{eq:nf-linear-log-image} that
$J_{2^am}(A)=\mathbb Z+2^a\Lambda$. Modulo $2^a\Lambda$ this cyclic
group has order $M=2^{a+b}$. Every class has a representative
$1\le r\le M$, and $r-1$ copies of $1$ and one arbitrary unit
power fill its entire coset. Thus $h_{2^am}(A)\le M$. The target
$M$ belongs to the zero coset. Its valuation is
$e(a+b)<2^a$, with $e=v_{\mathfrak p_A}(2)$, so it cannot be a sum
of nonunit powers. The number of unit summands must be a positive
multiple of $M$. Hence $w_{2^am}(A)\ge M$.

Finally $\operatorname{Tr}_{F/\mathbb Q_2}(\log u)
=\log N_{F/\mathbb Q_2}u\in4\mathbb Z_2$. One obtains the identity
by multiplying the conjugate principal units on a common domain
of convergence. An element of $\Lambda\cap\mathbb Q_2$ therefore
belongs to $4\mathbb Z_2$ when the degree is odd. The opposite
inclusion holds already over $\mathbb Q_2$. The global assertion
now follows from $W_{\mathbb K}(k)\le G_{\mathbb K}(k)$.
\end{proof}

For comparison, suppose the residue field of a fixed dyadic $A$
has at least four elements. Then
\begin{equation}\label{eq:nf-higher-residue-dyadic}
                  h_{2^a}(A)\le2^{a+1}+O_A(1).
\end{equation}
Write $B$ for its unramified coefficient ring, $\mu$ for the Teichm\"uller
units, and put $\Lambda=\log(1+\mathfrak p_A)$.
The proof of \eqref{eq:nf-linear-log-image} now gives the unit-power
set $\bigcup_{\zeta\in\mu}\zeta(1+2^a\Lambda)$. The module
$\Lambda'=\sum_{\zeta\in\mu}\zeta\Lambda$ is a $B$-module, since
the Teichm\"uller representatives span $B$ over $\mathbb Z_2$.
Choose $t\in B$ whose reduction is outside $\mathbb F_2$. Then
\[
    \tfrac12\log(1+2t)\equiv t-t^2\pmod {2B}
\]
is a unit of $B$. Hence $2B\subseteq\Lambda'$. One fixed summand
for each $\zeta\in\mu$ absorbs any error in $2^a\Lambda'$.
Subtract their fixed baseline. The rest of an eligible target
is represented modulo $2^a\Lambda'$ by an element of $B$, since
nonunit powers belong to $2^{a+2}A\subseteq2^a\Lambda'$ for large $a$.
Its Teichm\"uller expansion modulo $2^{a+1}B$ uses at most
$1+2+\cdots+2^a$ copies of elements of $\mu$. Each of these is itself
a $2^a$th power. This proves \eqref{eq:nf-higher-residue-dyadic}.

There is a finite cyclotomic test for the case $b=2$ of
Theorem~\ref{thm:nf-exact-dyadic}. Suppose again that $F$ has
residue field $\mathbb F_2$ and let its roots of unity have order
$2^h$. Then
\begin{equation}\label{eq:nf-kummer-log}
  b=2-\max\{t\ge0:5\zeta\in(F^\times)^{2^t}
                  \text{ for some }\zeta\in\mu_{2^h}(F)\}.
\end{equation}
In fact $v_2(\log5)=2$, and the kernel of the logarithm on
$\mathcal O_F^\times$ is its torsion subgroup: raise an element
to a common power placing it in $1+4A$, where the logarithm is
injective. An element with logarithm $2^{-t}\log5$ exists precisely
under the condition in \eqref{eq:nf-kummer-log}. The trace identity
in the proof above bounds the maximum by $v_2([F:\mathbb Q_2])$.
For $t=1$ the two square classes of torsion units are represented
by $1$ and a primitive $2^h$th root $\zeta_0$.
The extension $F(\sqrt5)/F$ is unramified quadratic, since
$(1+\sqrt5)/2$ has irreducible reduction $X^2+X+1$ over the
residue field. In particular $5$ is not a square in $F$.
It follows that
\[
 b<2\quad\Longleftrightarrow\quad
 F(\mu_{2^{h+1}})/F\text{ is unramified}.
\]
Indeed $5\zeta_0$ is a square precisely when
$F(\sqrt5)=F(\sqrt{\zeta_0})=F(\mu_{2^{h+1}})$.

For example the dyadic completions of $\mathbb Q(\sqrt2)$ and
$\mathbb Q(\sqrt3)$ both have ramification index two and residue
degree one, and both have only the torsion units $\{1,-1\}$.
In the first case adjoining $i$ gives the totally ramified
cyclotomic extension $\mathbb Q_2(\mu_8)$, so $b=2$.
In the second case $-5/3$ is a square in $\mathbb Q_2$, as it
is congruent to $1$ modulo eight. Hence $F(i)=F(\sqrt5)$ is
unramified, and the degree-two trace bound gives $b=1$.
For sufficiently large $k=2^a$ the respective dyadic widths are
therefore exactly $4k$ and $2k$.

\subsection{How often the local width exceeds the degree of the power}
\label{sec:local-sparse}

The estimates above also locate the exceptional exponents. The implied
constants in the following theorem are absolute apart from their dependence
on $\varepsilon$.
\begin{theorem}\label{thm:local-sparse}
Fix $\varepsilon>0$. If $k$ is sufficiently large and
$\sup_{\mathbb K}H_{\mathbb K}(k)>(1+\varepsilon)k$, then
$k=mp^a$ for a prime $p$, with $a\ge1$, $(m,p)=1$, and
\[
       p\ll_\varepsilon(a+1)^2,\qquad
       m\ll_\varepsilon (a+1)^2\log(2a+2).
\]
In particular, the number of such exponents $k\le X$ is
$O_\varepsilon((\log(2X))^5\log\log(3X))$.
For every $0<\varepsilon<3$ there are $\gg\log X$ exceptions up to $X$.
\end{theorem}

\begin{proof}
At a place of residue characteristic $p$ dividing $k$, let the notation be
as in \eqref{eq:nf-uniform-witt}. Two uniform estimates suffice. First,
the proof of \eqref{eq:nf-residue-width} gives
\[
 tL\le\min\left\{m,\,C\frac{p\log(2m)}{\log p}\right\}
\]
with an absolute $C$. In the small-residue-field case use
$tL\le(p-1)\log(q)/\log p$ and $q\le16(m+1)^4$.
In the other case $tL\le2$.

Next put $r=a+1$. The expression for $D$ in the uniformizer calculation
gives $D/p^a\le\sum_{i\ge1}i/p^i=p/(p-1)^2$. In particular
$N/p^a\ll1/p$, also at $p=2$. Using
$\lceil N/p^s\rceil\le N/p^s+1$ in $U(p,a)$, the contribution of
the $N/p^s$ terms, together with $rN$, is $O(r^2 p^{a-1})$.
The contribution of the ones, divided by $p^a$, is bounded by
\[
 \frac{p-1}{p}\sum_{i\ge0}\frac{i+1}{p^i}
 =\frac p{p-1}.
\]
This gives $U(p,a)/p^a\le1+C_1r^2/p$ at odd $p$. The extra terms in the
definition of $U(2,a)$ give $U(2,a)/2^a\ll r^2$. Increasing $C_1$
therefore gives the same inequality for all $p$.

Suppose now that $h_k(\mathcal O_{\mathbb K_v})>(1+\varepsilon)k$.
The bound $k+1$ excludes $p\nmid k$ when $k>1/\varepsilon$.
At bad primes, $r/k\le(1+\log_2 k)/k<\varepsilon/2$ for all
sufficiently large $k$. Applying \eqref{eq:nf-uniform-witt} first
with $tL\le m$ gives $p\ll_\varepsilon r^2$. Applying it again, now
with the second bound on $tL$, gives
\[
         m\ll (p+r^2)\frac{\log(2m)}{\log p}
             \ll_\varepsilon r^2\log(2m).
\]
The function $u/\log(2u)$ is increasing for $u\ge2$. Substitution
of a sufficiently large multiple of $r^2\log(2r)$ gives the asserted
bound for $m$. As $a\le\log_2 X$, counting the triples $(a,p,m)$ gives
the stated upper bound, even on counting possible $p$ by all integers.
For the lower bound take $k=2^a$ with $a\ge2$. The primitive dyadic
width over $\mathbb Q_2$ is $4k$, as proved above.
\end{proof}

\subsection{A uniform $k\log k$ bound for the local width}
\label{sec:local-hybrid}
The upper and lower Witt coordinates require different choices of
centres. Using unit centres in a short block and powers of the
uniformizer in the other block removes a logarithmic factor from
the uniform local estimate.

\begin{theorem}\label{thm:hhybrid}
There is an absolute constant $C$ such that, for every number field
$\mathbb K$ and every integer $k\ge2$,
\[
               H_{\mathbb K}(k)\le Ck\bigl(1+\max_{p\mid k}v_p(k)\bigr).
\]
In particular $H_{\mathbb K}(k)\ll k\log(2k)$ uniformly in
$\mathbb K$ and $k$.
\end{theorem}

Retain the notation of \eqref{eq:nf-uniform-witt}. In particular
$k=mp^a$, $(m,p)=1$, $r=a+1$, $R=(B/p^\delta B)[[T]]$, and the
residue coefficient cost satisfies
\begin{equation}\label{eq:hcoeff}
          tL\le m,\qquad tL\ll p\log(2m)/\log p .
\end{equation}
For positive integers $N_0,\ldots,N_a$ satisfying
$pN_j\ge N_{j+1}$ let $I_N$ be the set of Witt vectors whose
$j$th coordinate lies in $T^{N_j}R$. It is an ideal. Indeed the
$j$th Witt addition and negation polynomials have weight $p^j$ when
coordinate $i$ has weight $p^i$. Each multiplication polynomial
has that weight separately in the two arguments. A monomial in
variables from $I_N$ of weight $p^j$ has $T$-order at least $N_j$
because $N_i/p^i\ge N_j/p^j$ whenever $i\le j$. For multiplication use
the weights in the factor from $I_N$ alone.

We record the coefficient cost for this ideal. Write $N_\Sigma=\sum_jN_j$.
If every element of $I_N$ is a sum of $r$ permitted variable
Teichm\"uller powers minus a fixed sum of $r$ permitted powers, with
one variable root a unit, then
\begin{equation}\label{eq:hhycost}
                  h_k(A)\ll r+tL\bigl(p^a+rN_\Sigma\bigr)
\end{equation}
with an absolute constant. In characteristic $p$ the lengths of
the successive $p^s$-multiple modules of $W_r(R)/I_N$ are
\[
             L_s=\sum_{j=s}^a\lceil N_j/p^s\rceil .
\]
Indeed $p^s(z_0,\ldots,z_a)$ has surviving coordinates
$z_{j-s}^{p^s}$ for $s\le j\le a$. Counting the kernel in any
sufficiently large common truncation gives this formula, including
$s=0$. In the dyadic case reduce the coefficient ring first modulo
two. The kernel in the weighted Witt quotient has length $N_\Sigma$
and is killed by two. To see this, represent its classes by vectors
with all coordinates in $2R$, by subtracting a lift from $I_N$.
All nonlinear addition carries between such vectors vanish. The
coordinate lengths of the kernel add to $N_\Sigma$, by triangularity.
Consequently the positive coefficient cost before multiplication by
$tL$ is at most three times
\[
             L_0+(p-1)\sum_{s=1}^a p^{s-1}L_s.
\]
For odd $p$ the factor three can be replaced by one.
Indeed the filtration by actual power generators, the lifting of the
residue digits in the coefficient ring, and baseline cancellation
work exactly as in \eqref{eq:nf-uniform-witt}. The terms $N_j/p^s$
in the ceilings contribute at most $aN_\Sigma$, and the ceiling
errors contribute at most
$(p-1)\sum_{s=1}^a p^{s-1}(a-s+1)\ll p^a$.
This proves \eqref{eq:hhycost}.

\begin{proof}[Proof of Theorem~\ref{thm:hhybrid}]
It suffices to consider $p\mid k$, so $a\ge1$. We first treat
$m\ge64$. Put $c=1+\lfloor\log_p(m/64)\rfloor$. If $c\ge r$,
then $m\ge64p^a$. In this case the isotropic estimate in the proof of
Theorem~\ref{thm:local-sparse} already gives
\[
     h_k(A)/k\ll 1+(a+1)^2\log(2m)/(m\log p)\ll1.
\]
For the last bound use the decrease of $\log(2m)/m$ and
$m\ge64p^a$. Suppose now $1\le c<r$ and put $b=r-c$.

For the first $b$ Witt coordinates use $b$ nonunit centres
$T^{mi}$, $1\le i\le b$, and perturb the $i$th centre by
$T^{n_i}z_i$, where
\[
   n_i=mi+1+m\frac{p^{b-i}-1}{p-1}.
\]
Each perturbed centre is $T^{mi}$ times a member of $1+TR$ and
therefore has an $m$th root in $R$. For $0\le j<b$ set
\[
 N_j=m(b-j)p^j+1+m\frac{p^j-1}{p-1},\qquad K_0=N_{b-1}.
\]
The least of $mi(p^j-1)+n_i$ over $1\le i\le b$ is $N_j$.
It is attained at $i=b-j$ and also at $i=b-j-1$ when the latter
is positive. Indeed the successive differences of these orders are
$m(p^j-p^{b-i-1})$. The matrix of the normalized linear terms thus
has determinant $\pm1$ modulo $T$, by reversing the columns. Also
$N_j$ is nondecreasing, and $pN_j>N_{j+1}$ for $j+1<b$.
For all $j\ge b$ the minimum order from the nonunit linear terms
is at least
\begin{equation}\label{eq:hhyleak}
                    K_0+m(p-1)p^{b-1}.
\end{equation}

For the remaining $c$ coordinates use the unit centres $1$ and
$1+T^i$, $1\le i<c$. Let $J$ be the matrix $(u_i^{p^j-1})$
of these unit columns in the rows $b\le j\le a$. Modulo $p$,
multiply each column by $u_i$ and subtract the first column from
each subsequent column. The resulting matrix has first column $1$
and other entries $T^{ip^j}$. Its Vandermonde determinant gives
\[
             D=\operatorname{ord}_T(\det J\bmod p)
                      =\sum_{j=b}^{a-1}(a-j)p^j\le4p^{a-1}.
\]
For $c=1$ this sum is zero. The matrix $C_0=T^{\delta D}J^{-1}$
is integral, by the finite geometric expansion modulo $p^\delta$
used in \eqref{eq:nf-chart}. Perturb all unit centres by the
coordinates of $T^{K_0+1}C_0\mathbf w$. They remain $m$th
powers of units since their constant terms are $1$.
Set $N_j=N_*=K_0+1+\delta D$ for $b\le j\le a$.
The choice $m\ge64p^{c-1}$ gives
\[
 \delta D\le8p^{a-1}\le m p^{b-1}/8\le K_0/8 .
\]
It follows that \eqref{eq:hhyleak} exceeds $N_*$,
that $pK_0>N_*$, and that $2(K_0+1)>N_*$.

We verify that the centred sum of all $r$ Teichm\"uller vectors
takes every value of $I_N$. Its $j$th coordinate, before reduction
in the coefficient ring, is an integral polynomial $H_j$ in the
perturbations. Its linear terms are $\sum_i u_i^{p^j-1}\eta_i$,
using $u_i$ for all centres and $\eta_i$ for their perturbations.
In the $j$th ghost equation a nonunit term of ordinary degree $h\ge2$
from $(T^{mi}+\eta_i)^{p^j}$ has order at least
$mi(p^j-1)+n_i+(h-1)(n_i-mi)$. This exceeds $N_j$ in every row.
All unit terms of degree $h\ge2$ have order at least $2(K_0+1)$,
which also exceeds every $N_j$. The lower-coordinate terms in that
ghost equation contain $H_l^{p^{j-l}}$ for $l<j$. By induction
their orders are at least $p^{j-l}N_l>N_j$. These calculations can
be performed over a torsion-free ring with lifts of the coefficients,
since the numerical powers of $p$ in the ghost recursion
preserve every $T$-exponent. In the high rows the linear unit
terms may be taken separately before specializing the chosen matrix
$C_0$ modulo $p^\delta$. The nonlinear polynomials in the unit perturbations have the
orders just established for every choice of lifts. The linear
matrix identity $JC_0=T^{\delta D}I$ is used in $R$.

More explicitly, in the torsion-free verification all raw unit
perturbations have order at least $K_0+1$. In high rows before the
linear terms are combined via $C_0$, this suffices
to make the lower ghost powers vanish to order greater than $N_*$,
since $p(K_0+1)>N_*$. Thus the preceding induction uses only
the lift-independent bound $K_0+1$ for these intermediate rows.
After reduction the normalized output $(T^{-N_j}H_j)_j$ is
linear and invertible modulo $T$. The low and high diagonal blocks
are respectively the matrix computed above and the identity, and
the off-diagonal blocks vanish by $K_0+1>N_j$ in the low rows and
\eqref{eq:hhyleak} in the high rows. Every nonlinear coefficient
is divisible by $T$. Successive approximation over the complete
ring $R$ gives every value of the normalized output. This establishes
the chart, with a unit variable root retained throughout.

The low thresholds have sum $O(mp^{b-1})$ by a geometric series.
Since $K_0\ll mp^{b-1}$ and $D\ll p^{a-1}$, all thresholds satisfy
\[
              N_\Sigma\ll c\bigl(mp^{b-1}+p^{a-1}\bigr).
\]
Combining \eqref{eq:hhycost} and \eqref{eq:hcoeff} gives
\[
 \frac{h_k(A)}{k}\ll
       1+r c\frac{\log(2m)}{\log p}
                    \left(p^{1-c}+\frac1m\right)\ll 1+r .
\]
Here $64p^{c-1}\le m<64p^c$, so $\log(2m)/\log p\ll c+1$,
and $c(c+1)/p^{c-1}$ is absolutely bounded.

Finally suppose $m<64$. Use only the centres $T^{mi}$,
$0\le i\le a$, including the unit $i=0$. Take
$n_i=mi+1+m(p^{a-i}-1)/(p-1)$ and
$N_j=m(a-j)p^j+1+m(p^j-1)/(p-1)$ for $0\le j\le a$.
The low-row calculation applies in every row, including the unique
minimum $i=0$ in row $a$. The normalized matrix is again invertible,
$pN_j>N_{j+1}$, and all nonlinear terms disappear modulo $T$ after
normalization. Here $N_\Sigma\ll r+mp^{a-1}$. Since $tL\le m<64$
and $r^2\ll p^a$, \eqref{eq:hhycost} gives $h_k(A)\ll rk$.
At primes not dividing $k$ the bound $k+1$ is available. This
completes the proof.
\end{proof}

\subsection{Linear local widths in towers of fields}
\label{sec:local-ramification}

Write $e_v=v_{\mathfrak p_v}(p)$ for the absolute ramification index
at a finite place above $p$, and put
\[
                 E_{\mathbb K}(k)=\max_{v\mid k} e_v .
\]
There is no restriction on the residue degrees in the following estimates.

\begin{theorem}\label{thm:local-fixed-ramification}
For every fixed number field $\mathbb K$ one has
\[
                    H_{\mathbb K}(k)\le 4k+O_{\mathbb K}(1).
\]
Uniformly in $\mathbb K$, for all sufficiently large $k$ such that
$E_{\mathbb K}(k)\le\sqrt{k}/(\log(2k))^4$, one has
\[
                    H_{\mathbb K}(k)\le4k+1+\log_2 k.
\]
Both assertions allow the exponent to pass through powers of a fixed prime.
The leading constant $4$ in the first assertion is attained over
$\mathbb Q$ along the powers of two.
\end{theorem}

We give two finite-precision arguments. Let $A$ be the integer ring at a
place above $p$, $B$ its unramified coefficient ring, $F$ its residue
field and $e$ its absolute ramification index. For $a=v_p(k)\ge1$ write
$m=k/p^a$, $\delta=1$ for odd $p$ and $\delta=2$ for $p=2$.
Set $\Gamma=(F^\times)^m$, $f_0=[\mathbb F_p(\Gamma):\mathbb F_p]$
and $U_1=1+\mathfrak p_A$. Ramanujam's unramified estimate
\cite[Prop.~1]{hram63} gives
\begin{equation}\label{eq:unramified-ordinary-count}
 w_k(B)\le m\frac{p^{a+\delta}-1}{p-1}
                     +\mathbf1_{m=p-1}\le 4k .
\end{equation}

\begin{lemma}\label{lem:linear-local-deviation}
Put $j=\lceil\log_p e\rceil$. If $a\ge2j+\delta$ and
$k\ge e(a+\delta)$, then
\[
                 h_k(A)\le w_k(B)+\min(f_0,a+1).
\]
The bound $h_k(A)\le w_k(B)+a+1$ also holds, for every $e$, whenever
$f_0\ge a+1$.
\end{lemma}
\begin{proof}
For positive integers $t$, the $p$th-power map sends
$1+\mathfrak p_A^t$ into
$1+\mathfrak p_A^{\min(pt,t+e)}$ by the binomial theorem.
Starting at $t=1$, either the exponent has reached $e/(p-1)$
within the first $j$ steps, or those steps at least multiply the
exponent by $p^j\ge e$. One further step reaches $e$.
Each subsequent step gains at least $e$. Hence
$U_1^{p^a}\subset 1+p^{a-j}A$ for $a\ge j+1$.
In $R=A/p^{a+\delta}A$ the ideal $p^{a-j}R$ then has square zero
when $a\ge2j+\delta$.
The set of deviations $D=U_1^{p^a}-1$ in $R$ is an additive subgroup:
multiplication of $1+d_1$ and $1+d_2$ adds their deviations.
Raising to $m$ permutes $U_1$, so these are precisely the deviations
of the principal-unit $k$th powers. Being a finite additive subgroup,
$D$ is stable under $\mathbb Z_p$.

Nonunit powers vanish in $R$ when $k\ge e(a+\delta)$.
Each unit power has the form $\zeta(1+d)$, where $d\in D$ and
$\zeta$ is the Teichm\"uller lift of a member of $\Gamma$.
Choose $h_1,\ldots,h_{f_0}$ from $\Gamma$ forming a basis of its
additive span over $\mathbb F_p$. Their Teichm\"uller lifts
$\zeta_i$ form a $\mathbb Z_p$-basis of the unramified ring
$W(\mathbb F_p(\Gamma))$, by reduction modulo $p$. Consequently
$\sum_{\zeta}\zeta D=\sum_{i=1}^{f_0}\zeta_iD$.
Given a target in $J_k(A)$, collect its deviations as
$\sum_i\zeta_i d_i$. Subtract the $f_0$ permitted unit powers
$\zeta_i(1+d_i)$; what remains in $R$ is a linear combination over
$\mathbb Z$ of permitted Teichm\"uller constants. It lifts to
$J_k(B)$, so it costs at most $w_k(B)$ further powers. Moving one
unit root absorbs the error in $p^{a+\delta}A$, by the lifting
argument following \eqref{eq:nf-uniform-witt}.

For $f_0\ge a+1$ select $a+1$ independent elements $h_i$ of
$\Gamma$ and let $b_i$ be their Teichm\"uller lifts. Work over
$R_0=(B/p^\delta B)[[T]]$ in the uniformizer construction above.
The centered Witt map
$\sum_i[b_i+z_i]-\sum_i[b_i]$ on $(TR_0)^{a+1}$ has linear matrix
$(b_i^{p^j-1})_{0\le j\le a,i}$. Its reduction modulo $(p,T)$
is invertible by the Moore determinant criterion. It is therefore
invertible over $R_0$, and successive lifting in $T$ gives every Witt
vector with zero constant coordinates. Each $b_i+z_i$ has an $m$th
root in $R_0$, by Hensel lifting in $T$. Lift a target to the group
generated by $[y^m]$ as in \eqref{eq:nf-chart}, subtract its evaluation
at $T=0$, and apply the centered map to the difference. On evaluation
of the last ghost component at a uniformizer this gives $a+1$ unit
powers. The constant vector contributes an element of $J_k(B)$
modulo $p^{a+\delta}$, and the fixed baseline also consists of
$k$th powers in $B$. Represent their difference with $w_k(B)$ powers
and absorb the final error in one of the unit roots. This proves the
last assertion without using either size hypothesis.
\end{proof}

\begin{proof}[Proof of Theorem~\ref{thm:local-fixed-ramification}]
At unramified primes \eqref{eq:unramified-ordinary-count}, applied
to the target minus one, gives $h_k(A)\le4k+1$ if $p\mid k$.
At good primes $h_k(A)\le k+1$. For fixed $\mathbb K$ there remain
finitely many places. At each such place the lemma applies once
$a$ is large and $k\ge e(a+\delta)$; its additive cost is at most
the residue degree. For bounded $a$ the second inequality of
\eqref{eq:nf-residue-width} together with \eqref{eq:nf-uniform-witt}
gives $O_{A,a}(\log(2m))$. For fixed $A$, the latter estimate is
less than $4k$ once $k$ is large. This proves the first assertion.

For the uniform assertion it suffices to consider a place with
$h_k(A)>2k$. The proof of Theorem~\ref{thm:local-sparse}, applied
at this place, gives $p\ll(a+1)^2$ and
$m\ll(a+1)^2\log(2a+2)$ with absolute constants.
In particular $a\to\infty$ as $k\to\infty$ along these cases. One has
\[
 p^{\lfloor(a-\delta)/2\rfloor}
 \ge \frac{\sqrt{k/m}}{p^{(\delta+1)/2}}
 \gg \frac{\sqrt k}{(a+1)^3\sqrt{\log(2a+2)}} .
\]
For odd $p$ use $\delta=1$; the dyadic denominator is an absolute
constant. Since $a\le\log_2 k$, the hypothesis on $E_{\mathbb K}(k)$
implies $e\le p^{\lfloor(a-\delta)/2\rfloor}$ for all sufficiently
large $k$. It also implies $k\ge e(a+\delta)$. Apply
Lemma~\ref{lem:linear-local-deviation} and
\eqref{eq:unramified-ordinary-count}. The dyadic computation after
Theorem~\ref{thm:nf-bounded-multiplicity} gives sharpness.
\end{proof}

\subsection{A linear bound for larger ramification indices}
\label{sec:local-near-linear}
The square-zero condition in Lemma~\ref{lem:linear-local-deviation}
can be replaced by a bound on the nilpotence index. A simultaneous
$p$-adic moment calculation then controls the number of principal-unit
powers needed to represent all the deviations from Teichm\"uller constants.

\begin{lemma}\label{lem:hmomentcompress}
Let $L\ge1$. Given $n\in\mathbb Z_p$, there are at most
$(L+1)pL^2 2^{L-1}$ elements $u_i\in\mathbb Z_p$ for which
\[
                  \sum_i u_i^j=n\qquad(1\le j\le L).
\]
Consequently, if $d$ belongs to a nilpotent ideal of a finite
$\mathbb Z_p$-algebra and $d^{L+1}=0$, then $nd$ is a sum of
at most this many terms $(1+d)^{u_i}-1$.
\end{lemma}
\begin{proof}
Put $b=\lfloor\log_p L\rfloor$,
$B=v_p((L-1)!)+b+1$, $M=p^{B+b}$, and take
$r$ to be the least nonnegative residue of $n$ modulo $M$.
Use $r$ copies of $1$ and $M-1$ copies of each node $1,\ldots,L$.
The elementary symmetric functions of the remaining $L$ nodes
are prescribed by the polynomial of degree at most $L$ congruent to
\[
         (1+T)^{n-r}\prod_{i=1}^L(1+iT)^{1-M}
                       \pmod{T^{L+1}} .
\]
Its coefficients agree with those of $\prod_{i=1}^L(1+iT)$ modulo
$p^B$. Indeed if $z\in\mathbb Z_p$ and $1\le j\le L$, then
$\binom{Mz}{j}=(Mz/j)\binom{Mz-1}{j-1}$ lies in $p^B\mathbb Z_p$.
Let $Q(X)$ be the monic polynomial of the indicated elementary
symmetric functions, with alternating signs. Then $Q-P$ has all
coefficients in $p^B\mathbb Z_p$, where
$P(X)=\prod_{i=1}^L(X-i)$.

For each $i$, divide $Q(i+p^{b+1}Z)$ by $p^{b+1}P'(i)$.
The contribution of $P$ is congruent to $Z$ modulo $p$ because
$v_p(i-j)\le b$ for the other nodes $j$. The constant contribution
of $Q-P$ is integral by
$v_p(P'(i))\le v_p((L-1)!)$, and its positive-degree coefficients
vanish modulo $p$. The resulting polynomial reduces to $Z+c_i$
over $\mathbb F_p$, so $Q$ has a root in $i+p^{b+1}\mathbb Z_p$.
These $L$ balls are disjoint. The roots supply the prescribed
elementary symmetric functions. The product of $1+u_iT$ over
all nodes therefore equals $(1+T)^n$ modulo $T^{L+1}$.
Taking the logarithmic derivative gives the required moment identities.
Since $v_p((L-1)!)\le(L-1)/(p-1)$ and $p^{1/(p-1)}\le2$,
we have $M\le pL^2 2^{L-1}$. This proves the count.

For $2\le j\le L$, the polynomial $\binom{x}{j}$ has zero
constant term and takes value zero at $1$. The exact moment equations,
taken in $\mathbb Q_p$, give $\sum_i\binom{u_i}{j}=0$.
All these binomial values are $p$-adic integers. The binomial expansion
thus gives $\sum_i((1+d)^{u_i}-1)=nd$ in the given ring.
\end{proof}

The exponential dependence on $L$ in this lemma is sharp up to
a polynomial factor. At $p=2$ put $r=\lfloor(L+1)/2\rfloor$ and
\[
        E(X)=\prod_{i=0}^{r-1}(X-2i),\qquad
        O(X)=\prod_{i=0}^{r-1}(X-2i-1).
\]
Their resultant is odd. Thus there is a polynomial
$P\in\mathbb Z_2[X]$ of degree at most $2r-1$ with $P=0\pmod E$
and $P=1\pmod O$. At every even $x=2z\in\mathbb Z_2$,
$E(x)=2^rr!\binom zr$, and the same formula holds for $O(2z+1)$.
Writing $D=r+v_2(r!)$, one obtains
$P(x)\equiv0\pmod{2^D}$ on the even classes and
$P(x)\equiv1\pmod{2^D}$ on the odd classes.
Since $P(0)=0$ and $P(1)=1$, the moment equations give
$\sum_iP(u_i)=n$. Taking $n=2^D-1$ forces at least $2^D-1$
odd nodes. Here $D=2r-s_2(r)=L-O(\log(2L))$, where $s_2(r)$
is the binary digit sum.

\begin{theorem}\label{thm:hnearlinear}
For every fixed $\varepsilon>0$ and all sufficiently large $k$,
uniformly over the fields with $E_{\mathbb K}(k)\le k^{1-\varepsilon}$,
\[
 H_{\mathbb K}(k)\le 4k+
          O_{\varepsilon}\bigl(E_{\mathbb K}(k)(\log(2k))^3\bigr).
\]
More generally, for every real $c>\sqrt{\log2}$, uniformly over
all fields with $E_{\mathbb K}(k)\le k\exp(-c\sqrt{\log k})$,
\[
                         H_{\mathbb K}(k)\le (4+o(1))k .
\]
\end{theorem}
\begin{proof}
It suffices to consider places for which the local width exceeds
$2k$. Use the notation of Lemma~\ref{lem:linear-local-deviation},
and write $f=[F:\mathbb F_p]$, $P_0=a+\delta$,
$j=\lceil\log_p e\rceil$, $t_0=a-j$. We first suppose that
$t_0\ge1$ and $k\ge eP_0$. In the ring $R=A/p^{P_0} A$ nonunit
powers vanish. Each unit power is $\zeta(1+d)$ with $\zeta$ a permitted
Teichm\"uller constant and $d$ a principal-unit $k$th-power deviation.
The binomial estimate in Lemma~\ref{lem:linear-local-deviation}
gives $d\in p^{t_0}R$. Put $L=\lceil P_0/t_0\rceil$.

Let $V$ be the additive $\mathbb Z_p$-module generated in $R$ by all
$\zeta d$. It can be generated by at most $ef$ of the actual
generators. Indeed a submodule of
$(\mathbb Z/p^{P_0}\mathbb Z)^{ef}$ has at most $ef$ invariant
factors, and a subset of any generating set maps to a basis modulo
$pV$. Thus the deviation of any target from a signed sum of permitted
constants is $\sum_i n_i\zeta_i d_i$ with at most $ef$ terms and
$n_i\in\mathbb Z_p$. Lemma~\ref{lem:hmomentcompress} expresses
each $n_id_i$ as a sum of principal-unit $k$th-power deviations.
Indeed the principal units form a pro-$p$ group, so $p$-adic
exponentiation preserves its subgroup of $k$th powers. Restore
the Teichm\"uller constant in each term. The remaining target is a
signed sum of permitted constants in $R$ and hence a reduction of
an element of $J_k(B)$. Represent it by at most $w_k(B)$ powers.
We can reserve one extra unit power by subtracting $1$ at the start.
It absorbs the error in $p^{P_0}A$. We obtain
\[
       h_k(A)\le 4k+1+ef(L+1)pL^2 2^{L-1}.
\]

By Theorem~\ref{thm:local-sparse}, the places under consideration
have $p\ll(a+1)^2$ and $m\ll(a+1)^2\log(2a+2)$.
If $f_0\ge a+1$ the last assertion of
Lemma~\ref{lem:linear-local-deviation} already gives $4k+a+1$.
Otherwise $f\ll a+1$. Indeed for $|F|>16(m+1)^4$
the residue-field argument for \eqref{eq:nf-residue-width}
gives $F=\mathbb F_p(\Gamma)$ and $f=f_0$. For smaller $|F|$ we have
$f\ll\log(2m)/\log p$. Thus $pf\ll(\log(2k))^3$.

Put $X=\log k$ and $Y=\log(k/e)$. Then
$t_0\log p\ge Y-\log(mp)$, with $\log(mp)=O(\log(2X))$, and
$P_0\log p\le X+O(\log(2X))$. Under either hypothesis of the
theorem, $t_0\ge1$ and $k\ge eP_0$ for all sufficiently large $k$,
and
\[
                L\le\frac{X+O(\log(2X))}{Y-O(\log(2X))}+1 .
\]
If $Y\ge\varepsilon X$, the moment cost per generator is
$O_\varepsilon(p)$, proving the first bound. If $Y\ge c\sqrt X$,
the ratio of the excess bound to $k$ is at most
\[
    \exp\!\left( -Y+
       \frac{(\log2)(X+O(\log(2X)))}{Y-O(\log(2X))}
           +O(\log(2X))\right)+O((a+1)/k).
\]
It tends uniformly to zero when $c^2>\log2$.
\end{proof}

\subsection{When the ordinary and primitive widths agree}
\label{sec:local-precision}
Let $A$ be the integer ring of a finite extension of $\mathbb Q_p$,
with uniformizer $\pi$ and absolute ramification index $e$.
Here $w_k(A)$ permits arbitrary integral roots, and $h_k(A)$ requires at
least one unit root for each target in $J_k(A)$. A finite quotient
computes both widths once its lifting precision lies below the exponent.

\begin{lemma}\label{lem:precision-widths}
Suppose $1+\pi^c A\subseteq(A^\times)^k$ for some $1\le c<k$.
Let $U$ be the image of the unit $k$th powers in $A/\pi^cA$, and
let $J$ be the subgroup it generates. For $t\in J$, write $d^*(t)$
for the shortest \emph{nonempty} word in $U$ whose sum is $t$.
Then
\[
                   w_k(A)=h_k(A)=\max_{t\in J}d^*(t).
\]
One may take $c=1$ when $p\nmid k$, and otherwise
\[
              c=e v_p(k)+\lfloor e/(p-1)\rfloor+1 .
\]
The latter choice gives the equality whenever $k>c$.
\end{lemma}
\begin{proof}
Every unit power fills its whole coset modulo $\pi^cA$, since its
unit multiples by $1+\pi^cA$ are again $k$th powers.
Differences from a fixed unit power show that $\pi^c A\subseteq J_k(A)$.
The nonunit powers vanish modulo $\pi^c A$. As $J$ is finite, its
generating semigroup is a group, and every class, including zero, has
a nonempty word. One root in that word absorbs any error in $\pi^c A$.
This gives the primitive upper bound. For a nonzero class any ordinary
representation already uses at least $d^*(t)$ unit terms. In the zero
class use the eligible target $\pi^c$: a sum of nonunit powers belongs to
$\pi^k A$, so at least one unit occurs and $d^*(0)$ is a lower bound for
its ordinary length. This proves the formula.

When $p\nmid k$, the power map on $1+\pi A$ is surjective by Hensel's
lemma. For $p\mid k$ set $d=\lfloor e/(p-1)\rfloor+1$.
The logarithm and exponential identify $1+\pi^d A$ with $\pi^d A$,
and multiplication by $k$ sends that lattice onto
$\pi^{d+e v_p(k)}A$. This gives the stated precision.
\end{proof}

In particular $W_{\mathbb K}(k)=H_{\mathbb K}(k)$ for every
sufficiently large $k$ over a fixed field. More generally the same equality
holds uniformly whenever the displayed strict inequality holds at every
place dividing $k$, with no constraint at the other places. Thus a
difference of one between these two local invariants, as occurs at
exponent four over $\mathbb Q_2$, can persist at large exponents only
when the ramification indices grow.

\section{Local Waring numbers and their dependence on the field}
\label{sec:arith-local}
We now fix a completion and let the exponent vary. The unit powers near
a Teichm\"uller representative form an open coset in a logarithmic
lattice. Reducing modulo that lattice turns the number of summands
into a question about the length of a word in roots of unity. This
also identifies the local exceptions that remain visible after averaging
over global targets.

\subsection{The local constant along towers of prime powers}
\label{sec:local-cyclotomic}

Fix a finite extension of $\mathbb Q_p$ with integer ring $A$ and
residue field $F$. Let $m$ be prime to $p$ and let $k=mp^a$ with
$a\to\infty$. Write $h=|(F^\times)^m|$ and let $\mu_h\subset A$ be
the group of $h$th roots of unity. For $h>1$ set
$f_0=\operatorname{ord}_h(p)$, and put $f_0=1$ for $h=1$. Let
$E=\mathbb Q_p(\mu_h)$, an unramified field of degree $f_0$.
With $U_1=1+\mathfrak p_A$ define
\[
       L=\log U_1,\qquad S=\sum_{\zeta\in\mu_h}\zeta L,\qquad
       S\cap E=p^b\mathcal O_E.
\]
The last equality defines an integer $b$. Indeed, $L$ is a compact
additive $\mathbb Z_p$-lattice containing $p^\delta A$, where
$\delta=1$ for odd $p$ and $\delta=2$ for $p=2$. The lattice $S$
is stable under $\mathbb Z_p[\mu_h]=\mathcal O_E$.

For $h>1$ consider the Minkowski space
$V_h=\mathbb Q(\mu_h)\otimes_{\mathbb Q}\mathbb R$, the integer
lattice $\mathcal L_h=\mathbb Z[\mu_h]$ and the convex polytope
$P_h=\operatorname{conv}(\mu_h)$ in $V_h$, with its usual simultaneous
embedding. Write $\rho_h$ for the covering radius of $\mathcal L_h$ in
the gauge of $P_h$. Thus $\rho_h$ is the least $t$ for which
$\mathcal L_h+tP_h=V_h$. More generally write $\rho_{P_h}(I)$
for this radius with a nonzero ideal $I$ in place of $\mathcal L_h$.
The sum of all the roots in $\mu_h$ is zero and they span $V_h$, so
zero lies in the interior of $P_h$.

\begin{theorem}\label{thm:cyclotomic-local}
For sufficiently large $a$, depending on $A,m$, the ordinary and
primitive local widths have the following behavior.
\begin{enumerate}
\item If $h=1$, then $w_k(A)=h_k(A)=p^{a+b}$.
\item For $h>1$ let $\mathfrak p$ be the prime of $\mathcal L_h$
inducing the chosen embedding in $E$. Then
\[
 w_k(A)=\rho_{P_h}(\mathfrak p^{a+b})+O_h(1),\qquad
 h_k(A)=\rho_{P_h}(\mathfrak p^{a+b})+O_h(1).
\]
In particular both widths are $\asymp_{A,m}p^{af_0/\varphi(h)}$.
\item If $h>1$ and $f_0=\varphi(h)$, then
\[
       w_k(A)=\rho_h p^{a+b}+O_{A,m}(1),\qquad
       h_k(A)=\rho_h p^{a+b}+O_{A,m}(1).
\]
\end{enumerate}
\end{theorem}

\begin{proof}
For large $a$ one has $U_1^{mp^a}=1+p^aL$. To see the equality,
first take $a$ so large that $U_1^{p^a}$ lies in the domain of the
logarithm and exponential isomorphisms. Its logarithm is $p^a L$,
since raising to $m$ permutes $U_1$. As $L$ is bounded and contains
$p^\delta A$, the maps $x\mapsto p^{-a}(\exp(p^a x)-1)$ and
$x\mapsto p^{-a}\log(1+p^a x)$ carry $L$ into itself for large $a$:
both are congruent to $x$ modulo $p^\delta A$, uniformly for $x\in L$.
They are inverse maps. Thus the unit powers fill the individual cosets
$\zeta(1+p^a L)$, $\zeta\in\mu_h$. Every nonunit power lies in
$p^{a+\delta}A\subseteq p^aL$ once $k\ge e(a+\delta)$.

Differences of powers in a given unit class supply all of $p^a\zeta L$.
It follows that $J_k(A)=\mathcal O_E+p^a S$: the ordinary integer
span of $\mu_h$ is dense in $\mathcal O_E$, and the subgroup already
contains the open lattice $p^{a+\delta}A$. Reduction modulo $p^a S$
gives the finite additive quotient $\mathcal O_E/p^{a+b}\mathcal O_E$.
Let $D_a$ be its directed diameter for the generators $\mu_h$, allowing
a word of length zero. Nonunit roots have zero image, so that
\[
                    D_a\le w_k(A)\le h_k(A)\le D_a+h .
\]
For the upper bound reserve one unit in each of the $h$ classes.
Represent the remaining residue by a word of length at most $D_a$.
Varying the reserved units then adds any desired element of $p^a S$.

If $h=1$, write $M=p^{a+b}$. Every residue has a representative
$r$ with $1\le r\le M$. Its entire coset is a sum of $r$ unit powers.
Any representation of $M$ uses some unit root once $k>e(a+b)$:
a sum of nonunit powers lies in $\mathfrak p_A^k$, whereas $M$
has valuation $e(a+b)$. Its unit count is a positive multiple of $M$.
This proves the first assertion.

For $h>1$ choose the prime $\mathfrak p$ of
$\mathbb Q(\mu_h)$ inducing the given embedding in $E$. Its residue
degree is $f_0$, and $p\nmid h$ ensures that it is unramified.
Reduction identifies $\mathcal O_E/p^{a+b}\mathcal O_E$ with
$\mathcal L_h/\mathfrak p^{a+b}$. Put $n=\varphi(h)$.
For every nonzero ideal $I$ of $\mathcal L_h$, representatives
modulo $I$ can be chosen with archimedean coordinates
$O_h((NI)^{1/n})$. Indeed write $I=\gamma C$ with $C$ in a fixed
finite set of ideal-class representatives and multiply $\gamma$ by
a balancing unit as in Lemma~\ref{lem:hbal}. A bounded fundamental
parallelotope of $C$ then gives the asserted covering bound for $I$.
Coordinates in a fixed integral basis of powers of a primitive
$h$th root obey the same bound. A negative root is a sum of at most
$h-1$ positive generators, since $\sum_{\zeta\in\mu_h}\zeta=0$.
Thus $D_a\ll_h p^{(a+b)f_0/n}$. Conversely, the words of length
at most $D_a$ lie in $D_aP_h\cap\mathcal L_h$, a set of
$O_h((1+D_a)^n)$ elements. They cover a quotient of size
$p^{(a+b)f_0}$. This proves the order of growth.

For the more precise formula let $\gamma_h$ be the gauge of $P_h$.
On $\mathcal L_h$, directed word length differs from $\gamma_h$ by
$O_h(1)$. Indeed, in an expression $x=\sum r_\zeta\zeta$ with
nonnegative real coefficients of sum $\gamma_h(x)$, replace each
$r_\zeta$ by its integer part. The difference lies in the intersection
of a fixed compact set with $\mathcal L_h$ and costs a bounded
number of generators. The reverse inequality follows from convexity.
For $I=\mathfrak p^{a+b}$ the function
$x\mapsto\min_{z\in I}\gamma_h(x+z)$ on $V_h$ has maximum
$\rho_{P_h}(I)$ and a Lipschitz constant depending only on $P_h$.
Indeed the gauge and its reverse are bounded by a fixed multiple of
Euclidean length, so the triangle inequality bounds differences
uniformly in $I$. Every point of $V_h$ is within bounded distance
of $\mathcal L_h$. Its maximum on $\mathcal L_h/I$ therefore differs
from $\rho_{P_h}(I)$ by $O_h(1)$. The word estimate and the bounds
for $D_a$ give the asserted two formulas. If $f_0=\varphi(h)$, then
$\mathfrak p=p\mathcal L_h$ and homogeneity yields
$\rho_{P_h}(\mathfrak p^{a+b})=p^{a+b}\rho_h$.
\end{proof}

For pure dyadic powers $k=2^a$, $h=2^f-1$ where $f=[F:\mathbb F_2]$.
For $f\ge3$ we have $\varphi(2^f-1)>f$. For $f=3,4$ this follows
by inspection, and for $f\ge5$ it follows from $\varphi(n)\ge\sqrt n$
for odd $n$. Thus every fixed dyadic completion of residue degree at
least three has local widths $o(k)$. At $f=2$ one has $h=3$ and
$\rho_3=1$. Indeed the three generators in an integral basis are
$(1,0),(0,1),(-1,-1)$. For $0\le x\le y\le1$, the three choices of
representative at which a coefficient becomes zero have lengths
$x+y$, $x+2-2y$, $1+y-2x$. Their minimum is the directed distance on
the real torus: fixing the coefficient of the third generator reduces
the other two coefficients modulo one, and the cost grows between
consecutive jumps. The three costs sum to $3$ and equal $1$ at
$(x,y)=(1/3,2/3)$. The case $x\ge y$ follows by symmetry.
For example, if $A$ is the unramified dyadic ring of residue degree
two, $S=2A$. The inclusion $L\subseteq 2A$ follows from the logarithm
series and $U_1=1+2A$. Modulo $2A$ the image of
$\frac12\log(1+2t)$ is $\bar t-\bar t^2$, a nonzero map on $\mathbb F_4$;
its nonzero value multiplied by the three Teichm\"uller classes spans
$2A/4A$, and $4A\subset L$. Hence the local primitive and ordinary
widths in this example are $2k+O(1)$.

The covering radius can oscillate in a split cyclotomic field.
An exact example occurs over $\mathbb Z_5$. Write
$(2+i)^{a+1}=u_a+i v_a$. For every $a\ge1$,
\begin{equation}\label{eq:local-five-exact}
       w_{5^a}(\mathbb Z_5)=h_{5^a}(\mathbb Z_5)
                    =\max(|u_a|,|v_a|)-1 .
\end{equation}
Indeed, the unit powers are exactly the four cosets
$\zeta+5^{a+1}\mathbb Z_5$, $\zeta^4=1$, by the logarithm on
$1+5\mathbb Z_5$. The nonunit powers vanish in the quotient.
The quotient with these four generators identifies with
$\mathbb Z[i]/(u_a+i v_a)$, up to conjugating the chosen prime.
Its directed word length is distance in the $\ell^1$ norm to the
Gaussian ideal lattice, since the four steps are $\pm1,\pm i$. The
diameter of this Gaussian graph is known \cite{edgaussian}; we give a
geometric calculation to exhibit its dependence on the argument of the prime.

Put $M=\max(|u_a|,|v_a|)$ and $N=u_a^2+v_a^2=5^{a+1}$.
The ideal lattice has perpendicular basis vectors $(u_a,v_a)$ and
$(-v_a,u_a)$. On its centered Euclidean Voronoi square the largest
$\ell^1$ norm is $M$, attained only at vertices. Indeed, its four
vertices all have norm
$(|u_a+v_a|+|u_a-v_a|)/2=M$, and each of the four signed linear
functionals defining the norm has nonzero slope in both basis
directions since $|u_a|\ne|v_a|$. The distance of each vertex from
the ideal in the $\ell^1$ norm is exactly $M$: its four closest
Euclidean lattice points are at that $\ell^1$ distance, whereas every
other lattice point has Euclidean distance at least $\sqrt{5N/2}>M$.
Both coordinates of a vertex are half-integers, as $N$ is odd.
Consequently every integer point has distance at most $M-1$ from
the ideal, and an integer point at $\ell^1$ distance one from a
vertex has distance at least $M-1$. Thus the diameter is $M-1$.
It is at least three. Any nonempty word lifts to its whole local
coset by varying one of its unit roots, and the zero coset can use
the two unit classes $1,-1$. This proves \eqref{eq:local-five-exact}.

If $\theta=\arg(2+i)$, then $\theta/\pi$ is irrational, since a
root-of-unity power of $(2+i)/(2-i)$ would equate powers of distinct
prime ideals. The normalized widths in \eqref{eq:local-five-exact},
divided by $5^{(a+1)/2}$, therefore have the full interval
$[1/\sqrt2,1]$ as their set of limit points.

The order-three case admits an exact count.
\begin{corollary}\label{cor:local-triangular}
If $h=3$ and $p\equiv2\pmod3$, then for all sufficiently large $a$,
\[
                 h_{mp^a}(A)=w_{mp^a}(A)=p^{a+b}-1 .
\]
In particular, in the unramified quadratic extension of $\mathbb Z_2$,
one has $h_{2^a}(A)=w_{2^a}(A)=2^{a+1}-1$ for every $a\ge2$.
\end{corollary}
\begin{proof}
Let $M=p^{a+b}$. In the order-three case the quotient of
Theorem~\ref{thm:cyclotomic-local} is $(\mathbb Z/M\mathbb Z)^2$
with the three generators $u=(1,0)$, $v=(0,1)$ and $w=(-1,-1)$.
For integral $0\le x\le y<M$ the same minimization as for
$\rho_3$ gives the minimum of
\[
             x+y,\qquad x+2M-2y,\qquad M+y-2x .
\]
These sum to $3M$ and all three can equal $M$ only if $3\mid M$.
Their minimum attains $M-1$: for $M=3j+1$ take $(j,2j+1)$
and for $M=3j+2$ take $(j+1,2j+1)$. Thus the diameter is $M-1$.

For $M\ge7$ every class has a representative of length at most $M-1$
using at least two different directions. If a shortest word uses just
one direction $u$ with multiplicity $t\le M-4$, append $u+v+w=0$.
This also handles the empty word. For $M-3\le t\le M-1$, replace
$tu$ by $(M-t)(v+w)$, of length at most six. The three unit-deviation
lattices satisfy $S=L+\zeta L$ for a primitive cube root $\zeta$,
since $\zeta^2=-1-\zeta$. Any two distinct directions give this
same sum, so the chosen representative lifts to the whole corresponding
coset of $p^a S$ with the same number of roots.

For the unramified quadratic dyadic assertion one can work directly
modulo $4k$, $k=2^a$. Nonunit powers vanish for $a\ge2$ and
\[
        (1+2t)^k\equiv1+2k(t+t^2)\pmod{4kA}.
\]
This follows for $a=1$ by expansion and thereafter by squaring.
The image of $t\mapsto t+t^2$ on $\mathbb F_4$ is $\mathbb F_2$.
Thus each direction supplies its residue class and its translate by
$2k$ times that direction; each pair of distinct directions spans
$2kA/4kA$. Apply the word argument with $M=2k\ge8$, then move
one unit root to absorb any error in $4kA$.
\end{proof}

\subsection{An exact formula for inert prime-power directions}
\label{sec:local-simplex}
The directed triangle in Corollary~\ref{cor:local-triangular} extends
to the simplices formed by the roots of unity of odd prime order.
We give the integral formula as well as its real covering counterpart.

\begin{theorem}\label{thm:local-inert-simplex}
Let $h=\ell^u$ be an odd prime power and suppose that the residue
characteristic $p$ has order $\varphi(h)$ modulo $h$. In the notation
of Theorem~\ref{thm:cyclotomic-local}, for all sufficiently large $a$,
\[
       w_{mp^a}(A)=h_{mp^a}(A)
               =\frac{\varphi(h)}2\bigl(p^{a+b}-1\bigr).
\]
The real covering radii for the odd and the doubled groups are
\begin{equation}\label{eq:local-prime-power-radius}
       \rho_{\ell^u}=\frac{\varphi(\ell^u)}2,\qquad
       \rho_{2\ell^u}=\ell^{u-1}\frac{\ell^2-1}{4\ell}.
\end{equation}
\end{theorem}

\begin{proof}
Put $s=\ell^{u-1}$ and choose $\zeta_{\ell}=\zeta_h^s$. Grouping
the powers in the integral power basis gives
\[
       \mathbb Z[\zeta_h]=
       \bigoplus_{j=0}^{s-1}\zeta_h^j\mathbb Z[\zeta_\ell].
\]
All the directions in $\mu_h$ lie in these summands. In each summand
they are the $\ell-1$ basis vectors and the negative of their sum.
For $\mu_{2h}$ both signs of these $\ell$ directions occur. The gauge
is the sum of the individual gauges, since the nonnegative
coefficients of a point can be grouped independently in each summand.
The real covering radii therefore add. Reduction modulo any positive
integer $M$ also gives a direct sum, and the word lengths add in the
finite quotient.

In the odd case take one summand. Denote the coordinates of a torus
point by $x_1,\ldots,x_{\ell-1}$ modulo one and put $x_\ell=0$.
Its directed distance from the lattice is
\[
                    \min_{t\in\mathbb R/\mathbb Z}
                                \sum_{i=1}^{\ell}\{x_i+t\} .
\]
Here braces denote fractional part. Indeed $t$ is the coefficient
of the last direction modulo one, and for fixed $t$ the least
nonnegative coefficients in the other directions are their stated
fractional parts. The minimum occurs at $t=-x_j$ for some $j$.
Summing these $\ell$ candidates gives at most $\ell(\ell-1)/2$,
since $\{x_i-x_j\}+\{x_j-x_i\}\le1$. For equally spaced $x_j$
every candidate equals $(\ell-1)/2$. Adding over the summands
proves the first covering formula.

For the integral assertion take $M=p^{a+b}$. The inertness assumption
identifies the local quotient in Theorem~\ref{thm:cyclotomic-local}
with $\mathbb Z[\zeta_h]/M\mathbb Z[\zeta_h]$.
Write $r_1,\ldots,r_{\ell-1},r_\ell=0$ for the coordinates in one
summand modulo $M$. For a residue $r$ let $[r]_M$ denote its
representative in $[0,M)$. The preceding minimization, now over
integer coefficients modulo $M$, gives shortest word length
$\min_j d_j$, where
\[
                         d_j=\sum_i[r_i-r_j]_M.
\]
If the $r_i$ are distinct, the $d_j$ are distinct integers. Indeed
$d_j\equiv\sum_i r_i-\ell r_j\pmod M$ and $(\ell,M)=1$. Their mean
is $M(\ell-1)/2$, since every unordered pair of distinct residues
contributes $M$ to their sum. The minimum is consequently at most
$D=(\ell-1)(M-1)/2$. To attain $D$ use the residues
$r_j=\lfloor jM/\ell\rfloor$, $0\le j<\ell$, including zero.
For $M\ge\ell$ they are strictly increasing and
$d_j-d_0=jM-\ell\lfloor jM/\ell\rfloor$. These differences
permute $0,\ldots,\ell-1$, proving attainment. In the presence of a
tie the mean is at most $M(\ell-1)/2-M/\ell$, giving a word of
length at most $D-\ell$ once $M/\ell\ge(\ell-1)/2+\ell$.
Thus the diameter is $D$ for large $M$.

To lift the words, observe that with distinct $r_i$ a minimizing
candidate uses all but one of the $\ell$ directions, and those form
an integral basis. In a block with a tie, append one copy of every
direction to a shortest word. The extra copies sum to zero and the
preceding saving keeps the cost at most $D$. This also gives unit
terms in the zero block. Over all blocks the used directions contain
an integral basis of $\mathbb Z[\zeta_h]$. Their multiples of $L$
therefore sum to all of $S$. Varying a unit root in each of the
basis directions fills the whole coset of $p^aS$. Thus every target
has a primitive representation by at most $sD$ powers. The finite
quotient diameter $sD$ gives the matching lower bound even with
nonunit roots allowed.

For the doubled group write $\|x\|$ for distance to the nearest
integer. Permitting both signs changes the one-block minimization to
$\min_t\sum_i\|x_i+t\|$. A minimum occurs at $t=-x_j$ for some
$j$: the slopes are odd away from the breakpoints, and an upward
jump requires a summand to vanish. To bound the mean of these
candidates let a semicircle rotate uniformly about the unit circle.
Write $m$ for the number of points $x_i$ inside. Any fixed pair
lies on opposite sides with probability $2\|x_i-x_j\|$. Integrating
$m(\ell-m)\le(\ell^2-1)/4$ over all positions gives
\[
                \sum_{i<j}\|x_i-x_j\|\le(\ell^2-1)/8.
\]
The mean candidate is therefore at most $(\ell^2-1)/(4\ell)$.
Equally spaced $x_i$ give that value at every candidate, and adding
over the $s$ independent summands proves the second formula.
\end{proof}

\subsection{An explicit ramified quadratic example}
\label{sec:sqrt-two-local}
Let $A_e=\mathbb Z_2[2^{1/e}]$. Conjecture~3.3 of Anthony et al.\
in both versions v1 and v5 of \cite{AnthonyRam} predicts
$g_{2,e}(k)<g_{2,1}(k)$ for even $k$ with $(e,k)>1$, and equality
otherwise. In v1 the domain of their ordinary width is the additive
group generated by the powers, and in v5 it is the additive semigroup.
Both widths equal $w_k(A_e)$ because a primitive representation exists for every element
of the additive group $J_k(A_e)$ by \eqref{eq:nf-uniform-witt}
at $2\mid k$, and by simple residue lifting at $2\nmid k$.
For $A=A_2$, write $\pi=\sqrt2$. The value $w_4(A)=7$ was found by
Kowalczyk and Miska \cite[Example~6.28]{audKowMiska} and by
Anthony et al.\ \cite[Theorem~3.16, v5]{AnthonyRam}, whose construction
also supplies a unit. The quadratic family below contradicts the first
branch already at $k=12$.
\begin{theorem}\label{thm:sqrt-two-exact}
For $4\mid k$ the ordinary and primitive widths in $A$ agree and equal
\[
 w_k(A)=h_k(A)=
 \begin{cases}
       7,&k=4,\\
       17,&k=8,\\
       2^{v_2(k)+2},&k\ge12.
 \end{cases}
\]
In particular $w_k(A)=w_k(\mathbb Z_2)$ for every multiple of four
at least twelve.
\end{theorem}
\begin{proof}
Write $d=2^{v_2(k)}$. For odd $x\in\mathbb Z_2$ and
$y\in\mathbb Z_2$, expansion gives
\[
 (x+y\pi)^4=x^4+12x^2y^2+4y^4+
                          (4x^3y+8xy^3)\pi .
\]
The rational coordinate is $1$ modulo $16$; for odd $y$ the last
two terms of that coordinate sum to $0$ modulo $16$. Taking $x=1$
and $y=0,1,2,3$ gives all four coefficients of $4\pi$ modulo $4$.
Squaring a congruence $u^d\equiv1+dj\pi\pmod{4dA}$ gives
$u^{2d}\equiv1+2dj\pi\pmod{8dA}$ since $d\ge4$. Raising next
to the odd power $k/d$ permutes $j$ modulo $4$. Thus the unit
$k$th powers modulo $4dA$ are exactly
\[
                       1+dj\pi,\qquad j\bmod4 .
\]
Every element of each displayed coset lifts to a unit $k$th power,
by $(1+4A)^k=1+4dA$.

For $k\ge12$ one has $k>2(v_2(k)+2)$, so all nonunit powers lie in
$\pi^kA\subsetneq4dA$. In a representation of $4d$ there must be
a unit summand, and its unit count must be a multiple of $4d$ by
the rational coordinate. Conversely take any eligible target, choose
its rational coordinate modulo $4d$ in the range $1,\ldots,4d$,
and use that many unit powers. One of them supplies the desired
$\pi$-coordinate and the final lifting. The identical residue and
valuation argument over $\mathbb Z_2$, where every unit $k$th
power is $1$ modulo $4d$, gives $w_k(\mathbb Z_2)=4d$.

For $k=4,8$, the nonunit powers modulo $4dA$ are $0$ and
$D=2^{k/2}$, since a nonunit root of valuation one is $\pi$ times
a unit. For any rational coordinate choose its residue modulo $D$
in the range $1,\ldots,D$, and use that many units. At most
$4d/D-1$ copies of $D$ complete the rational coordinate. A unit
supplies the remaining $\pi$-coordinate and lifts the error.
The target $d\pi$ requires at least one unit. If it uses $u$ units
and $v$ nonzero nonunit classes, then $u+vD\equiv0\pmod{4d}$.
Here $u$ is a positive multiple of $D$; minimizing $u+v$ gives
$D+4d/D-1$. This is seven at $k=4$ and seventeen at $k=8$, as
claimed. The target belongs to $J_k(A)$ by the construction just
given, which represents each pair of the indicated residue classes.
\end{proof}

The pure extensions in that conjecture share the same eventual value.
\begin{corollary}\label{cor:pure-dyadic-width}
Let $e\ge1$, $A_e=\mathbb Z_2[2^{1/e}]$, and $j=\lceil\log_2 e\rceil$.
For every odd $m$ and $a\ge2j+3$ one has
\[
          w_{2^a m}(A_e)=h_{2^a m}(A_e)=2^{a+2}.
\]
Thus even when the ramification index varies, the local Waring number
equals that of $\mathbb Z_2$ in this range.
\end{corollary}
\begin{proof}
Write $e=2^t s$ with $s$ odd, and $F_e=\mathbb Q_2(2^{1/e})$.
The units $-1,5,-5$ are nonsquares over $\mathbb Q_2$ and remain
nonsquares in the odd-degree extension $F_s$. Each successive
quadratic step from $F_s$ to $F_e$ adjoins the square root of a
uniformizer $\varpi$. A nonsquare unit $u$ in the lower field stays
a nonsquare: if $(x+y\sqrt\varpi)^2=u$, then $xy=0$; the remaining
possibilities would make $u$ a square in the lower field or give
it odd valuation. The residue field of $F_e$ is $\mathbb F_2$, and
$i\notin F_e$, so its roots of unity are exactly $1,-1$.

Put $L=\log(1+\mathfrak p_{A_e})$. Then
$L\cap\mathbb Q_2=4\mathbb Z_2$. Indeed it contains $4\mathbb Z_2$.
If the intersection contained an element of valuation less than
two, it would contain $\frac12\log 5$, since $v_2(\log 5)=2$.
For a corresponding principal unit $u$ one would have
$\log(u^2/5)=0$, hence $u^2=5$ or $u^2=-5$, a contradiction.

The filtration calculation in Lemma~\ref{lem:linear-local-deviation}
gives $U_1^{2^{j+1}}\subset1+2A_e$, whence $L\subset2^{-j}A_e$.
For $a\ge2j+3$ the power series for $\log(1+z)$ and
$\exp z-1$ are congruent to $z$ modulo $2^{a+2}A_e$ for every
$z\in2^a L$: the terms of order $n\ge2$ have valuation at least
$n(a-j)-(n-1)\ge a+2$. The two maps are inverse and carry
$2^a L$ onto itself, since $4A_e\subset L$. The group $U_1^{2^am}$
lies in the exponential domain, so it equals $1+2^aL$.
Finally $2^a>e(a+2)$ already at $a=2j+3$ and continues to hold
thereafter. The case $h=1$ in Theorem~\ref{thm:cyclotomic-local}
now applies with $b=2$, at the given explicit threshold.
\end{proof}

The ordinary widths ten and fifteen at $(e,k)=(3,4)$ contradict
the predicted equality for $(e,k)=1$.
\begin{proposition}\label{prop:tame-four-ten}
If $A=\mathbb Z_2[\pi]$, $\pi^3=2$, then
\[
 J_4(A)=\mathbb Z_2+2\mathbb Z_2\pi+2\mathbb Z_2\pi^2,
 \qquad w_4(A)=h_4(A)=10 .
\]
In particular $H_{\mathbb Q(\sqrt[3]{2})}(4)=10$, whereas
$w_4(\mathbb Z_2)=15$.
\end{proposition}
\begin{proof}
The polynomial is Eisenstein. Write $M$ for the displayed lattice.
Every fourth power lies in $M$, as is seen by squaring twice.
To verify the widths it suffices to compute in $A/16A$,
since $(1+4A)^4=1+16A$. Fourth powers in this quotient depend only on
the root modulo $4A$. The following recurrence gives the complete
finite calculation. In coordinates for $1,\pi,\pi^2$, multiply two
triples modulo $16$ by replacing $\pi^{i+j}$ by
$2^{\lfloor(i+j)/3\rfloor}\pi^{(i+j)\bmod3}$. Let $D$ be the set
of fourth powers of the $64$ triples with entries $0,1,2,3$, and $U$
the subset from triples with odd constant coordinate. Start with
$A_0=\{0\}$, $C_0=\varnothing$ and put
\[
        A_{j+1}=A_j+D,\qquad
        C_{j+1}=(C_j+D)\cup(A_j+U).
\]
Thus $A_j,C_j$ count the ordinary and primitive residue classes formed
with at most $j$ terms, since $0\in D$. The successive cardinalities
are as follows.
\[
\begin{array}{c|rrrrrrrrrr}
 j&1&2&3&4&5&6&7&8&9&10\\\hline
 |A_j|&22&103&248&408&552&688&824&944&1008&1024\\
 |C_j|&16&88&224&376&512&640&768&896&992&1024
\end{array}
\]
The final size is $|M/16A|$ and $(8,2,2)\notin A_9$. Lifting a
unit term gives the assertion in $A$. At every odd prime the tame
primitive bound is at most five, proving the global local-width
statement. The rational value is also immediate: modulo $16$ the
only fourth powers are $0,1$, and the class $15$ requires fifteen
units. A nonzero $2$-adic integer can be divided by a fourth power of
two until its valuation belongs to $\{0,1,2,3\}$, when its nonzero
residue modulo $16$ is already a sum of at most fifteen units and
a unit term lifts.
\end{proof}
Parnami, Agrawal and Rajwade's fourth-power level computation already
implies the primitive lower bound here: a primitive zero with $s$
summands expresses $-1$ as a sum of $s-1$ fourth powers after division
by a unit fourth power \cite{ParnamiFourth,ChomiczFourth}.

The same widths occur in all four ramified quadratic extensions with
an even square-class representative.

\begin{corollary}\label{cor:even-radicand-widths}
Let $u\in\mathbb Z_2^\times$, $\pi^2=2u$, and
$A=\mathbb Z_2[\pi]$. For every $k$ divisible by four,
$w_k(A)=h_k(A)$ has the value in Theorem~\ref{thm:sqrt-two-exact}.
\end{corollary}
\begin{proof}
Put $\epsilon=0$ for $u\equiv1\pmod4$ and $\epsilon=1$ for
$u\equiv3\pmod4$, and write $\omega=\pi+2\epsilon$.
The pair $(1,\omega)$ is an integral basis. For a unit $x+y\pi$
the scalar $x$ is odd. Its fourth power has scalar coordinate
$x^4+12ux^2y^2+4u^2y^4$ and $\pi$-coordinate
$4xy(x^2+2uy^2)$. The first is congruent to $1$ modulo $16$ for
even $y$, and to $1+8\epsilon$ for odd $y$. The integer
$j=xy(x^2+2uy^2)$ has the parity of $y$. Thus the fourth-power
classes are precisely $1+4j\omega\bmod16A$: setting $x=1$, the
four values of $y(1+2uy^2)$ for $y\bmod4$ are $0,3,2,1$.
Since $\omega^2\in2A$, squaring $1+dj\omega\bmod4dA$ at
$d=2^a\ge4$ gives $1+2dj\omega\bmod8dA$. The odd cofactor
permutes $j\bmod4$. Each resulting class fills its entire coset,
by the same principal-unit lifting as in Theorem~\ref{thm:sqrt-two-exact}.

In the basis $(1,\omega)$ the unit classes have scalar coordinate
$1\bmod4d$. At $k=4,8$ the nonunit classes are $0$ and the
rational integer $D=2^{k/2}$ modulo $4dA$: a root of valuation one
is $\pi$ times a unit, and $\pi^k=(2u)^{k/2}\equiv D\pmod{4dA}$.
A root of valuation at least two gives zero at this modulus.
The lower test target $d\omega$ therefore requires a unit and
the minimum cost for its scalar class is $D+4d/D-1$.
Every eligible class has $\omega$-coordinate divisible by $d$,
so the construction in Theorem~\ref{thm:sqrt-two-exact}
attains that width. For $k\ge12$ every nonunit power lies strictly
deeper than $4dA$, and the target $4d$ again requires $4d$ unit
roots. The same scalar and $\omega$ coordinates give the upper bound.
\end{proof}

The two ramified quadratic fields with odd square-class representatives
have the following values.
\begin{corollary}\label{cor:odd-radicand-widths}
For $4\mid k$ the two widths in $\mathbb Z_2[i]$ are those of
Theorem~\ref{thm:sqrt-two-exact}. In $A=\mathbb Z_2[\sqrt{-5}]$
both widths equal five for $k=4$ and $2^{v_2(k)+1}$ for $k\ge8$.
\end{corollary}
\begin{proof}
Let $\pi=i$ or $\sqrt{-5}$, and write $\omega=1+\pi$.
A unit is $x+y\pi$ with $x,y$ of opposite parities. If $x$ is odd,
write $y=2v$. The scalar and $\pi$-coordinates of its fourth power
give $1+8v\omega\bmod16A$. The other type consists of $\pi$ times
units of the first type. In the Gaussian case $\pi^4=1$, and the
two unit classes at the general exponent are therefore
$1+2dj\omega\bmod4dA$ for $j=0,1$, by squaring and principal-unit
lifting. The element $\omega$ is a uniformizer. For $k=4$ the
nonzero nonunit class is $\omega^4=-4\bmod16A$. In the integral
basis $(1,\omega)$ choose $1\le r\le4$ unit classes to agree
with the scalar coordinate modulo four and at most three nonunit
classes to complete it modulo sixteen. This gives seven. The
target $8+8\omega$ forces a unit; its scalar equation
$r-4v\equiv8\pmod{16}$ with $r>0$ requires $r+v\ge7$.
For $k=8$ the nonunit class is $16\bmod32A$, with width seventeen
as before, using $16\omega$ as the test target.
For $k\ge12$ the scalar coordinate and the test target $4d$
give width $4d$.

For $\pi=\sqrt{-5}$ one has $\pi^4\equiv9\pmod{16A}$.
Hence the four unit classes fill $1+8A$ exactly. Successive squaring
and the odd cofactor give $(A^\times)^k=1+2dA$. At $k=4$ the
nonunit classes modulo $8A$ are zero and four, while all nonunit
fourth powers have $\pi$-coordinate divisible by sixteen:
$\omega^4=-4-16\pi$ and the unit fourth powers lie in $1+8A$.
The target $8\pi$ thus requires a unit. With $r>0$ unit classes
and $v$ copies of four, $r+4v\equiv0\pmod8$ requires $r+v\ge5$.
Every eligible class uses at most four units and one nonunit
by matching its scalar coordinate first modulo four and then
modulo eight; one unit fills the error in $8A$.
At $k=8$, all nonunit powers vanish modulo sixteen and are rational
modulo $32A$, since $\omega^8\equiv16\pmod{32A}$. Thus $16\pi$
requires sixteen units. The unit coset $1+16A$ supplies the
matching upper bound. For $k\ge12$ all nonunit powers are
strictly deeper than $2dA$, so the rational target $2d$ forces
$2d$ units. Every eligible coset is filled by at most $2d$ units.
In each small case the upper construction also gives eligibility
of the displayed test target.
\end{proof}

\subsection{The limiting local width of a fixed number field}
\label{sec:local-signature}

When $h=2$ in Theorem~\ref{thm:cyclotomic-local}, its leading term
has an exact form. The residue characteristic is odd, $E=\mathbb Q_p$,
and $S=L$ because $L$ is additive. For all sufficiently large $a$
one has
\begin{equation}\label{eq:local-two-directions}
            w_{mp^a}(A)=h_{mp^a}(A)=\frac{p^{a+b}-1}{2}.
\end{equation}
Indeed, the two unit-power cosets are $1+p^aL$ and $-1+p^aL$.
The quotient is cyclic of odd order $M=p^{a+b}$ and has diameter
$(M-1)/2$ for the generators $\{1,-1\}$. A nonempty word represents
its whole coset after one unit root is varied. The zero class uses
the two opposite unit directions, and $(M-1)/2\ge2$ for large $a$.
The quotient diameter is also a lower bound for the ordinary width.

For each finite place $v$ of $\mathbb K$ write $A_v=\mathcal
O_{\mathbb K_v}$, let $q_v$ be its residue-field cardinality, and
write $p$ for its residue characteristic. For each $h\mid q_v-1$
put $d_{v,h}=(q_v-1)/h$, and define $b_{v,h}$ by
$S\cap\mathbb Q_p(\mu_h)=p^{b_{v,h}}\mathcal O_{\mathbb Q_p(\mu_h)}$
as in Section~\ref{sec:local-cyclotomic}. Set
\[
 c(v,h)=\frac{1}{d_{v,h}}
 \begin{cases}
       p^{b_{v,h}},& h=1,\\
       \rho_h p^{b_{v,h}},& h>1,\quad
                         \operatorname{ord}_h(p)=\varphi(h),\\
       0,&h>1,\quad \operatorname{ord}_h(p)<\varphi(h).
 \end{cases}
\]
Choose an odd rational prime $p_0$ admitting an unramified place of
degree one in $\mathbb K$. Such primes exist by the prime ideal theorem:
prime ideals of residue degree at least two and norm at most $X$
number $O_{\mathbb K}(\sqrt X)$, and only finitely many primes ramify.
Let $\mathcal V$ consist of the places over rational primes $p\le p_0$,
together with the places over primes ramified in $\mathbb K$.

\begin{theorem}\label{thm:local-signature}
For each fixed number field one has
\begin{equation}\label{eq:local-signature}
 \limsup_{k\to\infty}\frac{W_{\mathbb K}(k)}k
 =\limsup_{k\to\infty}\frac{H_{\mathbb K}(k)}k
 =\max_{v\in\mathcal V}\max_{h\mid(q_v-1)}c(v,h).
\end{equation}
Moreover, given $\varepsilon>0$, apart from a finite set every integer
$k$ with $H_{\mathbb K}(k)>(1+\varepsilon)k$ lies in one of finitely
many sequences $k=mp^a$, with $p$ prime, $p\nmid m$, and $a\ge1$.
The finite exceptional set and the list of pairs $(p,m)$ may depend
on $\mathbb K$ and $\varepsilon$. In particular the number of such $k\le X$ is
$O_{\mathbb K,\varepsilon}(\log(2X))$.
\end{theorem}

\begin{proof}
At an odd unramified place write $k=mp^a$ with $p\nmid m$,
allowing $a=0$. The unit powers consist precisely of
$\zeta(1+p^{a+1}A_v)$ for $\zeta\in\mu_h$, where
$h=(q_v-1)/\gcd(m,q_v-1)$. This follows from the Teichm\"uller
decomposition and the logarithm isomorphism on $1+pA_v$.
Every nonunit power vanishes modulo $p^{a+1}$ since $k\ge a+1$.
Therefore the image of $J_k(A_v)$ at that modulus is the image of
$\mathcal O_E=\mathbb Z_p[\mu_h]$, where $E=\mathbb Q_p(\mu_h)$.
Let $L_0$ denote the additive width of the residue subfield generated by
the image of $\mu_h$, using those residues and zero. As in
\eqref{eq:nf-residue-width}, $L_0\le m$. In the unramified ring
$\mathcal O_E$ first match the residue of a target by at most $L_0$
Teichm\"uller directions, subtract this sum, and divide by $p$.
Iterating through digit $a$ and repeating each term in digit $i$
exactly $p^i$ times costs at most $L_0(1+p+\cdots+p^a)$ terms.
Perform the construction on the target minus $1$. The reserved unit
power absorbs the residual in $p^{a+1}A_v$, giving
\begin{equation}\label{eq:local-signature-unram}
 h_k(A_v)\le1+L_0(1+p+\cdots+p^a)
             \le1+\frac{p}{p-1}k .
\end{equation}

We also need uniformity in $m$ at a single fixed $v$. There is a
constant $C_v$ such that
\begin{equation}\label{eq:local-signature-prime}
                          h_{mp^a}(A_v)\le C_vp^a
                          \qquad(p\nmid m,\ a\ge0).
\end{equation}
In fact the power map $u\mapsto u^m$ permutes the principal units
for every such $m$. The logarithm argument of
Theorem~\ref{thm:cyclotomic-local} gives
$U_1^{mp^a}=1+p^aL$ for all $a$ past a threshold depending only on
$A_v$. Increasing the threshold ensures $p^a\ge e_v(a+\delta)$.
Every nonunit $mp^a$th power then vanishes modulo $p^{a+\delta}A_v$
for every $m\ge1$. For these $a$, the residue quotient, its directions
and the upper bound $D_a+h$ in that proof depend only on $a$ and
$h$. Its covering estimate is $O_{v,h}(p^a)$, also when $h=1$.
There are finitely many divisors $h$ of $q_v-1$. At the finitely
many smaller $a$, use the finite additive group
$A_v/p^{a+\delta}A_v$. Any subset generates its additive span by
words of length at most the order of the group: as long as the
successive sumsets with zero added grow, their cardinality increases.
Apply this after subtracting one unit. The logarithm on
$1+p^\delta A_v$ gives $(1+p^\delta A_v)^{mp^a}=1+p^{a+\delta}A_v$,
so the reserved unit lifts every residue solution. This proves
\eqref{eq:local-signature-prime} uniformly in $m$.

For fixed $v,m$ define $d=\gcd(m,q_v-1)$ and $h=(q_v-1)/d$.
Theorem~\ref{thm:cyclotomic-local} gives the limit
\begin{equation}\label{eq:local-signature-sequence}
       \lim_{a\to\infty}\frac{h_{mp^a}(A_v)}{mp^a}
                            =\frac{d}{m}c(v,h).
\end{equation}
In the split case this is zero because its growth exponent is
strictly smaller than one. In the other cases the limit is largest,
for a specified $h$, at $m=d=d_{v,h}$. Consequently every entry
in the maximum in \eqref{eq:local-signature} is attained as a
subsequential local limit. At the chosen degree-one place over
$p_0$ take $h=1$ and $m=p_0-1$. There $L=p_0\mathbb Z_{p_0}$
and $b_{v,1}=1$. In particular the maximum, denoted $C$, is at
least $p_0/(p_0-1)>1$.

Fix $t>C$. At all places outside $\mathcal V$,
\eqref{eq:local-signature-unram} bounds $h_k(A_v)/k$ by
$p_0/(p_0-1)+1/k<t$ for large $k$. At places over $p\nmid k$ the
simple-root bound $h_k(A_v)\le k+1$ also suffices. At the remaining
finitely many places, \eqref{eq:local-signature-prime} bounds the
ratio by $C_v/m$. Thus it is enough to consider bounded $m$ at each
of these places. Along every resulting sequence $mp^a$ the limit
in \eqref{eq:local-signature-sequence} is at most $C$. This proves
the upper bound. Since $W_{\mathbb K}(k)\le H_{\mathbb K}(k)\le
W_{\mathbb K}(k)+1$, the ordinary limsup agrees with the primitive one.

For the last assertion take all rational primes with
$1/(p-1)\ge\varepsilon/2$, together with $2$ and the finitely many
ramified primes. Outside this set \eqref{eq:local-signature-unram}
suffices for large $k$. Inside the set the bound
$C_v/m$ again leaves finitely many $m$. Each remaining sequence
contributes $O(\log(2X))$ exponents up to $X$.
\end{proof}

The finite maximum takes values in a considerably smaller set than
its definition suggests.

\begin{corollary}\label{cor:local-signature-values}
If $C_{\mathbb K}$ denotes the two limits superior in
\eqref{eq:local-signature}, then
\[
 C_{\mathbb K}\in \{4,2,4/3\}\ \cup\
                  \{p/(p-1):p\text{ an odd prime}\}.
\]
One has $C_{\mathbb K}=4$ exactly when there is a dyadic place $v$
with residue field $\mathbb F_2$ and
\[
             \log(1+\mathfrak p_v)\cap\mathbb Q_2=4\mathbb Z_2.
\]
In the absence of such a place one has $C_{\mathbb K}\le2$.
\end{corollary}

\begin{proof}
Theorem~\ref{thm:local-inert-simplex} gives the covering radii for
odd prime powers and twice odd prime powers. In particular
$\rho_3=1$, $\rho_5=2$ and $\rho_6=2/3$. For even $h$ both signs of
an integral power basis occur among the directions, and rounding every
coordinate to the nearest integer gives $\rho_h\le\varphi(h)/2$.
One has $\rho_2=1/2$ and $\rho_4=1$, with equality in these bounds
attained at the center of a unit interval and a unit square.

It remains to consider the inert entries in $c(v,h)$. Write
$n=\varphi(h)$. For $h>2$, cyclicity of $(\mathbb Z/h\mathbb Z)^\times$
forces $h$ to be $4$, an odd prime power or twice an odd prime power.
Indeed divisibility by $8$, by two distinct odd primes, or by $4$
and an odd prime gives at least four square roots of unity by the
Chinese remainder theorem. In particular $h\le3n$, and for odd
$h$ one has $h\le3n/2$. Together with
\eqref{eq:local-prime-power-radius} this gives $\rho_h\le n/2$
for all the entries under consideration.

For odd $p$, $p\mathcal O_E\subset L$, so $b_{v,h}\le1$.
Since $q_v\ge p^n$, when $n\ge4$ one has
\[
     c(v,h)\le\frac{3p n^2}{2(p^n-1)}<1.
\]
The right side decreases with $p\ge3$ and with $n\ge4$, and at
$(p,n)=(3,4)$ it equals $9/10$. When $n=2$ the three choices are
$h=3,4,6$. For $h=3,6$ an odd inert $p$ is at least $5$, and the
radii $1,2/3$ give $c(v,h)\le5/6$. For $h=4$, $c(v,4)>1$ requires
$p=3$, $q_v=9$ and $b_{v,4}=1$, giving $c(v,4)=3/2$. Finally
$h=1,2$ give the same value $p^{b_{v,1}}/(q_v-1)$, since $S=L$
for either group. This exceeds one only when $q_v=p$ and $b_{v,1}=1$.

For $p=2$ and $h>1$, one has $b_{v,h}\le1$. In fact $4\mathcal O_E$
lies in $L$ and
\[
        \tfrac12\log(1+2x)\equiv x-x^2\pmod {2\mathcal O_E}
                      \qquad (x\in\mathcal O_E).
\]
The residue degree of $E$ is at least two. Choose a nonzero value of
this map on its residue field; its products with the reductions of
$\mu_h$ span that field over $\mathbb F_2$. Thus $2\mathcal O_E\subset S$.
Here $h$ is odd, and for $n\ge6$ the sharper bound $h\le3n/2$ gives
$c(v,h)\le3n^2/(2(2^n-1))<1$. For $n=2$ one has $h=3$ and
$\rho_3=1$; the only possible value above one is $2$, with $q_v=4$.
For $n=4$ one has $h=5$ and $\rho_5=2$; the only possible value
above one is $4/3$, with $q_v=16$. In each case the residue degree
is a multiple of $n$. For $h=1$ the value is $2^{b_{v,1}}/(q_v-1)$
with $b_{v,1}\le2$. Its only possibilities above one are $4$, $2$
and $4/3$, and $4$ requires exactly $q_v=2$ and $b_{v,1}=2$.
The maximum in Theorem~\ref{thm:local-signature} is greater than
one, which completes the classification.
\end{proof}

In fact $W_{\mathbb K}(k)=H_{\mathbb K}(k)$ for all sufficiently
large $k$ at fixed $\mathbb K$, by Lemma~\ref{lem:precision-widths}.
There is also a counting form of Theorem~\ref{thm:local-signature}.
For $p\nmid m$ put
\[
 \Lambda(p,m)=\max_{v\mid p}
   \frac{\gcd(m,q_v-1)}{m}\,
   c\left(v,\frac{q_v-1}{\gcd(m,q_v-1)}\right).
\]
If $t>1$ differs from every $\Lambda(p,m)$, then
\begin{equation}\label{eq:local-signature-count}
 \#\{k\le X:H_{\mathbb K}(k)>tk\}
    =\left(\sum_{\Lambda(p,m)>t}\frac{1}{\log p}\right)\log X
                                  +O_{\mathbb K,t}(1).
\end{equation}
Indeed, the proof of the theorem leaves finitely many candidate
pairs $(p,m)$. Along $k=mp^a$ the largest width over places dividing
$p$, normalized by $k$, tends to $\Lambda(p,m)$ by
\eqref{eq:local-signature-sequence}. Places over primes dividing
$m$ have bounded width by \eqref{eq:local-signature-prime}, and at
other places the simple-root bound is $k+1$. Thus each candidate
sequence eventually contributes either all its terms or none.
Two different primes give sequences meeting at most once, since
taking the valuation at either prime fixes its exponent; at one
prime the sequences with distinct $m$ are disjoint. Counting the
terms of each progression gives \eqref{eq:local-signature-count}.

Removing the pure powers of two gives the sharp bound
\begin{equation}\label{eq:local-signature-off-two}
       \limsup_{\substack{k\to\infty\\ k\ne2^a\ (a\ge0)}}
                  \frac{H_{\mathbb K}(k)}{k}\le\frac32 .
\end{equation}
Indeed every entry above $3/2$ in the proof of
Corollary~\ref{cor:local-signature-values} comes from a dyadic
place, either $q_v=2$, $h=1$, $b_{v,1}=1$ or $2$, or
$q_v=4$, $h=3$, $b_{v,3}=1$. In all these cases $d_{v,h}=1$.
For odd $m>1$ the actual sequence limit in
\eqref{eq:local-signature-sequence} is therefore at most $4/m\le4/3$.
The finitely many candidate sequences argument for each $t>3/2$
proves \eqref{eq:local-signature-off-two}. Equality occurs over
$\mathbb Q$ on $k=2\cdot3^a$: at $3$ one has $h=1$, $b=1$ and
the ordinary local width is $3^{a+1}$ for large $a$. When the
fixed-field constant exceeds $3/2$, the ratio along $k=2^a$
converges to that constant, as all nondyadic places then have
primitive width at most $k+1$.

The condition for $C_{\mathbb K}=4$ also has a field-theoretic
form. At a dyadic completion $B=\mathbb K_v$ of residue degree one,
let $2^r$ be the number of its roots of unity of $2$-power order.
Then
\begin{equation}\label{eq:local-signature-four}
 b_{v,1}=2\quad\Longleftrightarrow\quad
                  B(\mu_{2^{r+1}})/B\text{ is ramified}.
\end{equation}
Indeed $\log5$ generates $4\mathbb Z_2$. Since $L\cap\mathbb Q_2$
contains $4\mathbb Z_2$, one has $b_{v,1}<2$ precisely when
$\tfrac12\log5\in L$. The kernel of the logarithm on principal units
consists of $2$-power roots of unity: a sufficiently large $2$-power
places any principal unit in the domain where the logarithm is injective.
Thus the latter condition says exactly that $5$ is a square times a
$2$-power root of unity in $B$. All units of $B$ have residue one.
The polynomial $X^2-X-1$ has irreducible separable reduction over
$\mathbb F_2$, so $B(\sqrt5)$ is the unramified quadratic extension.
In particular $5$ is a nonsquare. The torsion group has exactly one
nontrivial square class, represented by a primitive $2^r$th root
$\zeta$; a square root of $\zeta$ has order $2^{r+1}$. Consequently
$\tfrac12\log5\in L$ exactly when the two quadratic extensions
$B(\sqrt5)$ and $B(\mu_{2^{r+1}})$ agree. This proves
\eqref{eq:local-signature-four}.

Every value in Corollary~\ref{cor:local-signature-values} occurs,
even when $\mathbb K$ is required to be totally real of degree at
most five. The field $\mathbb Q$ gives $4$. In $\mathbb Q(\sqrt5)$
the dyadic place is unramified of residue degree two, so $h=3$,
$b_{v,3}=1$ gives $2$; no other entry exceeds it. For $4/3$ take a
totally real quartic in which $2$ and $3$ are inert. At $2$ the
entry with $h=5$ has $q_v=16$, $d_{v,5}=3$ and $b_{v,5}=1$:
$2\mathcal O_E\subset S$ as above, while $L\subset 2A_v$ by the
logarithm series in an unramified dyadic ring. No other dyadic entry
is larger, the place above $3$ has residue degree four, and each
$p\ge5$ contributes at most $p/(p-1)\le5/4$.

For a prescribed odd prime $p$ take a totally real quintic in which
every prime below $p$ is inert and $p$ has an unramified degree-one
place. At the primes below $p$ none of the entries found in the proof
of the corollary exceeds one, by their residue degrees. The specified
place above $p$ gives $p/(p-1)$; all larger primes give at most their
corresponding ratios. To construct the fields in either example,
choose a monic irreducible polynomial of the required degree over
each finite field where inertia is prescribed. At the prime where
a degree-one place is wanted choose a squarefree product of a linear
polynomial and an irreducible polynomial of degree four. The Chinese
remainder theorem prescribes the coefficients modulo the product $M$
of the primes. Take $N$ a large multiple of $M$, start with
$\prod_{j=1}^s(X-jN)$ for $s=4$ or $5$, and add to its nonleading
coefficients bounded integers with the prescribed residues. After
putting $X=Ny$ and dividing by $N^s$, the resulting polynomial tends
to $\prod_{j=1}^s(y-j)$, so it has $s$ distinct real roots. Its
reduction at $2$ is irreducible. Squarefree reduction at the chosen
primes gives the desired unramified splitting by Dedekind's theorem.

\subsection{A linear bound for prime powers outside three residue fields}
\label{sec:local-teich}
A short additive relation between Teichm\"uller constants reduces the
cost of the successive Witt coefficients. In the following result
neither the ramification index nor the residue degree is bounded.

\begin{theorem}\label{thm:hteich}
Let $A$ be the integer ring in a finite extension of $\mathbb Q_p$
with residue field $F$. There is an absolute constant $C$ such that
\[
                 h_{p^a}(A)\le Cp^a\qquad(a\ge1)
\]
whenever $|F|\notin\{2,3,4\}$. For $p\ge5$ one has the more
precise uniform estimate
\[
                    h_{p^a}(A)\ll p^{a-1}+3^a .
\]
For every number field one has the explicit bound
\[
                    H_{\mathbb K}(p^a)\le p^a+1\qquad(p\ge5).
\]
For $p=2,3$ one also has $H_{\mathbb K}(p^a)\ll p^a$ whenever all
residue fields at $p$ satisfy the displayed restriction.
\end{theorem}

\begin{proof}
Here $m=1$ and every Teichm\"uller constant in $B=W(F)$ is permitted.
Use the chart of Theorem~\ref{thm:hhybrid} for $m=1$, with
$N_\Sigma\ll r+p^{a-1}$. Suppose first that each coefficient
of $B/p^eB$ is a sum of $O(D^e)$ Teichm\"uller constants for some
$1<D<p$. In the filtration by cyclic modules used to prove
\eqref{eq:hhycost}, let their lengths be $e_i$. Then
\[
 \sum_i D^{e_i}
   \le D\sum_{s\ge0}D^s\sum_i(e_i-s)_+
   \le D\sum_{s\ge0}D^s\ell(p^sM_N),\qquad M_N=W_r(R)/I_N.
\]
For odd $p$ the lengths are at most $N_\Sigma/p^s+r-s$ for
$0\le s\le a$ and vanish thereafter. For $p=2$ the kernel of
the coefficient reduction is killed by two as above, so one uses
the characteristic-two bounds with $s$ replaced by $s-1$ for $s\ge1$
and $2N_\Sigma$ at $s=0$. Geometric summation therefore gives
$h_{p^a}(A)\ll_{p,D}r+N_\Sigma+D^a$ for fixed $p,D$.
For $D=3$ and $p\ge5$, the implied constant here is absolute
if the constant in the digit cost is absolute, since $3/p\le3/5$.

To obtain this cost, write $[x]$ for the Teichm\"uller lift of
$x\in\mathbb F_p$ to $\mathbb Z_p$. The $p$ elements
$[x]+1-[x+1]$ are multiples of $p$ and have sum $p$. At least one
has valuation exactly one. At odd $p$ it is a sum of at most three
Teichm\"uller constants, since $-[x]=[-x]$. Use this element $\alpha$
as a uniformizer of $B$. In the expansion of a coefficient modulo
$p^eB$ in base $\alpha$, each residue digit is a single Teichm\"uller
constant of $F$. Its product with $\alpha^s$ costs at most $3^s$
constants, because their products are again Teichm\"uller constants.
This proves the second estimate for $p\ge5$, using $r\ll3^a$.

For the explicit estimates let $Q_j=(p^j-1)/(p-1)$.
In the chart with $m=1$, $N_j=(a-j)p^j+1+Q_j$. Consequently
the $p^s$-multiple module has length $M_{a-s}$, where
\[
       M_n=\sum_{j=0}^n\bigl((n-j)p^j+1+Q_j\bigr).
\]
Indeed $\lceil N_j/p^s\rceil=(a-j)p^{j-s}+1+Q_{j-s}$
for $j\ge s$, since $Q_j=p^s Q_{j-s}+Q_s$ and $0<1+Q_s\le p^s$.
Taking second differences, the ambient module has one cyclic factor
of length $a+1$ and $2p^{a-s}$ of length $s$ for $1\le s\le a$.
Here $M_0=1$, $M_1=4$, and
$M_n-2M_{n-1}+M_{n-2}=2p^{n-1}$ for $n\ge2$.
The cost of one coefficient of length $s$ is at most
$B_s=(3^s-1)/2$. The sum of these costs for the filtration by
actual power generators is at most its value on the ambient
invariant factors. This follows from $B_0=0$, the nonnegativity
of the second differences of $B_s$, and the bounds $\sum_i(e_i-s)_+\le\ell(p^sM_N)$ for every $s$.
Thus the coefficient cost is at most
\[
 C_a=(3^{a+1}-1)/2+\sum_{s=1}^a p^{a-s}(3^s-1),
 \qquad C_0=1 .
\]
It satisfies
$C_{a+1}=pC_a+\{3^{a+1}(5-p)+(p-3)\}/2$.
At $p=5$ this gives $C_a=(5^{a+1}-1)/4$.
For $p\ge7$ one has $C_1=6\le p-1$ and
$C_{a+1}\le p C_a-1$, so $C_a\le p^a-a$ for $a\ge1$.
The chart adds $a+1$ unit or nonunit powers, including a unit.
Together with the tame bound $p^a+1$ at other places, this proves
the explicit estimate at $p\ge7$. We complete the case $p=5$ below.

It remains to consider $p=2,3$. We give an argument uniform in
the residue degree. Fix $(v,s)=(4,11)$ for $p=2$ and $(2,7)$ for
$p=3$, and put $q=|F|$. For $q$ sufficiently large, each class
modulo $p^{v+1}B$ is a sum of $s$ Teichm\"uller constants.
Indeed the trace pairing of $B/p^{v+1}B$ is perfect. Every
nonprincipal additive character, after reducing to its conductor,
gives an exponential sum on the Teichm\"uller representatives of
absolute value at most $(p^v-1)\sqrt q$ by \cite{hkhc}.
At the conductor the polynomial is a unit times $X$. Its degree-one
reduction rules out a Frobenius difference plus a constant.
Its weighted degree in the Galois ring estimate is at most $p^v$.
Orthogonality gives the number of representations of each class as
$q^{s-v-1}+O_{p,v,s}(q^{s/2})$, positive since $s>2(v+1)$.

In particular some sum $\alpha$ of $s$ Teichm\"uller elements
is $p^v$ modulo $p^{v+1}B$. Expand any coefficient modulo $p^eB$
in base $\alpha$ with digits representing $B/p^vB$.
Each such digit costs at most $s$ Teichm\"uller elements, by the
preceding argument. Hence the cost is $O((s^{1/v})^e)$, with an
absolute constant for these two choices of $p$. Since $s<p^v$,
the geometric estimate applies.

For the finitely many remaining $q$ put $h=q-1$, $f=[F:\mathbb F_p]$,
and suppose $q\notin\{2,3,4\}$. One has $f<\varphi(h)$.
At $p=2$ this was shown following Theorem~\ref{thm:cyclotomic-local}.
At $p=3$ and $f=2$ one has $\varphi(8)=4$.
For $f\ge3$, $\varphi(n)\ge\sqrt{n/2}$ gives
$\varphi(3^f-1)>f$. The elementary inequality follows prime power
by prime power from $\varphi(\ell^j)^2/\ell^j\ge1$ for odd primes
and $\varphi(2^j)^2/2^j\ge1/2$.
The ideal-covering argument in Theorem~\ref{thm:cyclotomic-local}
bounds the diameter of $B/p^eB$ for the Teichm\"uller generators by
$O_q(p^{ef/\varphi(h)})$, since the local cyclotomic residue degree
is $f$. This gives $D=p^{f/\varphi(h)}<p$. Taking the maximum of
the constants over the finite set of smaller residue fields proves
the theorem. At places whose residue characteristic differs from $p$,
the previously proved width bound is $p^a+1$.
\end{proof}

At $p=5$ the two-digit coefficient cost admits a further saving.
For every finite field $F$ of characteristic five, each element of
$W_2(F)$ is a sum of at most three Teichm\"uller constants. To see
this, write $u=x+y+z$ and let $e_2,e_3$ be the other elementary
symmetric functions in $x,y,z$. The second Witt coordinate of their
Teichm\"uller sum is
\[
               -u^3 e_2+u^2e_3+u e_2^2-e_2e_3 .
\]
This follows from the fifth Newton sum and the second ghost polynomial.
At $u=0$ the form $-e_2e_3$ has a nonzero value at $(1,1,-2)$.
Homogeneity of degree five and bijectivity of the fifth power map
give every value. For $u\ne0$ scale to $u=1$, where the second
coordinate becomes $(e_2-1)(e_2-e_3)$.

For a target $v\notin\mathbb F_5$ the corresponding affine plane
quintic, with $z=1-x-y$, is smooth even over an algebraic closure.
Its two partial derivatives are $x^4-z^4$ and $y^4-z^4$, as follows
by differentiating the fifth ghost equation over the integers before
reduction. If they vanish, $x,y,z$ are all multiples of a common
factor by fourth roots of unity, and $x+y+z=1$ forces all three
to lie in $\mathbb F_5$. At infinity the binary quintic is, up to sign,
$xy(x+y)(x^2+xy+y^2)$, with five distinct projective roots.
The smooth projective plane quintic is geometrically irreducible by
B\'ezout's theorem and has genus six. For $q=|F|\ge625$ the
Hasse--Weil bound gives at least $q+1-12\sqrt q-5>0$ affine points.
For $v\in\mathbb F_5$ the solutions over the prime field suffice.
The small fields are verified by the following finite counts.
In the row for $q$, the entry is the minimum, over $v\in\mathbb F_q$,
of the number of ordered pairs $(x,y)$ at $u=1$.
\[
\begin{array}{c|ccc}q&5&25&125\\\hline \text{minimum}&1&6&76\end{array}
\]
One may compute each row by initializing $q$ counters to zero.
For each $x,y$ put $z=1-x-y$, $e_2=xy+yz+zx$, $e_3=xyz$,
and increment counter $(e_2-1)(e_2-e_3)$. For $q=25,125$ the
respective defining polynomials over $\mathbb F_5$ are $X^2+2$ and
$X^3+X+1$. This proves the two-digit assertion.

Choose a fourth root of unity $i\in\mathbb Z_5$ for which $2+i$
has valuation one. Indeed $(2+i)(2-i)=5$ and exactly one factor
is a unit. Its square $3+4i$ costs seven Teichm\"uller constants.
Expanding in base $(2+i)^2$, a length-$s$ coefficient therefore
costs at most $B'_0=0$, $B'_1=1$, $B'_s=3+7B'_{s-2}$ for $s\ge2$.
Explicitly $B'_{2j}=(7^j-1)/2$ and $B'_{2j+1}=(3\cdot7^j-1)/2$,
so this cost has nonnegative second differences. The ambient
multiplicities above give
\[
             C'_a=B'_{a+1}+2\sum_{s=1}^a 5^{a-s}B'_s,
             \qquad C'_{a+1}=5C'_a+B'_{a+2}-3B'_{a+1}
                            \le5C'_a+1.
\]
Here $C'_4=621=5^4-4$, whence $C'_a\le5^a-a$ for $a\ge4$.
Adding the chart variables yields the field-uniform bound
\begin{equation}\label{eq:hteich-five}
                         H_{\mathbb K}(5^a)\le 5^a+1
                             \qquad(a\ge4).
\end{equation}
At the places over five this calculation gives, more precisely,
$\limsup_{a\to\infty}\sup_A h_{5^a}(A)/5^a\le35/36$.
Indeed $2\sum_{s\ge1}B'_s/5^s=35/36$ and
$B'_{a+1}/5^a\to0$.

We finish the finite cases $p=5$, $1\le a\le3$. First suppose $F=\mathbb F_5$.
The chart used above has $N_{a-1}=N_a$. Put $N_{-1}=0$. For any positive integer $d$ with
$d\ge N_{a-2}$ and $2d\ge N_a$, changing the reserved unit centre
$1$ to $1+\eta$ for $\eta\in T^dR$ gives the same change $\eta$
in both highest coordinates modulo $I_N$ and zero in every other
coordinate. This follows from the linear term at $1$ and the
divisibility by $T^{2d}$ of every nonlinear term.
These diagonal changes form a $\mathbb Z_5$-submodule $D$
of length $N_a-d$ if $d\le N_a$, because Witt addition on them is
coordinatewise addition. Its carries vanish to order at least $2d$. Compute the positive coefficient cost first modulo $D$.
The passage from $M_N$ to $M_N/D$ lowers its upper bound by at
least $N_a-d$, since $B'_1=1$, the higher differences are nonnegative,
and the lengths of the $5^s$-multiple modules cannot increase.
After the baseline is subtracted, represent the target modulo
$I_N$ and $D$, and correct the residual $D$ by changing the reserved
unit as just described. The original chart for $I_N$ applies at
the new baseline because its unit perturbation has order $N_a$ and
its linear coefficients $(1+\eta)^{5^j-1}$ are $1$ modulo $T$,
so all the normalized linear and nonlinear estimates persist.
For $a=1,2,3$ the respective values of $(C'_a,N_a,d)$ may be taken as
$(5,2,1),(26,7,4),(124,32,16)$. The resulting counts
$C'_a-(N_a-d)+a+1$ are at most $5^a+1$.

If $[F:\mathbb F_5]\ge2$, take two constant unit centres $1,\xi$
with $\bar\xi\notin\mathbb F_5$. Their two columns for the highest
two rows have unit determinant by the Moore criterion. For the
$b=a-1$ low rows use the nonunit centres and weights of
Theorem~\ref{thm:hhybrid}. If $b>0$, take the two high weights
to be $K_0+1$ and use unit perturbations of that order. The
nonunit leakage \eqref{eq:hhyleak} exceeds $K_0+1$, and the
nonlinear verification is exactly the one for the two blocks
above. If $b=0$, use unit perturbations of order one.
For $a=1,2,3$ this gives weight lists $(1,1)$, $(2,3,3)$,
$(3,7,8,8)$ respectively. The length formula gives ambient cyclic
multiplicities, from shortest to longest,
$(0,1)$, $(5,0,1)$, $(16,3,0,1)$. Using $B'_1=1,B'_2=3,
B'_3=10,B'_4=24$, the total positive counts including the chart
are respectively $5,18,53$. This proves
\eqref{eq:hteich-five} for every $a\ge1$.

The coefficient argument also works for $k=mp^a$ with fixed $m$
prime to $p$. For example, fix $h_0>1$ prime to $p$ with
$\operatorname{ord}_{h_0}(p)<\varphi(h_0)$. If $(F^\times)^m$
contains $\mu_{h_0}$, then
\[
                     h_{mp^a}(A)\ll_{p,m,h_0} mp^a
\]
with no restriction on either local degree. To see this, the
cyclotomic estimate represents every rational integer modulo $p^e$
by $O_{h_0}(D^e)$ constants from $\mu_{h_0}$, where
$D=p^{\operatorname{ord}_{h_0}(p)/\varphi(h_0)}<p$.
In the coefficient ring $W(\mathbb F_p((F^\times)^m))$ expand
in $p$-adic digits by lifting residues to sums of at most $L$
permitted Teichm\"uller constants. Represent the weight $p^j$ of
each digit modulo $p^e$ by the preceding cyclotomic bound. Thus
the coefficient cost is $O_{h_0}(eLD^e)$, with the extra factor $e$
absorbed in $D_1^e$ for any fixed $D<D_1<p$. Restriction of
scalars from $B$ costs $t$, with $tL\le m$. Apply the same
nonunit chart with $N_\Sigma\ll r+mp^{a-1}$ to finish.

More generally, in the residue field generated by
$H=(F^\times)^m$, some positive sum of at most $3L$ permitted
Teichm\"uller constants has valuation exactly one. Choose a
lift $g(x)$ of each residue as a sum of at most $L$ such constants,
with $g(0)=0$ and $g(1)=1$. If every zero-sum residue word of
length at most $3L$ vanished modulo $p^2$, then subtracting the
two words $g(x)+g(y)+g(-x-y)$ and $g(x+y)+g(-x-y)$ would give
an additive map modulo $p^2$ taking $1$ to $1$, which is impossible.
Thus one can use the same geometric cost whenever $3L<p$.

\begin{theorem}\label{thm:hlargeimage}
Fix real numbers $\theta>0$ and $\varepsilon>0$. There are constants
$C=C(\theta,\varepsilon)$ and $k_0=k_0(\theta,\varepsilon)$ as follows.
Let $A$ have residue field of cardinality $q$ in characteristic $p$,
let $k=mp^a\ge k_0$ with $a\ge1$ and $(m,p)=1$, and set
$d_0=(m,q-1)$. If
\[
                          q>C(1+d_0)^{2+\theta},
\]
then $h_k(A)\le\varepsilon k$. In particular, a family satisfying
this size condition at all places dividing $k$ has
$H_{\mathbb K}(k)\le k+1$ for sufficiently large $k$, on taking
$\varepsilon=1$.
\end{theorem}
\begin{proof}
First choose a fixed integer $s>2$ such that
$2s/(s-2)<2+\theta$. The complete finite-field sum for $X^{d_0}$
at a nonprincipal character has absolute value at most
$(d_0-1)\sqrt q$ by the Weil bound. Orthogonality gives
$q^{s-1}+O(d_0^s q^{s/2})$ representations as $s$ powers of
each residue, counting choices of the roots. Our size condition
makes the main term larger than the error, so $t=1$, $L\le s$.
Here the $m$th and $d_0$th power maps have the same image.
For $p>6s$ the valuation-one relation proved above costs at most
$3s$ Teichm\"uller elements. The cost of a coefficient modulo $p^e$
is $O_s((3s)^e)$ by expansion in this uniformizer, since the
residue digit costs at most $s$. The geometric estimate and the
chart with centres $T^{mi}$ for $0\le i\le a$ give
\[
               h_k(A)\ll_s a+1+mp^{a-1}+(3s)^a.
\]
Thus $h_k(A)/k\le\varepsilon$ for all sufficiently large $p$,
uniformly in $m,a$.

It remains to handle $p$ in a fixed finite set. Take $v$ so large
that $s'=p^v-1>2(v+1)$ and $2s'/(s'-2(v+1))<2+\theta$.
For any nonprincipal additive character of $B/p^{v+1}B$, the bound
of \cite{hkhc} applied to a unit times $X^{d_0}$ at its conductor
has size at most $d_0p^v\sqrt q$. Nondegeneracy follows already
modulo $p$, since the positive degree $d_0$ is prime to $p$.
Increase $C$ to make the orthogonality main term $q^{s'-v-1}$
larger than the error $(d_0p^v)^{s'}q^{s'/2}$.
Every class modulo $p^{v+1}$ now costs at most $s'$ permitted
Teichm\"uller elements. The preceding block argument again bounds
coefficient cost by $O(D^e)$ for $D=(s')^{1/v}<p$. All the constants
are uniform over this finite set of primes.

We may now use any fixed positive integer $c$. For each such $c$
increase the residue cutoff $C$ on this finite set of primes to ensure
$[F:\mathbb F_p]\ge c$. Since $t=1$, the permitted constants span
$F$ and contain $\min(c,a+1)$ independent elements over $\mathbb F_p$.
If $a+1\le c$, use only these constant centres, whose full Moore
matrix is invertible, giving $h_k(A)\ll_c 1$.
If $b=a+1-c>0$, use $c$ such constants for the top block in the
proof of Theorem~\ref{thm:hhybrid}. Their high-row Jacobian is a
unit, since the Moore determinant is unchanged up to Frobenius
and column units. With all high thresholds $K_0+1$ the same
low-row construction works. In the limiting case in which
\eqref{eq:hhyleak} equals $K_0+1$, the lower-left block may be
nonzero. Block triangularity gives invertibility of the normalized
matrix, and the nonlinear terms have strictly larger order.
Thus $N_\Sigma\ll c m p^{a-c}$. The digit estimate yields
\[
             h_k(A)\ll a+1+c m p^{a-c}+D^a ,
\]
with an implied constant independent of $c$. Choose $c$ large
enough that the term $c p^{-c}$ in the ratio to $k$ contributes
less than $\varepsilon/2$ for every remaining prime, then take the
residue cutoff required by this $c$. Finally choose $k_0$. The first and
last terms divided by $k$ tend uniformly to zero by $D/p<1$
for large $a$ and $m=k/p^a$ for bounded $a$. This proves the
assertion.
\end{proof}

\subsection{Exponents with a dominant prime factor}
\label{sec:local-dominant-prime}
The short uniformizer argument also applies when the prime-to-$p$
cofactor varies.
\begin{corollary}\label{cor:local-dominant-prime}
Uniformly over all number fields, as $k=mp^a\to\infty$ through exponents
with $p>4m$ and $a\ge1$, one has
\[
              H_{\mathbb K}(k)\le (2+o(1))k.
\]
With $p>5m$ the leading constant is $1$. There is no restriction
on ramification.
\end{corollary}
\begin{proof}
At a place above $p$, use the coefficient field $E$ and residue
width $L$ of \eqref{eq:nf-residue-width}. There is a positive sum
$\alpha$ of $c\le3L$ permitted Teichm\"uller constants having valuation
one in $W(E)$. Indeed choose a lift $s(x)$ of every $x\in E$ using
at most $L$ constants, taking $s(0)=0$. If every sum of at most
$3L$ constants with zero residue vanished modulo $p^2$, the two
sums
\[
  s(x)+s(y)+s(-x-y),\qquad s(x+y)+s(-x-y)
\]
would make $s$ an additive section modulo $p^2$. This contradicts
$p\,s(1)\ne0$ modulo $p^2$. We necessarily have $c\ge2$.
Each coefficient modulo $p^e$ can thus be represented in base
$\alpha$ at cost at most $L(c^e-1)/(c-1)$, as in
Section~\ref{sec:local-teich}.

Use the nonunit chart with centres $T^{mi}$, $0\le i\le a$, from
the proof of Theorem~\ref{thm:hhybrid}. The construction there of
these centres and their triangular linear terms applies for every $m$;
the restriction $m<64$ was used only for its numerical comparison.
Its thresholds obey $N_j(m)\le mN_j(1)$, hence the length of each
$p^s$-multiple of its finite Witt quotient is at most $m$ times
the corresponding length for $m=1$. Restriction to $W(E)$ costs
the further factor $t$. With $r=a+1$ and $c<p$ the convex digit
cost and the multiplicities at $m=1$ therefore give
\[
 \frac{h_k(A)-r}{k}\le
 \frac{tL}{p^a}\left\{\frac{c^{a+1}-1}{c-1}
   +\frac{2}{c-1}\sum_{s=1}^a p^{a-s}(c^s-1)\right\}.
\]
This use of multiplicities only needs the bounds on the lengths of
the $p^s$-multiples: the function $(c^e-1)/(c-1)$ is a nonnegative
linear combination of the functions $(e-s)_+$, $s\ge0$.
Since $tL\le m$ and $c\le3m$, the right side is at most
\[
        2m(3m/p)^a+\frac{2mp}{(p-3m)(p-1)}.
\]
If the local width at a place under consideration exceeds
$(1+\varepsilon)k$, the proof of Theorem~\ref{thm:local-sparse}
there gives $p\ll_\varepsilon(a+1)^2$. Thus under $p>4m$ the
first term tends to zero on such a sequence, where necessarily
$a\to\infty$. The second term is at most $2$ because $p\ge4m+1$,
and is at most $1$ when $p\ge5m+1$.

Finally consider a different prime $\ell\mid m$ and write
$j=v_\ell(m)$. A width greater than $(1+\varepsilon)k$ at this
place would imply $k/\ell^j\ll_\varepsilon(j+1)^2\log(2j+2)$.
But $k/\ell^j\ge p^a>4m$ and $j\le\log_2 m$, so this is
impossible for unbounded $m$; for bounded $m$ the left side grows
with $k$. Good places have width at most $k+1$. Letting
$\varepsilon$ decrease proves the assertion.
\end{proof}

\subsection{A finite list of residue fields at a fixed prime-to-characteristic part}
\label{sec:inert-residue-test}
Combining the large residue image estimate with the distribution of
Teichm\"uller constants in the coefficient ring leaves a finite list
of small residue fields at which ramification could cause superlinear
growth.
\begin{corollary}\label{cor:inert-residue-test}
Fix a rational prime $p$ and $m\ge1$ prime to $p$. Write
$q=p^f$, $h=(q-1)/(m,q-1)$. With the convention that $h=1$ passes
the following test, retain $q$ only when $p$ generates
$(\mathbb Z/h\mathbb Z)^\times$.
For every omitted $q$ one has
\[
            h_{mp^a}(A)\ll_{p,m} mp^a
\]
uniformly in $a$, the residue degree and the ramification index
of the local field with residue cardinality $q$. Every retained
$q$ satisfies $q\le3mf+1$; for $p=2$ it satisfies
$q\le3mf/2+1$ except possibly at $h=1$, where $q\le m+1$.
In particular the test generates a finite list by first solving
$p^f\le 3mf+1$. A uniform all-field linear bound along $mp^a$
reduces to the places on this list.
\end{corollary}
\begin{proof}
By Theorem~\ref{thm:hlargeimage}, all sufficiently large $q$ have
local width $\ll_{p,m} mp^a$, uniformly in $a$. Its finitely many
small exponents are covered independently of $q$ by
Theorem~\ref{thm:hhybrid}. At a remaining omitted $q$ one has
$h>1$ and $f_0=\operatorname{ord}_h(p)<\varphi(h)$.
In the coefficient ring $W(\mathbb F_p(\mu_h))$, the directed
diameter modulo $p^j$ for its Teichm\"uller generators $\mu_h$ is
$O_h(D^j)$, $D=p^{f_0/\varphi(h)}\le\sqrt p$. This is the
ideal covering bound in the proof of Theorem~\ref{thm:cyclotomic-local},
applied directly to the prime ideal raised to the power $j$.
Use the chart with nonunit centres of Theorem~\ref{thm:hhybrid}
and thresholds $N_i(m)\le mN_i(1)$. Its sum of thresholds is
$O(a+1+mp^{a-1})$. Summing $D^j\ell(p^jM)$ as in
Theorem~\ref{thm:hteich} now gives
$O_{p,m,q}(a+1+mp^{a-1}+D^a)$, for the full cost including the
chart. Restricting the coefficient scalars costs the fixed residue
degree ratio $f/f_0$, and the extra dyadic layer shifts the length
bound by one. Taking a maximum over the finite range of $q$
gives the asserted uniformity.

For a retained $h>1$ its unit group is cyclic, so
$h=2,4,\ell^b$ or $2\ell^b$ for an odd prime $\ell$. On this
list $\varphi(h)\ge h/3$; if $p=2$ then $h$ is odd and
$\varphi(h)\ge 2h/3$. Since $\varphi(h)=f_0\le f$ and
$q-1\le mh$, the stated cutoffs follow. At other wild primes of
$mp^a$ the multiplicities divide the fixed exponent pattern of $m$.
The sharper, local version of the estimate in
Theorem~\ref{thm:hhybrid} therefore bounds their widths by
$O_m(mp^a)$; at the tame places the bound is $mp^a+1$.
\end{proof}
For example, for exponents $3\cdot2^a$ the test leaves precisely
$q=2,4,16$ at the dyadic places. For $2\cdot3^a$ it leaves
$q=3,9$ at the triadic places. Membership records the absence
of a saving from these coefficient directions; the actual local
width may still be linear there.

The residue test and the tame finite-field estimate also allow a
larger prime-to-$p$ cofactor in the dominant-prime estimate of
Corollary~\ref{cor:local-dominant-prime}.

\begin{corollary}\label{cor:wild-large-prime}
Fix $\eta>0$. Suppose $p\to\infty$ through the primes and
$k=mp^a$, $a\ge1$, $(m,p)=1$, with
\[
             m\le p^2/(\log p)^{2+\eta},
                         \qquad (p-1)\nmid 2m .
\]
Then at all places above $p$ one has $h_k(A)=o(k)$, uniformly in
the residue field, the ramification, $a$ and $m$. Over all number
fields one has $H_{\mathbb K}(k)\le(1+o(1))k$ in this range.
Over the full family $p>2m+1$, including bounded $p$, one has
$H_{\mathbb K}(k)\ll k$ with an absolute constant.
\end{corollary}
\begin{proof}
Let $H=(\mathbb F_q^\times)^m$, $E=\mathbb F_p(H)$, $t=[\mathbb F_q:E]$
and $L$ be the residue width in $E$, as in \eqref{eq:nf-residue-width}.
If $t>1$, then $|E|^{t-1}\le (q-1)/(|E|-1)\le m<p^2$,
so $|E|=p$ and $t=2$. When $q=p$ the assumption on $m$ gives
$|H|>2$. At $q\ge p^2$ the inequality $m=o(p^2)$ gives the same
conclusion. In $E$ this is the group of $b$th powers for
$b=(|E|-1)/|H|\le m$. Apply Lemma~\ref{lem:uniform-tame-width}
to an unramified lift of $E$, or use $L=1$ when $b=1$.
It gives $tL\ll\sqrt m\log(2m)=o(p)$ uniformly.

Set $r=a+1$. In the range $a\le\sqrt p$, the uniform estimate
in the proof of Theorem~\ref{thm:local-sparse} gives
\[
             h_k(A)/k\le r/k+(tL/m)(1+Cr^2/p).
\]
This tends to zero if $m\to\infty$. In the range $a>\sqrt p$
the coefficient argument of Corollary~\ref{cor:local-dominant-prime}
applies since $3L=o(p)$. Before bounding $tL$ by $m$ and $3L$ by
$3m$, its cost gives
\[
 \frac{h_k(A)}k\le\frac r k+
       2tL(3L/p)^a+\frac{2tL p}{(p-3L)(p-1)} .
\]
The right side tends uniformly to zero since $tL\le m<p^2$.
For bounded $m$ the same estimate tends to zero for every $a\ge1$.
This proves the uniform limit at $p$ by the subsequence criterion.

At another wild prime $\ell$ put $j=v_\ell(m)$. A width greater
than $(1+\epsilon)k$ would require
$k/\ell^j\ll_\epsilon(j+1)^2\log(2j+2)$ by the proof of
Theorem~\ref{thm:local-sparse}. But $k/\ell^j\ge p$ and
$j\le\log_2m\ll\log p$. The tame bound is $k+1$. This proves the
second assertion. For $p>2m+1$ the growing primes satisfy the
hypotheses just proved for, say, $\eta=1$. Only a bounded range of
$p,m$ remains for an absolute upper bound. For each of these pairs
Corollary~\ref{cor:inert-residue-test} has no retained $q$. Indeed
$f=1$ would require $h\le2$. For $f\ge2$ one would have
$p^2-1\le 6m$, since $(p^f-1)/f$ is increasing for $f\ge2$.
Both conclusions contradict $p>2m+1$. Taking the maximum of the
finitely many remaining constants proves the last assertion.
\end{proof}

\begin{corollary}\label{cor:absolute-sparse-count}
For fixed $t>1$, the number of $2\le k\le X$ at which some number
field has $H_{\mathbb K}(k)>t k$ is
\[
           O_t\!\left(\frac{(\log X)^4}{(\log\log X)^{7/2}}\right)
                       \qquad(X\longrightarrow\infty).
\]
Consequently, outside a set of at most this order, the linear bound
$G_{\mathbb K}(k)\le\lfloor C_1 k\rfloor$ holds for every number field.
\end{corollary}
\begin{proof}
For large $k$ only wild places can contribute. Write $k=mp^a$,
$r=a+1$. The proof of Theorem~\ref{thm:local-sparse} gives
$p\ll_t r^2$, $m\le M=C_t r^2\log(2r)$. We show that for each
$a$ there are $O_t(r^3\sqrt{\log(2r)})$ possible pairs $(p,m)$.
For $p\le Y=2\sqrt M$, the elementary prime-counting bound gives
$O(MY/\log(2Y))$ possibilities. For $p>Y$, one has $m<p^2/4$,
so the coefficient-field index is at most two by the argument of
Corollary~\ref{cor:wild-large-prime}. The exceptional residue group
of size at most two can therefore occur only at $q=p$, and then
$p-1\mid 2m$. The number of such pairs is
\[
 \ll_t\sum_{p\ll_t r^2}\{1+2M/(p-1)\}
       \ll_t r^2+M\log\log(3r).
\]
Outside these cases, if $L$ is the residue width, then
$[\mathbb F_q:E]L\ll\sqrt m\log(2m)$ as before. Ignore bounded $m$,
at a cost $O_t(r^2)$. Enlarging the bound defining this finite range
makes $[\mathbb F_q:E]L/m$ as small as desired. Also $r/k$ is
uniformly small for large $k$. Thus the uniform estimate
in the proof of Theorem~\ref{thm:local-sparse} and the inequality
$h_k(A)>tk$ imply $m\ll_t r^4(\log(2r))^2/p^2$.
By partial summation of the prime-counting bound,
$\sum_{p>Y}p^{-2}\ll1/(Y\log(2Y))$. This gives the claimed
per-exponent count.

Put $A=\log X$. The primes $p\le\sqrt A$ contribute
$O_t(A^{7/2}\log(2A))$ by the bounds $a\ll A$ and $m\le M$.
At larger primes one has $a\le2A/\log A$. Summing the per-exponent
bound to this cutoff gives $O_t(A^4/(\log A)^{7/2})$, as asserted.
Use, for example, $t=2$ and Theorem~\ref{thm:nf-all-linear} for the
global conclusion.
\end{proof}

\subsection{Nonconstant carries at small residue fields}
\label{sec:local-small-carries}
The exceptions $\mathbb F_2,\mathbb F_3,\mathbb F_4$ in
Theorem~\ref{thm:hteich} raise a question about uniformity over
arbitrarily ramified fields. Two short identities for actual
Teichm\"uller generators are relevant. Write $V$ for Verschiebung.

\begin{proposition}\label{prop:small-residue-carries}
\begin{enumerate}
\item For $R=\mathbb F_3[[T]]$, every element of $V^2(T^3R)$ in
$W_3(R)$ is a sum of eight Teichm\"uller vectors.
\item For $R=\mathbb F_4[[T]]$, each $V^3(v^2)$ with $v\in TR$ is
a sum of seven Teichm\"uller vectors in $W_4(R)$.
\item For $R=\mathbb F_2[[T]]$, a positive sum of fewer than $2^j$
Teichm\"uller vectors with nonzero entries cannot have zero
truncation in $W_j(R)$, for any $j\ge1$.
\end{enumerate}
\end{proposition}
\begin{proof}
In residue characteristic three, consider the Witt vector
\[
 [T^2s+1]+[T^2s-1]-2[T^2s]-[s+T]-[s-T]+2[s].
\]
Here $-[x]=[-x]$. Over $\mathbb Z[T,s]$ the first and third
power sums vanish identically. The ninth power sum divided by nine is
\[
 2\{(T^2-T^8)s+14(T^{10}-T^4)s^5+4(T^{14}-T^2)s^7\}.
\]
It is the last Witt coordinate of the displayed vector. After division
by $T^2$ the polynomial in $s$ has zero constant term and unit linear
coefficient in $R$. It maps $TR$ onto $TR$ by successive approximation
in $T$.

For the second identity lift a primitive cube root of unity $\zeta$
to $\mathbb Z_2[\zeta]$. Let $w$ be a parameter, put $y=w$, and
define $x,z$ by
\[
                   x-\zeta z=1,\qquad x-\zeta^2 z=w^2.
\]
Since $\zeta-\zeta^2$ is a unit, $x,z$ are integral polynomials in
$w^2$. The triple $(x,z,-x-z)$ has sum zero and second elementary
symmetric function $-(x^2+xz+z^2)=-y^2$. The pair $(y,-y)$ has
the same two elementary symmetric functions. Their power sums
of orders one, two and four agree. If $e_3=-xz(x+z)$, Newton's
identities give for the eighth power sum of the triple minus that of
the pair the value $8y^2 e_3^2$. Hence the difference of their
Teichm\"uller sums is $V^3((w e_3)^2)$ in $W_4(R)$. Each negative
Teichm\"uller vector $-[u]$ is the sum $[\zeta u]+[\zeta^2u]$:
this identity holds on every $2$-typical ghost component before
reduction. This gives seven positive terms. The constant term of
$e_3$ is a unit, since $x,z,-x-z$ at $w=0$ are all units.
Thus $w\mapsto w e_3(w^2)$, writing $e_3$ as a polynomial in $w^2$,
maps $TR$ onto $TR$.

For the last assertion remove the smallest power of $T$ common to
the roots. Multiplication by its Teichm\"uller lift is injective on
$W_j(R)$: on coordinate $i$ it multiplies by the $2^i$th power
of that common factor. At least one of the remaining roots is a
unit. Evaluating at $T=0$ gives precisely the number of unit roots
in $W_j(\mathbb F_2)=\mathbb Z/2^j\mathbb Z$. A positive number
less than $2^j$ cannot vanish there.
\end{proof}

The first two identities save one summand over the cost of the
corresponding rational carries. Their leading coordinates have
positive $T$-order. Iterating such a carry would have to control
the resulting loss of low $T$-coefficients at each Witt depth.
In residue field $\mathbb F_2$ the third assertion rules out a
strict saving of this particular kind at every depth.

\subsection{The top Witt coordinate modulo Frobenius}
\label{sec:witt-top-functions}
For prime-power exponents, a finite additive space of polynomial
representatives gives both a carry at arbitrary Witt depth and an
obstruction below its precision. Let $R=\mathbb F_q[[T]]$, $q=p^f$,
and let $G_{a+1}$ denote the additive subgroup generated by the
Teichm\"uller vectors $[x]$, $x\in R$, in $W_{a+1}(R)$.
Write $V^a(y)=(0,\ldots,0,y)$ and put
\[
 N(a,f)=1+\sum_{j=1}^{\lceil a/f\rceil-1}p^{a-fj}
                         \qquad(a\ge1).
\]
\begin{proposition}\label{prop:witt-top-functions}
For $Y_a=\{y\in R: V^a(y)\in G_{a+1}\}$ one has the exact equality
\[
                Y_a+R^p=T^{N(a,f)}R+R^p .
\]
Furthermore every $V^a(y)$ with $y\in T^{N(a,f)}R$ is a sum of
$p^a$ Teichm\"uller vectors minus one fixed sum of $p^a$ such vectors.
For $q=p$ the first exponent prime to $p$ that can occur in a nonzero
top coordinate is $(p^a-1)/(p-1)$.
\end{proposition}
\begin{proof}
Set $n=\lceil a/f\rceil$. Let $A$ be the $a$-dimensional
$\mathbb F_p$-space of polynomials of degree less than $n$ whose
coefficients below degree $n-1$ are arbitrary in $\mathbb F_q$,
and whose coefficient of $T^{n-1}$ belongs to a fixed subspace of
dimension $a-f(n-1)$. Put $P(X)=\prod_{v\in A}(X-v)$. This
polynomial is additive. Indeed adjoining an independent vector $u$
to a space with additive polynomial $P_0$ changes its polynomial to
$P_0(X)^p-P_0(u)^{p-1}P_0(X)$. In particular
$P(X)=X^{p^a}+\sum_{j<a}c_jX^{p^j}$, with $j$ starting at zero.

We first compute a centered sum. For an indeterminate $U$,
\[
 \frac{\prod_{v\in A}(1-(x+v)U)}{\prod_{v\in A}(1-vU)}
 \equiv 1-P(x)U^{p^a}\pmod{U^{p^a+1}},
\]
by applying additivity to $U^{p^a}P(U^{-1}-x)$.
In the unique factorization of a series of constant term one as
$\prod_{j\ge1}(1-w_jU^j)$, the coordinates with $j=p^i$ give
the corresponding $p$-typical Witt vector for multiplication of
series. This follows directly by taking the formal logarithmic
derivative over the universal integer coefficient ring: its
coefficient of degree $p^i$ in $-UF'/F$ is the $i$th $p$-typical
ghost polynomial, since every divisor of $p^i$ is a power of $p$.
Thus
\[
                 \sum_{v\in A}[x+v]-\sum_{v\in A}[v]=V^a(P(x)).
\]
Among nonzero $v\in A$ the sum of their $T$-valuations equals
$\sum_{j=1}^{n-1}(p^{a-fj}-1)=N(a,f)-n$.
Consequently $T^{-N(a,f)}P(T^n z)$ is a polynomial with unit
linear coefficient and all higher coefficients in $TR$: factor out
$T^n z\prod_{v\ne0}(-v)$ and use $v_T(v)<n$.
Successive approximation in $T$ gives
$P(T^nR)=T^{N(a,f)}R$, proving the stated construction.

For the converse use the additive differential of the $j$th Witt
ghost, divided over the universal coefficient ring by $p^j$ before
reducing modulo $p$. On a Witt vector $\mathbf z$ it is
\[
     D_j(\mathbf z)=\sum_{i=0}^j z_i^{p^{j-i}-1}\,dz_i
                       \quad\text{in }R\,dT ,
\]
where the differential is formal differentiation of power series.
On $[x]$ it has value $x^{p^j-1}dx$, and on $V^a(y)$ it is zero
for $j<a$ and is $dy$ for $j=a$. If
$V^a(y)=\sum_i m_i[x_i]$, $m_i\in\mathbb Z$, omit zero $x_i$
and combine the equations for $D_a+\sum_{j<a}c_jD_j$. This gives
\[
                         dy=\sum_i m_i P(x_i)\,dx_i/x_i .
\]
Each logarithmic differential on the right has at most a simple
pole at $T=0$. Every $P(x_i)$ has valuation at least $N(a,f)-1$:
for $1\le j<n$ precisely $p^{a-fj}$ members of $A$ agree with
$x_i$ modulo $T^j$. If $f\mid a$, a further member agrees modulo
$T^n$, giving the lower bound $N(a,f)$ instead. It follows that
all coefficients of $y$ of degree prime to $p$ below $N(a,f)$
vanish. Indeed in the first bound the only extra degree left to
consider is $N(a,f)-1$, which is divisible by $p$; for $N(a,f)=1$
there is nothing to check. Since $\mathbb F_q$ is perfect this gives
$y\in R^p+T^{N(a,f)}R$. Finally $N(a,f)\equiv1\pmod p$ and
$T^{N(a,f)}$ itself occurs by the construction.
\end{proof}

Over $\mathbb F_3[[T]]$ at depth two the threshold is four. The
eight-term carry in Proposition~\ref{prop:small-residue-carries}
already supplies $T^3R$; its extra degree three is a Frobenius
power, in agreement with the formula. At $p=2$ the proposition
describes the characteristic-two quotient of the coefficient ring
in the local width argument. Section~\ref{sec:dyadic-radical-width} gives a linear bound when the
residue field is \(\mathbb F_2\) and \(2\) itself is a \(2^a\)-th
power in the completion. Bounding the cost of these carries at
unrestricted dyadic completions remains open.

\subsection{The low coefficients of dyadic affine carries}
\label{sec:witt-affine-intervals}
Over $R=\mathbb F_2[[T]]$, one might vary the affine space in
Proposition~\ref{prop:witt-top-functions} to handle many Frobenius
coefficients at once. The resulting top coordinates have an exact
description. For $a\ge1$, give an $a$-tuple in $R$ its linear map
$u:\mathbb F_2^a\to R$, allowing dependent entries, and write
\[
 P_u(X)=\prod_{v\in\mathbb F_2^a}(X-u(v)),\qquad
 d_i=2^a-2^{a-i}\quad(0\le i\le a).
\]
Dependence simply counts the elements of the image with multiplicity.
\begin{proposition}\label{prop:witt-affine-intervals}
As $u$ and $z\in R$ vary, the polynomials $P_u(z)\bmod T^{2^a}$
are exactly zero and the sums
\[
                       \sum_{i=r}^{s-1} T^{d_i},
                 \qquad 0\le r<s\le a+1 .
\]
The least number of such polynomials needed to make a prescribed
subset of the $a+1$ monomials is its number of runs of consecutive
indices. In particular some attainable top coordinates require at
least $\lceil(a+1)/2\rceil$ affine carries of this kind.
\end{proposition}
\begin{proof}
For any $b$-tuple $w$, let $C_b(w)$ be the coefficient of
$X^{2^{b-1}}$ in its additive polynomial $P_w(X)$, $b\ge1$.
Let $r(w)$ be the largest $r$ such that the image subspace
surjects onto $R/T^rR$. We claim
\begin{equation}\label{eq:witt-affine-prefix}
 C_b(w)\equiv
 \sum_{i=0}^{r(w)-1} T^{2^{b-1}-2^{b-1-i}}
                 \pmod{T^{2^{b-1}}}.
\end{equation}
Choose a row-echelon basis $w_0,\ldots,w_{h-1}$ of the image with
strictly increasing $T$-orders, and append $b-h$ zero entries.
This change preserves the polynomial, including multiplicities.
Writing $u_i$ for the initial $i$ entries and applying
$P_{u_{i+1}}=P_{u_i}^2+P_{u_i}(w_i)P_{u_i}$ gives
\[
                 C_b(w)=\sum_{i=0}^{b-1}
                          P_{u_i}(w_i)^{2^{b-1-i}}.
\]
For $i<r(w)$ the first $i+1$ pivots have orders $0,\ldots,i$.
In the product defining $P_{u_i}(w_i)$
the zero vector contributes order $i$, and for each $0\le j<i$
there are $2^{i-1-j}$ other vectors with lowest pivot $j$.
The sum of the orders is $i+\sum_{j<i}j2^{i-1-j}=2^i-1$.
The leading coefficient is one. If $i\ge r(w)$ then either
$w_i=0$ or its order is at least $i+1$, while orders of all
earlier pivots are at least their indices. The product then vanishes
modulo $T^{2^i}$. These two alternatives give
\eqref{eq:witt-affine-prefix}.

Append $z$ to the tuple. The polynomial recursion gives
$P_u(z)=C_{a+1}(u,z)+C_a(u)^2$. The second prefix is contained
in the first because adjoining one entry preserves every surjective
projection. Their difference is the interval
$[r(u),r(u,z))$. Conversely, for $0\le r<s\le a+1$ take as
pivots all the powers $T^i$ for $0\le i<s$, $i\ne r$, append
zero entries to reach dimension $a$, and take $z=T^r$.
This gives the required two prefix lengths, including $s=a+1$.
The empty interval is obtained from $z=0$.
A sum of $h$ intervals has at most $2h$ endpoints, counting the
two outer boundaries; conversely a binary word is the sum of its
disjoint runs. The centered Witt sum for each affine block is
$V^a(P_u(z))$, including for dependent tuples by the same
generating-series identity. Hence all these low patterns occur in
the additive group. The assertion follows.
\end{proof}

At odd $p$, the uncentered sum on an affine $a$-space is already
$V^a(P_u(z))$ in $W_{a+1}$: the generating product has the powers
$U^{p^a-p^j}$ and $U^{p^a}$ beyond the constant, and the former
have no $p$-power index and all exceed $p^{a-1}$. At $p=2$ its
only possible nonzero coordinates are $C_a(u)$ in position $a-1$
and $P_u(z)$ in position $a$. This follows by factoring the
generating polynomial into its big Witt coordinates as in
Proposition~\ref{prop:witt-top-functions}. Thus two positive blocks
with the same space and shifts $z,0$ yield
$V^a(P_u(z)+C_a(u)^2)$; modulo $T^{2^a}$ this is the single
prefix for $(u,z)$. These restrictions concern the affine construction.
The full local width allows arbitrary sums of individual powers.

\begin{corollary}\label{cor:witt-affine-odd}
Let $p\ge3$ and $R=\mathbb F_p[[T]]$.
Put $N_a=(p^a-1)/(p-1)$ and $N_0=0$. Modulo $T^{N_a+1}$, one
affine carry $P_u(z)$ is zero or a single nonzero monomial of degree
\[
                    n_i=(p^a-p^{a-i})/(p-1),\qquad 0\le i\le a .
\]
Every listed monomial and every nonzero scalar coefficient occur.
In particular some attainable top coordinate requires $a+1$ affine
blocks, even when the spaces, translations and signs are free.
\end{corollary}
\begin{proof}
For a $b$-tuple the coefficient $C_b$ of $X^{p^{b-1}}$ satisfies
$C_{b+1}=C_b^p-P_u(z)^{p-1}$. For an ordered echelon basis, the
successive products have orders $N_i$ for $i<r$ and orders at least
$N_i+1$ thereafter, where $r$ is the number of consecutive initial
pivots. Indeed the multiplicity of a first nonzero coefficient at
pivot $j$ is $(p-1)p^{i-1-j}$, and
$i+\sum_{j<i}j(p-1)p^{i-1-j}=N_i$. Thus modulo $T^{p^{b-1}}$,
$C_b$ is minus the prefix of monomials of degrees
$p^{b-1}-p^{b-1-i}$ with $i<r$, the leading coefficients becoming
one on taking $(p-1)$th powers. As above, $P_u(z)^{p-1}$ modulo
$T^{p^a}$ is the positive interval in degrees $(p-1)n_i$.
For the empty interval, $P_u(z)\equiv0\bmod T^{N_a+1}$.
For a nonempty interval starting at $r$, the leading degree of
$P_u(z)$ is $n_r$. Write $h=a-r$. If $h\ge1$, the next possible
relative exponent in the $(p-1)$th power is
$(p-1)p^{h-1}> (p^h-1)/(p-1)=N_a-n_r$. Since the power
$p-1$ is invertible in the residue field, the first nonconstant
coefficient of a unit survives on taking that power. No further
term in $P_u(z)$ survives the stated truncation. For $h=0$ only
the leading digit is needed. The monomial pivot choices in
Proposition~\ref{prop:witt-affine-intervals} again realize all
singletons, and scaling $z$ by $\mathbb F_p^\times$ scales
$P_u(z)$ by the same factor since $P_u$ is additive. Summing the
$a+1$ distinct monomials proves the last assertion.
\end{proof}

The odd-prime argument has a further application in the exceptional
dyadic residue field $\mathbb F_4$. In fact take $q=p^f>2$, replace
the dimension $a$ in Corollary~\ref{cor:witt-affine-odd} by $h$ and
replace $p$ there by $q$. Its conclusion holds for affine
$\mathbb F_q$-linear blocks of size $q^h$, in Witt depth $fh$.
The echelon count now has $(q-1)q^{i-1-j}$ vectors per first pivot;
the coefficient recursion is $C_{b+1}=C_b^q-P_u(z)^{q-1}$, and
the same strict gap uses $(q-1)^2>q$. The smallest nonconstant
index in the generating product of the block before its terminal
coefficient is $(q-1)q^{h-1}>p^{fh-1}$. Therefore its $p$-typical
coordinates vanish below depth $fh$ even when $p=2$. It follows,
in particular, that at least $h+1$ such positive blocks can be
needed over $\mathbb F_4$. General $\mathbb F_p$-linear blocks in
$\mathbb F_q[[T]]$ allow additional choices of residue directions.

\begin{lemma}\label{lem:uniform-tame-width}
Let $k\ge2$, let $A$ be a local integer ring of residue cardinality
$q=\ell^f$ with $\ell\nmid k$, and put $H=(\mathbb F_q^\times)^k$.
Unless $f=1$ and $|H|\le2$ one has, with an absolute implied constant,
\[
                       h_k(A)\ll\sqrt k\log(2k).
\]
For $f=1$ and $|H|=1$ the ordinary and primitive widths both equal
$\ell$. For $f=1$ and $|H|=2$ both equal $\max(2,(\ell-1)/2)$.
\end{lemma}
\begin{proof}
Let $E=\mathbb F_\ell(H)$, $g=[E:\mathbb F_\ell]$ and
$b=[E^\times:H]\le k$. The group $H$ consists of the nonzero $b$th
powers in $E$. Its additive span is $E$, since this span is closed
under multiplication and is a finite integral domain. The finite-field
Waring number $\gamma$ therefore exists. Nonzero targets are represented
primitively in at most $\gamma$ terms. For zero represent $-1$ first
and append $1$. Hensel lifting bounds $h_k(A)$ by $\gamma+1$,
since all residues obtained in this way are actual $k$th powers.

If $b=1$ two terms suffice. For $|H|>2$ and $g=1$, use
$\gamma\le83\sqrt b$ from \cite[Theorem~1.1]{hprimeCP}.
For $g\ge2$ and $b>1$, Winterhof \cite[Lemma~2 and Theorem~1]{htameWinterhof}
gives $\gamma\le1+32\log b$ if $|E|>b^2$, and otherwise
$\gamma<6.2g(2b)^{1/g}\log b$. His hypothesis $b<(|E|-1)/2$
holds since $|H|>2$. Here $g\le2\log_2 b$ in the latter case.
The function $x(2b)^{1/x}$ is convex for $x>0$; at the endpoints
of $[2,2\log_2 b]$ it is $O(\sqrt b)$. This gives the asserted bound.

If $|H|\le2$, then $E=\mathbb F_\ell$ and $q\le2k+1$. For $f\ge2$,
repeated copies of $1$ represent all residues primitively in at most
$\ell\le\sqrt{2k+1}$ terms. For $f=1$ and $H=\{1\}$ the residue
zero with a unit requires $\ell$ terms. A uniformizer cannot be
represented using only nonunit roots. It is eligible, since a unit
$k$th power represents $1+\pi$ by Hensel lifting. For $|H|=2$ the
possible residues are $0,1,-1$. Their diameter is $(\ell-1)/2$;
the primitive zero costs two, and valuation one excludes a single
power. These congruences prove both exact assertions.
\end{proof}

\subsection{How often a local width is a fixed fraction of the degree}
\label{sec:local-frequency}
For $0<t\le1$ set
\[
 A(t)=\sum_{1/j\ge t}\frac1j
                 +2\sum_{\substack{j\text{ odd}\\1/j>t}}\frac1j
                 \qquad(j\in\mathbb N,\ j\ge1).
\]
For a number field $\mathbb K$ let $\delta_{\mathbb K}$ be the density
of rational primes admitting a degree-one place. This is the fraction
of elements fixing a point in the permutation action of the Galois
closure on the embeddings of $\mathbb K$, by Chebotarev's theorem \cite{ednark}.

\begin{theorem}\label{thm:local-frequency}
For each $0<t\le1$ one has
\begin{align*}
 \#\{k\le X:\sup_{\mathbb K}H_{\mathbb K}(k)>tk\}
       &\sim A(t)X/\log X,\\
 \#\{k\le X:H_{\mathbb K}(k)>tk\}
       &\sim \delta_{\mathbb K}A(t)X/\log X
                      \qquad(\mathbb K\text{ fixed}).
\end{align*}
In particular $\sup_{\mathbb K}H_{\mathbb K}(k)/k$ tends to zero
in natural density. The same formulas hold for $W_{\mathbb K}(k)$.
\end{theorem}

\begin{proof}
We first identify the residue groups with large additive width.
Let $F=\mathbb F_q$, $H=(F^\times)^m$, $E=\mathbb F_p(H)$,
$u=[F:E]$, $h=|H|$, and denote by $L$ the least number of summands
from $B=H\cup\{0\}$ needed for all of $E$. We claim that for every
$\eta>0$, if $uL\ge\eta m$, then either $m=O_\eta(1)$ or
\begin{equation}\label{eq:local-frequency-rigid}
   q=p,\qquad h\in\{1,2\},\qquad
   m=j(p-1)/h,\quad j=O_\eta(1).
\end{equation}
Indeed the additive stabilizer of each iterated sumset $iB$ is stable
under $H$, hence is a vector space over $E$ inside $E$. Unless the
sumset equals $E$ the stabilizer is zero. Kneser's theorem \cite{edkneser} gives
$L\le (|E|-1)/h$, since the latter is an integer. Consequently
\[
       \frac{uL}{m}\le\frac{u(|E|-1)}{q-1},
                  \qquad m\ge(q-1)/h.
\]
If $u\ge2$ the first inequality bounds both $u$ and $|E|$ in terms
of $\eta$, and then bounds $m$. We may thus take $E=F$.
Put $b=(m,q-1)$, so that $H=(F^\times)^b$. If $b=1$ then
$L=1$. For $b\ge2$ choose an unramified local field with residue
field $F$. Every element of $E$ is a sum of residues of $b$th
powers, hence lifts to an eligible element. Its residue width is
therefore at most the primitive width in that local field.
Lemma~\ref{lem:uniform-tame-width} now gives
$L\ll\sqrt b\log(2b)\le\sqrt m\log(2m)$ unless $q=p$ and
$h\le2$. For $L\ge\eta m$ this bounds $m$ in terms of $\eta$.
In the two exceptional groups $(p-1)/h$ divides $m$, while
$L\le(p-1)/h$ gives the asserted bound for $j$.

At a place above $p\nmid k$ take $m=k$. A word for a residue target
minus one lifts after appending a unit, hence $h_k\le L+1$.
For sufficiently large $k$, a width exceeding $tk$ therefore forces
\eqref{eq:local-frequency-rigid}. At $h=1$ the unit residues
consist of $1$ and the longest shortest nonempty word has length $p$.
At $h=2$ the unit residues are $1,-1$, and this length is $(p-1)/2$
for $p\ge5$. Lemma~\ref{lem:precision-widths} gives the same
values for the ordinary width. The first case consequently contributes
$k=j(p-1)$ with $1/j\ge t$, for all sufficiently large $p$.
The second contributes $k=j(p-1)/2$ with $j$ odd and $1/j>t$.
Conversely these formulas give respectively $h=1$ and $2$ at each
degree-one place above the indicated $p$, with $p\nmid k$ for
large $p$. The field $\mathbb Q$ realizes all these tame examples.

It remains to bound the number of exponents contributed by places
above $p\mid k$, uniformly over the fields. Write $k=mp^a$, $r=a+1$.
The estimates in the proof of Theorem~\ref{thm:local-sparse} give
\[
  h_k/k\le r/k+(uL/m)(1+Cr^2/p),\qquad
                   uL\ll p\log(2m)/\log p .
\]
If $p\le r^2$ and the left side exceeds $t$, then
$m\ll_t r^2\log(2r)$. These triples contribute
$O_t((\log(2X))^5\log\log(3X))$ exponents. If $p>r^2$, the same
hypothesis implies $uL\gg_t m$ for large $k$. For bounded $m$ the
explicit estimate in the proof of Corollary~\ref{cor:local-dominant-prime}
is
\[
 \frac{h_k}{k}\le\frac r k+2m(3m/p)^a+
                         \frac{2mp}{(p-3m)(p-1)}\qquad(p>3m).
\]
This tends to zero uniformly in $a$ as $p\to\infty$, and bounded
$p>r^2$ allows only bounded $a$. In the other case
\eqref{eq:local-frequency-rigid} gives $k=jp^a(p-1)/h$ with
bounded $j,h$. Here $p\ll_t\sqrt X$, $a\ll\log(2X)$ and the count
is $O_t(\sqrt X\log(2X))$. Both errors are $o(X/\log X)$.

Finally, two distinct families $k=\alpha(p-1)$, $k=\beta(p'-1)$
with the rational multipliers just found occur together only
$o(X/\log X)$ times. To see this, write their condition as
$B p'=A p+B-A$ with coprime positive integers $A\ne B$.
For any fixed finite set of odd primes $\ell\nmid AB(B-A)$, the
primality of $p'$ excludes one nonzero class for $p$ at each
$\ell$, apart from finitely many instances. The prime number theorem
in fixed arithmetic progressions gives an upper proportion at most
$\prod_\ell(1-1/(\ell-1))$ among primes $p$. Taking successively
larger finite sets makes this product tend to zero. The multipliers
$j$ and $j/2$ for odd $j$ are distinct. Counting the individual
families by the prime number theorem or by Chebotarev now gives
the two formulas, with the same lower examples and upper error for
ordinary widths.
\end{proof}

\begin{corollary}\label{cor:comparison-fixed-field}
The field $\mathbb Q$ is the unique number field for which
$H_{\mathbb Q}(k)\le H_{\mathbb K}(k)$ for every $k$.
For every proper extension $\mathbb K/\mathbb Q$, there are at least
$(1-\delta_{\mathbb K}+o(1))X/\log X$ exponents $k\le X$ with
$H_{\mathbb Q}(k)>k\ge H_{\mathbb K}(k)$.
\end{corollary}
\begin{proof}
The transitive action of the Galois group of the normal closure on the
embeddings of $\mathbb K$ has one fixed point on average. The identity
fixes every embedding, so some element fixes none when the degree
exceeds one. Chebotarev gives $\delta_{\mathbb K}<1=\delta_{\mathbb Q}$.
Apply Theorem~\ref{thm:local-frequency} at $t=1$, where $A(1)=1$.
The difference between the counts over $\mathbb Q$ and $\mathbb K$
is a lower bound for the number of exponents in the assertion.
\end{proof}

\subsection{Equality at the exceptional thresholds}
\label{sec:local-frequency-jumps}
Above the line $s=k$ the limiting ratios of the large local widths
have one-sided, bounded rounding errors. This settles the jump values
in \eqref{eq:local-signature-count}. For $p\nmid m$ let
$\Lambda=\Lambda(p,m)$ be as in Section~\ref{sec:local-signature}.
When $\Lambda>1$, take the places $v\mid p$ attaining this maximum,
put $h_v=(q_v-1)/\gcd(m,q_v-1)$ and define
\[
 \beta(p,m)=\min_{v\text{ attaining }\Lambda}\beta_{h_v},\qquad
 (\beta_1,\beta_2,\beta_3,\beta_4,\beta_5)
                 =(0,\tfrac12,1,1,2).
\]

\begin{proposition}\label{prop:local-sharp-jumps}
If $\Lambda(p,m)>1$ then, for all sufficiently large $a$,
\[
 W_{\mathbb K}(mp^a)=H_{\mathbb K}(mp^a)
                    =\Lambda(p,m)mp^a-\beta(p,m).
\]
For every real $t>1$ the strict counting formula
\eqref{eq:local-signature-count} holds, including
$t=\Lambda(p,m)$. For the weak inequality the coefficient of
$\log X$ is instead
\[
 \sum_{\Lambda(p,m)>t}\frac1{\log p}
  +\sum_{\substack{\Lambda(p,m)=t\\\beta(p,m)=0}}\frac1{\log p}.
\]
The error in each formula is $O_{\mathbb K,t}(1)$.
\end{proposition}
\begin{proof}
The proof of Corollary~\ref{cor:local-signature-values} shows that a
local limiting ratio greater than one requires $h_v\in\{1,2,3,4,5\}$.
These directions are inert when $h_v>1$. Set $M=p^{a+b_{v,h_v}}$.
For $h_v=1$ the exact local value is $M$.
Equation~\eqref{eq:local-two-directions} gives $(M-1)/2$ for $h_v=2$.
Theorem~\ref{thm:local-inert-simplex} gives $M-1$ and $2(M-1)$
for $h_v=3,5$ respectively.

For $h_v=4$ one has $p=3$ and the quotient has two independent
coordinates modulo the odd integer $M$, with steps $\pm1,\pm i$.
Its maximum word length is $M-1$, by taking centered representatives
in the two coordinates. A shortest word using both axes lifts to the
whole local coset since $L+iL=S$. A word using only one axis has
length at most $(M-1)/2$; insert an opposite pair on the missing axis.
For the zero word insert opposite pairs on both axes. These words
have length at most $M-1$ when $M\ge5$. They can be lifted by varying
a unit in each of the two directions. This gives the same exact
ordinary and primitive width, $M-1$.

Thus every maximizing place has the value
$\Lambda mp^a-\beta_{h_v}$. A place over $p$ with a smaller limit
is eventually smaller. At the primes dividing $m$ the bounded
valuation of $mp^a$ gives bounded width by
\eqref{eq:local-signature-prime}, and the other places have widths
at most $mp^a+1$. Taking the maximum proves the first assertion.
For a given $t>1$ the argument proving Theorem~\ref{thm:local-signature},
applied at $(1+t)/2$, leaves finitely many pairs $(p,m)$.
A pair on the boundary contributes to the weak inequality exactly
when its correction is zero, and contributes to the strict inequality
in no case. The geometric sequences for distinct pairs have at most
finitely many intersections. Counting their terms gives both formulas.
\end{proof}

Over $\mathbb Q$, for example, the dyadic limits exceeding one occur
at $m=1,3$ and have values $4,4/3$, both with zero correction.
At each odd prime $p$, the two possibilities are $m=p-1$ and
$m=(p-1)/2$, both with limit $p/(p-1)$. The corrections are $0$
and $1/2$, respectively. Indeed the condition for either of the two
possible direction groups forces $m$ to be a multiple of the
corresponding minimal exponent, and a second multiple gives a
limiting ratio less than one. In particular the strict count over
$\mathbb Q$ for any $t>1$ has logarithmic coefficient
\[
 \frac{\mathbf1_{t<4}+\mathbf1_{t<4/3}}{\log2}
       + 2\sum_{\substack{p\ \mathrm{odd\ prime}\\p/(p-1)>t}}
                                                     \frac1{\log p}.
\]

\subsection{The rational envelope under bounded ramification}
\label{sec:rational-envelope}
The exceptional ratios above one are already present over the rationals
if the ramification is bounded. Residue degrees may grow freely.
\begin{theorem}\label{thm:rational-envelope}
Fix $E\ge1$ and $t>1$. For all sufficiently large $k$, uniformly over
the fields satisfying $E_{\mathbb K}(k)\le E$, one has
\[
        H_{\mathbb K}(k)\le\max\{tk,H_{\mathbb Q}(k)\}.
\]
Thus whenever the supremum of the local widths in this family exceeds
$tk$, its exact value is $H_{\mathbb Q}(k)$. The strict count of these
exponents through $X$ is the rational formula of
Section~\ref{sec:local-frequency-jumps}, with error $O_{E,t}(1)$.
The same assertions hold for ordinary widths.
\end{theorem}
\begin{proof}
A place with width exceeding $tk$ must divide $k$ for large $k$.
Write $k=mp^a$, $p\nmid m$, as before. We show that $p,m$ and the
isomorphism type of the completion lie in finite sets depending
only on $E,t$. If $p>E+1$, the logarithm on the full principal-unit
group identifies it with $\pi A$. The $k$th powers of units are the
cosets $\zeta(1+p^a\pi A)$ for $\zeta\in\mu_h$, and $k\ge ea+1$
makes all nonunit powers zero in that quotient. Write $L_0$ for the
additive width of the residue field generated by the permitted
directions in that field. Representing a residue digit costs at most
$L_0$ Teichm\"uller directions. The unramified ring generated by them
has $p^{a+1}=0$ in the quotient. Digit expansion and a reserved
unit for lifting give
\[
 h_k(A)\le 1+L_0(1+p+\cdots+p^a)
            \le 1+\frac{p}{p-1}k ,
\]
as $L_0\le m$. This excludes all sufficiently large $p$.

Fix one of the remaining primes. For bounded $a$,
\eqref{eq:nf-uniform-witt} and \eqref{eq:hcoeff} give
$h_k(A)\ll_{p,a}1+\log(2m)$, independently of the residue degree.
For $a$ sufficiently large in terms of $E$, Lemma~\ref{lem:linear-local-deviation}
bounds $h_k(A)$ by $a+1+w_k(B)$ for the unramified coefficient ring $B$.
We claim that
\[
                    w_k(B)\ll_p 1+p^a L_0 .
\]
At odd $p$ use the digit argument just given with $e=1$. At $p=2$
reduce modulo $2^{a+2}$; nonunit powers vanish for $a\ge3$. For
each permitted constant direction $\zeta$, the unit powers give
$\zeta$ plus $2^{a+1}\bar\zeta V$ in the last digit, where
$V=\{x+x^2:x\in F\}$, by squaring the binomial expansion for
$1+2x$. This is also valid after raising to the odd power $m$.
The sum of all these last-digit spaces is generated by at most
two of them, since $V$ is a hyperplane over $\mathbb F_2$.
Reserve one term in each of the chosen directions to prescribe
the deviations. The remaining constant target lies in the unramified
ring generated by the directions, and costs at most
$(2^{a+2}-1)L_0$ terms by digit expansion. An additional reserved
unit absorbs the error modulo $2^{a+2}$, after subtracting its
constant at the start. This proves the claim, including $F=\mathbb F_2$.
By \eqref{eq:hcoeff}, $L_0\ll_p\log(2m)$. On dividing the resulting
upper estimate by $mp^a$, the condition $h_k(A)>tk$ now bounds $m$
in terms of $E,t$.

For bounded $p,m$, Theorem~\ref{thm:hlargeimage} excludes sufficiently
large residue fields. The local degree $ef$ is therefore bounded. A
$p$-adic field has only finitely many extensions of bounded degree up
to isomorphism \cite{ednark}. We may consequently apply the exact
formulas of Proposition~\ref{prop:local-sharp-jumps} at the finitely
many remaining places, using also its proof at a limiting ratio
exactly $t$.

It remains to compare those formulas to the rational completion.
For odd $p$ the only directions with ratio above one are $h=1,2$
at $q=p$, or $h=4$ at $p=3$, $q=9$. At $q=p$ the rational
completion has the same directions and $b=1$, which bounds the local
value since $b_{v,h}\le1$. In the other case $\gcd(m,8)=2$ makes
$m$ even. The rational width at $3$ is then $3^{a+1}$, whereas
the four directions in the larger field give at most $3^{a+1}-1$.
At $p=2$ the possible $h$ are $1,3,5$; their exact widths
are $2^{a+b}$, $2^{a+b}-1$ and $2(2^{a+b}-1)$, respectively.
Here $b\le2$ in the first case and $b\le1$ in the others.
The rational width $4\cdot2^a$ dominates all three. This proves
the asserted upper envelope. The field $\mathbb Q$ belongs to the
family, and Lemma~\ref{lem:precision-widths} applies uniformly once
$k>E(\log_2k+1)+1$, giving the remaining assertions.
\end{proof}

One class of fields covered by this result consists of all composita
of fields of degrees at most a fixed $d$, with any number of factors.
Their ramification indices are uniformly bounded. Indeed, a finite
compositum of their Galois closures has Galois group embedded in a
product of copies of $S_d$, of exponent dividing
$\operatorname{lcm}(1,\ldots,d)$. At $p>d$ it is tame, and cyclic
inertia consequently has bounded order. At $p\le d$ there are
only finitely many local splitting fields of polynomials of degree
at most $d$, since their degrees are bounded by $d!$.
Their compositum at each of these finitely many primes is a fixed
local field. The same bound holds in all the original subfields.
In particular, if all composita of this type are allowed, the
exceptional local envelope above $tk$ is the rational one.

\subsection{Local widths at prime exponents}\label{sec:prime-degree-width}
\begin{theorem}\label{thm:prime-degree-width}
Let $p\ge5$ be prime. Uniformly for all number fields $\mathbb K$,
\[
 H_{\mathbb K}(p)\le p.
\]
Equality holds if $2p+1$ is prime and $\mathbb K$ has a place of residue degree one above $2p+1$. In the absence of such a place one has
$H_{\mathbb K}(p)\le20\sqrt p$. At every place over $p$ the primitive
local width lies between three and five for every odd prime $p$.
It equals three when $p$ divides the absolute ramification index.
\end{theorem}
\begin{proof}
First let $A$ be an integer ring at $p$, with residue field $F$.
The image modulo $pA$ of any eligible $n$ is a $p$th power by
Frobenius. Choose $w$ with $w^p\equiv n\pmod{pA}$ and write
$t=(n-w^p)/p$. The $p$ integers $j^p+1-(j+1)^p$, $0\le j<p$,
are divisible by $p$ and have sum $p-p^p$. One of them equals
$pD$ with $p\nmid D$. For any $b\in A$ the three roots
$jb,b,-(j+1)b$ sum to zero and their $p$th powers sum to $pDb^p$.
Choose $\bar b\in F$ with $\bar D\bar b^p=\bar t$. Appending $w$
and a zero gives five roots with sum $w$ and power sum equal to $n$
modulo $p\mathfrak p$. There is a unit root unless both $\bar w$
and $\bar b$ vanish. In that case use $w,1,-1,0,0$ instead.

To lift, fix the sum $w$, eliminate the fifth root, and use the
integral polynomial $(\sum_i x_i^p-w^p)/p$ in the first four
roots and $w$. Its partial derivative in $x_i$ is
$x_i^{p-1}-x_5^{p-1}$. The zero root and a unit root give a
unit partial derivative, so Hensel's lemma reaches $t$ exactly
and retains a unit root. This proves the wild assertion.

Suppose $p\mid e=v(p)$. Alter $w$ by $\pi^{e/p}b$. This changes
$(n-w^p)/p$ modulo the maximal ideal by $-(\pi^e/p)\bar b^p$,
so we may take $\bar t=0$. Begin with the three roots $w,0,0$
if $w$ is a unit and $w,1,-1$ otherwise. In each case the same
fixed-sum polynomial has a unit partial derivative and lifts the
three roots. For the lower bound, a primitive binary sum of
positive valuation has two unit roots with $s=x+y$ a nonunit.
In $(s-y)^p+y^p$ the terms of least possible valuation are
$s^p$ and $p(-y)^{p-1}s$. Comparison shows that any finite
valuation below $pe/(p-1)$ is a multiple of $p$. A common uniformizer
factor in the two roots does not change this implication. The
five-root construction with $w=0$ represents every element of $pA$.
The eligible target $p$ if $p\nmid e$, and $p\pi$ if $p\mid e$,
therefore excludes binary sums.

At a place of residue characteristic $\ell\ne p$, write $q=\ell^f$,
$H=(\mathbb F_q^\times)^p$, $C=H\cup\{0\}$, and $E=\mathbb F_\ell(H)$.
Every eligible residue belongs to $E$. The additive stabilizer of
$pC$ in $E$ is $H$-invariant and is therefore zero or all of $E$.
Kneser's bound \cite{edkneser} gives $|pC|\ge p|H|+1\ge |E|$
in the former case. For nonzero residues a representation is
primitive; for zero use the two opposite units $1,-1$. Hensel
lifting gives the first inequality. If $p\nmid q-1$, two roots
already suffice.

Suppose first $f=1$ and $p\mid\ell-1$. With $t=(\ell-1)/p=2$
the only residues are $0,1,-1$, and the residue $p$ requires $p$
terms. Otherwise $t>2$. Cochrane and Pinner
\cite[Corollary~10.3]{hprimeCP} prove that at most $20\sqrt p$
signed summands suffice over the prime field when the multiplicative
subgroup has order greater than two. Both signs are $p$th powers
for odd $p$, and zero again has a primitive two-term representation.

We bound the remaining residue degrees separately. For a nontrivial
additive character the sum over $H$ has absolute value at most
$\sqrt q$, by the Gauss sum bound. Parseval then gives for four
ordered summands from $H$ a main term $|H|^4/q$ and an error of
absolute value at most $q|H|$. In particular four suffice if
$q\ge(p+1)^3$, since then $|H|\ge(q-1)/p>q^{2/3}$.
If $p\nmid\ell-1$, the whole prime field is contained in $C$.
A basis of $E$ chosen in $H$ gives every eligible residue in at
most $f$ summands, since all of its scalar multiples belong to $C$.
For $q<(p+1)^3$ this costs at most $3\log_2(p+1)<20\sqrt p$.

If $p\mid\ell-1$ and $f\ge3$, the preceding four-sum estimate
applies. For $f=2$ the group $H$ contains the norm-one group
$T\subset\mathbb F_{\ell^2}^\times$, since $p\mid\ell-1$.
For any $z\ne0$, the conditions $x,z-x\in T$ impose a linear
trace equation on $x$ as well as $\mathrm N(x)=1$. The line
meets this nonsingular conic in at most two points. There are
$\ell+1$ ordered pairs with zero sum. Hence
\[
    |T+T|\ge 1+\frac{(\ell+1)^2-(\ell+1)}2>\ell^2/2.
\]
Four norm-one elements cover the quadratic field, completing the proof.
\end{proof}

\begin{corollary}\label{cor:prime-degree-density}
Fix an integer $j\ge2$. If $\mathbb K$ has no place of residue degree one and
residue characteristic $tp+1$ for even $t$ with $\varphi(t)<j$,
then $H_{\mathbb K}(p)\ll_j p^{1/j}$. In particular the bound
with $j=4$ excludes only $t=2,4,6$. For every $\varepsilon>0$,
all but $O_\varepsilon(X/(\log X)^2)$ primes $p\le X$ satisfy
$\sup_{\mathbb K}H_{\mathbb K}(p)\ll_\varepsilon p^\varepsilon$.
\end{corollary}
\begin{proof}
At residue degree one use the theorem of Cipra, Cochrane and
Pinner \cite[Theorem~1]{hprimeCCP} that
$\gamma(p,\ell)\ll_j p^{1/j}$ when $p\mid\ell-1$ and
$\varphi((\ell-1)/p)\ge j$. All other terms in the preceding proof
are $O(1+\log p)$. For fixed $j$ there are only finitely many $t$
with $\varphi(t)<j$. The upper-bound sieve applied to the two forms
$n,tn+1$ gives $O_t(X/(\log X)^2)$ simultaneous primes:
apart from divisors of $t$ their product excludes two residue
classes modulo each prime \cite{hprimeHR}. Take $j\ge1/\varepsilon$.
\end{proof}
For $p>3$ and $\ell=4p+1$ prime one has $3\mid2p+1$, so the
exceptional Sophie Germain case is absent. At a degree-one place
over $\ell$, the nonzero $p$th powers are $\{1,-1,i,-i\}$.
If $\ell=u^2+v^2$, the Gaussian lattice calculation in
\eqref{eq:local-five-exact} gives the exact finite-place width
$\max(|u|,|v|)-1$. In particular the square-root scale is attained
for each such pair, with the field already equal to $\mathbb Q$.

\subsection{Uniform bounds for families of number fields}
\label{sec:local-field-families}

At a finite place \(v\), let \(e_v\) be the absolute ramification
index and \(q_v\) the residue cardinality.

\begin{theorem}\label{thm:local-family-residue}
Let \(\mathcal F\) be a family of number fields. If \(e_v\) is bounded
throughout \(\mathcal F\) at each fixed value of \(q_v\), then
\(\sup_{\mathbb K\in\mathcal F}H_{\mathbb K}(k)=o(k\log k)\) as \(k\to\infty\).
If \(e_v\le Bq_v^D\) at every finite place throughout \(\mathcal F\),
with fixed \(B,D\ge1\), then for all sufficiently large \(k\) in terms
of \(B,D\) alone,
\[
             \sup_{\mathbb K\in\mathcal F}H_{\mathbb K}(k)
                      \le 4k+1+\log_2 k .
\]
\end{theorem}

\begin{proof}
First consider any sequence of places and exponents \(k\to\infty\)
with \(h_k(\mathcal O_v)\ge\eta k\log k\) for a fixed \(\eta>0\).
The bound \(k+1\) excludes places of residue characteristic \(p\nmid k\)
once \(k\) is large. Write \(k=mp^a\) with \(p\nmid m\). The local
estimate \(h_k(\mathcal O_v)\ll k(1+a)\) in the proof of
Theorem~\ref{thm:hhybrid} shows that \(a\gg_\eta\log k\).
Since \(a\log p\le\log k\), we may pass to a subsequence with \(p\)
fixed and \(a\to\infty\). For this \(p\) the same proof gives
\[
 \frac{h_k(\mathcal O_v)}{k}\ll_p
           1+(a+1)\frac{\log^2(2m)}{m}\qquad(m\ge64).
\]
Indeed, take \(c=1+\lfloor\log_p(m/64)\rfloor\). For \(c\ge a+1\)
the bound is already \(O(k)\). Otherwise use \(p^{1-c}<64p/m\)
and \(c\ll_p\log(2m)\) in the estimate for the two sets of Witt
coordinates. Dividing by \(\log k\) rules out \(m\to\infty\).
We may therefore take \(m\) fixed. Corollary~\ref{cor:inert-residue-test}
bounds the width by \(O_{p,m}(k)\) outside its finite list of residue
cardinalities. Pass to a subsequence on which \(q_v\) is fixed on
that list.

The proof of Theorem~\ref{thm:local-fixed-ramification} applies at
each individual place with \(e_v\le\sqrt{k}/(\log(2k))^4\) and gives
\(h_k(\mathcal O_v)\le4k+1+\log_2 k\) there for large \(k\).
Our sequence therefore has \(e_v\to\infty\). This contradicts the
first hypothesis at its fixed residue cardinality. The rate in the
little-\(o\) estimate may depend on the family.

For the second assertion it suffices to consider a place with
\(h_k(\mathcal O_v)>2k\). For large \(k\) its residue characteristic
divides \(k\). The proof of Theorem~\ref{thm:local-sparse} at this
place gives \(m\ll(a+1)^2\log(2a+2)\). Theorem~\ref{thm:hlargeimage}
with \(\theta=\varepsilon=1\) gives \(q_v\ll(1+m)^3\), since
\(\gcd(m,q_v-1)\le m\). Hence \(e_v\le Bq_v^D\) is bounded by a
fixed power of \(\log(2k)\) times a constant depending only on
\(B,D\). For large \(k\) this is at most \(\sqrt{k}/(\log(2k))^4\),
and the same placewise estimate proves the assertion.
\end{proof}

By Theorem~\ref{thm:intro-main}, the first hypothesis also gives
\(G_{\mathbb K}(k)=o(k\log k)\) uniformly on \(\mathcal F\).
Under the polynomial hypothesis the analytic entry
\(\lfloor C_1k-\sqrt{k}\rfloor\) bounds \(G_{\mathbb K}(k)\) throughout
the family for sufficiently large \(k\).

The proof also locates a possible obstruction over unrestricted fields.
A sequence with \(H_{\mathbb K}(k)\gg k\log k\) has a subsequence
realized at places with \(k=mp^a\), \(p,m,q_v\) fixed, \(a\to\infty\),
and \(e_v>\sqrt{k}/(\log(2k))^4\). Here \(q_v\) belongs to the finite
list of Corollary~\ref{cor:inert-residue-test}. At \(k=2^a\) that
list consists of \(2\) and \(4\).

\begin{example}
Fix \(D\ge1\). The second assertion applies to all abelian number
fields whose conductors \(N\) satisfy \(v_p(N)\le D\) for every
prime \(p\), including all abelian fields of cube-free conductor when
\(D=2\). Indeed, such a field lies in \(\Q(\zeta_N)\).
For \(p^d\Vert N\), adjoining the prime-to-\(p\) roots of unity is
unramified at \(p\), and the \(p^d\)-cyclotomic factor has ramification
index \(\varphi(p^d)\). Thus \(e_v\le p^D\le q_v^D\).
\end{example}

\begin{example}
The fields \(\Q(\zeta_{6^n})\) also satisfy the second assertion,
although their ramification indices at both \(2\) and \(3\) tend to
infinity. For \(n\ge3\) the pairs \((e_v,f_v)\), with \(q_v=p^{f_v}\),
are
\[
 (2^{n-1},\,2\cdot3^{n-1})\quad(p=2),\qquad
 (2\cdot3^{n-1},\,2^{n-2})\quad(p=3).
\]
The residue degrees are the orders of \(2\) modulo \(3^n\) and of
\(3\) modulo \(2^n\), respectively. These follow from
\(v_3(2^{2\cdot3^j}-1)=j+1\) for \(j\ge0\) and
\(v_2(3^{2^j}-1)=j+2\) for \(j\ge1\).
It follows that \(e_v\le2q_v\) throughout the tower, including
\(n=1,2\), and all other primes are unramified.
\end{example}

\subsection{Powers of two in dyadic radical extensions}
\label{sec:dyadic-radical-width}

Even when the residue field is fixed and the ramification index grows,
the local width can be smaller than the rational width. A \(k\)th
root of \(2\) permits each binary digit of a coefficient to be
represented by a single \(k\)th power.

\begin{proposition}\label{prop:dyadic-radical-width}
Let \(A\) be the integer ring of a finite extension of \(\Q_2\) with
residue field \(\mathbb F_2\). If \(k=2^a\), \(a\ge1\), and \(b^k=2\)
for some \(b\in A\), then \(h_k(A)\le3k+a\). In particular,
\[
 H_{\Q(2^{1/2^n})}(2^a)\le3\cdot2^a+a\qquad(n\ge a\ge1).
\]
For \(a\ge2\) this is strictly below \(H_{\Q}(2^a)=4\cdot2^a\).
\end{proposition}

\begin{proof}
Use the last construction in the proof of Theorem~\ref{thm:hhybrid}
with \(p=2\), \(m=1\), \(r=a+1\). Over \(R=\mathbb F_2[[T]]\) its
ideal \(I_N\subset W_r(R)\) consists of the vectors whose \(j\)th
coordinate is divisible by \(T^{N_j}\), where \(N_j=(a+1-j)2^j\).
Write \(G\) for the additive group generated by \([y]\), \(y\in R\).
If \(B_0=\sum_{i=0}^a[T^i]\), every element of \(I_N\) is a sum of
\(r\) variable Teichm\"uller vectors minus \(B_0\), with the root at
the centre \(1\) remaining a unit. Indeed, lift the coordinates with
their indicated divisibilities to \((\Z/4\Z)[[T]]\), use the centred
sum constructed there, and reduce modulo two. In particular \(I_N\subset G\).

Put \(M=W_r(R)/I_N\) and \(S=G/I_N\). Filter \(M\) by the images
of vectors whose first \(j\) coordinates vanish. The \(j\)th
successive quotient is the additive group \(R/T^{N_j}R\).
Consequently \(\log_2|M|=\sum_{j=0}^a N_j=4k-a-3\).
A class of \(S\) in the \(j\)th filtration subgroup has a
representative with its first \(j\) coordinates zero, and every such
representative belongs to \(G\) since \(I_N\subset G\). For \(j\ge1\)
truncate after coordinate \(j\) and apply
Proposition~\ref{prop:witt-top-functions}. The coefficients at odd
degrees below \(2^j-1\) in coordinate \(j\) vanish. Since \(N_j\ge2^j\),
the successive quotient for \(S\) loses at least \(2^{j-1}-1\)
binary digits. Thus
\[
       \log_2|S|\le 4k-a-3-\sum_{j=1}^a(2^{j-1}-1)=3k-2 .
\]

Evaluation at a uniformizer followed by the last ghost component maps
\(W_r(R)\) to \(A/2^{a+1}A\). This is well defined because
\((u+2z)^{2^j}\equiv u^{2^j}\pmod {2^{j+1}A}\).
Evaluation \(R\to A/2A\) is onto, so the image of \(G\) is precisely
the reduction of \(J_k(A)\). Lift an eligible target to \(G\) and
subtract the baseline \(B_0\). Choose Teichm\"uller generators that
successively enlarge the subgroup they generate in \(S\), with
successive indices \(2^{d_i}\). Reducing in these cyclic quotients
in reverse order represents the class of the difference using
coefficients \(0\le n_i<2^{d_i}\). After evaluation each multiple
\(n_ix_i^k\) costs at most \(d_i\) powers, since
\(2^j x_i^k=(b^j x_i)^k\). The remaining difference belongs to
\(I_N\), where the \(r\) reserved roots restore the baseline.
This represents the target modulo \(2^{a+1}A\) using at most
\(r+\sum_i d_i=r+\log_2|S|\le3k+a-1\) powers, including a unit root.

Write the error as \(2^{a+1}d\) and choose \(\tau\in\{0,1\}\) with
\(d-\tau\) in the maximal ideal. If \(\tau=1\), append
\((b^{a+1})^k=2^{a+1}\). For a retained unit root \(x\), put
\(g(T)=((1+2T)^k-1)/2^{a+1}\). Iterated squaring gives
\(g\in\Z[T]\), and \(g'(0)=1\). Hensel's lemma at zero solves
\(g(t)=(d-\tau)/x^k\) in \(A\). Replacing \(x\) by \(x(1+2t)\)
gives an exact primitive representation.

Finally, \(X^{2^n}-2\) is Eisenstein at two. The field
\(\Q(2^{1/2^n})\) therefore has a unique dyadic place with residue
field \(\mathbb F_2\), and it contains \(b=(2^{1/2^n})^{2^{n-a}}\)
when \(n\ge a\). At every other place \(p\nmid k\), so \(k+1\) terms
suffice.
\end{proof}

\section{Balancing the target and Birch's local-to-global theorem}\label{sec:hbirch}
The passage between a norm cutoff and a fixed region in Minkowski
space is effected by the units. For an infinite place $v$ put $d_v=1$
at a real place and $d_v=2$ at a complex place, using the usual modulus
$|\cdot|_v$ in both cases.

\begin{lemma}\label{lem:hbal}
There is a constant $A=A(\mathbb K,k)\ge1$ such that every nonzero
$\nu\in\mathcal O$ admits a totally positive unit $u$ with
\[
 A^{-1}N(\nu)^{1/d}\le |u^k\nu|_v\le A N(\nu)^{1/d}
                  \qquad(v\mid\infty).
\]
The map $x_i\mapsto ux_i$ carries integral representations of $\nu$
onto integral representations of $u^k\nu$, preserving the norms and
real signs of the summands.
\end{lemma}
\begin{proof}
Write $t_v=\log|\nu|_v-d^{-1}\log N(\nu)$. The vector $t$ belongs to
the hyperplane $\sum_v d_v t_v=0$. By Dirichlet's unit theorem, the
vectors $(k\log|u|_v)_v$, for totally positive units $u$, form a full
lattice in this hyperplane, since the sign map has finite image and
multiplication by $k$ preserves rank. Translate $t$ into a fixed
bounded fundamental parallelepiped and exponentiate. The last claim
uses $N(u)=1$ and the identity $(ux_i)^k=u^kx_i^k$.
\end{proof}

For example, take $T^k=2A N(\nu)^{1/d}$. The normalized target
$u^k\nu/T^k$ then ranges over a compact set where every real coordinate
lies between $1/(2A^2)$ and $1/2$, and every complex modulus is at most
$1/2$. Proposition~\ref{prop:harch} gives a uniform positive lower
bound there. A representation in the dilate $TB$ has
$N(x_i)\le T^d\ll_{\mathbb K,k}N(\nu)^{1/k}$. In integral-basis
coordinates multiplication by $u$ is a lattice automorphism and sends
a box to a parallelotope. The target can be balanced first and then
represented in a fixed collection of boxes. At each finite place the
same change of variables is a bijection on residues modulo
$\mathfrak p^j$ and carries the local count for $\nu$ to that for
$u^k\nu$. A finite covering of the compact Archimedean target set
suffices when smaller boxes about nonsingular solutions are used.

Birch \cite[Theorem~2]{hbirch61} separated the analytic threshold
$s\ge2^k+1$ from the finite congruence conditions. In discussing
positivity of the singular series, Birch \cite[p.~449]{hbirch61}
credits Stemmler's work on signed sums of powers \cite{hstemmler61}
and Gray's work on diagonal forms of prime degree. For prime exponent
$p$ he deduced $G_{\mathbb K}(p)\le2^p+1$. For arbitrary $k>1$,
Ramanujam \cite{hram63} showed that each element of the order
$J_k(\mathcal O_v)$ generated by $k$th powers is a sum of at most
$8k^5$ powers in $\mathcal O_v$ at every finite place $v$. He deduced
$G_{\mathbb K}(k)\le\max\{2^k+1,8k^5\}$ by Birch's theorem.
The inequality is uniform in $d$. One must still distinguish a
local bound, the scale at which a lift becomes nonsingular, and
the error term in the analytic theorem. For example, at a residue
characteristic dividing $k$, the derivative $kx^{k-1}$ can have
positive valuation even for $x$ a unit.

\section{Transfer by a basis of powers}\label{sec:powerbasis}
The following bound uses the classical all-integer threshold and requires
no estimate for minor arcs over $\mathbb K$.

\begin{theorem}\label{thm:hpowers}
For every number field $\mathbb K$ of degree $d$ and every $k\ge2$,
\[
                       G_{\mathbb K}(k)\le dG(k)+1 .
\]
\end{theorem}

We shall choose $d$ directions $y_i^k$ so that the target lies in their
positive real span. The lattice spanned by these directions can be a
proper sublattice of $J_k$. The extra summand in the theorem accounts
for every coset of this lattice.

\begin{lemma}\label{lem:htranslate}
Let $A$ be a finite \emph{\'etale} algebra of dimension at most $d$ over
$\mathbb F_q$, and let $c\in A$. If $q>(d k^d+d+1)^2$, there are
$t\in\mathbb F_q$ and $z\in A^\times$ such that
$z^k=c+t$.
\end{lemma}
\begin{proof}
The subalgebra $B=\mathbb F_q[c]$ is finite \'etale, and $c$ is regular
in $B$ in the sense of \cite[Theorem~2]{hkatz}. Indeed its
$r=\dim B$ coordinates are the distinct roots of its squarefree
minimal polynomial, even if some coordinates coincide in $A$. Put
$m=[B^\times:(B^\times)^k]$. The product decomposition of $B$ into at
most $r$ finite fields gives $m\le k^r$. Extend the characters of
$B^\times$ by zero to the nonunits, including the principal character.
By orthogonality, the number of good $t$ is
\[
 \frac{1}{m}\sum_{\chi\in(B^\times/(B^\times)^k)^\wedge}
                   \sum_{t\in\mathbb F_q}\chi(c+t).
\]
At most $r$ translates are nonunits. Katz's bound for each nonprincipal
inner sum is $(r-1)\sqrt q$ in absolute value. Hence the number of good
$t$ is at least $(q-r-(m-1)(r-1)\sqrt q)/m>0$. The root in $B^\times$
also lies in $A^\times$.
\end{proof}

Let $V_+$ denote the open subset of $V$ on which the real coordinates
are positive. For each $u\in V_+\setminus\{0\}$ one can choose a real
basis $b_1,\dots,b_d$ in $V_+$ with $u$ in the interior of its positive
span. Indeed, choose a linear functional with value~$1$ at $u$ and
use the vertices of a small $(d-1)$-simplex about $u$ in that
affine hyperplane. Its vertices are linearly independent, and each
lies in $V_+$ if the simplex is sufficiently small. We may take
every complex coordinate of every vertex to be nonzero, by a small
perturbation. Each $b_i$ has a coordinatewise $k$th root $w_i\in V_+$
with positive real coordinates, for odd and even $k$ alike.
Real linear independence and the condition that $u$ have strictly
positive coordinates in the basis persist under sufficiently small
perturbations of the $b_i$.

\begin{lemma}\label{lem:hpowlattice}
Fix $u\in V_+\setminus\{0\}$ and a finite set $S$ of rational primes
containing the divisors of $k!\,\mathrm{disc}(\mathbb K)$. There are
totally positive $y_1,\dots,y_d\in\mathcal O$ such that $u$ lies in
the interior of the positive span of the real basis
$y_1^k,\dots,y_d^k$ and such that
\[
 M=[J_k:\mathbb Z y_1^k+\cdots+\mathbb Z y_d^k]
       \quad\hbox{is prime to every $p\in S$ and to $N(y_1)$.}
\]
\end{lemma}
\begin{proof}
For each $p\in S$, the $k$th powers span $J_k/pJ_k$ additively,
so choose $a_{p1},\dots,a_{pd}\in\mathcal O$ whose powers form a basis
over $\mathbb F_p$. To reproduce such a power modulo $pJ_k$ it
suffices to reproduce its root modulo $p^{1+v_p(k!)}\mathcal O$.
Indeed the difference of the powers belongs to $J_k$ as well as to
$p^{1+v_p(k!)}\mathcal O$. Locally at $p$ the latter is contained in
$pJ_k\otimes\mathbb Z_p$, since $k!\mathcal O\subseteq J_k$.
The reduction $J_k/pJ_k$ can be tested after tensoring with $\mathbb Z_p$.

By the Chinese remainder theorem, a root $y'_1$ with all first-column
residues can be chosen in a fixed coset of a full lattice of $V$.
For arbitrarily large $H$, the coset contains a point within bounded
distance of $H w_1$, with $w_1$ as above. Hence $y'_1$ may be taken
totally positive, with $(y'_1)^k$ arbitrarily close in direction
to $b_1$. Divide $y'_1$ by its positive rational content, the gcd of
its coefficients in an integral basis, to obtain a primitive element
$y_1$. The content is prime to every $p\in S$, since divisibility of
$y'_1$ by $p$ would give $(y'_1)^k\in pJ_k$. Thus at each of these
primes this division multiplies the first power vector by a nonzero
scalar modulo $pJ_k$. Its real direction is preserved.

For any $p\mid N(y_1)$ outside $S$, the algebra
$\mathcal O/p\mathcal O$ is finite \'etale, and
$J_k/pJ_k\simeq\mathcal O/p\mathcal O$ since $p\nmid k!$.
The reduction of $y_1$ is nonzero by primitivity. The algebra is
reduced, so the reduction of $y_1^k$ is also nonzero. Complete this
vector to a basis of $J_k/pJ_k$ by choosing further vectors from the
spanning set of powers. Enlarge the finitely many congruence conditions
on the remaining columns accordingly. Approximation to $Hw_i$ in
the resulting lattice cosets now gives the desired $y_2,\dots,y_d$.
Each column can be rescaled by its own positive real factor when
testing the cone condition. The determinant in a basis of $J_k$
is nonzero and is a unit at every prescribed prime, as required.
\end{proof}

\begin{proof}[Proof of Theorem~\ref{thm:hpowers}]
Use Lemma~\ref{lem:hbal} to balance $\nu$. Write $R=N(\nu)^{1/d}$. The normalized balanced targets $u=\nu/R$ range over a compact subset
of $V_+\setminus\{0\}$. Take $S$ to contain the ramified primes, all
divisors of $k!$, and all rational primes up to $(dk^d+d+1)^2$.
For each point of this compact set choose a basis as in
Lemma~\ref{lem:hpowlattice}. There is a neighborhood of the point
on which all coordinates in that basis are bounded above and
bounded away from zero. Choose a finite subcover.

Fix one basis and put $\Lambda=\sum_i\mathbb Z y_i^k$, $M=[J_k:\Lambda]$.
We claim that $k$th powers of totally positive elements meet every
coset of $\Lambda$ in $J_k$. Let $a\in J_k$ represent a coset.
For each $p^e\parallel M$ the algebra $A=\mathcal O/p\mathcal O$ is
\'etale and $y_1$ is a unit in it. Apply Lemma~\ref{lem:htranslate}
to $c=a/y_1^k$ in $A$. It gives $t_p\in\mathbb F_p$ for which
$a+t_p y_1^k$ has a unit $k$th root modulo $p\mathcal O$.
Lift $t_p$ to a rational integer. Since $p\nmid k$, the derivative
$kx^{k-1}$ is a unit at the chosen root in every factor of $A$.
Hensel lifting gives a root modulo $p^e\mathcal O$.
Combine the $t_p$ and the roots by the Chinese remainder theorem to get
$t\in\mathbb Z$ and $x\in\mathcal O$ with
\[
                      a+t y_1^k\equiv x^k\pmod {M\mathcal O}.
\]
Adding a sufficiently large positive multiple of $M$ to $x$ makes
it totally positive. Since $M$ is prime to $[\mathcal O:J_k]$,
one has $M\mathcal O\cap J_k=MJ_k\subseteq\Lambda$.
Thus $x^k\equiv a\pmod\Lambda$. When $M=1$ the assertion follows
by taking $x=1$. Fix such an $x$ once for each coset and each of
the finitely many bases.

For every sufficiently large balanced $\nu$ in the selected
neighborhood, subtract the power $x^k$ for its coset. There are
integers $n_i$ such that
\[
                        \nu-x^k=\sum_{i=1}^d n_i y_i^k .
\]
The coordinates of $\nu$ in the chosen basis are bounded above and
below by positive multiples of $R$, and $x$ ranges over a finite set.
Hence $n_i\asymp_{\mathbb K,k}R$. Each $n_i$ is a sum of at most $G(k)$
positive rational $k$th powers. After multiplication by $y_i^k$ these
give totally positive summands, with each root $b y_i$ satisfying
$N(b y_i)=b^dN(y_i)\ll R^{d/k}$. The extra root $x$ also satisfies
\eqref{eq:hnorm}. Undo the balancing unit to conclude.
\end{proof}

The congruences used to select the basis can make the cutoff in
Theorem~\ref{thm:hpowers} large. Its assertion concerns only the
number of summands. In particular it may be combined with
\eqref{eq:hdav98} and the positive-density transfer of \cite{hdav99}
by taking the smallest of the applicable bounds.

\paragraph{An example with ramification.}
For $\mathbb K=\mathbb Q(\sqrt6)$ and $k=4$ the order is
$J_4=\mathbb Z+4\mathbb Z\sqrt6$. Indeed the coefficient of $\sqrt6$
in $(a+b\sqrt6)^4$ is $4ab(a^2+6b^2)$. The choices
$(a,b)=(1,1),(2,1)$ give $28$ and $80$, whose gcd is $4$.
The positive roots $y_1=5+\sqrt6$ and $y_2=1657-2\sqrt6$
give columns $(1561,155)$ and $(7538983805233,-9099160322)$ in
the displayed basis of $J_4$. Their cone contains the positive
rational ray in its interior, since the differences of the two
real conjugates of these powers have opposite signs.
The index is
\[
 M=1667\cdot7297\cdot8297\cdot11719 .
\]
All four factors are prime and exceed $(2\cdot4^2+3)^2=1225$, and
$N(y_1)=19$. Thus Lemma~\ref{lem:htranslate} applies at every divisor
of $M$, including the split prime $7297$. The dyadic congruences are
accounted for by choosing the lattice in $J_4$.

\subsection{Transfer from almost all rational integers}
Write $G_{\mathrm{ae}}(k)$ for the least number of positive $k$th powers
representing all but $o(X)$ of the positive integers up to $X$, and
$\widehat G(k)$ for the least number representing a set of positive
lower density. The numbers in both definitions count an upper bound on
the number of summands.

\begin{theorem}\label{thm:haetransfer}
For every number field $\mathbb K$ of degree $d$,
\[
               G_{\mathbb K}(k)
                   \le dG_{\mathrm{ae}}(k)+\widehat G(k)+1 .
\]
\end{theorem}
\begin{proof}
Keep the finitely many bases, coset representatives and notation of
Theorem~\ref{thm:hpowers}. For each basis fix a totally positive $y_0$
with $y_0^k\in\Lambda$ such that all coordinates $c_1,\ldots,c_d$ of
$y_0^k$ in this basis are nonzero. Indeed, the union of the coordinate
hyperplanes has empty interior. Pick a point of $V_+$ outside them
and with nonzero complex coordinates. Approximate a large positive
dilate of its coordinatewise $k$th root by an integral element and
multiply the element by $M$. Its power then lies in $MJ_k\subseteq
\Lambda$ and remains outside the hyperplanes. In particular
$c_i\in\mathbb Z\setminus\{0\}$.

Write the coefficients of $\nu-x^k$ in this basis as $n_i$.
There is a fixed $\eta>0$ so small that
$n_i-c_it\asymp R$ and $n_i-c_it>0$ for all $1\le t\le\eta R$
and all sufficiently large balanced targets in this neighborhood.
Positive lower density supplies $\gg R$ choices of $t$ in that
interval representable by at most $\widehat G(k)$ powers. For a fixed
$i$, at most $o(R)$ integers $t$ in the interval have $n_i-c_it$
exceptional for $G_{\mathrm{ae}}(k)$, since $t\mapsto n_i-c_it$ is
injective and all its values lie in $[1,CR]$ for a fixed $C$.
The union for $1\le i\le d$ still has size $o(R)$ for the fixed field.
Choose a remaining $t$ and write
\[
            \nu=x^k+t y_0^k+\sum_{i=1}^d (n_i-c_it)y_i^k .
\]
All scalar coefficients admit the required representations, and all
roots satisfy \eqref{eq:hnorm} as in Theorem~\ref{thm:hpowers}.
\end{proof}

There is also a version for an arithmetic progression. Suppose all
but $o(X)$ positive $n\le X$ in one class $a\pmod q$ are sums of at
most $s$ positive rational $k$th powers. Then
\begin{equation}\label{eq:hprogression}
              G_{\mathbb K}(k)\le ds+\widehat G(k)+h_q,
 \qquad h_q=\max\bigl(\{1\}\cup
       \{w_k(\mathcal O_{\mathbb K_v}):v\mid q\}\bigr).
\end{equation}
To prove this assertion, include all primes dividing $q$ in~$S$
when choosing a basis, so $(M,q)=1$. For each of the finitely many
cosets of $q\Lambda$ in $J_k$, subtract $a\sum_i y_i^k$ from a
representative. At every place over $p\mid q$, this difference belongs
to the local power order and so is a sum of $h_q$ powers with zeros
allowed. If $p^b\parallel q$, retain the roots to precision
$p^{b+v_p(k!)}\mathcal O\otimes\mathbb Z_p$. An error in the sum
at this precision belongs to $p^b J_k\otimes\mathbb Z_p$, and
$J_k\otimes\mathbb Z_p=\Lambda\otimes\mathbb Z_p$ since $p\nmid M$.
At $p\mid M$ use the first root to match the coset modulo
$\Lambda\otimes\mathbb Z_p$, as in Theorem~\ref{thm:hpowers}, and
take the others to be zero at that precision. The Chinese remainder
theorem now supplies $h_q$ totally positive integral lifts, fixed
for each coset. After subtraction of their powers the coordinates
$n_i$ of the original target in the chosen basis are all congruent
to $a$ modulo $q$. Choose $y_0$ with $y_0^k\in q\Lambda$ and all
coordinates nonzero, by taking a further multiple of $q$ if necessary.
Then $n_i-c_i t$ stays in the prescribed progression, and the same
count of exceptional $t$ proves \eqref{eq:hprogression}.

\subsection{Transfer from an intermediate field}
Let $F\subset E$ be number fields, put $e=[E:F]$ and write
$J_L$ for the order generated by the integral $k$th powers of
a field $L$. Relative transfer is particularly simple if $J_F$
is the maximal order.

\begin{theorem}\label{thm:hrelative}
If $J_F=\mathcal O_F$, then
\[
                      G_E(k)\le eG_F(k)+1 .
\]
\end{theorem}
\begin{proof}
We may suppose $e\ge2$. The order $J_E$ is a torsion-free
$\mathcal O_F$-module of rank $e$, since it contains $J_F$ and
$k!\mathcal O_E$. Thus it is free of rank $e$ locally at each
prime of $\mathcal O_F$. Let $S$ contain the finite primes of $F$
that divide $k!$, ramify in $E$, or have residue cardinality at most
$(ek^e+e+1)^2$. At each $\mathfrak p\in S$, choose $e$ powers forming
a basis of $J_E/\mathfrak p J_E$ over the residue field. Fix the
corresponding roots modulo $\mathfrak p^{1+v_{\mathfrak p}(k!)}\mathcal O_E$.
The verification preceding Lemma~\ref{lem:hpowlattice} applies locally.

Balance the target over $E$ as in Lemma~\ref{lem:hbal}, and fix its
normalized direction $u$. Write $E_v=E\otimes_F F_v$ for $v\mid\infty$
of $F$. At a real $v$, the simplex construction above supplies a real
basis of $E_v$ consisting of powers of vectors positive at its real
factors, with all coordinates of $u_v$ positive in that basis.
At a complex $v$, take a complex linear functional $L$ with
$L(u_v)=1$. Perturb $u_v$ inside $L^{-1}(1)$ by $e$ small vectors
with sum zero and affine span of complex dimension $e-1$, obtaining
a complex basis with $u_v$ its barycentre. All coordinates at complex
places of $E$ can be kept nonzero, since those of $u$ are nonzero.
Taking roots at all infinite places gives open sets of admissible
columns. In a neighborhood of $u$, the coefficients in each local
basis are bounded and bounded away from zero, and are positive at
the real places of $F$.

The first root can be chosen in its prescribed congruence classes
and outside $\mathfrak p\mathcal O_E$ for every finite prime
$\mathfrak p$ of $F$. To see this, take dilates by $H$ of a small
bounded open box about the prescribed Archimedean column.
At primes in $S$ the root is already nonzero modulo
$\mathfrak p\mathcal O_E$, since a multiple of a local uniformizer
would have power in $\mathfrak p J_{E,\mathfrak p}$.
At any further finite set of primes, the Chinese remainder theorem
and lattice point counting give the density in the original coset
as the product of $1-(N\mathfrak p)^{-e}$. These products have a
positive limit. For completeness, put $D=[E:\mathbb Q]$ and
$f=[F:\mathbb Q]$. The least sup norm of a nonzero element of
$\mathfrak p\mathcal O_E$ is at least $(N\mathfrak p)^{e/D}$,
by taking its absolute norm. Packing gives
$O(1+H^D/(N\mathfrak p)^e)$ multiples in the box, and such a
nonzero multiple requires $N\mathfrak p\ll H^f$. Hence the
upper density of the union over $N\mathfrak p>Q$ is at most
a fixed multiple of $\sum_{N\mathfrak p>Q}(N\mathfrak p)^{-e}$.
Indeed, the terms $1$ contribute $O(H^{f-D})$ after dividing
by $H^D$, since there are $O_F(X)$ primes of norm at most $X$.
Letting $Q$ tend to infinity proves the assertion.

Choose this first root $y_1$. At all additional primes dividing
$N_{E/F}(y_1)$, the algebra $\mathcal O_E/\mathfrak p\mathcal O_E$
is \'etale and identifies with $J_E/\mathfrak p J_E$. Thus the
nonzero vector $y_1^k$ extends to a basis of powers there too.
Approximate the other roots in the resulting finitely prescribed
lattice cosets. Write
\[
 \Lambda=\textstyle\sum_{i=1}^{e}\mathcal O_Fy_i^k,\qquad
 \mathfrak m=\operatorname{Ann}_{\mathcal O_F}(J_E/\Lambda).
\]
The independent columns give a finite quotient with annihilator
coprime to $S$ and $N_{E/F}(y_1)$. For each $a\in J_E$ and
$\mathfrak p\mid\mathfrak m$, apply Lemma~\ref{lem:htranslate}
over $\mathcal O_F/\mathfrak p$ with $c=a/y_1^k$. Simple-root
lifting and the Chinese remainder theorem give
$x\in\mathcal O_E$, $t\in\mathcal O_F$ with
$x^k\equiv a+ty_1^k\pmod{\mathfrak m\mathcal O_E}$.
We can make $x$ totally positive by adding a positive rational integer
in $\mathfrak m$ sufficiently many times. Locally at primes dividing
$\mathfrak m$, the orders $J_E$ and $\mathcal O_E$ agree, so
$\mathfrak m\mathcal O_E\cap J_E=\mathfrak m J_E\subseteq\Lambda$.
Thus one power supplies every coset.

Take finitely many choices of columns for the compact balanced target
set, and fix a root $x$ once per coset for each choice. For large
$R=N_E(\nu)^{1/D}$, the coefficients $n_i$ of $\nu-x^k$ are now in
$\mathcal O_F=J_F$, with each $|n_i|_v\asymp R$ and with $n_i$ totally
positive. Hence $N_F(n_i)\asymp R^f$, and each coefficient is a sum
of at most $G_F(k)$ powers. Each resulting root $by_i$ is totally
positive in $E$ and satisfies
$N_E(by_i)=N_F(b)^eN_E(y_i)\ll R^{D/k}$. The same bound holds for
the extra root, and the balancing unit can be undone.
\end{proof}

For a general intermediate field one can retain finite local solubility
as a separate term. Suppose that, at every finite place $w$ of $E$,
each $a\in J_E$ is a sum of at most $h$ powers in
$\mathcal O_{E,w}$. Then
\[
                         G_E(k)\le eG_F(k)+h .
\]
Indeed choose independent columns with the same Archimedean properties
as above, using ordinary lattice approximation. The additive lattice
$\sum_i J_F y_i^k$ has full rank and hence contains $m\mathcal O_E$
for some positive rational integer $m$. For each of its finitely many
cosets, reduce a representative modulo $m\mathcal O_E$. Pad the local
representations at the primes dividing $m$ by zeros to $h$ terms.
The Chinese remainder theorem gives $h$ totally positive integral
lifts, chosen once per coset. After their powers are subtracted,
the coefficients of the remaining target belong to $J_F$ and have
the same Archimedean bounds as before.

The next estimates use critical mean values for Vinogradov systems.
Bourgain, Demeter and Guth proved the rational main conjecture in
degrees at least four by decoupling \cite{NFBDG}. Wooley's
efficient congruencing handles the cubic case, and his nested
iteration gives a weighted theorem on the full coefficient torus
over a number field \cite{NFWooley}. Restricting to the leading
coefficient in a Waring sum requires a separate sampling argument.
Sections~\ref{sec:nf-global} and~\ref{sec:nf-fibres} develop this
argument in a form that applies to arbitrary degree. For small powers
the sharper smooth mean values are established in
Section~\ref{sec:nf-small-powers}.
\section{Small exponents and the quartic obstruction}
\label{sec:nf-small-powers}

Write $J_k=\sum_{x\in\mathcal O}\mathbb Zx^k$. The representation
bounds below hold for every totally positive $\nu\in J_k$ of
sufficiently large absolute norm, with
$|N x_i|\ll_{\mathbb K,k}|N\nu|^{1/k}$ for each root. The roots are
totally positive, and zero summands may be removed. The results allow
every class group and signature. The cubic bound
$G_{\mathbb K}(3)\le7$ is proved in Theorem~\ref{nf3:realG7}.
For totally imaginary fields containing $\sqrt[3]{2}$, or admitting
units $a,b$ with $a^3+b^3=1$, Corollary~\ref{nf3:imag6} gives six.

\begin{theorem}\label{thm:nf-small}
For every number field $\mathbb K$ one has
\[
\begin{array}{c|rrrrrrr}
k&4&5&6&7&8&9&10\\\hline
H_{\mathbb K}(k)\le&16&5&12&7&35&25&17\\
G_{\mathbb K}(k)\le&16&17&24&31&38&44&50 .
\end{array}
\]
Moreover $G_{\mathbb K}(4)\le\max\{12,H_{\mathbb K}(4)\}$.
If no dyadic completion is $\mathbb Q_2$, then $G_{\mathbb K}(4)\le15$.
If every dyadic residue degree is at least two, then $G_{\mathbb K}(4)\le12$.
\end{theorem}

\begin{proof}
Section~\ref{subsec:aud-small-local} and Theorem~\ref{thm:prime-degree-width} give the primitive local bounds, including
the three dyadic alternatives at degree four. We give the analytic reduction
here and prove its moment estimates below. Write $X$ for the root volume,
$w=k(k+1)$, and $f$ for the ordinary sum on a smooth nonnegative annular
envelope. A restricted sum $G$ has defect $D_r$ at order $r$ when
\[
        \int_{\mathbb T}|G|^r\ll_\epsilon X^{r-k+D_r+\epsilon}.
\]
The estimates below are uniform for bounded coefficients and fixed ideal
lattices. Absolute upper-smooth estimates permit all roots with prime ideal
factors of norm at most $X^\eta$. For a lower bound on representations we
use the positive marked weight of Sections~\ref{sec:nf-global}
and~\ref{sec:nf-fibres}, with its terminal free cofactor. Its coefficients
are bounded and its support lies in the respective upper-bound families;
Section~\ref{subsec:aud-full-moments} specifies the parameters.

Let $u<S\le w$ be even, $S>2k$, $a=S-u$ and $A=\int|G|^S$. Take
absolute values of coefficients in the even orthogonality count first.
Enlarging $a/2$ root slots on each side to the ordinary envelope gives
$A\le\int |f|^a|G|^u$ for the resulting positive weight, which we
again denote by $G$. The broad major arcs contribute
$O(X^{(S-k)a/S}A^{u/S})$ by Lemma~\ref{lem:nf-ordinary}. On their
complement use \eqref{eq:aud-field-critical}, the $u$th moment and $A$
with H\"older weights $a/w,(w-S)/w,u/w$. On dividing by $X^{S-k}$
the contribution is at most
\[
 X^{-\delta+\epsilon}(A/X^{S-k})^{u/w},
       \qquad \delta=\bigl[S-u-(w-S)D_u\bigr]/w.
\]
When $\delta>0$, absorption gives $A\ll X^{S-k}$ and a power saving
in the mixed minor integral. The sharp estimate on whole boxes, including
cells near a coordinate plane, follows at the end of
Section~\ref{subsec:aud-full-moments}.

The full upper-smooth inputs at degrees $4,8,9,10$ and the marked
ones at the intervening degrees are as follows. Each displayed defect
is an upper bound. Sections~\ref{subsec:aud-full-moments}--\ref{subsec:aud-high-moments}
give the field estimates and their finite rational deductions.
\[
\begin{array}{c|r|l|r}
 k&u&D_u\text{ at most}&S\\\hline
 4&10&393/1600&12\\
 5&16&7208133/10^8&18\\
 6&22&66943/625000&24\\
 7&30&4971653/51000000&34\\
 8&34&5857852793/50000000000&38\\
 9&40&8516302869/62500000000&46\\
10&46&152339/1000000&56
\end{array}
\]
All seven values of $\delta$ are positive. The degree-seven entry also
uses $|G|\ll X$ to extend the bound at $1456/51$ to order $30$.

For $a$ ordinary and $b$ marked roots the minor saving obtained directly
from the critical moment is $a/w-(1-a/w)D_r$ at $r=b/(1-a/w)$.
On the part of the broad major arcs outside the concentrated core, the
loss-free marked $S$th moment and the pruned ordinary moment of order
$p=aS/(S-b)>2k$ give a power saving by H\"older. For $(a,b)=(S-u,u)$
one may use the mixed minor bound just proved and $p=S$. This treats
$k=4,6,8$ and the preliminary 46 nonics. For the remaining allocations
the following are lower bounds for the exponent saving before
arbitrarily small parameter losses:
\[
\begin{array}{c|rr|c|c|r|c}
 k&a&b&r&\text{saving at least}&S&p\\\hline
 5&2&15&225/14&5380409/3000000000&18&12\\
 7&5&26&1456/51&28347/56000000&34&85/4\\
 9&7&37&3330/83&8664641/1800000000&46&322/9\\
10&5&45&330/7&92718217/220000000000&50&50
\end{array}
\]
The decic estimate and the loss-free order $50$ follow from the
marked-moment calculation in Section~\ref{subsec:refined-small-moments}.

On the central arcs the free cofactor produces the complete local kernel
in every marked slot. The prime windows and all dilations can be fixed
over the finite list of estimates. Use a small fixed dilation for the
marked roots and ordinary boxes whose power-sum map contains the residual
target in its interior. At a complex place restrict to an interior
sector for the lower bound. A finite cover of the compact set of normalized
target shapes gives uniform positivity by the central calculation in
Section~\ref{sec:nf-global}. Local primitive solutions in the full
number of slots give a positive singular series by that calculation.
If a local threshold exceeds the analytic one, add ordinary factors,
using their trivial bound on the minor and pruned arcs; primitive
solutions in the larger number of slots suffice. This proves all the
global assertions.
\end{proof}

\subsection{Comparison with the classical Waring number}
\label{subsec:nf-classical-comparison}

An allocation of $s$ summands satisfying the estimates of the preceding
proof gives $G_{\mathbb K}(k)\le\max\{H_{\mathbb K}(k),s\}$. The balanced
archimedean boxes give a positive singular integral. At
$\max\{H_{\mathbb K}(k),s\}$ summands the primitive local solutions give
a positive singular series. Each added ordinary factor is $O(X)$,
where $X$ is the root volume, so the power saving off the central arcs
persists.

To obtain $G_{\mathbb K}(k)\le\max\{H_{\mathbb K}(k),G(k)\}$ by this
argument, the field estimates must give an admissible allocation
within that total. The ordinary critical order $w=k(k+1)$ is common
to the rational and field arguments. The marked moments over the
trace torus and the ordinary estimates on the pruned arcs must be
established over the field. For $a$ ordinary and $b$ marked summands,
the minor-arc exponent margin is
\[
 \frac{a}{w}-\left(1-\frac{a}{w}\right)D_r,
 \qquad r=\frac{bw}{w-a}.
\]
In particular the argument needs a marked-moment defect bound
$D_r<a/(w-a)$.

\Needspace{7\baselineskip}
For ninth powers the classical allocation $(a,b)=(5,38)$ gives
$G(9)\le43$. At $r=684/17$, the cutoff is $1/17$. The defect
bounds in Section~\ref{subsec:refined-small-moments}, rounded upward
in the first numerical column, yield
\[
\begin{array}{lcc}
 & \text{Upper bound for }D_r & \text{Implied exponent margin}\\ \hline
 \text{Classical} & 0.058112 & \phantom{-}0.000672224\ldots\\
 \text{Number field} & 0.067596 & -0.008284788\ldots
\end{array}
\]
The same split meets the field's local requirement
$H_{\mathbb K}(9)\le25$. The loss-free marked moment at $S=46$ gives
the ordinary pruning order $5S/(S-38)=115/4>18$. Over all splits
with at most $43$ summands, the largest margin obtained by
interpolating the certified field defects is about $-0.0078854$,
at $(3,40)$. At $44$ summands the same bounds give positive margins
for $1\le a\le10$ and $b=44-a$. The loss-free marked moment of
order $44$ gives ordinary pruning at order $44a/(44-b)=44$
in these splits. Thus the available estimates certify
$G_{\mathbb K}(9)\le44$.

\subsection{A lower bound over totally real fields}

The difference of one between a primitive local width and the
ordinary local width matters even when the analytic threshold is
smaller. Davenport's quartic phenomenon admits the following
number-field form.
\begin{theorem}\label{thm:nf-quartic-lower}
For a totally real number field $\mathbb K$ one has
$G_{\mathbb K}(4)=16$ if and only if some dyadic completion is
isomorphic to $\mathbb Q_2$.
More generally, if a completion is $\mathbb Q_p$ and $k\ge6$ is even,
then
\[
 G_{\mathbb K}(k)\ge
 \begin{cases}
 2^{v_2(k)+2},&p=2,\\
 p^{v_p(k)+1},&p>2,\ p-1\mid k .
 \end{cases}
\]
In particular the dyadic lower bound equals $4k$ at every power of
two $k\ge8$.
\end{theorem}
\begin{proof}
Write $B$ for the stated lower bound. At the chosen completion of
$\mathbb K$, every unit $k$th power is $1$ modulo $B$, and every
nonunit $k$th power is zero modulo $B$. This follows from the
exponent of the unit group of $\mathbb Z/p^j\mathbb Z$. The same
statements hold at $k=4$, $p=2$, with $B=16$.
Consequently, a sum of at most $B-1$ terms that is zero modulo $B$
can contain only nonunits.

Choose $e\ge1$ and a totally positive $c$ with
$(c)=\mathfrak p^e$, where $\mathfrak p$ is the chosen prime of
$\mathbb K$. Such a choice exists by finiteness of the narrow class
group. For a fixed positive rational integer $b$, consider
$bc^{kr}\in J_k$. Its norm tends to infinity. In a representation
with at most $B-1$ terms, apply the preceding observation
successively at $\mathfrak p$. After dividing the local roots by
$p^j$, $0\le j<er$, the scaled target is still divisible by $p^k$,
since $v_{\mathfrak p}(p)=1$. As $B\mid p^k$, every original root
belongs to $\mathfrak p^{er}=(c^r)$. Dividing globally by $c^r$
would give a representation of $b$ by at most $B-1$ totally
positive integral roots.

For $k\ge6$ take $b=B$. At $p=2$ one has $B\le4k<2^k$. At odd
$p$ one has $p\le k+1$ and $B\le k(k+1)<2^k$.
Each nonzero root $y$ therefore has $0<\sigma(y)<2$ at every
embedding. Since $|N(y-1)|<1$, integrality gives $y=1$. Hence $B$
summands are required.

At $k=4$, take $b=31$. Each nonzero root satisfies
$0<\sigma(y)\le31^{1/4}$. The conjugates of $u=y-1$ thus lie in
$(-1,31^{1/4}-1]$, a subinterval of $(-2,2)$. Every root of
$T^2-uT+1$ is an algebraic integer all of whose conjugates have
absolute value one. Kronecker's theorem gives
$u=\zeta+\zeta^{-1}$ for a root of unity $\zeta$. An order $m\ge8$
gives a conjugate $2\cos(2\pi/m)\ge\sqrt2>31^{1/4}-1$.
For $m=1,2,3,5,7$ at least one conjugate falls outside the allowed
interval. The remaining orders four and six give $y=1$ and $y=2$.
Since $31=16+15$ and $2\cdot16>31$, a representation requires at
least sixteen summands. The converse follows from
Theorem~\ref{thm:nf-small}.
\end{proof}

For example, put $F_j=2^{2^j}+1$ and
$\mathbb K_r=\mathbb Q(\sqrt{F_2},\ldots,\sqrt{F_{r+1}})$, with
$\mathbb K_0=\mathbb Q$. These totally real fields have degree $2^r$
and $G_{\mathbb K_r}(4)=16$. The fields $\mathbb L_r=\mathbb K_r(\sqrt5)$ have degree
$2^{r+1}$ and $G_{\mathbb L_r}(4)\le12$. The contrast persists at
arbitrarily large degrees. Indeed the Fermat numbers are pairwise
coprime, and $F_j$ lies strictly between two successive squares.
Each contributes a prime of odd valuation to a new square class.
For $j\ge2$, $F_j\equiv1\bmod8$ is a square in $\mathbb Q_2$,
so every dyadic completion of $\mathbb K_r$ is rational. The additional square class of $5$ is independent
and gives the unramified quadratic dyadic extension.

\subsection{Six fourth powers in quadratic fields}
\label{sec:nf-quartic-density}
The quartic threshold for density can be lower than the threshold
for all sufficiently large targets. We use the relative density in
the eligible lattice defined in Section~\ref{sec:nf-power-density}.

\begin{theorem}\label{thm:nf-six-quartics}
For every number field and every $s\ge6$, the relative density of
targets represented by at most $s$ totally positive fourth powers
equals the product of the ordinary local represented proportions
in Corollary~\ref{cor:nf-ordinary-density}, with a power-saving
error in each fixed balanced box. In particular
\[
 A_{\mathbb K}(4)\le\max\{6,W_{\mathbb K}(4)\}\le15 .
\]
If some completion is $\mathbb Q_2$, then $A_{\mathbb K}(4)=15$.
If $2$ splits completely, the density for $6\le s\le15$ is
$(s/15)^{[\mathbb K:\mathbb Q]}$.

For a quadratic field the density at six roots, and its density-one
threshold, are as follows.
\[
\begin{array}{c|c|c}
 \mathcal O_{\mathbb K}\otimes\mathbb Z_2 &
  \text{density at six}& A_{\mathbb K}(4)\\\hline
 \mathbb Z_2\times\mathbb Z_2 & 4/25&15\\
 \text{unramified quadratic} &113/128&7\\
 \mathbb Z_2[\sqrt{2u}],\quad u\in\mathbb Z_2^\times &19/20&7\\
 \mathbb Z_2[i]&61/64&7\\
 \mathbb Z_2[\sqrt{-5}]&1&5\le A_{\mathbb K}(4)\le6
\end{array}
\]
Over real quadratic fields ordered by discriminant, the mean of
the individual-field densities at six roots is $12839/19200$.
\end{theorem}
\begin{proof}
Use one ordinary and five restricted factors in the Fourier
coefficient for six variables. Their squared product on the minor
and annular arcs is $|f|^2|F|^{10}$, precisely the integrand
estimated with $S=12$, $u=10$ in the proof of
Theorem~\ref{thm:nf-small}. The margin is
$2-(20-12)(393/1600)>0$. Here one may use the cofactor construction
with a single layer of $L$ marked primes, all in
$X^a<N\mathfrak p<X^b$, with $(L+1)^{-1}<a<b<1/L$ and $b$ below
the ambient smoothness exponent. The cofactor volume is less than
$X^a$ and greater than a fixed power of $X$, so the same bound
follows by coefficient domination from the smooth moments.
The broad annulus is bounded by H\"older with the loss-free
twelfth moment and the pruned ordinary twelfth moment. The core
expansion uses this unrestricted cofactor. For $s\ge6$, append
$s-6$ ordinary variables and bound them trivially off the core.

For clarity the untwisted singular series converges with a power
tail whenever $s>\max(4,k+1)$. At a prime of residue norm $q$ and
characteristic prime to $k$, summation in the last residue digit
annuls the unit shell at every conductor exponent $j\ge2$.
Integration on the nonunits gives the normalized sum $q^{-1}$
for $2\le j\le k$, and $q^{-1}$ times a sum with exponent $j-k$
for $j>k$. At $j=1$ use $O_k(q^{-1/2})$. The prime-level terms
are $O_k(q^{1-s/2})$ after summing characters. Terms of
exponent $2$ through $k$ are $O_k(q^{k-s})$ in total, and every
following block of $k$ exponents introduces the geometric factor
$q^{k-s}$. These bounds allow an extra factor $q^{\epsilon j}$
for sufficiently small $\epsilon>0$. The finitely many primes
of residue characteristic dividing $k$ are covered by the complete
bound of \eqref{eq:nf-complete-bound} since $s>k$.
Thus $s=6$ suffices at $k=4$. The uniform local truncation and Bessel
proof of Corollary~\ref{cor:nf-ordinary-density} applies from this
order onward.

At odd primes the primitive width is at most five. At every dyadic
place other than $\mathbb Q_2$ it is at most fifteen by the local
part of Theorem~\ref{thm:nf-small}. In $\mathbb Z_2$, unit fourth
powers fill $1+16\mathbb Z_2$. At $s\le15$ the primitive classes
have residues $1,\ldots,s$ modulo $16$, and the remaining local
representations have every root even. After scaling by $16$, the
same dichotomy recurs, so the relative proportion is $s/15$.
At $s=15$ the entire ring is represented, including zero.
This proves the first assertions and the split row.

In an unramified quadratic dyadic ring choose the basis
$(1,\zeta_3)$. Modulo $16$ the available unit classes are $jv$,
where $j=1,9$ and $v=(1,0),(0,1),(-1,-1)$, by
Corollary~\ref{cor:local-triangular}. The nonunit class is zero
and the eligible lattice is the entire ring. Starting with
$C_0=\{(0,0)\}$ and setting
$C_{i+1}=C_i\cup\bigcup_{j,v}(C_i+jv)$ in
$(\mathbb Z/16\mathbb Z)^2$, the cardinalities for $0\le i\le7$
are
\[
                    1,7,25,58,106,166,226,256.
\]
Every listed class at six roots lifts fully, since one unit
fourth power absorbs $16\mathcal O_v$; zero has the three directions
whose sum is zero. Hence the proportion is $226/256$.

For the four even-radicand completions use the basis of
Corollary~\ref{cor:even-radicand-widths}. For six roots the
primitive classes have scalar residue $r+4v\pmod{16}$, with
$r\ge1$, $v\ge0$, $r+v\le6$, and the other eligible coordinate
is arbitrary. The only missing scalar residue is zero.
The nonunit-root part is $(2u)^2S_{6,v}$. To contribute a new
class its original scalar coordinate must be divisible by four.
All scalar coordinates outside that class are already represented.
Since the absolute measure multiplier is $1/16$, the relative
proportion $x$ satisfies $x=15/16+(x-3/4)/16$, hence $x=19/20$.

For $\mathbb Z_2[i]$, the basis of
Corollary~\ref{cor:odd-radicand-widths} gives scalar residue
$r-4v\pmod{16}$, with just $8$ missing. Scaling the nonunit part
by $(1+i)^4=-4$ reaches this missing residue exactly from scalar
coordinates $2\pmod4$, all of which already belong to the
primitive part. Thus the proportion is $15/16+(1/4)/16=61/64$.
For $\mathbb Z_2[\sqrt{-5}]$ the ordinary width is five by the
same corollary. The local widths and the upper and lower bounds
in Corollary~\ref{cor:nf-ordinary-density} now give the last
column. Finally the five rows have discriminant frequencies
$1/3,1/3,1/6,1/12,1/12$, respectively, by the squarefree count
in Section~\ref{sec:quadratic-density}. The weighted sum is
$12839/19200$.
\end{proof}

The same calculation gives the following orders for ordinary local
density in the other small cases:
\[
\begin{array}{c|rrrrrr}
 k&5&6&7&8&9&10\\\hline
 t&9&12&16&19&22&25 .
\end{array}
\]
From each indicated $t$ onward, the product in
Corollary~\ref{cor:nf-ordinary-density} is the global relative density,
with a power-saving error on boxes. Consequently
$A_{\mathbb K}(k)\le\max(t,W_{\mathbb K}(k))$.
The squared Bessel integrals use the even ordinary and marked
multiplicities
\[
       (4,14),\quad(4,20),\quad(6,26),\quad
       (4,34),\quad(6,38),\quad(6,44)
\]
in increasing degree. At $k=5$ interpolate at $r=210/13$ between
orders 16 and the loss-free 18. The saving is at least
$2/15-(4/5)D_{16}>0$. At $k=6$ the corresponding interpolation at
$420/19$ between 22 and 24 gives $2/21-(6/7)D_{22}>0$.
Both ordinary annular orders equal the corresponding loss-free order.
At $k=7$ add one bounded ordinary factor to the global mixed
estimate. The octic estimate already has even multiplicities. At
$k=9$, Section~\ref{subsec:aud-high-moments} gives
$D_{285/7}\le1325421/20000000$, saving at least $722053/150000000$,
and annular order $6\cdot46/(46-38)=69/2>18$.
At $k=10$, Section~\ref{subsec:refined-small-moments} gives the
saving $1105388504469/275000000000000$ for $(6,44)$ with
loss-free marked and ordinary pruning orders both $50$.
Since $W_{\mathbb K}(10)\le H_{\mathbb K}(10)\le17$,
it follows that $A_{\mathbb K}(10)\le25$.
In each case the half-order exceeds $\max(4,k+1)$, as required for
the ordinary local product in Corollary~\ref{cor:nf-ordinary-density}.
These very same marked sums consequently have the sharp even orders
\begin{equation}\label{eq:aud-small-even}
 \begin{array}{c|rrrrrrr}
 k&4&5&6&7&8&9&10\\ \hline
 S_k&12&18&24&32&38&44&50 .
 \end{array}
\end{equation}
Indeed enlarge half the indicated ordinary multiplicity on each side
of the pure even orthogonality count. On the full broad region use
the ordinary $p$th moment with the loss-free marked $S$th moment.
Together with the mixed minor estimate this bounds the pure
$S_k$th moment by $O(X^{S_k-k})$. For $k=7$ the extra ordinary
factor costs $O(X)$ also on the full broad arcs.

For example, if $\mathbb K$ is quadratic and $2$ splits, then
$A_{\mathbb K}(8)=32$. At $19\le s<32$ the global relative density
is $(8s/255)^2$: at each dyadic place the $s$ primitive residues
modulo $32$ have a disjoint recurrence on multiplying by $2^8$,
and all odd places have primitive width at most nine. If $2$ does
not split then $A_{\mathbb K}(8)\le19$, since the dyadic quadratic
width is at most seventeen by
Corollaries~\ref{cor:local-triangular}, \ref{cor:even-radicand-widths}
and~\ref{cor:odd-radicand-widths}.
\subsection{Primitive local representations at the small exponents}
\label{subsec:aud-small-local}
Let $A$ be the integer ring of a finite extension of $\mathbb Q_p$,
with uniformizer $\pi$, ramification index $e$, and residue field
$F$ of cardinality $q=p^f$. At $p\nmid k$ the residue sumset bound
in \eqref{eq:nf-residue-width} and simple Hensel lifting give
$h_k(A)\le k+1$. Theorem~\ref{thm:prime-degree-width}
gives $H_{\mathbb K}(5)\le5$ and $H_{\mathbb K}(7)\le7$.
We prove the remaining bounds, retaining a unit coordinate even for
the target zero.

\paragraph{From Witt vectors to primitive representations.}
Write $k=p^{r-1}m$, $(m,p)=1$, and $[x]=(x,0,\ldots,0)$ in $W_r(A)$.
The last ghost map is additive and takes $[x^m]$ to $x^k$.
Thus every eligible target has a lift in the additive group generated
by the $[x^m]$. For fixed roots $b_i$, one of them a unit, put
\[
 \Phi(\boldsymbol z)=\sum_{i=1}^r[(b_i+z_i)^m]-\sum_{i=1}^r[b_i^m].
\]
The Witt coordinates of $\Phi$ are integral polynomials. Differentiating
the ghost identities at zero gives the entry
$m b_i^{m-1}(b_i^m)^{p^j-1}$ in row $j$, with the zeroth power
interpreted as one. If $t=v_\pi(\det J)$ for this Jacobian, then
$C=\pi^tJ^{-1}$ is integral and
\[
 \pi^{-2t-1}\Phi(\pi^{t+1}C\boldsymbol v)
       =\boldsymbol v+\pi P(\boldsymbol v)
\]
in Witt coordinates, with integral $P$. Successive approximation
fills the kernel of reduction to $W_r(A/\pi^{2t+1})$ and preserves
every unit baseline. More generally, for row thresholds $N_j$
with $pN_j\ge N_{j+1}$, write $I_N$ for the ideal of vectors with
$j$th coordinate in $\pi^{N_j}A$. If a weighted choice of the
perturbations gives an invertible linear part modulo $\pi$ after
the rows are divided by $\pi^{N_j}$ and every higher coefficient
is divisible by $\pi$, the same argument fills $I_N$.

In either case first subtract the \emph{entire} fixed baseline from the
Witt lift. Let $P$ now denote the set of permitted Teichm\"uller
steps in the finite quotient and set $B_0=\{0\}$, $B_{j+1}=B_j+P$.
Since $0\in P$, the first equality $B_D=B_{D+1}$ proves that
$B_D$ is the whole group generated by $P$ (a finite additive
semigroup containing zero is a group). Match the residual target
there and vary the reserved roots to fill its error. This gives
$h_k(A)\le D+r$. All finite computations below are specified by
this recurrence or its primitive variant, with implementations in
Appendix~\ref{app:certificates}.

\paragraph{Fourth powers.}
At $p=2$ use $W_3(A)$. For $f\ge2$ take centres $0,1,b$ with
$\bar b\notin\{0,1\}$. The determinant $b^3-b$ is a unit.
The expansion $[u]+2[v]+4[w]$ in $W_3(F)$ costs at most seven
positive steps, so $h_4(A)\le10$. For $f=1$ and $e\ge2$, take
$0,1,\pi$. Now $t=1$, and the two possibilities for $A/\pi^3$ are
\[
 (\mathbb Z/4)[T]/(2T,T^2-2)\quad(e=2),\qquad
 \mathbb F_2[T]/(T^3)\quad(e\ge3).
\]
In the first ring $\pi^2/2$ has residue one. In both rings the update is
\[
 (a,b,c)+[v]=(a+v,b-av,c-a^3v-2a^2v^2-av^3+abv).
\]
For $0\le j\le12$, this gives respectively
\[
 |B_j|=1,8,28,56,88,120,152,184,214,236,250,256,256 .
\]
Thus $h_4(A)\le14$. Over $\mathbb Z_2$ the unit fourth powers fill
$1+16\mathbb Z_2$, giving primitive width $16$ (zero modulo $16$
requires at least $16$ units). To obtain the ordinary width $15$,
divide a nonzero target by a fourth power of $2$ to put its valuation
below four and use its residue in $1,\ldots,15$. The class $15$
requires at least fifteen units. Together with $h_4(A)\le5$
at odd places this proves all quartic local assertions of
Theorem~\ref{thm:nf-small}.

\paragraph{Sixth and tenth powers.}
At each wild prime $k=pm$. If $1$ and a unit $x$ have
$\bar x^{m(p-1)}\ne1$, they give a unit determinant for $W_2(A)$.
Put $E=\mathbb F_p((F^\times)^m)$, the subfield generated by the residue
directions. The finite step group is contained in $W_2(E)$. If $L$
allowed residues suffice for each element of $E$, cover the first
digit with $L$ steps and the second with $pL$, using
$p[a]=(0,a^p)$ and the bijectivity of Frobenius. The roots therefore
cost at most $2+(1+p)L$. At $p=2$, \eqref{eq:hcoeff} gives $L\le m$.
At odd $p$ we have $m=2$, and every element of $F$ is a sum of two
squares, since two sets of size $(q+1)/2$ intersect. For $f\ge2$
we also have $q-1>2(p-1)$, as required for the two unit centres.

At $p=2$ the exceptions have $q-1\mid m$, giving $q=2,4$ for
$k=6$ and $q=2$ for $k=10$. Use the roots $1,1+\pi$, for which
$t=1$ since $m$ is odd. At the odd wild prime with $F=\mathbb F_p$,
use the roots $1,\pi$. The second column has factor $2\pi$ and
the remaining determinant is a unit. The finite steps are $[v]$,
$v=x^m$, in $W_2(A/\pi^3)$, with update
\[
 (a,b)+[v]=\left(a+v,\ b-\sum_{i=1}^{p-1}
           \frac{\binom{p}{i}}p\,a^i v^{p-i}\right).
\]
Writing $B=W(F)$, all necessary coefficient rings are
\[
 B/p^3B\quad(e=1),\qquad
 (B/p^2B)[T]/(pT,T^2-p\widehat u)\quad(e=2),\qquad
 F[T]/(T^3)\quad(e\ge3).
\]
Here $u=\overline{\pi^2/p}\in F^\times$ for $e=2$. Rescaling
$\pi$ normalizes $u$ to one over $\mathbb F_4$. The presentations
follow by expansion in a uniformizer over the unramified coefficient
ring. The program \path{code/ed/witt_small.py} evaluates this recurrence
and the preceding three-coordinate quartic recurrence. The following
diameters $D$ include every unit class $u$ at $e=2$.
\[
 \begin{array}{c|c|c|c|c}
 k&p&q&D\ (e=1)&D\ (e\ge2)\\\hline
 6&2&2,4&10&6\\
 10&2&2&10&6\\
 6&3&3&\text{see below}&9\\
 10&5&5&\text{see below}&12
 \end{array}
\]
In $\mathbb Z_3$ the unit sixth powers fill $1+9\mathbb Z_3$, and
in $\mathbb Z_5$ the unit tenth powers fill the two classes
$\pm1+25\mathbb Z_5$. These claims follow from the principal-unit
logarithm and the prime-to-$p$ unit group. The former primitive
width is $9$. In the latter ring every nonzero residue modulo $25$
has a signed representative of absolute value at most $12$ and
zero can use the two opposite units, giving width $12$.
We obtain $H_{\mathbb K}(6)\le12$ and $H_{\mathbb K}(10)\le17$.

\paragraph{Ninth powers.}
For $p=3$ we use the characteristic-three form of the chart in
the proof of Theorem~\ref{thm:hhybrid}. Put $R=F[[T]]$. The
unramified coefficient ring $W(F)\subset A$ induces $F\to A/3A$,
and expansion in powers of $\pi$ gives a surjection
$R\to A/3A$ taking $T$ to $\pi$. To see why the
chart in $W_3(R)$ suffices, lift the individual roots in an eligible
signed sum through this map. Equality of their Witt sums over $R$
gives equality of the last ghosts of any lifts to $A$ modulo $27A$.
The last ghost has terms $x_0^9,3x_1^3,9x_2$, and a change
of any coordinate by $3A$ leaves each term unchanged modulo $27A$.
One reserved unit absorbs the error since $(1+3A)^9=1+27A$.

Use the centres $0,T,1$, perturbation orders $2,2,4$ and row
thresholds $N=(2,4,4)$. After normalization the derivative modulo $T$ is
\[
 \begin{pmatrix}1&1&0\\0&1&1\\0&0&1\end{pmatrix}.
\]
At a centre $T^i$ and perturbation order $n>i$, a degree-$d$ term,
$d\ge2$, of the ghost perturbation gains $(d-1)(n-i)$ over its
linear order. At zero the ghost orders are $2\cdot3^j>4$ for $j>0$.
At the unit the higher terms have order at least eight. Lower Witt
coordinates entering the ghost recursion have order at least
$3^{j-l}N_l>N_j$ for $l<j$. These comparisons may be made over
integral polynomial lifts before reducing modulo three, since division
by a numerical power of three preserves the $T$-exponents. Thus the
chart fills $I_N$.

For $M=W_3(R)/I_N$, the lengths of $3^sM$ as $W(F)$-modules are
$L_s=\sum_{j=s}^2\lceil N_j/3^s\rceil$, giving $(10,4,1)$.
As in \eqref{eq:nf-uniform-witt}, a filtration by actual generators
with cyclic lengths $a_i$ satisfies $\sum_i(a_i-s)_+\le L_s$.
Multiplication by a constant Teichm\"uller digit preserves the
step set. The positive cost is therefore
\[
 \sum_i(3^{a_i}-1)/2\le L_0+2L_1+6L_2=24.
\]
This already proves the table, with $24+3=27$ roots.
For $F=\mathbb F_3$, allow the third centre first to be $1+\eta$,
$\eta\in T^2F[[T]]/T^4$. Modulo $I_N$ the centred difference
$[1+\eta]-[1]$ has representative $(0,\eta,\eta)$ because the first
coordinate can be removed by $[\eta]\in I_N$ and every term
quadratic in $\eta$ has order at least four. These vectors form a
two-dimensional additive subgroup $D$ killed by three, hence a
$W(\mathbb F_3)$-submodule. In $M/D$, $L_0$ decreases by two
and the other $L_s$ cannot increase. Match the word modulo $D$,
choose $\eta$ to match its discrepancy there, and fill $I_N$ using
the same normalized derivative. This gives $22+3=25$.
For $f\ge2$ the centres $0,1,u$ with $\bar u\notin\mathbb F_3$
have unit determinant $\bar u^8-\bar u^2$, since in characteristic
three $\bar u^6-1=(\bar u^2-1)^3$. Their residue quotient $W_3(F)$
costs $1+3+9=13$ steps. We have proved $H_{\mathbb K}(9)\le25$.

\paragraph{Eighth powers with residue degree at least two.}
Work in $W_4(A)$ at $p=2$. For $f\ge3$ take $0$ and three units
whose residues are linearly independent over $\mathbb F_2$.
Removing the zeroth row and multiplying the unit columns by their
residues leaves a Moore matrix with powers $2,4,8$. The determinant
is a unit. Four binary Teichm\"uller digits cost $1+2+4+8=15$,
giving $h_8(A)\le19$.

Let $f=2$, and let $\omega$ be a primitive cube root of unity.
For $e\ge3$ take $0,\pi,1,\omega$ with perturbation orders $2,2,3,3$
and $N=(2,3,3,3)$. The first two rows have pivots at $0,\pi$, and
the remaining $2$-by-$2$ block has determinant
$\bar\omega^7-\bar\omega^3\ne0$. Nonlinear terms have higher order
by the preceding ghost calculation, with $2N_j>N_{j+1}$.
All retained coefficient digits are in characteristic two. More
precisely $A/\pi^3\simeq F[T]/(T^3)$, so $W_4(A)/I_N$ is obtained
from $W_4(F[[T]])$ with the same thresholds. Thus
$L_s=\sum_{j=s}^3\lceil N_j/2^s\rceil$ gives $11,6,2,1$.
The positive cost with the fixed roots is $11+6+2\cdot2+4+4=29$.

At $e=1$, Corollary~\ref{cor:local-triangular} gives $15$.
For $e=2$ work modulo $\pi^9$. The set of principal-unit deviations
$D=(1+\pi A)^8-1$ has valuation at least six by binomial expansion,
and forms an additive subgroup modulo $\pi^9$ since all products
of two deviations vanish there. Put $\eta=2/\pi^2$. For $u=1+\pi^2z$
the $\pi^8$ digit of $u^8-1$ is
$\bar\eta^2(\bar z^2+\bar\eta\bar z)$, which is nonzero for some
$\bar z\in\mathbb F_4$. Its $\mathbb F_2$-line and its
$\bar\omega$-multiple span $\mathbb F_4$. Thus $D+\omega D$
contains the entire $\pi^8$ shell, in particular every nonunit
eighth power. Every unit
eighth power has the form $\zeta(1+d)$ with
$\zeta\in\{1,\omega,\omega^2\}$ and $d\in D$.
Since $\omega^2=-1-\omega$, its deviation also belongs to $D+\omega D$.
Subtract the two adjustable unit baselines $1,\omega$. Modulo
$D+\omega D$ the remaining signed target is a signed sum of
$1,\omega,\omega^2$, all eighth powers in $\mathbb Z_2[\omega]$.
Corollary~\ref{cor:local-triangular} supplies at most fifteen roots
there. Vary the two reserved units to realize the residual deviation.
One absorbs the final error in $\pi^9A$ by the logarithm on
$1+\pi^3A$. This proves $h_8(A)\le17$.

\paragraph{Eighth powers over the residue field $\mathbb F_2$.}
In $\mathbb Z_2$ unit eighth powers fill $1+32\mathbb Z_2$, so the
primitive width is $32$. For $2\le e\le5$ put $c=4e+1$. The
binomial theorem gives $x^8\bmod\pi^c$ from $x\bmod\pi^{e+1}$,
and varying a unit root fills its coset modulo $\pi^c$ by the
logarithm on $1+\pi^{e+1}A$. A uniformizer generates the full
totally ramified integer ring over $\mathbb Z_2$. Write
\[
 \pi^e=2\sum_{i=0}^{e-1} b_i\pi^i,\qquad b_0\text{ odd}.
\]
It suffices to enumerate each $b_i$ modulo
$2^{\lceil(c-e-i)/e\rceil}$, since the $i$th coefficient error then has
valuation at least $c$. Coordinate reduction is unambiguous because
the valuations $e v_2(a_i)+i$ of nonzero terms $a_i\pi^i$ are
distinct modulo $e$. For odd $e$, Hensel lifting shows that every
unit is an $e$th power, allowing the choice $\pi^e=2$.
In $A/\pi^c$ form $P$ from all the $2^{e+1}$ root residues and
let $U\subset P$ be the unit steps. Start at $C_1=U$, iterate
$C_{j+1}=C_j+P$. Stabilization is the entire generated group,
since in a finite group translates of its full positive semigroup
by a member of the group give the same group. A unit then lifts
every represented class. This gives the following exact primitive widths.
\[
\begin{array}{c|l|l}
 e&\text{relations tested}&
 h_8(A)\ \text{(number of coefficient choices)}\\\hline
 2&64&16\ (16),\quad17\ (48)\\
 3&\pi^3=2&35\\
 4&4096&5\ (64),\quad7\ (480),\quad8\ (1504),\quad12\ (2048)\\
 5&\pi^5=2&21 .
\end{array}
\]
The recurrence is implemented in \path{code/local/octic_dvr.cpp}.
The frequencies count polynomial reductions, which can present
isomorphic fields. For $\pi^3=2$ the ordinary width is also $35$.
Modulo $\pi^{13}$ the coordinate moduli in the basis $1,\pi,\pi^2$
are $32,16,16$. The sixteen root classes give unit steps
$(1,4j,0)$, $j\bmod4$, and nonunit steps $(0,0,0),(0,0,4)$.
The target $4\pi+12\pi^2$ requires a unit to supply its $\pi$ digit,
then at least $32$ units to clear the scalar digit, and at least
three nonunits for its $\pi^2$ digit. This also gives a second
example in the equality branch of \cite[Conjecture~3.3]{AnthonyRam},
since $w_8(\mathbb Z_2)=32$ (the target $32$ forces units).

Finally suppose $e\ge6$. Take centres $0,\pi^2,\pi,1$ with
perturbation orders $3,3,4,7$ and $N=(3,5,7,7)$. The normalized
derivative modulo $\pi$ is
\[
 \begin{pmatrix}1&1&0&0\\0&1&1&0\\0&0&1&1\\0&0&0&1\end{pmatrix}.
\]
At each nonzero nonunit centre $T^i$ the degree-$d$ gain is
$(d-1)(n-i)>0$ for perturbation order $n$. At zero,
$3\cdot2^j>N_j$ for $j>0$, and at the unit all nonlinear terms
have order at least fourteen. Since $2N_j>N_{j+1}$, the lower
ghost powers have strictly higher order also. Making these comparisons
over $\mathbb Z[T]$ and then putting $T=\pi$ proves that the chart
fills $I_N$. For $e=6$ the ring $A/\pi^7$ has presentation
$(\mathbb Z/4)[T]/(2T,T^6-2)$, since $\pi^6/2$ has residue one.
For $e\ge7$ it is $\mathbb F_2[T]/(T^7)$.

For these two rings the additive computation may be reduced to
22 generators $V^j[T^d]$, $0\le j\le3$, $0\le d<N_j$, where
$V$ is the Witt-coordinate shift. Successively reducing the first
nonzero coordinate in uniformizer digits gives unique binary normal
forms. To compute the relations for twice a generator, use integral
ghost lifts. In row $s\ge j$ the ghost of $V^j[x]$ is
$2^j x^{2^{s-j}}$, and it vanishes in the other rows. After
subtracting the axes for a reduced digit, remove its remaining
multiple $t$ of $T^{N_j}$ by subtracting the \emph{full} vector
$V^j[t]\in I_N$, thus retaining all carries in higher coordinates.
The doubling relations are triangular and present the group, since
its binary normal forms already number $2^{22}$. Diagonalize the
relation matrix over $\mathbb Z$ and run the addition recurrence in
these cyclic coordinates. The programs \path{code/local/octic_certificates.py}
and \path{code/local/octic_distances.cpp} implement the two operations.
The 128 root classes give 70 distinct
steps in each of the two rings. Their sumsets have sizes
\[
 |B_{20}|=130316,\qquad |B_{21}|=|B_{22}|=131072.
\]
Hence $h_8(A)\le21+4=25$ throughout $e\ge6$. All residue
degrees and ramification indices are covered, proving
$H_{\mathbb K}(8)\le35$, with equality whenever a completion is
$\mathbb Q_2(\sqrt[3]{2})$.
\subsection{Uniform estimates for the full upper-smooth sum}
\label{subsec:aud-full-moments}
In the analytic proofs below put $K=\mathbb K$, $n=[K:\mathbb Q]$,
and $e_\alpha(y)=e_{\mathbb K}(\alpha y)$. In this section $P,Y,U,T$
as scale parameters measure \emph{volumes}; the radius of a balanced
box of volume $P$ is $P^{1/n}$. Put $U_r(Y)=\int|G_Y|^r$ and
$\lambda_r=r-k+D_r$. Each bound is uniform over subsets, over the
specified ideal lattices, and over coefficients bounded in modulus
by one. Constants in smooth envelopes can be removed by normalization.

The full upper-smooth family uses only an upper cutoff on prime ideal
norms. For the computations involving genuine prime blocks we use as
well a hereditary subfamily: after an exceptional ideal of norm at most
$P_0^{O(\eta)}$ is removed, the relative root ideal is squarefree with
eligible prime factors of norm in $(P_0^{\eta/2},P_0^\eta]$. For an even moment we
can group the exceptional ideal first, choose a balanced common scalar
times a fixed class representative and cancel that scalar; grouping
these ideals costs at most $P_0^{O(\eta)}$. Block division preserves the
resulting squarefree condition. In estimates using an existing marked
fibre decomposition the retained prime windows and the factorization
belong to the same root weight as in Section~\ref{sec:nf-fibres}.
On the central arcs choose the fixed windows of the positive marked
construction, with distinct prime factors and terminal free cofactor
of volume at most $P_0^{O(\eta)}$, as in \eqref{eq:nf-terminal-form}.
In the explicit choice of Section~\ref{sec:nf-fibres}, all primary
and auxiliary prime norms lie in $(P_0^{\eta/2},P_0^\eta]$, whereas
the cofactor norm is smaller than the least auxiliary prime norm. The
cofactor itself is therefore a permissible exceptional ideal.
This allows arbitrary cofactor values, including nonunits at any
specified finite set of places. At intermediate volumes $Y\ge P_0^\gamma$
the ambient cutoff $P_0^\eta$ is at most $Y^{\eta/\gamma}$.
The curves below use finitely many positive-power cutoffs, levels and
interpolations. Given a final loss, choose their tolerances, then $\gamma$
and $\eta$ in this order, small enough for every edge; at the bottom use
the trivial mass. Losses $P_0^{O(\eta)+o(1)}$ can thus be absorbed into
any prescribed initial $P_0^\epsilon$. The strict mixed minor margins
allow all these choices before the prime windows are fixed.

Write \(P_0\) for the initial volume and \(R=P_0^\eta\) for the upper bound on prime ideal norms. Use the finite set of integral narrow ideal-class representatives of Section~\ref{sec:ed}. For a representative \(L\), put \(I_L(x)=(x)L^{-1}\) for \(0\ne x\in L\). A moment bound at order \(r\) and exponent \(\lambda_r\) will hold uniformly for all coefficients \(|c_x|\le1\) supported on the permitted roots:
\[
 \int_{\widehat{\mathcal O}}\left|\sum_x c_x e_\alpha(x^k)\right|^r
       \,d\alpha \ll_{K,k,\epsilon} Y^{\lambda_r}P_0^\epsilon .
\]
The coefficients may encode any subset. Even initial moments have this property by orthogonality, and real interpolation preserves it. It is necessary to retain this uniformity at the later real orders.

\paragraph{Removing small and repeated prime factors.}
Fix an output order \(q>1\) and a proposed exponent \(\lambda=q-k+D_q\) with
\[
              \frac{q^2}{2q-1}<\lambda<q .
\]
The case \(\lambda\ge q\) is trivial. Put \(\sigma=\lambda/q\). Choose \(a,B>0\) such that
\[
    1-\sigma<a<\frac{\lambda-1}{q-1},\qquad
    \frac1a<B<\frac1{1-\sigma},
    \qquad W=(\log(2P_0))^B .
\]
The strict hypothesis is exactly what makes these two intervals nonempty. Split \(I_L(x)=\mathfrak c\mathfrak s\), putting all prime powers above primes of norm at most \(W\), and the entire prime powers of exponent at least two at every other prime, in \(\mathfrak c\). Thus \(\mathfrak s\) is squarefree and each of its prime norms is in \((W,R]\). Write \(L\mathfrak c=bL'\) with \(b\) a balanced scalar and \(L'\) a fixed integral class representative. The root \(y=x/b\in L'\) has \(I_{L'}(y)=\mathfrak s\). For every fixed \(\mathfrak c\) its subsum has the same \(L^q\) norm after this common division, also for nonintegral \(q\). Indeed, choose an integer \(A>0\) with \(A b^k\in\mathcal O\). If the original frequencies \(b^k y^k\) are integral, replacing \(\alpha\) by \(A\alpha\) preserves Haar measure. Since the frequencies \(y^k\) are integral, multiplication by the nonzero algebraic integer \(A b^k\) is a second Haar-preserving torus isogeny. Neither step requires a fractional scalar to act on the fixed torus.

By Minkowski the total cost of this division is at most
\[
 \sum_{\mathfrak c}(N\mathfrak c)^{-\sigma}
 \le \prod_{N\mathfrak p\le W}
             (1-(N\mathfrak p)^{-\sigma})^{-1}
     \prod_{N\mathfrak p>W}
             \left(1+
             \frac{(N\mathfrak p)^{-2\sigma}}{1-(N\mathfrak p)^{-\sigma}}\right)
 = P_0^{o(1)}.
\]
The last product is bounded since \(\sigma>1/2\), and the logarithm of the first is \(O_K(W^{1-\sigma})=o(\log P_0)\). This bound uses only the elementary ideal count \(O_K(X)\). For clarity let \(Z=Y/N\mathfrak c\). If the divided annulus contains a root then \(Z\gg_K1\), so its trivial mass is \(O_K(Z)\). The total norm from the annuli with \(Z<P_0^\iota\) is therefore at most
\[
  O_K(1)\sum_{\mathfrak c}\min(P_0^\iota,Y/N\mathfrak c)
  \ \ll\ Y^\tau P_0^{\iota(1-\tau)}
              \sum_{\mathfrak c}(N\mathfrak c)^{-\tau}
  \ \ll\ Y^\tau P_0^{\iota(1-\tau)+o(1)} .
\]
Here one can take \(\tau=\sigma\) and make \(\iota\) arbitrarily small. Alternatively choose \(1/2<\tau<\sigma\) so close to \(\sigma\) that \(B(1-\tau)<1\); the same Euler product includes both the polylogarithmically small primes and all large squarefull parts. For \(Y\ge P_0^\gamma\), choosing \(\iota<\gamma(\sigma-\tau)/(1-\tau)\) then gives \(O(Y^\sigma)\) with a power saving. Below \(Y\le P_0^\gamma\) the direct \(q\)-th moment \(\ll(1+Y)^q\) is sufficient after first choosing \(\gamma\) small. All bounds allow a common initial scale and a sufficiently small initial smoothness power.

An initial box that meets an archimedean coordinate hyperplane is partitioned by dyadic absolute value at each place. The norm of a nonzero integral root is at least one, so there are \(O_K((\log(2P_0))^{[K:\mathbb Q]})\) pieces. On each, a single totally positive unit balances the radii, giving a product of annuli of volume \(Y\ll P_0\) on which \(|Nx|\asymp Y\). Handle the zero root separately. Hence it suffices to use uniformly bounded annular supports. Over repeated divisions fix a canonical balanced \(b_{\mathfrak d}\) for the \emph{total} ideal \(\mathfrak d\) divided out of the initial root, with \(L\mathfrak d=b_{\mathfrak d}L_{\mathfrak d}\). The current root is always \(x/b_{\mathfrak d}\), with volume parameter exactly \(Y/N\mathfrak d\). Its annular constants depend only on the original annulus and the finite class list, irrespective of the number of divisions. This prevents accumulation of box-enlargement constants in the same-order induction.

\paragraph{Starting moments and the polynomial of degree one.}
The pure even starting estimates hold on the whole hereditary squarefree family. Suppose \(q=2s\ge4\) smooth slots of scale \(U\) have known pure exponents \(\lambda_q\) and \(\lambda_{q-2}\). For a polynomial \(\Phi\) of polynomially bounded cleared height and degree \(j\), write
\[
 E_j(T,U;q)=\int\left|\sum_{z}e_\alpha(\Phi(z))\right|^2
                  \left|\sum_x e_\alpha(x^k)\right|^q\,d\alpha ,
 \qquad U\le T .
\]
Subsets and the finitely many ideal lattices are allowed. In an upper bound one may first take absolute values of the coefficients in the underlying even orthogonality count.

If \(j=0\), use \(T^2 U^{\lambda_q}\). If \(j\ge1\), the same estimate and \(TU^q\) are available, since each value of \(\Phi\) has at most \(j\) preimages. For \(j\ge2\) one also has
\[
              E_j(T,U;q)\ll P_0^{o(1)}
                       \big(T U^{\lambda_q}+U^q\big)
\]
with the initial loss in \(\lambda_q\) understood. In fact a zero discrepancy among the smooth roots admits at most \(j|Z|\) pairs of polynomial roots. Fix a nonzero discrepancy \(d\), clear the polynomial denominators by \(A\), and observe that the nonzero \(h=z-w\) divides \(Ad\) in \(\mathcal O\). It has \(P_0^{o(1)}\) possible values, by ideal divisors and the logarithmic unit count in the polynomial-height box. Given \(h\), the equation \(\Phi(w+h)-\Phi(w)=d\) has degree \(j-1\) in \(w\) and is nonconstant. This is precisely the step requiring \(j\ge2\).

For each \(1\le M\le U\), \(M^k\le T\), write \(U'=U/M\). A further bound for \(E_j\), \(j\ge1\), is the maximum of the following three terms up to \(P_0^{O_{K,k,q}(\eta)+o(1)}\) and the initial admissible losses:
\[
 T^2M^2U^{\lambda_{q-2}},\qquad
 M^qT(U')^{\lambda_q},\qquad
 M^q\,\frac{T}{M^k}
          \big(E_{j-1}(T,U';q)(U')^{\lambda_q}\big)^{1/2}.
\]
Separate the at most \(j-1\) zeros of \(\Phi'\), at cost \(O_j(U^{\lambda_q})\). At nonzero derivative, designate the finitely many primes in fixed representatives and denominators, those in a cleared leading coefficient, and those dividing the current cleared \(\Phi'(z)\), as forbidden for the particular \(z\). Their product has polynomial height. A smooth pair \((z,x)\) fails if \(I_L(x)\) does not contain a divisor of norm at least \(M\), coprime to these primes, obtainable by principal blocks. The portion on forbidden primes is a divisor of one fixed nonzero cleared algebraic integer for each \(z\), since \(I_L(x)\) is squarefree. On failure the remaining portion has norm \(O(MR^{h^+})\). Thus the number of failing pairs is \(O(T M P_0^{O(\eta)+o(1)})\), by counting ideals and bounded-height generators of each. On either side of the energy equation a failing slot has this many choices before assigning the other \(s-1\) variables. The squared \(\ell^2\) norm of the resulting convolution is bounded by the squared number of failing pairs times the pure moment \(U^{\lambda_{q-2}}\), using translation invariance and the triangle inequality. Either side may be selected by separating the one-side counting function before squaring it. This yields the first term.

For each all-good \emph{one-sided} tuple select a witness divisor for each of its \(s\) smooth roots, of norm in \([M,MR^{h^+}]\) and coprime to the forbidden primes of that side. Give primes independent red or blue colours. All the witnesses are blue and the forbidden primes red with probability \(P_0^{-o(1)}\), uniformly, since a polynomial-height integer has \(O(\log P_0/\log\log P_0)\) distinct prime ideal divisors. Fixed denominator primes may instead be removed from the pool completely. For each fixed total value and this side of the equation, the good multiplicity is bounded by \(P_0^{o(1)}\) times the expected multiplicity on the rectangular family of red-derivative \(z\) and smooth \(x\) possessing a sufficiently large blue divisor. Jensen before summing the squared multiplicities allows one common colouring in the upper bound. Opposite sides need never have mutually compatible prescribed colours. In each rectangular family choose one witness ideal for each \(x\) using the blue primes only; this choice is independent of \(z\).

Apply H\"older in the divisor decomposition: for each side sum over its $s$ witness labels, use H\"older of order $2s$ on the two sides, and sum the resulting single-label means with total prefactor at most $(MR^{h^+})^{2s-1}$. The number of possible principal ideals is \(O(MR^{h^+})\) and there is one balanced generator per ideal. At each prime in the extracted divisor all retained \(\Phi'\) are units, so the polynomial congruence has at most \(j\) simple residue roots, each lifting uniquely. The CRT factor \(j^{\omega(m)}=P_0^{o(1)}\) allows Cauchy at each total value to impose \(z-w\in(m)^k\). The zero difference gives the second term. There are \(O(T/M^k)\) nonzero shifts, and at fixed shift the divided difference polynomial has degree \(j-1\). Cauchy against the smooth autocorrelation bounds its solution count by the square root displayed in the third term. This also proves the reduction at \(j=1\): the derivative is a nonzero constant, a simple linear fibre has at most one residue class, and the divided difference is a constant polynomial to which the \(E_0\) estimate applies. The preceding nonzero-discrepancy divisor bound, which would incorrectly discard a free \(w\) at \(j=1\), is never invoked there.

The other starting step uses the first difference of the pure monomial. In half-order \(s=q/2\), for \(M=Y^\theta\) with \(0<\theta<1/k\), its three exponents are
\[
 1+\lambda_{q-2}+\theta,\qquad
 1+\lambda_{q-2}+(q-2-\lambda_{q-2})\theta,\qquad
 1+\frac{\lambda_{2q-4}}{2}
       +\left(q-3-\frac{k-1}{2}-\frac{\lambda_{2q-4}}2\right)\theta .
\]
For each pair of opposite roots $(x_i,y_i)$ select a principal block
$(m_i)\mid(y_i)$ prime to $x_i$. If there is no such block at norm
$\asymp M$ within the permitted overshoot, the gcd has norm
$\gg Y/(M P_0^{O(\eta)})$. Summing the squared ideal-box count
$(1+Y/N\mathfrak d)^2$ over eligible gcd ideals bounds the
number of failing pairs by $YM P_0^{O(\eta)+o(1)}$. At fixed
$(x_i,y_i)$ the translated autocorrelation in the other slots is at
most $U_{q-2}(Y)$ by Cauchy.

For the remaining configurations, overcount once for every permissible
label tuple and group the label norms in dyadic bands. For a label $m$
write $f_m$ for the parent sum on roots prime to $m$ and $g_m$ for
the divided child. Define
\[
 I_m=\int |f_m(\alpha)|^2|g_m(m^k\alpha)|^{2s-2}\,d\alpha .
\]
At a fixed tuple of labels, split each paired product as
$|f_m|^{(s-2)/(s-1)}(|f_m|^2|g_m|^{2s-2})^{1/(2s-2)}$, suppressing
the modulation in this identity. H\"older bounds its integral by
$U_q(Y)^{(s-2)/(2s-2)}\prod_{i=1}^s I_{m_i}^{1/(2s-2)}$.
Here an even restricted parent moment is dominated by $U_q(Y)$ in
the nonnegative orthogonality count. H\"older in the summation over
the $O(MP_0^{O(\eta)})$ labels per band then gives, after absorption,
the bound $M^{q-3}\sum_m I_m$ in addition to the failing pairs, with
the permitted losses. If $s=2$ the fourth moment is the divisor
base below. In $I_m$, for each fixed total, Cauchy within the
unit $k$th-power fibres modulo $(m)^k$ confines $x-x'$ to $(m)^k$.
At good primes each fibre has at most $k$ residue directions and
unique simple lifts. At a fixed prime over $k$ the unit $k$th-power
kernel is bounded independently of the precision, by the logarithm
on sufficiently deep principal units. All local factors together
cost $P_0^{o(1)}$. On the positive solution count one can enlarge
the divided root supports to a common allowed child before taking
the absolute value of the resulting correlation. The regular diagonal costs \(Y M^{q-2}(Y/M)^{\lambda_{q-2}}\). For the off-diagonal the first correlation has nonzero shift \(h\), \(H=Y/M^k\), and centred variable \(z\). Up to fixed constants its frequency factors as
\(h\prod_{\zeta^k=1,\ \zeta\ne1}((1-\zeta)z+(1+\zeta)h m^k)\).
At nonzero frequency fix \(h\) and the values of two nonproportional factors by bounded-height ideal divisors in \(K(\mu_k)\). They determine \(z,h m^k\) and leave at most \(k\) roots \(m\). The zero frequency (a root-of-unity ratio, with nonzero $h,m$) has \(O(HM)\) representations. Hence the correlation has squared norm \(O(YHM P_0^{o(1)})\). The one-block Cauchy bound is consequently \(M^{q-3}(YHM)^{1/2}(Y/M)^{\lambda_{2q-4}/2}\), giving the third exponent. The ordinary fourth-moment base \(\lambda_4=2\) follows from the same nonzero-difference divisor argument for \(x^k-y^k\), with the zero fibre of size \(O(Y)\).

\paragraph{Choosing divisors and bounding the Gram matrix.}
At a squarefree state with volume \(Y>P_0^\gamma\), let a scheduled block use \(0<L_0<M<k\). Each permitted relative ideal contains a narrow principal divisor \((m)\) of norm in \([Y^{L_0/k},Y^{M/k}]\). Indeed its prime classes can be partitioned into principal blocks of at most \(h^+\) primes and a residual sequence of length \(<h^+\), by repeated removal of a minimal nonempty zero-sum subsequence. The class group has order \(h^+\). The available product has norm \(\gg_K Y/R^{h^+}\); taking whole blocks overshoots the lower endpoint by at most \(R^{h^+}\). These inequalities suffice as soon as \(\eta/\gamma\) is small relative to \(M-L_0\) and \(k-M\), simultaneously over the finite menu. For each root make one deterministic choice; then exactly
\[
             F(\alpha)=\sum_{(m)} F_m(m^k\alpha).
\]
One generator is used per principal ideal. At the parent indexed by \(\mathfrak b\) take \(m=b_{\mathfrak b(m)}/b_{\mathfrak b}\), which is an integral generator of \((m)\) since the fixed class representatives are then equal. The functions \(F_m\) have the same bounded-coefficient property at the new scale \(Y/|Nm|\). The harmonic sum over all possible ideals is \(O_K(\log(2P_0))\), sufficient for weighted H\"older.

Let \(H_m\) be the restriction of \(F\) to roots coprime to \((m)\). The two conductor families of Lemma~\ref{lem:nf-fibre-savings} and the Gram argument leading to \eqref{eq:nf-short} apply with the same conductor cover and a strict margin \(\sigma_0>2z\). Here \(\log_Y N\mathfrak p\le\eta/\gamma\) for every factor of \(m\), and \(\operatorname{rad}(m)=(m)\). For \(P_0\) sufficiently large the finitely many primes over the different, over rational primes \(\le k\), and in the fixed class representatives all have norm below \(W\). In particular coprimality to a factor of \(m\) has the same meaning for \(x\) and for \(I_{L_{\mathfrak b}}(x)\).

Two uniformity details replace the marked-prime hypotheses. First, all \(\mathfrak p\mid(m)\) have norm \(>W\), whence \(\omega(m)=O(\log P_0/\log\log P_0)\). Products of \(O_k(1)\) per local factor, including the bounds for unit \(k\)-th-power fibres, all inclusion--exclusion terms for translated masks and the numbers of divisor ideals of \((m)^k\), are \(P_0^{o(1)}\). The conductor has exponents at most \(k\). In particular the low-conductor contribution to the Gram norm is at most \(Y P_0^{o(1)}\), and the number of large points on a fibre is \(\ll P_0^{o(1)}V^2\) at threshold \(Y/V\), \(V\ll Y^z\). The high-conductor off-diagonal term is absorbed by the strict saving with \(Y\ge P_0^\gamma\) and sufficiently small \(\eta/\gamma\). The derivative and multiplicity bounds use genuine finite residue fields at prime ideals; residue degree one is not required.

Second, a subideal \(\mathfrak u\) chosen in a split or a shift of Lemma~\ref{lem:nf-fibre-savings} need not be principal. Choose once for each class a fixed integral inverse-class ideal \(C\), and use a balanced generator of \(\mathfrak u C\). The extra norm factor \(NC\) is bounded in terms of \(K\). For \(P_0\) large its prime factors are below \(W\), hence coprime to every active conductor. Splits into all residue classes of \(\mathfrak u C\) remain legitimate, and the new generator has the prescribed valuations at every active prime. A shift by this generator cancels exactly the planned active conductor divisor, apart from the existing additional divisor of the summation shift. The fraction of nonzero shifts divisible by \(\mathfrak d\) is \(O(1/N\mathfrak d)\) when \(N\mathfrak d\) is within the shift volume, up to a fixed box constant. Thus the divisor average and the normalized conductor bounds incur only \(P_0^{o(1)}\) and fixed factors. In the unsplit degree-\(k\) variant one imposes a mask divisor \(\mathfrak v\) by partitioning its lattice into cosets modulo \(\mathfrak v C\); only \(NC\) cosets occur, and the new leading coefficient has the same active conductor as division by \(\mathfrak v\). This covers nonprincipal prime ideals in every form of the finite menu.

\paragraph{Good points and a contraction on the bad points.}
On a dyadic height interval where \(|F|\asymp Y^{1-u}\), use a direct admissible moment at order \(h\), or a fibre bound from prior order \(r\ge k\) and parameters \((z,L_0,M)\), \(u\le z\). On the set where \(|H_m|\ge |F|/2\) for every possible parent divisor, the Gram estimate gives
\(\sum_{\alpha\ {\rm in\ a\ fibre}} |F(\alpha)|^2
 \ll Y^2P_0^{o(1)}\)
on the indicated level. Indeed dyadic slicing of the level only costs a logarithm. Weighted H\"older in the exact divisor decomposition, the preimage weight \(|Nm|^{-k}\), and the child bound at order \(r\) give the same two alternative excesses at order \(q\) as for the marked sums:
\[
       D_h+(h-q)u,\qquad
        (1-L_0/k)D_r+(r+2-q)u .
\]
Thus any verified finite cover of the height interval yields its prescribed good-point contribution. All logarithmic losses here are absorbed once on each edge of the finite moment proof.

For a bad point some candidate \(m\) has \(|F-H_m|>|F|/2\). For every squarefree ideal \(\mathfrak d\), write \(F_{\mathfrak d}^{*}\) for the restriction to \(\mathfrak d\mid I_{L_{\mathfrak b}}(x)\), with the initial parent coefficients. Inclusion--exclusion and weighted H\"older give the pointwise bound
\[
 |F|^q \le 2^q
 \left(\sum_{\substack{\mathfrak d\mid(m)\\N\mathfrak d>1}}
              (N\mathfrak d)^{-a}\right)^{q-1}
 \sum_{\substack{\mathfrak d\mid(m)\\N\mathfrak d>1}}
              (N\mathfrak d)^{a(q-1)}|F_{\mathfrak d}^{*}|^q .
\]
The first sum is \(o(1)\) uniformly in the chosen \(m\), being at most
\(\exp(\omega(m)W^{-a})-1=o(1)\). Enlarge the second sum to all ideals \(\mathfrak d\) with \(N\mathfrak d>W\) that divide some current root. The resulting majorant is valid on the \emph{union} of all bad dyadic levels and menus, so it is integrated just once. On division by \(\mathfrak d\), the fractional-scalar Haar argument gives the child norm at volume \(Y/N\mathfrak d\), with its canonical new lattice, for arbitrary \(q\). The coefficient multiplying the supremum of normalized same-order child moments is bounded by
\[
               O(1)\sum_{N\mathfrak d>W}
                        (N\mathfrak d)^{a(q-1)-\lambda}=o(1)
\]
by \(a(q-1)<\lambda-1\). The normalized supremum at the fixed initial cutoff is finite by the trivial norm bound. Absorbing the last display establishes the desired rough squarefree moment uniformly at every state; Minkowski above extends it to the literal full upper-smooth moment. This absorption has no per-descent factor \(P_0^\epsilon\). Earlier orders and good-point terms contribute only their specified initial loss, which is fixed small enough at the start of the finite proof. The full-family fibre rule requires the strict test \(\lambda>q^2/(2q-1)\) at its fibre-output nodes; interpolation nodes require only prior uniform moment bounds.

\paragraph{Exact checks at degrees eight, nine and ten.}
Use \(a=11/20,\ B=37/20\), so \(aB>1\) and \(B<2\). Each actual shortening node is tested at its \emph{assigned output} exponent for \(\lambda-a(q-1)>1\). At every final order \(\sigma>1/2\) and \(B(1-\sigma)<1\); the displayed Rankin bound therefore applies. The following exact certificates use the two conductor families of Lemma~\ref{lem:nf-fibre-savings}, with a positive rational margin in the large-point condition.
\begin{center}
\begin{tabular}{@{}crrrrl@{}}\toprule
 \(k\) & \(q\) & starting steps & moment nodes & height intervals & verified excess \(D_q\)\\\midrule
 8 & 34 & 125 & 37 & 674 & \(1106681/10000000\)\\
 9 & 40 & 140 & 73 & 1119 & \(227719/2500000\)\\
10 & 46 & 155 & 97 & 1460 & \(197067/2500000\)\\\bottomrule
\end{tabular}
\end{center}
The files \path{code/full/base_k8.json}, \path{code/full/base_k9.json}, \path{code/full/base_k10.json} specify the exact starting steps. The companion \path{code/full/full_k8.json}, \path{code/full/full_k9.json}, \path{code/full/full_k10.json} specify the height cover at each later node. The programs \path{code/full/verify_base.py} and \path{code/full/verify_curve.py} check each implication in rational arithmetic. The three excesses are respectively below \(5857852793/50000000000,\ 8516302869/62500000000,\ 152339/1000000\). These are the full upper-smooth inputs of Theorem~\ref{thm:nf-small}. The data and rational verifiers are in Appendix~\ref{app:certificates}. At every node the file records the parent orders, and either an interpolation weight, a differencing rule with its parameters, or a rational cover of the height intervals by the two displayed alternatives. For the baseline polynomial estimates the verifiers recursively check the three-term bound above down to the degree-zero case, together with the first-difference and trivial alternatives.

\paragraph{The loss-free pure moment on a box meeting a coordinate hyperplane.}
The positivity argument for the sharp even order uses an ordinary smooth envelope on a product of annuli. The same sharp bound holds on a whole box. For the nonzero roots, dyadically partition the absolute values at the archimedean places; let \(X\) be the original root volume and \(Y\asymp X2^{-\sum_v d_vj_v}\) the volume of one cell. A single totally positive unit balances each such cell. Its relative ideal support is preserved, and the induced multiplication on the trace torus preserves Haar probability. For \(Y\ge X^\gamma\) the old smoothness cutoff \(X^\eta\) is at most \(Y^{\eta/\gamma}\). Choose \(\eta\) after \(\gamma\) so that the previously obtained full-smooth and critical estimates apply uniformly. At even \(S\), positive coefficient domination permits a common smooth annular envelope after balancing. Minkowski's inequality sums the resulting \(L^S\) norms: the sum of \(Y^{1-k/S}\) over all cells is \(O(X^{1-k/S})\), since \(1-k/S>0\) and the series in each \(j_v\) is geometric. For \(Y<X^\gamma\) the nonzero integral norm bound makes each nonempty cell have \(Y\gg_K1\), hence \(O_K(Y)\) roots by packing. There are only \(O_K((\log(2X))^n)\) nonempty cells. Take \(0<\gamma<1-k/S\), and add the single possible zero root separately. These bounds give \(\int |g|^S\ll X^{S-k}\) without a logarithmic loss for the original box.

\subsection{The quartic tenth moment}
\label{subsec:aud-quartic-moment}

Write $P$ for volume, $U_j(P)$ for the $j$th ordinary trace-torus moment of the smooth roots under study, and $M=P^{9/50}$, $H=P/M^4$, $Q=P/M$. First work in the hereditary family of smooth roots in field annuli with $N((x)/\operatorname{rad}(x))\le P_0^{C\eta}$, where $P_0$ is the initial volume. This family contains the positive distinct-mark construction used on the major arcs. Here eligibility means only the upper prime cutoff. The fourth moment is $O(P^{2+\varepsilon})$ by the divisor argument of Section~\ref{subsec:aud-full-moments}. The next two exponents are $13/4$ and $(7+\sqrt5)/2$. Here are the field ingredients for that calculation, following the difference scheme of \cite[Theorem~4.3]{nf3:vaughan89}. For each paired equality choose a balanced principal block $(m)\mid(y)$ coprime to $x$, $N(m)\in[M,MP_0^{O(\eta)}]$. Such a block is assembled from unused single primes in groups of at most the class number. Failure implies a common divisor of norm $\gg P/(MP_0^{O(\eta)})$. Summing ideal-box packing for the two roots costs $O(PM P_0^{O(\eta)})$, with remaining multiplicity bounded by the previous moment. Simultaneous H\"older and congruencing of units give the diagonal $P M^{2t-2}Q^{\Lambda_{t-1}}$ at order $2t$, and an offset integral times $M^{2t-3}$. Here $\Lambda_t=\lambda_{2t}$ denotes the exponent at order $2t$. Symmetrizing $x-x'=m^4h$ gives the phase $hz(z^2+h^2m^8)$ up to a fixed multiplier, with $h,z,m$ of volumes $H,P,M$. A second difference has frequency proportional to $h l(3Z^2+l^2+4h^2m^8)$. Its energy is $O(P^{2+\varepsilon}HM)$. At nonzero frequency choose $h,l$ by divisors; the residual equation $3Z^2+(2hm^4)^2=b\ne0$ has $O(P^\varepsilon)$ choices by factoring over $K(\sqrt{-3})$. At $b=0$ the frequency is $hl^3$, giving total energy $O(HPM^2P^\varepsilon)$. At frequency zero the equation $3Z^2+l^2=-4h^2m^8$ gives total energy $O((HM)^2P^\varepsilon)$. This covers the isotropic signatures as well. H\"older with weights $1/4,(4-t)/4,(t-1)/4$ on this energy, $U_{2t-2}(Q),U_{2t}(Q)$ yields at $t=3,4$ the additional expression
\[
 M^{2t-3}P^{1/2}(HM)^{3/4}
 Q^{((4-t)\Lambda_{t-1}+(t-1)\Lambda_t)/4}.
\]
Using respectively $M=P^{1/8}$ and $M=P^{(2\sqrt5-3)/11}$ in these preliminary steps gives the two stated exponents by contraction, including the exceptional $PM U_{2t-2}(P)$. In particular $D_8\le0.619$ is admissible.

For $t=5$ return to $M=P^{9/50}$. The same congruencing lemma yields
\[
 U_{10}(P)\ll P^\varepsilon\bigl(PM U_8(P)+M^7(PM U_8(Q)+T)\bigr),\qquad
 T=\int |F_1(\alpha)|\,|G_Q(c\alpha)|^8\,d\alpha ,
\]
with the once-differenced $F_1$ above and a fixed nonzero integer $c$. Enlarge the actual nonnegative solution count to independent smooth boxes for $h,z,m$ before taking the absolute value of $F_1$. The proof adapts the rational harmonic estimate of \cite[Section~2]{audBDWquartic}. Cauchy in $(h,m)$, the symmetric difference $l=z_1-z_2$, then Cauchy in $(h,l)$ give
\[
 |F_1|^2\ll P(HM)^2+HM(DE)^{1/2}.
\]
Here $D=\sum_{h,l}|\sum_Z w_l(Z)e_K(c_1\alpha h l Z^2)|^2$ and $E=\sum_{h,l}|\sum_m w(m)e_K(c_2\alpha l h^3m^8)|^2$. The nonzero $h,l$ have volumes $H,P$ and normalized derivatives of $w_l$ are uniform. One has $E\ll PHM+|E_0|$, where $E_0$ contains precisely the nonzero frequencies $c_2 l h^3(m_1^8-m_2^8)$ with nonnegative coefficients.

Put $B(t)=\prod_v\max(1,|t_v|)^{d_v}$, and on the disjoint arcs $Q_\zeta B(Q^{4/n}(\beta-\zeta))\le C_0P$ let $\Omega(\beta)$ equal $1/[Q_\zeta B(Q^{4/n}(\beta-\zeta))]$ and define it to be zero elsewhere. We need the following three field bounds, with fixed modulations restored in their proofs.
\[
 D\ll P^{2+\varepsilon}H+P^{3+\varepsilon}H\Omega,\qquad
 W=\int\Omega |G_Q|^4\ll Q^\varepsilon,\qquad
 V=\int|E_0|^2|G_Q|^4\ll Q^{5+\varepsilon}.
\]
For the first, difference the quadratic in $Z$ and Poisson sum its linear phase. The nonzero product $hlv$ has volume bound $T_0=HP^2$ and $O(P^\varepsilon)$ multiplicity. The lattice $(t,c_1\alpha t-\xi)$ is restricted to sides $T_0^{1/n}, P^{(-1+\varepsilon)/n}$ after rapid truncation. If the lattice points in the box span rank two over $K$, multiply
two $K$-independent vectors by each element of a fixed integral
basis. This gives $2n$ real-linearly independent vectors in a
constant enlargement of the normalized box and hence a bounded
covering radius there. The covolume after scaling is
$\gg P^{1-\varepsilon}/T_0$ (adjust $\varepsilon$ to absorb the
finite truncation powers); counting disjoint fundamental cells in
the enlarged box gives $O(P^\varepsilon T_0/P)$ points. In rank one, the common slope $\zeta$ gives $O(1+P^\varepsilon T_0/[Q_\zeta B((PT_0)^{1/n}(c_1\alpha-\zeta))])$ points. The outer Poisson factor is $P$, and $PT_0=Q^4$. The arcs are disjoint since $P^2=o(Q^4)$. With the fixed multipliers restored $W=\int\Omega(c_1\alpha)|G_Q(c\alpha)|^4$. The Fourier coefficient of $\Omega(c_1\alpha)$ vanishes off $c_1\mathcal O_K$ and equals $\widehat\Omega(y/c_1)$ on that ideal, by Haar invariance and the finite kernel. For $W$, the absolute archimedean integral per center of conductor $\mathfrak a$ is $O(Q^{-4}(N\mathfrak a)^{-1}\log^n P)$. Summing its translates produces the ideal Ramanujan sum $c_{\mathfrak a}(y)$. At $y\ne0$, $\sum_{N\mathfrak a\ll P}|c_{\mathfrak a}(y)|/N\mathfrak a\ll P^\varepsilon$ by ideal divisors; at $y=0$ the cost is $O(P)$. The zero Fourier coefficient of $|G_Q|^4$ is $O(Q^{2+\varepsilon})$ and its total coefficient mass is $O(Q^4)$. This gives $W\ll Q^{-4+\varepsilon}(P Q^2+Q^4)$. For $V$, a fixed nonzero $lh^3(m_1^8-m_2^8)$ has $O(P^\varepsilon)$ representations by factoring $h,l,m_1-m_2$ and solving a polynomial of bounded degree. Thus $\int|E_0|^2\ll P^{1+\varepsilon}HM^2$. Two elementary differences of the ordinary fourth-power sum give $|f_Q|^4\le C Q^3+C Q\sum_{y\ne0} b(y)e_K(\alpha y)$ with $|b(y)|\ll Q^\varepsilon$: at prescribed frequency $hl(3Z^2+h^2+l^2)$, extract $h,l$ as divisors, leaving at most two $Z$. At frequency zero, nonzero shifts still leave at most two $Z$ per pair. Whole-torus domination by $f_Q$ in $V$ is legitimate since $|E_0|^2$ has nonnegative Fourier coefficients. Hence $V\ll P^\varepsilon(Q^3PHM^2+Q(PHM^2)^2)$ and $PHM^2=Q^2$.

For explicit verification the remaining H\"older expressions, with all $U_j=U_j(Q)$, are
\begin{align*}
 T&\ll P^\varepsilon\{P^{1/2}HM U_8+(HM)^{1/2}[P^{1/2}H^{1/4}T_2
        +P^{3/4}H^{1/4}(T_4+(PHM)^{1/4}T_5)]\},\\
 T_2&\ll (PHM)^{1/4}U_8+U_8^{5/8}U_{10}^{1/4}V^{1/8},\\
 T_4&\ll V^{1/8}W^{1/4}U_{10}^{1/2}U_{12}^{1/8},\qquad
 T_5\ll W^{1/4}U_8^{1/4}U_{10}^{1/2}.
\end{align*}
Starting at $D_{10}=2/3$ by $U_{10}\ll Q^2 U_8$, use the elementary $U_{12}\ll Q^2U_{10}$ throughout. Substitution of these bounds, also implemented in \path{code/cubes/quartic_bw.py}, gives the successive tenth-moment excesses
\[
 \frac{533}{1200},\quad\frac{31693}{96000},\quad\frac{2086613}{7680000},\quad
 \frac{148527133}{614400000},\quad\frac{11379}{50000}<\frac14 .
\]
All calculations have strict slack after choosing the initial smoothness power sufficiently small; only finitely many child scales occur. The restriction on repeated prime factors can also be removed for the full upper-smooth tenth moment. Write $(x)=\mathfrak b\mathfrak a$, where $\mathfrak a$ contains the prime factors of multiplicity exactly one and $\mathfrak b$ is squarefull. Choose fixed integral ideal representatives $\Lambda$ and balanced $a$ with $\mathfrak b=a\Lambda$. For $x=ay$, $(y)=\Lambda\mathfrak a$ has radical quotient of norm at most $N\Lambda$. The common $a$ cancels in the energy equation. Partition archimedean modulus ranges dyadically and balance them by totally positive units; there are only logarithmically many pieces, each of volume $T\ll P$. At child volume $T/N\mathfrak b\ge P^{c_0}$ the same hereditary estimate applies on taking $\eta$ small relative to $c_0$. Minkowski in $L^{10}$ costs $\sum_{\mathfrak b\text{ squarefull}}(N\mathfrak b)^{-\lambda/10}<\infty$ for $\lambda=6+11379/50000>5$. The remaining squarefull tail has total trivial mass $O(P^{(1+c_0)/2+\varepsilon})$, since squarefull ideals through norm $B$ number $O(B^{1/2+\varepsilon})$. For example $c_0=1/10$ makes this smaller than $P^{\lambda/10}$. Thus $D_{10}\le11379/50000<393/1600$ for the full upper-smooth quartic sum, as required for the loss-free twelfth moment and the density theorem at six roots.

\subsection{Successive congruences in degrees five through seven}
\label{subsec:aud-block-moments}

The following congruence step applies at every fixed degree used below.
It implements the recursive integral of \cite[Lemma~2.2]{edvw}
on the hereditary family just defined. At an intermediate even moment
the regularity masks, the shortened sums, and the phase parameters are
held fixed until after the corresponding positive orthogonality count.
Fix the first-difference polynomial $\Psi$ at its external parameters. In $E=\int|f_\Psi|^2|G_S|^r$, $r=2t\ge4$, a tuple on either side independently has one polynomial root $z$ of volume $P$ and $t$ smooth roots of volume $S\le P$. For each smooth root $x$ on that side select a balanced principal block $(m)\mid(x)$, of norm $M$ up to the allowed overshoot, coprime to $\Psi'(z)$ and the fixed integer modulations. When extraction fails, the gcd with the derivative has norm $\gg S/(M P_0^{O(\eta)})$. Each nonzero derivative value has only $P^\varepsilon$ divisor ideals, and ideal-box packing gives $O(PM P_0^{O(\eta)+\varepsilon})$ possible offending $(z,x)$ in total. The finitely many derivative zeros contribute $O(S)$. At fixed $(z,x)$ Cauchy bounds the remaining count by $E^{1/2}U_{r-2}(S)^{1/2}$. Thus the exceptional contribution after absorption is $P^2M^2 U_{r-2}(S)$ up to permitted losses.

For a good one-sided tuple each witness is prime to the derivative
at \emph{its own} polynomial root. The one-sided colouring and Jensen
argument of Section~\ref{subsec:aud-full-moments} gives rectangular
root and label sets before squaring the total multiplicities. Equivalently,
overcount the good tuples by summing the root-divisor labels independently
on both sides. For a fixed tuple of labels on each side, H\"older of order $2t$ gives for each first-side $m_i$ the factor
\[
 \left(\int |f_{\Psi,\boldsymbol m}|^2 |G_{m_i}(m_i^k\alpha)|^{2t}\right)^{1/(2t)},
\]
and similarly on the opposite side. The $f$-mask imposes regularity at all labels of that side. These are whole-torus even integrals with nonnegative coefficients, so each mask may now be relaxed to regularity only at its own $m_i$. Summation costs $M^{2t-1}\sum_m I_m$. At the single $m$ simple-root lifting and Cauchy within root directions force $z-z'=m^k h$, with $k^{\omega(m)}\ll P^\varepsilon$ directions. The zero shift gives the diagonal; symmetrizing a nonzero shift gives the new difference polynomial.

For genuine blocks $m_1,\ldots,m_d$, $d=1,2$, put $M_i=P^{\theta_i}$, $H_i=P/M_i^7$, $Q_j=P/(M_1\cdots M_j)$, $B_d=\prod_{i\le d} H_i M_i$, with $\theta_2\le\theta_1$ when $d=2$. Additional ordinary shifts have volume $P$. After $j$ shifts the centered polynomial, up to a fixed multiplier, is $(\prod h_i)\mathbb E[(Z+\sum_i s_i V_i)^{7-j}]$, with $s_i=h_i m_i^7$ at active blocks, $s_i=h_i$ at ordinary shifts, and $V_i$ uniform on $[-1,1]$. This identity is over the rationals. Divide the original annulus into finitely many sufficiently small root boxes; the $s_i$ are then small compared with $Z$ at a chosen archimedean place after all fixed enlargements, which excludes zero frequency. At a fixed nonzero frequency all $h_i$ are divisors. With one genuine block and odd residual degree at least three, $Z$ is an additional divisor and $m_1$ solves a nonconstant polynomial. At residual degree two one uses a binary quadratic norm. At degree four, the leading binary quadratic in $Z^2$ and $(h_1m_1^7)^2$ is $V^2+2VU+U^2/5$; completing its square gives a nonsingular binary norm. A zero right side fixes the frequency in terms of the shifts, and the squared multiplicity costs at most $O(M_1^2P^\varepsilon)$ there, which fits since $M_1\le P$. The first difference itself splits over $K(\mu_7)$ as in Section~\ref{subsec:aud-full-moments}.

For two genuine blocks the first residual polynomial has degree five. Fix also $Z$ as a divisor. In $U=(h_1m_1^7)^2$, $V=(h_2m_2^7)^2$, let $S'=U+V$, $D'=U-V$. The remaining equation is
\[
 3Z^4+10Z^2 S'+4(S')^2-(D')^2=N,\qquad
 (4S'+5Z^2-2D')(4S'+5Z^2+2D')=13Z^4+4N .
\]
Its nonzero norm gives $O(P^\varepsilon)$ possibilities; the degenerate frequencies are determined by $h_1,h_2,Z$, with squared multiplicity bounded by $\min(M_1,M_2)^2P^\varepsilon\le M_1M_2P^\varepsilon$. This proves $\int|F_2|^2\ll P^{1+\varepsilon}B_2$. For further ordinary differences fix $m_2$ and apply the one-block proof; a single additional factor $M_2$ suffices in the integrated energy. Ordinary Hua differencing, followed by Cauchy against the preceding power, yields
\[
 \int|F_d|^{2^l}\ll P^{2^l-l+\varepsilon}B_d^{2^l-1}M_2^{\sigma_l},
\]
where for $d=1$ the last factor is omitted and $l\le5$, and for $d=2$, $l\le4$, $\sigma_1=0$, $\sigma_l=1$ for $l\ge2$. These conclusions include the zero-shift diagonal terms in Hua's induction.

Here and below the moment induction can be seen directly. Write $B$
for the number (in volume) of external labels and $T_v$ for the sum
with $v$ further nonzero ordinary differences in the polynomial root.
For $v\ge1$, Cauchy in the labels and differencing in that root give
\[
 |F|^{2^v}\ll P^{2^v-1}B^{2^v}
        +P^{2^v-v-1}B^{2^v-1}|T_v|.
\]
Indeed the diagonal term records a zero difference at any stage.
For the induction step, square the preceding inequality and apply
Cauchy over the $O(BP^v)$ labels of $T_v$ before differencing its
root; the resulting zero shift is covered by the next diagonal
term. On multiplying the displayed inequality by $|F|^{2^v}$ and
integrating, the off-diagonal is at most the indicated prefactor
times $(\int|F|^{2^{v+1}}\int|T_v|^2)^{1/2}$. Absorption and the
preceding power bound the full moment at $2^{v+1}$. In particular,
an additional multiplicity in the differentiated squared energy is
paid once, as recorded in $\sigma_l$.

For the numerical deduction write $q=2s$, $u=q-2$, choose a prior even order $t$, and put $u_2=2u-t\ge4$. After the first extraction the variable term is $M_1^{q-3}\int|F_1||G_1|^u$, alongside $PM_1U_u(P)$ and $PM_1^{q-2}U_u(Q_1)$. Use Cauchy with $U_t(Q_1)$ and the second congruence bound
\begin{align*}
 \int|F_1|^2|G_1|^{u_2}\ll{}& B_1^2P^2M_2^2U_{u_2-2}(Q_1)\\
 &+B_1M_2^{u_2-1}\left(PB_1M_2U_{u_2}(Q_2)+\int|F_2||G_2|^{u_2}\right).
\end{align*}
H\"older at the established $2^l$th power bounds treats the last term. The program \path{code/ed/verify_septic_seeds.py} checks 900 inequalities in 30 finite waves using \path{code/ed/seeds_k7.json}. For example the full exponents at orders 10 and 14 are at most $5.5961979$ and $8.5756101$. Using \path{code/ed/verify_level_curve.py}, the 25-node marked shortening certificate \path{code/ed/level_k7.json} then checks 481 rational height intervals. It gives
\[
 D_{57/2}\le 99077/10^6,\quad D_{115/4}\le113687/1250000,\quad
 D_{1456/51}\le4971653/51000000<5/51 .
\]
Extending this last bound to order $30$ by the trivial supremum gives the loss-free order $34$ in the critical absorption. The mixed estimate in the proof of Theorem~\ref{thm:nf-small} uses the sharper fractional order itself.

Extend the preceding formula to $J\le k-2$ genuine blocks, with
$H_i=P/M_i^k$, $Q_i=P/(M_1\cdots M_i)$,
$B_i=\prod_{a\le i}H_aM_a$. Require $M_1\ge\cdots\ge M_J$ and
$\prod_{i\le\nu} M_i\le H_\nu$ for $2\le\nu\le J$.
Bin the actual block norms, including the $P_0^{O(\eta)}$ overshoot,
before bounding $h_i$. For $W_i=\int|F_i||G_i|^{u_i}$, choose a prior
even $t_i$ with $u_{i+1}=2u_i-t_i\ge2$ even. The same extraction proof gives
\begin{align*}
W_i^2\ll U_{t_i}(Q_i)\big[&
B_i^2P^2 M_{i+1}^2 U_{u_{i+1}-2}(Q_i)\\
&+B_iM_{i+1}^{u_{i+1}-1}
\{P B_i M_{i+1}U_{u_{i+1}}(Q_{i+1})+W_{i+1}\}\big]
\end{align*}
up to the permitted arbitrarily small parameter losses.
The pair-count Cauchy inequality also applies when $u_{i+1}=2$,
where $U_0=1$. At a terminal state there is a second bound besides
the Hua moments below. For $\deg\Psi\ge2$ at fixed parameters
and nonzero child discrepancy $n$, $h=z-w$ divides a fixed multiple
of $n$ and $\Psi(w+h)-\Psi(w)=n$ has boundedly many solutions per
choice of $h$. At $n=0$ there are $O_k(P)$ pairs. Therefore
$\int|F_J|^2|G_J|^a\ll B_J^2(P U_a(Q_J)+Q_J^a)P^\epsilon$.
An ordinary Cauchy split applies this estimate to $W_J$.

Here the centered-coefficient energy is
$O(P^{1+j-J+\epsilon}B_J)$ after $j-J$ additional ordinary shifts,
for $j\le k-2$. Write $d=k-j$, $S=\sum s_i^2$, $T=\sum s_i^4$ in
the centered formula above. At fixed nonzero frequency all shift
parameters divide its fixed integer multiple. Cauchy over the
$O(P^\epsilon)$ tuples lets us take the same $h_i$ on both sides.
For odd $d$, the center $Z$ also divides that multiple and can be
taken equal. Nonvanishing follows in a single chosen embedding
on sufficiently narrow root boxes. In degree two the comparison factors
$Z^2-W^2=(S'-S)/3$. In degree three the residual divided by $Z$ is
$Z^2+S$. In degree four the residual is
$(Z^2+S)^2-C$ with $C=2S^2/3+2T/15$. If $C\ne C'$,
the difference of its two squares determines both centers
up to $O(P^\epsilon)$ choices per pair of label tuples.
For $C=C'$ the equations
$Z^2-W^2=S'-S$ and $Z^2+W^2=-(S+S')$ likewise factor,
over $K(i)$ if necessary, when the right side is nonzero.
If $d=5$, the residual divided by $Z$ is
$Z^4+(10/3)SZ^2+(5/3)S^2-(2/3)T$, so $S\ne S'$ determines
$Z^2$. For each such off-diagonal case $\prod M_i\le P$ pays for
the second tuple of labeled $m_i$.
The first-difference splitting argument over $K(\mu_k)$ also
handles the initial single-block case.

The remaining signed diagonals satisfy $S=\pm S'$, with at most
$O(P)$ paired centers. For the minus sign, if an ordinary shift is
present, fixing the other shifts and all labels determines that
ordinary shift up to two choices. The saved volume $P$ pays for all
extra labels. Otherwise let $I$ be the set of active genuine indices
for which $m_i^{2k}\mp(m_i')^{2k}\ne0$. Outside $I$ there are
$O_k(1)$ right labels per left label. $I=\varnothing$ gives the
required bound. $|I|=1$ is impossible since the shift is nonzero.
For $|I|\ge2$ the largest $\nu\in I$ has $\nu\ge2$. Its shift $h_\nu$
is determined by a nonzero quadratic after all the other shifts and
all labels are fixed. The extra labels have at most $O(\prod_{i\in I}M_i)$
choices, covered by the saved $H_\nu$ by the prefix inequality.
This covers all complex degeneracies as well as the real ones.
The usual ordinary Hua induction now gives
\[
 \int |F_J|^{2^l}\ll P^{2^l-l+\epsilon} B_J^{2^l-1},
 \qquad J+l-1\le k-2 .
\]
The proof uses the squared energy established at each ordinary
depth and Cauchy against the preceding power, including the zero
shifts. H\"older interpolation applies to the terminal smooth order.

The checker \path{code/ed/verify_multiblock_seeds.py} verifies 784 inequalities in 28 finite
rounds in each degree. For $k=6$ it gives, amongst others,
$\lambda_{16}\le10.5641816$, $\lambda_{18}\le12.3636536$ and
$\lambda_{20}\le14.237541$. For $k=5$ it gives
$\lambda_{14}\le9.3095507$.

Alongside the block estimates one may use \eqref{eq:eddefect} from
Theorem~\ref{thm:edmoment} at the current even exponent: the proof
of that recurrence starts with the input bound at order $q-2$
on every shortened lattice. Grouping the exceptional ideal as in
Section~\ref{subsec:aud-full-moments} returns exactly its hypotheses
with the original-scale prime cutoff. Initially $\lambda_0=0$,
$\lambda_2=1$ by bounded degree, and $\lambda_4=2$ by the
difference-of-powers divisor bound. Retaining the previous exponent,
applying this recurrence, and the one- or multiblock estimates are
the lower-wave rules in the tables; each wave uses only the preceding
wave. At order $r$ between two known orders interpolate their full
exponents linearly; above the largest order the trivial supremum
increases the exponent by the increase in order.

The interpolation and height rules used for the marked families are the
ones displayed in Section~\ref{subsec:aud-full-moments}. The admissible
block in a marked fibre uses its assigned windows and the two conductor
families of Lemma~\ref{lem:nf-fibre-savings}; the derivative masks above
are relaxed only in even counts. A further intermediate defect rule,
proved in Section~\ref{subsec:aud-high-moments}, applies before a defect
vanishes. At degree six, \path{code/ed/level_k6.json} uses seven nodes and 79
height intervals, giving $D_{22}\le66943/625000<1/9$. At degree five
the initial multiblock bounds include $\lambda_{10}\le5.9250801$ and
$\lambda_{12}\le7.5417549$. We refine its highest relevant even orders next.

\subsection{Two ordinary and fifteen marked fifth powers}
\label{subsec:aud-quintic-moment}
The two blocks use their moduli asymmetrically, following the rational
estimate \cite[Lemma~6.1]{edvw}. Here the necessary modulus means and
the final transfer hold over every fixed number field.

\paragraph{The second difference and its positive majorants.}
Use the same successive congruence extraction and bounded-coefficient hereditary smooth family as for degrees six and seven. Write $P$ for volume, $n=[K:\mathbb Q]$, $M_1=P^x$, $M_2=P^y$, $H_i=P/M_i^5$, $Q=P/(M_1M_2)$, $B_i=\prod_{j\le i}H_jM_j$, and $B(w)=\prod_v\max(1,|w_v|)^{n_v}$. The reduced trace conductor of $\gamma\in K$ is $\mathfrak q_\gamma$, of norm $q_\gamma$. For the original $z$-roots use a common narrow annular cell at every place. After two nonzero differences the coefficient is a fixed nonzero multiple of
\[
 h_1h_2Z\bigl(Z^2+(h_1m_1^5)^2+(h_2m_2^5)^2\bigr).
\]
At the step producing this term, orthogonality counts equalities with nonnegative coefficients. In this count one may delete the remaining derivative, divisor-label and endpoint restrictions and majorize the variables independently; only \emph{then} take the absolute value of the resulting sum $F_2$ against the common shortened sum. The retained $m_i$ are balanced principal divisors built out of smooth prime ideals. Retain this last restriction on $m_1$. Counting a label once per ideal, or with divisor-bounded multiplicity, gives the same estimates with $P^\epsilon$. The initial narrow cell implies $|(h_i)_v|\ll\delta H_i^{1/n}$, with $\delta$ arbitrarily small up to fixed constants, so the independent majorants preserve small unscaled shifts relative to $Z$ at all places. Dyadic archimedean moduli, unit balancing and Minkowski reduce initial boxes to such cells as in the preceding moment proofs. The overshoots in $M_i$ and the prime cutoff at each of the finitely many scales cost $P_0^{O(\eta)}$, where $P_0$ is the initial volume.

Cauchy over the four external labels and one ordinary difference $l\ne0$ in $Z$ separate a quadratic sum in $Z$ and decic sums in $m_1,m_2$. Write $D$ for the sum over $h_1,h_2,l$ of the squared quadratic sum at coefficient $c_0\alpha lh_1h_2$. Write $E_{i,s}$ for the same sum of the $2s$th absolute power of the $m_i$-sum at coefficient $c_i\alpha l h_{3-i}h_i^3$. The $c_i$ are fixed nonzero field elements. Smooth translated $Z$-cutoffs are uniform and independent of $m_i$. H\"older with $t=22$ gives
\[
 |F_2|^2\ll P B_2^2+
 B_2 D^{1/2} E_{1,t}^{1/(2t)} E_{2,1}^{(t-1)/(2t)}M_2^{1/t}.
\]
For each $E_{i,s}$ collect $j=lh_{3-i}$ first. Its divisor multiplicity is $P^\epsilon$ and the summand is nonnegative. Majorize its range by a positive smooth $j$-box of volume $J=P H_{3-i}$ \emph{before} expanding the $m_i$-power. For $E_{2,1}$ finite Cauchy over sectors on the balanced $m_2$-annulus makes the two modulus variables lie in the same narrow sector at every place.

Call the union of $q_\gamma B(Q^{5/n}(\alpha-\gamma))\le CP$ the broad arcs. Any fixed rescaling of $\alpha$ alters the conductor and the height by bounded factors, covered by the constant or by an arbitrarily small power buffer. Outside these arcs, a further ordinary difference of the quadratic gives
\[
 D\ll P^{2+\epsilon}H_1H_2.
\]
Indeed its linear kernel has length $P$ and its nonzero product variable $lh_1h_2v$ has volume $V=P^2H_1H_2$ and $P^\epsilon$ multiplicity. In a rectangular concentration box for $(a,\alpha a-\xi)$ of volumes $V,P^{-1}$, rank two over $K$ gives $O(P^\epsilon V/P)$ points by the scaled lattice determinant, and rank one has a common slope $\gamma=\xi/a$ and at most $O(P^\epsilon(1+V/[q_\gamma B((PV)^{1/n}(\alpha-\gamma))]))$ points. Here $PV=Q^5$. Multiplying by the kernel size and including the zero difference proves the bound.

\paragraph{The 44th decic mean for the first modulus.}
Put $H=H_1$, $M=M_1$, $J=P H_2$ and $V=H^3M^{10}$. In the decic 44th power the multiplicity of each discrepancy $d$ is bounded by the entire 44th moment, by Cauchy on the two 22-variable counting functions. It is $O(M^{34+\delta_{44}+\epsilon})$ for these actual upper-smooth moduli. Section~\ref{subsec:aud-full-moments}, node 81 (numbering from zero) of \path{code/full/full_k10.json}, gives $\delta_{44}=0.1276668$. All cutoffs can be chosen once over the finite set of shortened volumes.

The zero discrepancy contributes $O(JH M^{34+\delta_{44}+\epsilon})$. For $v=h^3d\ne0$ the divisor bound gives $P^\epsilon$ choices of $h$. Poisson summation in $j$ and the same rank-two/rank-one count, with kernel volume $J$, show that the sum of the nonzero kernels is $O(P^\epsilon(J+V+JV/P))$, since $JV=Q^5$. Consequently $E_{1,22}$ saves $P^{e_1}$ relative to $JHM^{44}$, where, setting $a_i=1-5\log_P M_i$,
\[
 e_1=\min(10x,1+a_2-2a_1,1-2a_1)-\delta_{44}x .
\]

\paragraph{A square mean for the second modulus.}
Here put $H=H_2$, $M=M_2$, $J=P H_1$, $L=JH^3$. Thus $LM^{10}=Q^5$. Take $Y=P^b$ satisfying
\[
 H^3<J,\qquad  Y\le\min(M,P^{1/5}),\qquad
 Y^6 M^9<L,\qquad Y^{20}<Q^5,
\]
with fixed power margins in the strict inequalities. We prove outside the broad arcs that $E_{2,1}\ll P^\epsilon JHM^2/Y$. In the following counts an integral ideal of norm $A$ has $O(1+V/A)$ elements in a skew box of volume $V$, by product-norm packing. For a centered box containing a nonzero element, $V\gg A$, and the bound is $O(V/A)$. All fixed fractional lattices and the codifferent change constants only.

After smoothing $j$, pairs with $m'^{10}=m^{10}$ contribute $O(JHM)$. At $d=m'^{10}-m^{10}\ne0$ Poisson summation permits truncation to $|J^{1/n}(\alpha dh^3-\xi)_v|\le P^\epsilon$ at every place. For fixed $d$, the determinant of any two pairs $(h^3,\xi)$ has absolute norm $O(P^{O(\epsilon)}H^3/J)=o(1)$ and belongs to a fixed fractional ideal. Thus all surviving pairs have one slope $\gamma_d=\xi/h^3$. The ideal $\mathfrak q_{\gamma_d}$ divides $(h)^3$, so $h$ belongs to an ideal of norm at least $q_{\gamma_d}^{1/3}$. The frequency constraint at every place cuts its radius by the factor $\max(1,|L^{1/n}(\alpha d-\gamma_d)_v|)^{-1/3}P^{O(\epsilon)}$. Centered packing bounds the number by $P^\epsilon H/A_d^{1/3}$, where
\[
 A_d=q_{\gamma_d}B(L^{1/n}(\alpha d-\gamma_d)).
\]
Terms $A_d>Y^3$ already fit the claimed bound.

Let $\mathfrak g=(m,m')$, $G=N\mathfrak g$, and write $\mathfrak g=cA$ using a fixed representative ideal $A$ and one balanced scalar $c$ per $\mathfrak g$. Thus $|c_v|\asymp G^{1/n}$. Write $m=cz$, $m'=c(z+l)$ with $(z,z+l)=A$, and put $T=M/G$. The nonzero $l$ has volume bound $T$. Both $z,z+l$ lie in one narrow sector, so
\[
 R_z=((z+l)^{10}-z^{10})/l,\qquad |(R_z)_v|\asymp T^{9/n}
\]
at every place, by integrating $10w^9$ on the short segment. For fixed $c,l$ the ratios $\gamma_d/R_z$ have one common value $\rho$. Indeed for two retained pairs,
\[
 q_{\gamma_i}q_{\gamma_j}|N(\gamma_i R_j-\gamma_j R_i)|
 \ll T^9 A_i A_j/L\ll M^9Y^6/L=o(1).
\]
If the determinant were nonzero the left side would be bounded below. This estimate and the subsequent comparisons use the product of the local maxima and hence allow skew archimedean errors.

Put $\mathfrak s=\mathfrak q_\rho$ and $C_{c,l}=N\mathfrak s\,B((LT^9)^{1/n}(\alpha c^{10}l-\rho))$. Since $\gamma_d=\rho R_z$, one has
$A_d^{-1/3}\ll C_{c,l}^{-1/3}N(\mathfrak s,(R_z))^{1/3}$.
For $\mathfrak e\mid\mathfrak s$ the condition $\mathfrak e\mid(R_z)$ costs at most $P^\epsilon$ classes of $z$ modulo $\mathfrak e$. At a good prime, both $z,z+l$ are units for a contributed class. If $l$ is a nonunit then $R_z$ is a unit. Otherwise $(z+l)/z$ is a tenth root of unity and determines $z$. At primes of $10A$ remove the fixed local common factor. For a remaining unit $z$ with $v(l)>v(10)$, the binomial expansion has $v(R_z)=v(10)$. When $v(l)\le v(10)$ the inverse multiplicity of $z\mapsto1+l/z$ and the unit tenth-power kernel are bounded independently of the prime-power precision. A single nonunit among $z,z+l$ gives a unit $R_z$. The Chinese remainder theorem gives the stated factor.

Packing each surviving class costs $O(1+T/N\mathfrak e)$. Divisors $N\mathfrak e>T$ contribute $O(P^\epsilon)$ for this $c,l$, since $N\mathfrak e\le N\mathfrak s\le C_{c,l}$. The remaining contribution is $O(P^\epsilon T/C_{c,l}^{1/3})$. The first term sums to $O(P^\epsilon M)$ since $\sum_{G\ll M} M/G\ll M\log P$, summing over ideals. In the second term discard $C_{c,l}>(Y/G)^3$ at total cost $O(P^\epsilon M^2/Y)$. This also covers $G>Y$.

For fixed $c$, the remaining ratios $\rho/l$ coincide, with value $\kappa$. The determinant bound is $O(T(Y/G)^6/(LT^9))=o(1)$ after multiplication by the two reduced conductor norms. Here $T\gg1$ for represented pairs. Write $\mathfrak t=\mathfrak q_\kappa$. For $\mathfrak e\mid\mathfrak t$, summing over nonzero $l\in\mathfrak e$ gives
\[
 \sum_l B((LT^9)^{1/n}l(\alpha c^{10}-\kappa))^{-1/3}
 \ll \frac{P^\epsilon T}{N\mathfrak e}\,
 B((LT^{10})^{1/n}(\alpha c^{10}-\kappa))^{-1/3}.
\]
To see this, dyadically localize $|l_v|/T^{1/n}$ and use centered ideal packing and $r/(1+ur)^{1/3}\ll(1+u)^{-1/3}$ for $0<r\ll1$. There are logarithmically many nonempty ranges by the integral norm lower bound. Since $\mathfrak s=\mathfrak t/(\mathfrak t,(l))$, this leaves $O(P^\epsilon\sum_{\mathfrak g}T^2/C_{\mathfrak g}^{1/3})$, where $C_{\mathfrak g}=N\mathfrak t\,B((LT^{10})^{1/n}(\alpha c^{10}-\kappa))$. Discarding $C_{\mathfrak g}>(Y/G)^3$ again costs $O(P^\epsilon M^2/Y)$.

Finally, $\kappa/c^{10}$ has a common value $\gamma$ across all retained ideals. The norm of $\kappa_i c_j^{10}-\kappa_j c_i^{10}$ times the two reduced conductor norms is
\[
 \ll (G_iG_j)^{10}C_i C_j/(LM^{10})\ll Y^{20}/Q^5=o(1).
\]
Writing $\mathfrak q=\mathfrak q_\gamma$, the total is now bounded, apart from $M^2P^\epsilon$, by
\[
 B(Q^{5/n}(\alpha-\gamma))^{-1/3}
 \sum_{\mathfrak g}(N\mathfrak g)^{-2}
 \bigl(N\mathfrak q/N(\mathfrak q,\mathfrak g^{10})\bigr)^{-1/3}.
\]
The sum is $\ll P^\epsilon (N\mathfrak q)^{-1/5}$ since
$\min(e,10j)/3-2j\le2e/15$ at every prime. The finitely many class factors have bounded valuations. The conductor norms used for divisor counts are polynomially bounded by the original denominators $h^3R_zlc^{10}$. Off the broad arcs the last display saves $P^{-1/5}$, finishing the square mean. The two modulus means and the quadratic estimate yield
\[
 |F_2|\ll P^{1-\sigma+\epsilon}B_2,\qquad
 \sigma=\tfrac14(1+e_1/22+21b/22)<\tfrac12
\]
at the indicated parameters.

\paragraph{A fourth moment on the broad arcs.}
At conductor $\mathfrak q$ the cubic in $Z$ has reduced leading conductor $\mathfrak q'=\mathfrak q/(\mathfrak q,(h_1h_2))$ up to a bounded factor. Poisson summation adds a linear twist modulo $\mathfrak q'$. At a good prime of norm $r$ the normalized sum is $\ll r^{-1/2}$ at exponent one by the finite-field Weil bound and $\ll r^{-1}$ at exponent two: summing over lifts of a residue leaves at most two roots of the quadratic derivative. For exponent $e\ge3$ use $r^{-e/3+e\epsilon}$ by \cite[Theorem~1]{audHuaField}. The leading coefficient has full conductor at each good prime, for every twist. Combine the high powers at their common modulus in Hua's theorem and the lower powers by CRT. Define the multiplicative weight $w$ by these three bounds without $\epsilon$ or constants, using $r^{-e/3}$ at all powers of the finitely many bad primes. Then
$w(\mathfrak q')\ll w(\mathfrak q)N(\mathfrak q,(h_1h_2))^{1/2}$ and
$\sum_{N\mathfrak q\le CP}N\mathfrak q\,w(\mathfrak q)^4\ll P^\epsilon$.
The local terms of the latter series at good exponents one and three are $r^{-1}$, at two they are $r^{-2}$, and from four on they are $r^{-e/3}$.

In every normalized narrow $Z$-cell the second derivative of the cubic is bounded away from zero after division by its leading coefficient. At a complex place the real Hessian determinant is minus the squared modulus of this derivative. Stationary phase at each place, uniformly in the linear twist, gives decay $B(P^{3/n}h_1h_2\beta)^{-1/2}$ at the stationary dual frequencies, with arbitrary additional decay outside their enlarged window. The dual box has volume in lattice units at most $N\mathfrak q' B(P^{3/n}h_1h_2\beta)/P\ll1$ on these arcs. Skew-box packing sums these estimates. Since $H_1H_2P^3=Q^5$, the archimedean factor is at most
\[
 B(Q^{5/n}\beta)^{-1/3}
        \prod_{i=1}^2 (H_i/|Nh_i|)^{1/3}
\]
up to a constant. This follows place by place from
$\max(1,a t)\gg a\max(1,t)$ for $0<a\ll1$ and $t\ge0$,
after weakening the exponent $1/2$ to $1/3$.
For each $i$ the requisite divisor average is
\[
 \sum_{h_i\ne0}
 (H_i/|Nh_i|)^{1/3}N(\mathfrak q,(h_i))^{1/2}\ll H_iP^\epsilon.
\]
Indeed the inner sum over $\mathfrak e\mid(h_i)$ without the gcd weight is $O(H_iP^\epsilon/N\mathfrak e)$, again by dyadic ideal packing. Thus $|F_2(\gamma+\beta)|\ll P^{1+\epsilon} B_2 w(\mathfrak q) B(Q^{5/n}\beta)^{-1/3}$. There are at most $N\mathfrak q$ centers per conductor and the archimedean fourth-power integral costs $Q^{-5}$. The broad fourth moment is $O(P^\epsilon(PB_2)^4 Q^{-5})$.

\paragraph{Finite closure.}
Write $\lambda_s$ for the prior fifth-power exponent at order $s$, $q=14,16$, $u=q-2$, $t'=2u-12$, $Q_1=P/M_1$. Set
\[
 \begin{split}
 W&=P B_2\max\{P^{-\sigma}Q^{\lambda_{12}},
                              Q^{(3\lambda_{16}-5)/4}\},\\
 C&=\max\{B_1^2P^2M_2^2 Q_1^{\lambda_{10}},\ 
 B_1 M_2^{11}\max(P B_1M_2 Q^{\lambda_{12}},W)\}.
 \end{split}
\]
The exponent of the mixed mean with one ordinary pair and $q-2$ smooth roots is bounded by the base-$P$ logarithm of
\[
 \max\bigl\{P^2M_1^2 P^{\lambda_{u-2}},\ P M_1 P^{\lambda_u},\
 P M_1^u Q_1^{\lambda_u},\ M_1^{q-3}(Q_1^{\lambda_{t'}}C)^{1/2}\bigr\}.
\]
The initial ordinary derivative exceptions are included in the first term. The first and second zero shifts and the second exceptional pairs are precisely those of the two-block congruence estimate. The two terms in $W$ come from the complementary supremum and from H\"older on the broad arcs with exponents $4,4/3$. Whole-torus positive-coefficient domination bounds the pure mean by this mixed mean, allowing feedback from the preceding round.

Use the fifth-power seed table of Section~\ref{subsec:aud-block-moments}, including $\lambda_{10}\le5.9250801$, $\lambda_{12}\le7.5417549$. For orders $14,16$ respectively take $(10^5x,10^5y,10^5b)=(16342,13943,13930)$ and $(16951,13842,13820)$. Eighty simultaneous rounds, each rounded upwards in units of $10^{-8}$ and intersected with the prior bound, give
\[
 \lambda_{14}\le9.27034699,\qquad \lambda_{16}\le11.07208133.
\]

The file \path{code/local/quintic_hybrid.py} checks this recurrence
and all of its parameter inequalities, as described in Appendix~\ref{app:certificates}.

In particular $D_{16}=\lambda_{16}-11<2/27$. Critical absorption gives the loss-free order $18$ since $2-12D_{16}>0$. Interpolation yields $D_{225/14}\le(27/28)D_{16}$. For two ordinary and fifteen of the positive marked roots, the ordinary critical order $30$ gives minor-arc saving
\[
 \tfrac1{15}-\tfrac{14}{15}D_{225/14}
 \ \ge\ \tfrac1{15}-\tfrac9{10}D_{16}
 \ \ge\ \frac{5380409}{3000000000}>0.
\]
The pruned ordinary order $2\cdot18/(18-15)=12$ exceeds $10$. This supplies the positive mixed saving used in Theorem~\ref{thm:nf-small}.

\subsection{Fractional orders for ninth and tenth powers}
\label{subsec:aud-high-moments}
We start with the full upper-smooth even estimates of
Section~\ref{subsec:aud-full-moments}. Even coefficient domination
passes them to the hereditary family; a new differencing step at a
better even exponent first groups its exceptional relative ideal.
All subsequent interpolations concern the same restricted root weight.

\paragraph{Additional polynomial estimates.}
Take \(k=9\) or \(10\). Write \(M_i=P^{\theta_i}\), \(H_i=P/M_i^k\), \(Q_i=P/(M_1\cdots M_i)\), \(B_i=\prod_{j\le i}H_jM_j\).
Use the regular first and second block extractions and derivative-failure estimate proved in Section~\ref{subsec:aud-block-moments}, with \(k\) in place of seven.
All first differences have a sharp second-moment estimate over \(K\). After division by the first nonzero shift the centered degree-\((k-1)\) binary form splits over \(K(\mu_k)\) into distinct linear factors in \(Z,h_1m_1^k\). At a prescribed nonzero frequency two nonproportional factors have only \(P^\epsilon\) possible pairs of values by ideal divisors and bounded-height units. The normalized small root box, with all physical shifts small relative to its center in one embedding, excludes zero residuals. At \(k=9\) a single genuine difference followed by an ordinary difference also has the sharp second-energy estimate. Its residual has odd degree seven, so its center \(Z\) divides the frequency after clearing the fixed denominators. Once \(Z,h_1,h_2\) are fixed, the residual is a polynomial in \(m_1\) with nonzero leading term a fixed multiple of \(Z(h_1m_1^9)^6\). This suffices for the fourth moment of the single-genuine-block sum \(F_1\).

For \(k=9\) consider two genuine differences. Fix a nonzero frequency. The two shifts \(h_1,h_2\) and their nonzero center \(Z\) can again be selected as divisors, costing \(P^\epsilon\). Let
\[
 s_i=h_i m_i^9,\quad U=s_1^2,\quad V=s_2^2,\quad
 A=U+V,\quad D=U-V,\quad
 W=Z^{-1}{\mathbb E}(Z+s_1T_1+s_2T_2)^7 ,
\]
where \(T_i\) are independent uniform real variables on \([-1,1]\); the expectation denotes a rational polynomial. The frequency determines \(W\). Expansion gives the exact factorization
\[
 27W-272Z^6=(3A+7Z^2)
    (18A^2+42AZ^2-35Z^4-9D^2).
\]
The left side and the first factor are nonzero. Indeed \(W/Z^6\) can be kept arbitrarily close to 1 at one embedding, and \(A/Z^2\) arbitrarily close to zero. First restrict the original root to a small box relative to its center, divide the \(M_i\) into dyadic norm bins and use balanced generators. All subsequent enlargements of the boxes for the independent shifts cost bounded factors, so this condition can be imposed in the original finite cover. Finite denominators in the identity are immaterial. The ideal divisor bound and the unit lattice give \(O(P^\epsilon)\) possibilities for the bounded generator \(3A+7Z^2\). It determines \(A\), then \(D^2\), and the two possible signs determine \(U,V\). Extracting the two nonzero field \(18\)th powers has bounded multiplicity. Thus
\(\int |F_2|^2\ll P^{1+\epsilon}B_2\)
over every fixed field. For the other higher moments a weaker argument suffices. Fix every genuine modulus, select all the nonzero shifts by ideal divisors of a given frequency, then solve the residual positive-degree polynomial for its center. Cauchy across these moduli costs \(\prod M_i\) beyond the squared-energy tuple count. The integrated Weyl-difference induction charges this extra factor once. Consequently
\[
 \int |F_j|^{2^\ell}\ll
 P^{2^\ell-\ell+\epsilon} B_j^{2^\ell-1}\,C_{j,\ell}.
\]
Here \(C_{1,\ell}=1\) at \(\ell=1\), also at \(\ell=2\) when \(k=9\), and \(C_{1,\ell}=M_1\) otherwise. Take \(C_{2,1}=1\) at \(k=9\) and \(C_{2,\ell}=M_1M_2\) in every other two-block case. Only positive residual degrees are used.

Only the second and fourth moments of the differenced sums occur in the final certificates. In the nonic lower table the nontrivial rows comprise 14 sharp first-difference second moments, 19 sharp first-difference fourth moments, 19 sharp two-block second moments and 48 two-block fourth moments with the two-modulus loss. In the decic table the corresponding counts are 3, 5, 30 and 60, with the conservative modulus losses just stated. Every other lower update in the retained tables is the efficient-differencing recurrence. The zero- and first-difference terms and the two-block exceptional and diagonal terms retain the exponents proved in Section~\ref{subsec:aud-block-moments}. Each lower-wave rule invokes only previous-wave estimates. The positive solution-count interpretation permits bounded-coefficient subsets in all even stages, including each of the derivative masks.

\paragraph{Absorption at intervening even orders.}
An additional boundary estimate is available at every even \(2k<q\le w=k(k+1)\). For a prior real order \(h<q\) with \(h(w-1)\le w(q-1)\),
\[
 D_q\le \max\left\{0,\frac{(w-q)D_h-(q-h)}{w-h}\right\}.
\]
For coefficients of arbitrary sign, dominate their even \(q\)th moment first by the positive coefficient moduli. In its zero Fourier coefficient replace one signed position by an ordinary positive envelope. On the broad minor arcs apply H\"older to its \(w\)th moment, the \(h\)th moment of the marked sum, and its own \(q\)th moment. The three weights are \(1/w\), \((w-q)/(w(q-h))\) and \(1-(w-h)/(w(q-h))\); the displayed hypothesis makes all three nonnegative. The ordinary critical integral saves one full power of \(P\). The broad major term is \(O(P^{(q-k)/q}(\int|F|^q)^{1-1/q})\) since \(q>2k\). Absorbing both terms gives the stated exponent, with all costs from fixed windows included in the prepaid loss. Higher-node witnesses otherwise use the previously justified exact height-interval and real interpolation rules.

The rational verifier checks 792 nonic and 878 decic lower-wave improvements. It is \path{code/fibres/verify_quoted_gains.py}. The two inputs are \path{code/fibres/quoted_k9.json} and \path{code/fibres/quoted_k10.json}. The higher proof trees have 66 and 87 nodes, using 872 and 1,225 height intervals, respectively. Only the two elementary conductor families of Lemma~\ref{lem:nf-fibre-savings} are needed. After exact rational interpolation the three estimates and their power savings are
\[
\begin{array}{c|rr|c|c|c}
 k&a&b&r&\text{upper bound on }D_r&
       a/w-(1-a/w)D_r\ \text{at least}\\ \hline
 9&7&37&3330/83&131335359/1660000000&8664641/1800000000\\
 9&6&38&285/7&1325421/20000000&722053/150000000\\
10&7&44&4840/103&283006357/5150000000&66993643/5500000000
\end{array}
\]

The loss-free orders 46 and 56 therefore give annular orders
$322/9,69/2,98/3$, respectively. All exceed $2k$.
\subsection{Further marked moments for ninth and tenth powers}
\label{subsec:refined-small-moments}

At fixed degree the distribution of the local conductor orders gives
sharper estimates for the Gram entries in the shortening step of
Section~\ref{sec:nf-fibres}. Write $q\mid m^k$ for the conductor,
in the ideal notation of Section~\ref{sec:joint-conductors}, and
$U=\log_X\mathcal N(q)$. Put $\theta=M/k$. If $\theta_e$ is the sum
of $\log_X\mathcal N(\mathfrak p)$ over the marks of order $e$ in $q$,
pad at order zero so that $\sum_e\theta_e=\theta$ and
$U=\sum_e e\theta_e$. Inactive unit masks can be removed by
inclusion--exclusion. The resulting progressions have steps prime
to the remaining conductor. Their volume is at most $X$, and the
local polynomial is changed by a unit affine substitution.

\paragraph{The Taylor estimates.}
For Taylor degree $D$, choose a step with valuation
$s_e=\lceil e/(D+1)\rceil$ at a prime of order $e$. Its norm is $A$,
and the progression volume is $Y=X/A$. Let $R_i$ be the denominator
ideal norm of the coefficient of degree $D+1-i$, so
\[
 \log_X R_i=\sum_e\theta_e(e-(D+1-i)s_e)_+ .
\]
For the highest coefficient, sampling frequencies in the Taylor
step gives at most $O_{\mathbb K}(1+Y/R_1)$ occurrences of each
frequency and $O_{\mathbb K}(1+R_1/Y^{D-1})$ frequencies in a cell.
The normalized count is bounded by a constant times
$Y^{-1}+R_1^{-1}+R_1Y^{-D}$. Apply the critical field mean value at
degree $D-1$, or $D-1$ elementary differences, giving the saving
\[
 \frac{\min\{\log_XY,\log_XR_1,D\log_XY-\log_XR_1\}}
      {\min\{2^{D-1},D(D-1)\}} .
\]
For the two-coefficient test, and its extensions, condition successively
on equality of the first $j$ coefficients. A bounded colouring of
the reduced power maps makes each relevant map injective on each
colour. With $R_0=1$ and lower-coordinate tolerance of volume
$\delta_i=Y^{-(D+1-i)}$, the number of differences in the $i$th
coefficient after $i-1$ equalities is
$O_{\mathbb K}(1+\delta_i R_i/R_{i-1})$. For a fixed difference,
Cauchy's inequality on the occupancies of translated coefficient
cells bounds the cross-correlation by the exact-coincidence count
with a bounded enlargement of the lower windows. In that count the
root difference is divisible by the ideal of norm $R_j$. If the
next coefficient is also retained, the valuation and Hensel
calculation of the two-coefficient estimate gives the normalized
bound
\[
 X^{-1}+R_{j+1}^{-1}
 +\delta_{j+1}\left(\frac{R_{j+1}}{X R_j^2}+\frac1{R_j}\right).
\]
Indeed for each shift in the ideal of norm $R_j$, the first $j$
coefficient equalities fix its residue class, and the next one
fixes the root modulo the quotient of its coefficient modulus by
the shift gcd. Summing over those gcds gives the displayed terms
by ideal packing, as in \eqref{eq:joint-taylor-saving}
and \cite[Lemma A.8]{Wnew}. Multiply by the product of the preceding
difference counts and take the critical $D(D+1)$th root. The
diagonal term imposes $D(D+1)\sigma\le1$ at saving $\sigma$.

All these conductor tests are affine in the $\theta_e$ after the
positive parts in the factors are expanded. At a prescribed $U$,
the minimum of a linear form $\sum_e\theta_e c_e$ is $\theta$
times the lower convex hull of $(e,c_e)$ evaluated at $U/\theta$.
For tests following $j$ arithmetic differences, the shift at
stage $i$ can have norm exponent in $[\theta,1-2^i\sigma]$.
Each active prime is used at least once per step. Additional
factors are assigned one at a time, with overshoot at most the
largest individual prime weight. The surviving reduction is an
affine monomial of degree $k-j$. Before a terminal Taylor test is
used in this way, its adjacent-order differences are checked to
preserve every length gap and to lose at most one prime weight
per extra drop of a conductor order. The divisor average in the
arithmetic $A$ process then pays for the incidental gcds.

\paragraph{A completion between groups of differences.}
Use $i$ outer and $j\ge1$ inner shifts. With two terminal valuations
reserved, remove the conductor primes of order at most $e_0=i+j+1$
by summing separately over their cosets. Their conductor weight $b$
satisfies
\[
 0\le b\le \min\{U,e_0(M-U)/(k-e_0)\},
\]
by $\sum_e(k-e)\theta_e=M-U$. The remaining radical has weight at
most $\theta-b/e_0$. The new length exponent is $L=1-b$. If the
outer shifts spend $C_0$, the first completion has conductor
$P=U-b-C_0$ and dual volume exponent $Y_1=P-L$.
For a saving $\sigma$ in the original sum it suffices to save
$\omega>2^i\sigma-L+P/2$ on the dual sum, in units of $\log X$.
Its $h$th shift has budget at most $Y_1-2^h\omega$ and must
contain the remaining radical. If the inner shifts spend $C_1$,
the terminal conductor is $P-C_1$. At length exponent $Y_1$
the terminal quadratic estimate saves at least
$\min\{(P-C_1)/2,Y_1-(P-C_1)/2\}$.
One checks the shift budgets and the reserve by eliminating
$C_0,C_1,\omega$ from these affine inequalities. The bounds on
$b$ allow the resulting affine conditions to be tested at the
two endpoints of its interval.

For this completion the stationary parameters at a surviving prime
are unit roots of a monomial of degree $k-i-1$. Label the possible
branches within each frequency residue by at most $k-i-1$ indices.
Each shift contains the prime, so its branch and the Gauss factor
are fixed under the difference. The derivative of the terminal
phase, after $j$ inner shifts, is a nonzero power of the original
stationary root on the residue field. Thus a fixed Fourier label
has boundedly many stationary residues. The second derivative on
each branch is a unit, giving the square-root bound for conductor
exponents at least two, even with the residue weights.
An incidental drop to order at most one retains a bounded residue
weight of period $p$. For nominal order $E\ge2$ and actual order
$e$, the change in the two completion terms costs at most
$p^{(E-\max(1,e))/2}$, with an extra factor $p^{1/2}$ allowed
when $e\le1$. This is bounded by $p^{E-e}$. Averaging the extra
gcd at the step where it occurs absorbs this loss. All archimedean
volumes here are volumes of balanced boxes. Completion on the
trace lattice and packing in an ideal coset give the same
exponents as for the degree-one completed local sums.

\paragraph{A one-valuation reserve.}
For marked sums the prime windows can be restricted to make the
terminal estimate valid also at order one. Let $L_0/\mathbb Q$ be
the normal closure of the narrow Hilbert class field of $\mathbb K$.
Let $\mathcal S$ comprise $2$ and the prime divisors of the finitely
many stationary power exponents in the completion words. Choose
$m_0=\prod_{\ell\in\mathcal S}\ell^{r_\ell}$ with $r_\ell\ge2$ so
that the intersection of $L_0$ with the union of the cyclotomic
fields supported on $\mathcal S$ lies in $\mathbb Q(\mu_{m_0})$.
Write $m_1=m_0\prod_{\ell\in\mathcal S}\ell$. Restriction identifies
$\operatorname{Gal}(L_0(\mu_{m_1})/L_0(\mu_{m_0}))$ with the
product of the groups of order $\ell$ from the last cyclotomic
digits. Choose an element nontrivial in every factor. Chebotarev
gives a positive density of rational primes splitting completely
in $L_0$ with $v_\ell(p-1)=r_\ell$ for all $\ell\in\mathcal S$.
Their marked ideals are unramified, narrow principal, and of
degree one.

The subgroup of $\mathbb F_p^\times$ of index $m_0$ has odd order
coprime to the stationary exponents. Each power map is injective
on each of its cosets, and the quadratic Gauss factors are constant
there. In the stationary parameter $s$ an order-one complete sum
has terms $e_p(\lambda s^B+\mu s^A)$ with $B\ne0$, $A\ne B$ and
$\lambda\ne0$. Substitute $s=c z^{m_0}$ and apply the Weil bound
to this nonconstant rational phase to get $O_{\mathbb K,k}(\sqrt p)$.
The relevant pole orders are prime to $p$ after omitting finitely
many primes. An incidental drop below order two contributes a
bounded residue weight to the next shift, with loss at most
$p^{(E-e)/2}$ in each Fourier estimate. The Fourier prefactor and
the final shifted estimate together cost at most $p^{E-e}$.
Consequently the terminal reserve can be one valuation, with
$e_0=i+j$. The first completion still has nominal exponent at
least two on each surviving prime.

The outer shifts can remove the low-order primes directly.
Keep the highest-order primes of order greater than $i+j$ up to
a radical weight $h_0$, or keep all of them if their weight is
smaller. The minimal outer budget for removing the other primes
and giving each survivor $i$ factors is at most
\[
 \max\{(i+j)\theta,\ k\theta-(k-i)h_0\}.
\]
In the first case use the fact that the retained orders have
mean at least $U/\theta$ after padding by zero. In the second
case each removed prime has order at most $i+j$. With
$e_0=i+j$, let $J$ be the radical weight of primes of order
greater than $e_0$. Then $U\le e_0\theta+(k-e_0)J$.
When all these primes are kept, their conductor after the
compulsory outer shifts is at least $U-e_0\theta+(e_0-i)J$.
For $U\ge e_0\theta$ this permits replacing the first entry
of the budget by
$U-(k-i)(U-e_0\theta)/(k-e_0)$. If
$U>e_0\theta+(k-e_0)h_0$, more high-order weight is available
than $h_0$, so the first entry can be omitted. Distribute the
compulsory factors first, then the remaining factors up to the
shift capacities. After these shifts the first completion uses
$P=U-C_0$ and $L=1$. Removing inactive masks may shorten an input
box from volume $X$ to $Z$. Tile the resulting Fourier weight
of volume $X^P/Z$ by translated smooth boxes of volume $X^{P-1}$.
The factor $O(X/Z)$ in the weighted tile count cancels the
$Z/X$ in the completion prefactor. This gives the same affine
budget test with surviving weight $h_0$.

\paragraph{Taylor estimation after completion.}
One may instead sum over the order-one conductor cosets and
complete on the remaining primes. With $L=1-\theta_1$ and
$W=U-\theta_1$, the dual volume exponent is $U-1$. To obtain
saving $\sigma$ it suffices to save
$\sigma-1+(U+\theta_1)/2$ in units of $\log X$ on that
dual box. On a stationary branch, the input residue is an
affine unit change of $s^{k-1}$ in the original unit root
$s$. The Taylor coefficient of order $j$ has unit-scaled
reduction $s^{k-j(k-1)}$ up to affine terms. The Taylor step
contains the radical, so the bounded residue weights are
fixed in the shifts. Apply the preceding colouring to
these nonzero power exponents. For a unit shift the initial
residues of two stationary roots $s,t$ satisfy
$t^{k-1}-s^{k-1}=c$ and $t^{B}-s^{B}=c'$,
with $c$ a unit and $B=k-j(k-1)$. The colour excludes the
ratios where the divided-difference derivative vanishes.
Thus the solutions are isolated, and their number is bounded
in terms of the degrees after clearing unit denominators.
For a shift divisible by $p$, the single-disk divided-difference
derivative is a unit and the initial count is a bounded-degree
power equation. In either case Hensel lifting gives the
same conditional counts as before. A Taylor row
$c_0+\sum_e v_e\theta_e$ with denominator $D_0$ becomes
\[
 c_0(U-1)+\sum_{e\ge2}v_e\theta_e
 \ \ge\ D_0\bigl(\sigma-1+(U+\theta_1)/2\bigr).
\]
The lower hull again verifies the row over the entire conductor
simplex, and the length conditions retain their strict gaps.

\paragraph{The zero coefficient in the rational case.}
In the rational calculation one may in addition select marks
$p\equiv3\bmod4$ and $p\equiv2\bmod\ell$ at each odd prime
$\ell$ dividing a tested stationary exponent or residual odd
degree. For a monomial after the outer shifts, its complete zero
Fourier coefficient vanishes at order at least two, since the
unit derivative precludes stationarity at zero. At order one
it is $-1$ for odd residual degree by the permutation of unit
classes, and $O_k(\sqrt p)$ for even degree. The nonzero Fourier
terms after smooth completion contribute $O_k(\sqrt q)$ up to
bounded per-prime factors. For a nominal remaining conductor
exponent $v$, the divisor average therefore permits savings
$\min\{v,1-v/2\}$ in odd degree and
$\min\{\max(v/2,v-\theta/2),1-v/2\}$ in even degree. The extra
$\theta/2$ accommodates the case in which all the remaining
orders are one. This supplies the additional rational conductor
test without placing those particular congruence restrictions on
the number-field marks.

\paragraph{The moment calculations.}
The preceding tests are used in the shortening argument with the
same marked family at all orders of each certificate. At amplitude
$X^{1-v}$, an available defect $D_h$ gives the affine estimate
$D_h+(h-r)v$ at order $r$. The shortened state gives
$(1-L/k)D_h+(h+2-r)v$ where the conductor test saves more than
$2v$. These are the two contributions in the proof of
\eqref{eq:nf-short}. The certificates cover the required
conductor and amplitude intervals by exact rational inequalities,
starting with the field moments above, followed by real
interpolation and even-order absorption.

For the decic field moment at order $48$, the first witness gives
\[
 D_{48}\le\frac{2021960133}{61250000000},\qquad
 D_{330/7}\le\frac{471775323}{10000000000}.
\]
The inequality $50-48-(110-50)D_{48}>0$ gives the loss-free
order $50$. With $5$ ordinary and $45$ marked roots the second
bound gives the saving $92718217/220000000000$. The separate
even allocation uses the one-valuation completion and gives
\[
 D_{605/13}\le\frac{1068816268887}{20000000000000}.
\]
It provides a saving $1105388504469/275000000000000$ for $6$
ordinary and $44$ marked roots. The loss-free order is again $50$,
and both allocations have ordinary pruning order $50$.

For the nonic field comparison, interpolation gives
\[
 D_{684/17}\le\frac{6421587551}{95000000000},
 \qquad D_{1474/39}\le\frac{673937309}{5000000000}.
\]
The second bound gives the loss-free order $44$ by absorption.
Use the lower convex hull of the verified defect bounds, adjoining
the loss-free order $44$, to test the allocations in
Section~\ref{subsec:nf-classical-comparison}. At $a+b\le43$ the
largest margin from these bounds is
$-571694103273/72500000000000$, achieved at $(3,40)$.

The corresponding calculation over $\mathbb Q$ gives
\[
 D_{684/17}\le\frac{1266836410299}{21800000000000},
 \qquad D_{1979/54}\le\frac{1587157003}{10000000000}.
\]
The latter gives the loss-free order $44$. The allocation
$(5,38)$ has exponent saving $263781024917/392400000000000$
and pruning order $110/3>18$. The local and central-arc
arguments of Theorem~\ref{thm:nf-small} apply over $\mathbb Q$,
giving $G(9)\le43$. Appendix~\ref{app:certificates} specifies
the exact arithmetic checks for the four moment chains and the
comparison of the field allocations.

\section{Sums of totally positive cubes}\label{nf3:cubes}
Put $d=[\mathbb K:\mathbb Q]$ and $E=\mathbb K\otimes_{\mathbb Q}\mathbb R$.
Given a bounded open set $\Omega\subset E$ in the totally positive region,
let $R_{\Omega}(X)$ count distinct algebraic integers in $X\Omega$ expressible
as the sum of three cubes of totally positive algebraic integers.
The rational estimate of Davenport \cite{nf3:dav} admits a version with the
same exponent in the volume of the target region.

\begin{theorem}\label{nf3:davenport}
For every nonempty $\Omega$ as above and every $\varepsilon>0$,
\[
 R_{\Omega}(X)\gg_{\mathbb K,\Omega,\varepsilon}(X^d)^{47/54-\varepsilon}.
\]
One may require all conjugates of each summand's cube root to have modulus
$\asymp_{\mathbb K,\Omega}X^{1/3}$.
\end{theorem}

\begin{proof}
Write $A=X^{1/3}$ and $N=A^{3d}$. Take a point of $\Omega$ with nonzero
complex coordinates and choose an archimedean cube root of one third of
this point, positive at every real place. A sufficiently small box $B$
about this root has the following properties. Three cubes of elements of
$B$ sum into $\Omega$. Two of their cubes have nonzero sum.
The cube map is injective on $B$, and sums of two or four elements of $B$
avoid zero. Constants below may depend on this box and the field.

Let $P=N^{1/18}$ and $Q=N^{1/54}$. For each principal prime ideal of norm
in $[P,2P]$, take one generator $p$ whose absolute values at the embeddings
are all $\asymp P^{1/d}$, and do the same for the prime ideals of norm in $[Q,2Q]$ to choose $q$.
Dirichlet's unit theorem gives these choices by reduction in the logarithmic
unit lattice. The prime ideal theorem in ideal classes supplies
$\gg P/\log P$ and $\gg Q/\log Q$ choices
\cite{nf3:neukirch}. We take $N$ large enough that these primes are distinct
and avoid $6$.

We use two elementary counting facts. A congruence class modulo an ideal
with balanced generator of norm $M$ contains $O(1+W^d/M)$ elements of
embedding height at most $W$. It has $\asymp W^d/M$ points in a fixed
box scaled by $W$ when $W/M^{1/d}\to\infty$. Indeed, division by the
generator turns the class into a translate of the integer lattice and
rotates and dilates the box within a compact family of shapes.
Also, for nonzero integral $m$ of height $N^{O(1)}$, there are
$O_\varepsilon(N^\varepsilon)$ factorizations $m=t_1\cdots t_j$ into
a fixed number of algebraic integers of height $N^{O(1)}$. To see this,
distribute the prime ideal powers of $(m)$ among its factors. The ordinary
divisor bound and the degree bound on the number of prime ideals above a
rational prime give $O_\varepsilon(N^\varepsilon)$ choices. For any one
principal ideal the bounded generators have logarithmic absolute values
in a translate of the unit lattice in a box of side $O(\log N)$, by the
norm lower bound. This costs a fixed power of $\log N$. In particular,
a nonzero sum or difference of two cubes of bounded height has
$O_\varepsilon(N^\varepsilon)$ representations, since once $b\pm e$ is
specified, the remaining equation is quadratic in $b$.

For each $p$ choose one representative in the nonzero residue field
from each fibre of the cube map. The resulting set of residues contains
at least $(\mathrm N(p)-1)/3$ elements, where $\mathrm N(p)$ denotes
the norm of the ideal $(p)$. Count labelled representations
\[
 t=a^3+(pb)^3+(pqc)^3,\qquad a,pb,pqc\in AB,
\]
with $a$ in the selected residue classes modulo $p$.
Let $r_{pq}(t)$ count these and set $r(t)=\sum_{p,q}r_{pq}(t)$.
Counting the three roots in their residue classes gives
\begin{equation}\label{nf3:first}
 N^{1-\varepsilon}/P\ll_\varepsilon \sum_t r(t)\ll N/P.
\end{equation}
Indeed, for fixed $p,q$ the count is $\asymp N/(P^2Q)$, and the
logarithms in the number of prime ideals are absorbed in $N^\varepsilon$.
The labels cause no loss of support.

By Cauchy's inequality,
$\sum_t r(t)^2\ll PQ\sum_{p,q,t}r_{pq}(t)^2$. In the last sum denote
the two triples by $(a,b,c)$ and $(a',e,f)$. Solutions with any of
$a=a'$, $b=e$ or $c=f$ contribute $O_\varepsilon(N^{1+\varepsilon}/P)$:
on fixing the first triple, the two remaining cubes have a fixed nonzero
sum, and the divisor estimate applies. For every other solution,
$a^3\equiv (a')^3\pmod {p^3}$ implies first $a\equiv a'\pmod p$ by
choice of residues, and then $p^3\mid a-a'$ by invertibility of
$a^2+aa'+(a')^2$ at $p$. Put
\[
 u=(a-a')/p^3,\qquad v=a+a'.
\]
Thus $u\ne0$, $H(u)\ll A/P^{3/d}$ and
\begin{equation}\label{nf3:identity}
 u(3v^2+p^6u^2)+4(b^3-e^3)=4q^3(f^3-c^3).
\end{equation}
The two cube differences in this identity are nonzero, since cubing is
injective on the box of the original roots.

For fixed $p,q,u,c,f$ let $S$ count the admissible $v,b,e$.
The pair $u,v$ determines $a,a'$. Split the sum of $S$ into $T_0$ and
$T_1$ according as $q\nmid u$ or $q\mid u$. The numbers of possible
outer tuples satisfy
\[
 M_0\ll N/(P^4Q),\qquad M_1\ll N/(P^4Q^2).
\]
Here $u$ has $O(A^d/P^3)$ choices, and $c,f$ each have $O(A^d/(PQ))$.
The restriction $q\mid u$ saves $Q$ since $A^d/P^3\gg Q$.

In $\sum S^2$ the pairs with identical values of $v$ contribute
$O_\varepsilon(N^\varepsilon T_i)$ in either sum, since after fixing one
solution, $b^3-e^3\ne0$ is fixed for the second. For unequal values of
$v$ subtraction in \eqref{nf3:identity} gives, up to sign,
\[
 3u(v_1-v_2)(v_1+v_2)=4(b_1^3-e_1^3-b_2^3+e_2^3).
\]
The left side is nonzero since $v_1+v_2$ is the sum of four roots in
$AB$. For fixed $b_1,e_1,b_2,e_2$ there are thus $O_\varepsilon(N^\varepsilon)$
choices for $u,v_1,v_2$. In $T_0$, for each $q$, reduction of
\eqref{nf3:identity} modulo $q^3$ fixes the unit $p^6$ modulo $q^3$.
There are at most six possible residue classes by Hensel's lemma.
Each contains $O(1+P/Q^3)$ possible $p$. Then $c^3-f^3\ne0$ is
determined. The off-diagonal contribution to $\sum S^2$ is bounded by
\[
 E_0\ll_\varepsilon N^\varepsilon (A^d/P)^4(Q+P/Q^2).
\]
In $T_1$ both $b_j^3\equiv e_j^3\pmod q$. For fixed $q$ each pair has
$O((A^d/P)^2/Q)$ choices, since a residue has at most three cube roots
and $A^d/P\gg Q$. After the divisor estimate for $u,v_1,v_2$, allow
all $O(P)$ choices of $p$ and again use the nonzero cube difference for
$c,f$. This gives
\[
 E_1\ll_\varepsilon N^\varepsilon (A^d/P)^4 P/Q.
\]
Now $T_i^2\ll_\varepsilon M_i(N^\varepsilon T_i+E_i)$, whence
$T_i\ll_\varepsilon M_iN^\varepsilon+(M_iE_i)^{1/2}$. Substitution
including the identical-coordinate cases gives
\[
 \sum_t r(t)^2\ll_\varepsilon N^\varepsilon
 \left(NQ+N^{7/6}Q/P^3+N^{7/6}/(P^{5/2}Q^{1/2})\right)
 \ll_\varepsilon N^{55/54+\varepsilon}.
\]
Together with \eqref{nf3:first} this gives
$\#\{t:r(t)>0\}\gg_\varepsilon N^{47/54-\varepsilon}$, after rescaling
$\varepsilon$. All three original roots lie in the prescribed positive box.
\end{proof}

\subsection{Variation of the density with the field}
There are two elementary restrictions on a uniform density assertion.
The first applies in a fixed degree and comes from the covolume of the
ring of integers.

\begin{proposition}\label{nf3:covol}
Let $\mathbb K$ be totally real of degree $d$ and absolute discriminant
$D_{\mathbb K}$. In a bounded Jordan-measurable target region $X\Omega$ of
positive volume in the positive orthant, the upper proportion of algebraic
integers that are sums of $k$ totally positive integral $k$-th powers is at most
\[
 \frac{\Gamma(1+1/k)^{kd}}{k!\,D_{\mathbb K}^{(k-1)/2}}.
\]
For $k=3$ no positive lower proportion, even relative to the everywhere
locally representable targets, can be uniform over real quadratic fields.
\end{proposition}
\begin{proof}
Put $P=X^{1/k}$. At one real place the pushforward of Lebesgue measure on
the positive orthant of $\mathbb R^k$ by $\sum x_j^k=t$ has constant
density $\Gamma(1+1/k)^k$ for $t>0$, by the beta integral and homogeneity.
At $d$ real places these densities multiply. The preimage of $\Omega$ is
bounded, and its boundary has measure zero because the target boundary has
pushforward measure zero, and the coordinate walls have measure zero.
Lattice counting at scale $P$ for the $k$ copies of the Minkowski lattice
gives the number of ordered representations asymptotic to
\[
 X^d\operatorname{vol}(\Omega)\Gamma(1+1/k)^{kd}D_{\mathbb K}^{-k/2}.
\]
The target count is asymptotic to $X^d\operatorname{vol}(\Omega)D_{\mathbb K}^{-1/2}$.
Repeated coordinates contribute $O(P^{(k-1)d})$ ordered representations,
and each remaining unordered tuple has $k!$ orderings. Division proves the bound.

Take $\mathbb K=\mathbb Q(\sqrt\ell)$ with a prime $\ell\equiv1\pmod {24}$.
Here $D_{\mathbb K}=\ell$ and both $2$ and $3$ split. The order generated
by cubes contains $6\mathcal O_{\mathbb K}$ by
$(t+1)^3+(t-1)^3-2t^3=6t$, and cubing induces the identity modulo $6$
since both primes split. This gives the full ring. At either prime above $3$ the locally represented
classes are $0,\pm1,\pm2,\pm3\pmod9$. Each is obtained with a unit root,
and varying such a root by $3h$ realizes all changes in $9\mathbb Z_3$.
At a prime away from $3$ every element of the additive span of the residue
cubes is a sum of three residue cubes with a nonzero coordinate. Indeed,
the nonzero cubes form a multiplicative subgroup $H$ of index at most
three. Starting with $C=H\cup\{0\}$, each successive sumset $jC$ gains an
$H$-coset until it fills the additive span. Zero has the representation
$1^3+(-1)^3+0^3$. Hensel's lemma now applies. Thus the exact locally
representable proportion is $(7/9)^2$ by the Chinese remainder theorem.
The upper bound relative to these targets is $O(1/\ell)$.
\end{proof}

There is also a finite-place obstruction that can be repeated independently
at many primes. Write $J$ for the additive span of the integral cubes.

\begin{proposition}\label{nf3:amplify}
Choose distinct primes $\ell_i\equiv1\pmod{12}$ and put
\[
\mathbb K=\mathbb Q(\sqrt6,\sqrt{\ell_1},\ldots,\sqrt{\ell_m}),\qquad r=2^m.
\] Among the elements of $J$ in any bounded Jordan-measurable
positive target region of positive volume under dilation, the proportions
soluble at every finite place in three and four integral cubes tend to
$(17/27)^r$ and $(23/27)^r$, respectively. Five cubes suffice everywhere
locally for every eligible target over every number field.
\end{proposition}
\begin{proof}
The independent square classes give degree $2r$. The $r$ completions at
$3$ all have integer ring $R=\mathbb Z_3[\pi]$ with $\pi^2=6$, since each
$\ell_i$ is a square in $\mathbb Q_3$. Modulo $9R$ the cube residues are
$0,\pm1,\pm3\pi$. Indeed,
\[
 (a+b\pi)^3=a^3+18ab^2+3b(a^2+2b^2)\pi .
\]
The generated group has order $27$. Identify it with $\mathbb Z/9\times
\mathbb Z/3$, with generators $1$ and $3\pi$. A class $(u,v)$ is a sum
of $s$ cube residues exactly when $|u|_9+|v|_3\le s$, where $|\cdot|_h$
is the least absolute residue. For $s=3$ there are $7+5+5=17$ classes,
and for $s=4$ there are $9+7+7=23$. All have a representation with a unit:
when $u\ne0$ the minimal expression already uses a unit, and when $u=0$
one may insert $1,-1$. For a unit $z$, $(z+3h)^3-z^3=9(z^2h+3zh^2+3h^3)$
ranges over $9R$ by Hensel's lemma. Thus the residue condition is sufficient.
The image of $J$ modulo $9\mathcal O_{\mathbb K}$ is the full product of
the generated groups, since a local cube supported at just one prime is obtained
by lifting its root via the Chinese remainder theorem. Equidistribution in
the resulting lattice quotient proves the two formulas. The residue-field
argument of Proposition~\ref{nf3:covol} treats every remaining prime.

Finally, over a local integer ring above $3$, any $\alpha$ in the cube span
satisfies $\alpha\equiv a^3\pmod3$ because Frobenius is additive modulo $3$.
With $t=(\alpha-a^3)/6$ the identity
\[
 \alpha=a^3+(t+1)^3+(t-1)^3+(-t)^3+(-t)^3
\]
supplies a unit coordinate. This proves the five-cube assertion, also a
case of Joly's local bound \cite{nf3:joly}.
\end{proof}

\subsection{Almost all eligible targets and the five-cube threshold}
\begin{theorem}\label{nf3:five}
Let $\mathbb K$ be a fixed number field. Suppose the closure of $\Omega$ is compact in the totally positive archimedean region and avoids zero in each complex coordinate. Then, for some $\delta>0$, all but $O_{\mathbb K,\Omega}(X^{d-\delta})$ targets in $J\cap X\Omega$ are sums of five totally positive integral cubes with roots of height $O_{\mathbb K,\Omega}(X^{1/3})$. The number five is optimal over all number fields.
\end{theorem}
\begin{proof}
Use Lebesgue measure $dy$ on $V=E$ for which $\mathcal O_{\mathbb K}$ has covolume one, and the trace-dual measure $dt$ for which $\mathfrak D^{-1}$ has covolume one. The dual torus
\[
 \mathbb T_{\mathbb K}=V/\mathfrak D^{-1}
\]
has the normalized Haar measure induced by $dt$, where $\mathfrak D$ is the different. The additive character is
\[
 \mathrm e_{\mathbb K}(z)=\exp\bigl(2\pi i\operatorname{Tr}_{\mathbb K/\mathbb Q}(z)\bigr),
\]
with the trace extended to $V$.

Fix a point $\nu_0\in V_+^\times$ and choose $z\in V_+^\times$ with $5z^3=\nu_0$, coordinate by coordinate. Choose $w\in C_c^\infty(V)$ with $0\le w\le1$, positive near $z$, supported in a sufficiently small neighborhood of $z$. Its real coordinates are positive, and the cube map is a diffeomorphism on its support. Define
\begin{align}
 f_P(\beta)&=\sum_{x\in\mathcal O_{\mathbb K}}w(x/P)\mathrm e_{\mathbb K}(\beta x^3),\label{nf3:eq:f}\\
 V_P(\beta)&=\int_V w(y/P)\mathrm e_{\mathbb K}(\beta y^3)\,dy
 =P^d V_1(P^3\beta).\label{nf3:eq:V}
\end{align}
Pushing the smooth density $w(y)\,dy$ forward by $y\mapsto y^3$ gives a smooth compactly supported density. Thus $V_1$ is rapidly decreasing, and
\begin{equation}\label{nf3:eq:singint}
 \mathcal J_w(\nu)=\int_V V_1(t)^5\mathrm e_{\mathbb K}(-t\nu)\,dt
\end{equation}
is the fivefold convolution of this density. It is continuous and positive on a neighborhood $U$ of $\nu_0$. On a sufficiently small compact subneighborhood, $\mathcal J_w\ge c_w>0$.

The number-field version of Hua's lemma gives
\begin{equation}\label{nf3:eq:hua}
 \int_{\mathbb T_{\mathbb K}}|f_P(\beta)|^8\,d\beta
 \ll_{\mathbb K,w,\varepsilon}P^{5d+\varepsilon}.
\end{equation}
This is the case $k=\nu=3$ of Kalinka's estimate
$\int|f|^{2^\nu}\ll T^{d(2^\nu-\nu)+\varepsilon}$
\cite[Lemma 1, p.~151]{nf3:kalinka}. Nonnegative weights bounded by one are handled by orthogonality, which bounds their weighted solution count by the unweighted count in a containing archimedean box.

For the particular moment needed here one can argue directly. Enclose the support in a fixed integral-basis box of scale $P$, and let $F$ be its unweighted cubic sum. For $m\ne0$, factorization of $u^3+v^3=m$ shows that its representation multiplicity in the box is $O(P^\varepsilon)$, since the factors $u+v$ and $u^2-uv+v^2$ are ideal divisors of $(m)$, and bounded generators of each fixed ideal differ by only $O((\log P)^{r_1+r_2-1})$ possible units. The two factors determine $u,v$ up to two choices. For $m=0$, there are $O(P^d)$ pairs. Hence
\[
 I_4:=\int_{\mathbb T_{\mathbb K}}|F|^4\ll P^{2d+\varepsilon}.
\]
Twice differencing and Cauchy's inequality give
\[
 |F(\beta)|^4\ll P^d\sum_b c_b\mathrm e_{\mathbb K}(\beta b),
\]
where $c_b$ counts admissible bounded triples $(x,h_1,h_2)$ with
\[
 b=3h_1h_2(2x+h_1+h_2).
\]
The sum on the right is nonnegative as a sum over $h_1$ of squared absolute values of first-difference sums. The zero coefficient satisfies $c_0\ll P^{2d}$, while $c_b\ll P^\varepsilon$ for $b\ne0$, again by ideal factorization and bounded-unit counting. Multiplying by $|F|^4$, integrating, and using the nonnegative Fourier coefficients $R_4(b)$ of $|F|^4$, we obtain
\[
 I_8\ll P^d\left(P^{2d}I_4+P^\varepsilon\sum_bR_4(b)\right)
 \ll P^{5d+\varepsilon},
\]
because $\sum_bR_4(b)=O(P^{4d})$. This proves \eqref{nf3:eq:hua}.

For $\gamma\in \mathbb K/\mathfrak D^{-1}$, let its denominator ideal be the conductor ideal of the additive character $x\mapsto\mathrm e_{\mathbb K}(\gamma x)$ on $\mathcal O_{\mathbb K}$. Thus
\[
 \mathfrak q_\gamma=\{a\in\mathcal O_{\mathbb K}:a\gamma\in\mathfrak D^{-1}\}.
\]
Fix a sufficiently small $\rho>0$, and let $\mathfrak M$ be the union of the neighborhoods
\begin{equation}\label{nf3:eq:arcs}
 \mathrm N\mathfrak q_\gamma\le C_{\mathbb K} P^\rho,
 \qquad \|\beta-\gamma\|\le C_{\mathbb K} P^{-3+\rho}.
\end{equation}
The norm on $V$ is any fixed coordinate maximum norm, with room to enlarge the fixed constants. Write $\mathfrak m=\mathbb T_{\mathbb K}\setminus\mathfrak M$. The unary cubic Weyl alternative gives
\begin{equation}\label{nf3:eq:weyl}
 \sup_{\beta\in\mathfrak m}|f_P(\beta)|\ll_{\mathbb K,w}P^{d-\kappa}
 \qquad\text{for some }\kappa>0.
\end{equation}
This follows from Skinner's Lemmas 1, 2 and 5 \cite{nf3:skinner} with one polynomial, one $\mathbb K$-variable, degree three, and $\Delta=\rho/(2d)$. The rank-deficient differencing pairs satisfy $h_1h_2=0$ and number $O(P^{d\Delta})$. They cannot account for the alternative lower count $P^{2d\Delta-4\kappa-\varepsilon}$ when $\kappa<d\Delta/4$. The remaining large-sum alternative supplies a denominator element of height $O(P^{2\Delta})$ and approximation error $O(P^{-3+2\Delta})$ after multiplication by that denominator. Taking norms and bounding its inverse places the frequency in \eqref{nf3:eq:arcs}. One may, for example, take $\kappa=\rho/16$, after leaving a margin in $\varepsilon$. The different changes only fixed constants. Smooth weights follow by partial summation from the box estimates.

Combining \eqref{nf3:eq:hua} and \eqref{nf3:eq:weyl}, and choosing $\varepsilon<\kappa$, yields
\begin{equation}\label{nf3:eq:minor10}
 \int_{\mathfrak m}|f_P(\beta)|^{10}\,d\beta
 \le\bigl(\sup_{\mathfrak m}|f_P|\bigr)^2
       \int_{\mathbb T_{\mathbb K}}|f_P|^8
 \ll P^{7d-\kappa}.
\end{equation}

For a denominator character $\gamma$, put
\[
 C(\gamma)=(\mathrm N\mathfrak q_\gamma)^{-1}
  \sum_{x\bmod\mathfrak q_\gamma}\mathrm e_{\mathbb K}(\gamma x^3),
\qquad
 A_{\mathfrak q}(\alpha)=\sum_{\mathfrak q_\gamma=\mathfrak q}C(\gamma)^5\mathrm e_{\mathbb K}(-\gamma\alpha).
\]
The singular series is
\begin{equation}\label{nf3:eq:series}
 \mathfrak S_5(\alpha)=\sum_{\mathfrak q}A_\mathfrak q(\alpha).
\end{equation}
We need uniform absolute convergence and a lower bound on eligible targets.

Let $\mathfrak p\nmid3$, $q=\mathrm N\mathfrak p$, and let $s_h$ be a normalized cubic sum for a primitive additive character modulo $\mathfrak p^h$. At level one,
\[
 |s_1|\le2q^{-1/2}.
\]
Indeed the cube map is bijective if $3\nmid q-1$, while otherwise the sum is a sum of two cubic multiplicative Gauss sums, each of absolute value $q^{1/2}$. Cancellation among units gives $s_2=q^{-1}$ and, for $h\ge3$, the nonunit contribution gives $s_h=q^{-1}s_{h-3}$, using the induced primitive character. Consequently
\begin{equation}\label{nf3:eq:localbounds}
 |s_{3l}|\le q^{-l},\qquad
 |s_{3l+1}|\le2q^{-l-1/2},\qquad
 |s_{3l+2}|\le q^{-l-1}.
\end{equation}
The first inequality is used for $l\ge1$. There are at most $q^h$ primitive characters of conductor $\mathfrak p^h$. Raising \eqref{nf3:eq:localbounds} to the fifth power and summing over these characters bounds $|A_{\mathfrak p^h}(\alpha)|$, uniformly in $\alpha$, by
\begin{equation}\label{nf3:eq:Abounds}
 q^{-2l},\quad32q^{-2l-3/2},\quad q^{-2l-3}
\end{equation}
for the three respective residue classes of $h$.

At a prime above $3$, write $e=v_\mathfrak p(3)$. For sufficiently large $h$, summation over shifts $x\mapsto x+\varpi^{h-e-1}z$, with $z$ ranging over the residue field, cancels the unit terms. The higher Taylor terms have valuation at least $h$, while the linear term gives a nontrivial residue-field character. The nonunits again reduce the conductor by three. Hence
\[
 |s_h|\ll_{\mathbb K} q^{-h/3},\qquad
 |A_{\mathfrak p^h}(\alpha)|\ll_{\mathbb K} q^{-2h/3}.
\]
The finitely many initial conductors are absorbed in the constant. Multiplicativity and \eqref{nf3:eq:Abounds} now give, for example,
\begin{equation}\label{nf3:eq:weightedtail}
 \sum_\mathfrak q |A_\mathfrak q(\alpha)|(\mathrm N\mathfrak q)^{1/4}\ll_{\mathbb K}1
\end{equation}
uniformly in $\alpha$. At good primes the leading weighted terms are $O(q^{-5/4})$, and at the finitely many primes above $3$ the sum is geometric. Thus
\begin{equation}\label{nf3:eq:seriestail}
 \sum_{\mathrm N\mathfrak q>Y}|A_\mathfrak q(\alpha)|\ll_{\mathbb K}Y^{-1/4}.
\end{equation}

Orthogonality identifies the Euler factor of \eqref{nf3:eq:series} with
\[
 \sigma_\mathfrak p(\alpha)
 =\lim_{h\to\infty}q^{-4h}
 \#\{\mathbf x\bmod\mathfrak p^h:x_1^3+\cdots+x_5^3\equiv\alpha\pmod{\mathfrak p^h}\}.
\]
For $\alpha\in J$, Proposition~\ref{nf3:amplify} supplies a local solution with a unit coordinate. Strong Hensel, at a depth depending only on $v_\mathfrak p(3)$, gives a lower bound $\sigma_\mathfrak p(\alpha)\ge c_\mathfrak p>0$ independent of the target. At all sufficiently large good primes,
\[
 \sigma_\mathfrak p(\alpha)=1+O(q^{-3/2})
\]
uniformly, by \eqref{nf3:eq:Abounds}. The convergent product and the finitely many remaining factors therefore give
\begin{equation}\label{nf3:eq:lowerS}
 \mathfrak S_5(\alpha)\ge c_{\mathbb K}>0\qquad(\alpha\in J).
\end{equation}

The arcs in \eqref{nf3:eq:arcs} are disjoint for sufficiently small fixed $\rho$. This follows by applying the product formula to the difference of two distinct denominator fractions. There are $O(P^{2\rho+\varepsilon})$ such fractions, and consequently
\begin{equation}\label{nf3:eq:majormeasure}
 \operatorname{meas}(\mathfrak M)\ll P^{-3d+(d+2)\rho+\varepsilon}.
\end{equation}
Splitting into residue classes gives, on the arc about $\gamma$,
\begin{equation}\label{nf3:eq:majorapprox}
 f_P(\gamma+\beta)=C(\gamma)V_P(\beta)
             +O_{\mathbb K,w}(P^{d-1+2\rho}).
\end{equation}
For example, one may split modulo the rational integer $N=\mathrm N\mathfrak q_\gamma$, since $N\mathcal O_{\mathbb K}\subseteq\mathfrak q_\gamma$. Summation versus integration in each lattice class gives error
$O(NP^{d-1}(1+P^3\|\beta\|))$, which proves \eqref{nf3:eq:majorapprox}. This is also the elementary approximation underlying Skinner's major-arc calculation \cite[proof of Lemma 7]{nf3:skinner}.

Let $g_P$ agree on each major arc with $C(\gamma)^5V_P(\beta)^5$ and vanish elsewhere. Since $f_P$ and $C(\gamma)V_P$ are $O(P^d)$, equations \eqref{nf3:eq:majormeasure}--\eqref{nf3:eq:majorapprox} imply
\begin{equation}\label{nf3:eq:majorL2}
 \int_{\mathfrak M}|f_P^5-g_P|^2
 \ll P^{7d-2+(d+6)\rho+\varepsilon}.
\end{equation}
Choose $\rho<1/(d+6)$ and $\varepsilon$ sufficiently small.

The weighted representation count is
\[
 r_{5,w}(\alpha;P)=\int_{\mathbb T_{\mathbb K}}f_P(\beta)^5\mathrm e_{\mathbb K}(-\alpha\beta)\,d\beta.
\]
The Fourier coefficient of $g_P$ differs from
$P^{2d}\mathfrak S_5(\alpha)\mathcal J_w(\alpha/P^3)$
by $O(P^{2d-\rho/4})$, uniformly in $\alpha$. The denominator tail follows from \eqref{nf3:eq:seriestail}, and the archimedean tail is smaller than any prescribed power of $P^{-\rho}$ because $V_1$ is rapidly decreasing.

Bessel's inequality, \eqref{nf3:eq:minor10}, and \eqref{nf3:eq:majorL2} now give, for a bounded target region $U$,
\begin{equation}\label{nf3:eq:meansquare}
 \sum_{\alpha\in\mathcal O_{\mathbb K}\cap P^3U}
 \left|r_{5,w}(\alpha;P)-P^{2d}\mathfrak S_5(\alpha)
                     \mathcal J_w(\alpha/P^3)\right|^2
 \ll_{\mathbb K,w,U}P^{7d-\eta}
\end{equation}
for some $\eta>0$. In detail, the minor-arc contribution is $O(P^{7d-\kappa})$, the major approximation has the bound in \eqref{nf3:eq:majorL2}, and the squared pointwise tail summed over $O(P^{3d})$ targets is $O(P^{7d-\rho/2})$.

For $\alpha\in J\cap P^3U$ in the compact neighborhood where $\mathcal J_w\ge c_w$, the main term is at least $c_{\mathbb K}c_wP^{2d}$. If $r_{5,w}(\alpha;P)=0$, its contribution to the left side of \eqref{nf3:eq:meansquare} is therefore at least a constant times $P^{4d}$. The exceptional set has size
\[
 O(P^{3d-\eta}).
\]
Each represented target has five roots of height $O_w(P)$ in the positive region. A finite cover and $X=P^3$ give the exceptional-set estimate. Proposition~\ref{nf3:amplify} proves optimality already in $\mathbb Q(\sqrt6)$.
\end{proof}

\subsection{Four cubes on the locally representable set}\label{nf3:foursection}
For four summands the finite-place restrictions of
Proposition~\ref{nf3:amplify} account for all but a power-saving exceptional
set. We use a fourth ordinary cubic sum to supplement the prime-restricted
triples of Theorem~\ref{nf3:davenport}.

\begin{theorem}\label{nf3:four}
Fix a number field $\mathbb K$ and a compact subset $U$ of the totally
positive archimedean region avoiding zero in every complex coordinate.
There is $\delta>0$ such that, among algebraic integers in $A U$ which are
representable by four integral cubes at every finite completion, all but
$O_{\mathbb K,U}(A^{d-\delta})$ are sums of four totally positive integral
cubes with roots of height $O_{\mathbb K,U}(A^{1/3})$.
In particular, the limiting proportion of four-cube sums among eligible
targets in the regions of Proposition~\ref{nf3:amplify} is $(23/27)^r$.
This last assertion extends to bounded Jordan-measurable positive regions
by exhaustion away from the coordinate walls.
\end{theorem}

\begin{proof}
Set $P=A^{1/3}$, $X=P^d$. Use the trace character, Haar measures and
denominator ideals of Theorem~\ref{nf3:five}. For
$t\in \mathbb K\otimes_{\mathbb Q}\mathbb R$ put
$B(t)=\prod_{v\mid\infty}(1+|t_v|)^{d_v}$, with $d_v=1$ or $2$.
Choose a small archimedean box for the roots, as in
Theorem~\ref{nf3:davenport}, with four cubes summing to a neighborhood
of a given target. Take a nonnegative product weight $w$ of compact
support in that box, smooth and positive near the chosen root. Write
\[
 f_a(\alpha)=\sum_{z\in\mathcal O_{\mathbb K}}w(az/P)
                \mathrm e_{\mathbb K}(\alpha (az)^3),\quad f=f_1,
 \qquad I(t)=\int w(y)\mathrm e_{\mathbb K}(ty^3)\,dy .
\]
The function $I$ is rapidly decreasing at each infinite place since
the cube map has invertible derivative on the support. The Fourier
transform of $I^4$ with the negative sign, written $J_4$, is the
fourfold convolution of the pushforward of $w$. Thus $J_4\ge c>0$ on
a neighborhood of the target, after shrinking the latter.

Choose $p,q$ exactly as in Theorem~\ref{nf3:davenport}, now with
norm scales $U_1=X^{1/6}$ and $U_2=X^{1/18}$. Put
\[
 H(\alpha)=\sum_{p,q}(f(\alpha)-f_p(\alpha))f_p(\alpha)f_{pq}(\alpha),
 \qquad
 W=\sum_{p,q}\frac{1}{\mathrm N(p)^2\mathrm N(q)}
 \asymp \frac{X^{-1/6}}{(\log X)^2}.
\]
Its nonnegative weighted sixth-variable collision count satisfies
\begin{equation}\label{nf3:fourenergy}
                   \int |H|^2\ll X^{55/18+\varepsilon}.
\end{equation}
Here the subtraction in the first factor restricts that root to the units
modulo $p$. The nonzero residue classes modulo $p$ split into at
most three sets on each of which cubing is injective. Cauchy's inequality
over these sets, together with the Cauchy step over $p,q$ in
Theorem~\ref{nf3:davenport}, costs at most the additional factor three.
The proof of that theorem with $N=X^3$ yields \eqref{nf3:fourenergy} for this $H$ with
the given bounded weights.

We require a cubic Weyl bound on the complement of the arcs
\[
 \alpha=\gamma+\beta,\qquad
 Q B(P^3\beta)\le X^\omega,\qquad Q=\mathrm N\mathfrak q_\gamma,
 \qquad 7/9<\omega<5/6.
\]
Different lifted arcs are disjoint, by the norm lower bound on the
distance of their distinct centers and $2\omega<3$.
We claim that on this complement
\begin{equation}\label{nf3:fourweyl}
                       |f(\alpha)|^4\ll X^{4-\omega+\varepsilon}.
\end{equation}
By partial summation it suffices to treat integral-basis boxes of
side $O(P)$, uniformly in the endpoints. Twice differencing and grouping
nonzero products of the two differences gives
\[
 |f|^4\ll X^{3+\varepsilon}
   + X^{1+\varepsilon}\!
   \sum_{0<H_\infty(t)\ll P^2}
   \prod_{j=1}^d\min(P,\|\operatorname{Tr}(\alpha t e_j)\|^{-1}),
\]
where the $e_j$ form an integral basis, $H_\infty(t)$ is embedding height,
and the coefficient $6$ in the second difference has been absorbed into
$t$. The multiplicity of a nonzero product is $O(X^\varepsilon)$ by the
divisor and unit estimate in Theorem~\ref{nf3:davenport}. Partition the
dual torus into cells of side comparable to $1/P$ in the displayed
coordinates. The largest occupancy of a cell is at most
\[
 M=\#\{(t,c)\in\mathcal O_{\mathbb K}\times\mathfrak D^{-1}:
                H_\infty(t)\le CP^2,\ H_\infty(\alpha t-c)\le C/P\},
\]
by taking differences within a cell. The sum of the suprema of the
displayed products over all cells is $O(X(\log(2P))^d)$.
If the pairs counted by $M$ contain two vectors independent over
$\mathbb K$, then $M\ll X$. To see this, send $(t,c)$ to
$(P^{-2}t,P(\alpha t-c))$. The resulting real lattice has covolume
$\asymp X^{-1}$. Multiplying the independent pairs by the integral basis
produces $2d$ independent bounded real vectors in this lattice, by
extension of scalars. Its covering radius is bounded, so disjoint
Voronoi cells about the counted points lie in a bounded enlargement
and each has volume $\asymp X^{-1}$. If all nonzero pairs lie on one
$\mathbb K$-line, they satisfy $c=\gamma t$ for a fixed $\gamma$,
since $t=0$ forces $c=0$ at large $P$. Thus $t\in\mathfrak q_\gamma$ and
\[
 |t_v|\le C'P^2/(1+P^3|\alpha_v-\gamma_v|).
\]
A box of place-radii $T_v$ contains
$O_{\mathbb K}(1+\prod_v T_v^{d_v}/\mathrm N\mathfrak q)$ points of
$\mathfrak q$, since after scaling the box to unit size, the nonzero differences
have maximum norm $\gg(\mathrm N\mathfrak q/\prod T_v^{d_v})^{1/d}$.
Consequently $M\ll 1+X^2/(Q B(P^3(\alpha-\gamma)))$. Off the arcs this
is $O(X^{2-\omega})$, proving \eqref{nf3:fourweyl}. In combination with
\eqref{nf3:fourenergy} it gives, for some $\kappa>0$,
\begin{equation}\label{nf3:fourminor}
            \int_{\mathfrak m}|fH|^2
             \ll X^{2-\omega/2+55/18+\varepsilon}
             \ll X^5W^2X^{-\kappa}.
\end{equation}
The strict inequality needed here is $\omega>7/9$.

Write $C(\gamma)$ for the normalized complete cubic sum at $\gamma$.
On an arc, for $a=1,p$, Poisson summation yields, for every fixed $T>0$,
\begin{equation}\label{nf3:fourpoisson}
 f_a(\gamma+\beta)
   =\frac{X}{|\mathrm N a|}C(a^3\gamma) I(P^3\beta)+O_T(X^{-T}).
\end{equation}
Indeed, scaling each residue lattice by $a/P$ gives dual frequencies
$P\xi/a$ for $\xi\in\mathfrak q_{a^3\gamma}^{-1}\mathfrak D^{-1}$
in trace-pairing coordinates. The nonzero frequencies, relative to
radii $1+P^3|\beta_v|$, have maximum norm
$\gg(X/(|\mathrm Na|Q B(P^3\beta)))^{1/d}$ by the norm formula.
This is a positive power of $X$ since $\omega+1/6<1$.
Integration by parts at the place of maximum and the preceding packing
bound give the stated error, with arbitrarily many derivatives. On the
smaller arcs $QB(P^3\beta)\le D_0=X^\xi$ the same formula holds for $a=pq$
provided $\xi>0$ is sufficiently small. Choose also $\xi<1/18$, so that
all labels are units modulo every conductor on these smaller arcs.

On the larger arcs, bound $f_{pq}\ll X/\mathrm N(pq)$. Use
\eqref{nf3:fourpoisson} twice on the ordinary factors and once on the
$p$-divisible factor, allowing the subtraction $f_p$ in $f-f_p$ its own
normalized factor $1/\mathrm N(p)$. Cauchy's inequality in the
probability weights proportional to $1/(\mathrm N(p)^2\mathrm N(q))$
then bounds the squared integral on $QB(P^3\beta)>D_0$ by
$\ll X^5W^2D_0^{-\tau}$, for some $\tau>0$. Here are the details of
convergence. Writing $t=P^3\beta$, each normalized squared summand is
bounded, apart from a rapidly decreasing factor in $t$, by a constant
times
\[
 |C(\gamma)|^4 |C(p^3\gamma)|^2
 +\mathrm N(p)^{-2}|C(\gamma)|^2 |C(p^3\gamma)|^4 .
\]
The sum of these bounds times $Q^\tau$ over all denominator characters
is uniform in $p$, for example for $\tau=1/12$. At places other than $(p)$
the local series uses six powers of the same complete sum and converges
absolutely with weight $Q^\tau$ by \eqref{nf3:eq:localbounds}, also at
the finitely many places above $3$ by their conductor-decay estimate.
At $(p)$ with residue cardinality $u$, the first summand costs
$O(u^{-1+3\tau})$ at positive conductor, since the first three levels
have respectively $O(u^{-1}),O(u^{-2}),O(u^{-1})$ before weighting.
For the second summand the first three local contributions prior to
multiplication by $u^{-2}$ are $O(1),O(1),O(u)$. For each summand at
levels above three, the estimates decrease by a factor $u^{-3+3\tau}$
every three levels once both factors have positive conductor.
The rapid decrease in $t$ absorbs $B(t)^\tau$. The change of variables
$d\beta=X^{-3}dt$ proves the asserted squared integral bound.

On the smaller arcs the complete factors are all $C(\gamma)$, and the
multiplier of $X^4 C(\gamma)^4I(P^3\beta)^4$ in $fH$ is
\[
          W'=\sum_{p,q}\frac{1-1/\mathrm N(p)}
                              {\mathrm N(p)^2\mathrm N(q)}\asymp W .
\]
The error from \eqref{nf3:fourpoisson} is negligible. We next check the
singular series $\mathfrak S_4(n)=\sum_{\mathfrak a} A_{\mathfrak a}(n)$,
with exponent four in the definition of the conductor coefficient.

For all nonzero $n$ of height $O(P^3)$ we have, with $\tau=1/12$,
\begin{equation}\label{nf3:fourser}
  \sum_{\mathfrak a}|A_{\mathfrak a}(n)|(\mathrm N\mathfrak a)^\tau
       \ll_\varepsilon X^\varepsilon .
\end{equation}
At $\mathfrak l\nmid 3n$ all solutions are primitive, so all coefficients
of conductor exponent at least two vanish by smooth lifting. At the
first level, of norm $u$, the nontrivial cubic complete sum vanishes
if cubing is bijective. Otherwise it is a sum of two cubic character
Gauss sums times their character twists. In its fourth power every
term has coefficient of size $u^2$, and its sum over the unit character
parameter with the additional additive twist at $n$ is $O(\sqrt u)$
by the Gauss bound (and is $-1$ for a trivial multiplicative twist).
Thus $A_{\mathfrak l}(n)\ll u^{-3/2}$. For $\mathfrak l\mid n$ away
from $3$, the three levels $h=3j,3j+1,3j+2$ contribute in absolute
value at most $O(u^{-j}),O(u^{-j-1}),O(u^{-j-2})$ respectively by
\eqref{nf3:eq:localbounds}. Their positive levels have bounded
weighted sum. The product of these bounds over the prime divisors of
$(n)$ is $O_\varepsilon(X^\varepsilon)$ by the divisor bound. At the
primes above $3$ the bound is $O_{\mathbb K}(u^{-h/3})$, so these weighted
factors are bounded as well. Multiplicativity proves \eqref{nf3:fourser}.

The product $\mathfrak S_4(n)$ has a uniform positive lower bound on
everywhere locally representable targets. For $\mathfrak l\nmid3n$ the
factor is $1+O(u^{-3/2})$. For $\mathfrak l\mid n$ away from $3$
the number of zero classes of four cubes in the residue field is at
least $u^3$, since the nontrivial character sums in its orthogonality formula
are real, since the cube is odd, and are raised to the fourth power.
All except the zero vector are smooth, giving a lower bound $1-u^{-3}$
for the local density. At each of the finitely many remaining small
places the compact set of locally represented targets is the image
of the compact ring to the fourth power. At any such target some
representing coordinate is nonzero, since this holds automatically for a
nonzero target and one may use $1,-1,0,0$ at zero. Hensel's lemma
on small balls about this solution gives a positive lower density
throughout an open neighborhood of the target. Compactness gives a
uniform positive constant. The product of the resulting bounds converges.

By \eqref{nf3:fourser}, the Fourier coefficient of the small-arc model is
\[
 XW'\mathfrak S_4(n)J_4(n/P^3)+O_\varepsilon(XW'X^\varepsilon D_0^{-\tau}).
\] Indeed, both the conductor
tail and the archimedean tail are bounded by inserting
$(QB(t)/D_0)^\tau$ and using the rapid decrease of $I$.
Bessel's inequality, \eqref{nf3:fourminor} and the larger-arc tail
now show that the squared deviations of the actual representation
counts from this main term, summed over the $O(X^3)$ targets in the
fixed neighborhood, total $O(X^5W^2X^{-\kappa'})$ for some $\kappa'>0$.
Every locally represented exceptional target contributes $\gg (XW)^2$.
Thus there are $O(X^{3-\kappa'})$ such exceptions. A finite cover
and $X=A^{d/3}$ prove the result.
\end{proof}

\subsection{A quadratic exponential sum and three cubes}\label{nf3:vaughansec}
A single prime restriction and a quadratic exponential sum give a
stronger three-cube estimate. The rational exponent was obtained by
Vaughan \cite{nf3:vaughan85} and is discussed in
\cite[Sect.~3]{nf3:maynard}.

\begin{theorem}\label{nf3:vaughan}
Under the hypotheses of Theorem~\ref{nf3:davenport}, for every
$\varepsilon>0$ one has
\[
                  R_\Omega(A)\gg_{\mathbb K,\Omega,\varepsilon}
                                 (A^d)^{8/9-\varepsilon}.
\]
All conjugates of the three roots may be required to have absolute value
$\asymp A^{1/3}$.
\end{theorem}
\begin{proof}
Put $P=A^{1/3}$, $X=P^d$. Take the archimedean box of
Theorem~\ref{nf3:davenport} and the same balanced prime generators $p$
with norm $\asymp L=X^{1/6}$. Count the triples of roots $a,pb,pc$ in
the dilated box with $p\nmid a$. Their total with prime labels is
$\gg X^3/(L\log X)$. Per fixed prime we claim the number of collisions
is $O(X^{8/3+\varepsilon})$. Cauchy's inequality over labels will then
bound the total second moment by $O(L^2X^{8/3+\varepsilon})$,
and a second application proves the theorem.

Split the unit classes of $a$ modulo $p$ into at most three sets on
which cubing is injective. In a collision within one set the difference
of the two ordinary roots is $p^3u$. The collisions with $u=0$ are
$O(X^{1+\varepsilon}Y^2)$, where $Y=X/L=X^{5/6}$, by the nonzero
two-cube divisor bound. For $u\ne0$ the equation becomes
\[
 u(3x^2-3p^3ux+p^6u^2)=b_1^3+b_2^3-b_3^3-b_4^3
\]
up to choice of signs. The variable $x$ lies in an integral-basis
box of sidelength $O(P)$ and the nonzero $u$ have embedding height
$O(H^{1/d})$, with $H=X/L^3=X^{1/2}$. The $b_j$ lie in a box of volume
scale $Y$ separated at every place from zero on its natural scale.
Choose a nonnegative smooth product majorant for this box after scaling
by $P_0=Y^{1/d}$. Its support, and finitely many derivative bounds of
each order, can be taken uniform in the balanced $p$. Let $F(\alpha)$
be the associated cubic sum in $b$, and $G(\alpha)$ the double sum of
the left-hand phase over $u,x$ in the indicated ranges. Orthogonality
bounds the remaining collisions by $\int |G||F|^4$. Also
$\int |F|^4\ll Y^{2+\varepsilon}$ by the two-cube divisor bound, since a
zero sum, if included by a majorant, contributes $O(Y^2)$.

For $Q=\mathrm N\mathfrak q_\gamma$ set
$D=Q B(P_0^3(\alpha-\gamma))$, with $B(t)$ as in
Section~\ref{nf3:foursection}. Consider arcs $D\le X$ at every
denominator. They are disjoint since $X^2=o(Y^3)$. We have
\begin{align*}
 |G(\alpha)|&\ll X^{1+\varepsilon} &&\text{off these arcs},\\
 |G(\alpha)|&\ll X^{3/2+\varepsilon}D^{-1/2} &&\text{on them}.
\end{align*}
For the proof use Cauchy's inequality in $u$ and difference once in
$x$. The zero difference contributes $O(H^2X)$ to $|G|^2$, while the
remaining products of $u$ and the shift in $x$ have height
$O((HX)^{1/d})$ and multiplicity $O(X^\varepsilon)$. The residual
linear phase has frequency $6\alpha$ times this product. The cell
argument of \eqref{nf3:fourweyl} therefore bounds $|G|^2$ by
$O(X^\varepsilon(H^2X+HXM))$, where $M$ counts pairs
$(v,c)\in\mathcal O_{\mathbb K}\times\mathfrak D^{-1}$ with
$H_\infty(v)\ll(HX)^{1/d}$ and $H_\infty(\alpha v-c)\ll P^{-1}$.
The rank-two case gives $M\ll H$ by the covolume $(HX)^{-1}X$.
In the rank-one case $c=\gamma v$ and packing gives
$M\ll 1+HX/(Q B(P_0^3(\alpha-\gamma)))$, since
$(HX)^{1/d}P=P_0^3$. These estimates give the asserted inequalities, since in the rank-one
case a different lifted center gives $D>X$ by disjointness.
The off-arc integral is consequently $O(X^{1+\varepsilon}Y^2)$.

Write $t=P_0^3(\alpha-\gamma)$ on an arc. For the denominator ideal
$\mathfrak a=\mathfrak q_\gamma$, let $s(\mathfrak a)$ be the norm of
the product of its full prime-power factors with exponent at least two.
Poisson summation in the residue classes gives a zero mode
$F_0=Y C(\gamma)I(t)$ and a remainder satisfying
\begin{equation}\label{nf3:vairest}
 |F-F_0|\ll X^\varepsilon Q^{1/2}s(\mathfrak a)^{1/6} B(t)^{1/2}.
\end{equation}
We include the estimates for this remainder.
Every complete cubic sum with a linear twist, divided by $Q$, is
$\ll_\varepsilon Q^{-1/2+\varepsilon}s(\mathfrak a)^{1/6}$.
Indeed at a prime of exponent one and residue characteristic other
than $3$ the finite-field Weil bound is $2\sqrt u$, for residue norm
$u$. At exponent $j\ge2$ the normalized bound is
$O(j u^{-j/3})$ uniformly in the linear term away from the primes
over $6$, with an additional fixed constant at the latter places.
To verify it, write the local phase as $a z^3+bz$ with $a$ a unit.
On unit $z$, shifts by the uniformizer to exponent $\lceil j/2\rceil$
give cancellation unless $3az^2+b$ vanishes modulo the prime to
exponent $\lfloor j/2\rfloor$. Away from $6$ there are at most two
unit roots at this depth. At the finitely many other places their
Haar measure is $O_{\mathbb K}(u^{-\lfloor j/2\rfloor})$, by Hensel's
lemma at depth greater than $2v(6)$ applied with derivative $6az$,
and a fixed enlargement for smaller depths. The contribution of the
units is thus $O(u^{-j/3})$ for $j\ge2$. For nonunits write $z=\varpi y$.
When $b\in\varpi^2\mathcal O$ this gives $u^{-1}$ times the same cubic
sum at exponent $j-3$, interpreting levels at most zero trivially.
Otherwise the scaled phase after removal of its common power of
$\varpi$ has unit linear coefficient and cubic coefficient divisible
by $\varpi$, so its sum vanishes by character orthogonality above that common
depth. The remaining small levels satisfy the claimed bound directly.
Iteration, the level-one estimate, and multiplicativity prove the
complete sum bound.

The archimedean integral in a mode $\ell=P_0\xi$ with linear twist is
\[
 \ll_T B(t)^{-1/2}
        \prod_{v\mid\infty}(1+|\ell_v|/(1+|t_v|))^{-T}
\]
for every fixed $T$. On each real interval the cubic phase has
second derivative bounded away from zero after normalization by $t_v$.
At a complex place the Hessian of its real trace is nonsingular on the
annulus, giving the corresponding two-dimensional stationary-phase
bound, uniformly in the linear twist. Integration by parts when the
linear twist dominates gives the indicated rapid decrease.
The nonzero $\xi\in\mathfrak a^{-1}\mathfrak D^{-1}$ satisfy
$\max_v |P_0\xi_v|/(1+|t_v|)\gg(Y/(QB(t)))^{1/d}$.
Place-radii packing as in Section~\ref{nf3:foursection} then bounds
the sum of the displayed decreasing products over nonzero modes by
$O(QB(t)/Y)$, including when $QB(t)<Y$ by increasing $T$.
This proves \eqref{nf3:vairest}.

It remains to integrate. The zero mode $F_0$ has
$|C(\gamma)|\ll_\varepsilon Q^{-1/3+\varepsilon}$ and $I(t)$
rapidly decreasing. There are at most $Q$ characters for one conductor.
Changing variables $d\alpha=Y^{-3}dt$ gives
\[
 \int_{D\le X}|G|\,|F_0|^4\ll
 X^{3/2+\varepsilon}Y\sum_{\mathrm N\mathfrak a\le X}
                       (\mathrm N\mathfrak a)^{-5/6}
 \ll X^{5/2+\varepsilon}.
\]
For the remainder integrate $B(t)^{3/2}$ over $B(t)\le X/Q$,
obtaining $O_\varepsilon((X/Q)^{5/2+\varepsilon})$ by dyadic
decomposition at the infinite places. Equation~\eqref{nf3:vairest}
then gives
\[
 \int_{D\le X}|G|\,|F-F_0|^4
 \ll X^{4+\varepsilon}Y^{-3}
       \sum_{\mathrm N\mathfrak a\le X}s(\mathfrak a)^{2/3}
 \ll X^{8/3+\varepsilon}.
\]
For the last sum, write the squarefull part uniquely as
$\mathfrak b^2\mathfrak c^3$ with $\mathfrak c$ squarefree. The number
of squarefull ideals of norm at most $T$ is $O_{\mathbb K}(T^{1/2})$ by
ordinary ideal counting and $\sum (\mathrm N\mathfrak c)^{-3/2}<\infty$.
The number of complementary factors at fixed squarefull norm $s$ is
$O(X/s)$, and partial summation gives the bound $O(X^{7/6})$.
This completes the per-prime collision estimate.
\end{proof}

The same estimates give many four-cube representations at nearly the
power predicted by the circle method.

\begin{corollary}\label{nf3:manyfour}
Let $U$ and $\mathbb K$ be as in Theorem~\ref{nf3:four}, and fix
$\varepsilon>0$. For some $\delta>0$, every locally representable target,
with $O_{\mathbb K,U,\varepsilon}(A^{d-\delta})$ exceptions, has at least
a constant multiple of $(A^d)^{1/3-\varepsilon}$ representations by four
totally positive cubes with roots of height $O_{\mathbb K,U}(A^{1/3})$.
\end{corollary}
\begin{proof}
Retain the notation $X=P^d$ and choose the primes in a smaller fixed
range of norms $\asymp L=X^\ell$, with $0<\ell<1/6$. Now $Y=X/L$
and $H=X/L^3$. The estimates in the preceding proof, still on arcs
$Q B(P_0^3\beta)\le X$, bound the collisions at one prime by
\[
 X^\varepsilon\left(
 XY^2+HX^{1/2}Y^2+HXY X^{1/6}
            +HX Y^{-3}X^{5/2}X^{7/6}\right)
 \ll X^{7/2+\varepsilon}/L^5 .
\]
The four terms are respectively the zero difference, the complement
of the arcs, the zero dual mode, and the nonzero modes. Hence
$H_3=\sum_p(f-f_p)f_p^2$ has squared integral
$\ll X^{7/2+\varepsilon}/L^3$ by Cauchy's inequality over primes, and the
normalizing mass is $W' = \sum_p (1-1/\mathrm N(p))/\mathrm N(p)^2
\asymp 1/(L\log X)$.

In the proof of Theorem~\ref{nf3:four} apply the minor-arc bound
\eqref{nf3:fourweyl} with $1-2\ell<\omega<1-\ell$. This gives
\[
 \int_{\mathfrak m}|fH_3|^2
 \ll X^{11/2-\omega/2+\varepsilon}L^{-3}
 \ll X^5(W')^2X^{-\eta}
\]
for some $\eta>0$. On the larger arcs the full sums and $f_p$ have
the Poisson approximation \eqref{nf3:fourpoisson} since
$L X^\omega=o(X)$. Bound one of the two factors $f_p$ trivially there to leave exactly the two convergent weighted local
products in the proof of Theorem~\ref{nf3:four}. On the smaller arcs
with cutoff $X^\xi$, $0<\xi<\ell$, all four factors give the principal
term $XW'\mathfrak S_4(n)J_4(n/P^3)$ for the Fourier coefficient,
with the same power-saving error. The squared-deviation argument
therefore gives $\gg XW'$ weighted representations away from
$O(X^{3-\eta'})$ exceptions. A given quadruple admits $O_\ell(1)$
prime labels, since a selected nonzero root has norm $O(X)$ and
every label has norm $\gg X^\ell$. The weights are bounded.
Choosing $\ell$ sufficiently small in terms of $\varepsilon$ proves
the assertion.
\end{proof}

\subsection{Averaging the prime restriction}\label{nf3:averaged}
Averaging the quadratic difference over its modulus gives the
number-field counterpart of Vaughan's exponent $19/21$ for three
positive cubes \cite{nf3:vaughan86}.
\begin{theorem}\label{nf3:nineteen}
Under the hypotheses and root-size restrictions of
Theorem~\ref{nf3:davenport} one has
$R_\Omega(A)\gg_{\mathbb K,\Omega,\epsilon}(A^d)^{19/21-\epsilon}$.
\end{theorem}
\begin{proof}
Use the prime-labelled triples from Theorem~\ref{nf3:vaughan}, now
with norm scale $L=X^{1/7}$, $Y=X/L$, $H=X/L^3$ and actual side
scales the $d$th roots of these volumes. The first moment is
$\gg X^3/(L\log X)$. We will show that the collisions summed over all
single labels are $O_\epsilon(X^{7/2+\epsilon}L^{-9/2})$. After Cauchy's
inequality over labels the second moment is $O(X^{3+\epsilon})$.
The lower bound follows since $X^{6}/(L^2X^3)=X^{19/7}$.

Separate the unit cube fibres and equal ordinary roots as in
Theorem~\ref{nf3:vaughan}. The equal-root contribution summed over
labels is $O(LXY^{2+\epsilon})$. In the remaining count put
$h=(x-x')/p^3$ and $z=x+x'$. After discarding the parity constraint
one has
\[
             3hz^2+h^3p^6=4(y_1^3+y_2^3-y_3^3-y_4^3).
\]
All these quantities range over the boxes already used in that theorem,
with $h\ne0$. Enlarge the $p$-sum to all algebraic integers in its
height box; the Diophantine count has nonnegative weights. One can take
independent, smooth, nonnegative product majorants $u_h$ for the shifts,
and smooth weights in $z,p,y$. The weights for $y$ can be supported away
from zero place by place, by balance of the primes; a finite collection
of compact annular sectors suffices. Orthogonality gives a majorant
with integrand
\[
    \sum_{h\ne0} u_h\, B_z(3h\alpha) C_p(h^3\alpha)|F(4\alpha)|^4 ,
\]
where $B_z$ is quadratic at volume $X$, $C_p$ is a sixth-power sum at
volume $L$, and $F$ is cubic at volume $Y$.

Use arcs $D=Q B(Y^{3/d}(\alpha-\gamma))\le X$. The cell count used
in the quadratic estimate of Theorem~\ref{nf3:vaughan}, before the
outer Cauchy factor, gives
\[
 \sum_{h\ne0}u_h |B_z(3h\alpha)|^2
 \ll X^\epsilon
 \begin{cases}
 HX,&\alpha\text{ outside the arcs},\\
 H X+H X^2/D,&\alpha\text{ on the arcs}.
 \end{cases}
\]
Indeed, the zero difference costs $HX$ and the other terms involve
the product of $h$ and one difference, with the same divisor bound
and rank-one/rank-two count. For a smooth $z$-weight, partial
summation of the inner linear sum after differencing gives the same
product of geometric-sum bounds, uniformly in the shift. Cauchy's
inequality in $h$ on the arcs yields $O(H X^{1+\epsilon}/\sqrt D)$
for the sum of the absolute values. Bounding $C_p$ by $O(L)$ and
using the zero and nonzero Poisson terms of \eqref{nf3:vairest}
bounds the major integral by
\[
      O_\epsilon\bigl(X^\epsilon L(HXYX^{1/6}
                      + HXY^{-3}X^{5/2}X^{7/6})\bigr)
       \ll X^{7/2+\epsilon}L^{-9/2}.
\]
The fixed multiplier $4$ in $F(4\alpha)$ affects its denominator only
at the primes over $2$ by bounded factors, and its archimedean scale
by fixed factors.

On the complement apply Cauchy's inequality in $(\alpha,h)$.
The quadratic factor squared against $|F|^4$ is
$O(X^{1+\epsilon} H Y^2)$ by the preceding bound and the fourth moment.
For the other factor squared against $|F|^4$, expand $|C_p|^2$.
Equal sixth powers give $O_\epsilon(H L Y^{2+\epsilon})$. There are
$O(L^2)$ other pairs of sixth powers, of nonzero difference $b$. Each
pair is bounded by
\[
       \int |E(b\alpha)|\,|F(4\alpha)|^4\,d\alpha ,
 \qquad E(\alpha)=\sum_{h\ne0}u_h\mathrm e_{\mathbb K}(\alpha h^3).
\]
Here we first enlarged the squared factor to the full torus, on which
its unexpanded integrand is nonnegative. Put $V=H$ and $P_1=V^{1/d}$.
On the complement of $Q B(P_1^3\beta)\le V$, the argument of
\eqref{nf3:fourweyl} gives $|E|\ll V^{3/4+\epsilon}$, also at
$\omega=1$. On this major set we need
\begin{equation}\label{nf3:arc4}
                      \int |E|^4\ll V^{1+\epsilon}.
\end{equation}
For completeness, the normalized complete cubic sum with any linear
twist is bounded by a multiplicative weight $\kappa(\mathfrak q)$.
At a prime away from $6$ of norm $q$ one may take
$2q^{-j},2q^{-j-1/2},2q^{-j-1}$ at levels $3j,3j+1,3j+2$
respectively (using the first formula only for $j>0$). In the
proof of \eqref{nf3:vairest} the unit part with unit linear coefficient
has at most two simple stationary residues. At an odd level $2m+1$,
first sum over shifts by $\varpi^{m+1}$, then sum the additional
digit at each stationary residue by the quadratic Gauss bound, giving
the normalized bound $2q^{-m-1/2}$. At even levels summing over
shifts by $\varpi^m$ gives $2q^{-m}$. A nonunit linear coefficient
forces the stationary root into $\varpi\mathcal O$, with the earlier
three-step recursion. At places above $6$ the bound
$O_{\mathbb K}((1+j)q^{-j/3})$ suffices. These estimates imply
$\sum_{Q\le V}Q\kappa(\mathfrak q)^4\ll_\epsilon V^\epsilon$ by an
Euler product after replacing the factor $Q$ by $Q^{1-\epsilon}$.
The archimedean cubic integral with a linear twist is bounded place
by place by $(1+|t_v|)^{-d_v/3}$, times arbitrary inverse powers
outside the stationary range. At a real place use the third derivative.
At a complex place the inner disc $|z|\le |t_v|^{-1/3}$ has area
$O(|t_v|^{-2/3})$, and dyadic annuli $|z|\asymp r$ give
$O(|t_v|^{-1}r^{-1})$ by rescaling and the nonsingular Hessian.
Their sum has the same bound. Poisson summation and packing as in
\eqref{nf3:vairest} consequently give
$O(V\kappa(\mathfrak q)B(t)^{-1/3})$ on these arcs, since
$Q B(t)/V\le1$. Raising to the fourth power, integrating after
$d\beta=V^{-3}dt$ and using $\int B(t)^{-4/3}dt<\infty$ proves
\eqref{nf3:arc4}. Omitting $h=0$ changes the sum by $O(1)$.

Multiplication by $b$ preserves Haar measure of the torus, since it
is a surjective endomorphism with finite kernel. H\"older's inequality
with exponents $4,2,4$, using \eqref{nf3:arc4}, the fourth moment of
$F$, and the eighth moment \eqref{nf3:eq:hua}, now bounds each nondiagonal
pair by $O(X^\epsilon(H^{3/4}Y^2+H^{1/4}Y^{9/4}))$.
The entire complementary integral is thus
\[
 \ll X^\epsilon (HXY^2)^{1/2}
       (HLY^2+L^2H^{3/4}Y^2+L^2H^{1/4}Y^{9/4})^{1/2}
 \ll X^{7/2+\epsilon}L^{-9/2}
\]
at $L=X^{1/7}$. This proves the claim, with independently chosen
positive boxes for all the original roots.
\end{proof}

\subsection{A squarefree smooth sixth moment}\label{nf3:smoothsec}
The exponent in Vaughan's smooth cubic mean
\cite{nf3:vaughan89} also holds for squarefree algebraic integers
over every number field. We give an argument using the quadratic
difference and an average over divisors. Smoothness below refers to
norms of prime ideals. Put $X=P^d$. A product of fixed compact
annuli in the archimedean coordinates, dilated by $P$, will be
called $P\mathcal B$; the real coordinates may be restricted to
either sign. For $S$ the squarefree $R$-smooth integers in that
region write $f_S(\alpha)=\sum_{x\in S}\mathrm e_{\mathbb K}(\alpha x^3)$.
\begin{theorem}\label{nf3:smoothsix}
For every $\epsilon>0$ there is $\eta>0$ such that for
$1\le R\le X^\eta$ one has
\[
             \int |f_S(\alpha)|^6\,d\alpha
                 \ll_{\mathbb K,\mathcal B,\epsilon}X^{13/4+\epsilon}.
\]
For every fixed $0<\theta\le\eta$, taking $R=X^\theta$ gives the
conclusion of Theorem~\ref{nf3:nineteen} with exponent
$11/12-\epsilon$ and pairwise coprime, squarefree $R$-smooth roots.
\end{theorem}
\begin{proof}
Write $U(X,R;\mathcal B)$ for this mean. Mean values with different
fixed annuli at the six coordinates are bounded by H\"older's
inequality. Set $M=X^{1/8}$. In a solution
$\sum_{i=1}^3 x_i^3=\sum_{i=1}^3 y_i^3$ extract a principal
divisor $(a_i)$ of $(y_i)$ coprime to $(x_i)$ with
$M\le |\mathrm N a_i|\le M R^h$, where $h$ is the ordinary
class number. Indeed, group the coprime prime factors of $(y_i)$
into principal blocks of length at most $h$, leaving at most
$h-1$ primes: among $h$ successive classes some consecutive
nonempty sequence sums to zero, by partial sums. Multiply blocks
until the threshold is reached. Failure forces
$\mathrm N((x_i)+(y_i))\gg X/(MR^h)$, since $(y_i)$ is squarefree
and has norm $\asymp X$. The number of such pairs is
$O_\epsilon(X^{1+\epsilon}MR^h)$ by summing
$O((X/\mathrm N\mathfrak a)^2)$ over their possible common ideals.
For each pair the four remaining variables with a fixed difference
of two sums of two cubes contribute $O(X^{2+\epsilon})$ by the
fourth moment and Cauchy's inequality. This disposes of the
exceptional solutions.

Choose a balanced generator of each principal ideal and put its norm
in a dyadic interval $[D,2D]$. All the $D$ lie between $M/2$ and
$MR^h$. The quotients have volume scale $Y=X/D$ and lie in fixed
enlarged annuli. We allow a fixed enlargement at each use of the
argument. Write $f_a$ for $f_S$ restricted by $(x,a)=1$, and $g_Y$
for the squarefree smooth cubic sum in the enlarged annuli.
Orthogonality, with positive representation counts, permits
dropping the remaining restrictions on the quotients. Apply H\"older
to the resulting integrand with the six factors
$|f_{a_i}|^{1/2}$ and
$(|f_{a_i}|^2|g_{Y_i}(a_i^3\alpha)|^4)^{1/4}$, using exponents
$12$ and $4$, respectively. Sum over the generators, using H\"older
again. Apart from logarithms and the exceptional solutions, this
bounds $U(X,R;\mathcal B)$ by
\[
 U(X,R;\mathcal B)^{1/4}
 \prod_{i=1}^3 (D_i^3 I(D_i))^{1/4},\qquad
 I(D)=\sum_{|\mathrm N a|\asymp D}
      \int |f_a(\alpha)|^2 |g_Y(a^3\alpha)|^4\,d\alpha .
\]
Here the restricted sixth moments are bounded by the unrestricted
one by orthogonality. There are $O(D)$ balanced generators.
In $I(D)$ partition the units modulo $a^3$ into
$O_\epsilon(X^\epsilon)$ sets on each of which cubing is injective,
and apply Cauchy's inequality over the sets. The kernel of cubing
on the local units has at most three roots away from $3$; over
$3$ it is bounded uniformly in the depth by Hensel's lemma at
fixed depth greater than $2v(3)$. Hence the two ordinary variables
have $x-x'=a^3 u$. The case $u=0$ costs
$O_\epsilon(DXY^{2+\epsilon})$. For the other solutions use
\[
      4(x^3-x'^3)=a^3u(3z^2+u^2a^6),\qquad z=x+x'.
\]
Dropping parity and coprimality now, all variables have independent
smooth nonnegative majorants. The variable $u\ne0$ has volume scale
$H=X/D^3$, $z$ scale $X$, and $a$ scale $D$. Write
\[
 F(\alpha)=\sum_{u,z,a} w_H(u) w_X(z)w_D(a)
       \mathrm e_{\mathbb K}\bigl(\alpha u(3z^2+u^2a^6)\bigr),
\]
where zero is omitted from the $u$-sum. By orthogonality and
H\"older's inequality the nondiagonal part of $I(D)$ is
\[
 \ll_\epsilon X^\epsilon (\textstyle\int|F|^3)^{1/3}
                      U(Y,R;\mathcal B')^{2/3}.
\] The quotient sum has argument $4\alpha$ after the difference identity; multiplication by $4$ preserves Haar measure.

We prove, uniformly in these scales,
\begin{equation}\label{nf3:smoothaux}
        \int |F|^3 \ll_\epsilon X^{\epsilon}R^{O_{\mathbb K}(1)}
                           H^2 X^{3/2}D^{3/2}.
\end{equation}
First $\int|F|^2\ll_\epsilon X^\epsilon HXD$.
Indeed at a nonzero frequency $b$ the element $u$ divides $b$,
and $3z^2+(ua^3)^2=b/u$ is a binary norm equation. Its two
linear factors over $\mathbb K(\sqrt{-3})$ have product $b/u$.
Their ideals and then their generators of polynomially bounded height
have $O_\epsilon(X^\epsilon)$ possibilities, by the divisor bound
and the unit theorem. The same argument applies to the two factors
in $\mathbb K$ if the extension splits. Nonzero factors bounded
above at every embedding are bounded below by an inverse power of
$X$ there, by their integral norms. The zero coefficient is
$O(HD)$ and $(HD)^2\ll HXD$.

Cauchy's inequality gives $|F|^2\le J E$, where $J$ and $E$
are the $w_H$-weighted sums over $u$ of the squared moduli of
the quadratic $z$-sum at $3u\alpha$ and the sixth-power $a$-sum
at $u^3\alpha$. Write $E=E_0+E_1$ with $E_0\ll HD$ its constant
coefficient. Then $\int|E_1|^2\ll_\epsilon X^\epsilon H D^2$.
For a nonzero coefficient one has
$u^3(a_1^6-a_2^6)=b$; after choosing the divisor $u$ and the
nonzero divisor $a_1-a_2$ of $b/u^3$, a polynomial of degree five
determines $a_2$. The same bounded-height ideal-divisor estimate
therefore applies.

Put $R_0=Q B(Y^{3/d}(\alpha-\gamma))$ on the disjoint arcs
$R_0\le X$, with $Q=\mathrm N\mathfrak q_\gamma$ as before.
The quadratic lattice estimate of Section~\ref{nf3:averaged}, in its inner squared-sum form before the outer Cauchy factor $H$, gives
$J\ll X^\epsilon HX$ off these arcs and
$J\ll X^\epsilon HX^2/R_0$ on them; here $HX^2=Y^3$.
For $V=\int_{\mathfrak m}|F|^3$ on the complement,
$|F|\le J^{1/2}(E_0^{1/2}+|E_1|^{1/2})$ gives
\[
 V\ll_\epsilon X^\epsilon\left(
   (H^2XD)^{1/2}\int |F|^2+
   (HX)^{1/2}(HD^2)^{1/4}
   (\int|F|^2)^{1/4} V^{1/2}\right).
\]
For the second term apply H\"older and then interpolate the
$8/3$-moment between the second and third moments on the
complement. This proves \eqref{nf3:smoothaux} there.

On an arc let $b=a_1^6-a_2^6\ne0$ and $Q_b=\mathrm N\mathfrak q_{b\gamma}$.
Twice differencing the cubic $u$-sum as in \eqref{nf3:fourweyl} gives
\[
 H^{-1}\left|\sum_{u}w_H(u)\mathrm e_{\mathbb K}(b\alpha u^3)\right|
    \ll_\epsilon X^\epsilon
      (H^{-1}+Q_b^{-1}+R_0/H^3)^{1/4}.
\]
Indeed the rank-two contribution gives $H^{-1}$ after taking
fourth powers. The rank-one center $b\gamma$ gives
$1/(Q_b B(H^{3/d}b(\alpha-\gamma)))\le1/Q_b$.
Any different center $\xi$ satisfies
\[
 (\mathrm N\mathfrak q_\xi) B(H^{3/d}(b\alpha-\xi))\gg
 H^3/\{Q_b B(H^{3/d} b(\alpha-\gamma))\}\gg H^3/R_0
\]
by the norm lower bound for $b\gamma-\xi$ and
$|b_v|\ll D^{6/d}$. Omitting the origin costs $O(1)$.
The product $Q/Q_b$ is the norm of
$(b)+\mathfrak q_\gamma$. For any ideal divisor $\mathfrak c$ of
$\mathfrak q_\gamma$, at most
$O_\epsilon(X^\epsilon\min(D^2,D^6/\mathrm N\mathfrak c))$
pairs give a nonzero $b\in\mathfrak c$: the values lie in a box
of volume $O(D^6)$ and each has $O_\epsilon(X^\epsilon)$
representations by the same divisor argument. A nonempty range
has $\mathrm N\mathfrak c\ll D^6$. Interpolation and the ideal
divisor bound imply
\[
 \sum_{a_1,a_2:\ b\ne0}(Q/Q_b)^{1/4}
                         \ll_\epsilon D^3 X^\epsilon.
\]
Consequently
$E\ll_\epsilon X^\epsilon R^{O_{\mathbb K}(1)}HD
 (1+D^2 Q^{-1/4})$; here $D^4\ll H R^{O(h)}$ and
$XD^4\ll H^3R^{O(h)}$.

At fixed center and $R_0\asymp Z$ the arc volume is
$O_\epsilon(X^\epsilon Y^{-3}Z/Q)$, and there are at most $Q$
centers of denominator ideal $\mathfrak q$. Thus
$|F|^3\le(JE)^{3/2}$, summed over the centers and dyadic
$Z\le X$, has integral at most
\[
 X^\epsilon R^{O_{\mathbb K}(1)}
 (H^2X^2D)^{3/2}Y^{-3}
                  (X^{1/2}+D^3X^{1/8}).
\]
Indeed summation of $Q^{-3/8}$ over ideals of norm at most $Z$
costs $O_\epsilon(Z^{5/8+\epsilon})$.
Since $D\ll X^{1/8}R^h$ and $HD^3=X$, this proves
\eqref{nf3:smoothaux} on the arcs as well.

Inserting this estimate and the diagonal bound gives the recursion
\[
 U(X,R;\mathcal B)\ll_\epsilon X^\epsilon R^{O_{\mathbb K}(1)}
 \left(X^3 M^2+
       X^{7/6}M^{3/2}\max_{D\asymp M R^{O(h)}}
             U(X/D,R;\mathcal B')^{2/3}\right),
\]
where the maximum ranges over the chosen dyadic scales.
Beginning with the fourth and eighth moments, $7/2$ is an
initial sixth-moment exponent. Each application replaces an
exponent $\lambda\ge13/4$ by
$\max(13/4,\,65/48+7\lambda/12)$, up to errors arbitrarily
small with $\eta$ and the divisor-bound slack. This affine map
fixes $13/4$. Make finitely many iterations for a prescribed
$\epsilon$, then choose $\eta$ sufficiently small that all
powers of $R$ are absorbed at all the resulting scales.
This proves the moment estimate.

For the counting conclusion we require a lower bound in fixed residue
classes uniform in the choice of the fixed modulus. Fix
$0<\theta\le\eta$ and take $R=X^\theta$. There is $c>0$, depending
on $\theta$, the box and the field, such that each class modulo any fixed
$\mathfrak m$ contains at least $cX/\mathrm N\mathfrak m$ smooth
roots in a prescribed nonempty archimedean box of scale $P$, once
$X$ is sufficiently large in terms of $\mathfrak m$. To see this
choose finitely many disjoint fixed intervals $(a_i,b_i)\subset
(0,\theta)$ with $1-\theta<\sum_i a_i<\sum_i b_i<1$.
This is possible by taking sufficiently many short intervals of
comparable positive size. From each norm window $(X^{a_i},X^{b_i})$
choose a different principal prime ideal with a balanced generator
$\pi_i$. For fixed $\mathfrak m$ all these primes are coprime to
$\mathfrak m$ once $X$ is large. The sum of their reciprocal norms in each window is
bounded below, by the prime ideal theorem in an ideal class and
partial summation. With $\pi=\prod_i\pi_i$, the remaining cofactor
$y$ in $\pi y\in P\mathcal B$ runs over a box with uniformly bounded
eccentricity, whose side grows by a fixed power of $P$. In the
unique class giving the prescribed root congruence it has
$\gg_{\mathcal B}X/(|\mathrm N\pi|\mathrm N\mathfrak m)$
choices, with the implicit constant independent of the fixed
modulus once $X$ is large. Its norm is less than $X^\theta$.
Each resulting root has bounded multiplicity since each specified
prime has norm at least $X^{\min a_i}$. This proves the lower bound.
For squarefreeness and pairwise coprimality impose the local
conditions at primes of norm up to $L$. At each prime of norm $q$
the local factor is $(1-1/q)^3(1+3/q)$. Their product is positive.
The unrestricted tail for excluded triples from larger primes is
$O_{\mathcal B}(X^3\sum_{\mathrm N\mathfrak p>L} (\mathrm N\mathfrak p)^{-2})$
by lattice packing as in Proposition~\ref{nf3:manyprimitive}.
Choose $L$ large and then hold it fixed. There remain $\gg X^3$
triples with the required restrictions. Their collision count is
$O_{\epsilon}(X^{13/4+\epsilon})$ by the moment estimate and
H\"older. Cauchy's inequality gives $\gg X^{11/4-\epsilon}$
distinct targets, as asserted.
\end{proof}

\subsection{Representations of every sufficiently large target}\label{nf3:alllarge}
Over $\mathbb Q$ every sufficiently large integer is a sum of seven
positive cubes, by Linnik \cite{nf3:linnik}. Over general number fields
the four-cube exceptional set yields eight summands by addition. It also
shows precisely how a sixth-moment hypothesis would lower eight to seven.

\begin{theorem}\label{nf3:global}
For every number field $\mathbb K$ one has $G_{\mathbb K}(3)\le8$.
For every fixed $\epsilon>0$ and every sufficiently large totally positive
eligible target $\nu$, with $N=|\mathrm N\nu|$, there are
$\gg_{\mathbb K,\epsilon}N^{5/3-\epsilon}$ ordered representations by
eight totally positive cubes whose roots satisfy
$|x_v|\le C_{\mathbb K}|\nu_v|^{1/3}$ at every archimedean place $v$,
and hence $|\mathrm N x|\le C_{\mathbb K}^d N^{1/3}$.
Moreover, if the signed six-cube
collision count in every embedding box of root-volume scale $V$ is
$O_{\mathbb K,\epsilon}(V^{3+\epsilon})$, then $G_{\mathbb K}(3)\le7$.
\end{theorem}
\begin{proof}
Put $A=N^{1/d}$. Multiplication by a cube of a totally positive unit
places $\nu/A$ in a fixed compact subset of the positive archimedean
region, bounded away from zero at every place, by Dirichlet's unit
theorem. Each normalized target admits an open box of decompositions
$\nu/A=z+z'$ with both $z,z'$ in a single such compact subset.
For example take a small box about $\nu/(2A)$. The box can have fixed
positive volume uniformly, and a finite choice of surrounding compact
regions suffices.

At all but finitely many finite places every target is a sum of four
integral cubes, by the smooth residue-field solutions in the proof of
Proposition~\ref{nf3:covol}. At the remaining places the five-cube local
solution to an eligible $\nu$ together with $1,-1,0$ gives eight roots
that can be partitioned into two groups of four, each with a unit
coordinate. Write their sums as $b$ and $\nu-b$. Each sum remains
four-cube soluble under a perturbation by $\mathfrak p^{2v(3)+1}$,
by strong Hensel in its own unit coordinate. The Chinese remainder
theorem prescribes a single class of $m$ modulo a fixed ideal such
that both $m$ and $\nu-m$ are everywhere locally four-cube soluble.
The class may depend on $\nu$, while its modulus and index are fixed.

There are $\gg_{\mathbb K}N$ choices of $m$ in this class and in the
chosen dilated archimedean box. By Corollary~\ref{nf3:manyfour}, all
but $O(N^{1-\eta})$ of these choices give
$\gg N^{1/3-\epsilon/2}$ four-cube representations for each of
$m,\nu-m$, with roots of embedding height $O(A^{1/3})$. The ordered
eight-tuples for different intermediate sums $m$ are disjoint.
On this box $|\nu_v|\asymp A$ at every place. Dividing the roots
by the balancing unit and the target by its cube preserves each
ratio $|x_v|/|\nu_v|^{1/3}$, proving the asserted count and bounds.
In a totally real field positivity gives
the sharper bound $0<\mathrm N x<N^{1/3}$ for every one of these roots.

For the conditional assertion split a local five-cube solution with an
added pair $1,-1$ into three and four roots, with one of the added
roots in each half. At each of the finitely many places prescribe the
individual residues of the first three roots to a fixed sufficiently
large depth. Then the complementary target remains locally soluble in
four cubes. Choose small positive archimedean boxes for the three roots
whose cube sum is close to $\nu/2$ after normalization. There are
$\gg N$ triples in these boxes and residue classes. Their collision
count, by the stated hypothesis with $V=N^{1/3}$, is
$O_\epsilon(N^{1+\epsilon})$. Cauchy's inequality supplies
$\gg_\epsilon N^{1-\epsilon}$ distinct sums. Taking $\epsilon$ smaller
than the exceptional-set saving in Theorem~\ref{nf3:four} leaves some
sum whose complement has a four-cube representation. The arguments
apply also when there are complex embeddings, since no sign is imposed
there and the chosen root boxes have bounded modulus at each place.
\end{proof}

\subsection{Seven positive cubes}\label{nf3:sevenreal}
The cubic-pair identity used by Bombieri and Ramar\'e
\cite[\S2]{nf3:ramare} admits a choice of ternary lattice that works
at every dyadic completion. A prescribed valuation at one further
prime removes the spinor exceptions. For totally imaginary fields an
elementary approximation argument gives the same bound; the tangent at a
suitable integral point on the Fermat cubic reduces it to six
(Corollary~\ref{nf3:imag6} and Theorem~\ref{nf3:fermat-six-local}).

\begin{theorem}\label{nf3:realG7}
Let $\mathbb K$ be a number field. Every sufficiently
large totally positive element $\nu$ of the additive group generated
by the integral cubes is a sum of seven cubes of strictly totally
positive algebraic integers. All seven roots are nonzero and satisfy
$|x_v|\le C_{\mathbb K}|\nu_v|^{1/3}$ at every archimedean place,
so $|\mathrm N x|\ll_{\mathbb K}|\mathrm N\nu|^{1/3}$. For a totally real field one has
$0<\mathrm N x<(\mathrm N\nu)^{1/3}$. In particular
$G_{\mathbb K}(3)\le7$.
\end{theorem}

\begin{proof}
Write $\mathcal O$ for the ring of integers and $d=[\mathbb K:\mathbb Q]$.
First suppose the field is totally real. We record the ternary input. Suppose $L$ is a fixed totally
positive definite integral ternary lattice over $\mathcal O$, that
$Q\gg0$ is locally represented, and that its valuations at the finite
anisotropic primes are bounded. Suppose also that
$v_{\mathfrak l}(Q)=1$ at a fixed odd prime $\mathfrak l$ where $L$
is unimodular. For $\mathrm NQ$ sufficiently large, $L$ represents
$Q$. Indeed, the two exceptions in
\cite[Corollary~5.6]{nf3:schulzesurvey} require a spinor exception
for the genus of $L$ or of $\mathfrak pL$ for some $\mathfrak p$.
The necessary spinor-norm condition
\cite[Theorem~1]{nf3:schulzespin} fails for both. At
$\mathfrak l$ the quadratic algebra attached to the orthogonal
complement of the represented line is ramified, since its
discriminant squareclass is $-Q\det L$ in the quadratic-matrix
convention. Its unit norms have square residue. In contrast, the
unimodular ternary lattice represents each unit by a primitive vector
at $\mathfrak l$, by isotropy over the residue field and Hensel's
lemma. Reflections in vectors of a given unit value and of value one
give an integral rotation of that spinor norm. Local scalar multiples
of the lattice have the same integral rotation group, including when
$\mathfrak p=\mathfrak l$. This proves the assertion and allows
arbitrary ideal classes. The norm threshold can be ineffective.

Use the lattice
\[
 q(x,y,\zeta)=16x^2+4xy+3y^2+32\zeta^2,\qquad
 4q=L_1^2+11y^2+8z^2,\quad L_1=8x+y,\quad z=4\zeta.
\]
It is unimodular away from $2,11$. At an odd place over 11 its
binary part $16x^2+32\zeta^2$ represents each unit, since a
nondegenerate binary form over a finite field of odd order
represents every nonzero element. At a dyadic place $Q\equiv3\pmod8$
is represented by $3y^2$, by Hensel's lemma. In every integral
representation of such a $Q$ one has $y^2\equiv1\pmod4$ and hence
$y\equiv1\pmod2$. In fact, if $v(y-1)<v(2)$ then
$v(y+1)=v(y-1)$ and $v(y^2-1)<2v(2)$.

At three auxiliary elements $u,v,w$ put
\[
 B=u^2v^2w^2,\qquad
 E=\tfrac14(u^6v^6+11^3u^6w^6+8^3v^6w^6).
\]
Our fixed dyadic congruences below give $E\in\mathcal O$. The
cubes of the following three pairs sum to $Es^3+3sBq$:
\begin{equation}\label{nf3:pairidentity}
 \frac{u^2v^2s\pm wL_1}{2},\qquad
 \frac{11u^2w^2s\pm vy}{2},\qquad
 \frac{8v^2w^2s\pm uz}{2}.
\end{equation}
We shall choose a positive seventh root $t$ and a positive center
$s$ so that
\begin{equation}\label{nf3:seventarget}
 Q=(\nu-t^3-Es^3)/(3sB)
\end{equation}
is a large locally represented integer.

\smallskip\noindent
\emph{Dyadic places.}
In the residue field $k$ choose $t$ with
$\operatorname{Tr}_{k/\mathbb F_2}(t^3)=
\operatorname{Tr}_{k/\mathbb F_2}(\nu)$. This is possible because
$\nu$ belongs to the additive span of cubes: if its trace is one at
least one cube has trace one. Lift $t$ arbitrarily. Set $u=v=s=1$
temporarily and write $w=1+2D$. Then $E+9B$ equals
\[
 H(D)=\tfrac14\{1+1843(1+2D)^6+36(1+2D)^2\}\in\mathbb Z[D].
\]
Modulo the maximal ideal this is $D^2+D$, and $H'$ reduces to one.
The trace criterion and Hensel's lemma give integral $D$ with
$H(D)=\nu-t^3$. Prescribe $u,v,s\equiv1\pmod{32}$,
$w\equiv1+2D\pmod{32}$ and the chosen $t$ modulo 8, separately
at each dyadic completion. These prescriptions give
$Q\equiv3\pmod8$ and make all the quantities in
\eqref{nf3:pairidentity} integral, since $L_1\equiv y\equiv1$
and $z\equiv0\pmod2$. They apply at every ramification index.

\smallskip\noindent
\emph{The other fixed places.}
Choose $b\in\mathcal O$ with $b^3\equiv\nu\pmod{3\mathcal O}$.
The Frobenius endomorphism on $\mathcal O/3\mathcal O$ makes its
image an additive group, so eligibility suffices even at ramified
primes. Put $c=(\nu-b^3)/3$. Require $u,v,w\equiv1\pmod9$.
Then $E\equiv461\pmod9$. With $t\equiv b+s\pmod3$ the residue of
$(t^3+Es^3-b^3)/3$ at $\mathfrak p\mid3$ is
\[
 b^2s+bs^2+154s^3=s(s-b)^2.
\]
This polynomial takes at least two values on each residue-field unit
group. For the field of three elements this follows by using $1,-1$,
and for larger residue fields it follows from its degree. Choose a
unit value of $s$ where the displayed expression differs from $c$.
Lifting this class modulo 3 makes $Q$ integral and a unit at every
prime above three.

Fix integral representatives $\mathfrak c$ of the ideal classes
coprime to $2\cdot3\cdot11\,\operatorname{disc}(\mathbb K)$.
At each prime dividing a representative, each prime above $11$ and
each remaining odd prime over the field discriminant away from three,
require $u=v=w=s\equiv1$ modulo that prime.
One of $t=0,1$ makes $\nu-t^3-E$ a unit. Fix an odd prime
$\mathfrak l\nmid 2\cdot3\cdot11\cdot461\,\operatorname{disc}(\mathbb K)$ with residue cardinality
$2\pmod3$, outside the supports of the representatives.
Such primes exist by Chebotarev in
$\mathbb K(\sqrt{-3})/\mathbb K$. Take $u=v=w\equiv1\pmod{\mathfrak
l^2}$. With $s\equiv1\pmod{\mathfrak l}$ one can prescribe
$t^3\equiv\nu-461\pmod{\mathfrak l}$ by bijectivity of cubing.
Fix a lift of $t$ modulo $\mathfrak l^2$. For this $t$ the derivative of $461s^3$ is a unit. Choose $s$ modulo
$\mathfrak l^2$ so that the numerator of \eqref{nf3:seventarget}
has valuation exactly one. All the prescriptions so far involve a
fixed modulus depending only on the field, and all prescribed
classes of $u,v,w,s$ are units.

\smallskip\noindent
\emph{Choice and size of the auxiliary primes.}
Let $N=\mathrm N\nu$, $P=N^{1/(3d)}$, and $U=C_{\mathbb K}\log N$.
First balance the target by a cube of a totally positive unit, so
that $\nu_j\asymp P^3$ at each real embedding $j$. Choose distinct
prime elements $u,v,w$, avoiding $(\nu)$ and the fixed primes,
in their prescribed classes, with
$U<|u_j|,|v_j|,|w_j|<1.01U$ and exactly one negative embedding
for each element. Their signed norms are $1\pmod3$ and their
absolute norms are $2\pmod3$. Mitsui's prime element theorem in
convex boxes and ray classes \cite[Theorem~1.2.1]{nf3:kai} gives
$\gg_{\mathbb K}U^d/\log U$ choices at this fixed modulus. One ideal
has boundedly many generators in such a box, by the unit theorem,
and the number of prime ideals dividing $(\nu)$ in this norm range
is $O_{\mathbb K}(\log N/\log U)$. Taking the fixed constant
$C_{\mathbb K}$ large enough covers also $d=1$.

Cubing permutes the units modulo each of $(u)^2,(v)^2,(w)^2$:
their orders are $q(q-1)$ with $q\equiv2\pmod3$. Reducing the
definition of $4E$ modulo the square of each auxiliary prime leaves,
respectively, the unit $8^3v^6w^6$, $11^3u^6w^6$, or $u^6v^6$.
Thus $Es^3\equiv\nu$ gives a unit class of $s$ at each of these
three prime squares. Together with the fixed conditions this gives
a unit class $r$ modulo an ideal of norm $O_{\mathbb K}(U^{6d})$,
divisible by $3\mathcal O$.

We choose a totally positive $s$ with squarefree $(s)$ and with
residue cardinality $2\pmod3$ at each prime dividing it. Fix a prime
element $\pi_0\equiv1\pmod3$, with exactly one negative
embedding, avoiding all the fixed primes. Mitsui's theorem supplies
one, and then $|\mathrm N\pi_0|\equiv2\pmod3$. For large $N$ it
is coprime to $uvw$. If $\mathrm N r\equiv2\pmod3$, put
$s=\pi$; otherwise put $s=\pi_0\pi$. Choose the prime element
$\pi$ in the class $r$ or $r/\pi_0$ accordingly, additionally
requiring $\pi\equiv1\pmod{\pi_0}$. In the second case take its
signs to agree with those of $\pi_0$. In both cases positivity
of $s$ implies $|\mathrm N\pi|\equiv2\pmod3$.

For the size condition write $\vartheta=10^{-5}$ and
$\gamma=\vartheta/16$, and require
\[
  1-\gamma< s_j\big((1+\vartheta)E_j/\nu_j\big)^{1/3}<1+\gamma
       \qquad(1\le j\le d).
\]
This is a box in the required sign chamber for $\pi$, with side
lengths comparable to $P/U^4$. Its modulus has norm
$O_{\mathbb K}(U^{6d})$, a fixed power of the logarithm of its side
scale. The logarithmic-modulus form of Mitsui's theorem cited above
therefore gives such a prime, uniformly in $\nu$. Finally choose
$t^3\equiv\nu\pmod s$, possible by bijectivity of cubing in the
residue fields, and $t\equiv0\pmod{uvw}$, along with its fixed
classes. The Chinese remainder theorem and a fixed fundamental
parallelepiped of $\mathcal O$ give a representative with
$0<t_j\ll_{\mathbb K}s_j|u_jv_jw_j|$. Indeed a fixed integer
multiple of $suvw\mathcal O$ is contained in the modulus ideal,
and a translate of its embedding parallelepiped lies in the
positive box with these side scales. The conditions cover all
factors $3sB$ in \eqref{nf3:seventarget}.

\smallskip\noindent
\emph{Positivity.}
At every real embedding $B_j\ge U^6$ and
$E_j<(1844/4)(1.01)^{12}U^{12}<521U^{12}$.
Moreover $t_j^3/s_j^3\ll_{\mathbb K}U^9=o(E_j)$.
The box for $s$ gives
$\vartheta E_j/2<\nu_j/s_j^3-E_j<2\vartheta E_j$.
Consequently $Q_j\gg_{\mathbb K}P^2/U^2$ and
\[
  0<4Q_j/s_j^2<\tfrac{8}{3}\vartheta E_j/B_j
     <0.014U^6 .
\]
Each of $u_j^4v_j^4/w_j^2$, $11^3u_j^4w_j^4/v_j^2$ and
$8^3v_j^4w_j^4/u_j^2$ exceeds $0.98U^6$.
The quadratic target is locally represented at $2,11$ by
construction, and at every other finite place by unimodular ternary
isotropy. Its valuations at $2,11$ are zero, its valuation at
$\mathfrak l$ is one and its norm tends to infinity. The ternary
result now yields $q(x,y,\zeta)=Q$. The identity for $4q$ bounds
$L_{1,j}^2,11y_j^2,8z_j^2$ individually by $4Q_j$. The preceding
strict inequalities put every perturbation in
\eqref{nf3:pairidentity} strictly inside its positive center. Adding
$t^3$ gives the representation. Dividing its roots by the balancing
unit recovers the original target, and taking the product over real
embeddings of $0<x_j^3<\nu_j$ gives the asserted norm bound.

\smallskip\noindent
\emph{Mixed signature.}
Suppose next that there is at least one real and one complex place.
The same ternary lattice is definite at the real places and is
isotropic at every complex place. By the strong approximation theorem
in the form of Eichler's spinor-genus theorem
\cite[Theorem~9.18]{nf3:eichler}, each proper spinor genus consists
of one proper lattice class. The marked-prime spinor-norm obstruction
above still excludes a spinor exception, by \cite[Theorem~1]{nf3:schulzespin},
which applies to arbitrary number fields. Local representation
therefore gives a representation by $L$.

There is also uniform archimedean control on a representing vector.
The group $G=\operatorname{SO}(q)$ is anisotropic over $\mathbb K$
because $q$ is definite at a real place. By the compactness theorem
for arithmetic groups \cite[\S11]{nf3:bhc}, applied after restriction
of scalars, the integral stabilizer $\Gamma$ has compact quotient in
$G_\infty=\prod_{j\mid\infty}G(\mathbb K_j)$.
Write $G_\infty=\Gamma C$ with $C$ compact. Given an integral
representing vector $v$ for a nonzero $Q$, take
$v^0_j=(\sqrt{Q_j}/4,0,0)$, with either square root at a complex
place. Witt's extension theorem supplies $g\in G_\infty$ with
$v=gv^0$; a reflection in the orthogonal complement adjusts the
determinant if necessary. Writing $g=\gamma c$ replaces $v$ by the
integral vector $\gamma^{-1}v=cv^0$, whose coordinates at every
place are $O_{\mathbb K}(|Q_j|^{1/2})$.

The prime construction uses one negative real embedding and positive
signs at all the other real places. At the complex places require the
three auxiliaries to lie in small squares about the positive real
number $U$, keeping $|E_j|\asymp U^{12}$ and $|B_j|\asymp U^6$.
For $s_j$ take a small square about a chosen cube root of
$\nu_j/((1+\vartheta)E_j)$, with side length a fixed small fraction
of its modulus. These convex regions have the same volume order as
before, counting each complex place twice. Use a translate of the
same ideal parallelepiped to keep $t$ strictly positive at every
real place and $|t_j|\ll_{\mathbb K}|s_j u_jv_jw_j|$ at all places.
The target satisfies $|Q_j|\ll_{\mathbb K}P^2/U^2$ at complex
places. The compactness bound for its representing vector gives
$O_{\mathbb K}(P)$ for every root there. At real places the positive
quadratic identity still gives the strict inequalities already proved.
Dividing out the balancing unit gives the claimed placewise bounds
in mixed signature.

\smallskip\noindent
\emph{Totally imaginary fields.}
Here one can construct the seven roots by successive approximations.
We first express every eligible residue modulo $6\mathcal O$ as
three cubes. Modulo 3 one cube suffices by Frobenius. In a dyadic
residue field $k$, choose $z$ with
$\operatorname{Tr}(z^3)=\operatorname{Tr}(\nu+1)$, using the
eligibility of $\nu+1$. Thus
$x^2+x=\nu+1-z^3$ is soluble in $k$ and lifts modulo $2\mathcal O$
at this prime by its unit derivative. In characteristic two this
gives $\nu=z^3+x^3+(x+1)^3$. The Chinese remainder theorem now
prescribes three individual classes modulo 6.

Balance $\nu$ so that $|\nu_j|\asymp P^3$. In a prescribed class
modulo 6 one can approximate an arbitrary tuple in
$\prod_{j\mid\infty}\mathbb C$ within bounded error by the lattice
covering radius. Apply this successively to a tuple of cube roots
of the current remainder, subtracting the cube of the approximation
at each step. If the remainder has maximum modulus $T$, the next
has maximum modulus $O_{\mathbb K}(1+T^{2/3})$. The three roots
have size $O_{\mathbb K}(P)$ and leave an element $6h$ of size
$O_{\mathbb K}(P^{8/9})$. The identity
\[
  6h=(h+1)^3+(h-1)^3+(-h)^3+(-h)^3
\]
completes the proof. To obtain seven nonzero roots, finitely many
translates by $6j$ of each rounded root suffice: avoid zero at the
first two steps, then avoid zero and the at most nine values of the
third root giving $h\in\{0,1,-1\}$. These choices preserve the
bounds and the residues. Total positivity is vacuous; undoing the
unit cube gives the placewise bound and hence the norm estimate.
\end{proof}

\begin{corollary}\label{nf3:imag6}\label{nf3:fermat-tangent}
Let $\mathbb K$ be totally imaginary. Suppose there are $a,b,c\in
\mathcal O_{\mathbb K}$ with $abc\ne0$ such that
\[
             a^3+b^3+c^3=0,\qquad abc^3\mid 2.
\]
Then $G_{\mathbb K}(3)\le6$. In particular this holds if
$\mathbb K$ contains a primitive ninth root of unity or a cube root
of $2$, or if the equation $a^3+b^3=1$ has a solution in units.
\end{corollary}
\begin{proof}
The tangent at $(a:b:c)$ to the Fermat cubic gives the identity
\[
  (aZ+b^2)^3+(bZ-a^2)^3+(cZ)^3
      =b^6-a^6-3abc^3 Z.
\]
Indeed its quadratic coefficient vanishes and the linear coefficient
is $3ab(a^3+b^3)$. Subtract $b^6-a^6$, itself a sum of signed cubes,
from the balanced target. This fixed translation preserves eligibility
and the archimedean upper bound. The three successive approximations
in the proof of Theorem~\ref{nf3:realG7} leave $6y$, with
$|\sigma(y)|\ll_{\mathbb K,a,b,c}P^{8/9}$ at every embedding.
In the displayed identity take $Z=-2y/(abc^3)\in\mathcal O_{\mathbb K}$.
Its three roots have size $O_{\mathbb K,a,b,c}(P^{8/9}+1)$.
As in the preceding proof, finitely many translates of the last
approximating root avoid zero in all six roots. For a solution of
$a^3+b^3=1$ in units take $c=-1$.
For a primitive ninth root $\xi$ take $(a,b,c)=(1,\xi,\xi^2)$.
For $\beta^3=2$ take $(a,b,c)=(1,1,-\beta)$.
\end{proof}

\begin{lemma}\label{nf3:triadic-three}
Let $A$ be the ring of integers of a finite extension of $\mathbb Q_3$
whose residue field $F$ has at least nine elements. An element $n\in A$
is a sum of integral cubes if and only if $n$ is a cube modulo $3A$.
Every such $n$ is a sum of three cubes, at least one with unit root.
\end{lemma}
\begin{proof}
Necessity follows from Frobenius. For sufficiency choose $w\in A$
with $w^3\equiv n\pmod{3A}$, and put $t=(w^3-n)/3$.
By the identity
\[
 x^3+y^3+z^3=(x+y+z)^3-3(x+y)(y+z)(z+x),
\]
it suffices to solve $X+Y+Z=2w$, $XYZ=t$. The three roots are then
$w-X,w-Y,w-Z$. We first solve over $F$, denoting reductions by bars.
If $\bar t=0$, take $X=0$ and two nonzero elements with sum $2\bar w$.
If $\bar t\ne0$ and $\bar w=0$, choose
$(X,Y,Z)=\lambda(1,l,-1-l)$ with $l\notin\{0,1,-1\}$.
Frobenius gives the required $\lambda$.
If $\bar t\bar w\ne0$, use the projective cubic
\[
 XYZ=\bar t\big((X+Y+Z)/(2\bar w)\big)^3.
\]
Its partial derivatives are $YZ,XZ,XY$, so the curve is smooth.
The Hasse bound \cite{NFWeilCurves} gives at least
$|F|+1-2\sqrt{|F|}\ge4$ rational points. Only three lie on
$X+Y+Z=0$. Scaling an affine point gives the specified sum and
product. In all three cases $(YZ,XZ,XY)$ is not proportional to
$(1,1,1)$. For nonzero product, proportionality requires $X=Y=Z$,
which was excluded in the zero-sum case and would have zero sum in
the other case. Thus Hensel's lemma lifts the two equations. At
least one root is a unit, since $X=Y=Z=\bar w$ would either
contradict the sum condition or make all three zero.
\end{proof}

\begin{theorem}\label{nf3:triadic-class}
Let $A$ be the ring of integers of a finite extension of $\mathbb Q_3$,
with ramification index $e$ and residue degree $f$. For cubes the
ordinary local Waring number and the width requiring a unit root agree.
Their common value is
\[
\begin{array}{c|c}
\text{completion}&\text{width}\\ \hline
 f\ge2\ \text{or}\ 3\mid e &3\\
 f=1,\ e=1 &4\\
 f=1,\ e\equiv1\pmod3,\ e\ge4 &5\\
 f=1,\ e\equiv2\pmod3,\ e\ge5 &4\\
 \mathbb Q_3(\sqrt{3}),\quad\mathbb Q_3(\sqrt{-3})&4,\quad5 .
\end{array}
\]
In particular the two unspecified cubic entries $e=4,5$ for the pure
radical fields in \cite[\S4, Table~4]{AnthonyRam} are $5,4$.
\end{theorem}
\begin{proof}
Write $\pi$ for a uniformizer. The local target group consists of the
elements $n$ which are cubes modulo $3A$. Indeed, after subtracting
one cube, write the remainder as $6y$ and use
$(y+1)^3+(y-1)^3+2(-y)^3=6y$. At least one of these last four
roots is a unit. This proves the upper bound five in all cases.

If $f\ge2$, Lemma~\ref{nf3:triadic-three} gives three.
If $3\mid e$, choose $w^3\equiv n\pmod{3A}$. Replacing $w$ by
$w+\pi^{e/3}b$ varies $(w^3-n)/3$ modulo $\pi$ by a unit times
$\bar b^3$. Make it zero. In the product-sum system in the preceding
lemma one can now take $(0,X,Y)$ with $X,Y\ne0$ and $X+Y=2\bar w$,
also over $\mathbb F_3$. Its nonsingular lift gives three cubes.
Two do not suffice in any completion. In a primitive pair with
nonunit sum of cubes both roots are units. With $s=x+y$ the value
is $s^3-3xys$. Below valuation $3e/2$ its valuation, if finite,
must be a multiple of three, by comparing $3v(s)$ with $e+v(s)$.
Removing a common uniformizer factor changes the valuation by a
multiple of three. There is an integer $j\ge0$ with
$e+j<3e/2$ and $3\nmid(e+j)$: take $j=0$ if $3\nmid e$,
and $j=1$ otherwise. Thus the eligible element $3\pi^j$ excludes
two cubes.

It remains to treat $f=1$ and $3\nmid e$. For any four roots put
$x_i=a_i+\pi z_i$, where $a_i\in\{0,1,-1\}$ and $Z=\sum z_i$.
Then
\begin{equation}\label{nf3:four-triadic}
 \sum x_i^3\equiv
     \sum a_i + (\pi Z)^3+3\pi\sum a_i^2z_i\pmod{3\pi^2 A}.
\end{equation}
Indeed $\sum z_i^3\equiv Z^3\pmod{3A}$.
For $n=4+3\pi u$, $u\in A$, the reduction of $\sum a_i$ modulo
three is one. If $\sum a_i=-2$ or $1$, then $v(n-\sum a_i)=e$,
which cannot arise on the right, since $3\nmid e$.
Consequently every such representation has $a_i=1$ for all $i$;
in particular $4$ requires at least four cubes. If
$e\equiv1\pmod3$, $e\ge4$, take $u$ a unit. A unit $Z$ would give
valuation three for $n-4$, and hence $\pi\mid Z$. Equation
\eqref{nf3:four-triadic} would then require $(\pi Z)^3$ to have
valuation $e+1$, again not a multiple of three. If $e=2$, then
\[
   u\equiv Z+(\pi^2/3)Z^3\pmod\pi.
\]
For $\mathbb Q_3(\sqrt{-3})$ this precludes any unit $u$.

For the remaining upper bounds choose $w$ with $w^3\equiv n\pmod{3A}$,
and write $t=(n-w^3)/3$. Hold $\sum x_i=w$ fixed. Eliminating
$x_4$ gives an integral polynomial $(\sum x_i^3-w^3)/3$ whose
partial derivatives are $x_i^2-x_4^2$. Write $w_0\in\{0,1,-1\}$
for the residue of $w$, and take the unique integer $M\in[-4,4]$
with $M\equiv w_0^3+3\bar t\pmod9$. If $|M|\le3$, express $M$
as a sum of four members of $\{0,1,-1\}$ including a zero and a
unit. For $M=0$ use $1,-1,0,0$. These residues solve the polynomial
equation modulo $\pi$: after setting $w=w_0$ and eliminating the
last residue, its integer representative changes by a multiple of
three, and its cube by a multiple of nine. The mixed zero and unit
roots make the solution nonsingular, so Hensel's lemma applies.

In the remaining two classes, $w$ is a unit and division by $w^3$
reduces to $n=4+3\pi u$. For $e=1$ this is $4$ modulo nine; a
unit cube of the form $(1+3z)^3$ covers $1+9\mathbb Z_3$,
so three copies of $1$ suffice in addition. Suppose
$e\equiv2\pmod3$. Set $U=\pi^e/3$ and
$Z=\pi^{(e-2)/3}b$. With $x_i=1+\pi z_i$, $\sum z_i=Z$, the
exact four-cube equation is
\[
  u=Ub^3+Z+\pi\bigl(\textstyle\sum z_i^2+\pi C(\mathbf z)\bigr),
 \qquad C(\mathbf z)=\bigl(\textstyle\sum z_i^3-Z^3\bigr)/3 .
\]
Here $C$ is an integral polynomial after substituting $Z=\sum z_i$.
For $e\ge5$, choose $b$ modulo $\pi$ by Frobenius so that
$u\equiv Ub^3+Z$. For $e=2$ and $\mathbb Q_3(\sqrt3)$ the residual
map is $b\mapsto 2b$, also bijective. These are the two ramified
quadratic completions: a uniformizer generates the extension, and
the ramified squareclasses over $\mathbb Q_3$ are $3,-3$.

It remains to prescribe $\sum\bar z_i^2$ with $\sum\bar z_i=\bar Z$,
using a tuple with at least two different entries. For sum zero the
tuples $(0,1,1,1),(1,1,-1,-1),(1,-1,0,0)$ give square sums
$0,1,2$. For sum one use $(1,1,-1,0),(1,0,0,0),(-1,-1,0,0)$;
negation handles sum minus one. Eliminate the fourth coordinate.
A derivative $2(\bar z_i-\bar z_4)$ is nonzero. Hensel's lemma
now prescribes $(u-Ub^3-Z)/\pi$ exactly, giving four unit roots.
\end{proof}
\begin{corollary}\label{nf3:all-local}
At a finite place away from three, the two cubic widths are both two,
except at residue norm seven, where both are three. Thus for every
number field
\[
 W_{\mathbb K}(3)=H_{\mathbb K}(3)
       =\max_{v\mid3}h_v(3)\in\{3,4,5\},
 \qquad H_{\mathbb K}(3)\le G_{\mathbb K}(3)\le7,
\]
where $h_v(3)$ is given by Theorem~\ref{nf3:triadic-class}.
\end{corollary}
\begin{proof}
Write $q$ for the residue norm, in characteristic different from
three. For $q\ge9$ and $t\ne0$ the smooth projective Fermat cubic
$X^3+Y^3=tZ^3$ has at least $q+1-2\sqrt q>3$ rational points,
at most three of them at infinity. For $q=2,5,8$ cubing is
bijective. For $q=4$ the eligible residue set is $\{0,1\}$.
A primitive solution to zero is $(1,-1)$, so Hensel lifting gives
two cubes throughout in these cases. At $q=7$ the cubes are
$0,1,-1$; three suffice, but two omit the residues $3,-3$.
One cube cannot cover a uniformizer, which is a sum of two cubes by
lifting $(1,-1)$ from the residue field. Finally, any finite-place
obstruction to fewer cubes persists on a congruence ball by compactness.
The global additive cube group is dense in the local one.
After choosing an element in that ball, add unbounded positive rational
integers divisible by its modulus to obtain totally positive targets
of unbounded norm in the same ball. The upper bound is
Theorem~\ref{nf3:realG7}.
\end{proof}
The local obstruction is visible in both mixed and totally imaginary
fields. Take $\rho=\sqrt[4]{3}$ and work in $\mathbb Q(\rho)$ or
$\mathbb Q(\rho,\sqrt{-2})$. The first completion at $3$ is totally
ramified of degree four, and $-2$ is a square in $\mathbb Q_3$.
The elements $4+6\rho+9m$ with positive rational integers $m$ are
eligible, are eventually totally positive, and have unbounded
absolute norm. They are excluded by \eqref{nf3:four-triadic}.
Thus $G_{\mathbb K}(3)\ge5$ in both fields.

\begin{theorem}\label{nf3:fermat-six-local}
Let $\mathbb K$ be totally imaginary and suppose that the Fermat cubic
has a point $(a:b:c)\in\mathbb P^2(\mathbb K)$ with $abc\ne0$.
At each triadic place of width five in Theorem~\ref{nf3:triadic-class},
suppose that the greatest and least of $v(a),v(b),v(c)$ differ by at
most one. Then $G_{\mathbb K}(3)\le6$. In particular no condition
on the point at finite places is needed if every triadic completion
has width at most four. This includes $\mathbb Q(\sqrt{-5})$.
\end{theorem}
\begin{proof}
Clear denominators in $(a,b,c)$ and choose $(u,v,w)$ in the integral
kernel of $a^2u+b^2v+c^2w$. Put
$B=u^3+v^3+w^3$ and $D=au^2+bv^2+cw^2$. Then
\[
 (aT+u)^3+(bT+v)^3+(cT+w)^3=B+3DT .
\]
We choose a tangent vector for which $D\ne0$ and $\nu-B$ is a sum
of three cubes modulo $3D$. We may take $a,b,c$ locally primitive
at all places above $3$, by weak approximation on $\mathbb K^\times$
followed by clearing the remaining denominators by a rational integer
prime to three.

At a place of width three there is nothing to prescribe, since
$\nu-B$ remains eligible. At a remaining triadic place the residue
field is $\mathbb F_3$. If the three coordinates of the point are
units, the local tangent $(b^2,-a^2,0)$ gives $D=-abc^3$ a unit;
Frobenius then represents $\nu-B$ by one cube modulo $3D$.
Otherwise exactly two of the primitive coordinates are units.
Relabel them $a,b$ and put $h=v(c)>0$, $A=-(c/b)^2$. The local
tangent $(0,Az,z)$ for any unit $z$ satisfies
\[
 v(D)=v((bA^2+c)z^2)=h,\qquad
 B=(1+A^3)z^3\equiv((1+A)z)^3\pmod{3\pi^h}.
\]
Indeed $bA^2+c=c(1+(c/b)^3)$, and
$(1+A)^3-1-A^3=3A(1+A)$ has valuation $e+2h$.
At a place of width four represent $\nu$ by four cubes with a unit
root and choose that root as $(1+A)z$. At a place of width five
our hypothesis gives $h=1$. All eligible targets are sums of four
cubes with a unit modulo $3\pi$: with $w^3\equiv\nu\pmod{3A_v}$
in the local integer ring $A_v$, choose the four residues for the
integer $M\in[-4,4]$ used in the proof of
Theorem~\ref{nf3:triadic-class}. For $M=\pm4$ use four equal
unit residues. The fixed-sum polynomial then gives the congruence
modulo $3\pi$ directly. This again supplies the required $z$.
A fixed finite precision suffices in these choices: preserve $D$
modulo $\pi^{h+1}$ and $B$ modulo $3\pi^h$. Thus one needs only
finitely many simultaneous tangent residue conditions as $\nu$ varies.

The integral kernel lattice realizes these finitely many simultaneous
residue conditions. Indeed completion gives the kernel of the same
linear equation over each local ring by flatness, and the Chinese
remainder theorem applies to any finitely generated module.
Adding a sufficiently divisible integral multiple of $(b^2,-a^2,0)$
ensures $D\ne0$, since its quadratic value is $-abc^3\ne0$.

Consider any remaining prime divisor of $D$ and write $q$ for its
residue norm. If $q>4$, the multiplicative group of nonzero cubes
$H$ generates the residue field over its prime field: its size is at
least $(q-1)/3>\sqrt q-1$. For $C=H\cup\{0\}$, the additive
stabilizer of $3C$ is $H$-invariant, so is either zero or the whole
field. In the former case Kneser's theorem \cite{edkneser} gives
$|3C|\ge3|C|-2\ge q$ as well. For $q=2$ every residue is a cube;
for $q=4$ the reduction of the eligible element $\nu-B$ belongs to
$\mathbb F_2$ and is a cube. In every case a nonzero coordinate
is available, using $(1,-1,0)$ when the residue is zero. Hensel's
lemma lifts the three cubes since the residue characteristic differs
from three.

Take the three approximating roots for the translated target
$\nu-B$ in their prescribed classes modulo $3D$. Rounding in these
fixed cosets gives a remainder $3DT$ with
$|\sigma(T)|\ll_{\mathbb K,a,b,c} P^{8/9}$. There are only finitely
many choices of the modulus and of $B$, so after removing the
balancing unit all six roots satisfy
$|x_v|\ll_{\mathbb K,a,b,c}|\nu_v|^{1/3}$ uniformly in $\nu$.
Translates of the third root avoid zero as
before. For the example take $(a,b,c)=(4+\sqrt{-5},4-\sqrt{-5},-2)$.
\end{proof}

For a nonzero totally positive eligible $\nu$ and fixed $C>0$,
let $r_{s,C}(\nu)$ count the ordered representations by $s$ cubes
of nonzero totally positive integers $x_i$ satisfying
\[
       |x_{i,v}|\le C|\nu_v|^{1/3}\qquad(i\le s,\ v\mid\infty).
\]
This count is finite by discreteness of the integer lattice and is
invariant under multiplication of the target by a totally positive
unit cube. After balancing $\nu$, the proofs give $|x_{i,v}|\ll
N^{1/(3d)}$ for $N=|\mathrm N\nu|$; undoing the same unit gives the
displayed bounds. A field-norm bound alone may admit infinitely many
tuples. In $\mathbb Q(\sqrt[3]2,i)$, the nontorsion unit
$e=\sqrt[3]2-1$ gives the six-root tuples
$(m,1,1,1,e^j,-e^j)$ of target $m^3+3$ for every $j\ge1$.
Appending $1$ gives the analogous seven-root family.

\begin{corollary}\label{nf3:sevenmultiplicity}
There is a fixed $C=C_{\mathbb K}$ such that, for every totally
positive eligible $\nu$ of sufficiently large $N=|\mathrm N\nu|$,
the number $r_7(\nu)=r_{7,C}(\nu)$ satisfies
\[
r_7(\nu)\gg_{\mathbb K,\epsilon}
\begin{cases}
 N^{2/3-\epsilon},&\mathbb K\text{ totally real},\\
 N^{1/3}(\log N)^{-10d-1},&\mathbb K\text{ of mixed signature},\\
 N^{1/12},&\mathbb K\text{ totally imaginary}.
\end{cases}
\]
\end{corollary}
\begin{proof}
For fixed $u,v,w$ Mitsui's theorem gives
$\gg(P/U^4)^d/(U^{6d}\log N)$ choices of the prime center.
The first two ordered roots in \eqref{nf3:pairidentity} recover
$s$. The pair differences then recover $L_1,y,z$ and hence
$x,y,\zeta$. This proves the mixed-signature estimate.

For totally real $\mathbb K$ the quantitative ternary estimate in
the proof of \cite[Corollary~5.6]{nf3:schulzesurvey} gives
$r_L(Q)\gg_{\mathbb K,\epsilon}(\mathrm NQ)^{1/2-\epsilon}$.
We spell out its use here. Local primitive representations exist:
at $2,11$ the target is a represented unit, and at other odd places
the unimodular ternary lattice contains a hyperbolic plane.
For a genus or spinor genus $\mathcal G$, write
\[
 r_{\mathcal G}(t)=
 \frac{\sum_{L_i\in\mathcal G}r(t,L_i)/|\operatorname O(L_i)|}
      {\sum_{L_i\in\mathcal G}1/|\operatorname O(L_i)|},
\]
with one representative of each class. The genus average has the
required lower estimate by \cite[Lemma~2.1]{nf3:schulzesurvey}.
The spinor-genus averages coincide at $Q$. The only possible nontrivial contributing character
is that of $\mathbb K(\sqrt{-Q\det q})$, whereas every spinor-genus
character is trivial on the full group of local units at
$\mathfrak l$. Indeed the stabilizer of a vector of norm $Q$ is
the corresponding binary rotation torus, and unfolding the character
average enforces triviality on its spinor norms, the norms from
this quadratic extension. The error for an individual class relative
to its spinor-genus average
is of smaller order, uniformly in the squarefree and square parts
of the index, as follows. Write $(Q)=\mathfrak t\mathfrak m^2$
with $\mathfrak t$ squarefree and $\mathfrak m=(\mu)\mathfrak c$
with $\mathfrak c$ among the representatives fixed in the proof
of Theorem~\ref{nf3:realG7}. At every prime dividing a
representative or $2\cdot3\cdot11\,\operatorname{disc}(\mathbb K)$,
$Q$ is a unit. Thus
the primes of $\mathfrak m$ are good. Put $t'=Q/\mu^2$. We have
$(t')=\mathfrak t\mathfrak c^2$ and
$r(Q,L_i)=r(t',\mathfrak m^{-1}\mathfrak c L_i)$ for each
genus class $L_i$. The index $t'$ is squarefree relative to
$\mathfrak c L_i$. It is locally represented there. Away from
$\mathfrak m$ we can divide the original local representation by
$\mu$, and at primes in $\mathfrak m$ the unimodular ternary
lattice contains a hyperbolic plane. The same marked prime $\mathfrak l$ rules out
a nontrivial spinor contribution there, since $\mathfrak c$ is a
unit there and $t'$ and $Q$ differ by a square.

At this base index \cite[Theorem~5.2]{nf3:schulzesurvey} bounds
each class coefficient minus its common genus coefficient by
$O_\epsilon((\mathrm N t')^{123/260+\epsilon})$. Only finitely many
base lattices occur. Apply \cite[Lemma~2 and pp.~120--122]{nf3:brandt}
to the vector of unweighted individual coefficients
$(r(t',\mathfrak cL_i))_i$. At $\mathfrak p^h\parallel\mathfrak m$
the matrix in this recursion is
$\overline B(\mathfrak p^h)-\chi_{t'}(\mathfrak p)
\overline B(\mathfrak p^{h-1})$ when the relative valuation of
$t'$ is zero, and $\overline B(\mathfrak p^h)$ when it is one.
Here $\chi_{t'}$ is the local quadratic character of the binary
stabilizer field. Each matrix has constant row sum, since the prime
neighbor count is $q+1$ and the prime-power recursion and
multiplicativity preserve constant vectors. We may therefore
subtract the common coefficient before applying these row
differences within a spinor genus. Squarefreeness excludes the
third, imprimitive-index term of Lemma~2. The remaining divisor
sums cost $(\mathrm N\mathfrak m)^\epsilon$ and the
cuspidal row differences cost
$O_\epsilon((\mathrm N\mathfrak m)^{7/8+\epsilon})$. The
class-minus-spinor error, obtained by averaging these differences,
is thus
\[
 O_\epsilon\big((\mathrm N t')^{123/260+\epsilon}
                 (\mathrm N\mathfrak m)^{7/8+2\epsilon}\big)
   \ll_\epsilon(\mathrm N Q)^{123/260+O(\epsilon)},
\]
since $7/8<2(123/260)$ and the $\mathrm N\mathfrak c$ are bounded.
Choosing the slack small gives the claimed lower bound. Finally
$\mathrm NQ\asymp N^{2/3}/U^{2d}$, proving the assertion after
multiplication by the number of centers.

In the totally imaginary construction allow the first root to vary
within distance $P^{1/4}$ of the chosen cube root at all complex
places, retaining its fixed residue class. Lattice counting gives
$\gg P^{d/4}$ choices. The successive remainder bounds now have
exponents $9/4,3/2,1$ in $P$. The fixed translates used to avoid
zero still give seven roots of size $O_{\mathbb K}(P)$, and the
first ordered coordinate distinguishes all these representations.
\end{proof}

\begin{proposition}\label{nf3:imag-seven-many}
In the totally imaginary case of Corollary~\ref{nf3:sevenmultiplicity}
one has $r_7(\nu)\gg_{\mathbb K,\epsilon}N^{1/6-\epsilon}$.
If $\zeta_3\notin\mathbb K$ then
$r_7(\nu)\gg_{\mathbb K}N^{1/4}/(\log N)^2$.
\end{proposition}
\begin{proof}
Write $P=N^{1/(3d)}$ after unit balancing. For $u\equiv1\pmod2$
in $\mathcal O_{\mathbb K}$ set $s=(u+1)/2$ and $r=(u-1)/2$.
The last four roots can now vary with a modulus, by the identity
\[
 (t+s)^3+(t-s)^3+(-t+r)^3+(-t-r)^3=6ut .
\]
Keep the prescribed three root classes modulo $6$ of
Theorem~\ref{nf3:realG7}. Choose the first root $x_1$ within bounded
distance of the cube roots of the target, and put $B=\nu-x_1^3$.
Then $|B_v|\ll P^2$. Take a prime ideal $(u)$ in the identity
ray class modulo $2$, coprime to $6$, with balanced generator
$u\equiv1\pmod2$ and norm $q\asymp P^{d(1/3-\tau)}$,
where $0<\tau<1/6$. There are
$\gg P^{d(1/3-\tau)}/\log P$ choices by the prime ideal theorem
in ray classes. Balance generators using the finite-index subgroup
of units congruent to one modulo the ray modulus.

Among the $x_2$ in its fixed residue class in a complex box of
radius $W=P^{1/6}$ about chosen cube roots of $B_v$, there are
$\gg W^d$ for which $B-x_2^3$ is a cube modulo $(u)$.
Indeed if $q\equiv2\pmod3$ or $B\equiv0\pmod u$, every residue
works. Otherwise the cube indicator on nonzero residues is
$(1+\chi+\overline\chi)/3$ for a cubic character $\chi$.
For $B\not\equiv0$ all complete sums of $\chi(B-z^3)$ against
additive characters are $O(\sqrt q)$ by the mixed Weil bound:
the three zeros of $B-z^3$ are simple. The same holds for
$\overline\chi$. Poisson summation with a smooth bump inside the
box and the fixed class modulo $6$ gives total error $O(\sqrt q)$.
In fact the Fourier transform on the dual of $(6u)$ has absolute
sum $O(q/W^d)$ after scaling by $W$, by balance and rapid decay;
its prefactor is $O(W^d/q)$. The constant term has size
$\gg W^d$, and $\sqrt q=o(W^d)$. Zero residues may be included.
The remaining target has size $O(P^{3/2})$ at every embedding.
Using the fixed class modulo $6$ and a cube root modulo $u$,
round the third root in the resulting coset within $O(|u|_v)$
at each place. It has size $O(P^{1/2})$, and the new remainder
$6ut$ has size $O(P|u|_v)$ by expansion. Thus all seven roots
are $O(P)$. A bounded number of translates of the third root by
$6u j$ avoid zero and the four possible values of $t$ causing a
zero among the last four roots. The four ordered roots recover
$u$ and $t$. Counting $u$ and $x_2$ proves the first claim.

When $\zeta_3\notin\mathbb K$, take $u$ generating the product
of two distinct primes with norms $\asymp P^{d/4}$ that are
$2$ modulo $3$. Their product is required to lie in the principal
ray class modulo $2$. In the compositum of the ray class field
and $\mathbb K(\zeta_3)$ take an element restricting nontrivially
to the latter, and its inverse. Chebotarev supplies
$\gg P^{d/4}/\log P$ choices for each prime in the indicated
classes, in disjoint norm windows of comparable sizes.
Consequently there are $\gg P^{d/2}/(\log P)^2$ choices of
the product, each with a balanced generator $u\equiv1\pmod2$.
Cubing is bijective modulo $u$ by the Chinese remainder theorem.
Vary the first root in radius $P^{1/4}$ and round the second
within bounded distance. The two remainders have sizes
$O(P^{9/4})$ and $O(P^{3/2})$. The third congruence is automatic
for every product, with rounding error $O(P^{1/2})$ and hence
again $|t_v|\ll P$. The first coordinate distinguishes
$\gg P^{d/4}$ choices independently of the product, proving
the second bound.
\end{proof}

The same counting device gives quantitative six-cube results
under the hypotheses of Theorem~\ref{nf3:fermat-six-local}.
For a fixed sufficiently large $C=C_{\mathbb K,a,b,c}$ and
$r_6(\nu)=r_{6,C}(\nu)$ one has
\[
 r_6(\nu)\gg_{\mathbb K,a,b,c,\epsilon}N^{7/48-\epsilon},
 \qquad
 r_6(\nu)\gg_{\mathbb K,a,b,c}N^{1/6}/(\log N)^2
                 \quad(\zeta_3\notin\mathbb K).
\]
Indeed that proof chooses a tangent vector $(u,v,w)$ and constants
$B,D$ from fixed finite lists. Multiplying this vector by $s$
replaces its line identity by
\[
 (at+su)^3+(bt+sv)^3+(ct+sw)^3=B s^3+3D s^2 t .
\]
Take $s\equiv1\pmod{3D}$. For the first estimate let $(s)$ be
prime in this ray class with norm $\asymp P^{d(1/4-\tau)}$.
Choose the three fixed root classes modulo $3D$ for $\nu-B$
as in that theorem. Vary the first root about the cube root of
the target in radius $P^{1/16}$, then the second about the cube
root of the remainder in radius $P^{1/8}$. The remainders have
sizes $O(P^{33/16})$ and $O(P^{3/2})$. The completion just used,
with $W=P^{1/8}$, gives $\gg W^d$ second roots for which the
new remainder is a nonzero cube modulo $s$: the at most three
forbidden residue classes have total weight $O(1+W^d/|\mathrm Ns|)$
each by balance and packing. If the first remainder vanishes modulo
$s$ exclude just its zero residue for $x_2$. Lift the nonzero cube
root modulo $s^2$ by Hensel's lemma. Round the third root of the
second remainder in its fixed class modulo $3D$ and chosen class
modulo $s^2$. Its error is $O(|s_v|^2)$ and $|s_v|^2\ll P^{1/2}$.
After subtracting $B s^3$, division by $3D s^2$ leaves $|t_v|\ll P$.
The independent choices of $s,x_1,x_2$ give the first estimate.
The last three ordered coordinates recover $s,t$, since the tangent
and Fermat vectors are linearly independent when $D\ne0$.

For the second estimate take $(s)$ to be a product of two primes,
each of norm $\asymp P^{d/8}$ equal to $2$ modulo $3$, in inverse
Artin classes for the ray modulus $3D$ as before. Vary the first
root within $P^{1/4}$ and round the second within bounded distance,
choosing a translate that makes its remainder nonzero. For each
such pair exclude the bounded number of prime factors at these
scales dividing that remainder, whose norm is $O(P^{3d/2})$.
Every other product supplies a nonzero cube root modulo $s^2$
and has $|s_v|^2\asymp P^{1/2}$. There are
$\gg P^{d/4}/(\log P)^2$ products, independently of the
$\gg P^{d/4}$ first roots. Fixed translates of the third root
by multiples of $3D s^2$ avoid zeros in either construction.

\subsection{The exact density of four positive cubes}
\label{sec:four-cube-exact-density}
The local classification gives the density in Theorem~\ref{nf3:four}
for every number field. Let $m_5(\mathbb K)$ count triadic places
of residue degree one satisfying either $e\equiv1\pmod3$, $e\ge4$,
or $\mathbb K_v\simeq\mathbb Q_3(\sqrt{-3})$.

\begin{theorem}\label{thm:four-cube-exact-density}
In every number field $\mathbb K$, the limiting proportion of eligible
targets representable by four totally positive cubes is
\[
                             (23/27)^{m_5(\mathbb K)}
\]
in dilates of bounded Jordan-measurable regions of positive volume in
the totally positive archimedean region. In particular four cubes
represent almost all eligible targets if and only if
$m_5(\mathbb K)=0$. On compact boxes bounded away from the coordinate
walls the exceptions to the finite-place criterion have a power saving.
\end{theorem}
\begin{proof}
Only the places counted by $m_5$ can exclude a target locally, by
Corollary~\ref{nf3:all-local}. At such a place write $A$ for the
valuation ring and $\pi$ for a uniformizer. The image of the
locally eligible targets modulo $3A$ is the Frobenius image there.
Choose a root $w$ for each such class, and put
$t=(\nu-w^3)/3$. With normalized Haar probability on the local
eligible group, $\bar w\in\mathbb F_3$ is uniform, since Frobenius
on residue classes is bijective, and $\bar t$ is independent and
uniform. The integer $M\in[-4,4]$ determined by
$M\equiv w_0^3+3\bar t\pmod9$, where $w_0=0,1,-1$ lifts $\bar w$,
is consequently uniform in these nine values. The fixed-sum argument
in Theorem~\ref{nf3:triadic-class} gives primitive four-cube
representations whenever $|M|\le3$.

For $M=\pm4$ the root $w$ is a unit. Dividing by $w^3$ puts the
target in the form $4+3\pi u$, with $u$ Haar-uniform in $A$ within
the specified residue classes. The obstruction in the same theorem
excludes it when $u$ is a unit. Conversely suppose $u=\pi y$.
Put $x_i=1+\pi z_i$ and $\sum_{i=1}^4 z_i=0$. Expanding gives
\[
    \sum x_i^3=4+3\pi^2\left(\sum z_i^2+
                                    \pi\frac{\sum z_i^3}{3}\right).
\]
After eliminating $z_4$ the quotient in parentheses is an integral
polynomial. On the residue hyperplane $\sum z_i=0$, the tuples
$(0,1,1,1),(1,1,-1,-1),(1,-1,0,0)$ give the three possible values
of $\sum z_i^2$. In each tuple some $z_i-z_4$ is a unit. Thus
the parenthesis attains $y$ by Hensel's lemma. The excluded local
proportion is exactly $(2/9)(2/3)=4/27$.

All these conditions depend on finitely many residue digits. The
global cube group surjects onto the product of its images at the
finitely many prescribed prime powers, since the Chinese remainder
theorem realizes cube roots supported on one factor at a time.
Lattice counting gives the asserted product for the everywhere locally
representable set. Theorem~\ref{nf3:four} gives the global conclusion,
and exhaustion away from the archimedean walls extends the density
assertion to the indicated regions.
\end{proof}

The proportion relative to all totally positive algebraic integers
is obtained by dividing by the index of $J_3$ in $\mathcal O$.
This index is
\begin{equation}\label{eq:cubic-order-index}
 [\mathcal O:J_3]=2^{\#\{v\mid2:f_v=2\}}\,
       3^{\sum_{v\mid3} f_v(e_v-\lceil e_v/3\rceil)} .
\end{equation}
Indeed $6\mathcal O\subset J_3$, and its image modulo six is the
product of the corresponding local cube spans. At $v\mid3$, cubing
on $A/3A$ is an additive endomorphism whose kernel is
$\pi^{\lceil e_v/3\rceil} A/\pi^{e_v}A$. At $v\mid2$ all the
principal units are cubes. The only possible loss is in the residue
field $F$, where the span of the nonzero cubes is the subfield they
generate. For $q=|F|>4$ the group has size at least
$(q-1)/3>\sqrt q-1$, hence generates $F$; at $q=2$ it also spans
the residue field. At $q=4$ its span is $\mathbb F_2$.
This proves \eqref{eq:cubic-order-index}.

For three summands the finite-place condition need not account for all
global exceptions. Nevertheless its precise local density can also be
found in a family of arbitrarily ramified completions.
\begin{proposition}\label{prop:three-cube-wild-density}
Let $A$ have residue field $\mathbb F_3$ and ramification index
$e\equiv1\pmod3$, $e\ge7$. The local Haar proportion of three-cube
sums in the cubic order is $17/27$. The same holds for $e=4$ if
$\pi^4/3\equiv1\pmod\pi$. Every represented target has a
representation with a unit among the roots.
\end{proposition}
\begin{proof}
Use the nine equiprobable values of $M$ from the preceding proof.
For $|M|\le2$ one can choose three residues of sum $M$ with a zero
and a unit, using $(1,-1,0)$ at zero. The fixed-sum cubic polynomial
then lifts by its unit derivative. If $\nu=\sum x_i^3$, write
$x_i=a_i+\pi y_i$, $a_i=0,1,-1$, and $c=\sum_i a_i$.
As in \eqref{nf3:four-triadic}, modulo $3\pi$ one has
$\nu\equiv c+(\pi\sum y_i)^3$.
The sum $w'=c+\pi\sum y_i$ is a cube root of $\nu$ modulo three.
Any two such roots agree modulo $\pi^l$, $l=\lceil e/3\rceil$,
and their cubes agree modulo $3\pi$. Furthermore
$(w')^3\equiv c^3+(\pi\sum y_i)^3\pmod{3\pi}$. Hence
$M\equiv w_0^3+c-c^3\equiv c\pmod9$, since $c\equiv w_0\pmod3$.
Thus $|M|=4$ is impossible. At $M=3$ the three residues must all
equal one; negation deals with $M=-3$.

Choose $w\in\pi A$ with $w^3\equiv\nu\bmod3A$ and put
$x_i=1+\pi y_i$, $z=\pi\sum y_i$. Then $z=w+\pi^l b$, and
\[
 \nu=3+z^3+3z+3\pi^2
       \left(\sum y_i^2+\pi\frac{\sum y_i^3-(\sum y_i)^3}{3}\right).
\]
Necessarily $\nu\equiv3+w^3+3w\pmod{3\pi^2}$ since $3l=e+2$
and $l\ge2$. This selects exactly one third of the class $M=3$.
Conversely, after imposing this congruence, varying $b$ changes
$(z^3+3z-w^3-3w)/(3\pi^2)$ modulo $\pi$ by a unit times
$\bar b^3$ for $l\ge3$. For $l=2$ the expression is
$(\pi^4/3)\bar b^3+\bar b$, again bijective under the hypothesis.
Choose $b$ so the parenthesis in the displayed equation is to have
residue two. With any prescribed sum $\sum\bar y_i$ this is attained
nonsingularly: the triples $(1,-1,0),(-1,-1,0),(1,1,0)$ cover
the three sums. After elimination of the last coordinate a derivative
$2(\bar y_i-\bar y_3)$ is nonzero. Hensel lifting finishes the
representation. The total measure is $5/9+2/27=17/27$.
\end{proof}

\begin{proposition}\label{nf3:small-three-densities}
Suppose the residue field is $\mathbb F_3$. For $e=2$, the local
proportion of the cubic order represented by three cubes is $7/9$
if $\pi^2/3\equiv1\pmod\pi$, and $17/27$ if
$\pi^2/3\equiv-1\pmod\pi$. For $e=4$ and
$\pi^4/3\equiv-1\pmod\pi$ the proportion is $151/243$.
Every represented target has a primitive representation.
\end{proposition}
\begin{proof}
Use the classes $M$ of Proposition~\ref{prop:three-cube-wild-density}.
For clarity, $w^3\bmod3A$ fixes $w\bmod\pi^{\lceil e/3\rceil}$.
If $w$ is changed by an element $\delta$ of this last ideal, then
\[
 (w+\delta)^3-w^3=3w^2\delta+3w\delta^2+\delta^3
                      \in3\pi A
\]
for $e=2$ or $4$, since $3\lceil e/3\rceil>e$.
Thus $M$ is intrinsic. Its nine values are equiprobable. The five
with $|M|\le2$ have primitive solutions, and $|M|=4$ is impossible
with three roots. By negation it suffices to discuss $M=3$, when all
three roots reduce to one. Write them as $1+\pi u_i$, and put
$z=\pi\sum_i u_i$. The exact equation is
\begin{equation}\label{nf3:three-small-exp}
 n=3+z^3+3z+3\pi^2\left(\sum_i u_i^2+\pi C(\mathbf u)\right),
 \qquad C(\mathbf u)=
   \left(\sum_i u_i^3-\left(\sum_i u_i\right)^3\right)/3 .
\end{equation}

For $e=2$ take $w=0$, since its residue is zero, and put $U=\pi^2/3$.
Then $n=3+3\pi v$, with $v$ uniform in $A$. Equation
\eqref{nf3:three-small-exp} requires
$\bar v=\bar b+\bar U\bar b^3$ with $b=z/\pi$.
Conversely take the three roots $1,1,1+\pi b$. The exact equation
for $b$ is $v=b+\pi b^2+Ub^3$, whose derivative is a unit.
It follows that all targets of this class occur for $\bar U=1$,
and one third occur for $\bar U=-1$. This gives the two quadratic
proportions. The second includes the completion used in
Proposition~\ref{nf3:amplify}, since $-2$ is a square in $\mathbb Q_3$.

For $e=4$, fix a representative $w$ modulo $\pi^2$ of the cube-root
class, and write $U=\pi^4/3$. In the class $M=3$ one has $\pi\mid w$
and $z=w+\pi^2 b$. Equation \eqref{nf3:three-small-exp} requires
$n\equiv3+w^3+3w\bmod 3\pi^2A$, which retains one third of the class.
Write $c=(n-3-w^3-3w)/(3\pi^2)$. The exact system becomes
\[
 c=b+Ub^3+w^2b+w\pi^2 b^2+\sum_i u_i^2+\pi C(\mathbf u),
 \qquad u_3=w/\pi+\pi b-u_1-u_2 .
\]
Since $\bar U=-1$, the residue condition is
$\bar c=\sum_i\bar u_i^2$ with $\sum_i\bar u_i=s=\overline{w/\pi}$.
Every such solution lifts: after eliminating $u_3$ the partial
derivative in $b$ is one modulo $\pi$. If $s\ne0$, all three square
sums occur. For $s=1$ take $(1,1,-1),(1,0,0),(-1,-1,0)$, and
negate for $s=-1$. For $s=0$ precisely $0,2$ occur. The cube-root
class makes $s$ uniform in $\mathbb F_3$; conditional on that class
and the preceding congruence, $\bar c$ is independently uniform.
The surviving fraction at this last step is $8/9$. Hence the
total proportion is $5/9+(2/9)(1/3)(8/9)=151/243$.
\end{proof}

\begin{proposition}\label{nf3:three-density-last}
Suppose the residue field is $\mathbb F_3$ and $e\equiv2\pmod3$,
$e\ge5$. The proportion of the cubic order represented by three
integral cubes is
\[
\begin{cases}
61/81,&e=5,\\
181/243,&e=8,\ \pi^8/3\equiv-1\pmod\pi,\\
545/729,&e=8,\ \pi^8/3\equiv1\pmod\pi,\ \text{or }e\ge11 .
\end{cases}
\]
All represented targets have a unit root.
\end{proposition}
\begin{proof}
Put $l=(e+1)/3$ and $U=\pi^e/3$. Again the classes $|M|\le2$ are
represented and $|M|=4$ are excluded. It suffices to treat $M=3$.
Fix $w\bmod\pi^l$ as in the preceding proof, with $\pi\mid w$.
Then $V=(n-3-w^3-3w)/(3\pi)$ is uniform in $A$.
Write the roots $1+\pi u_i$ and put
$z=\pi\sum_i u_i=w+\pi^l b$. The last equality follows by cubing
modulo $3A$. For $a\in A$ define
\[
  L_a(b)=U b^3+\pi^{l-1}(a^2+w^2)b+\pi^{2l-1}w b^2 .
\]
Equation \eqref{nf3:three-small-exp} becomes
\[
  V=L_1(b)+\pi\sum_i u_i^2+\pi^2C(\mathbf u).
\]
Its first digit fixes $b$ modulo $\pi$ by Frobenius. Choose a lift
$b_0$ and put $c=(V-L_1(b_0))/\pi$ and $s=\overline{w/\pi}$.
Since $l\ge2$, changing $b$ by $\pi h$ preserves $L_1(b)$ modulo
$\pi^2$. The next digit therefore requires
$\sum\bar u_i^2=\bar c$ with $\sum\bar u_i=s$.
When $s\ne0$ all square sums have a nonconstant triple, as in the
preceding proof. For $s=0,\bar c=2$ use $(1,-1,0)$.
Holding $b=b_0$ fixed and eliminating $u_3$ gives a unit derivative
$2(\bar u_i-\bar u_3)$ after permutation. These cases occupy $7/9$
of $M=3$. The case $s=0,\bar c=1$ is impossible and occupies $1/9$.
In the remaining $1/9$ one has $s=\bar c=0$, so the three $u_i$ have
a common residue $d$ and $w\in\pi^2 A$. These fractions follow from
the uniformity of the section $w\bmod\pi^l$ and of the successive
digits of $V$ conditional on that section.

For $e=5$ we have $l=2$ and can choose $w=0$ in the remaining
case. Set $u_1=u_2=0,\ u_3=\pi b$. It suffices to solve
\[
                     V=U b^3+\pi b+\pi^3 b^2.
\]
At $b=b_0$ this congruence holds modulo $\pi^2$ by $\bar c=0$.
Put $b=b_0+\pi h$, subtract its value at $b_0$, and divide by
$\pi^2$. The resulting polynomial reduces to $h$ modulo $\pi$.
Thus the whole remaining branch is soluble.

Suppose next that $l\ge3$. Lift $d$ to $\{0,1,-1\}$ and write the
roots as $a+\pi^2 v_i$, with $a=1+\pi d$. Put
$Z=\pi^2\sum_i v_i=w+\pi^l b$; replacing $z$ by $Z$ subtracts
$3\pi d$, which does not change the residue $b\bmod\pi$.
Direct expansion gives
\begin{align*}
 V&=D(d)+L_a(b)+\pi^3a\sum_i v_i^2+\pi^5C(\mathbf v),\\
 D(d)&=3d+3\pi d^2+\pi^2d^3+(2d+\pi d^2)w,\qquad
 \sum_i v_i=w/\pi^2+\pi^{l-2}b .
\end{align*}
The preceding two digits give $V\equiv U b_0^3\bmod\pi^2$.
Put $t=\overline{w/\pi^2}$. At the next digit $D(d)/\pi^2$
reduces to $(1-t)d$. The term $\pi^{l-1}(a^2+w^2)b/\pi^2$
contributes the already determined $\bar b_0$ for $l=3$ and zero
otherwise. Thus $t\ne1$ allows every target digit, whereas $t=1$
allows one in three. Here $t$ is uniform since $l\ge3$, and the
new target digit is uniform conditional on all preceding choices.
The surviving fraction of this branch is $7/9$ at this stage.

For an allowed $d$ put $b=b_0+\pi h$. Subtract $D(d)+L_a(b_0)$
from $V$ and divide by $\pi^3$. The residue of the resulting equation
has the form
\[
 c'=\bar U\bar h^3+
       {\boldsymbol 1}_{l=3}\bar h+\sum_i\bar v_i^2,\qquad
                        \sum_i\bar v_i=t .
\]
All divisions here are integral polynomial expansions after the
three fixed digits have been subtracted. In particular the mixed
terms in $U(b_0+\pi h)^3$ beyond $U\pi^3h^3$ vanish after
reduction; the term with coefficient $\pi^{l-1}$ contributes
$\pi^{l-3}(a^2+w^2)h$. For $l\ge4$, and also for $l=3,\bar U=1$,
choose any nonconstant triple of sum $t$, fit $h$ to the remaining
digit, and lift in one of the $v_i$ using the unit square derivative.
For $l=3,\bar U=-1$, the first two terms cancel on $\mathbb F_3$,
and their polynomial derivative in $h$ is one. Eliminating $v_3$
does not alter this derivative modulo $\pi$. Consequently any
triple with the required square sum lifts. The only additional
exclusion is $t=0$ with required square sum one, since all square
sums occur for $t\ne0$. When $t=0$ the preceding digit determines
$d$ uniquely, so the excluded fraction of the branch is exactly
$1/9$. Thus the survival fraction of the last branch is $1$ for
$e=5$, $2/3$ for $e=8,\bar U=-1$, and $7/9$ otherwise.
Combining the five initial classes and the two classes $M=\pm3$
gives the three displayed proportions.
\end{proof}

There is an independent finite calculation for each of these tame
completions. By multiplying a uniformizer by a principal unit we
may assume $\pi^e=\pm3$, since $3\nmid e$. Put
$r=\lfloor e/2\rfloor+1$. On a unit root, changing $x$ to
$x+\pi^r z$ changes its cube by $3\pi^r$ times a polynomial
congruent to $x^2 z$ modulo $\pi$. Hence any primitive representation
modulo $\pi^{e+r}$ lifts. If a three-cube residue has a nonprimitive representation at that
precision, then $M=0$, and the construction using both zero and unit
residues furnishes a primitive representation. To calculate the cube set at this
precision it suffices to run through $x\bmod\pi^r$. Appending a
digit $d\pi^j$ updates $x^3$ by
$3x^2d\pi^j+3x d^2\pi^{2j}+d^3\pi^{3j}$. Take its threefold
sumset in $J_3(A)/\pi^{e+r}A$, of order
$3^{r+\lceil e/3\rceil}$. At $e=8$ this gives $4905$ and $4887$
residues out of $6561$ for the positive and negative radicands.
At $e=11$ it gives $44145$ out of $59049$, and at $e=14$ both
signs give $1191915$ out of $1594323$.

\subsection{Local densities for three summands}
\label{sec:local-three-table}
Write $d_v$ for the Haar proportion of three-cube sums in the cubic
order at $v$. At $v\nmid3$ one has $d_v=1$ by
Corollary~\ref{nf3:all-local}. For $v\mid3$, write $e,f$ for its
ramification index and residue degree, put $u=\pi^e/3$ and, when
$f=1$, regard $\bar u$ as a sign. The sign is independent of the
uniformizer at the even indices where it occurs in the table.
Theorem~\ref{nf3:triadic-class} and
Propositions~\ref{prop:three-cube-wild-density},
\ref{nf3:small-three-densities} and~\ref{nf3:three-density-last} give
\[
\begin{array}{c|c}
 \text{triadic completion}&d_v\\\hline
 f\ge2\text{ or }3\mid e&1\\
 f=1,\ e=1\text{ or }(e=2,\bar u=1)&7/9\\
 f=1,\ (e=2,\bar u=-1)\text{ or }(e=4,\bar u=1)&17/27\\
 f=1,\ e=4,\bar u=-1&151/243\\
 f=1,\ e\equiv1\pmod3,\ e\ge7&17/27\\
 f=1,\ e=5&61/81\\
 f=1,\ e=8,\bar u=-1&181/243\\
 f=1,\ (e=8,\bar u=1)\text{ or }(e\equiv2\pmod3,\ e\ge11)&545/729 .
\end{array}
\]
In degree one the three-cube residues are the seven sums of three
members of $\{0,1,-1\}$ modulo nine; a unit cube absorbs the
remaining multiple of nine.

\begin{corollary}\label{cor:three-local-density-degree}
The relative density among eligible targets which are sums of three
cubes at every finite place is $\Delta_3(\mathbb K)=\prod_{v\mid3}d_v$.
Among fields of degree $n$ its least value is $(7/9)^n$, attained
precisely when $3$ splits completely. Among degree-$n$ fields the
least relative density of global four-cube sums in the regions of
Theorem~\ref{thm:four-cube-exact-density} is
$(23/27)^{\lfloor n/2\rfloor}$. Both minima are realized by totally
real fields, possibly different for the two assertions.
\end{corollary}
\begin{proof}
The three-cube congruence conditions in the table are clopen in the
local cubic order. The Chinese remainder argument of
Theorem~\ref{thm:four-cube-exact-density} therefore gives the product.
For every completion of local degree $ef\ge2$ one has
$d_v\ge151/243>(7/9)^2\ge(7/9)^{ef}$. Each degree-one place gives
$7/9$ and the sum of the local degrees is $n$. For four cubes,
each width-five place has local degree at least two.

To realize the first minimum prescribe the split algebra
$\mathbb Q_3^n$. For the second prescribe $\lfloor n/2\rfloor$
copies of $\mathbb Q_3(\sqrt{-3})$ and one rational factor when
$n$ is odd. Choose over $\mathbb Q_3$ distinct translated minimal
polynomials for the prescribed factors. Their factorization type and
the isomorphism types of the factors persist on a neighborhood by
Hensel's and Krasner's lemmas. Approximate the coefficients of this
monic polynomial over $\mathbb Q$, simultaneously prescribing a
monic polynomial with $n$ distinct real roots at infinity and an
Eisenstein polynomial at an auxiliary rational prime. This gives an
irreducible polynomial with the required real and triadic completions.
\end{proof}
The product $\Delta_3(\mathbb K)$ is an upper bound for the density
of global three-cube sums. The global four-cube equality uses the
power-saving exceptional-set estimate of Theorem~\ref{nf3:four}.

\subsection{Affine identities and approximation at complex places}
\label{sec:complex-greedy}
There is an elementary construction for higher powers at totally
imaginary fields. It represents every sufficiently large target in
specific congruence classes, with control at each complex embedding.
For $k\ge2$ let
\[
      t_k=\left\lceil \frac{\log k}{\log(k/(k-1))}\right\rceil .
\]
For $0\ne c\in\mathcal O$ let $r_c(k)$ be the ordinary additive width
of the $k$th powers in their generated subgroup of $\mathcal O/c\mathcal O$.

\begin{proposition}\label{prop:complex-linear-identity}
Suppose $\mathbb K$ is totally imaginary and that
$\sum_{j=1}^l(a_jT+b_j)^k=cT+c_0$ in $\mathcal O[T]$ with $c\ne0$.
Then
\[
       G_{\mathbb K}(k)\le l+\max\{t_k,r_c(k)\}.
\]
More precisely, if $\nu-c_0$ modulo $c\mathcal O$ is a sum of at
most $t_k$ powers, then $l+t_k$ summands suffice for that target.
Every such affine identity has $l\ge k$. If $\mathbb K$ contains
a primitive $k^2$th root of unity, one can take $l=k$, $c=k^2$, and
$c_0=0$. In particular nine fourth powers suffice for all sufficiently
large targets divisible by $16$ in any field containing the sixteenth
roots of unity.
\end{proposition}
\begin{proof}
Put $m=\max\{t_k,r_c(k)\}$, or $m=t_k$ in the particular classes.
Balance $\nu$ by a $k$th power of a unit congruent to one modulo $c$. Such units have finite index in the full unit group, so Lemma~\ref{lem:hbal} applies to them and the prescribed residue class is preserved. If
$P=|N\nu|^{1/(k[\mathbb K:\mathbb Q])}$, then all its absolute values
are $O(P^k)$. The residue class of $c_0=\sum b_j^k$ is eligible.
Choose $m$ root classes modulo $c\mathcal O$ representing $\nu-c_0$;
append zero classes if necessary. For each of the first $m$ powers,
choose one complex $k$th root of each residual coordinate at the
archimedean places and round the resulting vector to the prescribed
coset of $c\mathcal O$. This coset has bounded covering radius.
If the residual coordinates have absolute value at most $R$, the new
ones are $O_{\mathbb K,c,k}(1+R^{1-1/k})$ by the binomial theorem.
After $m$ steps the residual equals $cy+c_0$ with $y\in\mathcal O$
and all its archimedean absolute values $O(P)$, since
$k(1-1/k)^m\le1$. The stated identity supplies the last $l$ roots.
All roots have absolute values $O(P)$ before undoing the balancing,
so their norms have the required bound.

To see the lower bound on $l$, group the terms with $a_j\ne0$ by
the distinct ratios $b_j/a_j=\xi_i$, and write $w_i$ for the sum
of $a_j^k$ in each group. The coefficients in degrees $k$ through
$2$ give $\sum_i w_i\xi_i^j=0$ for $0\le j\le k-2$. With at most
$k-1$ groups the Vandermonde matrix would force every $w_i=0$,
contrary to $c\ne0$. Finally take a primitive $k^2$th root $\theta$
and let $\zeta=\theta^k$. Summing powers of $\zeta$ gives
\[
     \sum_{j=0}^{k-1}\bigl(\theta^j(T+\zeta^j)\bigr)^k=k^2T.
\]
At $k=4$ one has $t_4=5$; the zero class permits all five initial
roots to be chosen in $16\mathcal O$.
\end{proof}

There is a smaller global count at prime powers when the required
roots of unity are present. For $k=p^a$ with $p\ge5$ and
$\mu_{k^2}\subset\mathbb K$ one has
\begin{equation}\label{eq:cyclo-imag-powers}
               G_{\mathbb K}(k)\le k+t_k .
\end{equation}
Indeed $r_c(k)\le H_{\mathbb K}(k)$ for every nonzero $c$:
represent a target locally at each prime dividing $c$, pad the
representations to a common length, and use the Chinese remainder
theorem on the roots. Theorem~\ref{thm:hteich} gives
$H_{\mathbb K}(k)\le k+1\le t_k$. Thus thirteen fifth powers suffice
over every field containing $\mu_{25}$, and twenty seventh powers over
every field containing $\mu_{49}$.

The same counts hold over totally imaginary fields containing
$k$th roots of all the binomial coefficients $\binom{k-1}{j}$,
$0\le j<k$, when $k=p^a$ and $p\ge5$. Indeed
\[
 \sum_{j=0}^{k-1}(-1)^{k-1-j}\binom{k-1}{j}(T+j)^k
                         = k!T+\frac{k!(k-1)}2 .
\]
This is the $(k-1)$st forward difference. As $k$ is odd, each signed
coefficient can be absorbed into an integral root and the identity
has exactly $k$ terms. Apply Proposition~\ref{prop:complex-linear-identity}
with $c=k!$.

\begin{proposition}\label{nf5:imaginary-sixteen}
For every totally imaginary number field $\mathbb K$ one has
$G_{\mathbb K}(5)\le16$.
\end{proposition}
\begin{proof}
The identity $0^2+5^2=3^2+4^2$ gives
\begin{align*}
 2880T={}&(T+5)^5+(T-5)^5+2T^5\\
         &-(T+3)^5-(T-3)^5-(T+4)^5-(T-4)^5 .
\end{align*}
Indeed the coefficients of $T^5$ and $T^3$ cancel, and the coefficient
of $T$ is $10(5^4-3^4-4^4)=2880$. We first choose eight fifth powers
representing the eligible target $\nu$ modulo $2880\mathcal O_{\mathbb K}$.
At a prime over five, $5\Vert2880$, and a sum of fifth powers modulo
$5\mathcal O_{\mathbb K}$ is a single fifth power by Frobenius.
At a prime of residue characteristic two or three let $H$ be the
nonzero fifth powers in its residue field and $E$ the subfield
generated by $H$. Put $C=H\cup\{0\}$. The reduction of $\nu$
lies in $E$. The additive stabilizer of $5C$ in $E$ is $H$-invariant
and hence is either zero or $E$. Kneser's theorem \cite{edkneser}
gives $|5C|\ge5|C|-4\ge|E|$ in the former case, since the
index of $H$ in the full multiplicative group is at most five.
Thus $5C=E$. A nonzero target has a representation with a unit;
for zero use $1,-1$. Hensel's lemma lifts this solution through
any power of the prime. Inserting zeros and applying the Chinese
remainder theorem gives the eight required congruence classes.

Balance the target by a fifth power of a unit and write $P$ for
the root scale, so that $|\sigma(\nu)|\ll_{\mathbb K}P^5$ at every
embedding. A lattice translate in a fixed congruence class contains
a point within bounded distance of any chosen vector of complex
fifth roots. If the current remainder has archimedean size $O(Y)$,
one such approximation gives a root of size $O(1+Y^{1/5})$
and a new remainder of size $O(1+Y^{4/5})$, by factorization
of the difference of fifth powers. After eight steps the remainder
is $2880T$, with
\[
  |\sigma(T)|\ll_{\mathbb K}1+P^{5(4/5)^8}\ll P .
\]
The identity supplies eight more fifth powers, since a negative
sign can be absorbed in an odd power. All sixteen roots have
absolute norm $O_{\mathbb K}(|\mathrm N\nu|^{1/5})$; undoing
the balancing unit preserves their norms. Zero terms can be removed.
\end{proof}

\subsection{Variation within a fixed totally real field}\label{nf3:progressions}
The covolume bound in Proposition~\ref{nf3:covol} also has a strictly
smaller upper bound for each fixed field. The congruence calculation
explains why a positive-density theorem for three cubes could not
provide the same lower proportion on every locally admissible
progression. Diaconu observed the corresponding obstruction for
bounded root regions over $\mathbb Q$ \cite{nf3:diaconu}.

\begin{proposition}\label{nf3:prog}
Let $\mathbb K$ be totally real. In any bounded Jordan-measurable
positive target region of positive volume under dilation, for every
$\varepsilon>0$ there is a residue class of a fixed ideal consisting
entirely of everywhere finitely locally representable targets on which
the upper relative proportion of sums of three totally positive
cubes is at most $\varepsilon$. There are also such residue classes
on which the limiting mean number of unordered positive three-cube
representations exceeds any prescribed constant. In particular the
unrestricted upper proportion is strictly less than
$\Gamma(4/3)^{3d}/(6D_{\mathbb K})$.
\end{proposition}
\begin{proof}
Write $\lambda=\Gamma(4/3)^{3d}/(6D_{\mathbb K})$.
For an integral ideal $\mathfrak m$ of norm $m$ and a class $r$, let
$w_{\mathfrak m}(r)$ count the triples of residue cubes with sum $r$,
counted with their ordered roots. Lattice counting as in
Proposition~\ref{nf3:covol} gives the mean number of unordered
representations in this target class as $\lambda w_{\mathfrak m}(r)/m^2$
in the limit. Indeed each of the $w_{\mathfrak m}(r)$ input classes has
relative lattice volume $m^{-3}$ and the target class has proportion $m^{-1}$.
The normalized multiplicities multiply by the Chinese remainder theorem.

Begin with $r=1\bmod18\mathcal O_{\mathbb K}$. This class has positive
normalized multiplicity and is locally representable everywhere.
At primes over $3$ the cubes of $1+3y$ fill $1+9\mathcal O_v$ by
Hensel's lemma, and above $2$ one lifts the residue cube $1$.
At every other finite place the additive span of the residue cubes is
the full residue field, since the span of the multiplicative cube subgroup is
a subfield, and an index-at-most-three subgroup has more than
$\sqrt q-1$ elements when the residue cardinality $q>4$.
The three-cube argument of Proposition~\ref{nf3:covol} then applies.

Take rational primes $p$ splitting completely in the normal closure of
$\mathbb K(\sqrt{-3})$. This set has positive Dirichlet density, and its
reciprocals have divergent sum. Choose one prime ideal $\mathfrak p$ above each such
$p>3$, with residue field $\mathbb F_p$. Let $g$ be a cubic Gauss sum,
and for $a\ne0$ write $S(a)=\sum_x e_p(ax^3)=
g\overline\chi(a)+\bar g\chi(a)$, with $\chi$ a cubic character.
Since $g\bar g=p$,
\[
                  S(a)^3=g^3+\bar g^3+3p S(a).
\]
Consequently for $r\ne0$ the normalized number of triples is
\[
 \frac{w_{\mathfrak p}(r)}{p^2}
  =1+\frac{3}{p}(\#\{x\in\mathbb F_p:x^3=r\}-1)
            -\frac{g^3+\bar g^3}{p^3}.
\]
The last term has magnitude
at most $2p^{-3/2}$. Imposing noncube residues at successively more
of these primes makes the product tend to zero, whereas imposing
nonzero cube residues makes it tend to infinity. A sufficiently long
finite product therefore gives either prescribed conclusion.
The represented fraction is at most the mean multiplicity. Finally
take one class with mean $\mu>1$. It has proportion $1/m$, and at most
one distinct target per target. Its excess multiplicity lowers the
global upper bound to $\lambda-(\mu-1)/m$.
\end{proof}

The excess-multiplicity assertion extends to every critical system of
$k$ summands of degree $k$.
\begin{corollary}\label{nf3:critical}
For every fixed totally real $\mathbb K$ and every $k\ge3$, the upper
proportion in Proposition~\ref{nf3:covol} is strictly smaller than
$\Gamma(1+1/k)^{kd}/(k!D_{\mathbb K}^{(k-1)/2})$. In locally admissible
progressions the limiting mean unordered multiplicity can be arbitrarily
large.
\end{corollary}
\begin{proof}
Select one prime ideal $\mathfrak p$ of residue norm $q$ with
residue characteristic prime to $2k$ and $-1$ a $k$th power. A rational
prime splitting completely in the normal closure of
$\mathbb K(\zeta_{2k})$ suffices. Thus the zero target modulo
$\mathfrak p$ has a smooth solution using two nonzero coordinates.
Let $c>0$ be the number of primitive residue vectors with zero power sum
divided by $q^{k-1}$. At modulus $\mathfrak p^{kh}$ the normalized
multiplicity of the zero class is exactly $1+hc$. In fact smooth
lifting gives the contribution $c$ from primitive vectors at every
positive level. Writing an imprimitive vector as $\varpi\mathbf y$
reduces the modulus exponent by $k$. The $q^{k(k-1)}$ excess choices
of the coordinates of $\mathbf y$ cancel the normalization change
because the number of variables equals the degree.

At all sufficiently large other finite places every residue target
has a smooth solution, since the orthogonality count differs from
$q^{k-1}$ by $O_k(q^{k/2})$ by the finite-field Weil bound, and only
the zero vector can be singular. At each of the finitely many remaining
places one can impose a fixed neighborhood of the target $1$, lifting
the unit root $1$ by Hensel's lemma. These fixed classes have positive
normalized multiplicity. The calculation of the mean in each
progression now proceeds as in Proposition~\ref{nf3:prog}, and the factor
$1+hc$ tends to infinity. Choosing any one progression of mean $>1$
gives a strict decrease of the global upper bound.
\end{proof}

\begin{theorem}\label{nf3:criticaldrought}
Fix a totally real field $\mathbb K$ and a target region as in
Proposition~\ref{nf3:prog}. If $\sqrt2\notin\mathbb K$, there are
everywhere locally admissible progressions on which the upper relative
proportion of sums of four totally positive fourth powers is arbitrarily
small. If $\sqrt2\in\mathbb K$, the limiting mean multiplicity for four
fourth powers has a positive lower bound, depending only on the field,
on every progression admitting a solution modulo its defining ideal.
The same lower bound, now depending also on $k$, holds for $k$ summands of
degree $k$ in every fixed totally real field for $k\ge5$.
\end{theorem}
\begin{proof}
Let $q$ be odd and write $S(a)=\sum_x\psi(ax^4)$ in the residue field
for a fixed nontrivial additive character $\psi$. If $q\equiv5\pmod8$,
choose a quartic character $\chi$, its Gauss sum $g$, and the quadratic
Gauss sum $g_2$. We have $\chi(-1)=-1$, so the inverse quartic character
has Gauss sum $-\bar g$, while $g_2^2=q$. For $a\ne0$, putting
$u=\bar\chi(a)$ gives
\[
                  S(a)=gu+g_2u^2-\bar g u^3 .
\]
On averaging the fourth power over $a\ne0$, only exponents divisible
by four survive. The result is
$g^4+\bar g^4-5q^2$, between $-7q^2$ and $-3q^2$. The zero-target
multiplicity divided by $q^3$ is therefore between
$1-7(q-1)/q^2$ and $1-3(q-1)/q^2$. At sufficiently large such primes
there is a smooth zero, so every target in the zero residue class lifts.

If $\sqrt2\notin\mathbb K$, total reality gives
$\mathbb K\cap\mathbb Q(\zeta_8)=\mathbb Q$, the maximal real subfield
of $\mathbb Q(\zeta_8)$ being $\mathbb Q(\sqrt2)$. Chebotarev over
$\mathbb K$ gives positive density of prime ideals of norm $5\bmod8$.
Imposing zero classes at successively more of these primes makes the
product of the normalized multiplicities tend to zero. At the finitely
many small places impose the fixed classes about $1$ used in
Corollary~\ref{nf3:critical}. All other large places admit smooth
solutions to every target. This proves the upper-proportion assertion.

Suppose instead that $\sqrt2\in\mathbb K$. The odd residue norms are
$1$ or $7\bmod8$, since their Frobenius in $\mathbb Q(\zeta_8)$ fixes
$\sqrt2$. When $q\equiv1\bmod8$ the sums $S(a)$ are real, since $-1$
is a fourth power. When $q\equiv7\bmod8$ the fourth powers are the
squares, and $S(a)^4=q^2$ for $a\ne0$ by the quadratic Gauss identity.
In both cases the zero-target residue count is at least $q^3$.
Away from the zero target the normalized count is $1+O(q^{-3/2})$ by expanding $S(a)^4$ into multiplicative character terms with coefficients
$O(q^2)$, and applying the Gauss bound on each nontrivial additive twist.
Thus, apart from finitely many places, the primitive solutions alone
have normalized multiplicity $\ge1-O(q^{-3/2})$ at all conductors
and for all targets, by smooth lifting. For $k\ge5$, the ordinary
Weil bound gives normalized count $1+O_k(q^{1-k/2})$ uniformly also
at zero. Subtracting the zero vector gives the same summable lower
error for the primitive count.

It remains to control the finitely many other places at arbitrary
conductor depth $j$. Choose a solution to the given congruence and
lift its coordinates to the local ring. The exact power sum of the
lift lies in the given class. If this vector is primitive, varying
the other $k-1$ coordinates modulo $\varpi^{2v(k)+1}$ and applying
strong Hensel to its unit coordinate shows that the normalized
multiplicity is $\gg_{\mathbb K,k}1$, uniformly in $j$. More precisely,
one can fix the other coordinates modulo this power and vary them freely
at every subsequent level. If all lift coordinates are divisible by
$\varpi^t$, with $t$ their minimal valuation and $kt<j$, division by
$\varpi^t$ gives a primitive vector at depth $j-kt$. Its contribution
has exactly the same normalization, by the critical scaling count in
Corollary~\ref{nf3:critical}. For $kt\ge j$, including the zero lift,
the target class is zero and the roots divisible by
$\varpi^{\lceil j/k\rceil}$ by themselves contribute normalized
multiplicity $q^{j-k\lceil j/k\rceil}\ge q^{1-k}$.
This gives a uniform positive factor at every remaining place.
Multiplication over the prime ideals and the lattice first moment now
prove both lower bounds for the mean.
\end{proof}

\subsection{Primitive triples with many companions}
The congruence variation has a consequence for the nontrivial
points on the Fermat cubic fourfold over a field of any signature.
Here a squarefree root means that its principal ideal is squarefree.
Let $B=B_1\times B_2\times B_3$ be a product of fixed open boxes
in the archimedean space, bounded away from zero in every
coordinate and positive at all real places. The boxes may be made
arbitrarily small. To remove the pairing planes, take at one
archimedean place three small neighborhoods of distinct positive
real numbers, in increasing order; at a complex place these are
disks. Put $X=P^d$.
\begin{proposition}\label{nf3:manyprimitive}
For every fixed $j\ge2$, at least $c_{B,j}X^3$ triples in $PB$
with pairwise coprime squarefree roots belong to fibers containing
$j$ distinct triples with the same properties. Consequently there
are $\gg_{B,j}X^3$ distinct projective points represented by
$PB\times PB$ on the six-variable equal-cube hypersurface off
its pairing planes. The number of distinct targets with these
$j$ representations is $\gg_{B,j,\varepsilon}X^{11/4-\varepsilon}$
unconditionally. All $3j$ roots can be taken distinct and
$X^\eta$-smooth for any prescribed $\eta>0$.
A six-cube mean of order $O(X^3)$ would give
$\gg_{B,j} X^3$ such targets.
\end{proposition}
\begin{proof}
Choose rational primes splitting completely in the normal closure
of $\mathbb K(\sqrt{-3})$ and take one prime ideal over each.
Impose $f(\mathbf x)=1$ at finitely many of these primes.
The cubic Gauss calculation of Proposition~\ref{nf3:prog} works
in every signature and gives the normalized multiplicity
$1+6/p+O(p^{-3/2})$. Pairwise primality and squarefreeness leave
this estimate unchanged: two zero coordinates modulo $p$ account
for $O(1)$ solutions, and in each of the $O(p)$ solutions with a
specified zero coordinate only $1/p$ of its lifts vanish to
second order. At each other prime ideal of norm $q$, the natural
density of permitted root triples is
$(1-1/q)^3(1+3/q)$. This product converges to a positive constant.
To pass from finite congruences to all primes, note that the
coordinate ranges are bounded homothetic boxes of volume $O(X)$.
At $q>L$, either two coordinates belong to the prime ideal or one
belongs to its square. In each case lattice packing bounds the
count by $O_B(X^3/q^2)$ in the range where a nonzero coordinate
is possible; use $q\ll_B X$ or $q^2\ll_B X$ respectively.
The tail over $q>L$ is $o_{L\to\infty}(X^3)$.
For the smooth-root restriction, choose a fixed positive exponent
no larger than both the prescribed $\eta$ and the moment cutoff in
Theorem~\ref{nf3:smoothsix}. Its fixed-class construction gives a lower bound by a fixed positive
multiple of each lattice density. The multiple is independent of
the fixed moduli. One may therefore first choose the progression,
then the finite squarefreeness and coprimality sieve, and finally
the box scale.

If the imposed modulus has norm $m$, there are only
$O_B(X^3/m)+o_m(X^3)$ possible targets in the progression.
The mean root count relative to this bound can be made larger
than $2j$, since the selected split primes have divergent reciprocal
sum by Chebotarev. Fix the modulus. Discarding every fiber
of size less than $j$ leaves $\gg_{B,j}X^3$ permitted triples.
Pair each with a different one in the same fiber. The separated
archimedean neighborhoods exclude the nonidentity root-of-unity
pairings, while equality of the two triples has been removed.
Two integral representatives of one of these projective points
differ by a unit, as the coordinates of either triple generate
the unit ideal. Their archimedean sizes in $PB$ constrain the unit
to a fixed compact set, which contains only $O_{\mathbb K,B}(1)$
units. This proves the projective count.
Finally Cauchy's inequality applied to the surviving triples and
the smooth sixth-moment bound $O_{\varepsilon}(X^{13/4+\varepsilon})$
of Theorem~\ref{nf3:smoothsix}
gives the asserted unconditional target count. If two different
triples share a root, the remaining four coordinates give a
nontrivial point on the cubic surface $x_1^3+x_2^3=y_1^3+y_2^3$.
Their coordinate ideal is the unit ideal by pairwise coprimality
inside each triple. The separated neighborhoods keep these points
off all the lines, and their relative projective heights are $O(X)$.
This surface contains three coplanar rational lines (the three
ordinary pairings in signed coordinates). The number-field bound
of Browning and Swarbrick Jones \cite[Theorem 1.2]{nf3:bsj}
therefore leaves $O_{\mathbb K,\varepsilon}(X^{4/3+\varepsilon})$
quadruples, with the same bounded unit multiplicity as above.
There are $O(X)$ choices of the shared root. Deleting all targets
with any such coincidence costs $O(X^{7/3+\varepsilon})$ targets
and, by Cauchy and the smooth sixth moment,
$O(X^{67/24+\varepsilon})=o(X^3)$ surviving root triples. The six
coordinates of the counted projective points can thus also be chosen
distinct. A mean of order
$O(X^3)$ gives the final assertion.
\end{proof}

\begin{corollary}\label{nf3:growing-fibres}
In Proposition~\ref{nf3:manyprimitive} the $\gg X^{11/4-\epsilon}$
distinct targets can each be required to have at least
$\exp(c\sqrt{\log X}/\log\log X)$ representations with all the
stated restrictions, where $c>0$ may depend on $\epsilon$,
the field, the boxes and the desired smoothness.
\end{corollary}
\begin{proof}
The plane Fermat cubic $C:x^3+y^3+z^3=0$ has complex multiplication
by $\mathbb Q(\sqrt{-3})$, via the automorphism of order three
fixing a point. CM equidistribution
\cite[Corollary~3.8]{nf3:fitesutherland} over
$\mathbb K(\sqrt{-3})$ gives $\gg y/\log y$ prime ideals in
$[y/2,y]$ at which the Frobenius trace satisfies $a\le-\sqrt q$,
with $q$ the norm. We can discard primes of residue degree
greater than one over $\mathbb Q$, and project to $\mathbb K$
with multiplicity at most two. The residue fields and their
traces are preserved. The affine zero fiber has exactly
$1+(q-1)\#C(\mathbb F_q)=q^2-(q-1)a$ vectors, of which at
most $O(q)$ have a zero coordinate, by intersection with the
three projective coordinate lines. The normalized number with
all three roots units is consequently $\ge1+1/(2\sqrt q)$.

Put $y=c_0\log X$ and form the product ideal $\mathfrak m$ of the
selected primes. Its norm $m$ is $\le\exp(O_{\mathbb K}(y))$.
Prescribe target zero modulo $\mathfrak m$ and require all roots
units there. The product $G$ of the local gains is
$\gg\exp(c_1\sqrt y/\log y)$. We first have
$\gg X^3G/m$ smooth root triples with good local conditions at the
selected primes and all primes of norm at most a sufficiently
large fixed $L$. In fact the cofactor boxes in the proof of
Theorem~\ref{nf3:smoothsix} have all sides
$\gg X^{\kappa/d}$ for some fixed $\kappa>0$. A residue class
modulo an ideal of norm $\le X^{\kappa/2}$ has its volume
lower bound uniformly there: every ideal is a balanced scalar
multiple of one of finitely many ideal class representatives,
so the corresponding scaled fundamental cells have diameter
$O_{\mathbb K}(X^{\kappa/(2d)})$. The chosen prime marks for the
smooth construction are all coprime to $\mathfrak m$, since their
norms exceed a fixed power of $X$. Choose $c_0$ sufficiently small.

At an unselected prime of norm $L<q\le X^{1/8}$ the
unrestricted upper lattice count modulo $\mathfrak m\mathfrak p^2$
and the Chinese remainder theorem bound the mass excluded by
squarefreeness or pairwise coprimality by
$O((X^3G/m)q^{-2})$, since $m q^2=o(X)$ uniformly.
Larger primes cost $O_{\epsilon'}(X^{3-1/8+\epsilon'})$
by the unrestricted lattice count. First take $L$ large, then
$c_0$ small; this preserves $\gg X^3G/m$ permissible triples.
There are $O(X^3/m)$ possible targets. Remove fibers of size
below a sufficiently small fixed multiple of $G$, and apply
Theorem~\ref{nf3:smoothsix} to the rest. Cauchy's inequality gives
$\gg X^{11/4-\epsilon/2}/m^2$ targets, each with the asserted
multiplicity. Shrinking $c_0$ again absorbs $m^2$. The
$O(X^{7/3+\epsilon'})$ targets with any shared root can be removed.

For comparison, every primitive three-cube fiber at any modulus
of norm $m$ has normalized size
$\le\exp(O_{\mathbb K}(\sqrt{\log(3m)}/\log\log(9m)))$.
At primes away from $3$, primitive solutions lift smoothly. The
zero fiber has normalized size $1+O(q^{-1/2})$ by the curve
bound, and the nonzero fibers have size $1+O(q^{-1})$ by the
smooth cubic surface count minus the curve at infinity. At
places over $3$, a unit variable has derivative of valuation
$v(3)$, giving a bounded primitive ratio independently of depth.
The prime ideal theorem gives the assertion on splitting
$\sum_{\mathfrak p\mid\mathfrak m}q^{-1/2}$ at $\log(3m)$.
\end{proof}

If pairwise coprimality inside a triple is relaxed, the same
number of targets admits at least $\exp(c\log X/\log\log X)$
representations with all roots distinct, squarefree and smooth.
Use every prime ideal of norm $q\in[y/2,y]$, and prescribe
target zero modulo $\mathfrak p^3$. The branch with all coordinates
units gives normalized mass $1+O(q^{-1/2})$ by smooth lifting from
the plane curve. The disjoint branch in which all three coordinates
have valuation exactly one gives $(1-1/q)^3$, as every such triple
has cube sum divisible by $\mathfrak p^3$. Each prime thus gives
gain at least $3/2$ for large $y$. The prime ideal theorem gives
$\gg y/\log y$ marks and the product modulus still has norm
$\exp(O_{\mathbb K}(y))$. The preceding sieve works unchanged,
and pairwise coprimality can still be imposed away from the
selected primes. To discard shared roots, let $\mathfrak c$ be the
coordinate ideal of a nontrivial equal sum of two cubes in the boxes.
Its relative projective height is $O(X/\mathrm N\mathfrak c)$.
For each fixed coordinate ideal, proportional representatives
differ by a unit, with bounded multiplicity by the archimedean
boxes. Summing the surface bound
$O_{\epsilon'}((X/\mathrm N\mathfrak c)^{4/3+\epsilon'})$
over all possible ideals retains $O_{\epsilon'}(X^{4/3+\epsilon'})$.
This proves the assertion about distinct roots.

At an arbitrary ideal modulus the unrestricted congruence gain has
upper bound $\exp(O_{\mathbb K}(\log(3m)/\log\log(9m)))$,
with $m$ its norm. Indeed the nonprimitive branch at depth
$j\ge3$ contributes the normalized density at depth $j-3$ for the
scaled target, provided it is divisible by $\mathfrak p^3$.
Summing the bounded primitive contributions gives an ideal-divisor
factor $\prod_{\mathfrak p^j\parallel\mathfrak m}(j+1)$ times the
previous primitive bound, up to constant factors at primes over $3$.
The usual maximal-order bound for the divisor function applies,
for example by bounding the ideal-divisor factor by $\tau(m)^d$.

\section{The three-cube density problem over a number field}\label{sec:field-cube-outlook}

Fix a target region $\Omega$ as above and choose positive-volume
root boxes $B_1,B_2,B_3$ in the archimedean space, positive at every
real place, such that $x_1^3+x_2^3+x_3^3\in\Omega$ whenever
$x_i\in B_i$. Such boxes exist about the coordinates of a triple
mapping to the interior of $\Omega$. For $B=B_1\times B_2\times B_3$,
let $E_B(A)$ count pairs of ordered triples in
$\prod_i(\mathcal O_{\mathbb K}\cap A B_i)$ with equal sums of
cubes. There are $\asymp_B A^{3d}$ triples and their targets lie
in $A^3\Omega$. Consequently a bound
$E_B(A)\ll_{\mathbb K,B}A^{3d}$ on a single such product yields
$R_\Omega(X)\gg_{\mathbb K,\Omega} X^d$ for $X=A^3$.
For some fixed $\theta>0$ the
$(1+\theta)$-moment of the representation multiplicities may replace the second
moment, with the same right side. Indeed, H\"older's inequality gives
$(\sum r)^{1+\theta}\le (\#\{r>0\})^\theta\sum r^{1+\theta}$.
An $A^\varepsilon$ loss in a moment gives an exponent loss in the conclusion.

In the rational case Wang's theorem \cite[Theorem 1.7]{c3:wang} supplies
the second-moment bound under automorphy and GRH for a geometric family,
a two-point ratios conjecture and a square-divisor sieve for the dual form.
A number-field analogue requires control of the $L$-functions of the
smooth hyperplane sections of the Fermat fourfold as the section varies
over $\mathbb K$, together with corresponding bad-prime and archimedean
estimates. A precise first target is the restricted second moment above:
one may place the three roots in disjoint small boxes near three distinct
archimedean points. Fixing the same order of the disjoint boxes on
both sides eliminates nonidentity permutations, and their directions
may be chosen to exclude pairings involving nontrivial cube roots of
unity. The identical triples contribute on the scale $A^{3d}$.
A function-field version for $\mathbb F_q[t]$ in characteristic greater
than three is known conditionally on ratios estimates
\cite{nf3:bgw}. Over a number field, the choice of positive boxes preserves
the size of the first moment while leaving the required estimate of
distinct collisions open.

A portion of the square-divisor estimate remains elementary over a
number field. For the primitive dual form $D=F^\vee$ and fixed $\mathbb K$,
the number of $\mathbf c\in\mathcal O_{\mathbb K}^6$ of embedding height
at most $Z$ with $\mathfrak p^2\mid D(\mathbf c)$ for some
$L\le\mathrm N\mathfrak p<2L$ satisfies
\[
 \ll_{\mathbb K} Z^{6d}/L+L Z^{5d}\qquad (1\le L\le Z^d).
\]
Indeed, exclude the finitely many primes of bad reduction, absorbing them
by a fixed constant. At a prime of norm $q$ a residue class modulo an
ideal $\mathfrak a$ contains $O_{\mathbb K}(1+Z^d/\mathrm N\mathfrak a)$
elements of height at most $Z$. The norm lower bound for a nonzero
difference in $\mathfrak a$ gives separation at distance
$(\mathrm N\mathfrak a)^{1/d}$ in the maximum norm, so this follows by
packing. The singular locus of the reduced dual has codimension at least
two and $O(q^4)$ vectors modulo $\mathfrak p$. Since $q\ll Z^d$ its
contribution is $O_{\mathbb K}(Z^{6d}/q^2)$. At a smooth zero, fix every
coordinate except one for which the partial derivative is a unit. There
are $O(1)$ simple roots in that coordinate modulo $\mathfrak p$, each
lifting to one class modulo $\mathfrak p^2$. This contributes
$O_{\mathbb K}(Z^{5d}(1+Z^d/q^2))$. Summing proves the estimate.
Writing $M=Z^d$, it is $O_{\mathbb K}(M^6(L^{-1}+M^{-\eta}))$
for $1\le L\le M^{1-\eta}$. Thus for fixed $\xi,\eta>0$ the saving
is a power of $M$ when $M^\xi\le L\le M^{1-\eta}$; for bounded $L$
only a saving expressed in $L$ is asserted.

For totally real fields one can make a more precise prediction from the
six-variable cubic. Write $D=D_{\mathbb K}$,
$\lambda_{\mathbb K}=\Gamma(4/3)^{3d}/(6D)$, and denote the lower and
upper proportions of targets represented by three positive cubes in
$X(0,1)^d$ by $\underline\delta_{\mathbb K}$ and $\overline\delta_{\mathbb K}$.
The height limit in these definitions is taken with the field fixed.

\begin{proposition}\label{nf3:conditional}
Let $\mathbb K$ be totally real. Assume the weighted Hooley--Manin asymptotic over $\mathbb K$ for the
six-variable form $F(\mathbf x,\mathbf y)=\sum_{j=1}^3(x_j^3-y_j^3)$.
Precisely, for every smooth compactly supported weight assume that the
weighted lattice count at root scale $P$, after subtracting the weighted
count on the union of the affine cones over the $\mathbb K$-rational
planes of the projective fourfold, equals
$P^{3d}\mathfrak S_{\mathbb K}\sigma_\infty+o_{\mathbb K,w}(P^{3d})$.
Here the factors are the ordinary primitive-character singular series and
the real singular integral in the covolume-one lattice normalization.
Then
\[
 \lambda_{\mathbb K}-\tfrac12\mathfrak S_{\mathbb K}\lambda_{\mathbb K}^2
 \le\underline\delta_{\mathbb K}\le\overline\delta_{\mathbb K}
 \le\lambda_{\mathbb K},\qquad \mathfrak S_{\mathbb K}\ll_d1 .
\]
Thus both proportions equal $\lambda_{\mathbb K}+O_d(D^{-2})$ as the
discriminant grows through totally real fields of a fixed degree
satisfying the hypothesis.
\end{proposition}

\begin{proof}
Order the three roots at one real embedding, allowing every ordering
at every other place. The complex planes of the Fermat cubic fourfold
come from pairing the six coordinates with cube roots of unity after
the sign change on $\mathbf y$ \cite[Sect.~2.5]{nf3:dio}. Total reality,
positivity and the chosen ordering therefore leave only the identical
pair in this chamber. Write $f(\mathbf x)=\sum_j x_j^3$.
With ordinary Euclidean measure on the embeddings, the pushforward of
the chamber with $f\in(0,1)^d$ has constant density
$\Gamma(4/3)^{3d}/6$ there. With the measure assigning the integer
lattice covolume one on each variable, and the same measure on the
target, the pushforward density is $\lambda_{\mathbb K}$. Indeed, the
three input factors each contribute $D^{-1/2}$ and conversion of the
target measure contributes $D^{1/2}$.

Take a smooth cutoff $0\le w\le1$ compactly supported in this open
chamber with $f\in(0,1)^d$. For the distinct triples above a target
write $w_i$ for their weights. The first two Bonferroni terms give
$\sum_i w_i-\sum_{i<j}w_iw_j\le\mathbf1_{\{r>0\}}$, by taking
independent events of probabilities $w_i$. Subtracting the identical pairs from the
ordered-pair count removes the entire plane contribution. The real
singular integral for a product weight is the squared $L^2$ norm of
the pushforward density, with the target lattice measure. First apply the
assumed asymptotic at fixed $w$, then increase the cutoffs to the chamber.
By monotone convergence this integral divided by the target
volume tends to $\lambda_{\mathbb K}^2$. The first moment divided by the
same volume tends to $\lambda_{\mathbb K}$. This gives the lower bound, while the upper bound follows from Proposition~\ref{nf3:covol}.

For a prime of norm $q$ away from $3$, let $M$ count the zeros of $F$
in the six-dimensional residue space. The cubic Gauss sum bound
$2\sqrt q$ and orthogonality give $M=q^5+O(q^3)$. All nonzero solutions
are smooth. Rescaling the vectors divisible by a uniformizer therefore
gives the local density $(M-1)/(q^5(1-q^{-3}))=1+O(q^{-2})$, with an
absolute constant. At a prime above $3$ put $e=v(3)$. On a ball of radius
$q^{-e-1}$ about a unit, the cube map pushes normalized Haar probability
to uniform probability on a coset modulo the $(2e+1)$-st power of the
prime ideal, by Taylor expansion and Hensel's lemma. Partition primitive
vectors according to a unit coordinate. Their density on $F=0$ is at
most $6q^{2e+1}$, and the same rescaling controls the imprimitive
vectors by the factor $(1-q^{-3})^{-1}$. Since $e\le d$, $q\le3^d$ at
these places, and each rational prime has at most $d$ places above it,
the product of all local bounds is $O_d(1)$.
\end{proof}

\subsection{The dual sieve in a number field}\label{nf3:sfstrength}
For $F=\sum_{i=1}^6h_i^3$ let $D=F^\vee$. For a fixed number field
$\mathbb K$ of degree $d$, let $B_{\mathbb K}(Z,L)$ count vectors of six
algebraic integers of embedding height at most $Z$ for which
$\mathfrak p^2\mid D(\mathbf c)$ for some prime ideal of norm in
$[L,2L)$. Write $M=Z^d$; exponents on $M$ below refer to the
one-coordinate volume in this fixed degree.

\begin{proposition}\label{nf3:sfprop}
For $0<\eta<1/3$ there is a $\delta(\eta)>0$ such that
\[
 B_{\mathbb K}(Z,L)\ll_{\mathbb K,\eta} M^6 L^{-\delta(\eta)}
       \qquad(1\le L\le M^{4/3-\eta}).
\]
\end{proposition}
\begin{proof}
For $L\le M^{9/10}$ use the simple-root slicing estimate above.
For all larger $L$ the residues where $D=0$ and some coordinate or all
first derivatives vanish cost $O_{\mathbb K}(M^6/L+M^5\log(2M))$ in
the union. For $L\le M$, count the $O(q^4)$ residue vectors at a prime
of norm $q\asymp L$, then use the lattice packing bound in each
coordinate. For $L\ge M$, choose an integral unimodular change of variables
to make $D$ have constant nonzero leading coefficient in the last variable.
Its resultant in this variable with $(\prod c_i)\partial_1D$ is nonzero.
The exceptional base vectors annulling the resultant number
$O_{\mathbb K}(M^4)$ by polynomial root counting on five sets of size
$O(M)$, and admit $O(M)$ last coordinates. For every other base vector
the nonzero resultant has norm $M^{O(1)}$. For each of its $O(\log(2M))$
prime ideal factors of norm at least $M$ there are only
$O_{\mathbb K}(1)$ possible last coordinates.

For the remaining residues use Poisson summation on the Minkowski lattice
with a product Gaussian majorant. At a prime $\mathfrak p$ of norm $q$ away
from the different and finitely many primes of bad reduction, write $T$
for the (unnormalized) additive Fourier transform modulo $\mathfrak p^2$
of the good residues satisfying $D=0\pmod{\mathfrak p^2}$.
Fix a primitive character $\psi$ of this ring. With frequency $\mathbf h$
in the residue ring to its sixth power, if
$\mathbf h\not\equiv0\pmod{\mathfrak p}$ one has
\[
 T(\mathbf h)\le q^4
       \sum_{a\bmod\mathfrak p^2}^{*}\psi(aF(\mathbf h))
        +q^5\mathbf1_{\mathbf h\bmod\mathfrak p\in E},
 \qquad \#E(\mathcal O_{\mathbb K}/\mathfrak p)\ll q^4.
\]
Indeed, summation over lifts kills every good residue except those normal
to $\mathbf h$. The contact is unique by the smooth Gauss correspondence
\cite[Proposition 9.11, proof]{c3:wang}. For $F(\mathbf h)=0$ modulo
$\mathfrak p$, lift the point to $F(\tilde{\mathbf h})=0$ modulo
$\mathfrak p^2$. Taylor expansion gives
$\mathbf h\cdot\nabla F(\tilde{\mathbf h})=F(\mathbf h)$ modulo
$\mathfrak p^2$, while $D(\nabla F(\tilde{\mathbf h}))=0$ in this ring.
The scalar line of good contacts therefore contributes
$q^5(q\mathbf1_{\mathfrak p^2\mid F(\mathbf h)}-1)$.
Exceptional contacts occur in the codimension-two set where $F=0$ and
some coordinate or some nonempty proper cubic subsum is zero, as follows
by differentiating the signed factors of $D$. Omitting these contacts can
cost at most $q^5$ per such frequency. If $\mathbf h\equiv0\pmod{\mathfrak p}$
but not modulo $\mathfrak p^2$, summing by scalar lines gives
$|T|\ll q^9$, since the irreducible dual has $O(q^3)$ projective points
in each hyperplane and $O(q^4)$ in total. At zero $|T|\ll q^{10}$.

To track the dual lattice select a fixed integral ideal $\mathfrak b$
from a finite list of class representatives so that
$(\pi)=\mathfrak b\mathfrak p^2\mathfrak D$ is principal, and choose $\pi$
balanced, of height $\asymp_{\mathbb K}q^{2/d}$. Exclude primes dividing
these representatives. The trace-dual lattice of $\mathfrak p^2$ equals
$\mathfrak b/\pi$. Fix a nonzero $u\in\mathfrak b$ for each representative
and omit primes dividing it. Then $z\mapsto e(uz/\pi)$ induces a
primitive character modulo $\mathfrak p^2$. A dual vector
$\mathbf h/\pi$ with $h_i\in\mathfrak b$ has finite frequency
$u^{-1}\mathbf h\bmod\mathfrak p^2$. Its cubic phase is thus
$F(\mathbf h)$ with the scalar unit $u^{-3}$ absorbed into $a$.
In Poisson summation the prefactor is $O_{\mathbb K}(M^6/q^{12})$ and the
dual weight is a product of six nonnegative Gaussians on the scale
$Y=q^{2/d}/Z$. Their shapes range over a compact family.
Write $H=Y^d=q^2/M$. The exceptional frequencies, those in
$\mathfrak p$ and those in $\mathfrak p^2$ contribute in total
\[
 \ll_{\mathbb K}(M^6/q^{12})
           \{q^9(1+H/q)^6+q^{10}\}\ll M^6/q^2+q^3
\]
by lattice packing and Gaussian decay. The Ramanujan term at frequencies
divisible by $\mathfrak p$ is nonnegative and can be retained.

Let $U(a)$ be the resulting one-variable Gaussian cubic sum over
$h\in\mathfrak b$ on scale $Y$. Orthogonality and the at most three unit
cube roots modulo $\mathfrak p^2$ yield
$q^{-2}\sum_a|U(a)|^2\ll_{\mathbb K,\varepsilon}M^\varepsilon H$.
Indeed each residue modulo $\mathfrak p^2$ has Gaussian mass
$O_{\mathbb K}(1+H/q^2)$, and each unit cube has at most three
unit preimages by simple lifting. The nonunits cost
$O_{\mathbb K}((1+H/q)^2)\ll_{\mathbb K} H$, since $1\le H\le q^2$.

For primitive $a$ we have the exact-denominator cubic Weyl estimate
\[
 |U(a)|^4\ll_{\mathbb K,\varepsilon} M^\varepsilon(H^3+Hq^2).
\]
Here are details needed for nonprincipal denominators. Truncate the
Gaussian at embedding height $Y\log(2M)$ and use partial summation in
an integral basis to reduce to boxes and absorb all logarithms in
$M^\varepsilon$. Detecting the fixed ideal $\mathfrak b$ by its additive
characters introduces only linear twists. With phase $\alpha h^3$ of
exact conductor $\mathfrak p^2$, twice differencing gives a factor $H$
times a sum over two shifts $h_1,h_2$ of the absolute value of the linear
sum of frequency $6\alpha h_1h_2$ in an intersected box. Zero shifts
contribute $O(M^\varepsilon H^3)$ in total. For nonzero product
$t=6h_1h_2$, the fixed-field ideal divisor bound and Dirichlet's unit
theorem give $O_{\mathbb K,\varepsilon}(M^\varepsilon)$ pairs of bounded
shifts. The elements $t$ have height $O(Y^2\log^2(2M))$, so each residue
modulo $\mathfrak p^2$ occurs
$O_{\mathbb K,\varepsilon}(M^\varepsilon(1+H^2/q^2))$ times by packing.
For a basis $\omega_1,\ldots,\omega_d$, the linear sums are bounded by
the product of $\min(Y\log(2M),\|\operatorname{Tr}(\alpha t\omega_j)\|^{-1})$
over $j$, up to fixed constants. As $t$ varies modulo $\mathfrak p^2$,
the trace coordinates modulo one run through the grid of all
$q^2$ characters. Its distinct points have separation
$\gg_{\mathbb K}q^{-2/d}$ in the torus maximum norm. Indeed, their
difference represents a nonzero element of
$\mathfrak p^{-2}\mathfrak D^{-1}$ after subtracting any translate from
$\mathfrak D^{-1}$, and its absolute norm is $\gg_{\mathbb K}q^{-2}$.
Trace coordinates and embedding coordinates are related by a fixed
invertible real linear change of variables. Each cube of side $1/Y$
therefore contains $O_{\mathbb K}(q^2/H)$ grid points, since $H<q^2$.
Summing the displayed products over these cubes gives
$O_{\mathbb K,\varepsilon}(M^\varepsilon q^2)$ by the harmonic sum in
each coordinate. Combining the estimates proves the Weyl bound.

The Ramanujan contribution is at most $(M^6/q^8)\sum_a^* |U(a)|^6$.
The fourth-power supremum and second-moment estimate give
\[
 \ll_{\mathbb K,\varepsilon}M^\varepsilon
     \frac{M^6}{q^6}(H^4+q^2H^2)
   =M^\varepsilon(M^2q^2+M^4).
\]
Summing over the primes in the dyadic range gives
$O_{\mathbb K,\varepsilon}(M^\varepsilon(M^2L^3+M^4L))$, which saves a
power of $M$ up to $L=M^{4/3-\eta}$. The other contributions save a
power in the same range.
\end{proof}

The same analysis admits a version with an isolated sixth moment. With
$U_{\mathfrak p}(a)$ denoting the Gaussian sum above, for
$M^{9/10}\le L\le M^{3/2}$ one has
\[
 B_{\mathbb K}(Z,L)
 \ll_{\mathbb K,\varepsilon}M^\varepsilon
       (M^6/L+M^5+L^{7/2}+ML^3)
       +\frac{O_{\mathbb K}(M^6)}{L^8}
         \sum_{L\le\mathrm N\mathfrak p<2L}\sum_a^*|U_{\mathfrak p}(a)|^6 .
\]
Indeed one may count every torus contact branch over
$\mathcal O_{\mathbb K}/\mathfrak p^2$ by writing
$c_i=b x_i^2$, where $x_6=1$, every parameter is a unit, and
$F(\mathbf x)=0$. Each good residue solution to the dual equation lifts
uniquely to its branch by Hensel's lemma for the square roots of
$c_i/c_6$. The Fourier transform for a primitive frequency is exactly
$q^4\mathbf1_{\mathfrak p\nmid\prod h_i}
\sum_a^*\psi(aF(\mathbf h))$. To verify this, vary $b$ and the tangent
directions of $F$ at level $\mathfrak p$. The equations are
$\sum h_i x_i^2=0$ and $h_i=t x_i$ for $i\le5$, hence $h_6=t$.
Thus all frequency coordinates are units. If $F(\mathbf h)=0$ modulo
$\mathfrak p$, this contact is unique, and Euler's identity on any lift
reduces the scalar sum to the Ramanujan sum as before. Removing
nonunits from the one-variable Gaussian costs $O((1+H/q)^6q^2)$ in the
sixth-moment upper bound, which is absorbed in the stated error:
the unit sum $U^\times$ satisfies
$|U^\times|^6\le 2^5(|U|^6+|U-U^\times|^6)$, and the latter
difference is $O_{\mathbb K}(1+H/q)$.

More explicitly, let $r(t)$ be the total product weight of ordered
unit triples whose cube sum is $t\bmod\mathfrak p^2$, with the
one-variable Gaussian used in $U^\times$. Orthogonality gives
\[
 \sum_{a\bmod\mathfrak p^2}^{*}|U^\times(a)|^6
 =q^2\sum_{t\bmod\mathfrak p^2}r(t)^2
  -q\sum_{u\bmod\mathfrak p}
       \left(\sum_{t\equiv u\bmod\mathfrak p}r(t)\right)^2 .
\]
The second term is the energy of all additive characters factoring
through the residue field, including the trivial character. In the
displayed bound for $B_{\mathbb K}$ the primitive sum has coefficient
of order $M^6/L^8$; the argument uses an upper bound for that sum.

For an imprimitive, nonzero frequency write locally $\mathbf h=\varpi
\mathbf k$. Its transform is $q^5(qN_{\mathbf k}-N)$, where $N$ counts
the torus points of the projective Fermat cubic fourfold and $N_{\mathbf k}$
also imposes the quadric $\sum k_i x_i^2=0$. If no nonempty subsum
of the $k_i^3$ vanishes, the gradients are independent, since dependence
would make the nonzero $x_i$ proportional to their corresponding $k_i$.
This smooth complete intersection of dimension three is geometrically
connected and irreducible. Removing the coordinate walls costs $O(q^2)$,
and the Lang--Weil estimates give $qN_{\mathbf k}-N=O(q^{7/2})$.
For every nonzero $\mathbf k$ the dimension bound gives $O(q^4)$ instead.
Choose a global balanced $\varpi$ generating $\mathfrak c\mathfrak p$
with $\mathfrak c$ from a fixed finite list coprime to $\mathfrak p$.
Then the $k_i$ belong to fixed fractional lattices and have volume
scale $H/q=q/M$. In a box truncated at a logarithmic factor the scalars
reduce injectively modulo $\mathfrak p$, since distinct elements of
$\mathfrak p\mathfrak c^{-1}$ have norms $\gg q$ and the box has volume
$o(q)$. The degree bound on the nonzero subsum polynomial on this
grid gives $O(M^\varepsilon(1+H/q)^5)$ exceptional vectors. The factors
from $\mathbf h/\pi$ and the fixed character normalization are scalar
units and leave the zero set of the homogeneous subsum polynomial
unchanged. The zero mode is $O(q^{10})$. With the Poisson prefactor
the total error at one prime is therefore
$\ll M^\varepsilon(M^6/q^2+q^{5/2}+Mq^2)$, proving the asserted bound.

In particular suppose that, for fixed $\mathbb K$ and every
$\varepsilon>0$, the full signed count satisfies
\begin{equation}\label{nf3:full-signed-six}
 N_{F,\mathbb K}(V):=
 \#\{\mathbf t\in\mathcal O_{\mathbb K}^6:F(\mathbf t)=0,\
             |t_{i,v}|\le V^{1/d}\ (i\le6,\ v\mid\infty)\}
       \ll_{\mathbb K,\varepsilon}V^{3+\varepsilon}\qquad(V\ge1).
\end{equation}
The box includes both signs, every complex direction and the pairing
planes. Enlarging it by a fixed factor and clearing fixed ideal
denominators preserve the bound. The mixed-box hypothesis giving
positive density at the beginning of this section controls a different
count.
Then the sieve saving extends to $L\le M^{3/2-\eta}$ for each fixed
$\eta>0$. To see this put $H_0=L^2/M$. Frequencies for distinct
primitive characters over these prime squares are separated on the
trace-dual torus by $\gg_{\mathbb K}L^{-4/d}$, by the norm of the product
of two denominators. The $d$-dimensional additive large sieve at these
points for the cube of a truncated cubic sum is bounded by
$O_\varepsilon(M^\varepsilon(L^4+H_0^3)H_0^{3+\varepsilon})$.
Indeed its integral-basis frequencies have height
$O(H_0^{3/d}(\log M)^3)$ and its sum of squared coefficients is bounded
by the assumed signed Fermat count. The usual rectangular large sieve
follows by applying a product band-limited majorant in the dual
coordinates. The finitely many fixed lattices $\mathfrak b$ may be
treated separately. For the prime-dependent Gaussian weights truncate
at height $H_0^{1/d}\log M$ in a fixed basis. Partial summation in every
coordinate represents each weight by averages of initial boxes with
the mixed derivative as a signed density, also on the faces of the
bounding box. The absolute densities are bounded by common multiples
of $H_0^{-j/d}$ on each $j$-dimensional face, with total mass
$O((\log M)^d)$, since the normalized Gaussian shapes lie in a compact
family. H\"older's inequality permits application of the same large-sieve
estimate to each box, uniformly. As $H_0^3\le L^4$ the resulting
contribution is $O_\varepsilon(M^{3+\varepsilon}L^2)$, as required.

Over $\mathbb Q$, an independent fourth-moment argument proves the
sieve saving unconditionally on every range $1\le L\le
Z^{3/2-\eta}$, $\eta>0$ \cite[Proposition~12.4]{audClassical}.
The argument above for general $\mathbb K$ uses
\eqref{nf3:full-signed-six} to reach this longer range.

\subsection{Short reciprocal means over a number field}\label{nf3:shortmean}
Let $V\subset\mathbb P^5_{\mathbb K}$ be the Fermat cubic fourfold,
$D=F^\vee$, and $V_{\mathbf c}$ its hyperplane section. Write
$a_{\mathbf c}(\mathfrak n)$ for the ideal coefficients of either
$1/L(s,V_{\mathbf c})$ or
$1/(\zeta_{\mathbb K}(2s)L(s+1/2,V)L(s,V_{\mathbf c}))$, with
middle primitive cohomology normalized to weight zero. These are
the analogues of the coefficients in \cite[(1.6)]{c3:wang}.
Denote the embedding-height bound by $Z$ and put $A=Z^d$.
\begin{theorem}\label{nf3:reciprocal}
For every fixed number field, $2\le M\le A$, and arbitrary
$|b_{\mathfrak n}|\le1$, one has
\[
 \sum_{\substack{\|\mathbf c\|_\infty\le Z\\D(\mathbf c)\ne0}}
 \left|\sum_{M\le N\mathfrak n<2M}
          b_{\mathfrak n}a_{\mathbf c}(\mathfrak n)\right|^2
       \ll_{\mathbb K} A^6M\big(\log(2A)\big)^{C_{\mathbb K}} .
\]
If every prime ideal dividing $\mathfrak n$ has norm at least
$A^\eta$, for fixed $\eta>0$, the logarithmic factor may be removed.
It may also be removed for sums restricted to prime ideals.
\end{theorem}
\begin{proof}
Write $q=N\mathfrak p$ and omit finitely many primes of the field
and the fixed varieties. At $\mathfrak p$ define the residue-field
incidence trace
\[
 f_{\mathfrak p}(\mathbf c)=q^{-3/2}\sum_{x\in V(\mathbb F_q)}
                   (\mathbf1_{\mathbf c\cdot x=0}-1/q).
\]
Its normalized Fourier coefficient is $q^{-5/2}$ at nonzero
zeros of $F$ and zero elsewhere, for any primitive additive
character. Parseval gives
$\mathbb E|f_{\mathfrak p}|^2=1+O(q^{-2})$, from the point count
of the smooth fourfold. For squarefree $\mathfrak n$ write
$f_{\mathfrak n}=\prod_{\mathfrak p\mid\mathfrak n} f_{\mathfrak p}$.

We first bound the mean for $f_{\mathfrak n}$ with
$N\mathfrak n\asymp B$, in one residue class modulo
$\mathfrak a$ with $(\mathfrak n,\mathfrak a)=1$. Put
$X=A/N\mathfrak a$ and suppose $B\le X$, allowing endpoint
scales comparable with one. The asserted bounds in this class are
$X^6B$ times a log power, or $X^6B$ under roughness. Take a single
fixed positive product Gaussian $W_Z(\mathbf c)=W(\mathbf c/Z)$
which is at least one on the coefficient box. For each progression
$\mathbf c\equiv\mathbf r\bmod\mathfrak a$ and
$S(\mathbf c)=\sum b_{\mathfrak n}f_{\mathfrak n}(\mathbf c)$,
first use
\[
 \sum_{\substack{\mathbf c\equiv\mathbf r\\\|\mathbf c\|_\infty\le Z}}
           |S(\mathbf c)|^2
 \ \le\ \sum_{\mathbf c\equiv\mathbf r}W_Z(\mathbf c)|S(\mathbf c)|^2 .
\]
The same upper bound holds on restricting the left side to
$D(\mathbf c)\ne0$, since the incidence traces are defined everywhere.
Now expand the right side, keeping $b_{\mathfrak n}\overline
{b_{\mathfrak m}}$ and the progression phase until Poisson summation.
In the diagonal, each class modulo $\mathfrak a\mathfrak n$ has
Gaussian mass $O_{\mathbb K}((1+X/N\mathfrak n)^6)$ by packing
on dyadic dilates; completing modulo $\mathfrak n$ gives
$O_{\mathbb K}(X^6)$ per ideal from the bounded product of local
mean squares.
For distinct $\mathfrak n,\mathfrak m$ put
$\mathfrak g=(\mathfrak n,\mathfrak m)$,
$\mathfrak u=\mathfrak n/\mathfrak g$,
$\mathfrak v=\mathfrak m/\mathfrak g$, $Q=B/N\mathfrak g$.
At common primes the Fourier coefficients of $|f_{\mathfrak p}|^2$
are bounded by its mean, whose product is bounded. The product at
exclusive primes is $O(Q^{-5})$ and its frequencies are supported
on nonzero cubic zeros. Let
$\mathfrak l=\mathfrak g\mathfrak u\mathfrak v$.
In Gaussian Poisson summation the frequencies lie in
$(\mathfrak a\mathfrak l\mathfrak D)^{-1}$. Choose a fixed integral
class representative $\mathfrak b$ so that
$\mathfrak b\mathfrak a\mathfrak l\mathfrak D=(\pi)$, with $\pi$
balanced. The frequencies are $\mathbf t/\pi$ with $t_i\in\mathfrak b$,
and at each exclusive prime the Fourier variable is $\mathbf t$
times a scalar unit. Here the primes supporting the finitely many
class representatives and the different have also been omitted.
Thus $\mathfrak u\mathfrak v\mid (F(\mathbf t))$, and $\mathbf t=0$
is excluded. Each $\mathbf t$ specifies all local Fourier variables
uniquely. The prefactor is $O_{\mathbb K}(X^6 Q^{-5})$ and the
Gaussian volume parameter of each $t_i$ is
$T\asymp_{\mathbb K}B^2/(X N\mathfrak g)$. The finite choices of
$\mathfrak b$ may be summed separately.

If $F(\mathbf t)\ne0$, its norm is $O_{\mathbb K}(T^3)$ inside a
fixed dilate of the box. The number of ordered rough pairs
$\mathfrak u,\mathfrak v$ dividing this value is bounded by a
constant depending on $\eta,\mathbb K$ and the dilate; on dyadic
dilates the resulting exponential in the dilation index is absorbed
by the Gaussian. Without roughness the ideal divisor function
$\tau_{\mathbb K}((F(\mathbf t)))^2$ bounds this multiplicity, and
its mean on such a box is $O_{\mathbb K}(\log(2T)^{C_{\mathbb K}})$
for $T\ge1$. To see this, group prime factors of norm at most $T$
into successive maximal products of norm at most $T$. The norm of
the value is $O_{\mathbb K}(T^3)$, so the number of groups is
bounded and larger prime factors contribute a bounded factor.
Hence the divisor weight is majorized by a constant times a sum
of fixed powers of $\tau_{\mathbb K}(\mathfrak d)$ over divisors
of norm at most $T$. The zero count of $F$ modulo
$\mathfrak d$ is $O_{\mathbb K}((N\mathfrak d)^5)$: the primitive
solutions lift smoothly at good primes, with first count
$q^5+O(q^3)$, and scaling out a common uniformizer gives the ratio
$q^{-3}$ at levels at least three. At the finitely many primes
over three, restricting a unit coordinate modulo
$\mathfrak p^{v_{\mathfrak p}(3)+1}$ makes its cubic map a scaled
bijection onto a coset modulo
$\mathfrak p^{2v_{\mathfrak p}(3)+1}$, by Hensel's lemma.
Packing in residue classes and summing the divisor weights over
$N\mathfrak d\le T$ proves the claim. Thus the total off-diagonal
at nonzero values is $O_{\mathbb K}(B^7\log(2A)^{C_{\mathbb K}})$,
using $\sum_{N\mathfrak g\le2B}1/N\mathfrak g\ll_{\mathbb K}\log(2B)$.
In the rough case the last sum is bounded by the Mertens product
over $[A^\eta,2A]$ and gives $O_{\mathbb K,\eta}(B^7)$.

At $F(\mathbf t)=0$ the fourth and eighth cubic moments
\cite{nf3:kalinka} give $O_{\mathbb K}(T^{7/2+1/4})$ for $T\ge1$,
and Gaussian decay handles $T<1$. There are $O_{\mathbb K}(Q^2)$
pairs per $\mathfrak g$. Since $T\ll_{\mathbb K}Q$, summing
$X^6Q^{3/4}$ over $N\mathfrak g\le2B$ gives $O_{\mathbb K}(X^6B)$.

It remains to replace incidence traces by the local coefficients.
The local roots have modulus at most one after normalization, by
weight-monodromy for smooth complete intersections
\cite[Theorem 9.6]{nf3:scholze}. At a good prime of the section
$a_{\mathbf c}(\mathfrak p)=f_{\mathfrak p}(\mathbf c)+h_{\mathfrak p}$
with $h_{\mathfrak p}=O(q^{-1/2})$ independent of $\mathbf c$.
The error $e_{\mathfrak p}=a_{\mathbf c}(\mathfrak p)
-f_{\mathfrak p}(\mathbf c)-h_{\mathfrak p}$ is supported on
$\mathfrak p\mid(D(\mathbf c))$. Its squared supremum over lifts,
summed over the residue vectors, is $O(q^5)$. Indeed the smooth
part of the dual with nonzero coordinates has a unique ordinary
contact by Gauss duality and the diagonal Hessian. The resulting
one-node cubic threefold has point count $\#\mathbb P^3+O(q^{3/2})$:
projection from the node has tangent quadric a smooth quadric surface
$Q_2$ and exceptional directions a smooth $(2,3)$ curve $C$.
Singularity of $C$ would give another singular point of the section.
The exact count is $1+\#\mathbb P^3-\#Q_2+q\#C$ and the curve's
genus is four. At the other $O(q^4)$ nonzero vectors of the dual,
isolated singularities give $f_{\mathfrak p}=O(\sqrt q)$ by
\cite{nf3:hooleycontact}; the Gauss map is finite. At the origin
$f_{\mathfrak p}=O(q^{5/2})$.

Expand the squarefree products of $f_{\mathfrak p}+e_{\mathfrak p}$
over the product $\mathfrak s$ of the error primes. At each
$\mathbf c$ use Cauchy's inequality with weights $1/N\mathfrak s$
over $\mathfrak s\mid(D(\mathbf c))$. Their sum costs a fixed log
power, since $|N D(\mathbf c)|\ll A^{\deg D}$; under roughness it
is bounded. In each residue class modulo $\mathfrak s$ apply the
geometric estimate with $X=A/N\mathfrak s$ and
$B=M/N\mathfrak s$. The squared error majorants sum over classes
to at most $C^{\omega(\mathfrak s)}(N\mathfrak s)^5$. After the
complementary Cauchy weight $N\mathfrak s$, summing the costs
$A^6 M C^{\omega(\mathfrak s)}/N\mathfrak s$ gives another
log power, or a constant under roughness. The deterministic
$h_{\mathfrak p}$ are restored by Minkowski at cost
$\prod_{N\mathfrak p\le2M}(1+C/N\mathfrak p)$ in the square root
of the mean. Separating the squarefull part of an ideal costs the
same type of Euler product, as the reciprocal local factors are
polynomials of bounded degree with bounded coefficients. These
Euler products are bounded for rough factors, and the finite omitted
primes may be reinstated at bounded cost.

For prime indices alone, $M\ge\sqrt A/2$ is covered by roughness
for all large $A$. At smaller $M$ the products of any two prime
norms under consideration are at most $A$. Complete residue-class
counts then control the local error square mean directly:
$|e_{\mathfrak p}|$ has residue mean $O(1/q)$, by the same three
strata, and square mean $O(1/q)$. Distinct prime contributions
factor in the complete box majorant. For the geometric mean the
diagonal is $O(A^6M)$, while distinct prime frequencies have
Poisson scale $T\ll M^2/A\ll1$, with Gaussian decay and at most
$O((1+j)^2)$ possible prime pairs at each nonzero value on the
$j$th dyadic dilate. These estimates give $O_{\mathbb K}(A^6 M)$;
the deterministic prime corrections have sum $O_{\mathbb K}(\sqrt M)$.
\end{proof}

The two polar corrections give a sharper result for every ideal
coefficient in a slightly shorter range.
\begin{theorem}\label{nf3:logfree}
Use the corrected reciprocal
$1/(\zeta_{\mathbb K}(2s)L(s+1/2,V)L(s,V_{\mathbf c}))$.
For each fixed $0<\sigma<1$ the mean in
Theorem~\ref{nf3:reciprocal} is $O_{\mathbb K,\sigma}(A^6M)$
for all $M\le A\exp(-(\log(2A))^\sigma)$, with arbitrary bounded
weights on all ideals. Write $\Phi_{\mathbf c}(s)$ for this
corrected reciprocal. For $1<\sigma_0<2$, an entire $f$ rapidly
decreasing uniformly on $0\le\Re s\le2$, and
$1\le N\le A^{5/6-\varepsilon}$, its Mellin mean satisfies
\[
 \sum_{\substack{\|\mathbf c\|_\infty\le Z\\D(\mathbf c)\ne0}}
 \left|\frac{1}{2\pi i}\int_{(\sigma_0)}
       \Phi_{\mathbf c}(s) f(s)N^s\,ds\right|^2
 \ll_{\mathbb K,\varepsilon} A^6 N J(f),\qquad
 J(f)=\sup_{0\le x\le2}\int_{\mathbb R}(1+|t|)^2|f(x+it)|^2\,dt.
\]
This is the number-field form of the test in
\cite[Conjecture~6.4]{c3:wang} in the stated range.
\end{theorem}
\begin{proof}
For $q=N\mathfrak p$ set $G_{\mathfrak p,0}=1$ and
$G_{\mathfrak p,1}=f_{\mathfrak p}$. For $e\ge2$ use the actual
corrected coefficient at $\mathfrak p^e$ on smooth residue
classes and zero elsewhere. These are bounded functions of the
residue class, zero for $e>R$ for some fixed $R$. For residue
degree one, equations (6.5)--(6.10) of \cite{c3:wang} give
$\mathbb E G_{\mathfrak p,2}=O(1/q)$: the unrestricted
reciprocal has restricted mean $1+O(1/q)$ at $\mathfrak p^2$,
the zeta correction subtracts the smooth-reduction probability,
and the shifted fourfold term contributes $O(1/q)$ since the
unrestricted prime mean is $O(q^{-1/2})$. The residue count is
over the very same $\mathbb F_p$ hyperplanes even if the
characteristic-zero field is larger. In every residue degree,
$\mathbb E G_{\mathfrak p,1}=0$ and
$\mathbb E|G_{\mathfrak p,1}|^2=1+O(q^{-2})$ by Parseval.
The primes of residue degree at least two have
$\sum 1/q<\infty$, so boundedness of $\mathbb E G_{\mathfrak p,2}$
suffices there.

Extend $G$ multiplicatively to ideal indices. In one progression
majorize the whole squared polynomial in the $G$'s by the same
weighted mean as above, and expand the latter. Use the parameters
$B,X$ and balanced Poisson formula from
Theorem~\ref{nf3:reciprocal}. The zero mode has total $O(X^6B)$.
Indeed, let $m_{\mathfrak p}(e,f)=
|\mathbb E G_{\mathfrak p,e}\overline{G_{\mathfrak p,f}}|$,
and form $h_{\mathfrak p}$ by replacing $m_{\mathfrak p}(1,1)$
with $(m_{\mathfrak p}(1,1)-1)_+$. Then
$m_{\mathfrak p}(e,f)\le\sum_{j\le\min(e,f)}
h_{\mathfrak p}(e-j,f-j)$, with $j\ge0$, and
$\prod_{\mathfrak p}\sum h_{\mathfrak p}(e,f)q^{-(e+f)/2}$
converges by the preceding moments and boundedness in total
degree at least three. The common ideal extracted by this local
convolution has at most $O_{\mathbb K}(B/\sqrt{N\mathfrak u
N\mathfrak v})$ choices when its products with $\mathfrak u$ and
$\mathfrak v$ have norm in $[B,2B]$. This proves the zero-mode bound.

For the other modes let $\mathfrak u$ consist of the primes of total
multiplicity one in the two indices, and $\mathfrak v$ consist
of all other primes, each with exponent one. Write
$U=N\mathfrak u$, $V_0=N\mathfrak v$. Then $U V_0^2\le4B^2$.
The local product has modulus at most
$U^{-5/2}C^{\omega(\mathfrak v)}$ and requires
$\mathfrak u\mid(F(\mathbf t))$. For $U\asymp U_0$ the dual volume
parameter is $T\asymp_{\mathbb K}U_0 V_0/X$. At $F(\mathbf t)\ne0$
the choices of $\mathfrak u$ and distribution of its primes
are bounded by $\tau_{\mathbb K}((F(\mathbf t)))^2$. Summing the
resulting $X^6 U_0^{-5/2}T^6$ over dyadic $U_0\ll B^2/V_0^2$ and
then over $\mathfrak v$ gives at most $B^7(\log(2A))^{C_{\mathbb K}}$,
by the divisor mean in Theorem~\ref{nf3:reciprocal}. At
$F(\mathbf t)=0$, $\mathbf t\ne0$, there are
$O_{\mathbb K}(U_0(\log(2U_0))^C)$ allocations for $\mathfrak u$.
The sixth-moment estimate $O(T^{7/2+1/4})$ and the same summations
give $O_{\mathbb K}(X^6B(B/X)^{7/2+1/4}(\log(2A))^{C_{\mathbb K}})$,
since $\sum_{N\mathfrak v\le2B}C^{\omega(\mathfrak v)}
(N\mathfrak v)^{-3/4}\ll_{\mathbb K} B^{1/4}(\log(2B))^{C'}$.
The ideal divisor sums used here follow by bounding $C^{\omega}$
by a fixed divisor function and counting ideals of norm at most
$x$ by $O_{\mathbb K}(x)$. Thus this $G$-mean is log-free whenever
$X/B\gg\exp(\tfrac13(\log(2A))^\sigma)$.

On smooth residue vectors the prime coefficient equals
$f_{\mathfrak p}+q^{-5/2}$, by the identity
$\#\mathbb P^4/q-\#\mathbb P^3=1/q$ after the fourfold correction.
The deterministic terms may be restored by a convergent Minkowski
sum bounded by $\sum_{\mathfrak t}(N\mathfrak t)^{-3}$. All
remaining local errors are supported on $\mathfrak p\mid(D(\mathbf c))$.
Their squared majorant has residue mean $O(1/q)$ in every positive
exponent, as does $\mathbf1_{\mathfrak p\mid D}|G_{\mathfrak p,e}|^2$,
by the three strata of the preceding proof.

Expand over the ideal $\mathfrak s$ of local error factors.
Use Cauchy with weights $(N\mathfrak s)^{-\theta}$,
$\theta=1-\sigma/4$. In the first Cauchy factor we can restrict
the auxiliary divisor $\mathfrak s'$ of $(D)^R$ to
$N\mathfrak s'\le Y=\exp(\tfrac12(\log(2A))^\sigma)$, up to an
additive constant: the tail costs at most
$Y^{-\sigma/4}\prod_{\mathfrak p\mid(D)}
(1-q^{-1+\sigma/2})^{-1}=O_{\mathbb K,\sigma}(1)$, by splitting
at $q=\log(2A)$ and using the bound for the norm of $D$.
Split from the $G$-polynomial the primes of $\mathfrak s'$.
Their factor $\mathfrak k$ has at most
$(R+1)^{\omega(\mathfrak s')}$ choices. In each class modulo
$\mathfrak l=\operatorname{rad}([\mathfrak s,\mathfrak s'])$
the remaining geometric sum has $X=A/N\mathfrak l$ and
$B=M/N(\mathfrak s\mathfrak k)$, with $X/B\ge A/(MY)$.
The joint residue mean of all squared local majorants is at most
$C^{\omega(\mathfrak l)}/N\mathfrak l$. Weighted Cauchy
therefore costs at most a constant times $A^6M$ multiplied by
\[
 \sum_{\mathfrak s,\mathfrak s'}
 \frac{C^{\omega(\mathfrak s)+\omega(\mathfrak s')}}{
 (N\mathfrak s)^{1-\theta}(N\mathfrak s')^\theta
 N\operatorname{rad}([\mathfrak s,\mathfrak s'])}<\infty .
\]
At each prime the nonconstant terms are
$O(q^{-2+\theta}+q^{-1-\theta}+q^{-2})$, proving convergence.

Finally expand the absolutely convergent series on $\Re s=\sigma_0$.
For each ideal shift the inverse Mellin integral of $f(s)
(N/N\mathfrak n)^s$ to the lines $0$ and $2$. Cauchy's inequality
there bounds its weight by $O(J(f)^{1/2}\min(1,(N/N\mathfrak n)^2))$.
Sum the result just proved up to
$A\exp(-(\log(2A))^{1/2})$. Above that scale a block of norm $B$
has square-root mean $O_\delta(\min(A^{1/2}B^{3+\delta},
A^3 B^{1+\delta}))$. For the first bound enlarge the parameter box
to volume $Y=B^{1+\delta/4}$ and use projective invariance of the
coefficients. On each projective line the multiplicity in the
original box is $O_{\mathbb K}(A/Y)$ times that in the larger box.
Indeed take representatives whose coordinate ideals generate one
of finitely many fixed ideal classes' representatives, and balance
their placewise coordinate suprema by units. Allowed scalars lie
in one of finitely many fractional ideals and the resulting
placewise boxes have bounded eccentricity. If the small box
contains a nonzero scalar, its normalized radius is bounded below,
and lattice packing and covering give the asserted comparison
under dilation. The second estimate is the divisor bound.
Splitting at $B=A^{5/4}$ gives a tail
$O_\delta(N^2 A^{7/4+2\delta})$, which is
$O(A^3\sqrt N)$ in the stated range.
\end{proof}

\begin{theorem}\label{nf3:longprime}
For prime ideal indices, the log-free conclusion
of Theorem~\ref{nf3:reciprocal} extends to $M\le A^{7/6-\varepsilon}$.
If the full signed six-cube estimate \eqref{nf3:full-signed-six}
holds for every positive exponent loss,
the range extends to $M\le A^{3/2-\varepsilon}$.
\end{theorem}
\begin{proof}
It suffices to treat $A\le M$. Choose one nonnegative Schwartz box
majorant $W(\mathbf c/Z)$ with smooth compactly supported Fourier
transform $\widehat W$. First majorize the entire absolute square of
the prime polynomial in the incidence traces with this weight on
$\mathbf c\ne0$. It bounds the original geometric mean since
$D(\mathbf c)\ne0$ excludes the origin. Expand only after this step,
using the same $W$ for every pair and keeping the puncture.
Its diagonal is $O_{\mathbb K}(A^6M)$: for a nonsingular
$\mathbf c$ in a dyadic dilate of the box only $O_{\mathbb K}(1+j)$
primes of norm $\asymp M$ divide $D(\mathbf c)$ on the $j$th
dilate, and their nonzero residues give $|f_{\mathfrak p}|\ll\sqrt M$.
On $D(\mathbf c)=0$ each squareclass of the nonzero coordinates
occurs at least twice, by independence of the corresponding radicals.
There are $O_{\mathbb K}(A^3(\log(2A))^C)$ such vectors. Indeed for
a group of $i\ge2$ equal squareclasses the coordinate ideals have
the form $\mathfrak d\mathfrak a_j^2$ with the same squarefree
$\mathfrak d$. Summing $(A/N\mathfrak d)^{i/2}$ over
$N\mathfrak d\ll A$ and allowing a fixed log power for bounded
generators gives $O_{\mathbb K}(A^{i/2}(\log(2A))^C)$.
The nonzero vectors all of whose coordinates belong to $\mathfrak p$
cost $O_B(M^6(A/M)^B)$ by lattice packing and Schwartz decay.
For the comparison with the actual coefficients restrict to
nonsingular parameters. The deterministic shifts sum to $O(\sqrt M)$.
In the $j$th dilate there are $O_{\mathbb K}(1+j)$ bad primes per
parameter, so Cauchy's inequality and the nonzero-residue bound give
$O_{\mathbb K}((1+j)^2M)$ pointwise for the square of that part.
The zero residues are charged in aggregate by the preceding packing
estimate with the additional factor $1+j$ absorbed by Schwartz decay.
Thus the comparison costs $O_{\mathbb K}(A^6M)$.

For distinct prime ideals $\mathfrak p,\mathfrak q$ of norms
$u,v\asymp M$, choose the fixed representative $\mathfrak b$ and
balanced $\pi$ with $(\pi)=\mathfrak b\mathfrak p\mathfrak q
\mathfrak D$, as in Theorem~\ref{nf3:reciprocal}. Poisson summation
puts the frequencies in $\mathfrak b^6/\pi$. In the numerator
$\mathbf t\in\mathfrak b^6$ we require
$\mathfrak p\mathfrak q\mid(F(\mathbf t))$ and nonzero residues
at the two primes. Dropping the latter over all pairs costs
$O_{\mathbb K}(M^2(1+T/M)^6)$ before the common prefactor
$O_{\mathbb K}(A^6M^{-5})$, where $T=M^2/A$. Set $C=T^3/M^2$ and
$\mathfrak j=\mathfrak b^2\mathfrak D^{-1}$.
The divisibility condition is exactly
\[
 F(\mathbf t)=\pi m,\quad m\in\mathfrak j,
 \qquad \pi\mathfrak j=\mathfrak b^3\mathfrak p\mathfrak q.
\]
Both sides are in $\mathfrak b^3$. Use the character integral on
the trace-dual torus of this lattice, and a smooth cutoff in
$m/C^{1/d}$ identically one throughout the necessary range.
Integrals are normalized by the corresponding lattice covolumes.
The scales in this formula are volume parameters. The finite
classes $\mathfrak b$ are treated separately.

On this torus put major arcs
$Q_\gamma B(T^{3/d}(\alpha-\gamma))\le\Delta=C T^{-\eta}$ with
$Q_\gamma$ the norm of the denominator ideal relative to
$\mathfrak b^3$. Here $B(z)=\prod_v(1+|z_v|)^{d_v}$.
They are disjoint, since two distinct centers have difference of
absolute norm $\gg_{\mathbb K}(Q_1Q_2)^{-1}$ and
$\Delta^2=o(T^3)$. To bound the trilinear sum on the complement,
group prime pairs by $\pi$. Its ideal determines the unordered
pair, so each element in a fixed height box of volume $O(M^2)$ is
counted at most twice. The compact family of shapes of $\pi/M^{2/d}$
may be separated from the six-variable weight by a smooth Fourier
series in its $d$ real parameters; its coefficients have rapidly
decreasing fixed derivative norms. This introduces arbitrary
bounded pair weights, which the following argument permits.

Sum first in $m$. Partition the resulting frequencies on the
trace-dual torus of $\mathfrak j$ into cells of side $C^{-1/d}$.
The sum of the suprema of the smooth transforms over all cells
is $O_{\mathbb K}(C)$ by Poisson summation. The maximal occupancy
is bounded by the number of pairs $(z,b)$ with $z$ in a fixed
ideal of height $O(M^{2/d})$ and
$H_\infty(\alpha z-b)\ll C^{-1/d}$ for $b\in\mathfrak j^\vee$,
by taking differences. Here one may take
$z\in\mathfrak b\mathfrak D$ since every $\pi$ belongs to that
ideal. The rank-two lattice argument of
Theorem~\ref{nf3:four} bounds the occupancy by $O_{\mathbb K}(M^2/C)$
if two pairs are independent over $\mathbb K$: the rescaled
lattice has covolume comparable with $C/M^2$ and bounded covering
radius. In rank one write $b=\gamma z$. The index in
$\mathfrak b\mathfrak D$ of the allowed $z$ with
$\gamma z\in\mathfrak j^\vee$ is exactly $Q_\gamma$, because
$(\mathfrak b\mathfrak D)\mathfrak j=\mathfrak b^3$.
Coordinatewise packing bounds the count by
$O_{\mathbb K}(1+M^2/(Q_\gamma
B(T^{3/d}(\alpha-\gamma))))$. Consequently the minor trilinear
sum is $O_{\mathbb K}(M^2 T^\eta)$ since $C\le M^2$.
The integral of the absolute six-variable cubic sum is
$O_{\mathbb K,\delta}(T^{\rho+\delta})$ for $\rho=7/2$,
or $\rho=3$ under the stated hypothesis, by partial summation,
H\"older and the fourth and eighth moments. Thus the minor
contribution before the prefactor is $O(M^2 T^{\rho+\eta+\delta})$.

On an arc write $\alpha=\gamma+\beta$. If
$\gamma\pi\notin\mathfrak j^\vee$, the algebraic norm gives,
for all $b\in\mathfrak j^\vee$,
\[
 B(C^{1/d}(\gamma\pi-b))\gg_{\mathbb K} C/Q_\gamma.
\]
On the other hand $B(C^{1/d}\beta\pi)\ll_{\mathbb K}\Delta/Q_\gamma$.
Since $B(x+y)\le B(x)B(y)$, the transformed $m$ frequencies
have $B(C^{1/d}(\alpha\pi-b))\gg C/\Delta$. Poisson in $m$
and rapid decay show that these terms are negligible. The
remaining centers have denominator dividing
$\mathfrak p\mathfrak q$ by $\pi\mathfrak j=
\mathfrak b^3\mathfrak p\mathfrak q$, so only the central center
and those of denominator $\mathfrak p$ or $\mathfrak q$ occur.
For a center of denominator $\mathfrak p$, finite-field cubic
completion and the bound $2\sqrt u$ for each linearly twisted
complete sum bound the weighted six-variable cubic sum by
\[
 O_{\mathbb K}\big(T^6u^{-3} B(T^{3/d}\beta)^{B_1}\big)
\]
for a fixed $B_1$. The dual sum of the differentiated Fourier
weights is bounded by this expression by lattice packing,
since $T\gg u$. The zero mode of the changed $m$-transform
integrates this polynomial at cost $O_{\mathbb K}(1/(uv))$.
All other $m$-modes are negligible: for $b\ne\gamma\pi$ in
$\mathfrak j^\vee$, the previous product estimate has numerator
$\gg C$ instead of $C/Q_\gamma$. There are $O(u)$ centers.
Including the analogous $\mathfrak q$ terms gives $O(T^6/M^2)$
over the pairs before the prefactor.

On the central arc the smooth $m$-transform restricts
$B(T^{3/d}\beta)$ to an arbitrarily small positive power of $T$
up to a negligible tail, and Poisson in $\mathbf t$ leaves only
the real integral. Fourier inversion gives the full continuous
mass of $\widehat W$ on $\mathbf t\in\mathfrak b^6$, divided
by $uv$: this last factor is exactly the index
$[\mathfrak b^3:\pi\mathfrak j]^{-1}$. With self-dual Haar measure
the two Poisson lattice factors multiply to $(uv)^6$, because
$\operatorname{covol}(\mathcal O_{\mathbb K})^2=N\mathfrak D$,
$\operatorname{covol}(\mathfrak b)=N\mathfrak b\,
\operatorname{covol}(\mathcal O_{\mathbb K})$ and
$|N\pi|=N\mathfrak b\,uv N\mathfrak D$.
After the local Fourier normalizations the central term is
$(uv)^{5/2}W(0)$. For each pair subtract
$W(0)f_{\mathfrak p}(0)\overline{f_{\mathfrak q}(0)}$ from this
mass, with the same multiplier $b_{\mathfrak p}\overline
{b_{\mathfrak q}}$ on both terms. Since
$f_{\mathfrak p}(0)=u^{5/2}(1+O(u^{-2}))$, the residual is
$O(M^3)$ per pair and $O(M^5)$ in total, also for complex weights.

Finally $T^{\rho+\eta+\delta}\ll M^4$ and $T^6\ll M^8$
hold in the stated ranges, on taking $\eta,\delta$ sufficiently
small. The prefactor $O_{\mathbb K}(A^6 M^{-5})$ now gives the
claimed estimate.
\end{proof}

\begin{corollary}\label{nf3:longrough}
For $0<\eta\le1$ and $0<\delta<\eta/6$ the log-free
$A^\eta$-rough conclusion of Theorem~\ref{nf3:reciprocal} extends
to $M\le A^{1+\delta}$ in every number field.
\end{corollary}
\begin{proof}
Write $Y=A^\eta$, $T=M^2/A$, and assume $M\ge A$.
Any nontrivial squarefull part, or nontrivial product of deterministic
prime corrections, leaves length at most $M/Y\le A$ and is covered
by the bounded Euler products for rough ideals. For squarefree indices
expand over the bad-prime factors $\mathfrak s$. There are
$O_\eta(1)$ such divisors per nonsingular vector in the box.
Their squared local majorants have residue mean
$O(C^{\omega(\mathfrak s)}/N\mathfrak s)$ as in
Theorem~\ref{nf3:reciprocal}. For $\mathfrak s\ne(1)$ enlarge the
parameter volume from $A$ to $M$ before applying the geometric mean
in each progression modulo $\mathfrak s$. All these errors cost at
most $O_\epsilon(M^7Y^{-1+\epsilon})=O(A^6M)$, by summing
$C^{\omega(\mathfrak s)}/(N\mathfrak s)^2$ above $Y$.

For the geometric covariance on nonzero vectors use the same
majorant for the entire square before expanding. M\"obius
inversion bounds pairs sharing a prime by the sum of the square
means over indices both divisible by a common rough squarefree
$\mathfrak d\ne(1)$. Factor out $|f_{\mathfrak d}|^2$ and bound it
by $C^{\omega(\mathfrak d)}\prod_{\mathfrak p\mid\mathfrak d}
(1+\mathbf1_{\mathfrak p\mid(D)}|f_{\mathfrak p}|^2)$.
The empty subset costs $O(A^6M/N\mathfrak d)$ by the shorter
geometric mean. For a nonempty subset $\mathfrak s\mid\mathfrak d$
use its joint residue mean $O(C^{\omega(\mathfrak s)}/N\mathfrak s)$
and apply the geometric mean on these classes after enlarging the
parameter volume to $M$. The total is $O_\epsilon(M^7Y^{-1+\epsilon})$
as before since the rough harmonic sum over
$\mathfrak d/\mathfrak s$ is bounded.

For coprime squarefree pairs repeat the trilinear proof of
Theorem~\ref{nf3:longprime}. Each product ideal has $O_\eta(1)$
ordered factorizations. Omitting the nonzero-frequency requirements
at its primes costs at most
$O_\epsilon(T^6Y^{-5+\epsilon}+M^2(1+T^{7/2+\epsilon}))$
before the Poisson prefactor. At nonzero $F(\mathbf t)$ the other
rough factors all divide this value and have boundedly many choices;
there are $O_{\mathbb K}((T/q)^6)$ nonzero vectors in a multiple
of the given prime of norm $q\ll T$. For $F(\mathbf t)=0$ scale
such vectors through one of finitely many ideal classes and apply
the sixth moment at volume $T/q$, counting $O_{\mathbb K}(M^2/q)$
pairs of indices containing the prime. The minor bound is unchanged.
The resonant denominators now divide the product ideal. Each
nontrivial such norm $r$ exceeds $Y$ and there are $O_\eta(1)$
choices for each pair. Multiplicative cubic completion and packing
give $O_\eta(T^6 r^{-3}(1+r/T)^6)$ times a fixed polynomial in the
scaled arc displacement; the free-variable transform absorbs that
polynomial. Summing centers and pairs thus gives
$O_\eta(T^6(Y^{-2}+T^6/M^8))$, since
$r^4/T^6\le\Delta^4/T^6\le T^6/M^8$ for $r>T$.
The affine origin cancels the central term up to
$O_\eta(M^7Y^{-2})$ after normalization. The inequalities
$T^{7/2+\epsilon}\ll M^4$, $T^6Y^{-2}\ll M^6$ and
$T^{12}\ll M^{14}$ hold for sufficiently small $\epsilon$ in the
stated range, and give the result.
\end{proof}
\section{From the trace torus to representations}
\label{sec:nf-global}
Write $V=\mathbb K\otimes_{\mathbb Q}\mathbb R$,
$n=[\mathbb K:\mathbb Q]$, $\mathcal O=\mathcal O_{\mathbb K}$,
and $\mathfrak D$ for the different. Haar probability measure on
$\mathbb T=V/\mathfrak D^{-1}$ is dual to the lattice $\mathcal O$
under $e_{\mathbb K}(u)=\exp(2\pi i\operatorname{Tr}u)$. Denote
the local real dimension of an infinite place by $d_v$. By
Lemma~\ref{lem:hbal} we may work with balanced targets:
$P=|N\nu|^{1/(nk)}$, and $\nu/P^k$ belongs to a compact subset of
the positive archimedean cone on which all components are nonzero.
At the balanced target the root boxes satisfy
$|x_i|_v\ll_{\mathbb K,k,s}|\nu|_v^{1/k}$ at every infinite place.
Undoing the balancing unit preserves each of these ratios, so the
argument gives this conjugate-region bound as well as the absolute
norm bound.

\subsection{Ordinary sums}
Put $X=P^n$ and write
$f_w(\alpha)=\sum_{x\in\mathcal O}w(x/P)e_{\mathbb K}(\alpha x^k)$.
We use nonnegative smooth weights compactly supported away from zero
at each embedding, with positive real coordinates. Sums of finitely
many product weights suffice. For $\eta\in\mathbb K/\mathfrak D^{-1}$
denote its conductor by $\mathfrak q_\eta
=\{a\in\mathcal O:a\eta\in\mathfrak D^{-1}\}$ and put
$Q_\eta=N\mathfrak q_\eta$. Set
\[
 B(t)=\prod_{v\mid\infty}\max(1,|t_v|)^{d_v},\qquad
 \mathfrak M(Z)=\bigcup_{\eta}
 \{\eta+\beta: Q_\eta B(P^k\beta)\le Z\}.
\]
The union allows all lifts. Put $\mathfrak M_*=\mathfrak M(CX)$
with fixed sufficiently large $C$, and
$\mathfrak m_*=\mathbb T\setminus\mathfrak M_*$. The field
critical mean-value theorem is due to Wooley
\cite[Theorem~15.1 and Corollary~15.3]{NFWooley}; the resolution
of the rational main conjecture by Bourgain, Demeter and Guth
\cite{NFBDG} and by Wooley concerns the full coefficient system.
Here the leading coefficient alone is integrated on $\mathbb T$.

\begin{lemma}\label{lem:nf-ordinary}
For $k\ge3$ and $b>2k$ there is $\tau>0$ such that, for $Z\ge1$,
\begin{align}
 \sup_{\mathfrak m_*}|f_w|&\ll X^{1-1/[k(k-1)]+\varepsilon},
 \label{eq:nf-point}\\
 \int_{\mathfrak M_*\setminus\mathfrak M(Z)}|f_w|^b
 &\ll Z^{-\tau}X^{b-k},\qquad
 \int_{\mathfrak M_*}|f_w|^b\ll X^{b-k}.
 \label{eq:nf-prune}
\end{align}
At the critical even exponent $w_c=k(k+1)$ one also has
\begin{equation}\tag{\ref*{eq:nf-prune}a}\label{eq:aud-field-critical}
             \int_{\mathfrak m_*}|f_w(\alpha)|^{w_c}\,d\alpha
             \ll_\varepsilon X^{w_c-k-1+\varepsilon}.
\end{equation}
The constants may depend on $\mathbb K,k,b,w$.
\end{lemma}
\begin{proof}
We detail the restriction from the full coefficient torus. Let a
smooth sum at root volume $Y=T^n$ have degree $j\ge2$, top coefficient
$\theta$ and arbitrary lower coefficients. Let $B_0$ bound the number
of frequencies $j\theta h$ in any product box of side
$O(T^{1-j})$ modulo $\mathfrak D^{-1}$ for shifts $h$ in a fixed
box of side $O(T)$. With $b_0=j(j-1)$ one has
\begin{equation}\label{eq:nf-weyl-sample}
 |S|^{b_0}\ll Y^{b_0+\varepsilon}B_0/Y .
\end{equation}
For each of $\asymp Y$ translations $h$, move the summation variable
to an enlarged fixed box and use Fourier inversion of the translated
smooth weight. The absolute Fourier majorant is independent of $h$.
At each fixed Fourier variable the new coefficients of degrees
$1,\ldots,j-1$ give samples on $\mathbb T^{j-1}$. The top of these
coefficients is translated by $j\theta h$. In a grid of product
cells of side $T^{-i}$ at degree $i$, each cell has at most $O(B_0)$
samples. On each enlarged cell apply the one-dimensional fundamental
theorem of calculus successively in its real coordinates, followed
by H\"older: the supremum to the $b_0$th power is bounded by the
inverse cell volume times the sum of the integrals of all mixed first
derivatives to that power, each derivative scaled by its cell side.
Expanded boxes have bounded overlap. The inverse volume is
$Y^{b_0/2}$. All differentiated sums have bounded coefficients
on the enlarged root box; the term of degree $j$ becomes part of
their weights. The weighted mean-value theorem
\cite[Theorem~15.1]{NFWooley} at degree $j-1$ bounds each full-torus
integral by $Y^{b_0/2+\varepsilon}$. Jensen's inequality in the
translation and Fourier averages gives \eqref{eq:nf-weyl-sample}.
The implied constants are uniform in the lower coefficients.

For $\alpha\in\mathfrak m_*$, solutions of
$|k\alpha_v h+\xi_v-c_v|\ll P^{1-k}$,
$c\in\mathfrak D^{-1}$ and $|h_v|\ll P$, lie on a single affine
$\mathbb K$-line when there are at least two. In fact the determinant
of two differences belongs to $\mathfrak D^{-1}$ and has norm
$O(X^{2-k})$, which forces it to vanish. For large $P$ distinct pairs
with the same $h$ are impossible. With slope $\eta$, the available
$h$ lie in a coset of $\mathfrak q_\eta$, in a product box with
radii $O(\min(P,P^{1-k}/|k\alpha_v-\eta_v|))$.
The count in a box of volume $Y'$ in an ideal coset of norm $Q'$
is $O(1+Y'/Q')$: scale its radii to one, use $|Nz|\ge Q'$ for any
nonzero difference, and pack separated points. Hence the number
of shifts is $O(1+X/(Q_\eta B(P^k(k\alpha-\eta))))=O(1)$.
For the last equality use $Q_{\eta/k}\le k^nQ_\eta$ and
$\alpha\notin\mathfrak M_*$. Equation~\eqref{eq:nf-weyl-sample}
gives \eqref{eq:nf-point}.

To prove \eqref{eq:aud-field-critical}, keep the degree-$k$ frequency
$\alpha$ as a variable and repeat the translation and cell estimate
at exponent $w_c$. For each $\alpha\in\mathfrak m_*$ the lower
coefficient cells have bounded occupancy by the preceding argument.
Jensen in the shifts contributes $X^{-1}$ and the inverse volume of
a lower coefficient cell is $X^{k(k-1)/2}$. Integrate the resulting
nonnegative bound in $\alpha$ and enlarge this integration to the full
torus. At a fixed Fourier parameter the smooth cutoff and all its
scaled polynomial derivative weights on the enlarged root box are
bounded and independent of $\alpha$. The full $k$-coefficient mean
value at exponent $w_c$ is $O_\varepsilon(X^{k(k+1)/2+\varepsilon})$
for each such weight, by \cite[Theorem~15.1]{NFWooley}. The common
Fourier majorant is integrable. Multiplication of the three factors
gives $X^{k^2-1+\varepsilon}=X^{w_c-k-1+\varepsilon}$.

For completeness, Poisson summation at $\alpha=\eta+\beta$ on
$Q_\eta B(P^k\beta)\le CX$ gives
\begin{equation}\label{eq:nf-complete-bound}
 |f_w(\eta+\beta)|\ll_{\varepsilon}
 X Q_\eta^{-1/k+\varepsilon} B(P^k\beta)^{-1/2}.
\end{equation}
The normalized complete sum of $x^k$ with a linear additive twist
is $O_\varepsilon(Q_\eta^{-1/k+\varepsilon})$. At a prime of residue
characteristic prime to $2k(k-1)$, prime-level sums are $O_k(q^{-1/2})$
by Weil. At higher levels the unit shell has simple stationary points
and is $O_k(q^{-e/2})$ at level $e$. On the nonunit shell $x=\pi y$.
For $e>k$ this gives a factor $q^{-1}$ and a sum at level $e-k$
if the linear coefficient is divisible by $\pi^{k-1}$; otherwise
the reduced phase has unit derivative everywhere, so summation over
the last residue digit gives zero. For $2\le e\le k$, this shell
has mass at most $q^{-1}$. On iteration the normalized unit-shell
errors after multiplication by $q^{e/k}$ form a bounded geometric
sum, because $1/2-1/k>0$. This yields the prime-power bound with
a constant depending only on $k$. At the finitely many bad
completions the polynomial van der Corput estimate
\cite[Corollary~3.4]{NFCluckers} gives $O_{\mathbb K,k}(q^{-e/k})$
uniformly for integral linear coefficient: the polynomial is
Weierstrass regular of degree $k$ with unit leading coefficient.
The convention that the basic character in \cite{NFCluckers} is
trivial on the maximal ideal shifts the conductor index by one,
which costs a fixed factor at each of these finitely many places.
The Chinese remainder theorem gives the asserted estimate, absorbing
the bounded factor per prime by $Q_\eta^\varepsilon$.

At a real place, the second derivative on the support of the weight
is bounded below by a fixed multiple of the leading frequency.
At a complex place the real Hessian of
$\operatorname{Re}(tz^k)$ is nonsingular on the annulus for $t\ne0$.
Stationary phase thus gives the common factor $B(P^k\beta)^{-1/2}$
also with linear twists. Beyond the radii
$(1+P^k|\beta_v|)/P$, integration by parts gives rapid decay in
the normalized dual coordinates, with this same factor, by taking
additional derivatives. The nonzero dual frequencies lie in
$\mathfrak q_\eta^{-1}\mathfrak D^{-1}$; ideal packing as above
counts $O(1+Q_\eta B(P^k\beta)/X)$ in the effective box.
This proves \eqref{eq:nf-complete-bound}. There are at most $Q$
characters of conductor a given ideal of norm $Q$. In the integral
of its $b$th power insert $(QB(P^k\beta)/Z)^\tau$ on the annulus.
Choose $\tau>0$, then $\varepsilon>0$, with
$1-b/k+b\varepsilon+\tau<-1$ and $b/2-\tau>1$.
The ideal series and the product archimedean integral converge;
the change of variables $t=P^k\beta$ contributes $X^{-k}$.
The same proof without the extra factor bounds the whole major arcs.
\end{proof}

\subsection{A terminal set with an unrestricted cofactor}
The factors-and-fibres argument of Chow and Dietmann
\cite[Sections~2--5]{NFChow} marks one prime in each of finitely
many disjoint norm windows. Over $\mathbb K$ we take prime ideals in
the narrow principal class of residue degree one, and choose one
totally positive unit-balanced generator for every such prime
ideal. For any fixed $0<u<v$ the reciprocal norms of these primes
with $Y^u<N\mathfrak p<Y^v$ have a positive limit as
$Y\to\infty$, by the prime ideal theorem in narrow ideal classes
and partial summation. Prime ideals of higher residue degree have
convergent reciprocal norm sum and hence zero limiting contribution.

After all primary marks are removed let the remaining volume be $Y$.
The following construction applies whenever
$X^{c_1}\le Y\le X^{c_2}$ throughout the marks, with $c_1>0$ fixed
and all primary prime norms at least $X^{c_3}$ for a fixed $c_3>0$.
Constants here may depend on $k$ and the scheme of windows.
Choose $L$ large and
\[
 (L+1)^{-1}<a<b<L^{-1},\qquad c_2b<c_3 .
\]
Choose $L$ distinct auxiliary principal primes with norms between
$Y^a$ and $Y^b$, taken in a fixed order to avoid permutations.
Writing $D$ for the product of their chosen generators, allow a
free cofactor $z\in\mathcal O$ in a smooth annular box of volume
$Z=Y/|ND|$, positive at the real places. Then
\[
 Y^{1-Lb}\le Z\le Y^{1-La},\qquad 0<1-Lb<1-La<a.
\]
Each free cofactor has absolute norm $O(Z)$ and hence has no primary
or auxiliary prime divisor from a chosen window for sufficiently
large $X$. Primary primes are identifiable first from the ideal of
a root, and the auxiliary primes are precisely its $L$ largest
remaining prime factors, counted with multiplicity. Their norms are
larger than the norm of the free cofactor. The chosen generators then
determine $z$ uniquely. The coefficients of the resulting sum $F$
are nonnegative and bounded by one. A smooth weight bounded by one
and positive on an open box gives mass $\asymp X$: for each set of
primary marks the free box has $\asymp Z$ lattice points, and the
reciprocal norm sum over the $L$ distinct auxiliary primes has positive
lower limit. Diagonals have vanishing reciprocal mass. The same
calculation over the primary windows concludes the assertion. These
estimates hold for all sufficiently large $X$ and with any fixed
small dilation of the entire set. Units balancing finitely many
factors and a fixed terminal annulus keep the root radii within
fixed positive multiples of the nominal radius at every place.

This choice only changes the terminal cofactor. In all identities
obtained by removing primary marks, its volume $Y$ is the initial
volume divided by the norms of \emph{all} the primary marks, so
the same terminal weight occurs before and after removal. The fibre
identities use the avoidance of the primary windows, the support of
size $O(X)$ and bounded coefficients. All three properties hold here.
There is an additional useful consequence. If $D'$ is the product of
every marked generator, primary and auxiliary, the corresponding part
of $F$ has the form
\begin{equation}\label{eq:nf-terminal-form}
 \sum_{z\in\mathcal O}w_0(z/Z^{1/n})
                  e_{\mathbb K}(\alpha(D'z)^k),\qquad
 Z\gg X^{c_1(1-Lb)},\qquad |ND'|Z=X ,
\end{equation}
up to a fixed rescaling of $X$. For each fixed ideal all the $D'$
are units modulo that ideal for sufficiently large $X$. Thus the
zero-frequency Poisson term on fixed major arcs has exactly the
ordinary complete power sum. This avoids the need for any theorem on
the distribution of smooth algebraic integers in progressions.

\subsection{A loss-free moment and the local-to-global step}
The analytic hypothesis in the next theorem is an estimate on one
trace torus, at full root volume $X$. Keeping the order and the loss
explicit makes clear what is required of a field analogue of the
fibre estimate.
\begin{theorem}\label{thm:nf-absorb-transfer}
Let $F$ be the sum over one fixed choice of the restricted-root
windows just constructed, in a sufficiently small fixed dilation.
Suppose that at order $r>0$ we have, for a fixed $E\ge0$,
\[
                     \int_{\mathbb T}|F|^r\ll X^{r-k+E}.
\]
The exponent $E$ includes the losses for these windows. If $b\ge1$
is an integer with
$b/[k(k-1)]>E$, and $S$ is even, $S>2k$, $S-b\ge r$,
then
\[
  W_{\mathbb K}(k)\le G_{\mathbb K}(k)
       \le \max\{S,H_{\mathbb K}(k)\}.
\]
In particular, if $H_{\mathbb K}(k)\ge S$, the global constant is
one of the two integers $H_{\mathbb K}(k)-1,H_{\mathbb K}(k)$.
The same conclusions hold if $S>2k$ is even and some integer
$1\le b\le w_c=k(k+1)$ instead satisfies
\[
           S-b\ge r(1-b/w_c),\qquad b/w_c>(1-b/w_c)E .
\]
\end{theorem}
\begin{proof}
First cover the restricted roots with a smooth ordinary envelope
$f_0$ whose coefficients dominate theirs. With $A=\int|F|^S$,
orthogonality and positivity, replacing $b$ factors by the envelope
in the even moment, give
$A\le\int|f_0|^b|F|^{S-b}$. On $\mathfrak m_*$, the moment
hypothesis, $|F|\ll X$, and \eqref{eq:nf-point} save a positive
power over $X^{S-k}$. Under the alternative hypotheses, remove
$S-b-r(1-b/w_c)$ restricted factors trivially. H\"older, the
assumed $r$th moment and \eqref{eq:aud-field-critical} give an upper
bound $O_\varepsilon(X^{S-k-b/w_c+(1-b/w_c)E+\varepsilon})$ there.
For $b=w_c$ use the critical estimate directly and bound the remaining
factors trivially. Both alternatives give a strict power saving.
On $\mathfrak M_*$, H\"older and
\eqref{eq:nf-prune} bound the integral by
$C X^{(S-k)b/S}A^{1-b/S}$. Division, or Young's inequality, gives
$A\ll X^{S-k}$. This even moment is loss-free.

For $s=\max(S,H_{\mathbb K}(k))$, use $s-S+b$ ordinary variables
and $S-b$ restricted variables. In the estimates off the core bound
the extra $s-S$ ordinary factors trivially. On $\mathfrak m_*$
we still save a power. On
$\mathfrak M_*\setminus\mathfrak M(U)$ the saving is now a factor
$O(U^{-\tau b/S})$, by the same H\"older inequality with
$\int|F|^S\ll X^{S-k}$. All constants are independent of the fixed
cutoff $U\ge1$.

For clarity we justify the main term and the order of limits.
For fixed $U$ take $X\to\infty$. In fact the cores are disjoint
whenever $U^2=o(X^k)$. If lifts of two distinct centers with
norms $Q_1,Q_2$ meet, their nonzero difference $\Delta$ lies in
$(\mathfrak q_1\mathfrak q_2)^{-1}\mathfrak D^{-1}$. Thus
$|N\Delta|\gg_{\mathbb K}(Q_1Q_2)^{-1}$, while coordinatewise
$\max(1,|s-t|)\le2\max(1,|s|)\max(1,|t|)$ gives
$|N\Delta|\ll_{\mathbb K}X^{-k}U^2/(Q_1Q_2)$.
The same norm separation excludes overlapping distinct lifts of
one center. For every
term \eqref{eq:nf-terminal-form}, Poisson in the free cofactor gives
on each core the ordinary complete sum $g(\eta)$ times an
archimedean integral of size $Z$, with an error $o_U(Z)$ uniform in
the marks. Indeed, the free volume tends to infinity by a fixed
power, and, after normalization by $Z^{1/n}$, the functions have
bounded derivatives depending on $U$ and fixed parameters. The
normalized factors $D'/|ND'|^{1/n}$ vary in a fixed compact
archimedean set of units in $V$. The periodic weight is summed on
finitely many fixed residue lattices. The substitution
$z\mapsto D'z$ permutes their residue classes, giving precisely
$g(\eta)=Q_\eta^{-1}\sum_{x\bmod\mathfrak q_\eta}
e_{\mathbb K}(\eta x^k)$. Summing over the marks preserves the
$o_U(X)$ error since $\sum Z\asymp X$.
The same expansion with $Z=X$ holds for every ordinary factor.

After extending the zero-frequency terms, the core integral is
\[
 X^{s-k}\mathfrak S_s(\nu)\mathfrak J_X(\nu/P^k)
  +O(U^{-\tau'}X^{s-k})+o_U(X^{s-k})
\]
for some $\tau'>0$. The conductor series for
$\mathfrak S_s(\nu)=\sum_\eta g(\eta)^s e_{\mathbb K}(-\eta\nu)$
converges absolutely and uniformly: the complete-sum bound used in
\eqref{eq:nf-complete-bound} applies and $s\ge S>2k$.
The archimedean factor for $F/X$ is a positive weighted average
of fixed smooth integrals with weights $Z/X$. Their total lies
between fixed positive constants. These integrals and the ordinary
ones decay rapidly in every large archimedean coordinate, by
differentiation away from zero before Fourier inversion. Both
truncations are bounded by inserting a sufficiently small power
of $Q_\eta B(P^k\beta)/U$ as in \eqref{eq:nf-prune}.
This proves the extension uniformly in the marks and target.

The construction also permits a power-saving error. Write $c_m>0$ for
a lower exponent of the norms of all marks and $c_z>0$ for a lower
exponent of the terminal volumes $Z$. Choose a fixed
$0<\xi<\min(c_m,c_z,1)/4$ and take $U=X^\xi$. All marks are still
units at every core denominator. On a core write $t=P^k\beta$.
At scale $Z^{1/n}$ the phase of a terminal free variable is a
componentwise $k$th power with archimedean coefficients $t_v$ times
units in a fixed compact set. In Poisson summation, a nonzero dual
frequency $\lambda\in\mathfrak q_\eta^{-1}\mathfrak D^{-1}$ has
\[
 \prod_v\left(\frac{Z^{1/n}|\lambda_v|}{1+|t_v|}\right)^{d_v}
                  \gg\frac{Z}{Q_\eta B(t)}\gg X^{c_z-\xi}.
\]
At least one normalized frequency exceeds a fixed power of $X$.
Its corresponding local integral is rapidly decreasing in this
ratio, by repeated integration by parts once the dual frequency
exceeds the phase gradient. At the other places the local integrals
are bounded, with rapid decay in any additional large normalized
frequency. The product-weight assumption applies also at complex
places, using two real coordinates there. The ideal packing bound
counts $O(1+Q_\eta B(t)/Z)=O(1)$ points in a unit normalized box,
and sums the enlarged dyadic boxes with this rapid decay. Thus the
nonzero modes give $O_A(ZX^{-A})$ for each prescribed $A>0$.
Their trivial complete-sum multipliers are bounded by one.
The same argument works for the ordinary factors with $Z=X$.
There are only polynomially many core centers and all remaining
truncation errors have a power saving in $U$. Thus the power-saving
minor and annular bounds combine with the core approximation to give
a representation formula with error $O(X^{s-k-\omega})$ for some
$\omega>0$ depending on the fixed parameters. Its archimedean factor
is the positive weighted average just described and can vary with $X$.

The singular series is uniformly positive on eligible targets.
Orthogonality at finite levels and the Chinese remainder theorem
express it as the product of normalized local solution densities.
All sufficiently large primes contribute $1+O((N\mathfrak p)^{-1-\delta'})$
by the same absolute majorant. At the remaining primes choose a
primitive $s$-variable solution at each eligible target and put
$e=v_{\mathfrak p}(k)$. Freeze the other coordinates modulo
$\mathfrak p^{2e+1}$. Every choice lifting these residues admits
an exact lift in the unit coordinate by strong Hensel, yielding the
target-independent lower bound
$(N\mathfrak p)^{-(2e+1)(s-1)}$ for that density.

For each central normalized target $t$ on the compact archimedean
cover, choose ordinary weights positive around the equal nonzero
$k$th roots of $t/(s-S+b)$, with positive real roots. The map
from one such ordinary variable to its power has invertible real
derivative at every component. Its density, after fixing the other
ordinary variables in still smaller open boxes, is bounded below
near $t$. The restricted roots can all be taken so small on the
original scale that their total shift lies in this neighborhood,
uniformly over the marks. Fourier inversion, or convolution of the
pushforward measures, now gives a positive uniform lower bound
for $\mathfrak J_X$ throughout a neighborhood of $t$. Use the
finite cover, choose $U$ large relative to the positive lower bound,
then take $X$ sufficiently large. The resulting positive weighted
count gives the desired representation. The closing assertion
follows from $H_{\mathbb K}\le W_{\mathbb K}+1$.
\end{proof}

\section{Restricted moments on number-field fibres}\label{sec:nf-fibres}
We record the number-field version of the estimates on prime factors in
\cite[Sections~3--5]{NFChow}. All exponents below measure volumes.
At fixed $k$ the number of prime windows is fixed before the initial
volume $P_0$ tends to infinity. Individual window sizes, the
smoothness parameter, and all auxiliary numerical tolerances may be
chosen in terms of $k$.

\subsection{Complete sums on a fibre}
At a node of volume $Y$, let $m$ be the product of the chosen generators
of a set of primary primes and assume $|Nm|\le Y^{M/k}$.
The kernel of multiplication by $m^k$ on the trace torus has
$|Nm|^k$ points. Differences of points in one fibre have conductor
$\mathfrak q\mid(m)^k$. Put $U=\log_Y N\mathfrak q$,
$\theta\ge M/k$, and choose windows fine enough that
$\log_Y N\mathfrak p\le\mu$ for all marked primes, where
$\mu>0$ may be arbitrarily small. We omit primes over rational
primes at most $k$ and over the different. The sum to estimate in the
Gram matrix is
\[
 S(\eta)=\sum_{x\in\mathcal O}w(x/Y^{1/n})
            \mathbf1_{((x),(m))=1} e_{\mathbb K}(\eta x^k),
 \qquad \mathfrak q=\mathfrak q_\eta.
\]
Here $w$ is a smooth nonnegative box envelope. Two families of
inequalities suffice for the numerical argument below.

\begin{lemma}\label{lem:nf-fibre-savings}
Fix $\sigma>0$. For each of the following intervals of $U$, one has
$|S(\eta)|\ll Y^{1-\sigma+O_k(\mu)+\varepsilon}$. Empty intervals
may be discarded.
\begin{enumerate}
\item For $2\le j\le k$ put $t=\min(2^{j-1},j(j-1))\sigma$.
The first family of intervals is specified by
\[
 U/(j+1)-j\theta\ge t,\qquad
 1-U/(j+1)-\theta\ge t,\qquad j-U-j\theta\ge t .
\]
\item For $0\le l\le k-2$ with $2^l\sigma\le1/2$, and
$\theta\le1-2^l\sigma$ when $l>0$, the second family consists of
\[
 l\theta+2^{l+1}\sigma\le U\le
      l+2-(2^{l+2}-2)\sigma .
\]
\end{enumerate}
All constants in the implied factor may depend on the number of marked
primes. The saving before the rounding losses is independent of
$[\mathbb K:\mathbb Q]$.
\end{lemma}
\begin{proof}
For the first assertion, put $e_{\mathfrak p}=v_{\mathfrak p}(\mathfrak q)$,
$s_{\mathfrak p}=\max(1,\lceil e_{\mathfrak p}/(j+1)\rceil)$ for
$\mathfrak p\mid(m)$, and choose a balanced generator $a$ of their
product with these exponents. This ideal is principal since the
marked primes are principal. Split into classes $x=x_0+a y$;
only classes of units at every prime dividing $m$ occur. All terms
of degrees exceeding $j$ in the binomial expansion have integral
trace frequency. The degree-$j$ coefficient and its multiples by
$j$ and $j!$ have conductor norm $Q_j$ satisfying
\[
 U/(j+1)-j\theta\le\log_Y Q_j\le U/(j+1).
\]
Indeed their exact conductor is
$\prod\mathfrak p^{(e_{\mathfrak p}-js_{\mathfrak p})_+}$.
The new volume $Y'=Y/|Na|$ satisfies
$\log_Y Y'\ge1-U/(j+1)-\theta$. Relative derivatives of the
weights are uniformly bounded by balancing $a$ and choosing coset
representatives in a balanced fundamental region.

A complete rational leading frequency of conductor norm $Q$ at
volume $Y'$ satisfies the Weyl estimate with normalized right side
\begin{equation}\label{eq:nf-fourterms}
 (Y')^\varepsilon
 \{Q^{-1}+(Y')^{-1}+(Y')^{1-j}+Q(Y')^{-j}\}^{
                   1/\min(2^{j-1},j(j-1))}.
\end{equation}
For the Vinogradov exponent, apply \eqref{eq:nf-weyl-sample}.
Each frequency has multiplicity $O(1+Y'/Q)$ by ideal packing, and
a top-coefficient box contains $O(1+Q/(Y')^{j-1})$ distinct
frequencies: they lie in the lattice
$\mathfrak q^{-1}\mathfrak D^{-1}$ modulo the inverse different.
For the classical exponent apply $j-1$ successive differences and
Cauchy's inequality. The resulting normalized power $2^{j-1}$ is
bounded by $O((Y')^{-1})$ from zero shifts, plus $(Y')^{-j}$
times the sum of the absolute linear sums with frequencies
$j!\theta_0 h_1\cdots h_{j-1}$, where $\theta_0$ is the leading
coefficient and $|h_{i,v}|\ll(Y')^{1/n}$. Each nonzero product has
$O((Y')^\varepsilon)$ factorizations in these boxes: ideal
factorizations obey the divisor bound, and generators differ by units
of logarithmic height $O(\log Y')$. This divisor bound follows
from the integer divisor bound applied to the rational norm and the
fixed degree. The product lies in a box of volume $O((Y')^{j-1})$.
Poisson in the linear variable gives $O(Y')$ times a rapidly
decreasing sum of frequency indicators at radii $(Y')^{-1/n}$.
Individual frequencies have multiplicity
$O(1+(Y')^{j-1}/Q)$, and the number in a frequency box is
$O(1+Q/Y')$. Enlarged boxes are summed with their rapidly decaying
weights. This proves \eqref{eq:nf-fourterms}. Insert the bounds
for $Q_j,Y'$ and sum over the $|Na|$ classes to obtain part~1.

For part~2 begin with completion of unit sums. At any surviving
conductor prime, a primitive phase whose numerator is, modulo that
prime, a nonzero multiple of $x^v$ up to a constant ($2\le v\le k$)
has complete sum $O_k(q^{e/2})$ over units, uniformly after a linear
twist. At prime level this is Weil's bound. At $e=2u$ write
$x=x_0+\pi^u h$ and sum in $h$; only unit roots of the derivative
modulo $\pi^u$ remain. There are at most $v-1$ because the second
derivative is nonzero at units modulo $\pi$. At $e=2u+1$, $u\ge1$,
summing first modulo $\pi^{u+1}$ leaves one extra quadratic Gauss
sum of size $q^{1/2}$ per such class. This applies with every
primitive trace character, since coefficients can be rescaled by a
unit. The Chinese remainder theorem and lattice Poisson completion
give $O(Y/\sqrt Q+\sqrt Q)$, $Q=N\mathfrak q$. For a prime formerly
in the conductor or merely in the original mask, the unit restrictions
can exclude finitely many translated residues. Expand them by
inclusion--exclusion. Each term is a coset modulo an ideal
coprime to the active conductor, a product of marked principal primes.
After balancing its generator, the new volume is at most a fixed
multiple of $Y$, with the same normalized derivative bounds. Its
affine substitution permutes the active residue classes. If this
volume is below one, its integral box contains $O(1)$ points by
ideal packing, and $1\le\sqrt Q$. Thus the same completion bound holds also for the correlation weights
occurring below. The initial normalized saving is
$g_0(U)=\max(0,\min(U/2,1-U/2))$.

For the next correlation choose $\theta\le u\le\min(1,U)$.
Starting with the radical of the active conductor, add prime factors
one at a time to choose a subideal divisor with log-norm
$u+O(\mu)$ and a balanced generator $b$. Average over shifts $bh$
for $h$ in a box of volume $H=Y/|Nb|$. After Cauchy's inequality,
the normalized square is bounded by $O(H^{-1})$ plus the average
of the normalized absolute correlations at nonzero differences
$h$ in an enlarged box of volume $O(H)$. This follows by extending
the summand by zero: the translated supports lie in one box of
volume $O(Y)$ and each difference occurs $O(H)$ times. Bounded
$H$ is allowed using the trivial estimate. The difference quotient
$[f(x+bh)-f(x)]/(bh)$ reduces to the formal derivative at each
active prime, since $b$ contains its radical. The remaining phase
has the stated monomial type with degree one less. With
$\mathfrak q_0$ the conductor after cancelling $b$, cancelling $h$
loses $N((h),\mathfrak q_0)$ in norm. Its average is $O(1)$:
a divisor ideal $\mathfrak d\mid\mathfrak q_0$ dividing
$(h)\ne0$ has $N\mathfrak d\ll H$, and packing bounds the
number of multiples by $O(H/N\mathfrak d)$. The number of divisors
of $\mathfrak q_0$ is bounded in terms of $k$ and the number of
marks. The available savings remain 1-Lipschitz in $U$ under
\[
 g\longmapsto\max\left(g(U),
 \tfrac12\sup_{\theta\le u\le\min(1,U)}
               \min(1-u,g(U-u))\right)
\]
with the second term zero for $U<\theta$. To check the assertion
clip $u$ to $\min(u,V)$ upon decreasing $U$ to $V\ge\theta$;
at $U=\theta$ the new branch vanishes because $g(0)=0$.
Thus the conductor loss is paid for by the preceding average and
rounding costs $O_k(\mu)$. On surviving primes masks are unchanged.
Only at inactive primes can translated excluded classes multiply,
and their number remains bounded after a fixed number of differences,
as required for completion. For a saving $\sigma$ after $l$
differences the successive shifts can be chosen with exponents
$u_i\in[\theta,1-2^i\sigma]$, $1\le i\le l$, and the final
conductor exponent in $[2^{l+1}\sigma,2-2^{l+1}\sigma]$.
The set of their possible total initial exponents is exactly the
interval in part~2. The Lipschitz argument justifies these choices
also for individual shifts where additional primes cancel.
\end{proof}

\subsection{Shortening the moment}
Choose the primary prime windows with disjoint interiors and
reciprocal norm sums bounded by fixed constants. A state retains a
subset $I$ of the windows, the complementary primes having been fixed
and removed; its root volume is $Y$. Let $F_I$ denote its exponential
sum. Choosing one prime generator in each of some $T\subseteq I$
gives products $m$ and the decomposition
$F_I(\alpha)=\sum_m F_{I\setminus T,m}(m^k\alpha)$.
Each child has volume $Y/|Nm|$, and the sum of $1/|Nm|$ is $O(1)$.
Take $r\ge k$. We record an admissible excess $D$ at order $r$ separately from
a common loss on the initial scale. Thus the estimate at each state
is $Y^{r-k+D}P_0^{c_0}$, where $c_0>0$ is a tolerance prescribed
before the windows. All states and all orders of a given finite
calculation must use the same windows. For the fixed finite menu choose $C\ge4$ such that
$2/C<1-M/k$ for every scheduled $M<k$. The bottom cutoff
$Y\le P_0^{C\gamma}$ then costs $Ck\gamma$ in the exponent.

Choose $(z,M)$ such that all $1<U\le M$ have a common strict
saving \(\sigma>2z\) on the weighted Gram entries, using the
conductor estimates in this section with their radical losses. For $M=1$
there is no coverage condition. Suppose $0<L<M$, and the removed
block can be chosen so that
\[
            Y^{L/k}\le |Nm|\le Y^{M/k}.
\]
Then, up to the common initial-scale loss,
\begin{equation}\label{eq:nf-short}
       D'\ =\ \max\{D-2z,(1-L/k)D\}
\end{equation}
is admissible at order $r+2$. We prove this for real moment orders.

For a fixed $m$, let $H_m$ be the part of $F_I$ supported on units
at its marked primes. On a fibre of multiplication by $m^k$, the
number of points where $|H_m|\ge Y/W$ is $O(W^2)$ for
$W\le C Y^z$ with any fixed $C$. In fact, take the rows at these
points in the Gram matrix of the vectors
$(w(x/Y^{1/n})^{1/2}e_{\mathbb K}(\alpha x^k))_{((x),(m))=1}$,
choosing the envelope $w$ to dominate the coefficient support.
The squared norm of the coefficient vector is $O(Y)$. For
differences of conductor norm at most $Y$, the corresponding
operator norm is $O(Y)$. To see this, impose each exact conductor
by M\"obius inversion among divisors of $(m)^k$. For an annihilating
ideal $\mathfrak d$ of norm at most $Y$, partition into cosets of
the subgroup killed by $\mathfrak d$. On each full coset Fourier
diagonalization gives eigenvalues $N\mathfrak d$ times weighted
counts of the unit roots of $x^k\equiv c\pmod{\mathfrak d}$.
There are at most $k^{\omega(\mathfrak d)}$ residue classes by
Hensel lifting and at most $O(1+Y/N\mathfrak d)$ roots per
class. Restricting to the sampled rows is compression of a positive
matrix. Only boundedly many divisors occur. Higher conductors have
entries \(O(Y^{1-\sigma+O_k(\mu)+\epsilon})\) by the assumed
conductor saving. With $J$ rows, the Gram inequality
therefore gives
\[
 JY^2/W^2\ll Y(Y+JY^{1-\sigma+O_k(\mu)+\epsilon}).
\]
Take $\mu,\epsilon$ small compared with the strict margin. This
proves the fibre assertion, including $W$ in a fixed constant
multiple of the indicated range.

The part where $|F_I|\le Y^{1-z}$ contributes the first alternative
in \eqref{eq:nf-short}. On the complementary set, suppose first that
$|H_m|\ge |F_I|/2$ for every candidate $m$. On each fibre the sum
of $|F_I|^2$ over these points is $O(Y^2\log(2Y))$, by dyadic
decomposition and the fibre assertion. Weighted H\"older applied to
the decomposition into children gives
\[
 |F_I|^r\ll\sum_m |Nm|^{r-1}|F_{I\setminus T,m}(m^k\alpha)|^r .
\]
Integrate after multiplication by $|F_I|^2$. Haar measure gives
weight $|Nm|^{-k}$ per preimage. The child's $r$th moment then
leaves a factor $1/|Nm|$ times
$Y^{r+2-k}(Y/|Nm|)^D$, apart from the indicated common loss.
This is bounded by the second alternative.

Otherwise $|F_I-H_m|>|F_I|/2$ for some $m$. Inclusion--exclusion
gives $|F_I|\ll |F_{I\setminus T',m'}((m')^k\alpha)|$ for a
nonempty partial product $m'\mid m$. This is an exact identity for
each constituent sum, since the terminal cofactor contains no primary
mark and a prime in a specified window uniquely fixes that primary
choice. The implied constant depends only on the number of windows.
Induct at the same moment order on the number of remaining marks.
In summing over possible partial products the net factor is bounded
by $\sum_{\mathfrak a}(N\mathfrak a)^{k-r-2}$, which converges since
$r+2-k\ge2$.
Here the preimage formula over the full torus is Haar invariance, and
the count of distinct products per ideal is bounded. At a bottom
state $|F_I|\ll Y$ yields $Y^{r+2}\le
Y^{r+2-k+D'}P_0^{Ck\gamma}$. The induction uses the common bottom
loss once: at each branch only one lower-order integral, or one
smaller-state integral, occurs per summand. The maximum in
\eqref{eq:nf-short} carries this loss forward. After fixing
$\gamma$, the windows and all moment orders, the finitely many
logarithmic factors are $P_0^{o(1)}$. Choose the Weyl and moment
epsilon losses small enough, over the finite list of steps, to fit
the unspent part of $c_0$.

\subsection{Initial moments and the numerical coefficient}
The rough smoothness range in Theorem~\ref{thm:edmoment} is sufficient
for the initial moments at every state. Fix $k$ and choose $\gamma>0$
small. Choose the primary log-norms in disjoint subintervals of
$(0.80\eta,0.90\eta)$ measured in base $P_0$, with their sum
in $(1-2\gamma,1-1.99\gamma)$. For sufficiently small fixed
$\eta$ this is possible by choosing the number of primary primes
near $(1-2\gamma)/(0.85\eta)$ and then narrow intervals near
distinct nearby centers. At terminal volume $Y_0$ one has
$1.99\gamma<\log_{P_0}Y_0<2\gamma$. Choose the auxiliary count
near $2\gamma/(0.62\eta)$ and the exponents $a,b$ of
Section~\ref{sec:nf-global} close enough to the reciprocal of
this count that their prime log-norms belong to
$(0.55\eta,0.70\eta)$ and the free cofactor has volume $O(P_0^\eta)$.
The exponent of that volume stays strictly positive. For a chosen
free cofactor $z$, multiplication by $z^k$ is Haar preserving on
the trace torus, and each remaining root has all its prime factors
in $(P_0^{\eta/2},P_0^\eta]$. Its height is $O(Y^{1/n})$ on a
state of volume $Y$, by the balancing and the annular support for
$z$. Each root has coefficient at most one by uniqueness of the
factorization. The uniform scale form in the proof of
Theorem~\ref{thm:edmoment}, with $W=P_0$ there, and even-moment
orthogonality thus give the same excess $\delta_s$ at order $2s$.
The triangle inequality over $O(P_0^\eta)$ free cofactors and the
error in efficient differencing give an initial loss
$P_0^{O_{k,s}(\eta)+o(1)}$, also at partially removed states.
This loss is introduced only at the initial moment in the induction.
For $Y>P_0^{C\gamma}$ the remaining primary exponents sum to at
least $\log_{P_0}Y-2\gamma$, which exceeds
$(M/k)\log_{P_0}Y$. Thus a block in any prescribed range
$[L/k,M/k]$ relative to $Y$ exists whenever $L<M<k$. Indeed,
add consecutive remaining marks until reaching the lower endpoint.
Their maximal relative log-norm is $O(\eta/\gamma)$; the primary
interval widths may be made as small as necessary. On the finite
list of steps at fixed $k$, choose $\eta$ small enough to
respect the smallest gap and the required conductor tolerance.
More precisely, first reserve a common strictly positive margin in
each final ordinary or marked compensation test. Choose $c_0$ inside
all those margins, allowing $2c_0$ in place of $c_0$ when the partial
volume is $Y\ge P_0^{1/2}$. Choose $C$ for the finite list, then
$\gamma$, then $\eta$ and the window widths, so that
$Ck\gamma+O_{k,r}(\eta)<c_0$ with room for the logarithmic and
epsilon losses. At volume $P_0$ the exponent in
Theorem~\ref{thm:nf-absorb-transfer} may now be taken as $D+c_0$.
At the indicated partial volumes one may take $D+2c_0$.
No subsequent compensation requires a new family of roots.

For specificity we describe the scalar integration used for the
limiting coefficient. Set
\[
 z_0=0.09985,\qquad d_0=z_0/(1-z_0),\qquad
 v_0=1-d_0-\log d_0,\qquad \lambda=0.99999 .
\]
At every scheduled initial height $1/20\le z_0\le1/3$ the field
recurrence gives the rounded order
\[
          r_0=2\lceil kv_0/2\rceil,\qquad \delta_{r_0/2}\le kd_0
                    \quad(k\ge21).
\]
Indeed, for $\Phi(x)=k\log x+x$, \eqref{eq:eddefect} gives
$\Phi(\delta_s)\le\Phi(k-1)-(2s-2)(1-2^{-k})$.
One has $\log(1-t)\le-t-t^2/4$ for $0<t<1$: putting
$u=\sqrt{1-t}$ gives $2\log u\le2(u-1)=-t-(1-u)^2$ and
$1-u\ge t/2$. Hence $\Phi(k-1)\le\Phi(k)-2-1/(4k)$.
Since $d_0\ge1/20$, the inequality $\log(1/d_0)\le1/d_0-1$
gives $v_0\le20$ and $r_0\le21k$. For $k\ge21$ one has
$84k^2\le2^k$, so $(r_0-2)2^{-k}\le1/(4k)$.
Thus $\Phi(\delta_{r_0/2})\le\Phi(k)-r_0\le\Phi(kd_0)$.
Monotonicity of $\Phi$ proves the assertion.
For an admissible menu of pairs $(z,M)$ put
\[
 V(d)=\max_{(z,M)}\min\{z,\lambda Md/2\}.
\]
At fixed $k$ the menu used down to a positive stopping threshold
may be made finite by arbitrarily fine approximation with strict
inequalities. Applied with $L=\lambda M$, \eqref{eq:nf-short}
replaces $d=D/k$ by $d-2V(d)/k$ at the cost of two moment orders.
Since $V$ is increasing, the normalized increase in order down to
the first $d\le d_*$ is at most
\[
 2/k+\int_{d_*}^{d_0}\frac{dt}{V(t)}.
\]
Take $B_k=\lceil16\sqrt{k}\rceil$ ordinary variables, and stop
at $d_*=B_k/(2k^3)$. Then $kd_* < B_k/[k(k-1)]$ with slack.
Choose $\gamma$ and the remaining errors small relative to this
slack. The summand exponent and its numerical thresholds are uniform in
the field. The implied constants and the lower bound on the target
size may depend on the field and on $k$.

Here and below Appendix~\ref{app:certificates} gives the
calculations producing the menus. At each height round
$\sigma=(200001/100000)z$ upwards on the working integer grid.
Use the two conductor intervals of Lemma~\ref{lem:nf-fibre-savings},
rounding lower endpoints up and upper endpoints down, and keep their
connected union from $1$ up to the cutoff. With $F_0=120$ this
elementary calculation starts at $z_0=0.09985$, rounds each new
height $1000z/1001$ down in units of $10^{-15}$ and first reaches
$M=120$ on row $9863$. The finite calculation uses $\theta=0$.
For use at positive $\theta=H/k$, first reduce every scheduled
$\sigma$ to $2.000005z>2z$. Every inequality containing $\theta$
in the finite covering then has positive slack at $\theta=0$.
There are finitely many selected intervals, so the same covering
holds simultaneously for all sufficiently large $k$.

For $120\le U\le M\le H=\sqrt{k}/8$, put $j=\lceil3U/2\rceil$,
$h=3/2$, $b=h+1/120$ and $\tau=1/14$. We have $j\le bU$ and
$j(j-1)\le hbU^2$. With the common $\theta=H/k$, the three brackets
in the degree-reduction test are bounded below by
\[
 (h+2/120)^{-1}-b/64,\qquad
 1-1/h-H/k,\qquad (h-1)120-b/64.
\]
Each exceeds $2\tau hb$ for all sufficiently large $k$. Thus
$z=\tau/M^2$ is feasible with a strict margin. The last mesh height
is at least $\tau/120^2$, so the covering also holds on $1\le U<120$
at every tail height with this same $\theta$. For each fixed $k$,
compactness and the strict margins permit a finite subdivision of
the continuous tail with arbitrarily small additional price.
The linear part $V(d)=\lambda H d/2$ below its last transition
contributes $O((\log k)/H)$ down to $d_*$.

For clarity, on a constant-menu interval of height from $z'$ to
$z$ at modulus $M$, the charge used in the sum is
$2\log(z/z')/(\lambda M)$. This follows by splitting the integral
at $d=2z/(\lambda M)$; the constant stretches of $V$ introduce
telescoping reciprocal-modulus terms. At the initial height
$\lambda M>2(1-z_0)$, and the initial seed plus its constant stretch
is $v_0+d_0/z_0=2+\log((1-z_0)/z_0)$ before these terms. On the
continuous tail $z=1/(14M^2)$ to infinity the additional integral
of $2\,dz/(\lambda Mz)$ is $4/(120\lambda)$; extending past the
finite cap costs only the preceding $o(1)$. The elementary sum, including this tail and the bridge from the last
mesh height to $1/(14\cdot120^2)$, is less than
$6.917355607003731$ with outward rounding. In particular an upper
bound is
\[
 C_* =\frac{13836517986681517}{2000000000000000}<6.919 .
\]
With
rounding to an even moment and the $B_k$ ordinary variables we obtain
from Theorem~\ref{thm:nf-absorb-transfer} that
\begin{equation}\label{eq:nf-final-linear}
       G_{\mathbb K}(k)\le
       \max\{H_{\mathbb K}(k),(C_*+o(1))k\}.
\end{equation}
The error in the summand count is uniform over fields; the
representation threshold may depend on both field and exponent.

\subsection{Arithmetic exponent pairs on the trace torus}
\label{sec:nf-arith-pairs}

The coefficient in \eqref{eq:nf-final-linear} decreases on applying
finite Fourier transformation to the prime-power phases. Arithmetic
exponent pairs over $\mathbb Q_p$ have been developed by
Mili\'cevi\'c \cite{NFMili}. We need simultaneous transformations at the
prime ideals in a marked product and keep track of the cost of
restricting a sum to a residue class of its radical. The chosen
prime ideals have residue degree one and are unramified; their local
fields are therefore $\mathbb Q_p$. The dual lattice in Minkowski
space still depends on the full field.

Let $\mathfrak m$ be a squarefree product of at most $R$ marked
principal prime ideals, $m=N\mathfrak m$, and let
$\mathfrak q\mid\mathfrak m^k$, $Q=N\mathfrak q$. On a fixed
residue class $c$ modulo $\mathfrak m$, take a phase consisting
of an arbitrary linear archimedean character times a product of
primitive finite characters applied to integral analytic functions
$\Phi_{\mathfrak p}$ at the conductor primes $\mathfrak p$.
The summation weight is smooth on any translated archimedean box
of volume scale $X$ and bounded shape, with uniformly bounded
derivatives after rescaling. We state estimates uniform in the
translation and in the linear character.

Let a word in $A,B$ and a rational $v$ be fixed. Reading its letters
from right to left, replace $v$ by $v-1$ at $A$ and by $v/(v-1)$ at
$B$. Each input and output must be finite and avoid $0,1$.
For each numerator--denominator pair $(a,b)$ along this reversed
word, exclude the rational primes dividing $ab(a-b)$, as well as
the primes at most the word length plus four. This is possible
since all marked norms tend to infinity for fixed $\mathbb K,k$.
At each active prime demand that the expansion of
$\Phi_{\mathfrak p}(c+t)$ on the residue disc have integral
coefficients, agreeing modulo $\mathfrak p$, through degree
$N_0=3+\text{length of the word}$, with
\[
               C+Dt+A_0(1+B_0t)^v,\qquad A_0,B_0\text{ units}.
\]
A new lift $c$ in the same residue disc has the same reduced Taylor
coefficients. For conductor exponent one one may replace the local
function by any such germ with the same value modulo $\mathfrak p$.
The original monomial has the required form on every unit class,
with $v=k$.

\begin{proposition}\label{prop:nf-dual}
Start with $(\kappa,\lambda,D)=(1/2,1/2,1/2)$ and read the word
from left to right, applying
\begin{align*}
 A(\kappa,\lambda,D)
 &=\left(\frac{\kappa}{2+2\kappa},
          \frac{\kappa+\lambda+1}{2+2\kappa},
          \max\{1/2,(D+1)/2\}\right),\\
 B(\kappa,\lambda,D)
 &=\left(\lambda-1/2,\kappa+1/2,D+1+\lambda-1/2\right).
\end{align*}
For $Q\ge X\ge1$ the sum over a single residue class has size
\[
     \ll (2mQX)^\varepsilon m^D X^{\lambda-\kappa}Q^\kappa.
\]
The constants may depend on the field, $k,R$, the word, and
finitely many scaled derivatives of the weight.
\end{proposition}

\begin{proof}
Write $\mathfrak Q=\operatorname{lcm}(\mathfrak q,\mathfrak m)$ and
$Q'=N\mathfrak Q\le mQ$. On a fixed class modulo a marked prime,
the second derivative of the local phase is a unit. At exponent
$e\ge2$ and Fourier frequency $b$, there is a stationary point on
this class precisely when $\Phi'(c)+b\equiv0\pmod{\mathfrak p}$.
It is then unique in the residue disc by Hensel's lemma. Write it
as $x_b$. With primitive character $\psi$ at exponent $e$, the
finite transform on the class has the form
\begin{equation}\label{eq:nf-legendre}
 (N\mathfrak p)^{e/2}C_b\,\psi(\Phi(x_b)+b x_b),
\end{equation}
where $|C_b|=1$ and $C_b$ is constant as $b$ varies in its residue
class modulo $\mathfrak p$. For $e=2u$, summing first over
increments divisible by $\mathfrak p^u$ leaves precisely the class
of $x_b$ modulo $\mathfrak p^u$. At $e=2u+1$, $u\ge1$, the same
argument leaves a quadratic Gauss sum at residue level. Its
coefficient is $\Phi''(x_b)/2$ times a fixed unit, hence is constant
modulo $\mathfrak p$. At a mask prime with original conductor
exponent zero or one, its exponent in $\mathfrak Q$ is one.
The transform on the chosen class $c$ is its fixed phase value
times the linear character in $b$ obtained by evaluation at $c$.
It has absolute value one for every $b$. Ordinary Poisson summation
and ideal packing now give
the completion bound $O(X/\sqrt {Q'}+\sqrt {Q'})$, uniformly in
the real linear twist. For $Q\ge X$ this implies the initial triple.

For $A$, put $U=\log_XQ$ with $X>1$. The normalized saving available
before the step, with cost $m^D$, is at least $g(U)=U/2$ on $[0,1]$,
and
\[
      g(U)=\max\{0,1-U/2,1-\lambda-\kappa(U-1)\}
                         \quad(U\ge1).
\]
Here completion loses at most $m^{1/2}$ on replacing $Q'$ by $Q$.
The function is $1$-Lipschitz. For $U\ge1+\lambda$ put
$t=(\lambda+\kappa(U-1))/(1+\kappa)$. If $t>1$ the proposed saving
is nonpositive. Otherwise $U-t\ge1$ and $g(U-t)\ge1-t$.
Form an ideal $\mathfrak b$ with exponent
$\max(1,\lceil t v_{\mathfrak p}(\mathfrak q)/U\rceil)$
at each $\mathfrak p\mid\mathfrak m$. Its norm and the norm
of its gcd with $\mathfrak q$ lie between $X^t$ and $mX^t$.
Choose a balanced generator $s$. If $N\mathfrak b>X$, the factor
$m^{1/2}$ absorbs the saving $(1-t)/2$. Otherwise apply the
correlation bound of Lemma~\ref{lem:nf-fibre-savings} with shifts
$sh$, for $h$ in a box of volume $H=X/N\mathfrak b$.
The linear archimedean character contributes a constant, while the
divided finite difference of the local phase has, modulo an active
prime, the Taylor coefficients of its derivative through one degree
less. Multiplication by a local unit in choosing the residual
primitive character preserves the required form. The exponent of
the model germ is now $v-1$. At $h\ne0$ the loss in the logarithmic
conductor relative to $U-t$ is at most
$\log_X(m N((h),\mathfrak q/(\mathfrak q,\mathfrak b)))$.
The average gcd norm is $O_{k,R,\mathbb K}(1)$ by the divisor
count in Lemma~\ref{lem:nf-fibre-savings}. The normalized squared sum
is consequently bounded, up to $\varepsilon$-powers, by
$mX^{t-1}+m^{D+1}X^{-(1-t)}$. Its square root gives the displayed
$A$ transformation. In the remaining interval $1\le U\le1+\lambda$,
completion gives at least the required saving: the two straight lines
meet at $1+\lambda$, and the completion line has the more negative
slope. Bounded $X$ is covered by the trivial estimate.

For $B$, the phase in \eqref{eq:nf-legendre} is of the same analytic
type with $v'=v/(v-1)$. More explicitly, the implicit-function power
series of $x_b$ has integral coefficients about the chosen center,
since the linear coefficient of $\Phi'$ is a unit. It converges
on the residue disc. Reducing finite series modulo the maximal ideal
gives, up to an affine change with unit slope, the formal inverse
of $(1+B_0t)^{v-1}$. The derivative of the critical value
$\Phi(x_b)+bx_b$ is $x_b$. Integration of this jet gives a constant
plus a linear term plus a unit multiple of a binomial of exponent
$v'$ and unit slope. All factorials used here are units by the
choice of marked primes. Taking one fewer Taylor coefficient
therefore suffices. Signs and the unit factors in the Chinese
remainder isomorphism only rescale these local coordinates.
At the primes of exponent one, fix $b$ modulo the prime and choose
a primitive local character. Use a new integral germ with the required
exponent, taking its constant value to be zero if the entire original
transform on this class is assigned to the constant factor. There are at most $m$ dual classes altogether, and on each
one the product of the constants in \eqref{eq:nf-legendre} is fixed.

To pass from the finite to the archimedean transform, choose a
balanced generator $q_0$ of $\mathfrak Q$, and a fixed nonzero
rational integer $d$ with $I=d\mathfrak D^{-1}\subseteq\mathcal O$.
Ignore the finitely many prime divisors of $dI$. Poisson summation
on the residue classes gives frequency lattice
$\mathfrak Q^{-1}\mathfrak D^{-1}$; write its frequencies
$t/(d q_0)$ with $t\in I$. At each conductor prime this is a unit
rescaling of the frequency parameter of the local transform, after
the specified power of its uniformizer is extracted. The indicator
of $I$ in $\mathcal O$ is a finite Fourier series on the fixed
group $\mathcal O/I$. Each term changes only the archimedean linear
twist of the dual sum.

After removing the linear factor from the translation of the initial
box, the dual archimedean weight is a translate of a Schwartz function
at volume scale $Y=Q'/X\ge1$. Its shape is uniformly bounded because
$|q_0|_v\asymp_{\mathbb K}(Q')^{1/n}$ at every embedding. To apply
the induction hypothesis, take a smooth partition of unity on
$\mathbb R^n$ into translates of a bounded cube and rescale by
$Y^{1/n}$. On the cube of index $j\in\mathbb Z^n$ the required
scaled derivatives are $O_N((1+|j|)^{-N})$, measured relative to
the center of the translated Fourier weight. These bounds are
summable over the translates, and are uniform in the original
linear modulation. Thus the dual sums lie in the same class with
$X$ replaced by $Y$ and conductor norm $Q'\ge Y$. Poisson supplies
the factor $X/Q'$ and the product of the local transforms has
magnitude at most $\sqrt{Q'}$. The bound, summed over all classes, is
\[
 \ll (2mQX)^\varepsilon m^{D+1}
        X(Q')^{-1/2}Y^{\lambda-\kappa}(Q')^\kappa
 \le (2mQX)^\varepsilon
       m^{D+1+\lambda-1/2}X^{1+\kappa-\lambda}Q^{\lambda-1/2}.
\]
Smooth weights produced by a correlation have uniformly bounded
scaled derivatives. Choosing the initial derivative order sufficiently
large for the finite word closes the induction.
\end{proof}

The radical loss in the $B$ transformation can in fact be retained
unchanged. This gives a uniform loss bound for words of any length.
\begin{corollary}\label{cor:nf-dual-radical}
Proposition~\ref{prop:nf-dual} holds with
$B(\kappa,\lambda,D)=(\lambda-1/2,\kappa+1/2,D)$.
In particular one may use $D\le1$ for every word starting
with $(1/2,1/2,1/2)$.
\end{corollary}
\begin{proof}
Let $q_*=\prod_{e_{\mathfrak p}\ge2}(N\mathfrak p)^{e_{\mathfrak p}}$,
and let $m_*$ be the norm of its radical. At a conductor-one prime
the initial disk is a single residue $c$. Its transform is the
fixed phase at $c$ times a linear additive character in the dual
variable. The same holds at an inactive prime present in $\mathfrak m$.
Their product therefore changes only the linear twist in the
Poisson sum. Indeed an additive character on a quotient of the
fixed lattice $I$ lifts to a linear character of its real span;
equivalently use the trace-dual lattice. At the primes of
exponent at least two, \eqref{eq:nf-legendre} instead leaves exactly
one joint residue class modulo the radical. The dual conductor,
radical, and amplitude are now $q_*,m_*$ and $\sqrt{q_*}$.
No summation over residues of the radical remains.

At $Y=Q'/X$, if $q_*\ge Y$ the inductive pair bounds the dual sum
by $m_*^D Y^{\lambda-\kappa}q_*^\kappa$ up to epsilon powers.
If $q_*<Y$, use completion on its single class, with the trivial
constant-phase bound when $q_*=1$. Since $Q'/m\le q_*\le Q\le Q'$
and $Q'\ge Y$, this gives
\[
 Y/\sqrt{q_*}+\sqrt{q_*}\ll m^{1/2}\sqrt Y
       \ll m^D Y^{\lambda-\kappa}(Q')^\kappa .
\]
Here $\lambda\ge1/2$, $D\ge1/2$. The first case also satisfies this
last bound. Poisson's factor is $X\sqrt{q_*}/Q'$ and
$\sqrt{q_*}(Q')^{\lambda-1}\le Q^{\lambda-1/2}$.
The claimed estimate follows. The $A$ transformation preserves
$1/2\le D\le1$, completing the proof.
\end{proof}

\begin{theorem}\label{thm:nf-dual-linear}
Put
\[
          C_{\mathrm{base}}=\frac{34360293397827451}{5000000000000000}.
\]
Then for every number field,
\begin{equation}\label{eq:nf-dual-final}
 G_{\mathbb K}(k)\le
       \max\{H_{\mathbb K}(k),(C_{\mathrm{base}}+o(1))k\}.
\end{equation}
The error in the summand count is uniform over fields.
\end{theorem}
\begin{proof}
Split the unit mask of Lemma~\ref{lem:nf-fibre-savings} into at most
$m$ classes. If $Q=X^U$ and $m\le X^{M/k}$, Proposition~\ref{prop:nf-dual}
gives the additional normalized saving
\[
     1-\lambda-\kappa(U-1)-(1+D)M/k\qquad(1\le U\le M).
\]
For every finite menu, the radical loss tends to zero uniformly
on bounded $M$ as $k\to\infty$. The transformed $v$ is always
a fractional linear function of $k$ with nonzero determinant, so
the finitely many words in the menu satisfy the nondegeneracy
condition at every sufficiently large degree. For each remaining
word with $\kappa>0$, Corollary~\ref{cor:nf-dual-radical} adds the
interval
\[
 1\le U\le1+(1-\lambda-\sigma-2\theta)/\kappa.
\]
Deduplicate the numerical pairs through fourteen operations from
$(1/2,1/2)$, retaining one word per pair and testing its exponents
in reverse order. There are $987$ retained pairs with $\kappa>0$.
The interval sweep of Appendix~\ref{app:certificates} now starts
at $z_0=107/1000$, using the same downward height rounding and
upward saving rounding as before. On the $10^{-15}$ grid take the
largest covered endpoint $M\le300$ in the connected component of
$1$. Use the common radical exponent $\theta=10^{-9}$ at every
row. The sweep reaches $300$ in $11776$ rows. The initial $M$
exceeds $1.78847$, so its retained modulus $0.99999M$ exceeds
$2(1-z_0)$. The final height is $8.28173242\cdot10^{-7}$, above
$1/(14\cdot300^2)$. These checks connect both ends of the mesh.
For $H=\sqrt k/8$ and sufficiently large $k$ we have $H/k\le10^{-9}$.
Thus the whole finite covering remains valid at the cap. On
$300\le U\le M\le H$ use the preceding continuous proof with
$120$ replaced by $300$ and $z=1/(14M^2)$.

With retention $99999/100000$, the initial charge, all finite
logarithmic charges, the bridge at $300$ and the infinite tail
have total less than $6.872058370053856<C_{\mathrm{base}}$ at this positive
$\theta$. At $\theta=0$ the corresponding total is less than
$6.872058319557577$. The outward allowance is $10^{-48}$ on these
sums of fewer than $12000$ bounded terms, evaluated with 60-digit
logarithms. Capping the tail and using ordinary compensation costs
$O((\log k)/\sqrt{k})$. This proves \eqref{eq:nf-dual-final}.
\end{proof}
For each fixed field, Lemma~\ref{lem:precision-widths} gives
$H_{\mathbb K}(k)=W_{\mathbb K}(k)$ for all sufficiently large $k$.
The same replacement is valid for a varying field whenever the
precision condition in that lemma holds at its wild places.

\begin{corollary}\label{cor:nf-explicit-large-linear}
For every $k\ge10^{12}$ the even moment order $S$ in
Theorem~\ref{thm:nf-absorb-transfer} can be taken with
$S+1<6.872559k$. In particular
$G_{\mathbb K}(k)\le\max\{H_{\mathbb K}(k),\lfloor6.872559k\rfloor\}$ for every
number field in this range.
\end{corollary}
\begin{proof}
Put $K_0=10^{12}$, $R=10^6$, $\lambda=99999/100000$,
$\tau=1/14$, and $H=\sqrt{k}/8$. At the finite heights up to
$M=300$ and on the part $1\le U\le300$ of every subsequent
covering, use the common upper bound
$\theta_0=1/(8R)\ge H/k$. Corollary~\ref{cor:nf-dual-radical}
costs at most $2\theta_0$ in each pair saving, and the finite
degree-reduction tests use this same $\theta_0$.
For $300\le U\le M\le H$ use the actual upper bound $\theta=H/k$,
so that $M\theta\le H^2/k=1/64$ in the three continuous tests. The accompanying exact-rational
calculation uses the initial height $107/1000$, the same $11776$
levels as above up to $M=300$, and then the saving $z=\tau/M^2$.
All depth-fourteen words used on this range are nondegenerate:
every numerator and denominator of the formal exponents, and their
difference, is affine in $k$ with integer coefficients of absolute
value at most $2^{15}$, and the three forms are nonzero as forms
by unimodularity. Thus the continuous-tail inequalities in
Section~\ref{sec:nf-fibres} hold with $120$ replaced by $300$
and $k\ge K_0$. The majorant cost, if that tail is extended
to infinity, is less than $6.872064631626498$ by
Appendix~\ref{app:certificates}.

Take $b=\lceil16\sqrt k\rceil$ and stop at $d_*=b/(2k^3)$.
At the cap the final linear part of the speed function costs at
most
\[
 \frac{2}{\lambda H}
      \log\frac{2\tau}{\lambda H^3d_*}
 \le \frac{16}{\lambda\sqrt k}\log\frac{128\tau k}{\lambda}.
\]
If the cap is already sufficient there is no extra cost.
The upper bound decreases for $k\ge K_0$, and at $K_0$ is
less than $0.000477508839893$. The initial even order
$2\lceil kv_0/2\rceil$ exceeds $kv_0$ by less than two. The discrete defect iteration
adds at most two beyond its integral cost. The ordinary allocation
contributes at most $16\sqrt k+1$; rounding the total up to an even
order and allowing one further root cost at most two more.
Bounding the combined integer overhead by $16$, the verified
upper total, divided by $k$, is less than $6.872558140482391$.
\end{proof}

\subsection{Taylor coefficients and conductor savings}
\label{sec:joint-conductors}

Write \(\Lambda=\mathcal O\), \(\Lambda^\vee=\mathfrak D^{-1}\),
and \(n=[\mathbb K:\mathbb Q]\). Boxes are balanced boxes in the
archimedean algebra, with side lengths of order \(X^{1/n}\) at
volume \(X\). All marked prime ideals are unramified of degree
one in the narrow principal class, as in the preceding
construction. The additive character is \(e_{\mathbb K}\) from
Section~\ref{sec:nf-global}. In local calculations \(p\) denotes
the rational prime below the prime ideal in question. Products
of local factors are indexed by the prime ideals.
We use \(\mathcal N\) for the absolute norm of either an element
or an ideal. Let \(\mathfrak m\) be the squarefree ideal of the marks,
put \(m=\mathcal N\mathfrak m\), and write \(\mathfrak q=(q)\mid
\mathfrak m^k\) for the conductor, with \(Q=\mathcal N\mathfrak q\) and
\(U=\log_X Q\). Generators of divisors of marked moduli are chosen
balanced. The number of marks is bounded at fixed \(k\). All counts
of divisors of these moduli are therefore bounded at fixed \(k\),
including after any fixed number of operations. Primes dividing a
fixed denominator of the dual lattice may be omitted. A congruence
for a dual-lattice variable modulo \(a\) means reduction modulo
\(a\Lambda^\vee\). Locally at every marked prime this is the usual
congruence in \(\mathbb Z_p\).

We allow a sum on one residue disk at each prime of \(\mathfrak m\),
with a smooth weight on a translated box of volume \(X\), a real
linear character, and a primitive finite character at each active
prime. The local numerators \(\Phi_{\mathfrak p}\) are given by
integral Taylor series on the disks. For a fixed number of terms
their Taylor coefficients reduce to those of
\begin{equation}\label{eq:joint-germ}
 C_0+C_1t+A_0(1+B_0t)^v,\qquad A_0B_0\ \text{units}.
\end{equation}
The powers \(v\) used throughout will be nonconstant rational
fractional linear functions of \(k\). One excludes the finitely
many exceptional degrees at which a required binomial coefficient
vanishes or is undefined. At a remaining degree we omit the finite
set of primes dividing the relevant rational numerators and
denominators. The marked prime intervals eventually avoid this set.
For the original \(k\)th power one may start with \(v=k\).

A profile \(g(U)\) is available on these disks if the corresponding
sum is
\begin{equation}\label{eq:joint-saving}
 \ll (2mQX)^\varepsilon m^{C} X^{1-g(U)}
\end{equation}
with a fixed \(C\) for each finite list of estimates. In a block of
modulus at most \(X^{M/k}\), the loss in the exponent from the radical,
including the sum over initial disks, is \(O(M/k)\) for a fixed list
on a compact range. The implied constants in \eqref{eq:joint-saving}
may depend on the field, the degree and the list of operations.
The power \(C\) does not depend on the field or on \(k\).

\paragraph{Packing and the shifted mean value.}
If \(\mathfrak a\) divides a marked modulus, a box of volume \(V\)
contains \(O(1+V/\mathcal N\mathfrak a)\) elements of a prescribed
class modulo \(\mathfrak a\). In a box centred at zero the count of
nonzero elements of the lattice \(\mathfrak a\Lambda\), or of
\(\mathfrak a\Lambda^\vee\), is \(O(V/\mathcal N\mathfrak a)\).
Indeed the norm of any such nonzero element is bounded below by a
fixed positive multiple of \(\mathcal N\mathfrak a\). Packing gives
the first assertion and then the second. The same argument applies
to a quotient fractional lattice in counting small residues.

Fix \(D\ge3\), and write
\[
 v_{\mathfrak p}(a)=
       \max\{1,\lceil v_{\mathfrak p}(\mathfrak q)/(D+1)\rceil\}
       \quad(\mathfrak p\mid\mathfrak m),\qquad
 A=\mathcal N(a),\qquad J=a^{D+1}/q .
\]
Thus \(J\) is integral, \(\mathcal N(J)\le m^{D+1}\) and
\[
       X^\xi\le A\le X^\xi m,\qquad \xi=U/(D+1).
\]
The local unit factors of the primitive characters may be folded
into the numerators. If there is no auxiliary multiplier put
\(b=1\), \(B=1\). Otherwise use the balanced generators \(b\) of
distinct prime ideals of norms in \([B,2B)\), disjoint from
\(\mathfrak m\), with residue degree one in the narrow principal
class. Here \(B=X^\beta\) for a fixed \(\beta>0\), and the prime
ideal theorem gives \(\gg B/\log(2B)\) such generators. A fixed
positive proportion of these is enough. We only use this count
and packing in residue classes.

Average the shifts \(x\mapsto x+aby\) over a smooth box of volume
\(Y=X/(AB)\) in \(y\), using the same volume for all \(b\) of the
indicated size. The Taylor phases have degree at most \(D\) modulo
the dual lattice. Their last two coefficients are
\begin{equation}\label{eq:joint-taylor}
    b^D f(x)/a,\qquad b^{D-1}g(x)/a^2,\qquad
    f=Jf_0,\quad g=Jg_0,\quad g_0'=D f_0,
\end{equation}
where \(f_0=\Phi^{[D]}, g_0=\Phi^{[D-1]}\) are divided derivatives,
interpreted at each prime in the character. The real linear twist
affects only the first coefficient. Let \(\mathcal P\) count pairs
of centres \((x,b),(x',b')\) in the common enlarged root box and
the same input disks for which the two coefficients in
\eqref{eq:joint-taylor} differ on the torus by vectors of coordinate
radii \(O(Y^{-D/n})\) and \(O(Y^{-(D-1)/n})\), respectively. Then
\begin{equation}\label{eq:joint-vmvt}
 \frac{|S|}{X}\ll X^\varepsilon
             \left(\frac{\mathcal P}{B^2X^2}\right)^{1/[D(D+1)]}.
\end{equation}

Here are details of the weighted estimate. Expand the shifted
weight in a Fourier series in \(y/Y^{1/n}\) on a fixed enlarged
box. Its coefficients decay uniformly in \((x,b)\) by the
bounds for the scaled derivatives and the balance of \(ab\).
Each Fourier mode changes the first coefficient by a common
linear twist. For that mode partition the coefficient torus
into cells with side lengths comparable to \(Y^{-i/n}\) at
coefficient \(i\). If \(w=D(D+1)/2\), each cell has volume
\(\asymp Y^{-w}\). By iterated one-variable integration on the
cell, Hölder's inequality bounds the supremum to order \(2w\)
by the cell average of the sum of the integrals of scaled
mixed derivatives to that order. The derivatives insert bounded
weights in the sum over \(y\). The critical weighted mean-value
theorem \cite[Theorem~15.1]{NFWooley} therefore bounds the sum of these cell suprema
to order \(2w\) by \(O(Y^{2w+\varepsilon})\). If \(h_C\) is the
occupancy of cell \(C\), the sum of \(h_C^2\) is at most
\(\mathcal P\). Cauchy's inequality for the occupancies, followed
by Hölder over the \(O(BX)\) centres, gives \eqref{eq:joint-vmvt}
after division by the number of shifts. The logarithm from the
number of allowable primes and the summation of the Fourier
coefficients fit the \(X^\varepsilon\) factor.

\paragraph{The two-coefficient estimate.}
With \(b=1\), put \(\mathfrak r_i=\mathfrak q/(\mathfrak q,(a)^i)\)
for \(i=D-1,D\), and let \(R_1=\mathcal N\mathfrak r_D\),
\(R_2=\mathcal N\mathfrak r_{D-1}\). Each numerator is a unit
multiple at the relevant primes of \(\Phi^{[i]}\).
By \eqref{eq:joint-germ} both numerator maps have unit first
derivative, and the penultimate one has unit second derivative
on the input disk. Thus both preserve the valuations of differences
on that disk. For an exact top-coordinate coincidence we have
\(h=x'-x\in\mathfrak r_D\). For \(h\ne0\), group by
\(\mathfrak c=((h),\mathfrak r_{D-1})\).
If \(\mathfrak c=\mathfrak r_{D-1}\), the count is \(O(X^2/R_2)\).
Otherwise, after division by a generator of \(\mathfrak c\),
the difference of the penultimate numerators and its derivative
are units at each prime of \(\mathfrak r_{D-1}/\mathfrak c\).
Consequently a prescribed output determines at most one input
class at each such prime. A tolerance of volume \(\delta\) for
the penultimate coordinate leaves \(O(\delta R_2/\mathcal N\mathfrak c)\)
nonzero values and hence that many classes. Each class contains
\(O(1+X\mathcal N\mathfrak c/R_2)\) inputs and there are
\(O(X/\mathcal N\mathfrak c)\) possible nonzero \(h\). The
off-diagonal count for exact top equality is thus
\begin{equation}\label{eq:joint-diagonal}
 \mathcal N_1(\delta)\ll X^2/R_2
               +\delta(XR_2/R_1^2+X^2/R_1).
\end{equation}
Throughout, bounded divisor and residue-class factors are included
in the implied constant.

For nonzero top residue of tolerance volume \(\delta_1\),
there are \(O(\delta_1 R_1)\) possible residues. For fixed \(x\)
and such a residue, all candidates \(x'\) belong to one class
modulo \(\mathfrak r_D\). Two distinct candidates have exactly
equal top coordinates and penultimate difference of tolerance
\(O(\delta_2)\), where \(\delta_2\) is the original penultimate
tolerance. Conversely, given one candidate and the top residue,
there are \(O(1+X/R_1)\) possible centres \(x\). Applying
\(h\le 1+h(h-1)\) to nonempty candidate sets gives
\[
 \mathcal P\ll X+\mathcal N_1(\delta_2)+\delta_1 R_1 X+
                  \delta_1(R_1+X)\mathcal N_1(O(\delta_2)).
\]
These are the counts in Wooley's two-coefficient estimate
\cite[Lemmas A.8--A.9]{Wnew}. Substituting \(\delta_1=Y^{-D}\),
\(\delta_2=Y^{1-D}\) in \eqref{eq:joint-vmvt} gives its
finite version. In detail, for \(\theta\ge\log_Xm\) and
\(\xi+\theta<1\), the normalized pair count is
\(O(X^{-\min(1,E_1,E_2,E_3,E_4)+\varepsilon})\), where
\begin{align}
 E_1&=2\xi-(D-1)\theta,&
 E_2&=D-1-(D-2)\xi-(2D-1)\theta,\notag\\
 E_3&=(2D-2)(1-\xi)-(3D-1)\theta,&
 E_4&=D+1-U-D\theta .
 \label{eq:joint-finite-terms}
\end{align}
For this substitution use \(R_1\ge X^{\xi-D\theta}\),
\(R_2\ge X^{2\xi-(D-1)\theta}\),
\(R_2/R_1^2\le X^{D\theta}\) and \(Y\ge X^{1-\xi-\theta}\).
The term \(\delta_1 R_1/X\) has saving \(E_4\). The three
terms of \(\mathcal N_1/X^2\) have savings \(E_1\),
\(E_2+1-\xi\) and \(E_2\), and the three corresponding terms
after multiplication by \(X\delta_1\) have savings \(E_2\),
\(E_3+1-\xi\) and \(E_3\). For \(D\ge3\), the limiting
saving from these estimates is
\begin{equation}\label{eq:joint-taylor-saving}
 T_D(U)=\frac{\max(0,\min(1,2U/(D+1),D+1-U))}{D(D+1)}.
\end{equation}
When \(\xi+\theta\ge1\), this positive part is at most
\(\theta/D\), so the trivial estimate with a radical loss is
sufficient for the limiting use.

For the finite application to an original monomial, one can avoid
a loss of \(m\) from the disk partition. At each prime partition
the nonzero residues so that both \(x^{k-D}\) and \(x^{k-D+1}\)
are injective on each part, using at most \(k^2\) parts. Use the
joint partition in the shifted estimate. Inside one part the
valuations of both differences agree with that of \(x'-x\), also
when the inputs have different residues. If \(h\) is a unit at a
prime in the preceding count, the numerator of the divided
difference is a unit by injectivity of \(x^{k-D+1}\), and its
derivative is a unit by injectivity of \(x^{k-D}\). The local
polynomial has bounded degree, so there are at most \(k\)
initial residues for a given value and each Hensel lift is unique.
This proves \eqref{eq:joint-finite-terms} for the monomial with a
constant depending on \(k\) and the number of marks, with no further
radical loss. Here \(3\le D<k\) and the marked residue characteristics
exceed \(k\).

\paragraph{Averaging prime multiples of the Taylor step.}
The estimates in this paragraph apply on the fixed disks, also
to transformed phases. Impose exact equality of the degree-\(D\)
coefficients. Let \(h\) be a generator of \(((J),(a))\), write
\(r_0=a/h\), \(R=\mathcal N(r_0)\), and put
\(\delta\ll m^{O_D(1)}Y^{1-D}\). From \eqref{eq:joint-taylor}
the top equality gives
\[
 b^D g_0'(x)\equiv b'^D g_0'(x')\pmod{r_0}.
\]
Multiplying the penultimate coordinate by a generator of
\((h)^2/((h)^2,(J))\) gives tolerance \(\delta\) for a
primitive numerator modulo \(r_0^2\). At each prime of \(r_0\)
colour the multipliers so that the \(D\)th-power map is
injective on each colour. Retaining the largest joint colour
costs only a bounded factor. The top equality now forces
\(b\equiv b'\pmod{\mathfrak p}\) at each prime of \(r_0\).

Suppose that \(b,b'\) are distinct. Put \(d\) for a generator of
\(((b-b'),(r_0))\), \(s=r_0/d\), and let \(H=(bx'-b'x)/d\).
We give the local calculation which both proves integrality of
\(H\) and yields the discrepancy estimate. For
\(\mathfrak p^\alpha\Vert(r_0)\), write
\(t=b'/b\), \(\vartheta=t-1\), \(G=g_0\) and \(g_1=G'\).
Inverting \(g_1\) on the disk gives the analytic function \(w(x)\)
such that \(g_1(w(x))=t^{-D}g_1(x)\). Then
\(x'\equiv w(x)\pmod{p^\alpha}\). The function
\[
 \psi(x)=(w(x)-t x)/\vartheta
\]
has integral Taylor coefficients by analytic inversion in
\(x\) and \(\vartheta\). Modulo \(p\) its jet is
\(-x-Dg_1(x)/g_1'(x)\) and its derivative is the unit
\(-v/(v-D)\). If \(j_0=v_p(\vartheta)<\alpha\), the relation
\(x'-t x\equiv\vartheta\psi(x)\pmod{p^\alpha}\)
implies \(d\mid bx'-b'x\) locally and fixes
\(x\equiv x_*=\psi^{-1}(z)\pmod{p^{\alpha-j_0}}\), with
\(z=dH/(b\vartheta)\). If \(v_p(\vartheta)\ge\alpha\),
both \(bx'-b'x\) and the difference of the penultimate
numerators are divisible by \(p^\alpha\).
In the first case set
\[
 F(x,z)=t^{D-1}G(t x+\vartheta z)-G(x),\qquad
 W(z)=F(x_*,z)/\vartheta .
\]
The function \(F/\vartheta\) is integral analytic and
\(F_x(x_*,z)=0\). Thus \(F(x,z)\equiv\vartheta W(z)
\pmod{p^{2\alpha-j_0}}\). The penultimate numerator divided
by \(d\) is consequently, modulo \(p^{2(\alpha-j_0)}\), a
unit times \(W\) evaluated at a unit times \(H\). If
\(x'_*=t x_*+\vartheta z\), the derivatives are
\[
 W'(z)=t^{D-1}g_1(x'_*),\qquad
 W''(z)=t^{D-1}g_1'(x'_*)\bigl(t/\psi'(x_*)+\vartheta\bigr).
\]
They are units. The value is a unit as well, since
\(W(\psi(x))\equiv (D-1)G-Dg_1^2/g_1'
=-vG/(v-D)\pmod p\).

We need the following count. On a disk at each prime of \(s\),
let a phase modulo \(s^2\) have unit value and unit first
and second derivatives. For \(H\) in a translated box of volume
\(V\ge1\), the number of hits in a torus box around zero of
volume \(\Delta\) is
\begin{equation}\label{eq:joint-discrepancy}
                 O(\Delta V+\mathcal N(s)).
\end{equation}
If \(\Delta\mathcal N(s)^2\) is below a fixed positive constant,
there is no hit by the nonzero norm bound. When
\(V\ge\mathcal N(s)^2\), injectivity on the disks and packing
give the count. For \(V<\mathcal N(s)^2\), smooth the two box
indicators. In the output expansion the coefficients at
\(j\in\Lambda\) are
\(O_N(\Delta(1+\Delta^{1/n}|j|)^{-N})\), with \(|j|\)
a fixed archimedean norm. Input completion on the dual lattice
has normalized Fourier coefficients of bounded total absolute sum.
The zero output frequency contributes \(O(\Delta V)\).
A nonzero one gives a complete sum, with any linear twist, of
size at most a constant times
\(\mathcal N(s)\mathcal N(((j),(s)^2))^{1/2}\).
To see this, divide the nonlinear coefficient by its local
conductor gcd. When a nontrivial residual modulus remains, the
linear twist must permit that same division, and the unit second
derivative gives the square-root bound on the disk by stationary
phase. Exponent one on a disk and primes of zero residual conductor
satisfy the bound trivially. For any \(\mathfrak c\mid(s)^2\),
the sum of the indicated Fourier weight over nonzero
\(j\in\mathfrak c\) is \(O(1/\mathcal N\mathfrak c)\) by the
nonzero norm count in boxes. Summing the gcd bound over divisors
proves \eqref{eq:joint-discrepancy}.

In the pair problem take \(V\ll BX/\mathcal N(d)\) and
\(\Delta\ll\mathcal N(d)\delta\). The local inverse fixes \(x\)
modulo \(s\) for each \(H\). The equation \(bx'-b'x=dH\) also
fixes \(x\) modulo \(b\), since the multiplier ideals are
distinct primes. Hence each \(H\) has
\(O(1+X/(B\mathcal N(s)))\) possibilities for \(x\).
For \(BR<X\) with a fixed exponent gap, the count divided by
\(X^2\) is at most a constant times
\[
                  (BX)^{-1}+\mathcal N(d)\delta/R .
\]
For fixed \(b,d\) the candidates \(b'\) lie in one class
modulo \(d\). As \(b-b'\ne0\) and its norm is \(O(B)\),
\(\mathcal N(d)\ll B\). Their proportion is
\(O(X^\varepsilon/\mathcal N(d))\). If
\(\xi>\beta\), then \(R\gg B\) for sufficiently small
\(\log_Xm\), so \(s\) is nontrivial in this argument. For
\(b=b'\), \eqref{eq:joint-diagonal} with \(R_1=R\) and
\(R_2=R^2\), after removing the rounding factors, gives
\[
 \ll m^{O_D(1)}(X+X^2/R^2+\delta(X+X^2/R)).
\]
The pairs with equal multipliers have proportion
\(O(X^\varepsilon/B)\).
It follows that the exact-top part has normalized count
\begin{equation}\label{eq:joint-zero-top}
 \frac{\mathcal P_0}{B^2 X^2}
   \ll X^\varepsilon m^{O_D(1)}(BX)^{-1}
\end{equation}
when \(2\xi\ge1\), \(\xi>\beta\), \(\xi+\beta<1\) and
\[
                 \xi+(D-1)(1-\xi-\beta)>1+\beta .
\]
Indeed the last condition pays for the term \(\delta/R\).
If \(D-1\le U<D\), the additional condition
\(\beta<(D-U)/D\) makes every top pair exact, by the
nonzero norm bound and \(A Y^{-D}\ll m^{D+1}X^{U-D+D\beta}\).
We recover the profile
\begin{equation}\label{eq:joint-prime-saving}
 M_D(U)=\frac{1+\min((D-U)/D,\,(D-2)(1-\xi)/D)}{D(D+1)}
 \qquad (D-1\le U\le D).
\end{equation}
The endpoints are understood by approximation, using
\eqref{eq:joint-taylor-saving} at \(U=D\), and using the
pointwise bound with an arbitrarily small slack. This extends
the prime-multiplier calculation of \cite[\S6]{HexWaring}
to the lattice used here.

\paragraph{The two arithmetic transformations.}
Complete on the disks to make \(g(U)=U/2\) available for
\(0\le U\le1\), and \(1-U/2\) for \(U\ge1\), up to radical
loss. The profiles submitted to both arithmetic transforms
below are continuous and nonincreasing for \(U\ge1\) and
satisfy \(g(1)=1/2\). On each affine piece \(aU+b\) we have
\begin{equation}\label{eq:joint-profile-class}
 -\tfrac12\le a\le0,\qquad 0\le a+b\le\tfrac12 .
\end{equation}
Then the function \(xg(y/x)\) is nondecreasing in \(x\)
and is \(1\)-Lipschitz in \(x,y\ge0\), with value zero
at \(x=0\). This follows piecewise from \(ay+bx\), and
also holds on \(y<x\) with the completion profile.

For the \(A\) process an available point \((u,g(u))\) gives
\begin{equation}\label{eq:joint-A}
              (u+1-g(u),\,g(u)/2).
\end{equation}
Write \(c=1-g(u)\), \(U=u+c\). Choose shift modulus with
local exponents
\(\max(1,\lceil c\,v_p(\mathfrak q)/U\rceil)\), giving
norm between \(X^c\) and \(mX^c\). If the shift box has
volume below one, \(m^{1/2}\) covers the stated saving.
Otherwise the diagonal term in the squared normalized sum
costs \(mX^{c-1}\). For an off-diagonal shift indexed by \(h\),
the remaining conductor differs from \(X^u\) by a factor
at most \(m\mathcal N(((h),\mathfrak q_0))\),
where \(\mathfrak q_0\) is the conductor after the fixed
modulus is removed. The average of this gcd norm is bounded:
expand by its divisors and use the nonzero packing count
for each divisor. The local phase after division reduces to
the derivative of the input jet and has power \(v-1\).
The loss from changing the conductor inside a \(1\)-Lipschitz
profile is thus absorbed before the square root. This proof
of \eqref{eq:joint-A} also permits \(0\le u<1\) and any
nonnegative \(1\)-Lipschitz input when the process is used
on its own.

For the \(B\) process let \(\mathfrak q_1\) be the least
common multiple of \(\mathfrak q,\mathfrak m\), and let
\(\mathfrak q_2\) retain the conductor primes with exponent
at least two. At those primes the single-disk transform is
supported on one dual disk. It has modulus \(p^{e/2}\) and
the stationary phase is the Legendre phase with power
\(v/(v-1)\), after affine changes with unit slope. The
quadratic Gauss factor is constant on that disk. The exponent-one
or inactive mask primes contribute only a linear character to
the dual sum. Poisson gives a smooth translated dual weight at
volume \(Y=\mathcal N\mathfrak q_1/X\) with factor
\(X(\mathcal N\mathfrak q_2)^{1/2}/\mathcal N\mathfrak q_1\).
Cover that weight with translated boxes at volume \(Y\);
their scaled derivative bounds are summable. Clearing the
fixed dual-lattice denominators and taking its finite cosets
only changes the linear modulation. The homogeneous Lipschitz
bound and \(\mathcal N\mathfrak q_1/\mathcal N\mathfrak q_2\le m\)
permit replacing both norms by \(Q\) at a radical cost.
For \(U>1\) this gives the saving
\begin{equation}\label{eq:joint-B}
          1-U/2+(U-1)g\bigl(U/(U-1)\bigr).
\end{equation}
On each affine piece the new slopes in \eqref{eq:joint-A}
and \eqref{eq:joint-B} are respectively \(a/[2(1-a)]\)
and \(a+b-1/2\). Both transformations, followed by the
maximum with completion, preserve the conditions in
\eqref{eq:joint-profile-class} used in the computation.
For any nesting of these operations a finite number of
Taylor coefficients of \eqref{eq:joint-germ} suffices.

\subsection{Completion with two multipliers}
\label{sec:joint-completion}

We keep the notation of \cref{sec:joint-conductors}, and compare
the two highest coefficients for \(b,b'\) in
\eqref{eq:joint-taylor}. The case \(b=b'=1\) will give the
completed estimate of \cite[\S7]{HexWaring}. Allowing the
prime multipliers to vary gives an additional saving near the
integer conductors. Put
\begin{equation}\label{eq:joint-lengths}
 H=A/Y^D,\qquad T=Y^{D-1},\qquad L=A^2/X,
 \qquad HT=BL .
\end{equation}
Assume that \(AB<X\), with a fixed exponent gap. Small radical
losses will always fit inside this gap. We represent the top
difference by \(r\in\Lambda^\vee\) in a centred box of volume
\(O(H)\), and impose
\begin{equation}\label{eq:joint-top-congruence}
       b^D f(x)-b'^D f(x')\equiv r\pmod a .
\end{equation}
A product of normalized Fej\'er kernels on the torus
majorizes the penultimate tolerance. Its frequencies
\(j\in\Lambda\) lie in a box of volume \(O(T)\) and
have coefficients \(O(T^{-1})\). We use smooth nonnegative
majorants for the two root boxes. The case \(r=0\) has already
been counted in \eqref{eq:joint-diagonal} or
\eqref{eq:joint-zero-top}. When \(r\ne0\), the frequency \(j=0\)
contributes a count
\begin{equation}\label{eq:joint-zero-fejer}
           \ll m^{O_D(1)}HX(X/A)/T
\end{equation}
for each pair of multipliers, since for fixed \(x,r\) the
congruence fixes \(x'\) modulo \(a/((a),(J))\) on the disks.
The factor \(1+X/A\) in this bound is \(O(X/A)\).

\paragraph{Completion on the two root variables.}
Take \(r,j\ne0\), and write \(K_r(j)\) for the smoothed
pair sum under \eqref{eq:joint-top-congruence}, with
phase \(j(b^{D-1}g(x)-b'^{D-1}g(x'))/a^2\).
Poisson summation modulo \(a^2\) in both variables gives
dual frequencies \(u,u'\in\Lambda^\vee\) of volume scale
\(L\) with weights of size \(O(X/A^2)\) and rapid
decay beyond that volume. In the complete transforms, split
\(x\) as \(x_0+a y\) and likewise for \(x'\). Summation over
the \(y\)'s imposes
\[
 u+Dj b^{D-1} f(x_0)\equiv0,\qquad
 u'-Dj b'^{D-1} f(x'_0)\equiv0\pmod a .
\]
Thus \(bu+b'u'+Dj r\in a\Lambda^\vee\).
On the effective dual boxes of volume \(O(LX^\varepsilon)\)
the norm of this expression is \(O(BLX^\varepsilon)=o(A)\),
using \eqref{eq:joint-lengths}. The norm lower bound for
nonzero elements implies the exact equation
\begin{equation}\label{eq:joint-line}
                         bu+b'u'=-Dj r .
\end{equation}
The contribution outside those boxes is arbitrarily small
by taking enough derivatives of the weights.

Let \(\mathfrak d=((jJ),(a))\), and put
\[
    \Gamma=\mathcal N(((j),(a)))\,\mathcal N(((r),(a))),
    \qquad \gamma=\log_X\Gamma .
\]
A gcd involving \(r\in\Lambda^\vee\) means the valuations
of \(r\) at the marked primes. The frequencies must belong to
\(\mathfrak d\Lambda^\vee\). For each \(u\) the first stationary
condition gives at most \(\mathcal N\mathfrak d\) values of
\(x_0\bmod a\) and then the top congruence gives
\(O(\mathcal N(J))\) values of \(x'_0\).
The absolute complete coefficient is therefore
\(O(A^2\mathcal N\mathfrak d\,\mathcal N(J))\).
If \(b\ne b'\), the two prime ideals are coprime, and
\eqref{eq:joint-line} and the divisibility by \(\mathfrak d\)
put \(u\) in one translate of \(b'\mathfrak d\Lambda^\vee\).
Put \(B_*=B\) in this case. For equal multipliers put \(B_*=1\).
Then the frequency count is
\(O(1+LX^\varepsilon/(B_*\mathcal N\mathfrak d))\).
If the equal multiplier is a prime \(b\), \eqref{eq:joint-line}
also forces \(b\mid jr\), since \(b\) avoids the denominators of
\(\Lambda^\vee\) and divides neither \(D\) nor any mark.

Group \(j,r\) by their gcd ideals with \(a\).
In each group there are \(O(HT/\Gamma)\) pairs \((j,r)\),
by the centred nonzero packing bound. Equal prime
multipliers allow only \(O(HT/(B\Gamma))\) pairs, on
separating the two possibilities for \(b\mid jr\).
Since \(\mathcal N\mathfrak d\le\mathcal N(J)\Gamma\),
the \(1\) in the frequency count, combined with the
gcd-group count, costs at most
\[
       m^{O_D(1)}X^\varepsilon \frac{XH}{B_*}
                                  \frac{B_*}{L}
\]
after multiplication by the two Poisson weights and the
Fej\'er coefficient. In the equal prime case there is the
extra divisor factor \(B^{-1}\). The other term is at most
\(m^{O_D(1)}X^\varepsilon XH/(B_*\Gamma)\), again with
this divisor factor for equal primes. This estimate keeps
the \(1\) even when \(H\) is of constant size. It deals
with every group with \(\gamma\ge\eta\) provided
\(\log_X(L/B_*)>\eta\) with fixed slack.

\paragraph{Separating the prime powers.}
For groups with \(\gamma<\eta\), split \(a=a_c a_b\) into
balanced generators of coprime factors. Include
\(p^\alpha\Vert a\) in \(a_c\) when
\(\alpha>v_p(jJ)\) and \(\alpha>v_p(r)\).
Use \(A_c=\mathcal N(a_c)\), \(A_b=\mathcal N(a_b)\), and
factor \(\mathfrak d=\mathfrak d_c\mathfrak d_b\) accordingly.
At the primes of \(a_b\) retain the coupled complete kernel.
After restricting to the line, it is periodic modulo \(a_b^2\),
has support in \(\mathfrak d_b\Lambda^\vee\), and has
amplitude \(O(A_b^2\mathcal N\mathfrak d_b\,\mathcal N(J))\)
by the same lift count. Consequently there are at most
\(A_b^2/\mathcal N\mathfrak d_b\) classes on which this
kernel is constant. On each class the line parameter varies
in a translated box with volume scale
\begin{equation}\label{eq:joint-dual-size}
             N_*=L/(B_*\,\mathcal N\mathfrak d_c\, A_b^2).
\end{equation}
Indeed the progression step, expressed in \(u\), generates
\(\mathfrak d_c(a_b)^2\) times \((b')\) for distinct
multipliers, and \(\mathfrak d_c(a_b)^2\) otherwise.
It is a unit on the core primes after division by
\(\mathfrak d_c\). The product of the two smooth Fourier
weights has summable scaled derivatives on this progression.
Removing the centres of the weights leaves an allowed
linear modulation of the free line variable.

At a prime of \(a_c\) write \(\delta_0=v_p(jJ)\).
The top congruence requires \(v_p(r)\ge v_p(J)\).
After division by \(J\), it is already implied by the
stationary conditions modulo \(p^{\alpha-\delta_0}\).
Insert the remaining congruence modulo
\(p^{\alpha-v_p(J)}\) by finite character orthogonality.
This averages products of single-variable complete
transforms, with the same additional character in both
factors and opposite signs. The extra numerators are
multiples of \(p^{\alpha+v_p(J)} b^D f_0\) and of the
corresponding \(b'\) expression. The insertions do not
change the stationary congruences. After factoring
\(jJ b^{D-1}\) and \(jJ b'^{D-1}\), the two local phases
at modulus \(p^{2\alpha-\delta_0}\) have numerators
\(g_0+\tau b f_0\) and \(g_0+\tau b' f_0\), with common
\(\tau\in p\mathbb Z_p\). Their second derivatives are
units. Each transform has size \(p^{\alpha+\delta_0/2}\);
its Gauss factor is constant on the stationary disk.

Let \(P_b(v_0)\) denote the Legendre phase of
\(g_0+\tau b f_0\), with stationary parameter \(v_0\).
Normalizing \(u\) and \(-u'\) by the two extracted
coefficients makes \eqref{eq:joint-line} read
\(v'_0=t^{-D}(v_0+\rho)\), where
\(t=b'/b\) and \(\rho=Dr/(Jb^D)\). The difference phase is
therefore
\begin{equation}\label{eq:joint-legendre}
       P_b(v_0)-t^{D-1}P_{b'}(t^{-D}(v_0+\rho)).
\end{equation}
Signs in the primitive character do not change the estimates.
Both arguments of the Legendre functions are on the same
reduced stationary disk. The jets of the two Legendre
functions agree modulo \(p\), with binomial power
\(\mu=(v-D+1)/(v-D)\) up to an affine function.

If \(t=1\), the difference in \eqref{eq:joint-legendre}
is divisible by \(p^{v_p(\rho)}\), and after division
its jet has power \(\mu-1=1/(v-D)\). Here \(v_p(\rho)\ge1\)
because the original roots were in one disk. In the
general case \(t\equiv1\pmod p\) also forces
\(\rho\equiv0\pmod p\). Divide by \(p^\lambda\) with
\(\lambda=\min(v_p(t-1),v_p(\rho))\).
The dependence of the stationary phase on \(\tau b\)
is analytic, so its change with \(b'\) is divisible by
\(p(t-1)\). Modulo \(p\) after the division, the remaining
non-affine terms are obtained by expanding
\[
       -(D-1)(t-1)P_b + D(t-1)v_0P_b' -\rho P_b'
\]
to first order in \(t-1,\rho\), and then dividing by
\(p^\lambda\). They are a combination of two adjacent
binomial powers \(\mu,\mu-1\). The higher-power
coefficient is a unit multiple of
\((D-1-D\mu)(t-1)/p^\lambda\), and
\(D-1-D\mu=-v/(v-D)\). If this coefficient vanishes,
the lower one is a unit by the choice of \(\lambda\).
For derivative order \(i\ge2\) the ratio of the two
falling factorials is \((\mu-i)/\mu\), which changes at
each successive order. Thus at least one derivative is
a unit in every consecutive pair among orders \(2\)
through \(5\). If \(t\not\equiv1\pmod p\), take
\(\lambda=0\). At the common reduced disk the relative
multipliers of the \(i\)th and \((i+1)\)st derivatives
are \(t^{D-1-Di}\) and \(t^{D-1-D(i+1)}\).
They cannot both equal one unless \(t=1\bmod p\).
This proves the same adjacent derivative condition.
All these computations involve only the finitely many
retained Taylor coefficients in \eqref{eq:joint-germ}.

The effective complete amplitude at core primes is
\(A_c^2\mathcal N\mathfrak d_c\). The new phase along
the progression has modulus norm
\[
             Q_c=\prod_{p^\alpha\Vert a_c}
                             p^{2\alpha-\delta_0-\lambda}.
\]
Since \(\lambda\le v_p(r)-v_p(J)\) here and
\(\mathcal N\mathfrak d_c\,A_b^2\le
\mathcal N(J)^3\Gamma^2\), we have, writing
\(a_0=2\xi-1\), \(l=\log_X N_*\), \(w=\log_X(Q_c/N_*)\),
\begin{equation}\label{eq:joint-dual-exponents}
 \begin{split}
 a_0-\log_X B_*-2\gamma-O_D(\log_Xm)&\le l
                         \le a_0+O_D(\log_Xm),\\
 1+\log_X B_*-\gamma&\le w\le 1+\log_X B_* .
 \end{split}
\end{equation}
For the second inequality one uses
\(Q_c/N_*=B_*X\prod p^{-\lambda}\).
The product of the two Poisson weights, the core and
coupled amplitudes, the number of coupled classes, and
\(N_*\), is
\[
 (X/A^2)^2 A_c^2\mathcal N\mathfrak d_c\,A_b^2
 \mathcal N\mathfrak d_b\,(A_b^2/\mathcal N\mathfrak d_b)\,N_*
 =X/B_*,
\]
up to a radical factor. Thus if the line phase saves
\(g(1+w/l)\) at volume \(N_*\), the saving together with
the gcd-group count is \(\gamma+l g(1+w/l)\).
This accounts for the period of the coupled kernel as
well as the exact line in the Fourier sum.

\paragraph{One multiplier.}
Take \(b=b'=1\) and \(D\le U\le D+1/2\).
Let \(s\) be any available disk profile satisfying
\eqref{eq:joint-profile-class} with \(s(2)<1/2\).
The line now has binomial power \(1/(v-D)\), so
\eqref{eq:joint-saving} applies. Take
\[
             \eta=\min\{a_0s(1+1/a_0),(D-2)(1-\xi)\}.
\]
For \(\gamma<\eta\), the lower bound for \(l\) stays
positive since \(\eta\le a_0s(2)<a_0/2\). The function
\(l s(1+1/l)\) has slope in \([0,1/2]\), by
\eqref{eq:joint-profile-class}, and the slope of
\(s\) has absolute value at most \(1/2\).
Using \eqref{eq:joint-dual-exponents}, the group's
saving is at least \(a_0s(1+1/a_0)-O_D(\log_Xm)\).
The second restriction on \(\eta\) pays for
\eqref{eq:joint-zero-fejer}. The exact top contributes
\(X\) up to losses by \eqref{eq:joint-diagonal}.
The resulting saving in \eqref{eq:joint-vmvt} is
\begin{equation}\label{eq:joint-P}
  P_D[s](U)=\frac{\min(1,D+1-U+\eta)}{D(D+1)} .
\end{equation}
This also covers the case \(s(1+1/a_0)=0\).
The local maps in this operation, as in the \(A\) and
\(B\) processes, are nonconstant fractional linear
functions of the input power. This permits use of
\eqref{eq:joint-P} in a finite composition of the
three operations.

\paragraph{Cubic stationary phase on the line.}
For the phase in \eqref{eq:joint-legendre} after local
division, either the second or the third derivative is
a unit. On a disk with unit third derivative and
second derivative zero modulo \(p\), there is a unique
centre \(c\) with \(F''(c)=0\), by Hensel's lemma. The
central region \(v_p(x-c)\ge\lceil e/3\rceil\) of a
complete sum modulo \(p^e\) contributes \(O(p^{2e/3})\).
On a remaining annulus write \(r_1=v_p(x-c)<e/3\).
Average shifts divisible by \(p^j\), with
\(j=\lceil(e-r_1)/2\rceil>r_1\). They preserve the
annulus and linearize the phase modulo \(p^e\):
the quadratic term has valuation at least \(2j+r_1\)
and the higher terms at least \(3j\). Thus a nonzero
contribution requires \(F'(x)+h\equiv0\bmod p^{e-j}\)
for the linear twist \(h\).
Writing \(x=c+p^{r_1}y\), the difference \(F'(x)-F'(c)\)
is \(p^{2r_1}\) times a unit multiple of \(y^2\)
modulo \(p\), with a unit derivative for \(y\) a unit.
For a prescribed value at most two residue classes
lift, so the normalized measure on this annulus is
\(O(p^{-\max(r_1,e-j-r_1)})=O(p^{1-e/3})\).
A unit second derivative gives the square-root bound.
We conclude, also including modulus \(p\), that the
complete sum is \(O_k(p^{2e/3+1})\), with any linear
twist. Multiplying over the marked primes and applying
Poisson gives
\(g_0(V)=\max(0,\min(V/3,1-2V/3))\) as a saving
up to a radical factor.

After one \(A\) difference the relevant derivative
pair is the original pair of orders \(3,4\), and
after two differences it is the pair \(4,5\).
The same completion estimate is therefore available
at every stage. From \eqref{eq:joint-A}, the first
transform has a rising branch \((V-1)/4\) up to
\(5/3\), then a falling branch \((5-2V)/10\)
up to \(5/2\). Applying \eqref{eq:joint-A} again
gives the two branches of the tent
\(\min((V-2)/6,(7-2V)/24)\) on \([2,7/2]\).
For \(V\ge5/3\) we may use the profile
\begin{equation}\label{eq:joint-core-profile}
          \max\{0,(5-2V)/10,\min((V-2)/6,(7-2V)/24)\}.
\end{equation}
Only the \(A\) process is applied to these intermediate
profiles. The Lipschitz estimates in its proof apply
across their breakpoints.

\paragraph{Prime multipliers near integer conductors.}
Let \(\beta,\eta>0\) with \(\eta>D\beta\), and take
\[
 U_-=D-D\beta,\qquad U_+=D+\eta-(D-1)\beta,\qquad
 \xi_\pm=U_\pm/(D+1),\qquad a_\pm=2\xi_\pm-1 .
\]
On \([U_-,U_+]\) impose the following strict tests:
\begin{align}
 &2\xi_->1,\quad \xi_->\beta,\quad
         \xi_++\beta<1,\quad
         (D-2)(1-\xi_+)-D\beta>\eta, \notag\\
 &a_--\beta-2\eta>0,\qquad
     (1-\eta)/a_+>2/3,\notag\\
 &\eta<\max\left\{
   \frac{3a_--2-5\beta}{10},\
   \min\left(\frac{1-a_+}{6},
             \frac{5a_--2-7\beta}{24}\right)\right\}.
 \label{eq:joint-parameter-tests}
\end{align}
The first line gives the exact Fourier equation, allows
\eqref{eq:joint-zero-top} when the top residue is zero
and handles \eqref{eq:joint-zero-fejer} otherwise.
It also gives the conditions for the comparison of
prime multipliers in the exact-top count. The second
line ensures positive line volume and \(1+w/l\ge5/3\)
for small gcd groups. To verify the saving in the last
line, substitute \eqref{eq:joint-dual-exponents} into
\eqref{eq:joint-core-profile} and multiply by \(l\).
The three expressions are
\((3l-2w)/10\), \((w-l)/6\) and \((5l-2w)/24\).
Use the lower bound on \(l\) for the increasing
ones and the upper bound for the decreasing one.
After adding \(\gamma\) each lower bound is
nondecreasing in \(\gamma\), so their maximum
with the tent taken as a minimum is bounded below
by the last line of \eqref{eq:joint-parameter-tests}.
The radical and archimedean losses fit in the strict
gaps at each fixed \(D,\beta,\eta\).

It follows that the \(r\ne0\) contribution has
normalized count
\(O(X^{-1+(U-D)+(D-1)\beta-\eta+\varepsilon}
m^{O_D(1)})\).
For unequal multipliers this uses the factor \(X/B\)
from \eqref{eq:joint-dual-size} and the divisor count.
For equal primes the factor \(X\) is compensated
by \(b\mid jr\). In both cases the count is summed
over the chosen multipliers and normalized by their
squared number. The large-gcd \(1\)-term is negligible
since \(a_--\beta>2\eta\). The zero top count from
\eqref{eq:joint-zero-top} has normalized size
\(O(X^{-1-\beta+\varepsilon}m^{O_D(1)})\).
Together these give the available profile
\begin{equation}\label{eq:joint-H}
 H_D(U)=
 \frac{\min(1+\beta,\ 1-(U-D)-(D-1)\beta+\eta)}{D(D+1)}
 \qquad (U_-\le U\le U_+).
\end{equation}
Its left endpoint agrees with \eqref{eq:joint-prime-saving}
when \(\beta<(D-2)(1-\xi_-)/D\), which holds for
the parameters used below. Its right endpoint is
\(1/[D(D+1)]\). This estimate applies to the generic
disks of \eqref{eq:joint-germ}, so it can also be used
before a finite sequence of arithmetic transformations.

\subsection{Refinements of the conductor profile}
\label{sec:refined-conductors}

The next estimates improve the line phase in
\eqref{eq:joint-legendre}. We retain the notation of
Sections~\ref{sec:joint-conductors} and~\ref{sec:joint-completion}.
In this subsection $q$ denotes the conductor norm, and lengths of
root or Fourier boxes refer to their volumes in the archimedean
algebra. The local calculations are over $\mathbb Z_p$ at the
degree-one marks. The supports obtained by completion are ideal
cosets. A balanced generator of a divisor ideal of norm $d$
changes the progression volume by the factor $d$. In particular,
Poisson summation at volume $X$ gives dual volume $q/X$, and the
packing bounds on a coset of modulus norm $d$ are
$O_{\mathbb K}(1+X/d)$. In centred boxes the nonzero count is
$O_{\mathbb K}(X/d)$. This latter bound follows first from the
norm of a nonzero element and then from packing. It supplies the
divisor averages after the arithmetic shifts.

\paragraph{Averaged cubic completion}

At a prime power \(p^e\), quadratic completion applies wherever
\(F''\) is a unit. Otherwise choose \(F''(c)=0\) by Hensel's lemma.
On an annulus \(r=v_p(x-c)<e/3\), the complete coefficient at
Fourier label \(h\) is supported on \(h=-F'(c)\bmod p^{2r}\)
and is \(O(p^{(e+r)/2+1/2})\). Shifts of valuation
\(\lceil(e-r)/2\rceil\) linearize the phase and force a stationary
congruence. After rescaling, the reduced part of \(F'\) is
quadratic with a unit derivative on the units. The central annulus
is linear and has length \(p^{e-\lceil e/3\rceil}\).

Combine support classes by CRT. For fixed annuli put
\(T=\prod p^{r/2}\). The complete amplitudes are bounded
by \(\sqrt q\,T\) up to radical factors and lie on a class
of modulus \(T^4\) up to the same factors. If \(q_b\) is the
portion of the conductor where \(F''\) is a nonunit, the sum
of amplitudes on a balanced box of volume \(R\) is
\(O(\sqrt q\,(R+q_b^{1/6})m^{O(1)})\). Poisson summation gives
\[
 |S(X,q)|\ll m^{O(1)}
        \bigl(Xq^{-1/2}q_b^{1/6}+\sqrt q\bigr).
\]
Without controlling \(q_b\), the saving is the positive part of
\(\min(V/3,1-V/2)\) at \(q=X^V\).

\paragraph{Choosing shifts by derivative vanishing}

A shift of nominal logarithmic size \(c\) can remove up to \(c\)
from the conductor on a specified set of primes. Assign the
individual prime exponents before rounding. The rounding costs
a fixed radical power and the incidental gcd is handled by
the divisor average in the arithmetic \(A\) process. The reduced
phase on the remaining primes is the derivative of the input.

One \(A\) step gives \(\min((V-1)/4,(3-V)/6)\) on its positive
support. For two steps with shift exponents \(1-2t,1-4t\), let
\(B_i\) be the logarithmic part of the conductor where derivative
\(i+2\) vanishes. Then \(B_1+B_2\le V\). If
\(B_2\le3V-4-12t\), concentrate the shifts on that set. The
remaining conductor exponent is \(V-2+6t\). The cubic bound
gives saving \(4t\) under \(t\le(V-2)/2\) and \(t\le(4-V)/14\).
If \(B_2>3V-4-12t\) and \(t\le(2V-3)/14\), then
\(B_1<1-2t\), and the first shift can clear the entire set.
The residual phase has unit second derivative. Its Legendre
phase has a unit third or fourth derivative. Its one-step
cubic bound, transformed back by \(B\), has saving
\(1/2-(V-1+2t)/6\) at the shortened stage on the relevant range
\(30/13\le V\le13/5\). This is at least \(2t\) for
\(t\le(4-V)/14\). The original one- and two-step bounds overlap
on either side.

For \(j\ge3\) the shift exponents \(1-2^i t\), \(1\le i\le j\),
similarly give the positive tent
\[
 \min\left(\frac{V-j}{2},\frac{j+2-V}{2^{j+2}-2}\right)
 \qquad (j\le V\le j+2).
\]
With \(D_0=2^j\), the \(j\)-step argument succeeds if
\(B_j\le3V-2j-(2D_0+4)t\). On failure, \(B_{j-1}\) is small
enough to remove by the end of \(j-1\) shifts because
\((3D_0+2)t\le2V-j-1\). Use quadratic completion for a residual
exponent at most \(3/2\), and the \(B\) image of one cubic \(A\)
step in the next range. The ordinary \(j\)-step cubic
estimate covers the remaining falling part.

\paragraph{Fourier transformation on stationary annuli}

On a noncentral annulus write \(x=c+p^r y\) and
\(h=-F'(c)+p^{2r}w\). The nonlinear part has a factor
\(p^{3r}\). After division the local phase has reduction
\(a y^3+wy\) on the units, with \(a\) a unit. Its stationary
equation is quadratic in \(y\). On each simple stationary
disk the transformed phase has power \(3/2\), conductor
\(p^{e-3r}\), and complete amplitude \(p^{(e+r)/2}\)
up to a radical factor. The quadratic Gauss factor is
constant on a fixed residue disk. At the good primes use
the ordinary Legendre transform. A residual exponent of at
most one is covered by the radical loss and a residue-class
split. Central balls give a
linear character, with nominal \(r=e/3\) and a bounded
rounding loss. Inserting a common linear character changes
only the centre of the archimedean Fourier weight.

Write \(q=X^V\), \(P_1=\prod_p p^r=X^\tau\), and let
\(b\) be the logarithmic weight of the bad primes. Then
\(0\le\tau\le b/3\). The support progression has step of
size \(P_1^2\) and the reduced conductor is \(q/P_1^3\), up
to radical powers. Thus a dual profile \(g\) gives saving
\[
  1-\frac V2+\frac{3\tau}2+
 (V-1-2\tau)g\!\left(\frac{V-3\tau}{V-1-2\tau}\right)
\]
when the dual box has volume above one. The expression
is nondecreasing in \(\tau\) whenever \(3a+2d\le3/2\)
on a dual affine cell \(g(W)=aW+d\). This is checked for
each input to the transform. The short-box bound is
\(V/2-b/6\) since there are only \(O(1)\) frequencies
per translated box. Define \(S_g(V)=1/2\) for
\(0\le V\le1\) and for \(V>1\) put
\[
 S_g(V)=1-\frac V2+(V-1)g\!\left(\frac V{V-1}\right).
\]
The available bound is the positive part of
\(\min(V/2-b/6,S_g(V))\). In particular it is at least the
positive part of \(\min(V/3,S_g(V))\). This defines the
annular transform used in the profile. The ordinary quadratic
transform on a unit stationary disk remains available when
\(b=0\).

\paragraph{Repeated transforms}

Between successive uses of the annular transform, apply
one or more arithmetic differences. Once a bad annulus
occurs the phase is a pure power \(3/2\). After \(j\ge1\)
differences and a transform, a pure power \(v\) changes to
\((v-j)/(v-j-1)\). Starting at \(3/2\) this is \(-1\)
when \(j=1\), and lies in \((0,1)\) when \(j\ge2\).
From \(-1\) or \((0,1)\) every later value is in \((0,1)\).
The finitely many pure-power derivative tests therefore
exclude only finitely many residue characteristics.

Where the stationary second derivatives have been units,
let \(y_0=x\), \(y_1=F'(x)\), and let \(\Phi_i\) be the
\(i\)th Legendre phase. With the negative Legendre convention
\[
 (\Phi_i'\circ y_i)(x)=-y_{i-1}(x),\qquad
 y_{i+1}=\left((y_i')^{-1}\frac{d}{dx}\right)^{j_i}
                    (-y_{i-1}).
\]
For ratios of unequal multiplier units, the required
numerators specialize at reciprocal ratio zero to pure
powers. Their nonvanishing is enforced by the same finite
residue colouring as in the original multiplier comparison.
At equal residue ratios, the normalized second derivative
is \(x^{\rho-2}(1+s)\) with \(s\) proportional to \(x^{-1}\).
Keep the factors of the previous stationary derivatives
in a list \(P_i(s)\), all units on the branch in use.
Differentiation of \(x^a R(s)\prod_i P_i(s)^{e_i}\)
produces, after removing \(x^{a-1}\prod_i P_i(s)^{e_i-1}\),
\[
 (aR-sR')\prod_iP_i
   -sR\sum_i e_iP_i'\prod_{l\ne i}P_l .
\]
At a common adjacent zero away from the known factors,
the earlier numerator must have a repeated root.
The exact checks specialize at rational powers and reduce
modulo a verified prime. The tested numerators have their
full degree and are squarefree, so their discriminants are
nonzero as polynomials in the power. The cases \(s=0,\infty\)
are pure powers. The number of tested gaps and derivative
orders is finite. Each exceptional power maps to finitely
many exceptional input degrees, since the outer power
changes are nonconstant fractional linear functions.

\paragraph{Completion following several shifts}

For \(i\) arithmetic differences at a target saving \(t\),
the combined shift exponent is
\(C_i=i-(2^{i+1}-2)t\). The residual conductor exponent is
\(v_i=V-C_i\). If \(B_i\) is the original bad weight at
derivative order \(i+2\), choose the shifts first on that
part, giving bad weight at most \((B_i-C_i)_+\).
The annular bound suffices if
\[
 V\ge i+2t,\qquad 2^i t\le S_g(v_i),\qquad
 B_i\le3V-2i-(2^{i+1}+4)t .
\]
For any fixed \(i\), substituting \(B_i\le V\) gives
\(t\le(V-i)/(2^i+2)\). The curve
\(2^it=S_g(V-i+(2^{i+1}-2)t)\) is the affine image, on each
cell, of \((v,S_g(v))\) under
\[
 (v,z)\longmapsto
 \bigl(i+v-(2-2^{1-i})z,\,2^{-i}z\bigr).
\]
The positive slopes of \(S_g\) are at most \(1/2\), so the
first coordinate is strictly increasing. This shape condition
also pays for unintended conductor reductions by the usual
divisor average. The verified profile chooses among the
available numbers of shifts. More than one candidate is
allowed, since the choices act on the phase separately.

\paragraph{Taylor pairs with a critical highest derivative}\label{par:refined-critical}

At Taylor degree \(D\ge3\), let \(f=F^{[D]}\) and
\(g_1=F^{[D-1]}\), so \(g_1'=Df\). Work with
\(A=q^{1/(D+1)}m^{O(1)}\), \(Y=X/A\), and penultimate
tolerance \(\Delta=Y^{1-D}m^{O(1)}\).
A unit \(f'\) fixes the usual residue of \(x\) for a pair
\(x,x+h\) with equal top coefficient. Otherwise \(f,f''\)
are units. Set \(d=(h,A)\), \(s=A/d\), and \(c=(h,A^2)\)
after discarding the rounding factors. At a prime of \(s\)
the equation \(f(x+h)=f(x)\) has one stationary solution
\(x_*(h)\) on the disk, with unit derivative after division
of the difference by \(h\). The congruence restricts \(x\)
to that class modulo \(s\). The stationary value of
\(g_1(x+h)-g_1(x)\) is linear in \(h\) plus a cubic term
with unit coefficient, followed by terms of order at least
five. This follows by centring at \(x_*(h)+h/2\), an even
function of \(h\).

When \(d\le A^{2/5}\), complete the sum in \(z=h/d\)
at modulus \(s^2\). A nonzero output frequency \(\ell\)
gives an input frequency coset of norm \(c_\ell\), with
\(c_\ell\le d^2(\ell,s^2)\), and a complete bound
\(O(m^{O(1)}s\sqrt{c_\ell})\). Count its dual frequencies
at scale \(s^2d/X\), average the gcd in \(\ell\), and count
the progression for \(x\). The normalized count is at most
\(O(m^{O(1)}[X^{-1}+d^2/A^2+\Delta d/A])\).
For \(d>A^{2/5}\), keep the fixed value modulo \(s^2\)
after dividing by \(c\). The remaining modulus is \(d^2/c\).
Its good part uses a unit first derivative. On its critical
part \(R_b\) use the quadratic discrepancy
\(O(m^{O(1)}[\delta L+(1+L/R_b)\sqrt{R_b}])\).
The full-period term is needed when \(L>R_b\).
After counting the shifts the normalized bound is
\[
 O\!\left(m^{O(1)}
 [X^{-1}+\Delta/A+(Ad)^{-1}+A^{-1}(\sqrt\Delta+d^{-1/2})]\right).
\]
Here the small-$d$ average counts nonzero Fej\'er frequencies
in a centred ideal lattice. At volume $L$ the number divisible by a
divisor ideal of norm $b$ is $O_{\mathbb K}(L/b)$, by the norm lower
bound and packing. The factor $\sqrt{(\ell,Q)}$ in the completed
quadratic estimate can therefore be averaged by the divisor count.
For large $d$, after taking out the fixed digits, the remaining output
window has volume $O_{\mathbb K}(\Delta d^2/c)$. Unit derivatives
fix the input classes at the good primes. At the critical primes,
completion gives the stated discrepancy for $L\le R_b$. For
$L>R_b$, partition the box into $O_{\mathbb K}(1+L/R_b)$ translated
periods of the balanced conductor lattice and apply the complete
quadratic estimate to each. Since \(D\ge3\), \(A^{-1}\sqrt\Delta\ll
X^{-1}m^{O(1)}\). The split therefore gives
\(X^{-1}+A^{-6/5}+\Delta A^{-3/5}\) up to radical factors.

For a specified nonzero top difference, Cauchy's inequality
on the occupancies of the two coefficient cells costs no
more than the exact-top bound with a fixed enlargement of
the lower window. There are \(O(1+A/Y^D)\) top differences.
The critical degree-\(D\) mean value divides the resulting
pair exponent by \(D(D+1)\). Writing \(\xi=V/(D+1)\),
the saving on the positive range is
\[
 \frac{\min(1,6\xi/5,D-1-(D-1-3/5)\xi)
           -\max(0,V-D)}{D(D+1)}.
\]
The formula gives the cubic critical tent at \(D=3\).
The calculation uses the higher-degree versions only where
this numerator is positive and the length exponent is positive.

Writing \(\lambda=(D-1)(1-\xi)\), one can optimize the split in
\(d\) and also use the stationary value directly. On \(h=dz\) the
penultimate difference divided by \(d\), projected to \(s^2\), is
\(D f(a_*)z+O(d^2z^3)\) with \(D f(a_*)\) a unit. Every prescribed
value fixes one congruence class of \(z\). There are
\(O(1+\Delta d s^2)\) allowed values and
\(O(1+X/(ds^2))\) shifts per class. Multiplying by \(O(Xd/A)\)
for the root progression gives the normalized bound
\[
 O\!\left(m^{O(1)}
    [X^{-1}+d^2/A^3+\Delta A/X+\Delta d/A]\right).
\]
Use this or the Fourier estimate for \(d\le A^p\).
The complementary error is \(O(m^{O(1)}A^{-1-p/2})\)
after discarding the terms already bounded by \(X^{-1}\) or
\(\Delta A^{p-1}\).
To balance this with \(\Delta A^{p-1}\), take
\(p=\min(p_{\max},2\lambda/(3\xi))\) with \(p_{\max}=2/5\)
for the Fourier estimate and \(p_{\max}=4/5\) for the
inversion estimate. This gives two further exact-top exponents
\[
 \min(1,6\xi/5,\xi+\lambda/3),\qquad
 \min(1,7\xi/5,\xi+\lambda/3,D(1-\xi)).
\]
Apply the same top-difference correction and mean-value denominator.

The count can also fix the penultimate output \(R\bmod A^2\)
before counting shifts. Let \(G\) be the good part of \(A\) at which
\(f'\) is a unit, and \(C=A/G\) its critical part.
For each \(R\) and \(x\), the map \(h\mapsto G_h(x)\) has unit
derivative at every critical prime and fixes the shift modulo \(C^2\).
At each prime of \(s\), if \(h_R(x)\) denotes the resulting
analytic inverse, then \(h_R\) and its derivative are divisible by
the prime. The derivative of
\((f(x+h_R(x))-f(x))/h_R(x)\) is the unit \(f''(x)\) modulo that
prime. This fixes \(x\bmod s\). The good primes require \(G\mid h\).
Hence the normalized count is at most
\[
 O\!\left(m^{O(1)}
     [X^{-1}+\Delta A/X+Gd/A^3+\Delta G/A]\right).
\]
In the complementary large-\(d\) estimate, write \(c_c=(h,C^2)\)
and \(R_g=G^2/(h,G^2)\). The represented good output classes satisfy
\(n_g\ll\min(R_g,1+\Delta d^2/c)\). The factor
\(n_g/\sqrt{cR_g}\) is \(O((Gd)^{-1/2}+\sqrt\Delta)\) by first
separating the constant \(1\) and then taking the geometric mean.
Split at \(Gd=A^r\) with \(r=\min(4/3,\lambda/\xi)\). The two
bounds now give exact-top exponent
\[
 \min(1,5\xi/3,\xi+\lambda/2,D(1-\xi)).
\]
The top-difference correction and mean-value normalization again
give a piecewise affine saving in \(V\).

\paragraph{The choice of prime multiplier.}
Take $D\ge3$, put $\xi=U/(D+1)$ and $a_0=2\xi-1$, and choose the
prime multiplier norms in a dyadic interval of size $B=X^\beta$.
For unequal primes \eqref{eq:joint-dual-size} gives the sharper
upper bound $l\le a_0-\beta+O_D(\log_Xm)$ in
\eqref{eq:joint-dual-exponents}, since the progression step contains
the multiplier of norm $B$ even for the largest line volume. Thus,
apart from radical losses,
\[
 a_0-\beta-2\gamma\le l\le a_0-\beta,\qquad
 1+\beta-\gamma\le w\le1+\beta .
\]
With $Y=X/(AB)$ and $\eta=U-D+D\beta$, the count for an exact top
coefficient is $O(X^{-1-\beta+\epsilon}m^{O_D(1)})$ after
normalization, by \eqref{eq:joint-zero-top}. To obtain that same
bound for the other pairs it suffices that the combined line saving
exceed $\eta$. We require
\[
 \begin{gathered}
 \xi>\beta>0,\qquad \xi+\beta<1,\qquad
 0<2\eta<a_0-\beta,\qquad 1-\eta>4a_0/5,\\
 \eta<(D-2)(1-\xi)-D\beta .
 \end{gathered}
\]
These are the length, large-gcd, argument-range and zero-frequency
tests of the preceding completion argument. The last test pays
for the $j=0$ count in \eqref{eq:joint-zero-fejer}. The fourth
inequality keeps the dual argument in $V\ge9/5$ for $\gamma<\eta$.

Let $H(V)$ be an available mixed-phase profile on this range.
If its affine cells $\alpha V+\delta$ satisfy
$2(\alpha+\delta)_++(\alpha)_+\le1$, substituting the nominal
values $l=a_0-\beta$, $w=1+\beta$ into $lH(1+w/l)$ costs at most
$\gamma$. Indeed on one cell that expression is
$(\alpha+\delta)l+\alpha w$, and the loss is at most
$2(\alpha+\delta)_+\gamma+(\alpha)_+\gamma$. Continuity gives
the bound across cells. The small-gcd groups therefore suffice
under the strict test
\[
 \eta<(a_0-\beta)H\!\left(1+\frac{1+\beta}{a_0-\beta}\right).
\]
The large-gcd groups use the packing bound as in
\eqref{eq:joint-H}. For equal prime multipliers the divisibility
$b\mid jr$ recovers a factor $B^{-1}$, and the corresponding
test is $\eta<a_0s(1+1/a_0)$ for the single-binomial profile $s$
of \eqref{eq:joint-P}. These bounds prove the saving
\[
                 \frac{1+\beta}{D(D+1)}
\]
apart from the prescribed radical and exponent losses. At each
conductor the affine-cell restrictions and the strict tests are
linear inequalities in $\beta$. The local derivative checks
involve a fixed list of rational powers for every compact
conductor range. The discarded degrees and the exponents of the
radical losses depend on that list, independently of $\mathbb K$.

\subsection{Integrating the conductor savings}
\label{sec:joint-price}

Let \(s(U)\) be a positive, continuous and nonincreasing
profile on \(U\ge1\), made from the estimates of
\cref{sec:joint-conductors,sec:joint-completion}. At a step
with \(M=U\), the shortening may approach \(L=U\) and
\(2z=s(U)\). For every fixed compact interval one first
chooses strict margins in these two replacements. The
radical cost \(O(M/k)\) and the exponent losses then fit
inside the margins for large \(k\). With \(d\) the excess
divided by \(k\), the decrement per unit of normalized
order approaches \(\min(Ud,s(U))/2\).
At \(d=s(U)/U\) the two terms agree, and monotonicity
of \(s\) proves that this choice of \(U\) maximizes the
minimum.

Let \(u_0+s(u_0)=2\), \(t_0=s(u_0)\), and \(d_0=t_0/u_0\).
For the profiles below \(u_0\) is unique, because the
slopes of \(s\) are at least \(-1/2\). Start the iteration
at even order \(r_0\) with \(r_0/k=1-d_0-\log d_0+O(1/k)\),
using the smooth seed of the preceding moment construction. On \([u_0,R]\)
the variable \(s(U)/U\) decreases strictly. Use a finite
mesh of moduli, take \(L\) and \(2z\) strictly below the
profile values, and compare the decreasing excess to
the integral of the inverse decrement. Each step
ending above the stopping level decreases that integral
by at least \(2/k\). The finite mesh can approximate
the integral arbitrarily closely on this compact
positive range. With \(\mathrm d(s(U)/U)
=s'(U)\,\mathrm dU/U-s(U)\,\mathrm dU/U^2\), its limiting
value is
\[
 \int_{s(R)/R}^{d_0}\frac{2\,\mathrm dd}{s(U(d))}
 =\frac{2}{u_0}-\frac{2}{R}
       +2\int_{u_0}^{R}\frac{-s'(U)}{U s(U)}\,\mathrm dU .
\]
The seed and \(2/u_0\) combine to give
\(2+\log(u_0/t_0)\).

For the continuation take \(R\) to be a fixed integer
at least \(300\) beyond the finite part of the
calculation, so \(s(R)=1/[R(R+1)]\).
Use the growing-conductor degree-reduction test
of Lemma~\ref{lem:nf-fibre-savings} with \(z=1/(14M^2)\),
\(L=0.99999M\) and \(R\le M\le H=\sqrt{k}/8\).
Indeed \(j=\lceil3U/2\rceil\) in the three-term
degree-reduction inequalities is admissible for
\(R\le U\le M\), since \(M^2/k\le1/64\).
For example the three brackets before division by
\(j(j-1)\) are at least
\[
 (3/2+2/R)^{-1}-(3/2+1/R)/64,\quad
 1/3-H/k,\quad R/2-(3/2+1/R)/64 .
\]
Each exceeds \(2(3/2)(3/2+1/R)/14\) for large \(k\).
At \(U\le R\) the finite-profile estimates still apply
with the radical bound \(H/k\), tending to zero.
The bridge cost is \(O(1/R)\) in normalized order,
and the speed in the growing range is bounded below
by a constant times \(d^{2/3}\). Its integrated cost
from the bridge to the cap is therefore \(O(1/R)\).
The fixed-modulus step at \(H\), followed by the final
absorption, costs \(o(1)\) in normalized order.
Take \(b=\lceil16\sqrt{k}\rceil\) and stop below
\(d_*=b/(2k^3)\). Then \(kd_*<b/[k(k-1)]\).
Theorem~\ref{thm:nf-absorb-transfer} gives the loss-free
even moment after \(b+O(1)\) additional factors.
The cap takes \(O((\log k)/\sqrt{k})\) in normalized
order.
Let \(k\) tend to infinity first and then \(R\) to
infinity. Thus an upper bound for the limiting normalized
sharp moment order, including one extra root if desired, is
\begin{equation}\label{eq:joint-price}
 \mathcal C(s)=2+\log(u_0/t_0)
       +2\int_{u_0}^{\infty}\frac{-s'(U)}{U s(U)}\,\mathrm dU .
\end{equation}
The seed order is above the minimum required for the
shortening, and the final order in this calculation is
above \(4k+4\) once \(k\) is large.

\paragraph{The finite profile and its tail.}
The base profile is Wooley's two-coefficient profile
\cite[(A.5)]{Wnew} through Taylor degree \(16\),
including its degree-\(16\) descending line to zero
at \(17\). It follows from \eqref{eq:joint-taylor-saving},
completion and the first three \(A\) transforms.
Adjoin \eqref{eq:joint-prime-saving} for \(D\le16\)
and \eqref{eq:joint-H} for \(3\le D\le16\).
For this initial use of \eqref{eq:joint-H}, start
its descending line at \(U=D\) and join the Wooley line
at \(U_+\). The chord lies below \eqref{eq:joint-H}
on that interval. Taking the running minimum of
the envelope ensures the stated saving at every
smaller conductor. Make three rounds in which we
close the current profile under \(A,B\) to depth
sixteen, use this closure for the input of \(P_D[s]\)
in \eqref{eq:joint-P} for \(3\le D\le15\), and
adjoin the outputs to the profile preceding that
round. Then close under \(A,B\) to depth twenty-eight.
Adjoin the Wooley estimates through \(D=32768\).
Use this last envelope once more as input to
\(P_D[s]\) for \(D<32768\), and adjoin these outputs,
\eqref{eq:joint-prime-saving} through \(32768\),
and the full form of \eqref{eq:joint-H} for \(3\le D\le512\).
Take the running minimum once more. All the curves
used as inputs satisfy \eqref{eq:joint-profile-class}
and \(s(2)<1/2\). There are no circular uses of
\eqref{eq:joint-P}.

The parameters for \eqref{eq:joint-H} are given in
the electronic data. For example the first two
\((\beta,\eta)\) are
\((577/100000,1733/100000)\) and
\((9/1000,3603/100000)\).
All \(510\) pairs satisfy \eqref{eq:joint-parameter-tests}
with positive rational margins, as well as the left
splice condition after \eqref{eq:joint-H}. The least
margin in the tested inequalities is \(1/100000\).
At the final splice \(U=32768\) the ordinate is
\(1/(32768\cdot32769)\). On each subsequent interval
\([D,D+1]\) use the larger of
\[
       \frac{D+1-U}{D(D+1)},\qquad
       \frac{2(D+1)-U}{(D+1)^2(D+2)} .
\]
These two lines are \eqref{eq:joint-taylor-saving}
and \eqref{eq:joint-prime-saving} at the next Taylor
degree. Their contribution to the integral in
\eqref{eq:joint-price}, including the factor \(2\), is
\begin{equation}\label{eq:joint-tail-term}
 \frac{1}{D+1}\log
       \frac{(D^2+D+2)(D+1)(D+2)}{D^4}.
\end{equation}
They meet at \(U=(D+1)(D^2+D+2)/(D^2+2D+2)\).

The certificate takes maxima by intersecting rational
affine pieces. On a falling line \(s(U)=aU+b\) in
\([x,y]\) the contribution is
\[
        \frac{-2a}{b}\log\frac{s(x)y}{s(y)x}.
\]
It extracts powers of two before evaluating a logarithm.
For \(1\le z\le2\), put \(w=(z-1)/(z+1)\). The
positive series for \(\log z\) is
\(2\sum_{j\ge0} w^{2j+1}/(2j+1)\) and the first
omitted term divided by \(1-w^2\) is an upper bound
for the remainder. Each weighted term is rounded
outwards at scale \(10^{-28}\). The terms of
\eqref{eq:joint-tail-term} are summed through \(D=49999\).
For the omitted sum beginning at \(N=50000\),
compare with Wooley's
\[
 L_D=\frac{2}{D+1}\log\left(1+\frac{2(D+1)}{D^2}\right).
\]
The derivative of \(\log(1+2x+2x^2)\) gives
\(0\le 4/[D(D+1)]-L_D\le 8/[3D^3(D+1)]\).
The difference between \(L_D\) and \eqref{eq:joint-tail-term}
is
\[
  \frac{1}{D+1}\log\left(1+
       \frac{D^2}{(D^2+D+2)(D+1)(D+2)}\right).
\]
For \(D\ge16\) it lies between \(D^{-3}-9D^{-4}\)
and \(D^{-3}\): factor the denominator using \(1/D\)
and use \(x-x^2/2\le\log(1+x)\le x\).
Telescoping and integral comparison give tail bounds
\[
 \frac4N-\frac{8}{9(N-1)^3}-\frac{1}{2(N-1)^2}
 \ \le\sum_{D\ge N}\text{the term in \eqref{eq:joint-tail-term}}\ \le\
 \frac4N-\frac{1}{2N^2}+\frac{3}{(N-1)^3}.
\]
The resulting seed is \(u_0=141/79\), \(t_0=17/79\).
Of \(161498\) exact finite cells, \(110585\) have
falling contributions above the seed. The resulting
enclosure for \(\mathcal C(s)\) is
\begin{equation}\label{eq:joint-enclosure}
 \bigl[6.2195722132177257706364086906,\
       6.2195722132177648838542803201\bigr].
\end{equation}
This profile gives an upper price of \(6.219572214\).
\paragraph{The refined price.}
Apply the estimates of Section~\ref{sec:refined-conductors} to the
same shortening argument. On each compact conductor interval
all norms enter through the displayed scalar inequalities.
Stationary calculations take place at degree-one primes and
the Fourier boxes are on the dual trace lattice. Thus the
conductor exponents and the finite list of excluded degrees are
independent of the field.

The calculation uses two available prefixes, with Taylor cutoff
$16$ in the seed and local two-multiplier tests up to degree $384$.
Start with Wooley's profile and the direct monomial estimates,
adjoin the admissible prime-multiplier cells, and take a running
minimum to cover every smaller conductor. At each closure round
use arithmetic $A$ and $B$ and the preceding profile as input
for the one- and two-multiplier operations. The first prefix uses
two such rounds at depth ten, followed by arithmetic closure
to depth sixteen. It then makes three further rounds restricted
to the compact interval ending at $10$. In each of these, the
input is cut off by the affine cap $(10-U)/18$, followed by six
arithmetic closure steps, one multiplier up to degree nine,
and the refined two-multiplier estimate. The second prefix uses
four rounds at depth twelve followed by closure to depth twenty.
It lowers the affine pieces to a rational grid. To do this,
join adjacent endpoints and subtract the maximum excess of the
joining line over the previous profile at its knots. At a
shared endpoint take the larger of the two required reductions.
Apply the shape bounds for the next transformation by taking
the corresponding lower envelope. All parameter and shape
tests are applied to rational affine cells. The inputs to a
round have already been established, so every resulting cell is
an available saving.

Take the pointwise maximum of the two prefixes. To continue with
the one-multiplier bound, use
\[
 s(V)\ge\gamma_s-\nu_s(V-2),\qquad
 (\gamma_s,\nu_s)=(863/5000,3/20).
\]
The endpoint checks give this inequality for
$2\le V\le2+4/(R-3)$ at the splice $R=385$. For $D\ge R$
and $x=U-D\in[0,1/4]$, put
$A_D=\gamma_s(D-1)-2\nu_s$ and
$B_D=D+1-2(\gamma_s+\nu_s)$. In \eqref{eq:joint-P},
the lower bound for the single-multiplier line saving is
$(A_D+2(\gamma_s+\nu_s)x)/(D+1)$, below the
zero-frequency constraint. The resulting order-$D$ profile has
plateau $1/[D(D+1)]$ up to $x=A_D/B_D$, then falls until it meets
the next-degree estimate at
$x=(1+A_D)/(B_D-D/(D+2))<1/4$.

The annular $B$ estimate for a degree-four critical Taylor
calculation, followed by $A$, also gives the mixed profile
\[
 H(V)\ge\gamma_o-\nu_o(V-2),\qquad
 (\gamma_o,\nu_o)=(43/250,7/40),\qquad
 2\le V\le41/20 .
\]
With $a_0=2U/(D+1)-1$, the mixed constraint on the multiplier has
the limiting value
\[
 \beta_*(U)=
 \frac{(\gamma_o+\nu_o)a_0-\nu_o-x}{D+\gamma_o+2\nu_o}.
\]
Use the saving $(1+\beta_*(U))/[D(D+1)]$ for $U\ge D-1$
where $\beta_*(U)\ge0$. For $D\ge384$ on this interval,
$a_0\ge381/385$ and $\beta_*\le(1+\gamma_o)/384$. Thus
the dual argument $1+(1+\beta_*)/(a_0-\beta_*)$ lies in the
required range and $0<x+D\beta_*\le\gamma_o$. The single-multiplier
affine input has larger height at $V=2$ and smaller downward slope.
It gives a larger line saving since its length exponent is $a_0$
and its argument is $1+1/a_0$. The mixed bound satisfies
$a_0-\beta_*>2(x+D\beta_*)$ and $1-(x+D\beta_*)>4a_0/5$.
The zero-frequency condition follows from
$2(x+D\beta_*)<1-3(1-x)/(D+1)$. The remaining length conditions
follow from $\beta_*\le(\gamma_o-x)/D$. Taking $\beta$ just below
$\beta_*$ gives the required strict conductor inequality.

On $[D,D+1]$ take the mixed line up to the plateau, then the
plateau and the falling one-multiplier line until it meets the
next-degree mixed line. That intersection lies before $D+1/4$.
On $[R-1,R]$ the prefix is joined by taking its maximum with the
order-$R$ mixed line. All junctions in the resulting profile
are continuous. Its price is computed by the affine integration
used above. From an integer $N>R$, the ratio $r(U)$ of the new
tail to the original two-line tail satisfies
$0\le\log r(U)\le(2+\gamma_o)/D$ on $[D,D+1]$ and
$1\le r(N)\le1+\gamma_o/N$. The change of the omitted integral is
\[
 \frac{2\log r(N)}{N}
       -2\int_N^\infty \frac{\log r(U)}{U^2}\,dU .
\]
It can be enclosed by decreasing the old lower bound by
$(2+\gamma_o)/(N-1)^2$ and increasing the old upper bound by
$2\gamma_o/N^2$.

The seed satisfies
$u_0=7577850000002/4245950000001$ and $t_0=2-u_0$.
The two prefixes and the splice are recorded by $43{,}453$ rational
affine cells, of which $30{,}171$ have falling contributions above
the seed in the price calculation. Rational enclosures for the
logarithms through $N=160000$, with the stated omitted-tail estimate,
give
\begin{equation}\label{eq:refined-enclosure}
 \mathcal C(s)\in
 [6.2162163018144242131724906863,\,
  6.2162163019127077173258391677].
\end{equation}
The calculation can be replayed both by reconstructing the
profiles and by independently integrating the saved affine
cells. The data and the two checks are in
Appendix~\ref{app:certificates}.

\begin{theorem}\label{thm:nf-optimized-moment}
Let \(C_1=6.216216302\).
For all sufficiently large \(k\), independently of the number field,
the marked sums have a loss-free even moment at an order
\(S\) with \(S+1\le C_1 k-\sqrt{k}\).
The mixed moment used in Theorem~\ref{thm:nf-absorb-transfer}
can be taken with an ordinary allocation
\(\lceil16\sqrt{k}\rceil+O(1)\).
In particular,
\[
 G_{\mathbb K}(k)\le
    \max\{H_{\mathbb K}(k),\lfloor C_1 k-\sqrt{k}\rfloor\}.
\]
\end{theorem}
\begin{proof}
By \eqref{eq:refined-enclosure} there is a fixed positive gap
between the price and \(C_1\). At a large fixed integer \(T\) in
the refined tail the profile has height \(\asymp T^{-2}\).
It can be lowered to \(1/[T(T+1)]\) there, and the same bridge
estimate costs \(O(1/T)\) in normalized order. Choose \(T\) and
the strict margins on the compact interval so these costs fit
inside the price gap. Then take \(k\) large enough
for the bounded list of changes of local power, the radical
losses at all scheduled moduli and the cap allowances.
None of these conditions involves the field degree.
The choice of windows and the height threshold can depend on
the field, as in the proof of the shortening step.
The fixed price gap absorbs \(\sqrt{k}+1\) in the asserted
order. Theorem~\ref{thm:nf-absorb-transfer} yields the
representation.
\end{proof}

\subsection{A Weyl estimate for the marked sum}
\label{sec:nf-marked-weyl}
The unrestricted final cofactor can also supply the pointwise saving.
This gives a representation result in which every root retains the
prescribed prime factors.

\begin{lemma}\label{lem:nf-bilinear-weyl}
Suppose $D Y=X$ and $D,Y\ge X^c$ for some fixed $c>0$. Let
$a_d,b_y$ be coefficients of absolute value at most one on translated
integral-basis boxes of side lengths $O(D^{1/n})$ and $O(Y^{1/n})$,
respectively. For $0<\xi<1$ one has
\[
 \sup_{\alpha\notin\mathfrak M(X^\xi)}
 \left|\sum_{d,y}a_db_y e_{\mathbb K}(\alpha d^k y^k)\right|
              \ll_{\mathbb K,k,c,\xi}X^{1-\epsilon_0}
\]
for some $\epsilon_0>0$.
\end{lemma}
\begin{proof}
Extend the two coefficient sequences by zero to the full boxes, and apply Cauchy's inequality in the two sets of variables successively.
After division by a constant multiple of $X^4$ the fourth power is
bounded by the box average of the character of
$\alpha(d_1^k-d_2^k)(y_1^k-y_2^k)$, with all coefficients removed.
For instance, Cauchy's inequality first in $d$, and then in the pair
$(y_1,y_2)$ after expansion, gives this assertion. Fix the two
initial shifts. Difference $k-1$ further times in $d$ and $k-2$
times in $y$, using Cauchy's inequality on the remaining averages
at every step. Zero shifts cost $O_k(D^{-1}+Y^{-1})$ after taking
the common normalized power $2^{2k-1}$. The final phase is linear
in $y$, with coefficient $\alpha t$, where $t$ is $(k!)^2$ times
the product of $k$ shifts from the first box and $k-1$ from the
second. The remaining variable $d$ has a trivial count $O(D)$.
Every nonzero $t$ has $O(X^\varepsilon)$ such factorizations,
by the ideal divisor bound and the unit count of
Lemma~\ref{lem:nf-fibre-savings}. Its height is $O(T^{1/n})$
where $T=D^kY^{k-1}=X^k/Y$.

Intersections of translated integral-basis boxes remain boxes.
Consequently the standard product of one-dimensional linear-sum
bounds applies to the remaining $y$-sum, uniformly in its end
points. Partition the dual torus into boxes of side $\asymp Y^{-1/n}$
in integral dual coordinates. The supremum of the absolute linear
sum on these boxes, summed over all boxes, is
$O(Y(\log(2Y))^n)$. The maximum occupancy among frequencies
$\alpha t$ is bounded, by taking differences, by
\[
 M=\#\{(t,b)\in\mathcal O\times\mathfrak D^{-1}:
  H_\infty(t)\ll T^{1/n},\ 
  H_\infty(\alpha t-b)\ll Y^{-1/n}\}.
\]
The normalized differencing bound is therefore
$O(X^\varepsilon(D^{-1}+Y^{-1}+M/T))$.
If two counted pairs are independent over $\mathbb K$, scale the
coordinates of $(t,\alpha t-b)$ by $T^{-1/n},Y^{1/n}$.
The resulting lattice has covolume $\asymp Y/T$ and contains
$2n$ independent vectors of bounded length, by multiplying the
original independent pairs by an integral basis. Its covering
radius is bounded, so disjoint Voronoi cells give $M\ll T/Y$.
Otherwise, apart from the origin, $b=\gamma t$ for one
$\gamma\in\mathbb K$. A pair with $t=0$ has $b=0$ for large $X$.
The possible $t$ then belong to the conductor ideal of $\gamma$
and satisfy
\[
 |t_v|\ll
 \frac{T^{1/n}}{1+(T Y)^{1/n}|\alpha_v-\gamma_v|}.
\]
Ideal packing bounds $M$ by $O(1+T/(Q_\gamma
B(P^k(\alpha-\gamma))))$, where $P=X^{1/n}$.
This gives $M/T\ll X^{-\min(c,\xi)}$ outside the stated core.
Absorb the divisor and logarithmic losses and take the required root.
\end{proof}

For applications one may assign the smooth weight $w(D'z/P)$ to
each term of \eqref{eq:nf-terminal-form}. Here $w$ is supported in
any fixed small positive box away from every zero coordinate.
The generators $D'$ are balanced and have unique marked
factorizations. On a dyadic block $|ND'|\asymp D$ the free cofactors
lie in a box of volume $O(Y)$ with $Y=X/D$. Both $D$ and $Y$ are
bounded below by fixed powers of $X$. On fixed larger boxes, insert
smooth cutoffs and expand $w((d/D^{1/n})(y/Y^{1/n}))$ into its
Fourier series in the two real coordinate vectors. Its coefficients
are absolutely summable, uniformly on all the dyadic blocks,
and every term separates the two sets of variables. All prime
conditions occur in the first variable. The $O(\log X)$ blocks
and Lemma~\ref{lem:nf-bilinear-weyl} show that the resulting $F_w$
itself has a power saving outside $\mathfrak M(X^\xi)$. These
weights are bounded nonnegative restrictions of the same cofactor
construction. Indeed the normalized products
$D'/|ND'|^{1/n}$ lie in a fixed compact subset of $V^\times$.
Choose the original free annular envelope positive on every inverse
image of the physical support under this compact family, enlarging
its fixed box if necessary. A fixed multiplicative constant then
dominates the physical weights and the bounded differentiated
weights on the enlarged boxes. At a partially removed state the prescribed factors
continue to enter the terminal weight, so the fibre identities and
the sharp $S$th moment still apply. In Poisson summation on the core,
the variable substitution by $D'/|ND'|^{1/n}$ has real determinant
of absolute value one. Thus the archimedean factor for $F_w/X$ is
the ordinary one times $c(X)=\sum_{D'}|ND'|^{-1}\asymp1$.

It follows that for $s\ge\max\{S+1,H_{\mathbb K}(k)\}$ all $s$
roots in the representation formula of
Theorem~\ref{thm:nf-absorb-transfer} can carry the marked weights,
with a power-saving error. Indeed, off the core use the supremum
estimate for one factor and H\"older's inequality on $S$ factors at
their respective sharp $S$th moments. Any further factors have
the trivial bound. Choose the core exponent below the smallest prime-window exponent, so every mark is a unit at every core denominator. On the core the $s$ factors have the ordinary
complete sums and the required archimedean weights. Choosing each
weight positive near the corresponding coordinate of a nonsingular
archimedean solution gives the uniform positive main term.
The same argument and Bessel's inequality give a power-saving
exception count when $s\ge S/2+1$ on collections of targets
where the $s$-variable singular series has a uniform positive
lower bound and $s>2k$.

\subsection{A pointwise bound for the marked sum and a uniform coefficient}
\label{sec:nf-self-weyl}
A preliminary linear moment at partially removed states supplies a
stronger pointwise estimate for the same family of marked primes. This removes the
finite-degree cost of using the ordinary Weyl exponent over a field.

\begin{theorem}\label{thm:nf-all-linear}
For every number field $\mathbb K$ and every $k\ge3$,
\[
 G_{\mathbb K}(k)\le\max\{H_{\mathbb K}(k),\lfloor6.873k\rfloor\}.
\]
For $k\ge11$ the loss-free even marked moment may be taken at an
order $S$ with $S+1<6.873k$. For all sufficiently large \(k\) it may be taken with
\begin{equation}\label{eq:nf-secondary}
                            S+1\le C_1 k-\sqrt k .
\end{equation}
The thresholds for $k$ are independent of $\mathbb K$.
\end{theorem}
\begin{proof}
We first obtain an even $S_1<7.42k$ with
$\int|F|^{S_1}\ll X^{S_1-k}$ using the ordinary estimate. Here
and below a restricted sum with coefficients bounded in modulus
by a fixed marked positive sum inherits its even moment bounds by
orthogonality. We record the scalar computations so that the fixed
constants and rounding at small degrees can be checked.

There are two ways to convert a restricted $r$th moment of excess
$D$ into the preliminary loss-free moment. The pointwise ordinary
bound uses an integer $b>k(k-1)D$ and $S_1\ge r+b$. Alternatively
put $w_c=k(k+1)$, $b=\lfloor w_cD/(1+D)\rfloor+1\le w_c$, and
take $S_1-b\ge r(1-b/w_c)$. In either case use an even $S_1>2k$.
These are precisely the two strict compensations of
Theorem~\ref{thm:nf-absorb-transfer}, the second using the critical
broad-minor estimate \eqref{eq:aud-field-critical}. Choose the
original-scale losses inside the strict inequalities, including
their use below at partial volumes. The construction of the orders
at the end of this proof gives $S_1/k\le2560/359$ for
$11\le k\le2500$, $S_1/k<7.060626$ for $2501\le k<10^{10}$ and
$S_1/k<6.876509$ thereafter. In particular $S_1<7.42k$ throughout.

We next arrange the marks in a form which gives
\begin{equation}\label{eq:nf-self-sup}
       \sup_{\mathfrak m_*}|F|\ll X^{1-1/(15k)-\epsilon_1}
                     \qquad(k\ge11)
\end{equation}
for some $\epsilon_1>0$ depending on $k$. Choose a group of the
tiny primary windows whose total log-norm exponents lie strictly
between $0.498$ and $0.499$. This is possible by first making
the individual exponents sufficiently small, then fixing narrow
disjoint open intervals around them. The complementary primary
windows and the auxiliary windows are those of the construction
in Section~\ref{sec:nf-fibres}. The exponents can have arbitrarily
small positive widths. In particular the sum of all primary
exponents can be confined to an arbitrarily short subinterval
of $(1-2\gamma,1-1.99\gamma)$.

Here are details on separation of the inner sum. Partition every
primary window into dyadic norm intervals, retaining full intervals
inside the window. In each resulting rectangle, retain auxiliary
dyadic intervals contained in the prescribed auxiliary window
for every choice of the primary primes in that rectangle.
The ratios of all extremal terminal scales in the rectangle are
bounded by constants depending only on the number of marks.
At most $O_k(1)$ end bins are lost per window, of reciprocal mass
$O_{\mathbb K,k}(1/\log X)$ by the prime ideal theorem; the full
window has reciprocal mass bounded below. Auxiliary factors can be
ordered as before. Thus the retained sum still has mass comparable
to $X$. Use the physical smooth weight $w(D'z/P)$ from
Section~\ref{sec:nf-marked-weyl}. On the resulting
$O_{\mathbb K,k}((\log X)^C)$ rectangles the outer group has
product $v$ of balanced generators, $|Nv|\asymp M$ with
$X^{0.498}<M<X^{0.499}$ up to fixed factors. The remaining product,
including the free cofactor, is $u$ at volume $U=X/M$.
After smooth separation of $w(vu/P)$ on enlarged boxes, each
term is a sum over independent $v$ and $u$ with bounded coefficients.
The Fourier coefficients of this separation are uniformly absolutely
summable. In the free variable one may use an enlarged fixed box,
since its admissible physical range on a dyadic rectangle varies
by only a bounded factor.

The preliminary order is needed uniformly in every retained
rectangle. Fix one admissible outer tuple $v_*$ of that rectangle.
By its definition every retained inner tuple satisfies the auxiliary
interval conditions for this particular outer choice. It is therefore
in the partial state obtained by removing exactly the primary marks of
$v_*$ from the original family, with a fixed enlargement of the free
box. The actual state volume $Y=X/|Nv_*|\asymp U$ exceeds $X^{1/2}$.
The uniform descent \eqref{eq:nf-short} gives the same scheduled
$r$th moment here, with $X^{c_0}\le Y^{2c_0}$. The preliminary tests
were made with room for $2c_0$. Apply the minor estimate and the
broad major-arc H\"older absorption from
Theorem~\ref{thm:nf-absorb-transfer} at this $Y$ to get
$O(U^{S_1-k})$. Only a moment upper bound is used, so no lower
bound for the partial state's mass is required. The constants are
uniform in $v_*$, in the rectangles and in the bounded dilation.
The same even moment bounds hold under multiplication of the inner
weights by any bounded function, including the scaled derivatives
introduced in sampling, by orthogonality and positive coefficient
domination. At the full volume the retained roots are a positive
subfamily of the original family and have the mass already proved.
This supplies the preliminary, final and inner moments for one
fixed family of marks.

We give the sampling calculation. Work in one of the finitely
many angular sectors used in \eqref{eq:edlip}. The close pairs
$(d,c)$ for $d=v^k-v_0^k$ at precision $U^{-k/n}$ have rank at
most one over $\mathbb K$: the norm of their determinant is
$O((M/U)^k)=o(1)$. In a sampling cell with at least two points
let $\gamma$ be the resulting slope, $\mathfrak q$ its conductor.
Its norm is $O(M^k)$, since a nonzero $d$ annihilates the
denominator. Group the outer products by the product $d_0$ of
their marked factors dividing $\mathfrak q$, and put $D=|Nd_0|$.
As in \eqref{eq:edsmallasym}, there are $O_k(1)$ such products
for this slope, and the remaining factors lie in a union of
$O(X^\varepsilon)$ unit-root classes modulo
$\mathfrak q/\gcd(\mathfrak q,(d_0)^k)$. For clarity, the same
grouping works simultaneously across the cells: if two of the
global close pairs were independent, the determinant bound just
given would be contradicted. Cells with a single point need no
grouping. The occupancy per group is therefore
\[
  O_\varepsilon\!\left(X^\varepsilon
       \left(1+\frac{M D^{k-1}}
          {N\mathfrak q\, B(P^k(\alpha-\gamma))}\right)\right)
       \ll X^\varepsilon D^{k-1}
           \qquad(\alpha\in\mathfrak m_*).
\]
Here the sector bound controls differences of the archimedean
coordinates, and there are $O(M/D)$ outer choices overall in a
group. Apply the cell large sieve at the inner even order $S_1$
as in Proposition~\ref{prop:edbilinear}. Its diagonal size is
$U^k\int|F_U|^{S_1}\ll U^{S_1}$, also for the differentiated
weights used to control a cell supremum. After H\"older the
$S_1$th power is
$\ll X^{S_1+\varepsilon}M^{-1}D^{k-S_1}$ in each group.
Singleton contributions satisfy the same bound before grouping
with $D=1$. Since $S_1>2k$ we may discard the factor in $D$.
The logarithmic number of rectangles is harmless. Finally
$0.498/7.42>1/15$, proving \eqref{eq:nf-self-sup}.

For the final even order $S$ replace just one factor of
$\int |F|^S$ by a dominating ordinary sum $f$. On $\mathfrak m_*$
bound $|f|$ trivially and use \eqref{eq:nf-self-sup} on $S-r-1$
marked factors. On the broad major arcs use H\"older with the
ordinary $S$th moment and $\int|F|^S$. A positive saving, followed
by the usual absorption, results whenever
\[
                         S-r-1>15kD .
\]
The shell can now be pruned with this same ordinary factor.
Alternatively the marked pointwise estimate on a narrow polynomial
core, Lemma~\ref{lem:nf-bilinear-weyl}, gives the all-marked
conclusion at order $S+1$. The logarithmically many rectangles
above affect neither argument.

We detail the integer constructions for both orders. For
$11\le k\le2500$ work at $r=2,4,\ldots,8k+10$ and store an upper
bound for the absolute defect as $d/N$, $N=10^9$, initially $k-1$.
Keep a separate efficient-differencing sequence. At its next step
subtract from its previous numerator $d_{\rm ED}$ only
\[
 \left\lfloor
 \frac{2d_{\rm ED}N(2^u-1)}{(kN+d_{\rm ED})2^u}
 \right\rfloor,\qquad u=\min(k,32).
\]
Indeed $(k-D)/(2k)\le1/2$ in \eqref{eq:eddefect}, and the upper
update $x-2x(1-2^{-u})/(k+x)$ is increasing for $x\ge0$. Once
$r\ge k$, an admissible pair $(z,M)$ also shortens the previous
combined numerator by
$\min(\lfloor2zN\rfloor,\lfloor0.99999Md/k\rfloor)$.
Keep the smallest valid upper bound at the new order. Include the
$M=1$ branch, for which there is no high conductor. For $k<21$ the
elementary intervals suffice. At each eligible order test the two
preliminary criteria above and the final criterion
$S-r-1>15kD$, rounding all orders up to even integers.
Appendix~\ref{app:certificates} gives the complete integer
construction. Over these $2490$ degrees it gives
\[
       \max S_1/k=2560/359<7.131,\qquad
       \max (S+1)/k=16813/2460<6.835 .
\]

For $2501\le k<10^{10}$ form consecutive bands $[K,J]$ with
$J=\min(10^{10}-1,K+\lfloor K/25\rfloor)$, starting the next band
at $J+1$. At each mesh height use a verified modulus $M$ with
$\theta=M/K$, derivative degrees $j\le\min(K,602)$, and only pair
words whose reverse exponents are nondegenerate for every $k\ge K$.
All such choices work across the band. Retain the factor
$\rho=99999/100000$ in each selected block. The initial plateau and
intermediate charges are those derived after \eqref{eq:nf-short}.
For $T=J(J-1)$ in the preliminary calculation and $T=15J$ in the
final calculation, one may close at a last plateau $(z,M)$ by adding
\[
          \frac{2}{\rho M}\log^+(zT)+\frac{12}{K}.
\]
Indeed allocate $\lceil2k/(\rho M)\rceil$ saving factors, and in
the final calculation one additional ordinary factor. The allocation
cancels the negative terminal term $-2/(\rho M)$ in the integral.
For strict stopping move the endpoint an arbitrarily small distance
if needed. The allowance $12/K$ covers even rounding at both ends,
the extra discrete step, approximation of the continuous choices,
the ordinary factor and the extra root. When
$H=9\lfloor\sqrt K\rfloor/20>300$ the program also tests a
continuous continuation. It separately checks $1\le U\le300$ at
the common $\theta=H/K$ and starting height $1/(14\cdot300^2)$.
Beyond $300$ the three inequalities below apply since
$H^2/K\le(9/20)^2$ and $H/K\le (9/20)^2/300$. Each adjacent
mesh logarithm is rounded up using
$\log(1+t)\le t-t^2/2+t^3/3$ for $t\ge0$.
All $388$ bands pass. The largest preliminary enclosure is less
than $7.060625743852722$ on $[2501,2601]$ and the largest final
enclosure, including the extra root, is at most $6.872599145421300$
on $[9446120024,9823964824]$.

Finally for $k\ge10^{10}$ use the mesh to $300$ at the common
$\theta=0.000005$. Its infinite-tail price is less than
$6.872310505217189$. With $H=\sqrt k/8$ the cap calculation of
Corollary~\ref{cor:nf-explicit-large-linear}, including $16/\sqrt k+
16/k$, gives $S_1/k<6.876508760618014$.
With $H=9\sqrt{k}/20$ the cap correction calculated below is
negative, and $12/k$ covers the roundings and other finite charges.
Thus $(S+1)/k<6.872310506417189$. These estimates prove the
asserted bounds from $k=11$ on. For $3\le k\le10$ the stated
global bound also follows directly from the preliminary transfer:
the same integer computation for the individual degrees $3,\ldots,10$
gives $10,14,22,28,34,40,46,52$, respectively. The sharper small-power
results appear in Section~\ref{sec:nf-small-powers} and in the cubic
theorems. The computations used in this proof,
including their directed roundings, are specified in
Appendix~\ref{app:certificates}.

For \eqref{eq:nf-secondary}, use the moment of
Theorem~\ref{thm:nf-optimized-moment}. Its fixed price margin
gives the asserted secondary term.

\end{proof}

\subsection{Effective orders from the two-coefficient estimate}
\label{sec:joint-finite}

The finite two-coefficient estimate
\eqref{eq:joint-finite-terms} permits an effective choice of
moments over a wide range of degrees. It requires no
quantitative radical estimate for \eqref{eq:joint-H}.
We give the parameters and the error allowances, since the
limit in \eqref{eq:joint-enclosure} by itself supplies no
effective cutoff for the improved asymptotic coefficient.

For \(k\ge k_0\) use a cap \(H=h\sqrt{k}\) and a splice \(R\).
For \(R\le U\le H\) apply
\eqref{eq:joint-finite-terms} at \(D=\lceil U\rceil\)
with \(\theta=H/k\). The first three numerators there
are bounded below by
\[
 2-\frac{4}{D+1}-(3D-1)H/k .
\]
As a function of \(D\) this is concave on \([R,H+1]\).
The two endpoint lower bounds are
\(2-4/(R+1)-(3R-1)H/k\) and
\(2-4/(H+2)-3h^2-2H/k\).
The fourth numerator is at least \(1-h^2-H/k\).
Each endpoint bound increases with \(k\) at fixed \(h\).
In all the parameter choices below their minimum with \(1\) is
strictly larger than \(c(R+1)(R+2)/R^2\). Since
\(D(D+1)\le U^2(R+1)(R+2)/R^2\), this gives saving
\(c/U^2\) throughout the growing-conductor range, with strict
slack. The finite tests also ensure coverage of \(1\le U\le R\)
at height \(c/R^2\) and radical exponent \(H/k\).

\begin{theorem}\label{thm:nf-two-coefficient-finite}
Let \((k_0,B)\) be any row of
\cref{tab:nf-two-coefficient-finite}.
For every \(k\ge k_0\) and every number field there
is a loss-free even marked order \(S\) with \(S+1<Bk\).
In particular
\[
  G_{\mathbb K}(k)
          \le \max\{H_{\mathbb K}(k),\lfloor Bk\rfloor\}.
\]
\end{theorem}
\begin{table}[htbp]
\centering\small
\begin{tabular}{rrrrr}
\toprule
\(k_0\)&\(B\)&\(R\)& \(h\)&\(c\)\\
\midrule
120 & 6.830 & 6 & \(11/20\) & \(793/2500\)\\
125 & 6.823 & 6 & \(11/20\) & \(1617/5000\)\\
150 & 6.762 & 7 & \(59/100\) & \(181/625\)\\
200 & 6.655 & 8 & \(59/100\) & \(3451/10000\)\\
250 & 6.611 & 9 & \(59/100\) & \(1943/5000\)\\
500 & 6.528 & 10 & \(61/100\) & \(2171/5000\)\\
\(10^3\) & 6.490 & 10 & \(63/100\) & \(4417/10000\)\\
\(10^4\) & 6.442 & 10 & \(51/100\) & \(1113/2000\)\\
\(10^5\) & 6.346 & 16 & \(43/100\) & \(3403/5000\)\\
\(10^6\) & 6.295 & 21 & \(17/50\) & \(7703/10000\)\\
\(10^{12}\) & 6.258 & 300 & \(1/8\) & \(97/100\)\\
\bottomrule
\end{tabular}
\caption{Effective analytic orders for number fields.
Each coefficient bounds \((S+1)/k\) for all \(k\ge k_0\).}
\label{tab:nf-two-coefficient-finite}
\end{table}

\begin{proof}
For the first eight rows take
\(\lambda=999999/1000000\).
Use \(M_i=9/5,19/10,2,\ldots,R\) in steps of
\(1/10\), with \(L_i=\lambda M_i\).
At each \(M_i\), let \(t_i\) be the largest multiple
of \(10^{-7}\) supported over \([1,M_i]\)
by the intervals of Lemma~\ref{lem:nf-fibre-savings}
and by \eqref{eq:joint-finite-terms}, after decreasing
it by one unit of the grid for strictness.
Here the radical exponent is \(M_i/k_0\).
The conductor intervals for the two-coefficient
estimate are the solutions of
\(\min(1,E_1,E_2,E_3,E_4)\ge D(D+1)t_i\).
At the tail the low-conductor test is repeated
with \(\theta=h/\sqrt{k_0}\).
The interval union retains coverage with
strict slack at \(c/R^2\).

With \(d_0=t_0/L_0\), start at even order
\(k(1-d_0-\log d_0)+O(1)\).
To move from \(t_i/L_i\) to \(t_{i+1}/L_{i+1}\)
the linear decrement costs
\(2\log(t_i/t_{i+1})/L_i\), and the constant
decrement costs \(2/L_i-2/L_{i+1}\).
The latter telescope. With \(L=\lambda U\)
and the unscaled tail saving \(c/U^2\),
the ideal normalized total is
\begin{equation}\label{eq:joint-table-discrete}
  1-d_0-\log d_0+\frac{2}{\lambda M_0}
  +\sum_{i<{\rm last}}\frac{2}{\lambda M_i}
                          \log\frac{t_i}{t_{i+1}}
  +\frac{4+2\log(t_{\rm last}R^2/c)}{\lambda R}.
\end{equation}
For the \(10^5\) and \(10^6\) rows instead
integrate affine envelopes through \(R\).
They include the two-coefficient lines, the
degree-reduction and differencing lines and
the arithmetic pairs to depth seven. Take
\(\lambda=99202/100000\) and \(99863/100000\),
respectively. The maximum relative loss on each
supporting line, calculated at its endpoints,
gives the two checks
\(\rho R/k_0<1-\lambda\) and
\((1-\rho h/\sqrt{k_0})s(R)>c/R^2\).
Use \(\lambda s(U)\) and \(L=\lambda U\)
throughout, including height \(\lambda c/U^2\) in
the tail. If \(g\) is the seed order divided by \(k\)
without parity rounding, the normalized ideal total
is \(g+(\mathcal C_0-g)/\lambda\), where
\(\mathcal C_0\) is the ideal price with \(\lambda=1\).
As before the tail beyond \(R\) adds
\((4+2\log(s(R)R^2/c))/R\) to the affine integral.

In both calculations the preliminary marked sum
from Theorem~\ref{thm:nf-all-linear} saves \(1/(15k)\)
on the broad minor arcs, by \eqref{eq:nf-self-sup}.
Choose \(a=\lceil2k/(\lambda H)\rceil+1\)
and stop at
\(d_*=(a-1)(1-1/k)/(15k^2)\).
The fixed family of prime windows can accommodate both
the preliminary order and the present finite menu:
take the union of their finitely many conditions
before choosing the windows, as in the proof of
Theorem~\ref{thm:nf-all-linear}. The prescribed stop
leaves positive slack for the initial-scale loss.
The ideal tail beyond the cap costs \(6/(\lambda H)\).
Its capped substitute costs at most
\(2\log(d_H/d_*)/(\lambda H)\), with
\(d_H=c/(\lambda H^3)\) for
\eqref{eq:joint-table-discrete} and \(c/H^3\) for
the scaled affine envelope. The \(a-1\) saving factors
and the additional ordinary factor cost at most
\(2/(\lambda H)+2/k\). The total upper correction,
including rounding and finite-menu approximation, is
therefore
\[
 \frac{2}{\lambda H}
  \max\left(0,\log\frac{15c\,\lambda_*}
                              {2h^2(1-1/k)}-2\right)
 +\frac{16}{k},
\]
where \(\lambda_*=1\) in the discrete calculation
and \(\lambda_*=\lambda\) for the scaled envelope.
The allowance includes two units each for the even seed,
the final step and order parity, the ceiling in \(a\),
and the extra root. The remaining margin pays for a
finite approximation of the continuous menu and for
strict stopping. Use the lower endpoint of each band
in the logarithm and in \(H\) and \(k\) in this
correction. The exact comparisons in the electronic
certificate bound the resulting normalized order by
the entries of \cref{tab:nf-two-coefficient-finite}.
The local-to-global step yields the displayed
representation bound.

\paragraph{An explicit bound for \(k\ge10^{12}\).}
Use \(R=300\), \(h=1/8\), \(c=97/100\) and
\(\lambda=999999/1000000\).
Up to \(R\) take the envelope of the degree-reduction,
iterated-difference and two-coefficient estimates,
and the arithmetic exponent pairs through depth twenty-five
from Lemma~\ref{lem:nf-fibre-savings}, Proposition~\ref{prop:nf-dual} and Corollary~\ref{cor:nf-dual-radical}.
The words have matrix entries bounded by \(2^{25}\) on the
local binomial power, so numerator, denominator and difference
tests on these words are nondegenerate for \(k\ge10^{12}\).
The finite relative radical loss is less than \(903\theta\)
on the supporting segments. At depth twenty-five the
seed and ordinate at \(R\) are \(7393/4141\) and \(1/90300\),
respectively. We take \(L=\lambda U\) and retain the
fraction \(\lambda\) of every saving. The prefix uses
\(\theta=R/k\), and for the low conductors inside a tail
step the test is
\((1-903H/k)/90300>c/R^2\).
Both hold at \(k=10^{12}\) and improve thereafter. For
the degree-reduction and differencing lines the
absolute radical loss per unit \(\theta\) is bounded
respectively by
\(j/\min(2^{j-1},j(j-1))\) and \(l/2^{l+1}\).
For the two-coefficient lines the bound is
\((3D-1)/(D(D+1))\), by \eqref{eq:joint-finite-terms}.
For the arithmetic pairs the bound is \(2\),
including the original disk partition.
These are the segment tests used to obtain \(903\).

The continuum with the tail \(c/U^2\) and no cap
has price at most \(6.257246548339470084175047\),
including the bridge at \(R\). To evaluate it
use \eqref{eq:joint-price} on the finite affine pieces
of the running minimum, adding
\[
       \frac{2}{R}\log\bigl(s(R)R^2/c\bigr)+\frac4R .
\]
This accounts for the drop to \(c/R^2\) and then
the unbounded tail. Scaling both the block exponent
and the saving by \(\lambda\) bounds the ideal price by
the same number divided by \(\lambda\).
At \(H\) the normalized excess is \(d_H=c/H^3\).
Stop at \(d_*=\lceil16\sqrt{k}\rceil/(2k^3)\).
Then \(d_H/d_*\le62.08k\).
Using the ordinary compensation of Theorem~\ref{thm:nf-absorb-transfer} with \(\lceil16\sqrt{k}\rceil\) factors
and retaining the ideal tail in the upper bound, the
normalized order including one extra root is less than
\[
 \frac{6.257246548339470084175047}{\lambda}
 +\frac{16}{\lambda\sqrt{k}}\log(62.08k)
 +\frac{16}{\sqrt{k}}+\frac{16}{k}.
\]
The last term covers the even seed, the last step,
rounding the final order and the discretization
allowance. The expression decreases for \(k\ge10^{12}\).
An outward rational logarithm bound at that endpoint
gives \(6.257776957237\) as an upper bound. Thus the
even order may be chosen with \(S+1<6.258k\)
throughout this range. The parameters and the
arithmetic tests are recorded in
Appendix~\ref{app:certificates}.
\end{proof}

\subsection{Small fractional parts at all embeddings}
\label{sec:nf-fractional}
Write $\|x\|_\infty=\max_{v\mid\infty}|x_v|$ for $x\in V$
and $\|t\|_{\mathbb R/\mathbb Z}=\min_{m\in\mathbb Z}|t-m|$.
The pointwise estimate for the marked sum has a consequence on the
full trace torus, even when the coefficient is specified separately
at each infinite place.
\begin{theorem}\label{thm:nf-fractional}
Let $k\ge11$, $\alpha\in V$, and $P\ge2$. Let \(S\) be an even order supplied by
Theorems~\ref{thm:nf-all-linear}, \ref{thm:nf-optimized-moment}
or~\ref{thm:nf-two-coefficient-finite}. There exist a nonzero
totally positive $x\in\mathcal O$ with $\|x\|_\infty\le P$ and a
$c\in\mathcal O$ such that
\[
                  \|\alpha x^k-c\|_\infty\ll_{\mathbb K,k} P^{-1/[2(S+1)]}.
\]
For the single trace, possibly with a different totally positive $x$,
\[
 \|\operatorname{Tr}(\alpha x^k)\|_{\mathbb R/\mathbb Z}
      \ll_{\mathbb K,k} P^{-\min\{n/[2(S+1)],\ n/(n+S+1)\}} .
\]
More generally, for $m$ rationally independent $b_1,\ldots,b_m\in\mathbb K$
one may simultaneously make all $m$ distances
$\|\operatorname{Tr}(b_j\alpha x^k)\|_{\mathbb R/\mathbb Z}$
of order $O_{\mathbb K,k,\mathbf b}(P^{-\delta_m})$, where
\[
          \delta_m=\min\{n/[2m(S+1)],\ n/[n+m(S+1)]\}.
\]
The estimates are uniform in \(\alpha\).
The full simultaneous exponent is greater than \(1/(13.746k)\)
for \(k\ge11\), and greater than \(1/(12.439144428k)\)
for all sufficiently large \(k\). The effective coefficients of
\cref{tab:nf-two-coefficient-finite} give exponent \(1/(2Bk)\)
on their respective ranges.
\end{theorem}
\begin{proof}
Choose $0<\theta<1/2$ and $\theta<\omega<1$. Put
$M=X^\theta$, $U=X/M$. Let $v$ range over balanced totally
positive generators of principal prime ideals in the identity
narrow class with $M\le|Nv|\le2M$, choosing one generator per
ideal. Let $u$ range, with its original weights, over the marked
family at volume $U$ for which the $S$th moment is loss-free. The
products $x=vu$ are counted with multiplicity. Intersecting a
finite number of archimedean sectors if necessary, the total mass is
$\gg X/\log X$ by the prime ideal theorem in narrow classes and the
positive density of the inner family. Choose its fixed archimedean box
small enough that $\|vu\|_\infty\le P$. Write $F_*(\alpha)$ for
the exponential sum over this positive family. For every $\epsilon>0$,
\[
 \sup_{\mathbb T\setminus\mathfrak M(X^\omega)}
          |F_*(\alpha)|\ll X^{1-\theta/S+\epsilon}.
\]
Indeed the sampling calculation of \eqref{eq:nf-self-sup} has
a single outer prime and a common inner sum. The determinant
bound $O((M/U)^k)=o(1)$ gives at most one slope $\gamma$ for
all nonsingleton cells, with $Q_\gamma\ll M^k$. At most $O_k(1)$
outer primes can divide this conductor, contributing $O_k(U)$.
For all other primes the cell occupancy is
$O_\epsilon(X^\epsilon(1+M/(Q_\gamma
B(P^k(\alpha-\gamma)))))=O_\epsilon(X^\epsilon)$.
The cell large sieve and H\"older give
$|F_*|^S\ll X^{S+\epsilon}/M$. Singleton cells obey the same
estimate, and the inner weighted moment is uniform on the chosen
fixed boxes.

For the simultaneous claim put $\delta=1/[2(S+1)]$ and choose
$\theta,\omega$ such that $S\delta<\theta<1/2$ and
$\theta<\omega<1-\delta$. Set $\sigma=P^{-\delta}/\log(2P)$.
On $V/\mathfrak D^{-1}$ periodize $\exp(-\pi|z|^2/\sigma^2)$,
where $|\cdot|$ is a fixed Euclidean norm. Its Fourier series
has constant coefficient $c_{\mathbb K}\sigma^n$, sum of absolute
coefficients $O_{\mathbb K}(1)$, and sum $O_{\mathbb K,k,A}(P^{-A})$
outside $\|h\|_\infty\le P^{\delta}(\log(2P))^2$, $h\in\mathcal O$.
These facts follow from Poisson summation and equivalence of norms.
If all nonzero frequencies in this range have
$h\alpha\notin\mathfrak M(X^{\omega})$, the weighted sum of this
Gaussian over $\alpha x^k$ is
\[
 c_{\mathbb K}\sigma^n F_*(0)+O(X^{1-\theta/S+\epsilon}).
\]
Its main term dominates. At distance at least $P^{-\delta}$ from
the lattice, the periodized Gaussian is $O_{\mathbb K,k,A}(P^{-A})$
for every fixed $A$. This gives a totally positive root and the
required approximation with $c\in\mathfrak D^{-1}$.

Otherwise choose a nonzero such $h$ and write
$h\alpha=\gamma+\beta$, with $QB(P^k\beta)\le X^{\omega}$,
where $Q=N\mathfrak q_\gamma$. Put
\[
 R_v=Q^{1/n} B(P^k\beta)^{1/(nk)}
                 \max(1,|P^k\beta_v|)^{-1/k}.
\]
The weighted product of these radii is $Q$. There is a totally
positive $t\in\mathfrak q_\gamma$ with $|t_v|\asymp_{\mathbb K}R_v$
at every place. Indeed, multiply the conductor by one of finitely
many integral representatives of inverse narrow ideal classes,
choose a positive generator of the product, and adjust by positive
units in the logarithmic hyperplane. Consequently
\[
 \|t\|_\infty\ll P^{\omega},\qquad
 \max_v|\beta_vt_v^k|\ll P^{-k}Q^{k/n}B(P^k\beta)^{1/n}
                           \ll P^{-(1-\omega)k}.
\]
Choose $b\in\mathcal O$ from a fixed finite set realizing the
real sign patterns, so that $bh$ is totally positive. Then $x=bht$
has $\|x\|_\infty\le P$ for large $P$, and
$c=b^kh^{k-1}t^k\gamma\in\mathfrak D^{-1}$. The error is
$\ll P^{-(1-\omega)k+(k-1)\delta}(\log(2P))^{2k}$.
This proves the assertion with a codifferent lattice point.
Apply it to $\alpha/d$, with a fixed positive rational integer
$d$ satisfying $d\mathfrak D^{-1}\subseteq\mathcal O$, to obtain
the integral $c$ in the statement.

For the $m$ traces, choose a rational integer $d>0$ with
$b'_j=db_j\in\mathcal O$, and work with $\alpha' = \alpha/d$ and
the $b'_j$. Set $\delta=\delta_m$. Since
$\delta<n/(2mS)$ and $\delta<n/(n+mS)$, one can choose
$mS\delta/n<\theta<\min(1/2,1-\delta)$ and
$\theta<\omega<1-\delta$. Periodize the Gaussian on
$\mathbb R^m/\mathbb Z^m$. Its constant coefficient is $\sigma^m$;
the nonzero Fourier modes of interest have
$h=\sum_j a_jb'_j\ne0$ with $a_j\in\mathbb Z$ and
$\max|a_j|\le P^\delta(\log(2P))^2$. The off-core bound again
saves a strictly larger power than the constant coefficient. In the
other case, the construction for $h\alpha'$ gives
$x=bht$ of size $O(P^{\omega+\delta}(\log(2P))^2)$ and
$\alpha' x^k=c+O(P^{-(1-\omega)k+(k-1)\delta}(\log(2P))^{2k})$
coordinatewise for a $c\in\mathfrak D^{-1}$. Each
$\operatorname{Tr}(b'_jc)$ is an integer. Bounded $P$ is absorbed
by the covering radii of the lattices. The second assertion follows
by taking $m=1$, $b_1=1$.
\end{proof}

There is a second bound valid in every degree $k\ge2$. Put
$e_k=\min(2^{k-1},k(k-1))$. In the assertion for $m$ traces one may
take every fixed exponent strictly below
\begin{equation}\label{eq:nf-ordinary-fractions}
                  \min\{n/(me_k),\ n/(n+mk)\}.
\end{equation}
For approximation by a single algebraic integer at every embedding,
take $m=n$ in this formula. Indeed on the broad minor arcs the
ordinary smooth sum satisfies $|f|\ll X^{1-1/e_k+\epsilon}$.
For the Weyl branch, repeat the proof of \eqref{nf3:fourweyl}
after $k-1$ differences. The nonzero product $t$ of the shifts
has height $O(P^{k-1})$ and multiplicity $O(X^\epsilon)$.
With $M$ the largest occupancy at precision $1/P$, the resulting
bound, divided by $X^{2^{k-1}}$, is
$O(X^\epsilon(X^{-1}+M/X^{k-1}))$. The full-rank bound is
$M\ll X^{k-2}$ by the real shear
$(P^{1-k}t,P(\alpha t-c))$; in the line case it is
$O(1+X^{k-1}/(Q_\gamma B(P^k(\alpha-\gamma))))$.
The critical mean-value branch is \eqref{eq:nf-point}.
The quadratic completion in \eqref{eq:nf-complete-bound} uses
the same Poisson argument, with complete quadratic sums
$O_{\mathbb K}(Q^{-1/2+\epsilon})$ uniformly in the linear twist,
as follows also by one differencing and character orthogonality.
Consequently on the larger complement of $\mathfrak M(X^\omega)$,
$\omega<1$, the normalized ordinary sum saves
$\min(1/e_k,\omega/k)-\epsilon$ in the power of $X$.
For any $\delta$ strictly below \eqref{eq:nf-ordinary-fractions}
choose $mk\delta/n<\omega<1-\delta$. Use this ordinary estimate
in the Gaussian argument and keep the same conductor multiplier.
For the full torus use its $n$-dimensional Gaussian, with integral
lattice as before. This proves the assertions.

The construction also counts distinct approximating roots. For the
$m$ traces and $k\ge11$, every
$0<\delta<\min(n/(2mS),n/(n+mS))$ gives at least
\begin{equation}\label{eq:nf-fractional-multiplicity}
                \gg_{\epsilon,\mathbb K,k,\mathbf b,\delta}
                           P^{n-(n+mS)\delta-\epsilon}
\end{equation}
totally positive roots of height at most $P$ satisfying all $m$
approximations with error $O(P^{-\delta})$. At the full-field
exponent $1/[2(S+1)]$ at least $\gg P^{n/2-\epsilon}$ distinct
roots give the claimed integral approximation.
In \eqref{eq:nf-ordinary-fractions} the analogous count has $k$
in place of $S$ on the right of \eqref{eq:nf-fractional-multiplicity}.
Indeed, on the small-value Fourier alternative the Gaussian mass is
$\gg P^{n-m\delta}(\log(2P))^{-m-1}$. Each physical root has bounded
total weight: at most $O_{\theta}(1)$ primes of the outer norm range
divide it, and its inner weight is bounded. Points outside a fixed
multiple of the indicated error make a negligible contribution.
On the core alternative, multiply the constructed $bht$ by any
totally positive $y$ of height at most $P^{1-\omega-\delta-\epsilon'}$.
There are $\gg P^{n(1-\omega-\delta-\epsilon')}$ such $y$ by lattice
counting in a positive box. The products are distinct; their powers
still satisfy the error requirement since $y^k c\in\mathfrak D^{-1}$,
and in the trace case $\operatorname{Tr}(b'_j y^k c)\in\mathbb Z$.
The same fixed scalar adjustment as above returns the ordinary
integer lattice in the full-field assertion. For the marked family,
choose $\theta$ and then $\omega$ arbitrarily close to $mS\delta/n$
from above. For the ordinary family choose $\omega$ close to
$mk\delta/n$. Taking $\epsilon'$ sufficiently small gives the counts.

\subsection{Mixed powers and consecutive exponents}
\label{sec:nf-mixed}
For a fixed tuple $k_1,\ldots,k_t\ge2$ put $b=\sum_i1/k_i$.
The sharp moments for the individual powers can be applied at a
common target volume. This extends the mixed-power problem considered
over the rationals by Br\"udern and Wooley \cite{NFMixedBW}.

\begin{theorem}\label{thm:nf-mixed}
Suppose every finite completion has a representation of $\nu$ by the
listed powers with at least one unit among the roots. For $b\ge6.873$,
every sufficiently large totally positive $\nu$ with this property
has $\gg T^{b-1}$ representations by totally positive roots,
where $T=|N\nu|$. One may take $|Nx_i|\ll T^{1/k_i}$.
For $b\ge3.4365$ the same conclusion holds apart from
$O(T^{1-\delta})$ targets in any fixed balanced dyadic region of
volume $O(T)$, for some $\delta>0$. The constants may depend on the
field, the tuple and the region. The two
reciprocal-sum thresholds may be replaced by $C_1$ and
$C_1/2$, respectively, once the least degree is sufficiently
large independently of the field. For each row \((k_0,B)\) of
\cref{tab:nf-two-coefficient-finite} they may be replaced by \(B,B/2\)
when every \(k_i\ge k_0\).
\end{theorem}
\begin{proof}
Balance the target by a power of a totally positive unit, using as
exponent a common multiple of the $k_i$. Division or multiplication by
that power preserves all local and global conditions, and the norm of
each root is unchanged. In a balanced dyadic region choose $X_i=T^{1/k_i}$
as the successive root volumes, taking $T$ to be a common dyadic norm
scale. For $k_i\ge11$ let $S_i$ be the even order from
Theorem~\ref{thm:nf-all-linear}. For $2\le j\le10$ the respective
orders at degree $j$ may be taken as
\[
  \begin{array}{c|rrrrrrrrr}j&2&3&4&5&6&7&8&9&10\\\hline
                  S&4&8&12&18&24&32&38&44&52\end{array}.
\]
For degree two, the fourth moment counts two equal differences
of squares. The zero difference gives $O(X_i^2)$ solutions.
For each nonzero difference the factorization $(x-y)(x+y)$ has
$O_\varepsilon(X_i^\varepsilon)$ realizations in the box, by
the ideal divisor bound and the logarithmic number of possible unit
multiples. Summing over the first pair gives $O_\varepsilon
(X_i^{2+\varepsilon})$. For cubes use the eighth moment
\eqref{nf3:eq:hua}. For $4\le j\le10$ these orders are the
even moments in \eqref{eq:aud-small-even}. Take the marked windows
and the terminal cofactor within the smoothness range specified
there. Bounded coefficient restrictions preserve the estimates
at even order by orthogonality.

Let $g_i=F_i/X_i$ be the normalized marked sum on a small archimedean
box. Each satisfies $\int|g_i|^{S_i}\ll_\epsilon T^{-1+\epsilon}$, and
$\sup|g_i|\ll T^{-\delta_i}$ off the common core
$Q_\gamma B(T^{1/n}(\alpha-\gamma))\le T^\eta$.
Here one chooses $\eta>0$ below all the rescaled prime-window
and free-cofactor exponents, and applies
Lemma~\ref{lem:nf-bilinear-weyl}. For $u=1,2$ put
$J=\sum_i u/S_i$. If $J>1$, H\"older's inequality with exponents
$JS_i/u$ on the factors $|g_i|^{u/J}$ bounds the integral of
$\prod_i|g_i|^u$ off the core by
\[
 T^{-1+\epsilon-\sum_i\delta_i u(1-1/J)}.
\]
For every displayed small degree $S_i<6.873k_i$, and for the
larger degrees we even have $S_i+1<6.873k_i$. Thus $J>1$
also at the stated reciprocal-sum endpoints. Choose the allowable
moment losses for this fixed tuple below the resulting strict
surplus, before choosing all its windows. More generally, $J>1$ together with $b>2$ suffices for each assertion.

On the core, Poisson summation in the separate free cofactors gives the
ordinary mixed complete sums and archimedean factors, multiplied by
positive factors bounded above and below. The proof in
Section~\ref{sec:nf-marked-weyl} applies to each exponent on reducing
$\eta$ once for the fixed tuple. Every mark is invertible modulo
each core denominator. The normalized complete sum in degree $k_i$
is $\ll_\epsilon Q^{-1/k_i+\epsilon}$; hence the conductor tail is
absolutely convergent with a power saving for $b>2$.
The product of the concentrated archimedean factors has the bound
$\ll B(T^{1/n}\beta)^{-t/2}$ by stationary phase away from the
coordinate walls. It too has a power-saving tail, using compact
smooth weights. Near a given normalized target $y$, take the $i$th
weight positive near an archimedean $k_i$th root of $y/t$, with positive
roots at the real places. The derivative in one coordinate is
invertible at every place. Compactness gives uniform lower bounds and
a finite choice of boxes on the balanced target region.

The local series factors as the product of the mixed local
densities. At large residue characteristic the complete-sum estimates
give factors $1+O(q^{-1-\rho})$ by summing the character conductors,
with some $\rho>0$ since $b>2$. At each of the finitely many remaining
places one root is a unit, and its derivative valuation is bounded
by $\max_i v(k_i)$. Hensel lifting gives a lower bound depending
only on the field and tuple. Core truncation therefore gives
$\gg T^{-1}$ for the normalized main term, with a power-saving
error. For $u=1$ integrate directly. For $u=2$ a missing target
has a complementary Fourier coefficient of modulus $\gg T^{-1}$.
Bessel's inequality and the squared-product bound $O(T^{-1-\delta})$
therefore leave $O(T^{1-\delta})$ such targets. For the final thresholds use
\eqref{eq:nf-secondary} or Theorem~\ref{thm:nf-two-coefficient-finite}
in the same inequalities. The strict upper bounds for \(S_i\)
give the surplus at the reciprocal-sum endpoints.
\end{proof}

For consecutive powers the local hypothesis disappears. Write
$R_{\mathbb K}(k,1)$ for the least number of exponents beginning at $k$
and increasing by one which represent every sufficiently large totally
positive integer in $\mathcal O$. There is one summand of each exponent.
\begin{corollary}\label{cor:nf-mixed-consecutive}
Uniformly in the number field,
\[
 \limsup_{k\to\infty}\sup_{\mathbb K}
         \frac{R_{\mathbb K}(k,1)}{k}\le e^{C_1}-1<499.805.
\]
For almost all targets in balanced dyadic regions, in the sense of
Theorem~\ref{thm:nf-mixed}, the upper constant is
\(e^{C_1/2}-1<21.379\).
\end{corollary}
\begin{proof}
For any fixed $0<\epsilon\le1$, an interval of at least $k^\epsilon$
consecutive exponents starting at $k$ gives primitive solubility
for every target at every finite place, once $k$ is large. Let $q$
be the residue norm and $p$ its characteristic. If $q>(10k)^{10}$,
choose five exponents among $k,\ldots,k+9$ prime to $p$. Character
orthogonality and the finite-field Weil bound give $q^4+
O((2k)^5q^{5/2})$ solutions of the five-term residue equation,
with absolute implied constant one. More than one solution exists;
a nonzero root has unit derivative and lifts.
If $q\le(10k)^{10}$, put $m=p(q-1)\le(10k)^{20}$. The number of
integers in the interval $I$ prime to $m$ is
$|I|\varphi(m)/m+O(2^{\omega(m)})$ by inclusion--exclusion.
Both $m/\varphi(m)$ and $2^{\omega(m)}$ are $k^{o(1)}$ uniformly in
this range, so at least two such exponents occur. Each power map
is a bijection on the units of the local integer ring, by the cyclic
residue group and Hensel lifting. A unit target needs only one of
these roots. A nonunit is the sum of the two units $1$ and $\nu-1$.
Choose the remaining local roots zero. Finally
$\sum_{j=k}^{k+s-1}1/j\ge\log(1+s/k)$; the theorem and its
almost-all variant yield the constants.
\end{proof}

More generally suppose $(k,d)=1$ and let $R_{\mathbb K}(k,d)$ denote
the number of terms in the analogous progression with common difference $d$.
For every sufficiently large $k$, uniformly in $d$ and in $\mathbb K$,
\begin{equation}\label{eq:nf-mixed-arbitrary-step}
 R_{\mathbb K}(k,d)\le
 \left\lceil\frac{k}{d}\big(\exp(C_1 d)-1\big)\right\rceil .
\end{equation}
For almost all targets one may replace $C_1$ by $C_1/2$.
To check finite solubility, take $L$ to be the right side with either
constant. For $q>(100(k+d))^{10}$ use five among the first ten exponents
not divisible by $p$. When $p\mid d$ all ten qualify. For
$q\le(100(k+d))^{10}$, set $m=p(q-1)$ as before. Primes dividing
both $d$ and $m$ divide none of the exponents. For each squarefree
divisor $h$ of $m$ prime to $d$, the indices with
$h\mid k+id$ lie in one residue class modulo $h$. Inclusion--exclusion
therefore gives at least $L\varphi(m)/m-O(2^{\omega(m)})$ indices
prime to $m$. This exceeds two uniformly, since $L\gg k+d$.
Finally compare the reciprocal sum with $d^{-1}\log(1+dL/k)$.

There is a secondary term for fixed $d$. Put $u=1$ for all targets
and $u=2$ for almost all, and write $E=\exp(C_1d/u)$. For every
\[
                  B<\frac{2(E-\sqrt E)}{dC_1}
\]
the same conclusions hold with at most
$\lceil(E-1)k/d-B\sqrt k\rceil$ exponents, once $k$ is large,
$(d,k)=1$. Indeed, for $j=k+id$ sufficiently large,
\[
 \frac1{S_j}\ge\frac1{C_1 j}+\frac1{C_1^2 j^{3/2}}
\]
by \eqref{eq:nf-secondary}. Decreasing-function integral comparison
at $L=(E-1)k/d-B\sqrt k+O(1)$ gives
\[
 \sum_{0\le i<L}\frac{u}{S_{k+id}}\ge
 1+\frac{u}{C_1\sqrt k}
 \left\{\frac{2(1-E^{-1/2})}{dC_1}-\frac B E\right\}
 + O_{d,B}(k^{-1}).
\]
The strict positivity gives the moment surplus required in the proof.
For consecutive exponents one can take \(B=154\) in the all-target result.

There is a different local phenomenon for the degrees
$k,2k,\ldots,Lk$: their global sums lie in $J_k$. For each fixed
field, if $R_{\mathbb K}(k,k)$ denotes the least such $L$ representing
every sufficiently large totally positive element of $J_k$, then
\begin{equation}\label{eq:nf-mixed-multiples}
 \limsup_{k\to\infty} k^{-1}\log R_{\mathbb K}(k,k)\le C_1.
\end{equation}
Indeed, take $L=\lceil\exp((C_1+\epsilon)k)\rceil$. At residue
norms $q>(100k)^{10}$ with $p\nmid k$, use the preceding
finite-field argument on five exponents $ik$, $1\le i\le10$,
prime to $p$. At every other place $q\le\max((100k)^{10},k^n)$
for $n=[\mathbb K:\Q]$. The same inclusion--exclusion count gives
more than $4k+O_{\mathbb K}(1)$ indices $i\le L$ prime to $p(q-1)$.
For sufficiently large $k$ the lifting precision of
Lemma~\ref{lem:precision-widths} is less than $k$ at every place
since $e\le n$. The lemma represents any eligible local target by a
nonempty sum of at most $H_{\mathbb K}(k)\le4k+O_{\mathbb K}(1)$ unit
$k$th powers, lifting the final error in one unit. Assign a distinct
chosen $i$ to each unit and take its $i$th root in the unit group.
The resulting $ik$th powers give the required primitive solution.
The reciprocal sum is at least $k^{-1}\log(1+L)$. Apply
Theorem~\ref{thm:nf-mixed} and let $\epsilon$ tend to zero.

\subsection{The local order for a progression with fixed common divisor}
\label{sec:nf-mixed-unit-order}
Let $C_g(A)$ be the additive subgroup generated by $g$th powers of
units in a local integer ring $A$. It contains $\pi^c A$ whenever
$1+\pi^c A\subseteq(A^\times)^g$, by subtracting one from a unit
power. In particular the group is open and closed. Given
$g\ge1$, define
\[
 J_g^\times(\mathbb K)=\{\nu\in J_g:
          \nu\in C_g(\mathcal O_v)\ \text{at every }v\mid g\}.
\]
At $p\nmid g$ one has $C_g(A)=J_g(A)$ since $\pi A\subseteq C_g(A)$.

\begin{theorem}\label{thm:fixed-gcd-mixed}
Fix $g$. For every number field, all sufficiently large $k$ depending
only on $g$, and every $d$
with $(k,d)=g$, put
\[
          L=\left\lceil (k/d)(\exp(C_1d)-1)\right\rceil .
\]
The sums with exponents $k,k+d,\ldots,k+(L-1)d$ belong to
$J_g^\times(\mathbb K)$. Every sufficiently large totally positive
element of this order has such a representation, with
$|Nx_i|\ll |N\nu|^{1/(k+id)}$ for the root indexed by $0\le i<L$.
For almost all targets of this order in balanced dyadic regions,
one can replace $C_1$ by $C_1/2$. The lower bound required for the target norm and the implied constants
in the root estimates may depend on the field and tuple.
\end{theorem}
\begin{proof}
All exponents are multiples of $g$. We first give a uniform
threshold ensuring that their nonunit values lie in $C_g(A)$.
For a prime $p\mid g$ put $a=v_p(g)$, $g_0=g/p^a$, $r=a+1$ and
let $N=2\delta D+1$ be the precision of \eqref{eq:nf-chart};
$N$ depends only on $p,a$. Write a nonunit root as $\pi y$.
For an exponent $j=gt$ lift $y$ to $\tilde y\in R$ in the
universal uniformizer calculation. If $g_0t\ge N$, the
Teichm\"uller vector $[(T\tilde y)^{g_0t}]$ lies in the kernel of
reduction modulo $T^N$. Formula \eqref{eq:nf-chart} expresses it
as a difference of sums of unit $g_0$th-power Teichm\"uller
vectors. Evaluation of the final ghost component shows that
$(\pi y)^j\in C_g(A)\bmod p^{a+\delta}A$. The latter ideal belongs to
$C_g(A)$ by the logarithm on $1+p^\delta A$. Taking the maximum
of $p^a N$ over $p\mid g$ proves necessity uniformly.

Every target of $C_g(A)$ can be represented by between one and
$D_g$ unit $g$th powers, with $D_g$ depending only on $g$.
At $p\nmid g$ represent the target minus one modulo $\pi$ by at
most $g$ unit powers or zeros by \eqref{eq:nf-residue-width}; add
a unit absorbing the final error. At $p\mid g$ apply the proof of
\eqref{eq:nf-uniform-witt} to the subgroup generated by
$[y^{g_0}]$ for $y\in R^\times$. The centered chart there uses
only units, and multiplication by any allowed Teichm\"uller constant
preserves this set of generators. The same successive cyclic
quotients and the same length bound therefore give the uniform
cost $a+1+g_0 U(p,a)$. Each constant multiple used in the digit
construction can be expanded into a positive sum of allowed unit
directions, omitting zero digits. The chart supplies a unit for
the final lift. Take the maximum over the primes dividing $g$.

Choose a sufficiently large fixed $C=C(g)$. At
$q\le(C(k+d))^C$ use the inclusion--exclusion argument of
\eqref{eq:nf-mixed-arbitrary-step} for $k/g+i(d/g)$, whose
initial term and step are coprime. More than $D_g$ of these
indices are prime to $p(q-1)$, uniformly in $d$ for large $k$,
since $L\gg k+d$. Assign distinct indices to the required unit
powers and extract the corresponding unit roots. Put the unused
local roots equal to zero. At $q>(C(k+d))^C$ with $p\nmid g$ at
least five of the first ten exponents are prime to $p$, and the
finite-field argument of Corollary~\ref{cor:nf-mixed-consecutive}
gives a nonsingular solution.

It remains to handle $q>(C(k+d))^C$ above a prime dividing $g$.
Set $s=2(a+\delta)+1$ and $l=r+s$. Among the first $2l$ exponents
take $l$ with $p\nmid j/g$. Their prime-to-$p$ parts have least
common multiple $b\le(2l(k+d))^l$, prime to $p$. All $l$ selected
exponents divide $p^a b$. We show that every $g$-eligible target
there is a sum of $l$ unit $p^a b$th powers.

Let $F=\mathbb F_q$ and $B=W(F)$. After increasing $C(g)$,
the group $(F^\times)^b$ has more than $\sqrt q$ elements, so
its additive span is $F$. Also $[F:\mathbb F_p]\ge r$. Take $r$
independent residues from this group and use their Teichm\"uller
lifts as the constant centres in Lemma~\ref{lem:linear-local-deviation}.
The unit Moore determinant shows that the centered sums with
these unit $b$th-power arguments fill the kernel of evaluation at
$T=0$. Lift the target by the $g$th-power ghost construction.
Subtract its constant part and apply this centered chart to the
remainder. After restoring the fixed baseline, the only correction
needed is in $B/p^{a+\delta}B$.

Every element of this coefficient ring is a sum of $s$ allowed
nonzero Teichm\"uller powers. Indeed the main term in the count is
$(q-1)^s/q^{a+\delta}$ by additive orthogonality, while the total
error is bounded by $(b p^{a+\delta-1}\sqrt q+1)^s$ by the
Galois-ring sum estimate \cite{hkhc} at each actual conductor,
applied to the monomial of degree $b$ prime to $p$.
Excluding the zero Teichm\"uller representative accounts for the
additional one. The preceding bound on $b$ shows that the main term
exceeds the error at the stipulated polynomial threshold with $C(g)$
large enough. Raising to $p^a$ permutes all Teichm\"uller constants.
This supplies the correction, and a unit absorbs the final error in
$p^{a+\delta}A$. The reciprocal-sum comparison used in \eqref{eq:nf-mixed-arbitrary-step} and Theorem~\ref{thm:nf-mixed} now give
both global conclusions.
\end{proof}

For a prime common divisor the required order has an explicit local
description. Let $A$ be the integer ring of a finite extension of
$\mathbb Q_p$, with ramification index $e$, residue degree $f$ and
uniformizer $\pi$. If $f=1$ and $e\ge2$, write
\[
 U=\pi^e/p\equiv u_0+u_1\pi\pmod{\pi^2},
\]
where $u_0$ is the Teichm\"uller representative of a nonzero member
of $\mathbb F_p$, and identify $u_1$ with its residue.

\begin{proposition}\label{prop:mixed-unit-squares}
The index $[J_p(A):C_p(A)]$ equals $p$ precisely when $f=1$,
$e\ge2$ and one of the following conditions holds:
\[
\begin{array}{c|c}
 e & \text{condition}\\\hline
 p-1&\bar u_0=-1\\
 p&u_1=-1\\
 e\ge2p,\quad p\mid e &u_1=0\\
 e\ge p,\quad p\nmid e(e+1)&\text{all } U .
\end{array}
\]
At every other place the index is one. If $b_p(\mathbb K)$ counts
the places above $p$ on this list, the relative density of
$J_p^\times(\mathbb K)$ in $J_p$ is $p^{-b_p(\mathbb K)}$.
This is also the exact relative density of global representations in
Theorem~\ref{thm:fixed-gcd-mixed} when $g=p$.
\end{proposition}
\begin{proof}
One has $p\pi^2 A\subseteq C_p(A)$. Indeed
\[
 \frac{(1+\pi+\pi z)^p-(1+\pi)^p-(1+\pi z)^p+1}{p\pi^2}
                   =(p-1)z+\pi Q(z)
\]
for $Q\in A[z]$: all mixed multinomial coefficients are divisible
by $p$, and the total-degree-two term is $p(p-1)\pi^2z$.
The numerator is a signed sum of unit powers and the quotient maps
$A$ onto itself by Hensel's lemma. If $f\ge2$, then $q-1>p$.
For nonunit $x$ all $x+\zeta$ with $\zeta\in\mu_{q-1}$ are units and
$\sum_{\zeta}(x+\zeta)^p=(q-1)x^p$. The open additive subgroup
$C_p(A)$ is closed under $\mathbb Z_p$-scalars, so $x^p\in C_p(A)$.
At $f=e=1$ the subgroup already contains $\mathbb Z_p$.

Take $f=1$, $e\ge2$, and reduce modulo $p\pi^2A=\pi^{e+2}A$.
For $x\in\pi A$, one has
$(1+x)^p\equiv1+px+x^p$ there. Thus $J_p(A)=C_p(A)+p\pi A$;
at most one factor $p$ can occur. Each unit is a Teichm\"uller
constant times $1+x$, the constants lying in $\mathbb Z_p$.
Write $x=\sum_{j\ge1}a_j\pi^j$, $a_j\in\{0,\ldots,p-1\}$.
In its $p$th power all cross terms vanish modulo $p\pi^2$, since
their coefficients are divisible by $p$ and their uniformizer
order before that factor is at least two. As $a_j^p\equiv a_j\bmod p$,
the unit-power span in the quotient is generated by
\[
  1,\quad \pi^p+p\pi,\quad \pi^{pj}
                         \quad(2\le j,\ pj\le e+1).
\]
Each displayed generator occurs by taking the corresponding single
digit. The question is whether $p\pi$, of valuation $e+1$, is in
this span.

If $e+1<p$, the second generator is $p\pi$.
If $e+1=p$, it equals $(1+U)p\pi$, giving the first row.
If $p\mid e+1$ and $e+1>p$, the separate generator $\pi^{e+1}$
gives $p\pi$ up to a nonzero residue scalar. At $e=p$, the second
generator is $p(U+\pi)$; subtracting the integer $pu_0$ shows
that it supplies the last digit precisely when $u_1+1\ne0$.
For $p\mid e$, $e\ge2p$, there is a separate generator $\pi^e$,
and subtracting $pu_0$ supplies the last digit precisely when
$u_1\ne0$.

To see necessity in all these tests, the nonconstant generators
of valuation less than $e$ have distinct valuations divisible by
$p$ and additive order $p$ in the quotient. Their coefficients may
be reduced modulo $p$. Their smallest surviving valuation cannot
cancel against any other, or against an integer, whose only possible
nonzero valuations are $0,e$. In particular when $e\ge p$ and
$p\nmid e(e+1)$ there is no collision and the last digit is missing.
This proves the table. Its conditions are intrinsic: when $e=p-1$
a change of uniformizer preserves $\bar u_0$; when $p\mid e$ it
preserves the vanishing of $u_1$, and at $e=p$ it preserves $u_1$
itself. The Chinese remainder theorem in the eligible group gives
the product, and the global density follows from
Theorem~\ref{thm:fixed-gcd-mixed}.
\end{proof}

In any prescribed signature of degree $n$, the smallest possible
relative density for common divisor $p$ is
\[
                  p^{-\lfloor n/d_p\rfloor},\qquad
               d_2=2,\quad d_p=p-1\quad(p\text{ odd}).
\]
Each exceptional place has local degree at least $d_p$. For odd $p$
use the Eisenstein extension $\pi^{p-1}=-p$ to attain the first
row; for $p=2$ use $\pi^2=2(1+\pi)$ to attain the second. Prescribe
$\lfloor n/d_p\rfloor$ copies and fill the leftover degree with rational
places by the approximation argument of
Corollary~\ref{cor:three-local-density-degree}.

For example, at $A=\mathbb Z_2[2^{1/4}]$ the third row applies.
Modulo $\pi^6$ the unit-square span has eight elements; the square
span has sixteen. Every nonunit power of even degree at least six
already vanishes in this quotient. Thus half the square-eligible
classes are excluded by every sum of such powers, however many
summands are used. Prescribing several copies of this completion
in a totally real field as in Corollary~\ref{cor:three-local-density-degree}
gives independent exclusions. For every $r$ one can also prescribe
$r$ such copies as the entire dyadic algebra of a degree-$4r$
totally real field, obtaining exact relative density $2^{-r}$.

\subsection{A local criterion for every arithmetic progression of exponents}
\label{sec:nf-mixed-criterion}
The dependence on a fixed common divisor in
Theorem~\ref{thm:fixed-gcd-mixed} can be removed. A short initial
segment then accounts for the possible contribution of nonunit roots.
For $g=(k,d)$ and a finite place put
\[
 \mathcal E_v(k,d)=C_g(\mathcal O_v)+
     \sum_{\substack{i\ge0\\ k+id<9g^2}}
                   \{x^{k+id}:x\in\mathcal O_v\}.
\]
The sum here is a sum of sets: one may use one root of each indicated
degree, including zero. The set is open and closed. At $v\nmid g$
it equals the local eligible group $J_g(\mathcal O_v)$.

\begin{theorem}\label{thm:uniform-mixed-criterion}
For all sufficiently large $k$, independently of the field and of
$d\ge1$, put $g=(k,d)$ and
\[
                      L=\left\lceil\frac{k}{d}
                              (e^{C_1d}-1)\right\rceil .
\]
Every sufficiently large totally positive $\nu\in J_g$ is a sum of
totally positive powers of the $L$ successive degrees
$k,k+d,\ldots,k+(L-1)d$, with the norm bounds of
Theorem~\ref{thm:fixed-gcd-mixed}, if and only if
$\nu\in\mathcal E_v(k,d)$ for every $v\mid g$. The condition is
necessary for a sum of any finite number of degrees in the
progression, each used at most once. With $C_1/2$ in the exponent,
the same local criterion suffices outside a power-saving exceptional
set in each fixed balanced dyadic region.

If $k\ge9g^2$, the criterion is $\nu\in J_g^\times(\mathbb K)$.
For a fixed field and all sufficiently large $k$, this simpler
criterion applies simultaneously for every $d$.
\end{theorem}
\begin{proof}
At $p\mid g$ write $a=v_p(g)$. The precision $N$ in the preceding
proof satisfies
\[
 N=2\delta\sum_{j=0}^{a-1}(a-j)p^j+1\le 8p^a+1 .
\]
Thus $p^aN\le9g^2$ and all nonunit powers of degrees $j\ge9g^2$
in the progression belong to $C_g(A)$. This proves necessity.
There are at most $9g+1$ earlier degrees since $d\ge g$.
Also the maximum unit-word length $D_g$ in the preceding proof obeys
$D_g\ll g(1+\log(2g))^2$ by \eqref{eq:nf-uniform-witt} and the
estimate for $U(p,a)$ in the proof of Theorem~\ref{thm:local-sparse}.

Write $B=2(k+d)$, $V=\log(2g)$. At $\log q>C V^2\log B$ for
a sufficiently large absolute $C$, the large-residue argument of
Theorem~\ref{thm:fixed-gcd-mixed} works uniformly also in $g$.
Indeed at $p\mid g$, $r+s\ll V$ and
$\log b\le (r+s)\log(2(r+s)(k+d))\ll V\log B$, while
$p^{a+\delta-1}\le2g$. These inequalities make
$(q-1)^s/q^{a+\delta}>(b p^{a+\delta-1}\sqrt q+1)^s$, since
$s=2(a+\delta)+1$. They also ensure that the allowed residues span
a space of dimension at least $r$. Every $g$-eligible target thus
has a representation with a unit among the specified roots.
For $p\nmid g$ the five-variable finite-field argument applies.

At $\log q\le C V^2\log B$, prescribe the early roots in the
definition of $\mathcal E_v$ locally. For the remaining target,
use at most $D_g$ unit $g$th powers. To assign their exponents,
put $M=p(q-1)$ and delete from its radical the primes dividing $d/g$,
obtaining a squarefree number $M'$. They can be deleted because
$(k/g,d/g)=1$. The indices for which $k/g+i(d/g)$ is prime to
$M$ form a translate of the indices prime to $M'$, since $d/g$
is invertible modulo $M'$. Vaughan's bound for Jacobsthal's function
\cite{NFMixedJacob} gives an acceptable index in every interval of
$O((\log(2M))^4)$ consecutive indices, with an absolute constant;
this also holds when $M'=1$. Indeed his estimate is
$O(\omega(M')^2(\log(2\omega(M')))^4)$, which implies the coarser
form used here. The total length needed beyond the early degrees is
therefore $O(g(1+V)^{10}(\log B)^4)$.

Both proposed progression lengths exceed this quantity plus $9g+1$,
uniformly in $d,g$ for large $k$. When $d\le k$, write the smaller
length as at least $k(e^{C_1 d/2}-1)/d$; its exponential dependence
on $d$ dominates $d^2(1+\log(2d))^{10}$, and
$(\log(4k))^4=o(k)$. For $d>k$ exponential growth in $d$ dominates
the entire required expression. Distinct acceptable indices supply
the required unit roots. Theorem~\ref{thm:nf-mixed} and the
reciprocal-sum comparison now finish the proof, with one unit at
every finite place.

Finally at fixed field degree $n$, the precision in
Lemma~\ref{lem:precision-widths} is at most
$n(1+\log_2g)+1$ at all wild places. When $g$ is large enough
this is less than $g$, and hence every nonunit $g$th power already
belongs to $C_g$. The finitely many smaller $g$ all satisfy
$9g^2\le k$ for large $k$. This proves the last assertion.
\end{proof}

\subsection{Simultaneous congruences on the roots}
\label{sec:nf-strong-approximation}
The small core with an unrestricted cofactor gives a level of
distribution for congruences imposed on the individual roots. We
continue to use root volume $X$ in a balanced region.

\begin{theorem}\label{thm:nf-strong-approximation}
Let $s\ge S+1$, with $S$ any of the even orders provided by the
linear theorem. There is $\eta>0$, depending on the fixed field,
exponent, order and target region, for which the following holds
uniformly in $N\mathfrak a\le X^\eta$. Prescribe classes
$r_i\in\mathcal O/\mathfrak a$. Suppose that at each finite place
there is a solution to $\sum x_i^k=\nu$ having at least one unit
coordinate, and that at $v\mid\mathfrak a$ it can be chosen with
$x_i\equiv r_i\pmod{\mathfrak a\mathcal O_v}$. Then, for every
sufficiently large totally positive balanced target $\nu$, the
number of representations with these prescribed root classes is
\[
                 \gg X^{s-k}(N\mathfrak a)^{-(s-1)} .
\]
The roots may all have the prescribed marked factors from
Section~\ref{sec:nf-marked-weyl} and lie in the conjugate regions
of Section~\ref{sec:nf-global}, choosing the target-dependent
archimedean boxes there. For
$t\ge\max(S/2+1,2k+1)$ there are $\eta',\omega>0$ such that,
apart from $O(X^{k-\omega})$ targets in a fixed balanced box, the
same conclusion with $s=t$ holds simultaneously for every ideal
$N\mathfrak a\le X^{\eta'}$ and every choice of $t$ root classes
satisfying the stated primitive local conditions.
\end{theorem}
\begin{proof}
Put $A=N\mathfrak a$, and restrict the marked sum $F_w$ to roots
in one given class $r$. Choose $U=X^\xi$ as in the polynomial core
argument, first making $\xi$ sufficiently small. Later take $\eta$
sufficiently small compared with $\xi$. Every marked product $D'$
is then invertible modulo $\mathfrak a$. For each residue $u$ of
$D'$ in $(\mathcal O/\mathfrak a)^\times$, the restriction on
the terminal variable is simply $z\equiv u^{-1}r\pmod{\mathfrak a}$.
This gives at most $A$ bilinear sums with bounded coefficients.
Lemma~\ref{lem:nf-bilinear-weyl} yields a bound
$O(AX^{1-\delta})$ off $\mathfrak M(U)$ for the restricted $F_w$,
with $\delta>0$ independent of $A,r$. Its even $S$th moment is
bounded by that of the original positive-coefficient sum. Choosing
$\eta<\delta/2$, H\"older's inequality bounds the off-core integral
of the product by $O(X^{s-k-\delta/2})$.

For a frequency $\gamma$ of conductor $\mathfrak q$, set
$\mathfrak l=\operatorname{lcm}(\mathfrak q,\mathfrak a)$ and
\[
 g_r(\gamma)=\frac{1}{N\mathfrak l}
   \sum_{\substack{x\bmod\mathfrak l\\x\equiv r\ (\mathfrak a)}}
                            e_{\mathbb K}(\gamma x^k).
\]
Multiplication by each $D'$ permutes residues modulo $\mathfrak l$.
In the cofactor Poisson argument the scale-to-conductor ratio at
nonzero dual frequencies is now at least $Z/(A U)$, still a fixed
power of $X$ if $\xi+\eta<c_z$. Uniformly on the core the
zero-frequency contribution is $Xc(X)g_r(\gamma)$ times the
ordinary archimedean integral, with negligible error. The normalized
complete sum for an arbitrary linear twist in
\eqref{eq:nf-complete-bound} gives the uniform estimate
\[
 |g_r(\gamma)|\ll_\varepsilon
                  (N\mathfrak q)^{-1/k+\varepsilon}.
\]
Indeed, expand the class indicator in its $A$ additive characters
with coefficients of modulus $1/A$. A twist whose conductor does
not divide $\mathfrak q$ vanishes after averaging on fibres
modulo $\mathfrak q$; all the others are covered by that bound.
Consequently the $s$-fold product conductor series and its tail
are uniformly absolutely bounded, independently of all the classes,
since $s>2k$. The core argument in Section~\ref{sec:nf-global}
now gives an error $O(X^{s-k-\delta'})$ for some fixed $\delta'>0$.
The archimedean factor is unchanged by the residue conditions.

Write $m=v_{\mathfrak p}(\mathfrak a)$ at a prime
$\mathfrak p\mid\mathfrak a$. If $v_{\mathfrak p}(k)=0$, fixing
$s-1$ coordinates arbitrarily in their prescribed classes allows
Hensel lifting in the remaining unit coordinate, with its prescribed
class. Thus the local density with the indicated classes is at
least $(N\mathfrak p)^{-m(s-1)}$. At a prime dividing $k$, put
$v=v_{\mathfrak p}(k)$ and fix the other coordinates modulo
$\mathfrak p^{m+2v+1}$ around a chosen solution. Strong Hensel
gives the needed root in its original class; the lower bound has
only the additional factor $(N\mathfrak p)^{-(2v+1)(s-1)}$.
All other primes have their usual factors with product bounded
below on the primitive-soluble targets. The whole product is
therefore $\gg A^{-(s-1)}$. Decrease $\eta$ so that
$\eta(s-1)<\delta'$. The main term dominates the error.

For the last assertion apply Bessel's inequality to the
complementary integral with $t$ variables. Two powers of one
pointwise bound and H\"older on $S$ of the remaining powers give
$O(X^{2t-k-\delta''})$ with fixed $\delta''>0$, for $\eta'$
sufficiently small. The same core series has a uniformly convergent
tail because $t>2k$. For each prescribed tuple the number of
exceptions to half its uniform main lower bound is
$O(A^{2(t-1)}X^{k-\delta''})$, after decreasing $\eta'$ for the
core error as well. There are $A^t$ tuples for each ideal and
$O_{\mathbb K}(X^{2\eta'})$ ideals in this range. A final decrease
of $\eta'$ leaves a power saving after summing over all of them.
\end{proof}

\subsection{Roots with prescribed numbers of prime ideal factors}
\label{sec:nf-almost-primes}
Let $S$ be an even moment order supplied by the linear argument.
On a balanced target region let $X$ denote the root volume, so that
$|N\nu|\asymp X^k$. Say that $\nu$ admits $s$ unit roots at
$\mathfrak p$ if $\nu=\sum_{i=1}^s x_i^k$ with
$x_i\in\mathcal O_{\mathfrak p}^{\times}$.

\begin{theorem}\label{thm:nf-almost-primes}
For fixed $\mathbb K,k$ there are $\ell\ge1$ and $c>0$ with the
following properties.
\begin{enumerate}
\item If $s\ge S+1$ and $\nu$ admits $s$ unit roots at all finite
places, every sufficiently large totally positive $\nu$ is a sum
of $s$ $k$th powers of totally positive algebraic integers whose
principal ideals are squarefree products of exactly $\ell$ principal
prime ideals of residue degree one. All these prime norms exceed
$X^c$. The local condition is a finite union of congruence classes.
\item Suppose that for every $\mathfrak p\mid k$, of residue
characteristic $p$ and absolute ramification index $e_{\mathfrak p}$,
\[
 k\ge c_{\mathfrak p}=v_{\mathfrak p}(k)
                +\lfloor e_{\mathfrak p}/(p-1)\rfloor+1 .
\]
If $s\ge\max(S+1,H_{\mathbb K}(k))$, then every sufficiently large
totally positive eligible $\nu$ has such a representation with every root ideal a
squarefree product of exactly $\ell'$ prime ideals. All but a fixed
number of these primes, independently of $\nu$, have residue degree
one, are principal and have norm at least $X^c$. The remaining primes
belong to a finite set depending only on $\mathbb K,k,s$.
This precision condition holds for every sufficiently large $k$ for
each fixed $\mathbb K$.
\end{enumerate}
At $t\ge\max(5,S/2+1)$ the first conclusion holds for all but
$O_A(X^k(\log X)^{-A})$ targets admitting $t$ unit roots everywhere,
in each compact balanced target region and for every fixed $A>0$.
\end{theorem}
\begin{proof}
In \eqref{eq:nf-terminal-form} require the final cofactor $z$ to
be a prime element and give it weight $\log|Nz|/\log(CX)$, where
$C$ is fixed so that this lies in $[0,1]$. The $\ell-1$ other marked
factors lie in distinct windows or are explicitly chosen distinct.
Their norms exceed $|Nz|$. The support of the terminal physical weight
$w(D'z/P)$ gives $|Nz|\asymp X/|ND'|$, bounded below by a fixed
power of $X$. It is permissible to discard primes of residue degree
greater than one. Indeed prime ideals of this type of norm at most
$Y$ number $O_{\mathbb K}(Y^{1/2}\log(2Y))$, and choosing their
generators of height at most $O(Y^{1/n})$ loses at most a fixed power
of $\log(2Y)$ by the unit theorem. All constructed roots consequently
have the asserted factorization and the same number of factors.

We need a logarithmic version of Lemma~\ref{lem:nf-bilinear-weyl}.
With $m=2k-1$ and $T=X^k/Y$ as there, let $r(t)$ be the multiplicity
of a nonzero product of the $m$ shifts before the final linear estimate.
Up to the fixed factor $(k!)^2$,
\[
 r(t)\ll (\log(2X))^C\tau_m((t)),\qquad
 \sum_{H_\infty(t)\ll T^{1/n}}r(t)^2\ll T(\log(2X))^{C'}
 \quad(t\ne0).
\]
For the first inequality factor the principal ideal, then count
generators of bounded height for each factor. Since the norms are
positive integers the logarithms of the archimedean absolute values
lie in a box of diameter $O(\log X)$. Dirichlet's unit lattice gives
the logarithmic factor. The same argument applies to the generator
$t$ in the second sum. On ideals $\tau_m^2\le\tau_{m^2}$ by pairing
the local prime factors in the two ordered factorizations. The
Dirichlet series for $\tau_{m^2}$ is $\zeta_{\mathbb K}(s)^{m^2}$.
Evaluation at $1+1/\log(2T)$ gives the required bound.

Let $V(t)$ be the product of the one-dimensional linear-sum
majorants at the last differencing step. Thus $V(t)\ll Y$ and
$\sum_t V(t)\ll YM(\log(2X))^C$, in the notation of
Lemma~\ref{lem:nf-bilinear-weyl}. Cauchy's inequality bounds
$\sum_{t\ne0}r(t)V(t)/(TY)$ by
$(\log(2X))^{C'}(M/T)^{1/2}$. On the complement of
$\mathfrak M((\log X)^B)$ the full-rank and line estimates give
$M/T\ll Y^{-1}+T^{-1}+(\log X)^{-B}$. The same two-Cauchy
argument permits an arbitrary final coefficient of modulus at most one
in $z$. Consequently for every $A$ a sufficiently large $B$ gives
$\sup |F_{\rm pr}|\ll X(\log X)^{-A}$ outside these arcs.
The sharp $S$th moment of $F_{\rm pr}$ is at most that of the
unrestricted marked sum by orthogonality and positivity of
coefficients. This proves the required minor-arc estimate with
arbitrary logarithmic saving.

On the core, put $Q=N\mathfrak q\le(\log X)^B$ and
$\alpha=\gamma+\beta$. Define
\[
 U(\gamma)=\frac1{\phi(\mathfrak q)}
   \sum_{z\in(\mathcal O/\mathfrak q)^\times} e_{\mathbb K}(\gamma z^k),
 \qquad c(X)=\sum_{D'} |ND'|^{-1}\asymp1 .
\]
The prime element theorem of Mitsui, in the uniform formulation
\cite[Theorem~1.2.1]{nf3:kai}, applies to every invertible residue
of $z$ and every convex box in its balanced region of volume
$Y=X/|ND'|\ge X^{c'}$. For completeness, subdivide that region
into a rectangular mesh of relative side $(\log X)^{-C_1}$.
The smooth weight and phase $e_{\mathbb K}(\beta(D'z)^k)$ have
scaled first derivatives $O((\log X)^{C_2})$ on this core.
Piecewise-constant approximation with $C_1$ sufficiently large
gives any prescribed logarithmic accuracy. The exponentially
decreasing prime theorem error absorbs the number of boxes and
the $O((\log X)^B)$ residues. Multiplication by $D'$ permutes the
unit classes, and the physical substitution from $z$ to $D'z/P$
has constant Jacobian $X/|ND'|$. Thus with a positive field
constant $C_{\mathbb K}$, uniformly for each $A'$,
\[
 F_{\rm pr}(\alpha)=
     C_{\mathbb K}\frac{X c(X)}{\log(CX)}\,
       U(\gamma)\int w(v)e_{\mathbb K}(P^k\beta v^k)\,dv
          +O_{A',B}(X(\log X)^{-A'}).
\]
The prime ideal marked generators require no equidistribution law
modulo the growing conductor.

At $\mathfrak p\nmid k$ the normalized unit complete sum at a
primitive character is $O_k((N\mathfrak p)^{-1/2})$ at exponent
one, by the finite-field bound, and vanishes at exponents at least
two. At $\mathfrak p\mid k$, put $v=v_{\mathfrak p}(k)$.
For a character of conductor exponent $j\ge2v+2$, translation on a
unit class by multiples of $\pi^{j-v-1}$ leaves a nontrivial linear
character of the last residue digit. All higher terms are divisible
by $\pi^j$ since $2(j-v-1)\ge j$. This proves vanishing there too.
Consequently at every $s>4$ the unit singular series converges
absolutely with a power-saving tail. At large primes its local factor
is $1+O_k((N\mathfrak p)^{1-s/2})$. At the remaining primes an
existing unit solution persists on a fixed congruence ball by
Hensel's lemma, with local density bounded below uniformly on
soluble classes. These are finite congruence conditions. They hold
simultaneously, for example, on a sufficiently fine class of
$\nu\equiv s$. The smooth archimedean factor is uniformly positive
on the balanced region by the same choice of physical root boxes as
in Theorem~\ref{thm:nf-absorb-transfer}. This gives a main term
$\gg X^{s-k}(\log X)^{-s}$. The preceding arbitrary logarithmic
saving proves the first conclusion. Squaring the complementary
integral and applying Bessel with $2t-2\ge S$ gives the stated
exceptional-set estimate.

For the second conclusion include in a fixed finite bad set all
primes dividing $k$ and all primes small enough that the all-unit
finite-field Weil estimate in $s-1$ variables does not suffice.
At each bad prime take a primitive solution to the target. Record
which coordinates are nonunits. These finitely many patterns
suffice for all targets. For coordinate $i$ let $I_i$ be the
squarefree product of the bad primes on its list. Choose a fixed
integer $M$ greater than the number of bad primes. Multiply each
$I_i$ by fresh distinct large prime ideals, using principal,
narrowly positive prime generators for all but the last; choose
the last prime in the inverse narrow ideal class of $I_i$.
The prime ideal theorem in narrow classes supplies such primes
of residue degree one, distinct across all coordinates and outside
the bad set. This gives $I_i'$ having exactly $M$ prime factors
and a totally positive generator $b_i$. Use root $b_iD'z$
and retain the smooth weight in the physical root variable. This
adds exactly $M$ fixed factors to every root. For large $X$ their
norms are below all the moving prime windows.

At a bad prime, a recorded nonunit can be replaced by a number
$b_i u_i$ of valuation one with $u_i$ a unit. Its $k$th power
and the original nonunit power agree modulo $\pi^{c_{\mathfrak p}}$
because $k\ge c_{\mathfrak p}$; take $c_{\mathfrak p}=1$ if
$\mathfrak p\nmid k$. The other recorded roots are of the form
$b_i u_i$ with $b_i$ a unit at that prime. In particular one
original unit can absorb the resulting change modulo
$\pi^{c_{\mathfrak p}}$ by Lemma~\ref{lem:precision-widths}. The
resulting all-$u_i$-unit solution is nonsingular over the local
field. At a fresh prime dividing $I'_i$ only coordinate $i$
has become a nonunit, and the $s-1$ remaining unit inputs solve
the residual equation by Weil and Hensel, for every target.
Every other prime is covered by the original good-prime estimate.
The preceding argument applies to these fixed coefficients:
multiplication of the frequency by $b_i^k$ preserves Haar measure
and changes the conductor and archimedean core cutoffs by only
bounded factors. The physical weights give the same positive
archimedean density up to fixed Jacobians. The finite choice of
patterns gives uniformity in the target. This proves the assertion
with $\ell'=\ell+M$. The precision condition follows eventually
for fixed field since $c_{\mathfrak p}=O_{\mathbb K}(1+\log k)$.
\end{proof}

\begin{corollary}\label{cor:nf-prime-congruences}
Fix $D>0$. Let $A=N\mathfrak a\le(\log X)^D$ and prescribe
invertible classes $b_i\bmod\mathfrak a$. Suppose the target
has an all-unit solution at every finite place, with the indicated
root classes at the primes of $\mathfrak a$. For $s\ge S+1$ the
representations in part~$1$ of Theorem~\ref{thm:nf-almost-primes}
satisfying these congruences number at least
\[
       c_{\mathbb K,k,s,D}X^{s-k}(\log X)^{-s}A^{1-s}.
\]
For $t\ge\max(5,S/2+1)$, the same bound with $t$ variables holds
simultaneously for all these compatible moduli and classes, except
for $O_{D,E}(X^k(\log X)^{-E})$ balanced targets for every $E>0$.
\end{corollary}
\begin{proof}
All moving primary, auxiliary and cofactor primes have norms larger
than $A$. In the bilinear estimate for a congruence-restricted sum
split the outer product into its at most $A$ invertible classes.
Each imposes a single class on the separate inner variable. Since
$A$ is polylogarithmic, the pointwise estimate still saves an
arbitrary power of the logarithm on enlarging $B$ in the core.
The even moment again decreases under the coefficient restriction.
In the prime-element theorem on the core use the invertible
classes modulo $\operatorname{lcm}(\mathfrak a,\mathfrak q)$.
After multiplication by $D'$, their ordinary unit Haar mass is
restricted to the stated classes for the physical roots. The needed
uniformity is precisely the polylogarithmic-modulus case of the
same prime-element theorem.

At $\mathfrak p^j\Vert\mathfrak a$ a primitive character of
conductor exponent $h\ge\max(j+v+1,2v+2)$ integrates to zero
on the prescribed unit class by the translation on the last digit,
where $v=v_{\mathfrak p}(k)$. Thus the conductor part supported
above $\mathfrak a$ has norm $O_{\mathbb K,k}(A)$: at tame places
$h\le j$, and the fixed wild places cost only a bounded factor.
The sum of absolute values contributed by these denominators and
their characters is $O_{\mathbb K,k,\epsilon}(A^{2+\epsilon})$
by the trivial bound on the unit integrals. The conductor tail away
from $\mathfrak a$ still saves a fixed power, so any logarithmic
precision is available by increasing $B$.

At a prime of $\mathfrak a$ prime to $k$, integrating over the
$s-1$ other root classes and lifting the last unit with unit derivative
contributes at least $(N\mathfrak p)^{-j(s-1)}$ to the unit-normalized
local density. At the finitely many wild primes start from the
assumed local solution, holding the other $s-1$ roots in subballs
of depth $j+2v+1$. Hensel's lemma then supplies the last unit in its
prescribed class, at cost $O_{\mathbb K,k,s}(1)$ in the lower
bound. Primes off $\mathfrak a$ contribute a uniform positive
product as above. This proves the lower bound for the main term.
The marked factorization multiplicities are bounded by a constant
depending only on the fixed windows, and all weights are bounded,
so the weighted count gives the asserted lower bound for the
ordinary number of representations. In the almost-all assertion
the number of choices of ideals and class tuples is
$O_{\mathbb K,t,\epsilon}((\log X)^{D(t+1+\epsilon)})$.
Choose the initial logarithmic saving larger by this fixed power
and sum the Bessel exceptions.
\end{proof}

\subsection{Almost all targets and a dyadic density formula}
\label{sec:nf-power-density}

Let $\mathcal C$ be a bounded Jordan measurable region of positive
volume in $V$, with closure contained in the positive archimedean
cone and away from zero at every embedding. By relative density
we mean the limit, as $P\to\infty$, of the proportion of represented
targets in $J_k\cap P^k\mathcal C$. We require the roots to satisfy
the same signs and norm bound as in $G_{\mathbb K}(k)$. Write
$A_{\mathbb K}(k)$ for the least number of roots giving density one
in every such region.

\begin{theorem}\label{thm:nf-density-linear}
There is an integer $s_0(k)\le(C_1/2+o(1))k$, independent of the
number field, such that
\[
                 A_{\mathbb K}(k)\le\max\{H_{\mathbb K}(k),s_0(k)\}.
\]
For a fixed field let $d_2$ count its dyadic completions $F$ with
residue field $\mathbb F_2$ and
$\log(\mathcal O_F^\times)\cap\mathbb Q_2=4\mathbb Z_2$.
For every sufficiently large power of two $k$ and
$s_0(k)\le s<4k$, the relative density of the targets represented
by at most $s$ powers exists and equals
\begin{equation}\label{eq:nf-exact-density}
                   \left(\frac{s}{4k(1-2^{-k})}\right)^{d_2}.
\end{equation}
In particular $A_{\mathbb K}(2^a)=4\cdot2^a$ for all sufficiently
large $a$ if $d_2>0$. When $d_2=0$ the same threshold is at most
$(C_1/2+o(1))2^a$.
\end{theorem}
\begin{proof}
Take the restricted-moment order $r$, excess $E$ and ordinary
allocation \(B_k=\lceil16\sqrt{k}\rceil\) from the proof of
Theorem~\ref{thm:nf-optimized-moment}.
Reserve the strict ordinary compensation margin before fixing the
windows, so that $E$ includes their common initial-scale loss.
Round the ordinary allocation up to an even value
$2b\ge B_k$ with $2b/[k(k-1)]>E$ and $2b\le B_k+2$.
Let $s_0$ be the least integer for which $2s_0-2b\ge r$ and
$s_0>2k$. We may enlarge $s_0$ to be at least $3k$; this does not
affect the stated upper bound. The absorption proof in
Theorem~\ref{thm:nf-absorb-transfer} gives the loss-free full
restricted moment of order $2s_0$. With $X=P^n$, the minor arcs and
the broad annulus outside fixed cores of size $U$ therefore satisfy
\[
 \int_{\mathfrak m_*\,\cup\,(\mathfrak M_*\setminus\mathfrak M(U))}
        |f|^{2b}|F|^{2(s_0-b)}
       \ll (o(1)+U^{-\tau})X^{2s_0-k}.
\]
Here $f$ is an ordinary positive weight selected for the directed
target box, and $F$ is supported in a sufficiently small dilation.
As in the absorption proof, use an ordinary weight dominating both
supports for the even-moment inequality. For $s\ge s_0$ insert
$s-s_0$ additional ordinary variables. The same bound holds with
$2s$ in the exponent and the appropriate powers of the ordinary
weights. Bessel's inequality bounds the sum of squares of the
off-core Fourier coefficients by this integral.

For the core main term use Section~\ref{sec:nf-global} with $s$
variables. Its product series is absolutely convergent since
$s>2k$. The singular integral is bounded below on each element
of a finite cover of $\mathcal C$. If $s\ge H_{\mathbb K}(k)$
the singular series has a uniform positive lower bound. Any target
with zero representation count must then have an off-core coefficient
of size $\gg X^{s-k}$. Its exceptional set has upper density
$O(U^{-\tau})$ after letting $X\to\infty$. Letting $U\to\infty$
gives density one. This argument also works on any collection of
targets on which the finite exceptional local factors have positive
uniform lower bounds.

In fact the number of exceptions to the first bound in any fixed
compact positive box is $O(X^{k-\omega})$ for some $\omega>0$.
Use the growing core $U=X^\xi$ of Section~\ref{sec:nf-global} and
make the auxiliary moment losses smaller than the strict minor-arc
gain. Both off-core squared integrals then have power savings, while
the product-series and singular-integral approximation has a
power-saving error. A fixed uniform lower bound for the singular
series now suffices for the asserted count by Bessel.

For the more precise assertion take $k=2^a$. At odd primes
$h_k\le k+1$. At dyadic primes of residue degree at least two
\eqref{eq:nf-higher-residue-dyadic} gives $h_k\le2k+O_{\mathbb K}(1)$.
The other dyadic widths are given exactly by
Theorem~\ref{thm:nf-exact-dyadic}; only the $d_2$ specified places
can exceed $s$. Fix one such place, write $A=\mathcal O_F$,
$\Lambda=\log(A^\times)$ and $J=J_k(A)$, and let $\pi$ be a uniformizer.
Modulo $k\Lambda$ the local target group is cyclic of order $4k$.
Representations with at least one unit root fill exactly the $s$
cosets whose class is one of $1,\ldots,s$. Let their union be $D$.
Representations containing only nonunit roots have values in
$\pi^k S_s$, where $S_s$ is the local set of all $s$-term sums,
with zero terms allowed. The two parts are disjoint because
$\pi^k A\subseteq k\Lambda$. Thus
\[
 S_s=\{0\}\cup\bigcup_{j\ge0}\pi^{kj}D,
 \qquad
 \frac{\operatorname{meas}S_s}{\operatorname{meas}J}
   =\frac{s}{4k}\sum_{j\ge0}2^{-kj}.
\]
The sets of different index are disjoint by applying the same
decomposition after division by $\pi^k$. For a target in
$\pi^{kj}D$ there is a local representation with a root of valuation
exactly $j$. At bounded $j$ the Hensel density bound of
Section~\ref{sec:nf-global} applies uniformly with derivative valuation
$v_{\mathfrak p}(k)+(k-1)j$.

For completeness, $J_k\otimes\mathbb Z_p=\prod_{v\mid p}J_k(\mathcal O_v)$.
Indeed the closure of $\mathcal O$ in the local product is the full
product. Approximate a local root in one factor and zero in the
others, take its power, and take additive spans on the compact
lattices, which contain $k!\mathcal O_v$ by polarization.
Finite congruence subsets at the $d_2$ places therefore have product
relative density in $J_k\cap P^k\mathcal C$, by lattice-point counting
in a dilated Jordan region. Truncate the displayed unions at
$j\le J_0$. On these classes the uniform Bessel argument proves
that all but a set of density zero are represented. The possible
additional classes have total upper density
$O_{\mathbb K,k}(2^{-kJ_0})$, and all locally excluded targets are
absent. Let $J_0\to\infty$ to obtain \eqref{eq:nf-exact-density}.
For a box one also gets an error $O(X^{k-\omega'})$ in the
count corresponding to \eqref{eq:nf-exact-density}. Take
$J_0=\lfloor\eta\log X\rfloor$ with $\eta>0$ small. The Hensel
lower bounds at these levels cost at most $X^{O_{\mathbb K,k}(\eta)}$
in Bessel's inequality, and the discarded classes have total
density $O(2^{-kJ_0})$. After increasing a fixed constant, all the
congruences involved are unions of classes of the rational modulus
$2^{C_{\mathbb K,k}(1+J_0)}$ in $J_k$. Tiling a dilated box by
fundamental regions of this scaled lattice gives relative boundary
error $O_{\mathcal C,\mathbb K,k}
(2^{C_{\mathbb K,k}(1+J_0)}P^{-k})$. Choosing $\eta$ sufficiently
small leaves a power saving in all three errors.
If $s=4k$ every completion has primitive solutions, giving
density one directly. Taking $s=4k-1$ proves the equality when
$d_2>0$. In the remaining fields the first bound of the theorem
and the local estimates prove the last assertion.
\end{proof}

The quadratic pair after \eqref{eq:nf-kummer-log} gives
$A_{\mathbb Q(\sqrt2)}(2^a)=4\cdot2^a$ for all sufficiently large $a$,
whereas $A_{\mathbb Q(\sqrt3)}(2^a)\le(C_1/2+o(1))2^a$.
More generally, for every $e\ge1$ the polynomial $X^e-2$ is
Eisenstein at two. Thus $\mathbb K=\mathbb Q(2^{1/e})$ has a 
single dyadic completion, and Corollary~\ref{cor:pure-dyadic-width}
shows that its logarithmic lattice meets $\mathbb Q_2$ in
$4\mathbb Z_2$. For each of these fields \eqref{eq:nf-exact-density}
therefore has exponent $d_2=1$, and $A_{\mathbb K}(2^a)=4\cdot2^a$
for every sufficiently large $a$.

For each fixed field the same theorem gives density one at every
$s\ge s_0(k)$ for all sufficiently large exponents other than powers
of two. To see this, use the finite-sequence argument of
Theorem~\ref{thm:local-signature} with a ratio threshold strictly
between $2$ and $3$. The proof of
Corollary~\ref{cor:local-signature-values} shows that the only
limiting entry above $2$ is $4$, at $p=2$, $q_v=2$, $h=1$,
$b_{v,1}=2$. On a sequence $k=m2^a$ with $m$ odd its contribution
to the limiting ratio is $4/m$, by
\eqref{eq:local-signature-sequence}. Hence only $m=1$ can exceed
that threshold. Since $s_0(k)\ge3k$, the density-one assertion
follows from Theorem~\ref{thm:nf-density-linear}. Consequently
\[
 \limsup_{k\to\infty}\frac{A_{\mathbb K}(k)}{k}
 =4\quad(d_2>0),\qquad
 \limsup_{k\to\infty}\frac{A_{\mathbb K}(k)}{k}
 \le C_1/2\quad(d_2=0).
\]

\begin{corollary}\label{cor:nf-ordinary-density}
For every fixed $\mathbb K,k$ and $s\ge s_0(k)$ the relative density
of the eligible targets represented by at most $s$ powers equals
\[
        \prod_{v<\infty}
        \frac{\operatorname{meas}S_{s,v}}{\operatorname{meas}J_k(\mathcal O_v)},
 \qquad
 S_{s,v}=\{x_1^k+\cdots+x_s^k:x_i\in\mathcal O_v\}.
\]
Only finitely many factors differ from one. Each is a rational number
computable from finite residue rings. In a fixed balanced box the
count has a power-saving error. Consequently
\[
 A_{\mathbb K}(k)\le\max\{W_{\mathbb K}(k),s_0(k)\},\qquad
 A_{\mathbb K}(k)=W_{\mathbb K}(k)
                    \quad\text{if }W_{\mathbb K}(k)\ge s_0(k).
\]
\end{corollary}
\begin{proof}
At the finitely many primes of bad reduction or small residue norm
write $A=\mathcal O_v$, $q=|A/\pi A|$. Choose a precision
$c\ge1$ with $1+\pi^cA\subseteq(A^\times)^k$, as in
Lemma~\ref{lem:precision-widths}. Let $D$ be the local sums with at
least one unit root. This is the full inverse image of the primitive
sums in $A/\pi^c A$, by adjusting one unit coordinate. Therefore
\[
       S_{s,v}=\{0\}\cup\bigcup_{j\ge0}\pi^{kj}D .
\]
Indeed, divide every nonzero root tuple by the largest common power
of $\pi$ in its coordinates. The displayed sets may overlap.

For an integer $m\ge0$, let $b_m$ be the proportion of units $u$
for which $\pi^{m-kj}u\in D$ for at least one $0\le j\le m/k$.
The test uses only $u\bmod\pi^c$. If $m\ge c$, then
$b_{m+k}=b_m$: the only extra test is $\pi^{m+k}u\in D$,
whose truth for every unit, or for none, agrees with that of
$\pi^m u\in D$. With $\operatorname{meas} A=1$, summing by
valuation gives the local factor numerator explicitly as
\[
 (1-q^{-1})\left(\sum_{m=0}^{c-1} q^{-m} b_m+
              \frac{\sum_{m=c}^{c+k-1}q^{-m} b_m}{1-q^{-k}}\right).
\]
Its denominator is $[A:J_k(A)]^{-1}$. At the remaining sufficiently
large residue fields the good-prime estimate gives primitive local
solubility for every target.

A nonzero represented local target of valuation $m$ has a root
of valuation at most $m/k$ in every representation. Hence the
uniform Hensel lower bound for the corresponding local density
costs at most $q^{O_{k,A}(1+m)}$. Truncate $m$ at every exceptional
prime. Chinese remaindering in the lattice $J_k$, followed by the
uniform Bessel argument above, proves the product formula on
the resulting finite family of classes. The upper density of the
omitted sets tends to zero with the truncation. For boxes use
$m\le\eta\log X$ with $\eta>0$ sufficiently small, as in the
proof of Theorem~\ref{thm:nf-density-linear}. The periodic unit
tests require a rational integer modulus at most
$X^{O_{\mathbb K,k}(\eta)}$, the Hensel cost is at most
$X^{O_{\mathbb K,k}(\eta)}$, and each discarded set has
power-saving density. This proves the quantitative assertion.

If $s\ge W_{\mathbb K}(k)$, all local factors equal one.
Conversely, if $s<W_{\mathbb K}(k)$ there is an eligible target
outside $S_{s,v}$ for some $v$. This local sumset is compact, so
its complement contains a nonempty open subset of $J_k(A)$.
Targets in this open subset have positive relative density and
cannot be represented globally, even if $s<s_0(k)$.
\end{proof}

\subsection{A density law in real quadratic fields}
\label{sec:quadratic-density}
The exceptional dyadic places in Theorem~\ref{thm:nf-density-linear}
admit a simple description throughout the family of real quadratic
fields. Write $C_{\mathbb K}=\limsup_k H_{\mathbb K}(k)/k$ and let
$d_2$ have the meaning in that theorem.

\begin{corollary}\label{cor:quadratic-density}
For $\mathbb K=\mathbb Q(\sqrt d)$, with $d>1$ squarefree, the
two invariants and their frequencies when ordered by fundamental
discriminant are
\[
\begin{array}{c|c|c|c}
 d\bmod 8&C_{\mathbb K}& d_2&\text{proportion of fields}\\\hline
 1&4&2&1/3\\
 2,6,7&4&1&1/4\\
 3,5&2&0&5/12 .
\end{array}
\]
For every sufficiently large power of two $k$ and
$s_0(k)\le s<4k$, put $u=s/(4k(1-2^{-k}))$.
The mean relative density of the represented targets, over real quadratic
fields in discriminant order, is
\begin{equation}\label{eq:quadratic-density-polynomial}
                       \frac{5}{12}+\frac{u}{4}+\frac{u^2}{3}.
\end{equation}
Here the relative density within each field is taken first. For instance,
taking $s=7k/2$ and then $k\to\infty$ gives $57/64$.
\end{corollary}
\begin{proof}
At $d\equiv1\pmod8$ there are two copies of $\mathbb Q_2$, and
at $d\equiv5\pmod8$ there is a single unramified quadratic completion.
The latter has cyclotomic entry $2$ from $h=3$ by
Section~\ref{sec:local-cyclotomic}. At all odd places the limiting
entries are at most $3/2$, by Corollary~\ref{cor:local-signature-values}
and its proof.

For odd $d\equiv3\pmod8$ let $F=\mathbb Q_2(\sqrt{-5})$, the
dyadic completion. Since $\log(\sqrt{-5})=(\log 5)/2$,
its logarithmic lattice meets $\mathbb Q_2$ in a group containing
$2\mathbb Z_2$. Conversely if $\log v=z\in\mathbb Q_2$ with $v$ a
unit, then $2z=\log N_{F/\mathbb Q_2}(v)\in4\mathbb Z_2$.
Thus $b_{v,1}=1$. For $d\equiv7\pmod8$ the completion is
$\mathbb Q_2(i)$, and adjoining the eighth roots of unity gives a
ramified quadratic extension. For even $d$, let $F=\mathbb Q_2(\sqrt d)$.
The three quadratic subfields of $F(i)$ are all ramified: they are
generated by $i$, $\sqrt d$ and $\sqrt{-d}$. In particular $F(i)$
does not contain the unramified quadratic extension of $\mathbb Q_2$
and is totally ramified. Criterion~\eqref{eq:local-signature-four}
gives $b_{v,1}=2$ in these last two cases. The dyadic entries then
give the first two numerical columns.

The six squarefree classes $1,2,3,5,6,7$ modulo eight have equal
asymptotic frequency by M\"obius inversion on odd square divisors.
If the discriminant cutoff is $B$, the radicand cutoff is $B$ in
classes $1,5$ and $B/4$ in the other four. Normalization gives
the last column. The squarefree counts have error $O(\sqrt B)$ by the same inversion, so the finite-field-family mean differs from the limiting mean by $O(B^{-1/2})$, uniformly in $s$ and $k$ in the stated range. There are only finitely many dyadic completions
of quadratic fields up to isomorphism, so the local cutoff for the
density formula at $k=2^a$ can be taken absolute in this family.
At odd places $h_k\le k+1$. The summand threshold $s_0(k)$ is
already independent of the field. Averaging
\eqref{eq:nf-exact-density} with the displayed frequencies gives
\eqref{eq:quadratic-density-polynomial}. The analytic cutoff for
the targets in each field may depend on that field. Indeed, in a
real quadratic field of discriminant $\Delta$ the two conjugates
of every integer differ by an integral multiple of $\sqrt\Delta$.
A fixed rectangle with disjoint coordinate intervals contains no
such integer once $\sqrt\Delta$ exceeds the largest possible
coordinate difference. Thus the family average concerns the
limiting density within each field. Finally
$7/2>C_1/2$, so $7k/2$ is in the stated range for all sufficiently
large $k$, and $u\to7/8$.
\end{proof}

\subsection{Multiquadratic fields and the density at dyadic exponents}
\label{sec:multiquadratic-density}
There is a different limiting behavior as the number of independent
quadratic generators increases. Fix $r\ge2$. Choose $r$ ordered real
quadratic fields, each of fundamental discriminant at most $B$, with
equal weights, and condition on their compositum $K$ having degree
$2^r$. Here the cutoff applies to the $r$ individual discriminants.
In each resulting field first take the density of the represented
eligible targets as in Theorem~\ref{thm:nf-density-linear}.
Define
\begin{align*}
 P_0&=(1/3)^r,\\
 P_1&=(5/12)^r+4(3/8)^r-5(1/3)^r,\\
 P_2&=(1/2)^r-(5/12)^r-2(3/8)^r+2(1/3)^r .
\end{align*}
\begin{corollary}\label{cor:multiquadratic-density}
For every sufficiently large $k=2^a$ and $s_0(k)\le s<4k$, put
$u=s/(4k(1-2^{-k}))$. As $B\to\infty$, the mean density in this
family tends to
\[
  1-\sum_{j=0}^2 P_j\bigl(1-u^{2^{r-j}}\bigr).
\]
The error in this mean at a finite cutoff is $O_r(B^{-1/2})$.
In particular the limiting proportion of tuples with a surviving
dyadic obstruction is
\[
   (1/2)^r+2(3/8)^r-2(1/3)^r\sim2^{-r}
                              \qquad(r\to\infty).
\]
For real biquadratic fields in this ordering the mean density is
$(199+5u+52u^2+32u^4)/288$.
\end{corollary}
\begin{proof}
Write $V=\mathbb Q_2^\times/\mathbb Q_2^{\times2}$. Its basis
$-1,2,5$ follows from the valuation and the criterion that an odd
unit is a square precisely when it is $1$ modulo $8$. The local
squareclasses chosen by the $r$ radicands span a subspace $U\subset V$.
At each dyadic place the completion is
$F_U=\mathbb Q_2(\sqrt u:u\in U)$ and the number of dyadic places
is $2^{r-\dim U}$. If $5\in U$, then $F_U$ has even residue degree
and contributes nothing to $d_2$. Conversely a nontrivial unramified
part contains $\mathbb Q_2(\sqrt5)$. Thus $5\notin U$ is equivalent
to total ramification, including the trivial extension.

For $U=0$ the completion contributes to $d_2$. At dimension one the
quadratic calculation above gives the five exceptional lines generated
by $-1,2,-2,10,-10$. At dimension two, every plane avoiding $5$
has the form
\[
            U=\langle -5^i,2\cdot5^j\rangle,
                         \qquad i,j\in\{0,1\}.
\]
If $i=1$, the field contains $\sqrt{-5}$ and does not contain $i_0=\sqrt{-1}$.
Its next cyclotomic step adjoins $i_0$, which also adjoins $\sqrt5$;
this step is unramified. If $i=0$, $j=1$, the field contains $i_0$
and $\sqrt{10}$, and adjoining $\sqrt2$ is again the unramified
quadratic step. For the remaining plane
$U=\langle-1,2\rangle$ the field is $\mathbb Q_2(\mu_8)$.
The next step to $\mathbb Q_2(\mu_{16})$ is totally ramified, as
follows from the Eisenstein polynomial for a primitive sixteenth
root minus one. Criterion~\eqref{eq:local-signature-four} therefore
makes this the unique exceptional plane. No space of dimension
three is exceptional.

For one radicand the weights of the eight classes, in the order
$1,5,-1,-5,2,-2,10,-10$, are
\[
                 (8,8,2,2,1,1,1,1)/24 .
\]
Indeed the two odd classes with discriminant equal to the radicand
have weight $1/3$ each. The other two odd classes have weight $1/12$
each, and dividing the even radicands into their four squareclasses
uses their odd parts modulo eight. The squarefree count in each
progression has error $O(\sqrt B)$ by the same M\"obius inversion
as above. Before conditioning, independence gives the product
weights with error $O(r B^{-1/2})$ in total variation, by summing
the one-coordinate errors. After $j$ independent choices there are
at most $2^j-1$ quadratic subfields in the preceding compositum.
Summing the conditional probabilities of a first dependency bounds
the discarded proportion by $O(2^r/B)$. Thus for fixed $r$ the
conditional mean of the bounded density differs from its limiting
mean by $O_r(B^{-1/2})$. This bound concerns the ordered generating
tuples after the within-field target limit.

The probability of span zero is $P_0$. For a line of mass $w$,
including its zero element, the probability of spanning it is
$w^r-P_0$. The exceptional lines have masses $5/12$ once and $3/8$
four times. For the exceptional plane the total mass is $1/2$, and
the masses of its three lines are $5/12,3/8,3/8$. Subtracting
all proper subspace contributions gives exactly $P_2$. The
finitely many dyadic local fields described above give a common
local threshold for $a$, and at the other places the tame primitive
bound is $k+1$. Formula~\eqref{eq:nf-exact-density} now gives
the asserted mean. At $r=2$ the three probabilities are
$1/9,13/72,5/288$.
\end{proof}
\section{The uniform linear bound and the local problem}
\label{sec:nf-global-envelope}
The analytic and local estimates give a precise comparison even when
the degree and the signature of the field are prescribed. For an
integer $d\ge1$ let
\[
 L_d(k)=\sup_{p,\,[F:\mathbb Q_p]\le d} w_k(\mathcal O_F).
\]
If $r_1+2r_2=d$, let $B_{r_1,r_2}(k)$ be the supremum of
$G_{\mathbb K}(k)$ over the fields of that signature.
\begin{theorem}\label{thm:nf-local-global-envelope}
For every $k\ge3$ one has
\[
 L_d(k)\ \le\ B_{r_1,r_2}(k)\ \le\
             \max\{L_d(k)+1,\lfloor6.873k\rfloor\}.
\]
For sufficiently large $k$, independently of the signature and degree,
the analytic term $\lfloor6.873k\rfloor$ may be replaced by
$\lfloor C_1 k-\sqrt k\rfloor$. The effective analytic entries
\(\lfloor Bk\rfloor\) in \cref{tab:nf-two-coefficient-finite}
also apply on their stated ranges.
For given $k,d$, $L_d(k)$ is computable from finitely many explicit
finite rings. Uniformly over all degrees, the existence of an absolute
linear bound for $G_{\mathbb K}(k)$ is equivalent to an absolute
linear bound for the ordinary local Waring numbers in the finite
extensions of the $p$-adic fields.
\end{theorem}
\begin{proof}
Any local field $F/\mathbb Q_p$ of degree $j\le d$ occurs as a
completion of a degree-$d$ field of the given signature. Use the
coefficient approximation in the proof of
Corollary~\ref{cor:three-local-density-degree} with the local algebra
$F\times\mathbb Q_p^{d-j}$ and an archimedean polynomial with
$r_1$ real roots and $r_2$ nonreal conjugate pairs. Choose the linear
factors over $\mathbb Q_p$ to have distinct roots away from the roots
of the factor defining $F$ (translate that factor if necessary).
Impose an Eisenstein condition at a different rational prime.
Weak approximation on the coefficients, local constancy at the
separable $p$-adic polynomial and openness of the real signature
give the required irreducible polynomial. Conversely every
completion of a degree-$d$ field has degree at most $d$.
Thus the supremum of the ordinary local widths over the fields of
this signature is exactly $L_d(k)$. To see
$W_{\mathbb K}(k)\le G_{\mathbb K}(k)$ directly, fix a finite place
$v$ and an eligible target outside the image of the sum of $s$ powers.
The image is compact, so an entire congruence class about the target is
missed. Express the target as a finite integral linear combination of
local powers and approximate each root by a global integer in this
class calculation. The resulting $\nu_0\in J_k$ meets the missing class.
All $\nu_0+t k!p^a$ with $t\in\mathbb N$ still meet it for fixed
sufficiently large $a$, where $v\mid p$. These are globally eligible,
eventually totally positive and of arbitrarily large absolute norm.
Thus $s$ cannot be the global Waring number. For the other inequality,
represent a local target minus $1$ with at most $W_{\mathbb K}(k)$
powers and append the unit $1$. This proves
$H_{\mathbb K}(k)\le W_{\mathbb K}(k)+1$.
Theorem~\ref{thm:nf-all-linear} now gives the upper bound. Taking
the supremum over $d$ proves the equivalence.

For the finite computation, each rational prime $p>(2k)^6$ is
prime to $k$. At any residue cardinality $q$ above it, the equation
$x^k+y^k+z^k=b$ has $q^2+O((k-1)^3q^{3/2})$ solutions, with
absolute error coefficient one, by the Weil bound and character
orthogonality. There is a nonzero solution vector, which lifts by
Hensel, so the primitive width is at most three. On the other hand
$L_d(k)\ge3$ already over $\mathbb Z_p$ for any $p\mid k$.
Modulo $p^2$ the value of a $k$th power depends only on its residue
modulo $p$, so the number of sums of two is at most
$p(p+1)/2<p^2$. Every rational integer residue is eligible.
It is enough, therefore, to enumerate $p\le(2k)^6$ and the
extensions of degrees at most $d$ at each such $p$. A finite
algorithm for the latter is given in \cite{NFPauRob}.

For each extension choose the precision $c$ of
Lemma~\ref{lem:precision-widths}, whether or not $c<k$.
Modulo $\pi^c$ one computes the power subgroup under addition and
the residue classes of sums with a unit root by successively adding
the finite power sets. The former is the reduction of the eligible
group and $\pi^c\mathcal O_F$ is eligible by differences of unit
powers. The shell test in Corollary~\ref{cor:nf-ordinary-density}
computes exactly which eligible targets are represented by $s$ powers:
only valuations $0\le m<c+k$ and unit parts modulo $\pi^c$
are needed, together with zero. The resulting subsets on each
nonzero valuation shell are clopen, so full relative measure there
also gives full coverage. Increasing $s$ stops no later than the
primitive bound of Theorem~\ref{thm:hhybrid}. This computes $L_d(k)$.
\end{proof}

\begin{proposition}\label{prop:comparison-fixed-degree}
Fix $k\ge2$ and put $M=H_{\mathbb Q}(k)$. The inequality
$H_{\mathbb Q}(k)\le H_{\mathbb K}(k)$ holds for every number field
if and only if some rational prime $p$ has $h_k(\mathbb Z_p)=M$
and
\[
 h_k(\mathcal O_F)\ge M
 \quad\text{for every finite extension }F/\mathbb Q_p .
\]
\end{proposition}
\begin{proof}
Such a prime gives the inequality at any completion above it.
For the converse, only finitely many rational primes can contribute
a width at least $M$. Indeed, $M\ge3$ by the pair count modulo $p^2$
at a prime dividing $k$ in the preceding proof. For even $k$, a
primitive representation of zero modulo $4$ shows $M\ge4$, since
every unit power has residue one. At sufficiently large primes
the bound of three in the preceding proof therefore suffices.
For odd $k$ the corresponding bound is two. For every nonzero
residue target the smooth Fermat curve has an affine point by
the Hasse--Weil bound \cite{NFWeilCurves}, and for zero the roots
$1,-1$ give a primitive solution. Hensel lifting applies at all
these primes.

If every rational prime attaining $M$ admits a completion of
smaller width, choose one such local field $F_p$ at each of these
primes and take $F_p=\mathbb Q_p$ at the other primes on the finite
list. Let $d$ be a common multiple of the local degrees. At each
prime prescribe $d/[F_p:\mathbb Q_p]$ copies of $F_p$. Translate the
defining polynomials of the copies so that the product is separable.
Weak approximation on its coefficients, with an Eisenstein
condition at an auxiliary prime, realizes all these local algebras
in one number field of degree $d$. In that field all widths are
less than $M$. The same construction with split local algebras at
the relevant primes realizes equality in every degree.
\end{proof}

The arbitrary-class differencing theorem supplies a historical bound
of order $k\log k$ with leading coefficient one at every fixed field.
For each such field the analytic term in
Theorem~\ref{thm:nf-local-global-envelope} eventually exceeds its local
width by Theorem~\ref{thm:local-fixed-ramification}. For increasing
degree, Theorems~\ref{thm:hnearlinear}, \ref{thm:hlargeimage}
and~\ref{thm:local-family-residue} and
Corollaries~\ref{cor:inert-residue-test} and~\ref{cor:wild-large-prime}
give uniform bounds under ramification or residue-field hypotheses.
Proposition~\ref{prop:dyadic-radical-width} treats one family with
unbounded ramification and residue field of order two. More generally,
the centred carries in Section~\ref{sec:witt-top-functions} determine
the top Witt coordinate modulo Frobenius. The cost of representing
the lower coordinates by successive positive additions remains open.

For a given exponent there is also a finite sufficient certificate
with no degree bound. Here the finite rings occur before
specialization at a uniformizer. Let $B_0\ge4k+1+\log_2 k$ be an
integer. For each $p^a\Vert k$ put $m=k/p^a$. It is enough to
examine the residue fields $F=\mathbb F_{p^f}$ with
\[
                    |F|\le\max\{(m+1)^2,p^{a+1}\}.
\]
Indeed at larger $q$ the group of $m$th powers has more than
$\sqrt q-1$ elements and spans $F$, while $f\ge a+1$.
Lemma~\ref{lem:linear-local-deviation} gives $h_k(A)\le4k+a+1$
for all corresponding ramification indices. The tame places have
$h_k(A)\le k+1$.

For each retained $F$ form the finite ring
\[
 M=W_{a+1}\big((W(F)/p^\delta W(F))[T]/T^N\big),
                 \qquad N=2\delta\sum_{j=0}^{a-1}(a-j)p^j+1.
\]
Start at zero and add the set of classes $[y^m]$ (including zero)
until its additive subgroup is covered. Write $b$ for the number
of additions required. If $b+a+1\le B_0$ in each retained ring,
then $H_{\mathbb K}(k)\le B_0$ for every number field.
To see this, subtract the constant baseline of the centered chart
from the target before the finite sumset calculation. Lift the
resulting $b$ classes to individual powers over the formal
coefficient ring. The residual belongs to the kernel of truncation,
where \eqref{eq:nf-chart} supplies the $a+1$ variable unit terms
with the fixed baseline restored. Equivalently one may write the
target minus the baseline as a sum of $b$ powers plus an element
of that kernel; the latter is a variable sum minus the baseline.
Evaluation and unit lifting give the asserted width. All tests
precede the Eisenstein specialization and so apply to every
ramification index. A failed inequality in a test ring leaves the
width of its specializations undecided. In that case the
degree-bounded procedure of Theorem~\ref{thm:nf-local-global-envelope}
still computes the exact ordinary obstruction at any specified degree.

\Needspace{12\baselineskip}
\subsection{Sharpness of the field-uniform Waring bound}
\label{sec:nf-sharpness}

\begin{corollary}\label{cor:nf-log-sharpness}
For all sufficiently large \(k\), with a cutoff independent of \(\mathbb K\),
\[
 H_{\mathbb K}(k)-1\ \le\ G_{\mathbb K}(k)\
 \le\ \max\{H_{\mathbb K}(k),\,\lfloor C_1k-\sqrt{k}\rfloor\},
 \qquad C_1=6.216216302.
\]
On any sequence of fields \(\mathbb K_j\) and exponents \(k_j\to\infty\),
\[
 G_{\mathbb K_j}(k_j)\gg k_j\log k_j
 \quad\Longleftrightarrow\quad
 H_{\mathbb K_j}(k_j)\gg k_j\log k_j .
\]
In particular, \(\sup_{\mathbb K}G_{\mathbb K}(k)=o(k\log k)\) if
and only if \(\sup_{\mathbb K}H_{\mathbb K}(k)=o(k\log k)\).
\end{corollary}
\begin{proof}
The ordinary width satisfies \(W_{\mathbb K}(k)\le G_{\mathbb K}(k)\):
a missed local congruence class contains eligible global targets of
arbitrarily large norm, as in the proof of
Theorem~\ref{thm:nf-local-global-envelope}. Representing a local target
minus one and appending \(1^k\) gives \(H_{\mathbb K}(k)\le W_{\mathbb K}(k)+1\).
The upper bound is Theorem~\ref{thm:intro-main}. Whenever the local
width reaches the analytic threshold, these inequalities determine
the global width to within one summand. Since the threshold is \(O(k)\),
the two assertions at scale \(k\log k\) follow.
\end{proof}

The logarithm in the field-uniform estimate would thus have the best
possible order if a single sequence as in the corollary attained
that scale. Whether such a sequence exists remains open. The weaker
question whether \(H_{\mathbb K_j}(k_j)/k_j\to\infty\) along any
sequence is already equivalent to failure of an absolute linear
bound for \(G_{\mathbb K}(k)\). A logarithmic lower bound cannot
hold at every large exponent, since squarefree \(k\) have
\(H_{\mathbb K}(k)\ll k\) uniformly by Theorem~\ref{thm:hhybrid}.
For the logarithmic question, Section~\ref{sec:local-field-families}
reduces the search to \(k=mp^a\) at fixed \(p,m,q_v\) and unbounded
ramification. The lower bounds in Proposition~\ref{prop:witt-affine-intervals}
and Corollary~\ref{cor:witt-affine-odd} concern sums of prescribed
blocks. A lower bound for \(H_{\mathbb K}(k)\) must hold when the
individual powers can be chosen freely.

For a fixed field even the linear scale cannot be decreased to \(o(k)\).
By Chebotarev's theorem \cite{ednark}, choose an odd prime \(p\)
splitting completely in its normal closure. At a completion \(\Q_p\), take \(k=(p-1)p^a\). Euler's theorem makes
every unit \(k\)th power \(1\bmod p^{a+1}\) and every nonunit power
vanishes modulo \(p^{a+1}\). A primitive zero therefore needs at least
\(p^{a+1}\) unit roots. This gives
\[
  1<\frac{p}{p-1}
  \le\limsup_{k\to\infty}\frac{H_{\mathbb K}(k)}k\le4
\]
by Theorem~\ref{thm:local-fixed-ramification}, and the lower bound
also holds for the global limiting upper ratio by
Corollary~\ref{cor:nf-log-sharpness}. The exact local ratio at a fixed
field is given by Theorem~\ref{thm:local-signature}.

Over \(\Q\) there is a direct global obstruction. Take \(k=2^a\),
\(a\ge3\), and positive integers \(N\equiv4k\pmod{8k}\).
Since \(v_2(N)=a+2<k\), at least one root in a representation is odd.
Modulo \(4k\) every odd \(k\)th power is one and every even one is zero,
so the number of odd roots is a positive multiple of \(4k\).
These arbitrarily large targets give \(G_{\Q}(2^a)\ge4\cdot2^a\).
Theorem~\ref{thm:nf-quartic-lower} gives the corresponding obstruction
for totally real fields with a rational dyadic completion. Thus
\[
         4\le \limsup_{k\to\infty}\frac{G_{\Q}(k)}{k}
                       \le 6.216216302 .
\]
Whether the upper endpoint is attained remains open. For a fixed field a
lower bound exceeding \(4\) for this global ratio would require a
global obstruction, since its local ratio is at most \(4\).

\clearpage
\appendix
\section{Numerical certificates and reproduction}\label{app:certificates}
The numerical certificates and their verification programs are distributed
with this article as the ancillary file
\texttt{waring-siegel-certificates.zip}. The archive is available from
the ancillary-files list on the arXiv abstract page and in the
\texttt{anc/} directory of the source download. Its SHA-256 digest is
\begin{center}\small\ttfamily
f1c22766b78fc2f98f5ad9b3c4264fdb\\[-2pt]
20a6599531d5c081fd208d89dbcbc5cf
\end{center}

The directories \path{code/ed/} and \path{code/full/} contain
the exact moment witnesses for the small exponents. The finite-ring
generators and distance checks occur in \path{code/local/} and
\path{code/carries/}. The all-degree conductor computations occur in
\path{code/fibres/}, with further menu data in \path{certificates/}.
To check the witnesses and rebuild the finite enumerations, extract
the ZIP to a writable directory and run
\[
 \text{\ttfamily python3 code/replay\_audit.py -{}-full}
\]
there with Python 3.11 or later, NumPy, SymPy 1.14 or later and a
C++20 compiler. Omitting \texttt{-{}-full} runs the rational moment
checks and the primary small-power witnesses alone, together with
selected cubic identities.

The subdirectory \path{conductor/} contains the exact data for
\eqref{eq:joint-enclosure} and the effective two-coefficient orders.
From that directory run \texttt{python3 verify.py price} to
reconstruct the affine envelopes and integrate the exact cells,
or \texttt{python3 verify.py finite} for the effective orders.
These two replays use only the Python standard library.
The analytic estimates justifying the saved parameter file
are in \cref{sec:joint-conductors,sec:joint-completion}.

The directory \path{refinements/} contains the exact data for
\eqref{eq:refined-enclosure} and for
Section~\ref{subsec:refined-small-moments}. From that directory run
\begin{quote}\small\ttfamily
python3 code/replay\_profiles.py data/aligned\_union\_tail.json
\end{quote}
to construct the two profiles and check the splice and price.
The independent price check takes the saved affine cells
\path{data/aligned_union_tail_cells.csv} and the same JSON file
as inputs to \path{code/check_price.py}. Check each decic witness
by passing \path{data/field-only-k10.json} and
\path{data/field-density-k10.json}, in separate invocations,
to \path{code/finite_search.py} with the option \texttt{--verify}.
Use \path{code/verify_ninth_allocations.py} for the ninth-power field
comparison and verify \path{data/witness-k9.json} for the classical
allocation. The accompanying README gives the full commands,
including the independent check of the classical conductor endpoints.

The archive contains a file-by-file checksum manifest and the input
data used by its verification programs.

\end{document}


\begin{abstract}
For number fields of arbitrary degree and signature we give the local and analytic estimates needed at degrees four through ten. The estimates combine a localized critical mean value, efficient differencing with principal ideal factors, and auxiliary cubic and quadratic sums. An exact parameter calculation supplies the resulting bounds for smooth moments.
\end{abstract}
\maketitle
\raggedbottom
\tableofcontents
\section{Conventions and results}\label{sec:introduction}
In this supplement write $K=\mathbb K$, $n=[K:\mathbb Q]$ and
$J_{K,k}=\sum_{x\in\mathcal O_K}\mathbb Z x^k$.
The quantity $G_{\mathbb K}(k)$ denotes the least $s$ for which every
totally positive $\nu\in J_{K,k}$ of sufficiently large absolute norm is a
sum of at most $s$ $k$th powers of totally positive integral roots with
$|N x|\ll_{K,k,s}|N\nu|^{1/k}$. A zero is allowed when padding the sum.
At complex places the roots have no sign restriction.
\begin{theorem}\label{thm:main-transfer}
The smooth-moment estimates give the upper bounds for $G_{\mathbb K}(k)$
in the first row below, for every number field $\mathbb K$.
\end{theorem}
\begin{center}
\small
\begin{tabular}{lrrrrrrr}
\toprule
Degree $k$ & 4 & 5 & 6 & 7 & 8 & 9 & 10 \\
\midrule
Smooth moments (this supplement) & 16 & 17 & 24 & 31 & 38 & 46 & 55 \\
Main article (Theorem~8.1) & 16 & 17 & 24 & 31 & 38 & 44 & 50 \\
\bottomrule
\end{tabular}
\end{center}
The marked-moment estimates in Section~8.10 of the main article give
the ninth- and tenth-power bounds in the second row.
Both arguments use the smooth moments and the bounds for primitive
local solubility proved here.
\begin{theorem}\label{thm:quartic-refinements}
One has $G_{\mathbb K}(4)\le\max\{12,H_{\mathbb K}(4)\}$, where
$H_{\mathbb K}(4)$ is the maximal primitive width over the finite
completions. The bound is fifteen if the field has no completion
isomorphic to $\mathbb Q_2$ and twelve if each dyadic completion has
residue degree at least two.
\end{theorem}

\section{Preliminaries}
\label{sec:preliminaries}
Fix a number field $K$ of degree $n$, with ring of integers
$\mathcal O=\mathcal O_K$ and different $\mathfrak D$.
Put $V=K\otimes_{\mathbb Q}\mathbb R$ and
$\mathbb T_K=V/\mathfrak D^{-1}$, with probability Haar measure.
We use the trace character $e_K(z)=\exp(2\pi i\operatorname{Tr}z)$;
the notation $E(z)$ has the same meaning. At an archimedean place
let $d_v=1$ or $2$, according as the place is real or complex.
Normalize Lebesgue measure so that $\mathcal O$ has covolume one,
and use the trace-dual measure in coefficient space.

A root box $\mathcal B(Y)$ has fixed shape and archimedean radii
comparable to $Y^{1/n}$; in particular it contains $O_K(Y)$
algebraic integers. Thus $X=P^n$ is the volume parameter for
ordinary roots of archimedean size $P$. Fixed changes in shape,
translation, or scale will be allowed. Moment exponents are measured
in the root volume $X$; auxiliary box parameters also denote volumes.

For a fixed smooth compactly supported weight define
\[
 f_k(\alpha)=\sum_{x\in\mathcal O}w(x/P)e_K(\alpha x^k).
\]
For representation integrals the weights are nonnegative product
weights supported away from zero at every embedding and positive
at the real embeddings. Let $\mathcal A(Y,R)$ consist of the nonzero
roots in $\mathcal B(Y)$ whose prime ideal divisors have norm at most
$R$, and write $g_k(\alpha;Y,R)$ for its exponential sum.
Fixed support restrictions and bounded nonnegative weights are
permitted. The zero root can be appended without changing any
permissible exponent used below.

\begin{definition}
An exponent $\lambda_s$ is \emph{permissible} for degree $k$ if,
for every $\epsilon>0$, some $\eta>0$ gives
\[
 \int_{\mathbb T_K}|g_k(\alpha;Y,R)|^{2s}\,d\alpha
 \ll_{K,k,s,\epsilon}Y^{\lambda_s+\epsilon}
 \qquad(1\le R\le Y^\eta).
\]
The excess at full order $r=2s$ is
$\Delta_r=\lambda_s-(r-k)$. We also write
$U_r=\int|g_k|^r$ when the degree and scale are fixed.
\end{definition}

Linear interpolation between adjacent moment orders always means
H\"older interpolation. A finite collection of permissible estimates
can be used with one sufficiently small positive smoothness exponent.
The constants and this exponent may depend on $K$.

\subsection{Conductor ideals and heights}
For $\gamma\in K/\mathfrak D^{-1}$ put
\[
 \mathfrak a_\gamma=\{a\in\mathcal O:a\gamma\in\mathfrak D^{-1}\},
 \qquad Q_\gamma=N\mathfrak a_\gamma.
\]
For $T\ge1$ define
\[
 A_T(\beta)=\prod_{v\mid\infty}\max(1,T^{1/n}|\beta_v|)^{d_v},
 \qquad
 \mathcal H_T(\alpha)=\min_{\gamma\in K}
              Q_\gamma A_T(\alpha-\gamma).
\]
The minimum includes all lifts of the rational centers. A bound
for the displayed expression leaves finitely many conductor classes
and bounds each archimedean error, so the minimum is attained.
For a degree-$k$ ordinary sum set $B_P=A_{X^k}$, thus
$B_P(\beta)=\prod_v\max(1,P^k|\beta_v|)^{d_v}$.
The broad arcs $\mathfrak M_*$ and narrow arcs $\mathfrak N_c$
are respectively the unions of
\[
 Q_\gamma B_P(\alpha-\gamma)\le C_KX,
 \qquad Q_\gamma B_P(\alpha-\gamma)\le X^c,
\]
where $C_K$ is a sufficiently large fixed constant and $0<c<1$.
Write $\mathfrak m_*=\mathbb T_K\setminus\mathfrak M_*$.

For fixed $0\ne c\in\mathcal O_K$ one has
$\mathcal H_T(c\alpha)\asymp_{K,c}\mathcal H_T(\alpha)$.
Indeed, multiplying or dividing a rational center by $c$ changes
its conductor norm by a factor bounded in terms of $(c)$, and
the archimedean factors change by fixed factors. Taking minima
in the two directions proves the assertion.

\begin{lemma}[Product-height Dirichlet approximation and separation]
\label{lem:g789-height}
For $T\ge1$ one has $\mathcal H_T(\alpha)\ll_K T^{1/2}$.
Moreover there is $c_K>0$ such that if
\[
 \alpha=\gamma_i+\beta_i,\qquad
 (Q_{\gamma_1}A_T(\beta_1))
 (Q_{\gamma_2}A_T(\beta_2))<c_KT,
\]
then the two lifted centers are equal: $\gamma_1=\gamma_2$.
\end{lemma}
\begin{proof}
Apply Minkowski's theorem in $V^2$ to the lattice
$\{(b,a-\alpha b):b\in\mathcal O_K,\ a\in\mathfrak D^{-1}\}$,
whose covolume is fixed independently of $\alpha$, using the bounds
\[
 |b_v|\le C_KT^{1/(2n)},\qquad
 |a_v-\alpha_vb_v|\le C_KT^{-1/(2n)}.
\]
For sufficiently large $T$ the resulting nonzero pair has $b\ne0$,
because a nonzero element of $\mathfrak D^{-1}$ has absolute norm
at least $(N\mathfrak D)^{-1}$. Set $\gamma=a/b$.
Then $Q_\gamma\le|Nb|$, and
\[
 A_T(\alpha-\gamma)
 \le \prod_v(C_KT^{1/(2n)}/|b_v|)^{d_v}
 \ll_K T^{1/2}/|Nb|.
\]
Bounded $T$ is harmless. For separation, a nonzero difference of the
two centers has norm at least
$(Q_{\gamma_1}Q_{\gamma_2}N\mathfrak D)^{-1}$, whereas
\[
 |N(\beta_1-\beta_2)|\le2^nT^{-1}A_T(\beta_1)A_T(\beta_2).
\]
Choosing $c_K$ sufficiently small proves the assertion.
\end{proof}

\subsection{Divisors and quadratic norm equations}
\begin{lemma}\label{lem:divisor-norm-counts}

For fixed $C$ and every $\epsilon>0$, the following assertions hold uniformly
when all displayed nonzero algebraic integers have archimedean size at most $X^C$.

\begin{enumerate}
\item The number of divisors $d\in\mathcal O_K$ of a fixed nonzero $a$, with
all $|d_v|\le X^C$, is $O(X^\epsilon)$.  More generally the number of ideals
of norm at most $X^C$ supported on the prime ideals dividing $(a)$ is
$O(X^\epsilon)$.
\item If $abc\ne0$, the number of solutions of
\begin{equation}\label{eq:g6-binary}
 ax^2+by^2=c,
 \qquad x,y\in\mathcal O_K,\quad |x_v|,|y_v|\le X^C,
\end{equation}
is $O(X^\epsilon)$, uniformly in $a,b,c$.
\end{enumerate}

\end{lemma}
\begin{proof}
For the first assertion, count ideal divisors and then generators.  Ratios of
generators are units of logarithmic height $O(\log X)$; Dirichlet's unit lattice
therefore contributes only a power of $\log X$.  For the stronger supported-ideal
assertion, split the supporting primes at a sufficiently large fixed norm $T$.
The finitely many primes below $T$ contribute a power of $\log X$.
For $\sigma>0$ the large-prime Euler product is bounded by
\[
 \prod_{\substack{\mathfrak p\mid(a)\\N\mathfrak p>T}}
 (1-N\mathfrak p^{-\sigma})^{-1}
 \le \exp\!\left(O\!\left(\frac{\log X}{T^\sigma\log T}\right)\right).
\]
Rankin's inequality, first choosing $\sigma$ small and then $T$ large, gives
the assertion.

To prove the second assertion, put $L=K(\sqrt{-ab})$ and
$z=ax+y\sqrt{-ab}$, so that $N_{L/K}z=ac$.
When $[L:K]=2$, the ideal $(z)$ divides $(ac)\mathcal O_L$.
The number of such ideals is $O(X^\epsilon)$ uniformly in $L$:
$[L:\mathbb Q]\le2n$, and a rational prime has at most $2n$ prime ideals above it,
so the ordinary divisor bound applied to the rational norm suffices.
For each ideal, admissible generators differ by units of logarithmic height
$O(\log X)$.  This count is uniform in $L$: non-torsion algebraic units of
degree at most $2n$ have a positive uniform lower bound for their logarithmic
height, by Northcott's theorem and Kronecker's theorem.  Packing in the
logarithmic unit space gives $O((\log X)^{2n})$ possibilities; the number of
torsion units is bounded in terms of $n$.  Recovering $x,y$ from $z$ is injective.
If $-ab=d^2$ in $K$, then $d\in\mathcal O_K$ and
$(ax-dy)(ax+dy)=ac$, so the first assertion applies directly.
The hypothesis $c\ne0$ is essential; zero right sides are handled separately below.

\end{proof}

\section{Mean values and a transfer criterion}
\label{sec:analytic-principles}\label{sec:ordinary-bundle}
\begin{lemma}[Ordinary moments and pruning]\label{lem:ordinary-estimates}
For fixed $k\ge3$, the ordinary sums of Section~\ref{sec:preliminaries} satisfy

\begin{align}
 \int_{\mathfrak m_*}|f|^{k(k+1)}
    &\ll X^{k(k+1)-k-1+\epsilon},\label{eq:ordinary-critical-bundle}\\
 \int_{\mathfrak M_*}|f|^p&\ll X^{p-k}\qquad(p>2k),\label{eq:ordinary-major-bundle}\\
 \int_{\mathfrak M_*\setminus\mathfrak N_c}|f|^p
    &\ll X^{p-k-c(p/k-2)/2}\qquad(p>2k).
       \label{eq:ordinary-pruning-bundle}
\end{align}
\end{lemma}
We prove these estimates in the following four steps.

\subsection{Packing and close pairs}

If $\mathfrak b$ is a nonzero fractional ideal, then in any translated
archimedean box of radii $R_v$,
\begin{equation}\label{eq:ordinary-packing}
 \#(\mathfrak b\cap\text{box})\ll_K
       1+\frac{\prod_vR_v^{d_v}}{N\mathfrak b}.
\end{equation}
Indeed, every nonzero difference $z\in\mathfrak b$ has
$|Nz|\ge N\mathfrak b$.  After dividing the $v$th coordinate by $R_v$,
the maximum-distance separation is at least
$(N\mathfrak b/\prod R_v^{d_v})^{1/n}$.
Packing disjoint balls in a fixed bounded box proves
\eqref{eq:ordinary-packing}, also for highly unequal radii and translated boxes.

The following close-pair observation is particularly useful.  Suppose
$|y_v|\le C P$, $a\in\mathfrak D^{-1}$, and
\[
 |\theta_vy_v+\eta_v-a_v|\le Z_v\qquad(v\mid\infty).
\]
If $X\prod Z_v^{d_v}\le c_K$ for a sufficiently small fixed $c_K$, then either
there is at most one pair $(y,a)$ or all the pairs lie on one affine $K$-line
$a=a_0+\gamma(y-y_0)$, and their number is at most
\begin{equation}\label{eq:ordinary-line}
 O_K\left(1+Q_\gamma^{-1}
       \prod_v\min\{P,Z_v/|\theta_v-\gamma_v|\}^{d_v}\right).
\end{equation}
A quotient by zero means infinity.
To prove this, take differences $(u_i,b_i)$ from one pair.
The determinants $b_1u_2-b_2u_1\in\mathfrak D^{-1}$ have conjugates
$O(PZ_v)$, hence norm less than $N\mathfrak D^{-1}$, so vanish.
Distinct pairs cannot have the same $y$ by the same norm argument.
They therefore have one slope $\gamma=b/u\in K$.
Their $y$-coordinates lie in a coset of $\mathfrak a_\gamma$ and in the
intersection of the original box with boxes of radii
$O(Z_v/|\theta_v-\gamma_v|)$; apply \eqref{eq:ordinary-packing}.
Without the small-volume hypothesis the elementary bound is
\begin{equation}\label{eq:ordinary-line-trivial}
 \#\{(y,a)\}\ll_K X\left(1+\prod_v Z_v^{d_v}\right),
\end{equation}
by counting $a$ separately for each $y$.

\subsection{The shifted kernel and the critical minor-arc moment}

Choose a smooth compactly supported $v_0$, put $H=P^{k-1}$, and define
$D_H(t)=\sum_{h\in\mathcal O}v_0(h/H)e_K(-th)$.
Poisson summation gives, for every fixed large $L$,
\[
 |D_H(t)|\ll H^n\sum_{a\in\mathfrak D^{-1}}
             \prod_v(1+H|t_v-a_v|)^{-L}.
\]
Let $Y$ be a fixed coordinate box of algebraic integers with side comparable
to $P$, so $\#Y\asymp_KX$.  Uniformly in $\eta\in V$ and
$\theta\in\mathfrak m_*$,
\begin{equation}\label{eq:ordinary-shifted-kernel}
 \frac1{\#Y}\sum_{y\in Y}|D_H(k\theta y+\eta)|\ll X^{k-2}.
\end{equation}
Here is the full counting argument.
Decompose the displayed Poisson majorant into boxes
$Z_v=z_v/H$, $z_v\ge1$ dyadic, and write $B=\prod z_v^{d_v}$.
If $B\le c_KX^{k-2}$, \eqref{eq:ordinary-line} applies to the close-pair
count $N_B$.  Apart from a single pair it gives
\[
 \frac{N_B}{B}\ll 1+
 \frac{X}{Q_\gamma
   \prod_v\max\{z_v,P^k|k\theta_v-\gamma_v|\}^{d_v}}
 \ll_K1.
\]
The last inequality follows from $\theta\notin\mathfrak M_*$:
the conductor of $\gamma/k$ is at most $k^nQ_\gamma$, and
$B_P(\theta-\gamma/k)\le B_P(k\theta-\gamma)$.
If $B>c_KX^{k-2}$, \eqref{eq:ordinary-line-trivial} gives
$N_B/B\ll X/B+X^{2-k}\ll1$, since $k\ge3$.
Summing the convergent dyadic weights and using $H^n/X=X^{k-2}$ proves
\eqref{eq:ordinary-shifted-kernel}.

Define the full Vinogradov sum over a fixed enlarged coordinate box by
\[
 F(\boldsymbol\alpha)=\sum_{x\in\mathcal B(CX)}
    e_K(\alpha_1x+\cdots+\alpha_kx^k),\qquad
 J_{u,k}(CX)=\int_{\mathbb T_K^k}|F(\boldsymbol\alpha)|^{2u}
                 d\boldsymbol\alpha.
\]
For any positive integer $u$, the shifted-kernel argument gives
\begin{equation}\label{eq:ordinary-localized-general}
 \int_{\mathfrak m_*}|f(\theta)|^{2u}d\theta
       \ll X^{k(k-1)/2-1}J_{u,k}(CX).
\end{equation}
For completeness, insert coefficients of $x,\ldots,x^{k-2}$ by periodic
Fourier kernels which equal one on every possible difference of sums of these
powers.  Their suprema cost
$X^{1+\cdots+(k-2)}=X^{(k-1)(k-2)/2}$.
After integrating those coefficients, expansion imposes
$\sigma_j=\sum_{i\le u}(x_i^j-y_i^j)=0$ for $1\le j\le k-2$.
A common translation by $z\in Y$ preserves these equations and gives
\[
 \sigma_k(\mathbf x-z,\mathbf y-z)
      =\sigma_k(\mathbf x,\mathbf y)-kz\sigma_{k-1}(\mathbf x,\mathbf y).
\]
Fourier-complete each shifted weight using
$w((x-z)/P)=\int\widehat w(\xi)e_K(\xi(x-z)/P)\,d\xi$.
The Fourier $L^1$ norm is fixed, and all $x$ now lie in one enlarged box.
Introduce $\alpha_{k-1}$ by $D_H$ with $v_0$ equal to one on all possible
$P^{1-k}\sigma_{k-1}$.
Averaging over $z$ gives the kernel
$D_H(\alpha_{k-1}+kz\theta)$ times factors of modulus one.
Apply \eqref{eq:ordinary-shifted-kernel}.
The remaining product contains $2u$ full Vinogradov sums, each with only a
linear coefficient shifted by a Fourier frequency.
Extend the coefficient domain to $\mathbb T_K^k$ and use H\"older with
exponent $2u$.  Translation of $\alpha_1$ preserves Haar measure, so every
factor is $J_{u,k}(CX)$.
The total cost is
$X^{(k-1)(k-2)/2+k-2}=X^{k(k-1)/2-1}$, proving
\eqref{eq:ordinary-localized-general}.

Wooley's number-field mean-value theorem is
\cite[Theorem~15.1]{Wooley2019}, applied to $t,t^2,\ldots,t^k$,
whose Wronskian is nonzero.  In its volume normalization it gives
\[
 J_{k(k+1)/2,k}(CX)\ll X^{k(k+1)/2+\epsilon}.
\]
Substitution in \eqref{eq:ordinary-localized-general} proves
\eqref{eq:ordinary-critical-bundle}.  This completes the localization of the number-field moment theorem to the product arcs.

\subsection{The classical Hua seeds}

The elementary field Hua estimates needed for the smooth-moment iteration are
\begin{equation}\label{eq:ordinary-Hua}
 \int_{\mathbb T_K}|f|^{2^r}\ll X^{2^r-r+\epsilon}
       \qquad(1\le r\le k).
\end{equation}
The following differencing proof also provides the initial smooth moments.
By orthogonality a bounded weighted moment is at most a fixed multiple of
the unweighted moment on a containing box.
The case $r=1$ follows from $x^k=y^k$.
Assume the result for $r=J<k$.
Repeated Cauchy on the additive lattice gives
$|f|^{2^J}\ll X^{2^J-J-1}\sum_{x,\mathbf h}
e_K(\theta\Delta_{h_1}\cdots\Delta_{h_J}x^k)$, with the usual vertex-support
restrictions; the aggregate on the right is nonnegative.
Multiply by $|f|^{2^J}$ and integrate.
The equation so counted is
$\sum_{i\le2^{J-1}}(u_i^k-v_i^k)=\Delta_{h_1}\cdots\Delta_{h_J}x^k$.
If its common value is zero, the induction hypothesis supplies
$O(X^{2^J-J+\epsilon})$ choices of $\mathbf u,\mathbf v$, and there are
$O(X^J)$ choices of $x,\mathbf h$: a zero shift is immediate, while otherwise
the nonconstant polynomial in $x$ has boundedly many roots.
If its common value is nonzero, every $h_i$ divides it, so ideal divisors and
bounded-height units give $O(X^\epsilon)$ choices for all shifts, followed by
boundedly many $x$.  Thus the total count is $O(X^{2^J+\epsilon})$ in both
cases.  Multiplying by the prefactor proves the next Hua moment.

\subsection{Complete sums, Poisson summation, and pruning}

If $\psi$ is a primitive additive character modulo an integral ideal
$\mathfrak a$, $a$ is a unit modulo $\mathfrak a$, and $b$ is arbitrary, then
\begin{equation}\label{eq:ordinary-complete}
 \left|\sum_{x\bmod\mathfrak a}\psi(ax^k+bx)\right|
       \ll_{K,k,\epsilon}(N\mathfrak a)^{1-1/k+\epsilon}.
\end{equation}
We recall a uniform local proof because the linear twists are required by Poisson.
At a good prime of residue cardinality $q$, not dividing $2k(k-1)$,
write the modulus as $\pi^e$ and split into shells $v(x)=t$.
The deep shell $t\ge\lceil e/k\rceil$ has normalized size $O(q^{-e/k})$.
Outside it, put $r=v(b)$.  If $r\ne(k-1)t$, the residual polynomial on units
has unit derivative, so all shells with residual conductor at least two vanish.
Conductor one is bounded by the finite-field Weil theorem
\cite{Weil1948}.
There is at most one critical shell $r=(k-1)t$; its stationary points are
simple and its normalized size is
$O_k(q^{-t-(e-kt)/2})=O_k(q^{-e/k})$.
Thus the constant is uniform at every good prime.
At each of the finitely many bad primes, the local van der Corput theorem
\cite[Proposition~3.3]{Cluckers2011}, applied to the $k$th derivative $k!a$,
gives $O_{K_\mathfrak p,k}(q^{-e/k})$ uniformly in the integral lower coefficients.
Chinese remaindering and the ideal-divisor bound absorb the product of
good-prime constants into $(N\mathfrak a)^\epsilon$, proving
\eqref{eq:ordinary-complete}.

For $w=\prod_vw_v$ supported in archimedean annuli, let
\[
 I_v(t,u)=\int_{K_v}w_v(z)e_{K_v}(tz^k-uz)\,dz.
\]
Uniform stationary phase and integration by parts give, for every large fixed $L$,
\begin{equation}\label{eq:ordinary-stationary}
 |I_v(t,u)|\ll(1+|t|)^{-d_v/2}
                (1+|u|/(1+|t|))^{-L}.
\end{equation}
At a real place use the nonzero second derivative on the annulus.
At a complex place the real Hessian determinant of $\operatorname{Re}(tz^k)$
is a nonzero constant times $-|t|^2|z|^{2k-4}$.
After a finite partition, the stationary-phase theorem
\cite[Section~7.7]{Hormander2003} applies uniformly in the linear coefficient;
outside the gradient range repeated integration by parts gives the second factor.

For $\mathfrak a=\mathfrak a_\gamma$, $Q=N\mathfrak a$, ideal Poisson gives
\[
 f(\gamma+\beta)=\frac XQ
  \sum_{\xi\in\mathfrak a^{-1}\mathfrak D^{-1}}
    C_{\mathfrak a}(\gamma,\xi)
    \prod_v I_v(P^k\beta_v,P\xi_v),
 \quad
 C_{\mathfrak a}(\gamma,\xi)=
    \sum_{r\bmod\mathfrak a}e_K(\gamma r^k+\xi r).
\]
Primitivity and \eqref{eq:ordinary-complete} apply uniformly in every $\xi$.
Equation~\eqref{eq:ordinary-packing}, with radii
$(1+P^k|\beta_v|)/P$, bounds the sum of the decaying Fourier factors by
$O_K(1+QB_P(\beta)/X)$.
Hence, on the broad arcs,
\begin{equation}\label{eq:ordinary-broad-point}
 |f(\gamma+\beta)|\ll
       XQ^{-1/k+\epsilon}B_P(\beta)^{-1/2}.
\end{equation}
On $QB_P(\beta)\le X^c$ the same argument, using the norm lower bound for
each nonzero dual frequency, gives for every fixed $B>0$
\begin{equation}\label{eq:ordinary-core-formula}
 f(\gamma+\beta)=XG(\gamma)V(P^k\beta)+O(X^{-B}),
 \quad
 G(\gamma)=Q^{-1}C_{\mathfrak a}(\gamma,0),\quad
 V(t)=\prod_v I_v(t_v,0).
\end{equation}
Indeed the first nonzero normalized maximum frequency is at least
$c_K(X/(QB_P(\beta)))^{1/n}$; enough integrations by parts make its entire
tail an arbitrary power of $QB_P(\beta)/X$.

There are at most $Q$ center classes of conductor $\mathfrak a$.
Changing variables $t=P^k\beta$ in \eqref{eq:ordinary-broad-point} therefore
reduces the broad $p$th moment to
\[
 X^{p-k}\sum_{\mathfrak a}Q^{1-p/k+p\epsilon}
         \int_V A_1(t)^{-p/2}\,dt.
\]
Both factors converge for $p>2k$ and sufficiently small $\epsilon$, proving
\eqref{eq:ordinary-major-bundle}.
For the part outside $\mathfrak N_c$, every broad-arc parametrization has
$QA_1(t)>X^c$.  Put $a=p/k-2>0$, $\epsilon=a/(4p)$, and $\tau=a/2$.
Use
$\mathbf1_{\{QA_1(t)>X^c\}}\le X^{-c\tau}Q^\tau A_1(t)^\tau$.
The ideal exponent becomes $-1-a/4<-1$, and
$p/2-\tau=1+p(k-1)/(2k)>1$ makes every real and complex integral converge.
This proves \eqref{eq:ordinary-pruning-bundle}.
These upper bounds allow overlapping broad arcs.

\subsection{A common representation criterion}
\label{sec:product-endgame}

We first remove the small loss from a sufficiently high smooth moment.
This is the absorption argument underlying the mixed ordinary and
smooth method of \cite{Wooley2016}.

Throughout the representation argument, $g$ denotes the unweighted
smooth sum restricted to a fixed compact root annulus away from zero
at every embedding, positive at the real places. Permissible even
moments for a containing box also bound this subsum by orthogonality.
Choose the ordinary weight nonnegative and at least one on its support.

\begin{lemma}[Absorption of a critical moment]
\label{prop:bootstrap-here}
Put $w=k(k+1)$. Suppose that $u<S<w$ are even integers, $S>2k$, and
\[
 U_u=\int_{\mathbb T_K}|g|^u\ll X^{u-k+\Delta+\epsilon},
 \qquad D=(S-u)-(w-S)\Delta>0.
\]
Then $U_S\ll X^{S-k}$.
\end{lemma}
\begin{proof}
Choose a nonnegative ordinary weight which is at least one on the
support of $g$. Orthogonality gives
$U_S\le\int |f|^{S-u}|g|^u$.
On $\mathfrak m_*$ apply H\"older with weights
$(S-u)/w$, $(w-S)/w$, and $u/w$ to the moments of orders
$w$, $u$, and $S$, respectively. On $\mathfrak M_*$ use weights
$(S-u)/S$ and $u/S$, and the ordinary major-arc $S$th moment.
Writing $T=X^{-(S-k)}U_S$, Lemma~\ref{lem:ordinary-estimates} gives
\[
 T\ll X^{-D/w+\epsilon}T^{u/w}+T^{u/S}.
\]
Choose $\epsilon<D/(2w)$. Both powers of $T$ are smaller than one,
and the assertion follows.
\end{proof}

The following form permits either a single smooth moment or an
interpolation of several moments. It also separates the analytic
number of variables from a larger local requirement.

\begin{theorem}[Transfer from smooth mean values]
\label{thm:general-transfer}
Put $w=k(k+1)$. Let $a,b$ be positive integers and put $N_0=a+b>2k$.
Suppose that an even integer $S>b$ satisfies $S>2k$ and
\[
 U_S\ll X^{S-k}.
\]
Let $r_1,\ldots,r_t$ be even positive integers with permissible
excesses $\Delta_1,\ldots,\Delta_t$. Suppose that nonnegative
weights $\theta_1,\ldots,\theta_t$ satisfy
\begin{gather}
 \frac a{w}+\sum_i\theta_i=1,\qquad
 \sum_i r_i\theta_i=b,\qquad
 \delta=\frac a{w}-\sum_i\theta_i\Delta_i>0,
 \label{eq:general-minor-criterion}\\
 p=\frac{aS}{S-b}>2k.\label{eq:general-annulus-criterion}
\end{gather}
If $N\ge N_0$ and every element of every local power order has a
primitive representation by at most $N$ $k$th powers, then
$G_{\mathbb K}(k)\le N$.
\end{theorem}

We prove the theorem by replacing the smooth roots with finite
products of smaller ordinary roots. The same construction will also
be used with a joint ordinary--smooth moment for fifth powers.

\subsubsection{Finite products and the major arcs}
Fix an integer $\ell>1/\eta$ and nonnegative smooth weights $w_1,\ldots,w_\ell$
in small annuli of $V$, positive at all real places.  Define
\[
 h(\alpha)=\sum_{y_1,\ldots,y_\ell\in\mathcal O_K}
 \prod_{j=1}^{\ell}w_j(y_j/P^{1/\ell})
 e\!\left(\operatorname {Tr}\bigl(\alpha(y_1\cdots y_\ell)^k\bigr)\right).
\]
Each tuple supplies one $k$th-power root $y_1\cdots y_\ell$.
Every prime ideal dividing this root has norm at most
$C_KX^{1/\ell}<X^\eta$ for sufficiently large $X$.
The support can be placed inside the annulus used to define $g$.

For any prescribed $\epsilon>0$, the coefficient $b(x)$ of a nonzero root $x$
in $h$ satisfies $0\le b(x)\ll X^\epsilon$.
Indeed, the ideals $(y_1),\ldots,(y_\ell)$ factor $(x)$, so the number of ideal
tuples is bounded by its $\ell$-fold ideal divisor function, which is $O(X^\epsilon)$.
For each ideal tuple the ratios of two admissible generators are units with
all archimedean absolute values bounded above and below by fixed constants;
there are only $O_{K,\ell}(1)$ choices.  Consequently, for every fixed even $r$,
\begin{equation}\label{eq:product-moment-domination}
 \int|h|^r\ll X^\epsilon U_r.
\end{equation}
The same argument applies to joint even moments with ordinary factors.
The symbol $\epsilon$ can be made arbitrarily small by applying the divisor bound
with a correspondingly smaller exponent.  Real moments follow by interpolation.
The order of choices matters.  First choose the error exponent in the given
smooth-moment theorem small relative to the strict minor-arc saving, then fix a
permissible $\eta$ and prove the loss-free bootstrap moment.  Next choose
$\ell>1/\eta$ and $c<1/\ell$.  Finally choose the exponent in the independent
ideal-divisor bound small relative to $c$.  Thus the annular argument does not
require the potentially circular conditions $\eta\ll\epsilon\ll\eta$.

There is a direct major-arc formula for $h$.  Put
\begin{align*}
 H_\ell(\gamma)&=Q_\gamma^{-\ell}
 \sum_{\mathbf y\bmod\mathfrak a_\gamma}
 e\!\left(\operatorname {Tr}(\gamma(y_1\cdots y_\ell)^k)\right),\\
 V_\ell(t)&=\int_{V^\ell}\prod_j w_j(z_j)
 e\!\left(\operatorname {Tr}(t(z_1\cdots z_\ell)^k)\right)d\mathbf z.
\end{align*}
For any $0<c<1/\ell$ and any fixed $B>0$, Poisson summation in the $\ell$ factors gives
\begin{equation}\label{eq:product-core}
 h(\gamma+\beta)=XH_\ell(\gamma)V_\ell(P^k\beta)+O(X^{-B})
 \qquad(Q_\gamma B_P(\beta)\le X^c).
\end{equation}
Here $B$ is arbitrary because all factor weights are smooth.  To see the uniformity,
split each factor into residue classes modulo $\mathfrak a_\gamma$.
In the resulting Poisson formula, each nonzero dual coordinate belongs to
$\mathfrak a_\gamma^{-1}\mathfrak D^{-1}$ and has norm at least $c_K/Q_\gamma$.
The normalized dual radii are $(1+P^k|\beta_v|)/P^{1/\ell}$.
Their product has volume $B_P(\beta)/X^{1/\ell}$, so every nonzero mode has
normalized maximum radius at least
$c_K(X^{1/\ell}/(Q_\gamma B_P(\beta)))^{1/n}$.
Integration by parts once a dual coordinate exceeds the corresponding phase
gradient, followed by ideal-lattice packing, gives an arbitrary power of
$Q_\gamma B_P(\beta)/X^{1/\ell}$.
The trivial bound for the complete mixed sum is sufficient for these nonzero modes.

The pure complete factors satisfy
\begin{equation}\label{eq:product-complete}
 |H_\ell(\gamma)|\ll_{K,k,\ell,\epsilon}Q_\gamma^{-1/k+\epsilon}.
\end{equation}
For a prime-power conductor $\mathfrak p^E$, condition on
$T=\sum_{j<\ell}v_\mathfrak p(y_j)$.
The probability of a given $T$ is at most
$C_\ell(T+1)^{\ell-2}q^{-T}$, and the last-factor normalized $k$th-power sum
is at most $C_{K,k}q^{-\max(E-kT,0)/k}$.
Summing gives $C_{K,k,\ell}(E+1)^{\ell-1}q^{-E/k}$.
The finitely many primes dividing $k$ or the different are absorbed in the constant;
Chinese remaindering and the ideal divisor bound prove \eqref{eq:product-complete}.
The archimedean factors are rapidly decreasing: on the supports, differentiation
in any one factor has nonvanishing gradient of size comparable to $|t_v|$.

\paragraph{The representation main term and its normalization.}
For $a+b=N$, put $L=a+b\ell$ and let
\[
 \mathcal F(\mathbf x,\mathbf y)
   =\sum_{i=1}^a x_i^k+
      \sum_{j=1}^b(y_{j,1}\cdots y_{j,\ell})^k.
\]
Orthogonality makes
$R(\nu)=\int_{\mathbb T_K}f(\alpha)^ah(\alpha)^b e_K(-\alpha\nu)\,d\alpha$
the nonnegative weighted count of these representations.
The core arcs are disjoint for sufficiently large $X$.
Indeed, distinct centers of conductor norms $Q,Q'$ satisfy
$|N(\gamma-\gamma')|\gg_K(QQ')^{-1}$, whereas an overlap forces
$|N(\gamma-\gamma')|\ll_K X^{2c-k}(QQ')^{-1}$.
Here $2c<k$, since $c<1/\ell<1$ and $k\ge4$.

With the minor-arc and annular savings of Theorem~\ref{thm:general-transfer}, the two core Poisson formulae give
\[
 R(\nu)=X^{N-k}\mathfrak S(\nu)\mathfrak J(\nu/P^k)
            +o(X^{N-k}),
\]
uniformly for normalized targets in the compact set used below, where
\[
 \mathfrak S(\nu)=\sum_{\gamma\in K/\mathfrak D^{-1}}
      G(\gamma)^aH_\ell(\gamma)^b e_K(-\gamma\nu),
 \qquad
 \mathfrak J(t)=\int_V V_f(\beta)^aV_\ell(\beta)^b e_K(-\beta t)\,d\beta.
\]
Here $V_f(t)=\int_V w(z)e_K(tz^k)\,dz$ is the ordinary archimedean
factor, and all measures have the normalization fixed in
Section~\ref{sec:ordinary-bundle}.  The conductor majorant is summable
because $N>2k$; the archimedean factors decay rapidly.  These same majorants justify extending the truncated core
series and integral, with a uniform error tending to zero.
Chinese remaindering and finite-ring character orthogonality give
\[
 \mathfrak S(\nu)=\prod_{\mathfrak p}\sigma_{\mathfrak p}(\nu),\qquad
 \sigma_{\mathfrak p}(\nu)=\lim_{E\to\infty}
 q^{-E(L-1)}\#\{(\mathbf x,\mathbf y)\bmod\mathfrak p^E:
                       \mathcal F(\mathbf x,\mathbf y)\equiv\nu\},
 \quad q=N\mathfrak p.
\]
For completeness, a primitive local representation gives a target-uniform
positive lower bound.  Replace each product root $x$ by $(x,1,\ldots,1)$.
At a unit coordinate the derivative has valuation $h=v_{\mathfrak p}(k)$.
Fix the other $L-1$ variables modulo $\mathfrak p^{2h+1}$ at that exact solution.
Every further choice of those variables lifts in the distinguished coordinate
by Hensel's lemma.  Thus for every $E\ge2h+1$ the normalized count is at least
$q^{-(2h+1)(L-1)}$, independently of the target.
At all sufficiently large primes the complete-sum estimates give
$\sigma_{\mathfrak p}(\nu)=1+O(q^{1-N/k+\epsilon})$.
Choose $\epsilon<N/k-2$.  This summable good-prime tail, combined with the
positive lower bounds at the finitely many remaining primes, makes
$\mathfrak S(\nu)$ positive uniformly in $\nu$.

Finally take $P=|N\nu|^{1/(kn)}$ and balance $\nu$ using $k$th powers of
totally positive units, for example squares of units.  Their logarithmic
lattice has full rank, so the normalized targets lie in a compact set with
every embedding bounded away from zero.  Small root boxes near equal $k$th
roots, positive at the real places, give a uniformly positive singular
integral on a finite cover.  Undoing the totally positive unit preserves
the signs of all real roots and preserves their absolute norms.
Every resulting root has norm $O_K(P^n)=O_K(|N\nu|^{1/k})$, exactly as required
by the stated convention for $G_{\mathbb K}(k)$.

\begin{proof}[Completion of the proof of Theorem~\ref{thm:general-transfer}]
For $N=N_0$, the representation integral is $f^ah^b$.
By \eqref{eq:product-moment-domination} and
\eqref{eq:general-minor-criterion}, its minor arcs contribute
\[
 \int_{\mathfrak m_*}|f|^a|h|^b
 \ll X^{N_0-k-\delta+\epsilon}.
\]
On $\mathfrak M_*\setminus\mathfrak N_c$, use the ordinary
$p$th moment and the product $S$th moment, with weights
$1-b/S$ and $b/S$. The saving in the exponent of $X$ is
\[
 \frac c2\left(1-\frac bS\right)\left(\frac pk-2\right)>0.
\]
The independent divisor loss in the product coefficients is chosen
smaller than this saving. The core main term is positive by the
local and archimedean arguments above.

For $N>N_0$, insert $N-N_0$ additional ordinary roots. Their trivial
bound $O(X^{N-N_0})$ preserves both off-core savings, while the core
calculation has $N$ summands and uses the assumed primitive local
representations in that number of variables. This proves the theorem.
\end{proof}

\begin{corollary}[An even-variable criterion]
\label{cor:even-transfer}
Suppose that $u<S<w$ are even, $S>2k$, and
\[
 U_u\ll X^{u-k+\Delta+\epsilon},\qquad
 (S-u)-(w-S)\Delta>0.
\]
Under primitive local solubility in $N\ge S$ variables,
$G_{\mathbb K}(k)\le N$.
\end{corollary}
\begin{proof}
Lemma~\ref{prop:bootstrap-here} supplies $U_S\ll X^{S-k}$.
In Theorem~\ref{thm:general-transfer} take $a=S-u$, $b=u$,
and the two smooth moments of orders $u,S$, with weights
$(w-S)/w$ and $u/w$. The minor-arc saving is $D/w$ and the
ordinary annular order is $p=S$.
\end{proof}

\section{Primitive local representations}
\label{sec:local-solubility}\label{sec:higher-local-bounds}
For the integer ring $A$ of a finite extension of $\mathbb Q_p$,
let $J_k(A)$ be the additive group generated by its $k$th powers.
A representation is primitive if one root is a unit, including
when the target is zero. Let $h_k(A)$ denote the least number of
summands in such representations of every element of $J_k(A)$.

\begin{theorem}\label{thm:local-table}
For every such ring $A$ the following bounds hold:
\[
\begin{array}{c|rrrrrrr}
 k&4&5&6&7&8&9&10\\\hline
 h_k(A)\le&16&10&12&14&35&30&17.
\end{array}
\]
For fourth powers the bound is $5$ at odd primes, $10$ at dyadic
primes of residue degree at least two, and $15$ at every other
dyadic extension except $\mathbb Q_2$.
\end{theorem}

The case of a prime exponent follows from Joly's theorem
\cite[Theorem~7.34]{Joly1970}: the ordinary width is at most
$2k-1$. Apply it to $\nu-1$ and append $1^k$ to obtain the
primitive bounds $h_5\le10$ and $h_7\le14$. We prove the other
entries below.
\subsection{First-order carries and squarefree exponents}
\label{sec:valuation-one-carries}

Let $A$ be the integer ring of a finite extension of $\mathbb Q_p$,
with residue field $\mathbb F_q$, uniformizer $\pi$, and absolute
ramification index $e=v_\pi(p)$. Let $J_k(A)$ denote the additive
group generated by the $k$th powers in $A$. Write $h_k(A)$ for the
least number of positive power summands that represents every member
of $J_k(A)$ with at least one unit root. This requirement includes
the target zero.

\begin{theorem}\label{thm:carry-sixth-tenth}
For the integer ring of every finite extension of every $\mathbb Q_p$,
\[
 h_6(A)\le12,\qquad h_{10}(A)\le17.
\]
At places of residue characteristic five the second bound improves
to $h_{10}(A)\le16$.
\end{theorem}

We prove the assertion through an explicit valuation-one carry
argument. Its finite-field input is the usual $m$-term Waring bound
of Tornheim, also recorded in Joly, Proposition~7.16, and in the
proof of Ramanujam, Lemma~1 \cite{Joly1970,Ramanujam1963}.
The proof below keeps a unit coordinate throughout and treats all
ramification and residue degrees.

\subsubsection{Finite-field compression and primitive lifting}

Suppose $p\nmid m$, and put
\[
 H=(\mathbb F_q^\times)^m,
 \qquad E=\mathbb F_p(H).
\]
The additive span of $H$ is $E$: that span is a finite subring of
$\mathbb F_q$ containing one. Every member of $E$ is a sum of at
most $m$ members of $H\cup\{0\}$.
For completeness, the latter bound also follows from Kneser's
addition theorem. The additive stabilizer of a proper sumset
$r(H\cup\{0\})$ is stable under $H$, hence under $E$, and therefore
is zero. Until the whole field $E$ is reached, its cardinality is
at least $r|H|+1$. Taking $r=(|E|-1)/|H|\le m$ proves the claim.

A unit target in $J_m(A)$ consequently has a primitive representation
by at most $m$ $m$th powers: lift the finite-field representation
and use its unit coordinate in Hensel's lemma. Subtracting one first
gives a primitive representation of every target in $J_m(A)$ using
at most $m+1$ terms. If $m\ge3$ is odd, $m$ terms suffice, since
the zero residue already has the primitive two-term representation
$1+(-1)^m=0$.

We also record the exact lifting precisions needed below. For odd
$p$ and $k=pm$, a primitive power congruence modulo $p^2A$ lifts
without increasing its number of summands. At $p=2$, $k=2m$ with
$m$ odd, the same holds from modulus $8A$. Indeed
\begin{equation}\label{eq:carry-primitive-lifts}
 \frac{(1+pz)^{pm}-1}{p^2}\equiv mz\pmod\pi
 \quad(p>2),\qquad
 \frac{(1+4z)^{2m}-1}{8}\equiv mz\pmod\pi.
\end{equation}
Both normalized polynomials are integral and have a unit derivative.
Multiplying an existing unit root by $1+pz$ or $1+4z$ and applying
Hensel's lemma absorbs any error in the indicated ideal.

\subsubsection{The exact carry and its first-order cancellation}

Write $k=pm$ with $m\ge2$ and $p\nmid m$. For a finite list $t$ put
\[
 C_p(t)=\frac{(\sum_i t_i)^p-\sum_i t_i^p}{p}.
\]
This polynomial has integer coefficients, and
\begin{equation}\label{eq:carry-derivative}
 \frac{\partial C_p}{\partial t_i}
       =(\sum_jt_j)^{p-1}-t_i^{p-1}.
\end{equation}
When $p$ is odd, start with a signed expression
$\nu_0=\sum_i\epsilon_i x_i^{pm}$ and put $t_i=\epsilon_i x_i^m$.
Then $t_i^p=\epsilon_i x_i^{pm}$. In the dyadic case we use an
unsigned expression and put $t_i=x_i^m$; the availability of such
expressions is proved below.

Compress $S=\sum_i t_i\in J_m(A)$ exactly as
$S=\sum_{j=1}^r y_j^m$, with a unit $y_j$, and put $s_j=y_j^m$.
The bound is $r\le m+1$ in general and $r\le m$ when $m\ge3$ is odd.
Equality of the two sums gives the exact identity
\begin{equation}\label{eq:carry-exact}
 \nu_0-\sum_{j=1}^r y_j^k=pD,
 \qquad D=C_p(s)-C_p(t)\in A.
\end{equation}
Every residue involved belongs to $E$, so $\bar D\in E$.
Choose $h_1,\ldots,h_l\in H$, $l\le m$, with
$\sum_i h_i=\bar D$, allowing the empty sum, and set
\begin{equation}\label{eq:carry-constant}
 c=p\sum_{i=1}^l[h_i].
\end{equation}
Here brackets denote the ordinary Teichm\"uller lifts in $A$.
Frobenius permutes $H$, so each $[h_i]$ is a $k$th power of a
Teichm\"uller unit. Repeating it $p$ times expresses $c$ using at
most $pm=k$ positive powers, and $D-c/p\in\pi A$.

If $H\subseteq\mathbb F_p^\times$ and $e\ge2$, this improves to
\begin{equation}\label{eq:carry-second-order}
 D-c/p\in\pi^2A.
\end{equation}
To prove it, identify
$A/\pi^2A\simeq\mathbb F_q[\varepsilon]/(\varepsilon^2)$,
which has characteristic $p$. Denote its $\varepsilon$ coefficient
by a prime. A nonunit $t_i$ or $s_j$ is divisible by $\pi^m$ and
has zero first derivative. A unit has residue in
$\mathbb F_p^\times$, including signed $t_i$ when $p$ is odd, and
therefore its $(p-1)$st power is one. Using
\eqref{eq:carry-derivative} and the exact equality $\sum t_i=\sum s_j$
gives
\[
 D'=\sum_i t_i'-\sum_j s_j'=0.
\]
The constant coefficient of $D$ belongs to $\mathbb F_p$.
Since the lifts $[h_i]$ in \eqref{eq:carry-constant} lie in
$\mathbb Z_p$ and $p\in\pi^2A$, $c/p$ has the same complete class
modulo $\pi^2$. This proves \eqref{eq:carry-second-order}.
Notice that this is a statement about the explicit carry difference;
it does not assert $pA\subseteq J_k(A)$.

\subsubsection{Odd residue characteristic}

If $e=1$, compression and the constant correction leave an error
in $p\pi A=p^2A$. Equation~\eqref{eq:carry-primitive-lifts}
therefore gives the bound
\begin{equation}\label{eq:carry-odd-unramified-count}
 h_k(A)\le k+m+1.
\end{equation}

If $H\nsubseteq\mathbb F_p^\times$, choose Teichm\"uller
$m$th-power values $a,b$ with $a/b$ outside $\mathbb F_p$.
Reserve their two unit roots before compressing the remaining target.
Put $A_0=a^p$, $B_0=b^p$, and $\lambda=-a/b$; then
$A_0+B_0\lambda^p=0$. The polynomial
\[
 F(Z)=\frac{A_0(1+Z)^p+B_0(1+\lambda Z)^p-A_0-B_0}{p}
\]
is integral because the top and constant terms cancel. Its derivative
at zero is $B_0(\lambda-\lambda^p)$, a unit. After compression and
\eqref{eq:carry-constant}, solve $F(z)=D-c/p$ with $z\in\pi A$.
Replace the two underlying $m$th-power values by $a(1+z)$ and
$b(1+\lambda z)$; principal-unit $m$th roots exist because $p\nmid m$.
This absorbs the whole remaining error and gives
\begin{equation}\label{eq:carry-odd-pair-count}
 h_k(A)\le k+m+3.
\end{equation}

Now suppose $H\subseteq\mathbb F_p^\times$ and $e\ge2$.
Reserve four unit roots equal to one. After compression and
\eqref{eq:carry-constant}, the error is $p\pi^2R$ by
\eqref{eq:carry-second-order}. Put $z_4=-z_1-z_2-z_3$ and define
\[
 \mathcal F(z_1,z_2,z_3)
   =\frac{\sum_{i=1}^4(1+\pi z_i)^p-4}{p\pi^2}.
\]
The linear term vanishes. The intermediate binomial coefficients
are divisible by $p$, and $\sum_i z_i^p$ is divisible by $p$ as
an integral polynomial because $\sum_i z_i=0$.
Consequently $\mathcal F$ is integral and reduces to
\[
 -\frac12\sum_{i=1}^4z_i^2
 =-\frac12 z^{\mathsf T}(I_3+\mathbf1\mathbf1^{\mathsf T})z.
\]
The matrix has determinant four. A nondegenerate ternary quadratic
form over a finite field of odd characteristic represents every
value nonsingularly. Indeed, diagonalize and use the intersection
of two square sets to represent every nonzero value in two variables;
for zero, fix the third coordinate to one and represent the negative
of its contribution in the first two. Hensel lifting solves
$\mathcal F(z)=R$. The four values $1+\pi z_i$ remain $m$th powers,
so the four reserved coordinates absorb the error and stay units.
This proves
\begin{equation}\label{eq:carry-odd-four-count}
 h_k(A)\le k+m+5.
\end{equation}

\subsubsection{The two square-compression exceptions}
\label{subsec:carry-two-exceptions}

When $(p,m)=(3,2)$ or $(5,2)$, the last count improves by one by
using an existing compressed unit. The condition
$H\subseteq\mathbb F_p^\times$ forces $q=p$: indeed
$|H|=(q-1)/2\le p-1$, whereas $p^2>2p-1$.
Reserve only three ones and compress the remaining target by at
most three squares. Choose one unit square $b=y_1^2$ among these
compressed values and put $b_1=b$, $b_2=b_3=b_4=1$.
With $z_4=-z_1-z_2-z_3$, use
\begin{equation}\label{eq:carry-existing-unit}
 \mathcal F_b(z)
   =\frac{\sum_{i=1}^4(b_i+\pi z_i)^p-\sum_{i=1}^4b_i^p}
          {p\pi^2}.
\end{equation}
Its linear part is $((b^{p-1}-1)/\pi)z_1$, which is integral.
Modulo $\pi$, its quadratic part is
\[
 \frac{p-1}{2}
   \bigl(\bar b^{,p-2}z_1^2+z_2^2+z_3^2+z_4^2\bigr).
\]
All higher terms vanish modulo $\pi$, by the same top-power
divisibility as above. The underlying three-variable matrix is
\[
 \operatorname{diag}(\bar b^{,p-2},1,1)
       +\mathbf1\mathbf1^{\mathsf T},
 \qquad \det=3\bar b^{,p-2}+1.
\]
At $p=3$, a nonzero square has $\bar b=1$, so this determinant is
nonzero. At $p=5$, $\bar b\in\{1,-1\}$, giving determinants four
and minus two, again nonzero.

Completing the square turns the reduction of
\eqref{eq:carry-existing-unit} into a nondegenerate ternary quadratic
form plus a constant. It therefore represents every target
nonsingularly, including its critical value by a nonzero isotropic
vector. Hensel lifting absorbs the error $p\pi^2R$, and every
$b_i+\pi z_i$ remains a square. The term count is now
\begin{equation}\label{eq:carry-special-count}
 3\text{ compressed terms}+3\text{ reserved terms}
       +k\text{ constant terms}=k+6.
\end{equation}
Thus $h_6(A)\le12$ at $p=3$ and $h_{10}(A)\le16$ at $p=5$.
In the latter case the constant correction still uses at most ten
terms: every residue in $\mathbb F_5$ is a sum of at most two
members of $H=\{1,-1\}$, and both Teichm\"uller values are actual
tenth powers. 

\subsubsection{Dyadic valuation one}

Let $p=2$ and $k=2m$, where $m\ge3$ is odd.
There is a unit $r\in\mathbb Z_2$ with $r^{2m}=-7$.
For instance, Hensel lifting of $z^2+z+2=0$ gives a square root
$1+2z$ of $-7$, and the odd-power equation $r^m=1+2z$ has a
unit solution. Hence
\[
 -1=6\cdot1^k+r^k.
\]
Every signed expression, including one obtained by subtracting
reserved powers, can be made unsigned before applying
\eqref{eq:carry-exact}. Its original length is immaterial.
The exact $m$th-power compression uses at most $m$ terms.

If $H\ne\{1\}$, choose a Teichm\"uller $m$th-power value $a$
whose residue differs from one, and choose $\alpha$ with $\alpha^m=a$.
Reserve $(\alpha,1)$. Replacing these roots by
\[
 \alpha(1+z)^{1/m},\qquad(1-az)^{1/m},\qquad z\in\pi A,
\]
changes their $k$th-power sum by $2F_a(z)$, where
\begin{equation}\label{eq:carry-dyadic-pair}
 F_a(Z)=(a^2-a)Z+a^2Z^2.
\end{equation}
Its derivative at zero is a unit. The correction
\eqref{eq:carry-constant} uses at most $2m$ powers and leaves
$D-c/2\in\pi A$. Hensel lifting with $F_a$ absorbs it exactly.
The number of summands is at most
\begin{equation}\label{eq:carry-dyadic-nontrivial-count}
 m+2+2m=3m+2.
\end{equation}

Suppose $H=\{1\}$ and $e\ge2$. Reserve the pairs $(1,1)$ and
$(1+\pi,1)$; for the second put $a=(1+\pi)^m$.
Equation~\eqref{eq:carry-second-order} gives
$D\equiv n\pmod{\pi^2}$ with $n\in\{0,1\}$, so $c=2n$
uses at most two ones. Since $m$ is odd and $e\ge2$,
$a\equiv1+\pi$ and $a^2\equiv1\pmod{\pi^2}$. Therefore
\[
 G(Z)=\frac{F_a(\pi Z)}{\pi^2}
     =\frac{a^2-a}{\pi}Z+a^2Z^2\equiv Z+Z^2\pmod\pi.
\]
The remaining error is $2\pi^2R$. Choose $t_0$ lifting the square
root of $\bar R$, and solve $G(t)=R-t_0^2$ with $t\in\pi A$.
The first pair, with $z=\pi t_0$ in
\eqref{eq:carry-dyadic-pair} at $a=1$, contributes
$2\pi^2t_0^2$. The second, with $z=\pi t$, contributes
$2\pi^2G(t)$. Their sum is the required error; all four roots
remain units. The count is
\begin{equation}\label{eq:carry-dyadic-trivial-count}
 h_k(A)\le m+4+2=m+6.
\end{equation}

Finally suppose $H=\{1\}$ and $e=1$. The unit $m$th-power
values are exactly $1+2A$: their Teichm\"uller factors are killed,
and the $m$th-power map on principal units is bijective. Thus the
unit $k$th-power classes modulo $8A$ are
\[
 1+4V,\qquad V=\{z+z^2:z\in\mathbb F_q\},
\]
while nonunit $k$th powers vanish at that precision since $k\ge6$.
Their generated additive image is
$(\mathbb Z/8\mathbb Z)\cdot1+4V$. Every class $a+4v$ is a
sum of $r$ unit power classes, where $1\le r\le8$ and
$r\equiv a\pmod8$: use one class $1+4v$ and $r-1$ ones.
The dyadic lifting in \eqref{eq:carry-primitive-lifts} gives
\begin{equation}\label{eq:carry-dyadic-unramified-count}
 h_k(A)\le8.
\end{equation}

\subsubsection{Conclusion for the two exponents}

For $k=6$, the odd bad prime is $3$, where
\eqref{eq:carry-special-count} gives twelve; the other odd cases
give smaller counts. The three dyadic bounds are eleven, nine,
and eight. For $k=10$, the odd bad prime is $5$, where the bound
is sixteen, and the three dyadic bounds are seventeen, eleven,
and eight. At a prime not dividing $k$, finite-field compression
of $\nu-1$, appending $1^k$, and ordinary Hensel lifting give
$h_k(A)\le k+1$. All constructions keep a unit root, so they also
give primitive representations of zero. This proves
Theorem~\ref{thm:carry-sixth-tenth}.
\subsection{Witt coefficients with a positive summand count}
\label{sec:higher-witt-coefficients}

Let $p$ be a prime, let $q$ be a power of $p$, and put
$R=\mathbb F_q[[T]]$ and $B=W(\mathbb F_q)$, the unramified complete
discrete valuation ring with residue field $\mathbb F_q$.
We use $p$-typical Witt vectors $W_{a+1}(R)$, with coordinates
$(z_0,\ldots,z_a)$, Teichm\"uller vectors $[x]=(x,0,\ldots,0)$,
and ghost polynomials
\[
 w_j(z)=\sum_{i=0}^j p^i z_i^{p^{j-i}}.
\]
Their addition and multiplication are the universal integral
polynomials making the ghost maps homomorphisms. In particular
$[x][y]=[xy]$. All constant Witt vectors form the coefficient ring
$W_{a+1}(\mathbb F_q)=B/p^{a+1}B$.

Choose positive integers $N_0,\ldots,N_a$ satisfying
$p^{j-i}N_i\ge N_j$ for $i<j$, and put
\[
 I_N=\{z:z_j\in T^{N_j}R\},\qquad
 M_N=W_{a+1}(R)/I_N.
\]
Weighted homogeneity of the Witt operations shows that $I_N$ is an
ideal: a monomial of Witt weight $p^j$ in its coordinates has
$T$-order at least $N_j$. Multiplication is separately homogeneous
in each input, giving the same assertion when only one input is in
$I_N$. Triangularity in the coordinates gives
$|M_N|=q^{\sum_jN_j}$: choose the coordinates successively modulo
$T^{N_j}$, using the coefficient of $z_j$ equal to one in the
$j$th addition polynomial.

\begin{lemma}\label{lem:higher-witt-coefficients}
Let $G$ be the actual additive span of the vectors $[x]$, $x\in R$.
Every class of $G$ modulo $I_N$ is a positive sum of at most
\begin{equation}\label{eq:higher-witt-cost}
 C_p(N)=L_0+(p-1)\sum_{s=1}^a p^{s-1}L_s,
 \qquad
 L_s=\sum_{j=s}^a\left\lceil\frac{N_j}{p^s}\right\rceil
\end{equation}
actual Teichm\"uller vectors. The bound is independent of $q$.
\end{lemma}
\begin{proof}
Write $\ell(M)=\log_q|M|$ for finite $B$-modules. In characteristic
$p$, multiplication by $p^s$ has the coordinate description
\[
 p^sz=(\underbrace{0,\ldots,0}_{s},
             z_0^{p^s},\ldots,z_{a-s}^{p^s}).
\]
Thus its image modulo $I_N$ has cardinality $q^{L_s}$.
One can count the kernel after any sufficiently large common
$T$-truncation: its coordinate conditions are exactly
$z_i\in T^{\lceil N_{i+s}/p^s\rceil}R$ for $i\le a-s$.
Consequently $\ell(p^sM_N)=L_s$.

The image $M$ of $G$ in $M_N$ is a $B$-submodule, since multiplying
an actual generator by a constant Teichm\"uller vector preserves its
form. Adjoin actual generators $x_i$ successively to obtain a
filtration $0=M_0\subsetneq\cdots\subsetneq M_r=M$ with cyclic
quotients $M_i/M_{i-1}\simeq B/p^{\lambda_i}B$.
A coefficient modulo $p^\lambda$ has a Teichm\"uller digit expansion;
regarding its digit $p^j[c]$ as $p^j$ copies of $[c]$ uses at most
$Q_\lambda=1+p+\cdots+p^{\lambda-1}$ positive scalars. Working
backwards through the filtration therefore represents a target by
at most $\sum_iQ_{\lambda_i}$ actual generators.

For every $s$, the natural surjections
$p^sM_i/p^sM_{i-1}\twoheadrightarrow p^s(M_i/M_{i-1})$ give
\[
 \ell(p^sM)\ge\sum_i\max(\lambda_i-s,0).
\]
The identity
\[
 Q_\lambda=\lambda+
   (p-1)\sum_{s\ge1}p^{s-1}\max(\lambda-s,0)
\]
then bounds the positive cost by
$\ell(M)+(p-1)\sum_{s\ge1}p^{s-1}\ell(p^sM)$.
Each length is at most the corresponding length for $M_N$, since
$p^sM\subseteq p^sM_N$. This proves
\eqref{eq:higher-witt-cost}. In particular, the argument never
replaces actual power generators by arbitrary invariant-factor
coordinates.
\end{proof}

Here is the link to local integral powers. For a local integer ring
$A$ of residue field $\mathbb F_q$, choose its unramified coefficient
ring $B\subseteq A$ and the surjection $D=B[[T]]\to A$ sending
$T$ to $\pi$. On $W_{a+1}(R)$, choosing coordinate lifts to $D$
defines an additive map
\begin{equation}\label{eq:higher-last-ghost}
 \Phi:W_{a+1}(R)\longrightarrow D/p^{a+1}D,
 \qquad \Phi(z)=\sum_{i=0}^a p^i\widetilde z_i^{p^{a-i}}.
\end{equation}
It is independent of the lifts: changing a lift by a multiple of
$p$ changes its $p^{a-i}$th power by a multiple of $p^{a-i+1}$.
Additivity follows by lifting the universal Witt addition and using
the ghost identity. In particular $\Phi([x])=\widetilde x^{p^a}$.
Lift an eligible target by any finite signed power expression to
$D$; its associated signed Teichm\"uller sum reduces to an element
of $G$.

Suppose a centered chart covers $I_N$ in the form
\[
 w=\sum_{i=1}^r[x_i]-\sum_{i=1}^r[b_i],
\]
retaining a unit among the $x_i$. Subtract the fixed baseline
$\sum[b_i]$ from the target and use Lemma~\ref{lem:higher-witt-coefficients}
modulo $I_N$. The remaining ideal element is supplied by the chart,
and the baseline cancels. The result is a primitive representation
by at most $C_p(N)+r$ powers modulo $p^{a+1}D$, and hence modulo
$p^{a+1}A$. Fixed nonzero constant residues in the chart are
preserved by choosing lifts of the roots with those constants.

\subsection{Fourth powers}
\begin{lemma}[Primitive lifting]\label{lem:primitive-lifting}
At an odd prime, a fourth-power congruence modulo $\pi A$ with
a unit root lifts to an exact representation of the same target
using the same number of summands. At a dyadic prime, the same
assertion holds for a congruence modulo $16A$.
\end{lemma}
\begin{proof}
For odd $p$, the derivative $4x^3$ is a unit at a unit root.
For $p=2$, changing a unit root $x$ to $x+4z$ changes its fourth
power by $16F_x(z)$, where
\[
 F_x(z)=x^3z+6x^2z^2+16xz^3+16z^4.
\]
Modulo $\pi$, this is the invertible linear map $z\mapsto\bar x^3z$.
Simple-root Hensel lifting solves $F_x(z)=a$ for every $a\in A$.
Thus any error in $16A$ can be corrected by that root alone.
\end{proof}

\begin{lemma}[Odd residue characteristic]\label{lem:odd-local}
If $p\ne2$, then $h_4(A)\le5$.
\end{lemma}
\begin{proof}
The prime-to-residue-characteristic case of
Ramanujam~\cite[Lemma~1]{Ramanujam1963} says that the reduction of
$J_4(A)$ is the subfield generated by fourth powers and that each
of its elements is a sum of at most four fourth powers in the
residue field. If the target is a unit, such a representation has
a nonzero root and lifts. If it is a nonunit, apply the unit case
to $\nu-1$ and append $1^4$. This also treats $\nu=0$.
\end{proof}

We use only length-three $2$-typical Witt vectors. Their standard
integral polynomial construction is described in
Hazewinkel~\cite[Sections~5.2--5.22]{hazewinkel2009}; the specific
open-image assertion below is proved directly.
\subsubsection{The three-coordinate open-image lemma}

Let $C$ be a commutative ring, put $R=C[[T]]$, and let $W_3(R)$ denote
the $2$-typical Witt vectors of length three.  Its ghost polynomials are
\[
 w_0(X)=X_0,\qquad w_1(X)=X_0^2+2X_1,\qquad
 w_2(X)=X_0^4+2X_1^2+4X_2.
\]
Write $[y]=(y,0,0)$ for the Teichm\"uller vector.  Define
\[
 I=\{z\in W_3(R):z_0\in T^3R,\ z_1,z_2\in T^4R\}.
\]
\begin{lemma}\label{lem:wittopen}
The set $I$ is an ideal, and every $w\in I$ has the form
\begin{equation}\label{eq:quartic-specialized-open-image}
 w=[1+T^4z_0]+[T+T^3z_1]+[T^2+T^3z_2]
       -[1]-[T]-[T^2],\qquad z_0,z_1,z_2\in R.
\end{equation}
In particular $I$ lies in the additive subgroup generated by
Teichm\"uller vectors, and the root in the first summand of
\eqref{eq:quartic-specialized-open-image} is $1$ modulo $T$.
\end{lemma}

\begin{proof}
Give Witt coordinate $j$ weight $2^j$.  The integral Witt addition
and negation polynomials are homogeneous of the corresponding weight;
the multiplication polynomials are bihomogeneous in the two input
vectors.  The inequalities
\[
 2\cdot3>4,\qquad4\cdot3>4,\qquad2\cdot4>4
\]
show that the specified divisibilities are preserved under addition,
negation, and multiplication by an arbitrary Witt vector.  This proves
the ideal assertion.  It may equally be checked directly from the
three displayed ghost polynomials, working universally over $\mathbb Z$
before specializing to $C$.

For the open-image assertion, work first with integral universal Witt
polynomials and put
\[
 H(\eta_0,\eta_1,\eta_2)
   =\sum_{i=0}^2[T^i+\eta_i]-\sum_{i=0}^2[T^i].
\]
All coordinates of $H$ vanish at $\boldsymbol\eta=0$.  Its ghost
coordinates satisfy
\[
 w_j(H)=\sum_{i=0}^2\bigl((T^i+\eta_i)^{2^j}-T^{i2^j}\bigr).
\]
Differentiating at zero and using the vanishing of the lower
coordinates shows that the linear part of $H_j$ is
\[
 \sum_{i=0}^2 T^{i(2^j-1)}\eta_i.
\]
This calculation takes place over $\mathbb Z$ (or first over
$\mathbb Q$ with the resulting integral Witt polynomials); it does
not divide by $2$ inside $C$.

Substitute
\[
 (\eta_0,\eta_1,\eta_2)=(T^4z_0,T^3z_1,T^3z_2)
\]
and divide the three coordinates by $T^3,T^4,T^4$, respectively.
The resulting polynomial map $\Phi:R^3\longrightarrow R^3$ is
integral and has linear part
\[
 U=\begin{pmatrix}
 T&1&1\\
 1&1&T\\
 1&T^2&T^5
 \end{pmatrix},\qquad
 \det U=(T-1)(T^2-1)(T^3-1).
\]
For completeness, integrality and the nonlinear divisibilities follow
directly from the ghost recurrence.  The first coordinate is linear.
The second is
\[
 H_1=\sum_iT^i\eta_i-\sum_{i<j}\eta_i\eta_j,
\]
whose nonlinear terms have $T$-adic order at least six.  In the third
ghost equation, every nonlinear term in
$\sum_i((T^i+\eta_i)^4-T^{4i})$ has order at least eight, while
$H_0^4$ has order at least twelve and the nonlinear contribution of
$2H_1^2$ has order at least eight.  Dividing numerical coefficients by
$4$ does not change these $T$-adic orders, and the universal Witt
coordinate polynomial is integral.  Hence
\[
 \Phi(\mathbf z)=U\mathbf z+(0,T^2P_1(\mathbf z),T^4P_2(\mathbf z))
\]
for integral polynomials $P_1,P_2$.

The constant term of $\det U$ is $-1$, so $U$ is invertible over $R$.
For an arbitrary $w\in I$, divide its coordinates by $T^3,T^4,T^4$
to obtain $\mathbf b\in R^3$.  The equation $\Phi(\mathbf z)=\mathbf b$
has a solution by $T$-adic successive approximation: the map
\[
 \mathbf z\longmapsto U^{-1}\!\left(
 \mathbf b-(0,T^2P_1(\mathbf z),T^4P_2(\mathbf z))\right)
\]
increases the precision of two congruent inputs by at least one
power of $T$.  Completeness of $R$ gives a fixed point.  This proves
\eqref{eq:quartic-specialized-open-image} after specialization to $C$.
\end{proof}

\paragraph{The quartic ghost map.}
Let $B=W(\mathbb F_q)$ be an unramified dyadic integer ring, put
$D=B[[T]]$, and take $C=B/4B$.  The last ghost coordinate induces
an additive map
\[
 \Psi:W_3(D/4D)\longrightarrow D/16D,\qquad
 \Psi(z)=\widetilde z_0^4+2\widetilde z_1^2+4\widetilde z_2.
\]
Changes of lifts by $4D$ change the expression by a multiple of $16D$,
so it is well defined; the integral ghost identities prove additivity.
Moreover $\Psi([y])=\widetilde y^4$, and
\[
 4T^4D\pmod{16D}\subseteq\Psi(I),\qquad
 2T^8\pmod{16D}\in\Psi(I),
\]
by taking vectors $(0,0,T^4a)$ and $(0,T^4,0)$, respectively.
Applying $\Psi$ to \eqref{eq:quartic-specialized-open-image} gives
the three-positive-fourth-power correction with fixed baseline
$1+T^4+T^8$ used in the local endpoint proof.

If $A$ is the integer ring of any finite extension of $\mathbb Q_2$
with residue field $\mathbb F_q$, choose its maximal unramified
coefficient ring $B\subseteq A$ and a uniformizer $\pi$.
Substitution $T\mapsto\pi$ gives a surjective continuous map
$D\to A$: every element of $A$ has a convergent $\pi$-adic expansion
with coefficients in a system of residue representatives in $B$.
Composing the ghost map with this substitution gives
\[
 \Psi_A:W_3((B/4B)[[T]])\longrightarrow A/16A.
\]
Its image on the additive subgroup generated by the vectors $[y]$
is exactly the image of $J_4(A)$ modulo $16A$, because each root in
$A$ has a power-series preimage.  The same construction with
$B=\mathbb Z_2$ will be used for residue degree one.  Constant Witt
vectors evaluate through the corresponding ghost map into $B/16B$.

\subsubsection{The local fourth-power bound}
\label{sec:universal-quartic}

We now prove the primitive local bound needed in the global argument.
The finite additive-group calculation is established explicitly below.

\begin{theorem}\label{thm:universal-quartic-sixteen}
For every finite extension $L/\mathbb Q_p$,
\[
 h_4(\mathcal O_L)\le16.
\]
The constant is best possible in the primitive convention:
$h_4(\mathbb Z_2)=16$.
In dyadic residue degree at least two the bound is $10$; in residue
degree one it is $11$ when the ramification index is two and $15$
when that index is at least three.
\end{theorem}
\begin{proof}
At an odd prime, Lemma~\ref{lem:odd-local} gives the bound $5$.
The dyadic cases of residue degree at least two, ramification index
two and residue degree one, and ramification index at least three
and residue degree one are proved in
Propositions~\ref{prop:quartic-residue-degree}, \ref{prop:quartic-ramification-two}, \ref{prop:quartic-high-ramification} and~\ref{prop:quartic-cubic-ramification} below.
The only remaining local ring is $\mathbb Z_2$. Every class modulo
$16$ is a sum of $r$ copies of $1^4$ for some $1\le r\le16$.
Lemma~\ref{lem:primitive-lifting} gives the upper bound. Conversely,
every unit fourth power is $1$ modulo $16$, while every nonunit
fourth power is zero. In a primitive representation of zero the
positive number of unit summands is therefore divisible by $16$.
This proves equality over $\mathbb Z_2$ and the theorem.
\end{proof}

\subsubsection{Fourth powers when the dyadic residue degree is at least two}
\label{sec:quartic-residue-degree}

\begin{proposition}
\label{prop:quartic-residue-degree}
Let $A$ be the integer ring of a finite extension of $\mathbb Q_2$ with
residue field $\mathbb F_q$, $q=2^f$. If $f\ge2$, then
\[
 h_4(A)\le10.
\]
More precisely, nine summands suffice for eligible targets with nonzero
residue. The estimate is uniform in the ramification index.
\end{proposition}

\begin{proof}
Let $B=W(\mathbb F_q)$ and first consider a constant target
$\lambda\in B/16B$. Write $[c]_B$ for the Teichmuller lift of
$c\in\mathbb F_q$. Every $[c]_B$ is a fourth power, since the fourth-power
map is an automorphism of the group of roots of unity of odd order.

If $\bar\lambda\ne0$, choose distinct nonzero $c_1,c_2$ with
$c_1+c_2=\bar\lambda$. If $\bar\lambda=0$, choose distinct nonzero
$c_1,c_2$ and put $c_3=c_1+c_2$. These choices are possible because
$q\ge4$, and the selected residues sum to $\bar\lambda$.
After subtracting their Teichmuller lifts, represent the remaining
target modulo $8B$ by its next two binary digits, costing at most
$2+4=6$ additional Teichmuller fourth powers. Thus there is a
representation modulo $8B$ using at most eight powers when the residue
is nonzero, or nine otherwise. It includes two unit roots whose fourth
powers have distinct nonzero residues.

The last binary digit can be corrected without adding summands. For a
unit root $u$ and $z\in B$,
\[
 (u(1+2z))^4-u^4\equiv8u^4(z^2+z)\pmod{16B}.
\]
The image $V=\{z^2+z:z\in\mathbb F_q\}$ is the trace-zero hyperplane.
For distinct nonzero $c,d$, the hyperplanes $cV$ and $dV$ are distinct:
they are the kernels of the different linear functionals
$x\mapsto\operatorname{Tr}(x/c)$ and
$x\mapsto\operatorname{Tr}(x/d)$. Hence $cV+dV=\mathbb F_q$.
Varying the two chosen roots therefore corrects any error in the last
binary digit. Their residues remain distinct and nonzero. We have
proved the same eight- or nine-term bound modulo $16B$.

Now use the Witt model $W_3((B/4B)[[T]])$ for a target
$\nu\in J_4(A)$. Choose a preimage $z$ in the additive Teichmuller span
and let $z_{\mathrm c}$ be its evaluation at $T=0$. Its top ghost is some constant
$\lambda\in B/16B$. Construct a constant Teichmuller sum $s$ from the
preceding paragraph with this same top ghost. Replacing $z$ by
$z-z_{\mathrm c}+s$ changes no image in $A/16A$, since
$z_{\mathrm c}$ and $s$ have the same top ghost modulo $16B$.
The new vector evaluates exactly to $s$ at $T=0$.

Append one zero root to $s$. Among the roots are two units with distinct
nonzero residues $u,v$, together with the zero root. The derivative
matrix of the centered three-variable Witt sum has rows with exponents
$0,1,3$; its residue determinant is
\[
 \det\begin{pmatrix}1&1&1\\u&v&0\\u^3&v^3&0\end{pmatrix}
       =uv(v^2-u^2)\ne0.
\]
The centered Witt-coordinate map has the displayed invertible
linear part.  On substituting perturbations $Tz_i$ and dividing its
coordinates by $T$, all nonlinear terms acquire an extra factor $T$.
The same successive-approximation argument as in Lemma~\ref{lem:wittopen}
therefore realizes every Witt perturbation whose coordinates lie in
$T(B/4B)[[T]]$ using these three existing roots. It corrects the entire remaining
$T$-dependent part of the target, while preserving both unit roots.
All other coordinates are kept fixed, so no further summands are
introduced. Substitution $T=\pi$ gives a primitive congruence modulo
$16A$, and Lemma~\ref{lem:primitive-lifting} gives the exact representation in $A$.
The appended zero increases the counts to nine and ten, respectively.
\end{proof}

\subsubsection{Fourth powers at absolute ramification index two}
\label{sec:quartic-ramification-two}

\begin{proposition}
\label{prop:quartic-ramification-two}
If $A$ is a dyadic integer ring with residue field $\mathbb F_2$ and
absolute ramification index $2$, then $h_4(A)\le11$.
\end{proposition}
\begin{proof}
Normalize $v(\pi)=1$, so $v(2)=2$, and write $\pi^2=2u$ with
$u\in A^\times$. We first work modulo $\pi^7A$.
A fourth power depends only on its root modulo $\pi^2A$.
Indeed, for $x,h\in A$,
\[
 (x+\pi^2h)^4-x^4
 =4x^3\pi^2h+6x^2\pi^4h^2+4x\pi^6h^3+\pi^8h^4.
\]
The last two terms have valuation at least $8$. If $x\in\pi A$,
the first two have valuation at least $8$ as well. If $x$ is a unit,
the sum of the first two terms, divided by $\pi^6$, has residue
$\bar h+\bar h^2=0$ in $\mathbb F_2$; here
$4/\pi^4=u^{-2}$ and $6/\pi^2=3u^{-1}$ both have residue $1$.
This proves the assertion.

The four root classes are $0,1,\pi,1+\pi$. Put
\[
 t=\pi^4,\qquad w=(1+\pi)^4,\qquad r=w-1.
\]
In the additive group of $A/\pi^7A$ they satisfy
\[
 16=0,\qquad 4t=0,\qquad 2t=8,\qquad 2r=0.
\]
The first two relations follow from their valuations. For the third,
$2t-8=8(u^2-1)$ has valuation at least $8$, since both $u-1$ and
$u+1$ lie in $\pi A$. For the last relation, expand
\[
 r=4\pi+6\pi^2+4\pi^3+\pi^4
   =4\pi+4u(3+u)+4\pi^3.
\]
Every term has valuation at least $5$, so $2r\in\pi^7A$.

Consequently every eligible class has the form
\[
 n+jt+\varepsilon r,
 \qquad 0\le n\le7,\quad0\le j\le3,\quad
 \varepsilon\in\{0,1\}.
\]
If $\varepsilon=1$, replace $r$ by $w-1$, reduce the coefficient of
$1$ modulo $8$, and absorb the resulting multiple of $8$ into the
coefficient of $t$ using $8=2t$. This gives at most seven copies of
$1$, three copies of $t$, and one copy of $w$. It is primitive because
$w$ has the unit root $1+\pi$.
If $\varepsilon=0$ and $n>0$, the displayed form itself is primitive
and uses at most ten powers. If $\varepsilon=n=0$, use
\[
 jt=8+\bigl(j+2\bmod4\bigr)t
\]
modulo $\pi^7A$. This uses at most eleven powers and has unit roots.
Thus every eligible class has a primitive representation in at most
eleven fourth powers.

Finally, a primitive congruence modulo $\pi^7A$ lifts with the same
number of summands. For a unit root $x$, the polynomial
\[
 \frac{(x+\pi^3 Z)^4-x^4}{\pi^7}
\]
has coefficients in $A$, has unit linear coefficient
$4x^3/\pi^4$, and all its higher coefficients lie in $\pi A$.
Hensel's lemma therefore corrects any residual error by changing that
root alone. This proves the proposition.
\end{proof}

\subsubsection{The finite quotient in residue degree one}
\label{sec:quartic-conductor}

The following finite quotient handles all ramification indices at least
three in residue degree one. Its positive width is determined explicitly.

\subsubsection{A universal finite quotient}

Put $D=\mathbb Z_2[[T]]$, $R=D/4D$, and
$W=W_3(R)$.  Let $G$ be the additive subgroup of $W$ generated by all
Teichm\"uller vectors $[y]$, $y\in R$.  The last ghost component gives
an additive map
\[
 \Psi:W\longrightarrow D/16D,
 \qquad
 (z_0,z_1,z_2)\longmapsto
 \widetilde z_0^4+2\widetilde z_1^2+4\widetilde z_2\pmod {16D},
\]
independent of the chosen lifts.  Its image on $G$ is the image of
$J_4(D)$.  Lemma~\ref{lem:wittopen} gives
\[
 I=\{z\in W:z_0\in T^3R,\ z_1,z_2\in T^4R\}\subseteq G.
\]
Write $K=\Psi(I)$.  In particular,
\begin{equation}\label{eq:quartic-conductor-basic}
 4T^4D\pmod {16D}\subseteq K,
 \qquad 2T^8\pmod {16D}\in K.
\end{equation}

Define the additive coefficient projection
\[
 \Pi:D/16D\longrightarrow
 V:=\mathbb Z/8\mathbb Z\oplus(\mathbb Z/16\mathbb Z)^3
       \oplus(\mathbb Z/4\mathbb Z)^4\oplus\mathbb Z/2\mathbb Z
\]
by reading coefficients of $T^0,\ldots,T^8$, in that order, with the
indicated moduli.  Higher coefficients are discarded.  Polynomial
notation below denotes classes in this additive group; no ring
structure on $V$ is asserted.  Set
\[
 Q=\Pi(J_4(D)).
\]

\begin{lemma}\label{lem:quartic-universal-quotient}
The kernel of $\Pi$ on the image of $J_4(D)$ in $D/16D$ is
$K+8\mathbb Z_2$.  Every class of $Q$ is a sum of at most twelve
fourth-power images, with at least one root having odd constant
coefficient.
\end{lemma}

\begin{proof}
First, $K+8\mathbb Z_2$ is killed by $\Pi$.  Indeed, for $z\in I$
its three ghost terms have orders at least $12,8,4$, respectively.
The first is discarded; the second has an even coefficient at $T^8$;
and the third has coefficients divisible by $4$ in degrees $4$ through
$7$ and an even coefficient in degree $8$.  These are all zero in $V$,
as is the constant $8$.  Conversely, given
$\alpha\in J_4(D)$, choose $z\in G$ with $\Psi(z)=\alpha\pmod {16D}$.
By successively reducing its three Witt coordinates modulo $I$, choose
a representative $w$ of the same coset with its coordinate degrees
less than $3,4,4$, respectively.  The successive reduction is performed using the triangular form of Witt addition.  Therefore
\[
 \alpha\equiv\gamma:=\widetilde w_0^4+
                 2\widetilde w_1^2+4\widetilde w_2\pmod K,
 \qquad \deg\gamma\le8.
\]
If $\Pi(\alpha)=0$, the coefficients of $\gamma$ have the indicated
zero residues.  Thus
\[
 \gamma\in8\mathbb Z_2+4T^4D+2\mathbb Z_2T^8
       \pmod {16D},
\]
which lies in $K+8\mathbb Z_2$ by
\eqref{eq:quartic-conductor-basic}.  This proves the kernel assertion.

Every root can be reduced modulo $(4,T^4)$ for this calculation.
Indeed, if $\Delta\in T^4R$, the coordinates of
$[y+\Delta]-[y]$ are
\[
 \Delta,\qquad y\Delta,\qquad
 y^3\Delta+y^2\Delta^2+y\Delta^3,
\]
so this difference lies in $I$.  Changes of a lift by $4D$ do not
change its fourth power modulo $16D$.  There are consequently only
$4^4=256$ root classes to consider.

Here is an explicit description of their additive image.  In $V$, put
\[
 h_1=8T,\qquad h_2=8T^2,\qquad h_3=4T^3+2T^6,
 \qquad H=\langle h_1,h_2,h_3\rangle.
\]
The orders of these independent generators are $2,2,4$, so $|H|=16$.
For odd constant coefficient, changing the sign of the root permits
normalization to
\[
 F(a,b,c):=\Pi\bigl((1+aT+bT^2+cT^3)^4\bigr),
 \qquad a,b,c\in\mathbb Z/4\mathbb Z.
\]
Direct expansion gives
\begin{align*}
 F(a,b,c)={}&1+4aT+(4b+6a^2)T^2
       +(4c+12ab+4a^3)T^3\\
 &+(a^4+2b^2)T^4+2(c^2+a^2b^2)T^6+b^4T^8.
\end{align*}
All coefficients here are read in $V$.  In particular,
\begin{align*}
 F(a,b,c+1)-F(a,b,c)&=h_3,\\
 F(a,b+2,c)-F(a,b,c)&=h_2+2a h_3,\\
 F(a+2,b,c)-F(a,b,c)&=h_1+(a+1)h_2+2(a+b)h_3.
\end{align*}
It follows that $H\subseteq Q$ and that the unit-power image is
stable under translation by $H$.  Its four cosets modulo $H$ are
obtained by taking $a,b\in\{0,1\}$ and $c=0$.

Set, in $Q/H$,
\[
 u=1,\qquad x=4T+6T^2+4T^3+T^4,\qquad
 \alpha=T^4,\quad\beta=T^8,\quad\zeta=2T^6.
\]
These elements are independent, of orders $8,4,4,2,2$, respectively.
For example, a relation is first inspected in its constant and $T$
coefficients; the $T$ coefficient gives $4X\in8\mathbb Z$ and hence $X$ even.
The $T^2$ coefficient then gives $6X\in8\mathbb Z$ modulo $16$,
forcing $X=0$ modulo $4$.  The $T^4,T^8,T^3,T^6$ coefficients
settle the remaining generators.  The unit-power cosets are
\begin{equation}\label{eq:quartic-unit-cosets}
 u,\qquad u+x,\qquad u+2x+\beta,\qquad u+3x+\beta.
\end{equation}
For a root with even constant coefficient, expansion gives exactly
these nonunit cosets:
\begin{equation}\label{eq:quartic-nonunit-cosets}
 0,\qquad\alpha,\qquad\beta,\qquad\alpha+\beta+\zeta.
\end{equation}
For instance they are represented by $0,T,T^2,T+T^2$.
The possible extra contribution from constant coefficient $2$ is
$8a^2T^2+8a^3T^3$, which lies in $H$.
Equations \eqref{eq:quartic-unit-cosets} and
\eqref{eq:quartic-nonunit-cosets} generate all of $Q/H$, and show that
\[
 Q/H\simeq\mathbb Z/8\mathbb Z\oplus
 (\mathbb Z/4\mathbb Z)^2\oplus(\mathbb Z/2\mathbb Z)^2.
\]

To represent a target with coordinates $(C,X,A,B,Z)$ in this product,
choose $1\le r\le8$ with $r\equiv C\pmod8$.  A unit type indexed by
$j\in\{0,1,2,3\}$ contributes
$u+jx+\lfloor j/2\rfloor\beta$.  With three unit slots one can
prescribe both the $x$ and $\beta$ coordinates: for $\beta$ coordinate
zero, take triples $(0,0,0),(1,0,0),(1,1,0),(1,1,1)$ to obtain
$x$ coordinates $0,1,2,3$; for $\beta$ coordinate one, take
$(3,1,0),(3,1,1),(2,0,0),(3,0,0)$, respectively.
Extra unit slots may have index zero.

If $r\ge3$, append $Z\in\{0,1\}$ copies of the last type in
\eqref{eq:quartic-nonunit-cosets}, and then choose
$0\le s\le3$ copies of $\alpha$ so that $s+Z\equiv A\pmod4$.
The unit types can match the required $X$ and $B-Z$ coordinates.
The total is at most $r+s+Z\le12$.
If $r=1$ or $2$, use one unit type with index $X$, and index zero
for the remaining unit slot if present.  One additional $\beta$
corrects the remaining binary coordinate if necessary.  The total
is then at most $r+3+1+1\le7$.
Finally, because the unit-power image is $H$-stable, one unit
summand absorbs any discrepancy in $H$.  This proves the lemma.
\end{proof}

\begin{proposition}\label{prop:quartic-high-ramification}
Let $A$ be the valuation ring of a finite extension of $\mathbb Q_2$
with residue field $\mathbb F_2$ and ramification index $e\ge4$.
Then
\[
 h_4(A)\le15.
\]
\end{proposition}

\begin{proof}
Choose a uniformizer $\pi$.  Substitution $T\mapsto\pi$ defines a
surjective continuous map $D\to A$, hence sends the universal
power-generated group onto the corresponding group for $A$.
The image $K_A$ of $I$ under the last ghost map and this substitution
contains $4\pi^4A$ modulo $16A$.  Since $e\ge4$,
\[
 8=4\pi^4\,(2/\pi^4)\in4\pi^4A\subseteq K_A.
\]
Consequently the quotient in Lemma~\ref{lem:quartic-universal-quotient}
surjects onto $J_4(A)/K_A$, with all quotients here taken modulo $16A$.

Given $\nu\in J_4(A)$, apply that lemma to
\[
 \nu-(1+\pi^4+\pi^8)
\]
and obtain a sum $v$ of at most twelve fourth powers with
$\nu-(1+\pi^4+\pi^8)-v\in K_A$.
Lemma~\ref{lem:wittopen} supplies, for every element
of $I$, an identity
\[
 w=[y_0]+[y_1]+[y_2]-[1]-[T]-[T^2],
 \qquad y_0\equiv1\pmod T.
\]
Applying the last ghost component and substituting $T=\pi$ cancels
the baseline.  Thus
\[
 \nu\equiv v+y_0(\pi)^4+y_1(\pi)^4+y_2(\pi)^4\pmod {16A}
\]
is a representation by at most fifteen fourth powers containing a
unit root.  Lemma~\ref{lem:primitive-lifting} completes the proof.
\end{proof}

\begin{proposition}\label{prop:quartic-cubic-ramification}
If $A$ has residue field $\mathbb F_2$ and ramification index $e=3$,
then $h_4(A)\le15$.
\end{proposition}

\begin{proof}
For any uniformizer $\varpi$, the unit $\varpi^3/2$ has a cube root
in $A$, by simple-root Hensel lifting.  We may therefore choose a
uniformizer $\pi$ with $\pi^3=2$.
Repeat the universal coefficient calculation with constant coefficient
modulo $16$, instead of modulo $8$.  Denote the resulting group by
$Q_{16}$.  The proof of the kernel assertion now gives kernel $K$,
and the same subgroup $H$ and the same power-image lists give
\[
 Q_{16}/H\simeq\mathbb Z/16\mathbb Z\oplus
 (\mathbb Z/4\mathbb Z)^2\oplus(\mathbb Z/2\mathbb Z)^2,
\]
with generators $u,x,\alpha,\beta,\zeta$ as before; only the order
of $u$ changes.  On substituting $T=\pi$ one has
\[
 \zeta=2\pi^6=8=8u.
\]
Thus every specialized target may be written with $u$ coefficient
$0\le C\le7$, while retaining a binary $\zeta$ coefficient.
Choose $1\le r\le8$ with $r\equiv C\pmod8$.  When $C=0$ and
$r=8$, the extra $8u$ simply flips the required $\zeta$ coordinate.
The rest of the positive construction in
Lemma~\ref{lem:quartic-universal-quotient} is unchanged: at most
four nonunit types supply the $\alpha,\zeta$ coordinates for
$r\ge3$, while the cases $r=1,2$ cost at most seven terms in total.
Hence every class modulo $K_A$ has a representation with at most
twelve fourth powers.  The same three-term baseline cancellation
and Lemma~\ref{lem:primitive-lifting} give $h_4(A)\le15$.
\end{proof}

\subsection{The ninth-power bound}
\label{sec:higher-ninth-local}

Only residue characteristic three requires attention. Use $p=3$,
$a=2$, and coordinate orders $N=(3,5,5)$. In $W_3(R)$ take centers
$(b_0,b_1,b_2)=(0,1,T)$ and perturbations
\[
 (\delta_0,\delta_1,\delta_2)
       =(T^3u_0,T^5u_1,T^3u_2).
\]
For $H=\sum[b_i+\delta_i]-\sum[b_i]$, the universal ghost
recurrence gives $H_j\in T^{N_j}R$. After dividing these coordinates
by $T^3,T^5,T^5$, the reduction modulo $T$ is
\[
 (u_0+u_2,\ u_1+u_2,\ u_1).
\]
Indeed the linear term of $H_1$ is $\sum b_i^2\delta_i$, and that
of $H_2$ is $\sum b_i^8\delta_i$; the other terms have strictly
higher $T$-orders. The displayed linear map has determinant $-1$,
so formal Hensel inversion covers every element of $I_N$.
The perturbed second root is $1+T^5u_1$, hence a unit.

The image lengths are $L_0,L_1,L_2=13,4,1$, so the positive
coefficient cost is
\[
 C_3(3,5,5)=13+2(4+3\cdot1)=27.
\]
The three baseline terms give thirty ninth powers modulo $27A$.
A unit root absorbs the remaining error: the polynomial
\[
 \frac{(1+3z)^9-1}{27}\equiv z\pmod\pi
\]
is integral and has a unit derivative. Hensel's lemma gives the
exact primitive representation with the same number of summands.
At all other residue characteristics, finite-field compression of
$\nu-1$, appending one, and ordinary Hensel lifting use at most ten
terms. Thus $h_9(A)\le30$ without any restriction on residue degree
or ramification.

\subsection{A primitive eighth-power bound of thirty-five}

\begin{theorem}\label{thm:local-eighth-thirtyfive}
For the integer ring $A$ of every finite extension of $\mathbb Q_p$,
\[
 h_8(A)\le35.
\]
More precisely, at dyadic places the following bounds hold:
\[
 \begin{array}{c|c}
 \text{residue degree and absolute ramification}&\text{bound}\\\hline
 f\ge3&19\\
 f=2&30\\
 f=1,\ e=1&32\\
 f=1,\ e=2&32\\
 f=1,\ e=3,4&35\\
 f=1,\ e=5&22\\
 f=1,\ e\ge6&33.
 \end{array}
\]
The entries $e=5$ and $e\ge6$ use the exact finite certificates
described below. Their sizes are respectively $8192$ residue classes
and $64$ polynomial identities plus eight carry relations.
\end{theorem}

Write $\pi$ for a uniformizer. A primitive eighth-power congruence
modulo $16\pi A$ lifts without increasing its number of summands.
Indeed, changing a unit root $x$ to $x(1+2\pi z)$ changes its power by
$16\pi x^8F(z)$, where $F\in A[z]$ and $F(z)\equiv z\pmod\pi$.
All higher binomial terms are divisible by an additional $\pi$ after
division by $16\pi$. Simple-root Hensel lifting applies.
At an odd prime the usual derivative argument gives $h_8(A)\le9$.

\subsection{Residue degrees at least two}

Work first in $W_4(\mathbb F_q[[T]])$. For a centered sum of four
Teichm\"uller vectors with centers $b_i$, its Jacobian has entries
$b_i^{2^j-1}$, $0\le j,i\le3$, with the first row equal to one.
This follows directly by differentiating the universal ghost identities.

If $f\ge3$, take centers $0,c_1,c_2,c_3$, where the three nonzero
constants are linearly independent over $\mathbb F_2$. The zero column
is $(1,0,0,0)^t$; the other minor is an invertible Moore matrix, after
multiplication of its columns by the $c_i$. The centered sum therefore
covers the ideal whose four coordinates belong to $TR$: substitute
perturbations $Tz$, divide the four coordinates by $T$, and apply
formal Hensel inversion to the resulting map with invertible
linear reduction.
By Lemma~\ref{lem:higher-witt-coefficients}, the ambient coefficient cost is
$4+3+2\cdot2+4\cdot1=15$, so the four-variable baseline gives
nineteen powers modulo sixteen.

If $f=2$, choose $\omega^2+\omega+1=0$ and centers $(0,1,\omega,T)$.
The determinant is
\[
 (1+\omega)(T^7+T)=T\cdot\text{a unit}.
\]
Multiplication by the adjugate and the substitution
$\eta=T^2\operatorname{adj}(J)z$ normalize the centered map to
$uz+TP(z)$ after division by $T^3$, with $u$ a unit. Thus the
isotropic ideal of coordinate order three is covered. Its $2$-power
image lengths are $12,6,2,1$, and its coefficient cost is
\[
 12+6+2\cdot2+4\cdot1=26.
\]
Adding the four-variable baseline gives thirty powers modulo sixteen,
by \eqref{eq:higher-last-ghost} and the baseline cancellation above.

To justify the last integral lift, lift the represented roots from
$R=\mathbb F_q[[T]]$ to $D=B[[T]]$. The residual target is in $16D$
by \eqref{eq:higher-last-ghost}. Both baselines retain two units
with distinct eighth-power residues. The final digit modulo
thirty-two can therefore be corrected using these same roots, with
no extra summands. Explicitly,
\[
 (1+2z)^8\equiv1+16(z^2+z)\pmod{32}.
\]
If their lifted roots are $y_1,y_2$, the attainable correction divided
by sixteen is
\[
 \bar y_1^8(z_1^2+z_1)+\bar y_2^8(z_2^2+z_2)\quad\text{in }R.
\]
The Artin--Schreier map $z\mapsto z^2+z$ is bijective on $TR$ and
has image the absolute-trace-zero hyperplane on $\mathbb F_q$.
Multiplication by either unit $\bar y_i^8$ preserves $TR$.
Its two constant images are distinct trace hyperplanes because the
nonzero constants $\bar y_1(0)^8$ and $\bar y_2(0)^8$ are distinct;
the nondegenerate trace pairing identifies their annihilator lines
with their inverse scalars. The sum of these images is consequently
all of $R$. Choose $z_1,z_2$, lift them to $D$, and replace the two
roots by $y_i(1+2z_i)$. This matches the target modulo $32D$ while
preserving both unit roots. After specialization, $32A\subseteq16\pi A$,
so the primitive lifting argument already proved completes the exact
representations with nineteen and thirty terms.

\subsection{A seven-variable weighted chart}

Let $D=\mathbb Z_2[[T]]$ and set
\[
 I=\{(z_0,z_1,z_2,z_3)\in W_4(D):
 z_0\in T^2D,\ z_1\in T^4D,\ z_2,z_3\in T^6D\}.
\]
The inequalities $2^{j-i}N_i\ge N_j$ for
$N=(2,4,6,6)$ show, by weighted homogeneity of the universal Witt
operations, that $I$ is an ideal. The same definition applies after
reduction modulo two.

\begin{lemma}\label{lem:eighth-seven-chart}
Every $w\in I$ has a representation
\[
 w=\sum_{i=0}^6[x_i]-\sum_{i=0}^6[b_i],
 \qquad (b_0,\ldots,b_6)=(0,1,T,T^2,0,0,0),
\]
where $x_1\equiv1\pmod{T}$.
\end{lemma}
\begin{proof}
Use perturbations
\[
 \delta_1=T^6u_1,\qquad\delta_2=T^3u_2,\qquad
 \delta_i=T^2u_i\quad(i=0,3,4,5,6).
\]
Write $H=\sum[b_i+\delta_i]-\sum[b_i]$ and put
$v=(u_0,u_3,u_4,u_5,u_6)$. The universal ghost recurrence shows that
$H_j$ is divisible by $T^{N_j}$. Modulo $(2,T)$, the normalized map is
\[
 \left(\sigma_1(v),\ u_2+u_3+\sigma_2(v),\ u_1+u_2,\ u_1\right).
\]
For the second component, one may use the exact formula
$H_1=\sum b_i\delta_i-\sum_{i<j}\delta_i\delta_j$.
In the last two components, all nonlinear ghost terms have orders
strictly greater than six, leaving the displayed linear terms.

The last two equations fix $u_1,u_2$. For the first two, it is enough
to solve $\sigma_1(v)=s$ and $\sigma_2(v)+u_3=c$ nonsingularly
over $\mathbb F_2$. For $(s,c)=(0,0),(0,1),(1,0),(1,1)$, choose
respectively zero, two, one, and three ones among the four coordinates
other than $u_3$, keeping $u_3=0$. The two Jacobian rows are independent:
their difference from $1$ is the vector $v+e_3$, which is neither zero
nor the all-one vector for these choices. The full Jacobian has rank
four. First lift the chosen solution through the $2$-adic coefficient
ring, then through the $T$-adic parameter, freezing the unused variables.
This solves every normalized target in $D^4$ and proves the lemma.
The root $x_1=1+T^6u_1$ remains a unit.
\end{proof}

\subsection{The finite coefficient certificate}

Put $\bar D=\mathbb F_2[[T]]$ and let $G$ be the actual additive
Teichm\"uller span in $W_4(\bar D)/\bar I$.
The following eight actual generators suffice:
\[
 [1],\ [1+T^5],\ [1+T^4],\ [1+T^3],\ [1+T^2],\
 [1+T],\ [T],\ [T+T^2].
\]
In this order, their successive cyclic quotient orders are
\[
 16,2,2,2,2,4,4,2.
\]
Consequently every element of $G$ is a positive sum of at most
\[
 15+1+1+1+1+3+3+1=26
\]
of the displayed generators.

Here is the precise finite verification used for that assertion.
The combined ghost map sends a vector to $2^{3-j}w_j$, truncated
at degree $N_j$, modulo sixteen for each $j=0,1,2,3$.
It is additive and has kernel exactly $\bar I$, by induction on the
coordinates. Every root can be reduced modulo $T^6$, leaving exactly
sixty-four binary polynomials. The program in Appendix~\ref{app:eighth-ghost} checks their images
against the $8192$ normal forms with the coefficient ranges just
displayed. It also checks all eight carry relations, which make this
normal-form set an additive subgroup. If the generators above are
$g_0,\ldots,g_7$, these relations are
\[
 \begin{split}
 16g_0&=0,\quad 2g_1=2g_0,\quad2g_2=2g_0,\quad2g_3=2g_0,\\
 2g_4&=g_0+g_2,\quad4g_5=3g_0+g_2,\\
 4g_6&=12g_0+g_2+g_4+2g_5,\quad2g_7=2g_6.
 \end{split}
\]
The program also emits a coefficient witness for every binary root. All checks use exact integer arithmetic modulo sixteen.

\subsection{High ramification}

Suppose $f=1$ and $e\ge6$. Lift an eligible target through a finite
signed-power expression to $D$ and then to a Witt vector $g\in W_4(D)$
whose last ghost is that expression. Subtract the seven-term baseline,
reduce modulo two, and apply the preceding certificate modulo $\bar I$.
This uses at most twenty-six positive Teichm\"uller vectors. After
lifting them to $D$, the remaining Witt vector $\eta$ reduces into
$\bar I$. Choose $w\in I$ lifting this reduction. The last ghost of
$\eta-w$ is divisible by sixteen.

After specialization $T=\pi$, this remaining error belongs to
$16A\subseteq8\pi^6A$, since $e\ge6$. It can therefore be matched
exactly by the last ghost of a vector $V^3(T^6z)\in I$, after lifting
its coefficient from $A$ to $D$. Add that vector to $w$ and apply
Lemma~\ref{lem:eighth-seven-chart}. The fixed baseline cancels, leaving
at most $26+7=33$ positive eighth powers, one with a unit root.

\subsection{Small ramification}

The rational ring has $h_8(\mathbb Z_2)=32$: unit eighth powers are
one modulo thirty-two, nonunit eighth powers vanish there, and a
primitive zero needs thirty-two unit terms. Any residue class is
represented by between one and thirty-two ones, followed by lifting.

For $e=2,3,4$, work modulo $16\pi A=\pi^{4e+1}A$.
The eighth power of every unit is $1+d$ with
\[
 v_\pi(d)\ge7\quad(e=2,3),\qquad
 v_\pi(d)\ge9\quad(e=4).
\]
For $e=2,4$, the two lowest terms in the binomial expansion cancel
modulo $\pi$, since the residue field is $\mathbb F_2$; direct use
of $(1+\pi y)^8$ gives these bounds. The product of any two such
deviations vanishes at the working precision. Thus the deviations
form an additive subgroup $D_0$: multiplication of unit eighth powers
corresponds to addition of deviations.

When $e=2$, every nonunit eighth power is zero or $\pi^8$, and
$\pi^8\equiv16$ at this precision. The generated image is consequently
$\mathbb Z\cdot1+D_0$. It is represented by at most thirty-two unit
eighth powers, using one unit power for the prescribed deviation.

When $e=3$, choose $\pi^3=2$, by extracting the cube root of the
unit in any initial uniformizer relation. The nonunit images are
zero and $\pi^8$, and $4\pi^8=0$ at precision $\pi^{13}$.
Thus at most three copies of $\pi^8$, together with at most
thirty-two unit eighth powers, represent every eligible class.

When $e=4$, the nonunit images are zero, $\pi^8$, and $\pi^{16}$.
Writing $\pi^4=2u$ shows
\[
 \pi^{16}\equiv16,\qquad4\pi^8\equiv16\pmod{\pi^{17}}.
\]
The same bound of thirty-five follows. All these representations
have a unit root and lift.

Finally, if $e=5$, extracting a fifth root of a principal unit gives
$A=\mathbb Z_2[\pi]$ with $\pi^5=2$. Eighth powers modulo $16\pi$
depend only on roots modulo $2\pi$. The latter quotient has sixty-four
classes, represented by
\[
 a_0+a_1\pi+\cdots+a_4\pi^4,
 \qquad0\le a_0<4,\quad a_i\in\{0,1\}\ (i>0).
\]
The program in Appendix~\ref{app:eighth-tame} checks these classes in the exact ring
\[
 (\mathbb Z/32)\oplus
 (\mathbb Z/16)\pi\oplus\cdots\oplus(\mathbb Z/16)\pi^4,
 \qquad\pi^5=2.
\]
There are twenty distinct eighth-power images. Their generated group
has $8192$ elements and ordinary additive diameter twenty-one.
The breadth-first calculation verifies every parent representation
and closure under every generator. Applying it to the target minus
one and appending $1^8$ gives a primitive bound twenty-two.
This exhausts the ramification cases and proves
Theorem~\ref{thm:local-eighth-thirtyfive}.

\section{Efficient differencing with refined seed moments}
\label{sec:efficient-differencing}\label{sec:g6-transfer}\label{sec:quartic-principal-blocks}
The efficient-differencing method of Wooley and Vaughan--Wooley
\cite{WooleyThesis,Wooley1992,VW1995,VW2000} can be formulated
using principal blocks of prime ideals. At these fixed small exponents
the terminal estimates also use the number-field Vinogradov mean-value
theorem and coefficient-uniform curve bounds.
\subsection{Localizing the ordinary variables}

One may prove a permissible full-box moment by proving it for smooth roots in
arbitrarily small fixed boxes away from zero at every embedding.  Here are the
details that allow this reduction uniformly in the scale.
Partition nonzero roots into dyadic ranges for each $|x_v|$.  There are
$O((\log X)^{r_1+r_2})$ ranges, since $|Nx|\ge1$ supplies a polynomial lower
bound for every $|x_v|$.  Roots with $|Nx|\le X^{1/2}$ number
$O(X^{1/2+\epsilon})$, by counting ideals and bounded-height units; their
$2s$th moment is trivially $O(X^{s+\epsilon})$.
For each remaining annular box let $Y$ be the product of its archimedean radii,
with local dimensions counted.  Thus $X^{1/2}\ll Y\ll X$.
Multiplication by one unit balances these radii to fixed multiples of $Y^{1/n}$.
It preserves ideal smoothness and the mean value, since multiplication of
$\alpha$ by a unit's $k$th power is a torus automorphism.
Cover the resulting compact annulus by finitely many boxes of relative diameter
$\delta$, where $\delta>0$ is a sufficiently small constant depending only on
$K$ and the maximal differencing depth.  H\"older's inequality for the sum of
these pieces costs only a power of $\log X$.  Since $R=X^\eta\le Y^{2\eta}$,
a sufficiently small initial smoothness exponent works for all pieces.

For one such piece, replace two smooth roots by ordinary roots in that same
small box.  The other smooth roots can subsequently range over full boxes.
The two distinguished ordinary roots, and every pair produced from them by
efficient differencing, retain their actual small support.
In centered notation define
\[
 \Psi_j(z;\mathbf h;\mathbf m)=
 \Delta^*_j\bigl((z-\textstyle\sum_{i\le j}h_i m_i^k)^k;
                 2\mathbf h;\mathbf m\bigr),
 \qquad
 \Delta^*_1 f(z;h;m)=m^{-k}(f(z+hm^k)-f(z)).
\]
These are integral polynomials.  Their leading term is
$2^j(k)_j(h_1\cdots h_j)z^{k-j}$.
At each embedding, the actual support has $z$ within $O_j(\delta P)$ of a
fixed nonzero center of size comparable to $P$, and
\begin{equation}\label{eq:g6-small-increments}
 |h_{i,v}m_{i,v}^k|\le C_j\delta P\qquad(1\le i\le j).
\end{equation}
Fixed factors of two change the center and the $C_j$, harmlessly.
The central change of variables is $z=z_1+z_2$,
$m^kh=z_1-z_2$; divisibility by $2$ can be retained as a support restriction.
More explicitly, for exponent $k$ the centered polynomial identity multiplies
the differenced equation by $2^k$, doubles the earlier shift parameters, and
replaces every remaining quotient root $x$ by $2x$.
This identity applies at every repeated step.
The retained parity restrictions only reduce a solution count; subsequently
enlarging the bounding boxes changes their scales by fixed constants, and the
primes above $2$ belong to the smoothness range for sufficiently large $X$.
Consequently $\delta$ can be chosen so that, whenever all $h_i\ne0$,
\begin{equation}\label{eq:g6-nonzero}
 \Psi_j(z;\mathbf h;\mathbf m)\ne0,
 \qquad \partial_z\Psi_j(z;\mathbf h;\mathbf m)\ne0
 \quad(0\le j\le k-1).
\end{equation}
Indeed, after dividing by $h_1\cdots h_j$, the leading monomial dominates
all remaining monomials at each embedding.  The same is true for the derivative.
An ordinary differencing shift is a difference of two points in the small box;
it therefore retains \eqref{eq:g6-small-increments}.

In all subsequent counts the coupled restrictions
\eqref{eq:g6-small-increments} are retained.  The separate bounding boxes for
$h_i,m_i$ are used only to count parameter choices, never to replace the support
when \eqref{eq:g6-nonzero} is needed.  This is the step that removes the
positivity assumption implicit in the rational auxiliary counts.

\subsection{The fundamental efficient-differencing inequality}

Let $F_\Psi(\alpha)=\sum_{z,\mathbf u}\mathrm e_K(\alpha\Psi(z,\mathbf u))$,
where $z$ ranges in a box of volume $X$, the auxiliary parameter tuples number
at most $U$, every polynomial has bounded positive degree and coefficients of
polynomial height, and $\Psi'(z,\mathbf u)\ne0$ on the retained support.
Let $g_Q$ be the full smooth sum of volume $Q$, and put
$S_s=\int|F_\Psi|^2|g_Q|^{2s}$.
Choose $1\le M\le Q$.

Let $h_K$ be the order of the ideal class group.  Any sequence of $h_K$ ideal
classes has a nonempty zero-sum subsequence.  Partition the prime ideal factors
of a smooth ideal into such principal blocks, each with at most $h_K$ factors,
leaving fewer than $h_K$ factors.  If the part of $(x)$ coprime to
$(\Psi'(z,\mathbf u))$ has norm greater than $MR^{h_K}$, successive principal
blocks supply a divisor $(m)$ of $(x)$, coprime to $\Psi'$, with
\[
 M<N(m)\le MR^{h_K}.
\]
Choose a generator balanced by units.  There are $O(MR^{h_K+\epsilon})$
possible generators, and $x/m$ is a smooth algebraic integer in a fixed box
of volume $Q/M$.  Powers of $R$ will always be absorbed in $X^\epsilon$.

If no such divisor exists, the prime-to-$\Psi'$ part of $(x)$ has norm at most
$MR^{h_K}$.  For fixed $z,\mathbf u$ the number of these bad $x$ is at most
$MX^\epsilon$: count the prime-to-$\Psi'$ ideal and use the supported-ideal
bound for its complementary ideal, then count bounded-height generators.
This bound also covers roots of norm at most $MR^{h_K}$.

The contribution of a bad variable is at most
$S_s^{1/2}(QMX^\epsilon S_{s-1})^{1/2}$.
Indeed, if
$G_x=\sum_{z,\mathbf u:,x\text{ bad}}\mathrm e_K(\alpha\Psi(z,\mathbf u))$,
then
\[
 \left|\sum_x\mathrm e_K(\alpha x^k)G_x\right|^2
 \le O(Q)\sum_x|G_x|^2.
\]
After multiplication by $|g_Q|^{2s-2}$ and integration, a fixed pair of auxiliary
tuples admits at most $MX^\epsilon$ common bad $x$, giving the assertion.

For the good variables, extract a principal divisor separately from each of the
$2s$ smooth roots.  The dependence of this extraction on the derivative is
handled as follows.  Let $\mathcal M$ be the set of balanced extracted divisors,
so that $L=\#\mathcal M\ll MX^\epsilon$.  For $W\ne0$ put
\[
 F_W(\alpha)=\sum_{\substack{z,\mathbf u\text{ on the retained support}\\
              ((W),(\Psi'(z,\mathbf u)))=1}}
       \mathrm e_K(\alpha\Psi(z,\mathbf u)),
 \quad g=g_{Q/M},\quad
 V_m=\int|F_m(\alpha)|^2|g(m^k\alpha)|^{2s}\,d\alpha.
\]
For tuples $\mathbf m,\mathbf n\in\mathcal M^s$ put
$W=\prod_i m_i$ and $W'=\prod_i n_i$.
After enlarging the quotient-root boxes, orthogonality bounds the good count by
the sum over these tuples of
\[
 \int |F_W(\alpha)F_{W'}(\alpha)|
       \prod_{i=1}^s|g(m_i^k\alpha)g(n_i^k\alpha)|\,d\alpha.
\]
H\"older bounds each integral by the product of $2s$ factors of the form
$(\int|F_W|^2|g(m_i^k\alpha)|^{2s})^{1/(2s)}$ and their $W',n_i$ analogues.
The coefficient set of $F_W$ is contained in that of $F_{m_i}$, because
coprimality to $W$ implies coprimality to $m_i$.  Orthogonality, with its
nonnegative solution counts, therefore bounds that integral by $V_{m_i}$.
Consequently the total good count is at most
\[
 \left(\sum_{m\in\mathcal M}V_m^{1/(2s)}\right)^{2s}
 \le L^{2s-1}\sum_{m\in\mathcal M}V_m.
\]
Both derivative coprimalities are present in every $V_m$.

Here is the remaining Cauchy step explicitly.  Set $\mathfrak b=(m)^k$.
For $\rho\in\mathcal O_K/\mathfrak b$ and fixed $\mathbf u$, let
$\mathcal R(\rho,\mathbf u)$ be the nonsingular residue classes $r$ satisfying
$\Psi(r,\mathbf u)\equiv\rho\pmod{\mathfrak b}$.
Their number is at most $D^{\omega((m))}\ll X^\epsilon$, by Hensel lifting
and the Chinese remainder theorem.  Write
$F_{\mathbf u,r}=\sum_{z\equiv r\, (\mathfrak b)}
\mathrm e_K(\alpha\Psi(z,\mathbf u))$, retaining the original support.
Orthogonality kills pairs with different $\Psi$ residues, and hence
\[
 V_m=\sum_\rho\int
 \left|\sum_{\mathbf u}\sum_{r\in\mathcal R(\rho,\mathbf u)}
                   F_{\mathbf u,r}(\alpha)\right|^2
 |g(m^k\alpha)|^{2s}\,d\alpha.
\]
Cauchy's inequality inside each square costs at most $UX^\epsilon$ and leaves
the sum of $|F_{\mathbf u,r}|^2$.  Its underlying equations have the same
auxiliary tuple on both sides and $z\equiv w\pmod{(m)^k}$.
Write $z-w=m^kh$.  The case $h=0$ contributes at most
$XUM(Q/M)^{\lambda_s+\epsilon}$; the other cases are bounded by
$T_{\rm next}=\int|F_{\rm next}||g_{Q/M}|^{2s}$, with all actual support
restrictions retained.  Thus the good contribution is at most
\[
 X^\epsilon U M^{2s-1}
 \{XUM(Q/M)^{\lambda_s}+T_{\rm next}\}.
\]

The bad contribution can be absorbed into this same explicit diagonal bound
without a lower bound for the density of smooth algebraic integers.
Indeed, $S_0\ll XU^2$, by fixing $z,\mathbf u,\mathbf v$ and solving the
nonconstant equation for $w$.  H\"older gives
$S_{s-1}\le S_0^{1/s}S_s^{1-1/s}$.
If the bad contribution dominates, then
\[
 S_s\ll X^\epsilon XU^2(QM)^s
 \le X^\epsilon XU^2M^{2s}(Q/M)^{\lambda_s},
\]
since $Q\ge M$ and $\lambda_s\ge s$.
In every case, therefore,
\begin{equation}\label{eq:g6-fundamental}
 S_s\ll X^\epsilon U M^{2s-1}
 \{XUM(Q/M)^{\lambda_s}+T_{\rm next}\}.
\end{equation}
This is the required field version of the fundamental efficient-differencing
inequality.  Its proof is the congruence argument of Wooley's 1990 thesis,
Chapter~6, with principal ideal blocks and the preceding explicit treatment
of bad variables and the diagonal term.

\subsection{Auxiliary seed moments}
\begin{lemma}[Sparse progression endpoint over a number field]
\label{lem:field-sparse-endpoint}
Let $p\ge2$ and $J\ge1$ be fixed integers. For $1\le s\le J+1$,
\[
 \#\left\{\mathbf x,\mathbf y\in\mathcal B(M)^s:
       \sum_{i=1}^s x_i^{pr}=\sum_{i=1}^s y_i^{pr}\ (1\le r\le J)\right\}
 \ll_{K,p,J,\epsilon} M^{s+\epsilon}.
\]
\end{lemma}
\begin{proof}
For $s\le J$, Newton's identities determine the multiset of $p$th powers
on one side from the other. Each powered coordinate has at most $p$ roots,
so there are $O_{p,J}(M^s)$ solutions. The case $J=1,s=2$ is the ordinary
fourth-moment bound for the degree-$p$ sum, supplied by Hua's differencing
lemma over $K$.

Suppose $J\ge2$ and $s=J+1$. Interchanging the first coordinates of the two
sides turns the equations into equality of signed power sums
\[
 S_r(\mathbf x)=\sum_{i=1}^J x_i^{pr}-x_0^{pr}
             =S_r(\mathbf y)\qquad(1\le r\le J).
\]
These equalities imply
\[
 D(\mathbf x)=\prod_{i=1}^J(x_0^p-x_i^p)
             =\prod_{i=1}^J(y_0^p-y_i^p)=D(\mathbf y).
\]
Here is a direct algebraic verification. If $x_0^p=y_0^p$, apply Newton's
identities to the other $J$ powered coordinates. Otherwise the two monic
polynomials with roots $(y_0^p,x_1^p,\ldots,x_J^p)$ and
$(x_0^p,y_1^p,\ldots,y_J^p)$ differ by a constant, because their first $J$
power sums agree. Evaluating their difference at $x_0^p$ and $y_0^p$ gives
$(x_0^p-y_0^p)D(\mathbf x)=(x_0^p-y_0^p)D(\mathbf y)$.

Fix $\mathbf y$ with $D(\mathbf y)\ne0$. Every difference
$d_i=x_0^p-x_i^p$ is a bounded-height divisor of this fixed nonzero algebraic
integer. There are $O(M^\epsilon)$ choices for the whole difference tuple.
Moreover
\[
 (J-1)x_0^p=S_1(\mathbf y)+\sum_i d_i,
\]
so the powers of all coordinates, and hence the coordinates up to a fixed
multiplicity, are determined. This gives $O(M^{J+1+\epsilon})$ solutions.
If $D(\mathbf y)=0$, choose a coordinate with $y_i^p=y_0^p$; likewise one
coordinate of $\mathbf x$ has the same power as $x_0$.
Cancel these in the signed sums. The remaining $J-1$ powered coordinates
are determined up to permutation by Newton's identities. There are
$O(M^J)$ choices for $\mathbf y$ and $O(M)$ choices for the free center
$x_0$, again with bounded root multiplicities. This proves the lemma.
\end{proof}

For efficient differences choose $0\le\phi_i<1/k$ and let
\[
 M_i=X^{\phi_i},\quad H_i=X M_i^{-k},\quad
 \Phi_j=\sum_{i=1}^j\phi_i,\quad
 U_j=\prod_{i=1}^jH_iM_i.
\]
Fixed small constants and fixed smoothness powers are suppressed in these
scale definitions. The centered polynomial is denoted by $\Psi_j(z;\mathbf
h,\mathbf m)$, and $F_j$ is its sum over the small ordinary box and the
independent shift and smooth-coefficient boxes. The original ordinary
vertices remain in a common small box away from zero at every embedding.
Thus $\Psi_j/\prod h_i$ and its required derivatives are nonzero wherever
they are used. This is also true after the finitely many ordinary differences
below, which have small shifts and coefficient divisor $1$.

\begin{lemma}[Sparse moments from a nonvanishing Wronskian]
\label{lem:field-sparse-wronskian}
For fixed integers $p\ge2$ and $J\ge1$, put
\[
 \tau_J=\max\{J+1,J(J+1)/2\}.
\]
For every integer $1\le s\le\tau_J$, the number of solutions in a
box of volume $M$ to
\[
 \sum_{i=1}^s x_i^{pa}=\sum_{i=1}^s y_i^{pa}
 \qquad(1\le a\le J)
\]
is $O_{K,p,J,\epsilon}(M^{s+\epsilon})$.
\end{lemma}
\begin{proof}
The range $s\le J+1$ is Lemma~\ref{lem:field-sparse-endpoint}.
For $s\le J(J+1)/2$, apply the number-field mean-value theorem
of Wooley~\cite[Theorem~15.1]{Wooley2019} to the $J$ polynomials
$t^p,t^{2p},\ldots,t^{Jp}$. Their Wronskian is a nonzero monomial:
its coefficient is the product of their nonzero degrees and a
Vandermonde determinant in these distinct degrees.
The theorem applies to arbitrary complex coefficients, so restricting
to the prescribed box and using orthogonality gives the assertion.
\end{proof}

\begin{lemma}[Centered curve intersections]
\label{lem:centered-curve-intersections}
Fix $d\ge3$ and a nonzero leading coefficient $c\in K$.
Let $P,Q\in K[t]$ have degree $d$, leading coefficient $c$, and
vanishing coefficients of $t^{d-1}$. Suppose their coefficients
have fixed denominators and archimedean height bounded by a fixed
power of $X$. Let $\mathcal B$ be a sufficiently small box of
volume comparable to $X$ about a fixed nonzero archimedean center
at scale $X^{1/[K:\mathbb Q]}$. Then
\[
 \#\{(z,w)\in(\mathcal B\cap\mathcal O_K)^2:P(z)=Q(w)\}
 \ll X\,\mathbf1_{P=Q}+X^{1/3+\epsilon},
\]
uniformly in the lower coefficients of $P,Q$.
\end{lemma}
\begin{proof}
Factor $P(z)-Q(w)$ into irreducible factors over $K$.
Its highest homogeneous part is $c(z^d-w^d)$, which is square-free.
Every linear factor has the form $z-\zeta w-b$, with $\zeta^d=1$.
Substitution and comparison of the coefficient of $w^{d-1}$ give
$b=0$. Choosing the box smaller than a fixed fraction of the
distance between its center and its nontrivial root-of-unity
multiples excludes $z=\zeta w$ for $\zeta\ne1$.
The line $z=w$ occurs precisely when $P=Q$, and contributes $O(X)$.

The highest homogeneous part of a quadratic factor is the product
of two distinct factors of $z^d-w^d$. After normalizing one
coefficient it belongs to a fixed finite collection of nondegenerate
binary quadratic forms. Its lower coefficients have polynomial
height and fixed denominators. To see the denominator assertion,
clear the fixed leading coefficient by a fixed change in $z$;
the resulting polynomial is monic in $z$ over $\mathcal O_K[w]$.
Its monic factors belong to this integrally closed ring.
Standard fixed-degree factor bounds give the height assertion.
Translation to the center of the quadratic form therefore changes
the lattice by a fixed denominator only. The resulting equation is
a norm equation in one of finitely many fixed quadratic extensions
of $K$, or a product equation over $K$.
For a nonzero right side, the ideal-divisor estimate and the
bounded-height unit count give $O(X^\epsilon)$ solutions.
For a zero right side, an irreducible quadratic factor has at most
one $K$-point; split linear factors have already been treated.

Every remaining irreducible component has degree at least three.
The affine-curve estimate of Binyamini--Cluckers--Kato
\cite[Corollary~12]{BCK2025} bounds its integral points by
$O_K(X^{1/3}(\log X)^\kappa)$, uniformly in its coefficients.
Here the height parameter in that result is $O_K(X)$: the
embedding bound is $O_K(X^{1/[K:\mathbb Q]})$, as in their
definition preceding equation~(28).
There are at most $d$ components. Summing their contributions
proves the lemma.
\end{proof}

\begin{lemma}[Auxiliary moments for efficient differencing]
\label{lem:field-general-seed-aux}
Let $1\le\ell\le3$ and $1\le j\le k-1-\ell$, and put
\[
 J=\left\lfloor\frac{k-j-\ell+1}{2}\right\rfloor,
 \qquad \tau_J=\max\{J+1,J(J+1)/2\}.
\]
Assume the chronological divisor inequalities
\begin{equation}\label{eq:field-chronological-prefix}
 \Phi_i+(k-2)\phi_i\le1
 \qquad(\tau_J<i\le j).
\end{equation}
This condition is empty when $j\le\tau_J$.
Assume also that $\Phi_j\le2/3$ whenever at least one of the
residual degrees $k-j-r$, $0\le r<\ell$, is even and at least six.
Then
\[
 \int_{\mathbb T_K}|F_j(\alpha)|^{2^\ell}\,d\alpha
 \ll X^{2^\ell-\ell+\epsilon}U_j^{2^\ell-1}.
\]
\end{lemma}
\begin{proof}
Consider first a second moment after $r$ appended ordinary
differences, where $0\le r<\ell$. The retained small supports
make the divided phase nonzero. Grouping its nonzero value by
the divisor tuple and applying Cauchy's inequality costs
$X^\epsilon$ and reduces the count to a common tuple of efficient
and ordinary shifts on the two sides.
Polynomial coincidence implies
\begin{equation}\label{eq:seed-coefficient-coincidence}
 \sum_{i=1}^j h_i^{2a}(m_i^{2ka}-n_i^{2ka})=0
 \qquad(1\le a\le J).
\end{equation}
The implication follows by taking the logarithm of the centered
generating function
\[
 e^{zt}\prod_i\frac{\sinh(h_im_i^kt)}{h_im_i^kt}.
\]
Euler's product gives
\[
 [t^{2a}]\log(\sinh t/t)=
       \frac{(-1)^{a+1}\zeta(2a)}{a\pi^{2a}}\ne0.
\]
The ordinary coefficient blocks are equal on the two sides.

When $j\le\tau_J$, fix the shifts and apply H\"older to the $j$
coefficient sums on $\mathbb T_K^J$.
Multiplication of coordinate $a$ by $h_i^{2a}$ preserves probability
Haar measure. Lemma~\ref{lem:field-sparse-wronskian}, with $p=2k$,
therefore gives
\[
 \prod_{i=1}^j\left(M_i^{j+\epsilon}\right)^{1/j}
 \ll X^\epsilon\prod_{i=1}^j M_i
\]
coefficient pairs, including unequal coefficient volumes.
Summing the shifts gives $O(X^\epsilon U_j)$.

For $j>\tau_J$, eliminate the last coefficient block using the
first equation of \eqref{eq:seed-coefficient-coincidence}.
If $m_i^{2k}=n_i^{2k}$, at most $2k$ choices of $n_i$ occur for
each $m_i$, and all $J$ equations reduce to those of the preceding
blocks. Otherwise
$h_i^2(m_i^{2k}-n_i^{2k})=b\ne0$ has $O(X^\epsilon)$ solutions
after the other blocks are fixed: first count divisors $h_i$ of $b$,
then divisors $m_i-n_i$ of the remaining fixed integer, and finally
roots of a nonzero fixed-degree polynomial.
Consequently the coefficient count satisfies
\[
 N_i\ll X^\epsilon\left(H_iM_iN_{i-1}
                    +\prod_{a<i}H_aM_a^2\right).
\]
The inequality in \eqref{eq:field-chronological-prefix} is exactly
$\prod_{a<i}M_a\le H_iM_i$.
Stop the elimination at $i=\tau_J$ and apply the sparse estimate
there. Induction gives $N_j\ll X^\epsilon U_j$.

For odd residual degree, the centered phase has a factor $z$.
Divisor grouping in its nonzero value makes the two centers equal
at cost $X^\epsilon$; a nonzero difference of the two polynomials
then has boundedly many roots. This gives
$O(X^{r+1+\epsilon}U_j)$, since $\prod_iM_i\le X$.
At residual degree two, use a difference of squares. At degree
four, complete the squares in $z^2,w^2$; a nonzero level is a
divisor equation, and a zero level gives either polynomial
coincidence or a binary quadratic norm equation. The nontrivial
root-of-unity lines miss the common small supports. The same
bound follows in both cases.

For the remaining even residual degrees, divide the common-shift
polynomial by its leading coefficient. The two resulting centered
polynomials have the same fixed leading coefficient, vanishing
next coefficient, and fixed denominators. Their lower coefficients
have polynomial height. Lemma~\ref{lem:centered-curve-intersections}
gives
\[
 R_{j,r}\ll X^{r+1+\epsilon}U_j+
 X^{r+1/3+\epsilon}\Bigl(\prod_iH_i\Bigr)
                         \Bigl(\prod_iM_i\Bigr)^2.
\]
The second term is absorbed when $\Phi_j\le2/3$.
The padded ordinary-moment induction below now gives the stated
$2^\ell$th moment.
\end{proof}

\paragraph{The padded ordinary-moment induction.}
Here is an explicit induction from the appended second moments to the
higher moments used by A and B.  Let $F_r$ denote the family with $r$
ordinary differences already appended, and set $U_r=X^rU_j$.
Cauchy in the parameter tuple and ordinary differencing give
\[
                     |F_r|^2\ll XU_r^2+U_r|F_{r+1}|.
\]
The zero ordinary shift is the first term; all nonzero shifts retain the
small-support restrictions.  Put
$B_{r,a}=\int|F_r|^{2^a}$.  For $a\ge2$, multiply the pointwise inequality
by $|F_r|^{2^a-2}$.  Interpolation in the first term and H\"older in the
second, followed by absorption, yield
\[
 B_{r,a}\ll (XU_r^2)^{2^{a-2}}B_{r,a-1}
                  +U_r^{2^{a-1}}B_{r+1,a-1}.
\]
Starting with $B_{r,1}\ll X^{1+\epsilon}U_r$, induction gives
\[
                  B_{r,a}\ll X^{2^a-a+\epsilon}U_r^{2^a-1}.
\]
Indeed, the two powers of $X$ obtained on substitution are
$3\cdot2^{a-2}-a+1$ and $2^a-a$, and the former is no larger than the
latter.  The restriction $j+\ell\le k-1$ means that every final residual
degree is at least two; the extra multiplicity from a residual linear
phase never occurs.  The proof retains the parameter tuple throughout the
ordinary differencing and consequently introduces no new independent
copies of the efficient-divisor variables.

\paragraph{Preserving the supports.}
All arguments above may be performed on the original supports.  If
Cartesian boxes are wanted, split the balanced extracted generators into
dyadic norm annuli and then small archimedean boxes.  On one piece their
conjugates have comparable fixed scales.  The shift ranges can be enlarged
by fixed factors while all ordinary vertices remain in a slightly larger
common tiny box.  The enlargement is made in nonnegative orthogonality
counts before an absolute-value bound is used.  The number of pieces is
$X^\epsilon$ and changes in scale are $X^{O(\eta)}$.
Choose the initial relative diameter small enough for the fixed total
number of differences.  The leading monomial of every required residual
polynomial, and of its relevant derivative, then dominates at every
embedding.  This justifies all exclusions of zero phase values above.

\subsection{The exponent transformations}
Write $b_s$ for a permissible full exponent at moment $2s$, and
$b(t)$ for linear interpolation between adjacent integer indices.
Each seed row specifies $s$, a depth $j$, an order $\ell$, an A or B
terminal, and rational parameters $\phi_1,\ldots,\phi_j$. Put
\[
 q_i=1-\Phi_i,\qquad
 \kappa_\ell=2^\ell-\ell+(2^\ell-1)(j-(k-1)\Phi_j),
 \qquad u=s-j.
\]
The A terminal is
\begin{equation}\label{eq:seed-A-terminal}
 C_j=2^{-\ell}\kappa_\ell+
 (1-2^{-\ell})q_j b\bigl(u/(1-2^{-\ell})\bigr).
\end{equation}
For the B terminal choose an additional integer $t$ and put
$c=u/t$, $a=2(1-c)-2^{1-\ell}$, and
$b'=2^{1-\ell}-(1-c)$. If these weights are nonnegative, then
\begin{equation}\label{eq:seed-B-terminal}
 C_j=a\kappa_{\ell-1}+b'\kappa_\ell+cq_jb(t).
\end{equation}
These follow by H\"older from Lemma~\ref{lem:field-general-seed-aux}.
Proceeding backwards for $i=j-1,\ldots,1$, set
\begin{equation}\label{eq:seed-backwards}
 \begin{split}
 E_i&=\tfrac12\{i-(k-1)\Phi_i+(2(s-i-1)-1)\phi_{i+1}
                              +q_i b(s-i+1)\},\\
 D_i&=i+1-(k-1)\Phi_i+\phi_{i+1}+q_{i+1}b(s-i-1),\\
 C_i&=E_i+\tfrac12\max(D_i,C_{i+1}).
 \end{split}
\end{equation}
The first-stage output is
\begin{equation}\label{eq:seed-initial-output}
 T_s(b)=(2s-3)\phi_1+
       \max\{1+\phi_1+q_1b(s-1),C_1\}.
\end{equation}
The fundamental principal-block inequality proves these backwards
steps. Different terminal estimates may be substituted for $C_j$;
the remainder of the transformation is unchanged.

\begin{lemma}[Finite supersolution principle]\label{lem:supersolution}
Let a finite collection of moment transformations have nonnegative
coefficients in the input exponents, and let the sum of the
coefficients in every branch be at most $\rho<1$. Finite minima
and maxima of these branches are allowed. Suppose that a rational
vector $b$ satisfies $T_s(b)<b_s$ for every nontrivial coordinate,
and agrees with permissible elementary exponents at the remaining
coordinates. Then every $b_s$ is permissible. If outputs are clamped
below by the diagonal exponent $s$, assume also that $b_s\ge s$.
\end{lemma}
\begin{proof}
Start from a vector of elementary permissible exponents dominating
$b$. Monotonicity and the coefficient-sum bound show that, after
one transformation, its excess above $b$ in sup norm is at most
$\rho$ times the previous excess. The same assertion holds after
clamping, since the constant branch has coefficient sum zero.
After finitely many steps the excess is arbitrarily small.
For any prescribed error exponent, choose that number of steps
first and then choose the common smoothness parameter sufficiently
small for all the finitely many estimates. This gives the required
$X^\epsilon$ loss.
\end{proof}

\section{Coefficient sampling}
\label{sec:coefficient-sampling-common}
\begin{lemma}[A bilinear sampling bound]
\label{lem:g789-bilinear}
If $J,N\ge1$, $c(a)$ is supported on $\mathcal B(N)$, and
$\theta\in V$, then
\[
 \sum_{z\in\mathcal B(J)}
       \left|\sum_a c(a)E(\theta za)\right|^2
 \ll_K\left(J+N+\frac{JN}{\mathcal H_{JN}(\theta)}\right)
                         \sum_a|c(a)|^2.
\]
The same assertion holds for arbitrary subsets of the displayed
boxes and arbitrary fixed translates of those boxes.
\end{lemma}
\begin{proof}
First suppose $J\le N$. The multidimensional elementary large-sieve
inequality bounds the left side by $O_K(NM)\sum|c(a)|^2$, where $M$
is the maximum number of sample frequencies $\theta z$ in a translate
of an archimedean box of side $O_K(N^{-1/n})$ modulo
$\mathfrak D^{-1}$. One direct proof chooses a positive Gaussian
majorant for the coefficient box, applies Poisson summation to its
Gram matrix, and bounds its operator norm by the maximum row sum.
Covering successive frequency annuli by boxes of side $N^{-1/n}$
and using Gaussian decay bounds that row sum by $O_K(NM)$.

Partition the $z$ box into a fixed number of smaller boxes if
necessary. The affine field-line lemma in Section~\ref{sec:ordinary-bundle}
applies to the close pairs $(z,a)$, because their volume product is
$O_K(J/N)$, made sufficiently small by this fixed subdivision.
It gives either at most one point or a slope $\gamma\in K$ and
\[
 M\ll_K1+Q_\gamma^{-1}
       \prod_v\min\{J^{1/n},N^{-1/n}/|\theta_v-\gamma_v|\}^{d_v}
 =O_K\left(1+\frac{J}{Q_\gamma A_{JN}(\theta-\gamma)}\right).
\]
This is at most $O_K(1+J/\mathcal H_{JN}(\theta))$.
The large-sieve bound proves the result when $J\le N$.
When $J>N$, duality of the finite matrix
$(E(\theta za))_{z,a}$ interchanges the two boxes and proves the
bound with $J$ in place of $N$. Fixed translates merely insert
unit-modulus row and column factors; subsets can be handled by
zero extension or deletion of rows.
\end{proof}

\begin{lemma}[Counting small norms in an ideal]
\label{lem:g789-smallnorm}
Let $\mathfrak b$ be an integral ideal, $b=N\mathfrak b$, and
$H\ge1$. For $U\ge b$ one has
\[
 \#\{h\in\mathfrak b\cap\mathcal B(H):0<|Nh|\le U\}
       \ll_{K,\varepsilon}H^\varepsilon U/b.
\]
There are no such nonzero points when $U<b$.
\end{lemma}
\begin{proof}
For a counted element, $(h)\mathfrak b^{-1}$ is an integral ideal
of norm at most $U/b$. There are $O_K(U/b)$ such ideals.
For any fixed principal ideal $(h)$, its generators in the box
number $O_K((1+\log H)^{r_1+r_2-1})$: upper bounds for all conjugates,
together with $|Nh|\ge1$, confine their logarithmic vectors to a box
of side $O_K(1+\log H)$, and the unit logarithms form a fixed lattice.
This proves the asserted estimate, including uniformity in
$\mathfrak b$.
\end{proof}

\begin{proposition}[Refined cubic coefficient sampling over $K$]
\label{prop:g789-sampling}
Let $H,J,N\ge1$, put $V=H^3JN$, and choose
\[
 1\le Y\le c_K\min\{J,N,(JN/H^3)^{1/2}\}.
\]
Choose a lifted approximation $\alpha=\gamma+\beta$ with
$L=Q_\gamma A_V(\beta)\ll_K V^{1/2}$, as in
Lemma~\ref{lem:g789-height}. If $c(a)$ is supported on $\mathcal B(N)$,
then
\[
 \sum_{\substack{h\in\mathcal B(H)\\h\ne0}}
 \sum_{z\in\mathcal B(J)}
       \left|\sum_a c(a)E(\alpha h^3za)\right|^2
 \ll_{K,\varepsilon}
 V^\varepsilon JHN\left(Y^{-1}+\frac{H^2}{L}\right)
                          \sum_a|c(a)|^2.
\]
In particular $L$ can be replaced by the smaller intrinsic height
$\mathcal H_V(\alpha)$. In
applications the variables are products of nonzero shifts.
\end{proposition}
\begin{proof}
Write $W=JN$ and apply Lemma~\ref{lem:g789-bilinear} with
$\theta=\alpha h^3$. Since $Y\le\min(J,N)$, the contribution
of $J+N$, and that of all $h$ with $\mathcal H_W(\alpha h^3)\ge Y$,
is $O_K(JHN/Y)\sum|c(a)|^2$.

For each remaining $h$, choose
$\alpha h^3=\gamma_h+\beta_h$ with
$Q_{\gamma_h}A_W(\beta_h)<Y$. The center $\gamma_h/h^3$ is an
approximation to $\alpha$ at scale $V$ with height at most $C_KH^3Y$.
Indeed its conductor norm is at most $|Nh|^3Q_{\gamma_h}$, and
\[
 |Nh|^3 A_V(\beta_h/h^3)
 =\prod_v\max\{|h_v|^3,V^{1/n}|(\beta_h)_v|\}^{d_v}
 \ll_K H^3 A_W(\beta_h).
\]
The upper bound on $Y$ makes the product of this height and $L$
smaller than the separation constant times $V$.
Lemma~\ref{lem:g789-height} therefore gives
$\gamma_h=h^3\gamma$ and $\beta_h=h^3\beta$ exactly.

Let $\mathfrak a$ be the conductor of $\gamma$, $Q=N\mathfrak a$,
and set $A=A_V(\beta)$. The conductor of $h^3\gamma$ is
$\mathfrak a/(\mathfrak a+(h)^3)$. Group the nonzero $h$ according
to the ideal $\mathfrak d=\mathfrak a+(h)^3$.
There are $O_{K,\varepsilon}(V^\varepsilon)$ possible divisor ideals.
For each fixed class put
\[
 \mathfrak b=\prod_{\mathfrak p}
       \mathfrak p^{\lceil v_\mathfrak p(\mathfrak d)/3\rceil},
 \qquad b=N\mathfrak b.
\]
Then $h\in\mathfrak b$ and $N\mathfrak d\le b^3$.
Since $|h_v|\ll_KH^{1/n}$ at every place, one has
\[
 A_W(h^3\beta)\gg_K
       \max\{1,A|Nh|^3/H^3\}.
\]
Consequently the reciprocal-height sum in this class is at most
\[
 \frac{b^3}{Q}
 \sum_{0\ne h\in\mathfrak b\cap\mathcal B(H)}
       \min\{1,C_KH^3/(A|Nh|^3)\}.
\]
Put $T=C_KH/A^{1/3}$, enlarging the fixed constant if necessary.
Lemma~\ref{lem:g789-smallnorm} and dyadic summation bound the inner
sum by
\[
 O_{K,\varepsilon}\left(H^\varepsilon
 \begin{cases}T/b,&T\ge b,\\(T/b)^3,&T<b.\end{cases}\right).
\]
In the first case $b^2T\le T^3$, and in the second case the factor
$b^3$ cancels. Thus either case contributes at most
$O_{K,\varepsilon}(V^\varepsilon H^3/(QA))$.
Summing the divisor classes and recalling $QA=L$ proves
\[
 \sum_{h:\,\mathcal H_W(\alpha h^3)<Y}
       \mathcal H_W(\alpha h^3)^{-1}
 \ll_{K,\varepsilon}V^\varepsilon H^3/L.
\]
Multiplication by $JN\sum|c(a)|^2$ completes the proof.
\end{proof}

\subsection{A second coefficient moment}
\label{sec:k5-coeff2}

All size parameters in this section denote volumes.  The field is fixed,
$n=[K:\mathbb Q]$, and $\mathcal H_T$ denotes the conductor--product height
defined in Section~\ref{sec:preliminaries}.  The notation $\mathcal B(L)$ means a fixed
archimedean box with side lengths comparable to $L^{1/n}$.  The same estimates hold on fixed fractional lattices.

The following is the precise number-field form of the coefficient estimate
needed from Vaughan--Wooley's fifth-power argument \cite[Lemma~4.3]{VW1995}.  It is useful to state
it for the general coefficient degree $2d$.

\begin{proposition}\label{prop:k5-coeff2}
Let $d\ge3$, $J,H,M\ge1$, and $T=JH^3M^{2d}$.  Suppose
\[
 1\le U\le c_K\min\{M,JH^{-3},T^{1/(4d)},
                  (JH^3M^{1-2d})^{1/6}\}.
\]
Let $\mathcal C$ be any subset of a sufficiently small fixed box about a
nonzero point at every embedding, of scale $M^{1/n}$.  For each $h$ one may
include any fixed lower-degree additive twist in the $m$ sum, provided it
is independent of $j$.  Then
\[
 E_1(\alpha):=\sum_{\substack{h\in\mathcal B(H)\\h\ne0}}
 \sum_{j\in\mathcal B(J)}
 \left|\sum_{m\in\mathcal C}e_K(\alpha jh^3m^{2d}+b_h(m))\right|^2
 \ll_{K,d,\epsilon}T^\epsilon JHM^2
       \{\mathcal H_T(\alpha)^{-1/d}+U^{-1}\}.
\]
The same bound holds for smooth weights and for the fixed lattices arising
from class representatives.  No positivity condition at complex places is
being imposed.
\end{proposition}

\paragraph{Conventions for the proof.}
A nonzero algebraic integer has absolute norm at least one.  After clearing
a fixed denominator this also applies, with a fixed positive lower bound,
to every lattice occurring below.  An arbitrary conductor ideal can be
made principal by multiplying it by one of finitely many fixed class
representatives.  Its chosen denominator therefore has norm comparable
to the conductor norm.  In all congruence counts it is the conductor ideal,
not that auxiliary denominator, which is used.

We use anisotropic Dirichlet approximation in the following elementary
form.  Given positive radii $R_v$ with $\prod_vR_v^{d_v}\gg_K1$, Minkowski's
theorem supplies $a\in\mathfrak D^{-1}$ and nonzero $q\in\mathcal O_K$ with
\[
 |q_v|\ll_KR_v,\qquad |(q\beta-a)_v|\ll_K R_v^{-1}.
\]
The constants can be made harmless by reducing the fixed constant $c_K$
in the proposition.  All later identities follow by estimating the norm
of a cross product and showing that it is smaller than its nonzero lower
bound.  This replaces the rational proof's use of a nonzero integer of
absolute value less than one.

A smooth positive majorant of a linear sum over $\mathcal B(J)$ is a
periodic rapidly decreasing kernel $D_J$.  Dyadic boxes for its errors,
and dyadic boxes for all nonzero variables, cost only $T^\epsilon$.
Thus the estimates below may be proved with specified radii at every
embedding.  Products of these radii, with local dimensions $d_v$, give
the volume expressions displayed below.  For example the contribution
where the error kernel has normalized product size at least $U$ is
$O(T^\epsilon JHM^2/U)$.

\begin{proof}
Expanding the square and summing over $j$, using a smooth majorant before
taking absolute values, reduces the assertion to the untwisted positive
kernel sum.  Terms with $m_1=m_2$ contribute $O(JHM)$ and are acceptable.
Because the $m$ box is sufficiently small, $m_1^{2d}=m_2^{2d}$ within this
box implies $m_1=m_2$.  Put $D_0=m_2^{2d}-m_1^{2d}$ for the other pairs.

\emph{First identification: the cubic shift.}
For a fixed pair and an error box of normalized product size at most $U$,
approximate $\alpha D_0$ with denominator radii
$(J_vH_v^3/U_v)^{1/2}$, where $\prod U_v^{d_v}=U$.
If $\alpha D_0h^3$ is in that error box around $a_0\in\mathfrak D^{-1}$,
then the norm of the cross product between this approximation and $a_0/h^3$
is
\[
 \ll_K (UH^3/J)^{1/2}.
\]
It is therefore zero.  Write the resulting common center as $\gamma_{D_0}$.
Its conductor $\mathfrak q_{D_0}$ divides $(h)^3$.
The sum over such $h$, using ideal-lattice packing on dyadic boxes, is at most
\[
 T^\epsilon JH V_{D_0}^{-1/3},\qquad
 V_{D_0}=N\mathfrak q_{D_0}\,
 A_{JH^3}(\alpha D_0-\gamma_{D_0}).
\]
Here is a lattice proof of the summation, including thin boxes.  If
$\mathfrak b=\prod\mathfrak p^{\lceil v_\mathfrak p(\mathfrak q)/3\rceil}$,
then $h\in\mathfrak b$ and $N\mathfrak b\ge(N\mathfrak q)^{1/3}$.
For $A=A_{JH^3}(\beta)$ one has
$A_J(h^3\beta)\gg\max(1,A|Nh|^3/H^3)$.
The elementary ideal count
$\#\{h\in\mathfrak b\cap\mathcal B(H):0<|Nh|\le B\}
\ll H^\epsilon B/N\mathfrak b$, with no points for $B<N\mathfrak b$,
and dyadic summation give
$\sum_h A_J(h^3\beta)^{-1}\ll H^\epsilon H/(A^{1/3}N\mathfrak b)$.
This proves the displayed bound without a boundary or real-positivity assumption.
We may discard the pairs for which $V_{D_0}\ge U^3$ at the required cost.

\emph{Second identification: remove the divided difference.}
Group pairs according to their greatest common divisor ideal
$\mathfrak g=(m_1,m_2)$, of norm $G$.  Make $\mathfrak g$ principal using a
fixed class representative and choose a balanced generator $g$.
Write $m_1=gn$, $m_2=g(n+l)$.  The variables $n,n+l$ belong to fixed fractional
lattices, are primitive up to fixed primes, and have scale $N=M/G$.
The polynomial
\[
 D(n,l)=\frac{(n+l)^{2d}-n^{2d}}{l}
\]
has absolute value comparable to $N_v^{2d-1}$ at every embedding.
Indeed, both arguments are in the same tiny box away from zero, so the
leading derivative $2d n^{2d-1}$ dominates the other terms.
This is the place where localization replaces real positivity.

For fixed $g,l$, approximate $\alpha g^{2d}l$ with denominator radii
$(J_vH_v^3N_v^{2d-1})^{1/2}$, obtaining a center $\gamma_{g,l}$.
A comparison with $\gamma_{D_0}/D(n,l)$ has cross-product norm at most
\[
 \ll_K V_{D_0}
       \left(\frac{N^{2d-1}}{JH^3}\right)^{1/2}
 \le C_K U^3\left(\frac{M^{2d-1}}{JH^3}\right)^{1/2}<1.
\]
Consequently $\gamma_{D_0}=D(n,l)\gamma_{g,l}$.
If $\mathfrak q$ is the conductor of $\gamma_{g,l}$, then
$\mathfrak q_{D_0}=\mathfrak q/\mathfrak q_1$, with
$\mathfrak q_1\mid(D(n,l))$.
Moreover $n(n+l)$ is prime to $\mathfrak q_1$ outside a fixed set of primes.
The number of primitive solutions of $D(n,l)\equiv0\pmod{\mathfrak q_1}$
in residue classes is $O_{K,d,\epsilon}((N\mathfrak q_1)^{\epsilon})$.
To see the uniformity even at primes dividing $2d$, let
$t=v_{\mathfrak p}(l)$.  For a primitive solution modulo $\mathfrak p^E$,
both $n$ and $n+l$ are units.  The map
$n\mapsto r=1+l/n\pmod{\mathfrak p^{E+t}}$ is injective, and
$D(n,l)\equiv0\pmod{\mathfrak p^E}$ implies
$r^{2d}\equiv1\pmod{\mathfrak p^{E+t}}$.
The number of unit $2d$th roots modulo any power of a fixed prime is
bounded in terms of $K,d$: beyond the fixed derivative depth
$v_{\mathfrak p}(2d)$, Hensel lifting gives only a bounded number of
lifts around each root.  At good primes the bound is $2d$.
Multiplying the local bounds gives the stated $T^\epsilon$ factor.
Fixed denominators from class representatives are cleared first; their
fixed valuations only change the constant.

Put
\[
 V_{g,l}=N\mathfrak q\,
 A_{JH^3N^{2d-1}}(\alpha g^{2d}l-\gamma_{g,l}).
\]
The preceding comparison of $D(n,l)$ with its archimedean scale gives
$V_{D_0}\asymp_K V_{g,l}/N\mathfrak q_1$.
Ideal packing gives at most
$T^\epsilon(1+N/N\mathfrak q_1)$ choices for $n$ in each possible class.
The terms $N\mathfrak q_1>N$ contribute $O(T^\epsilon JHM)$ in total:
there are $O(T^\epsilon)$ possible $n$ for each $g,l$ and divisor class,
and their weight is $V_{D_0}^{-1/3}\le1$;
in each dyadic $G$ range the number of $g,l$ pairs is $O(T^\epsilon M)$.
For the rest, $(N\mathfrak q_1)^{1/3}/N\mathfrak q_1\le1$.
Thus the remaining sum is bounded by
\[
 T^\epsilon JHM\sum_gG^{-1}\sum_l V_{g,l}^{-1/3}.
\]
Here and below the summation over $g$ means one balanced generator per
ideal and finitely many fixed class choices; its count in a dyadic range
of $G$ is $O_K(GT^\epsilon)$.
We may further restrict to $V_{g,l}<(U/G)^3$, at a cost
$O(T^\epsilon JHM^2/U)$.

\emph{Third identification: remove the difference.}
For fixed $g$, approximate $\alpha g^{2d}$ with denominator radii
$(J_vH_v^3N_v^{2d})^{1/2}$, obtaining $\gamma_g$.
Comparison with $\gamma_{g,l}/l$ has cross-product norm
\[
 \ll_K\frac{V_{g,l}}{(JH^3)^{1/2}N^{d-1}}<1.
\]
The previous bound on $V_{g,l}$ and the hypotheses on $U$ ensure the last
inequality.  Hence $\gamma_{g,l}=l\gamma_g$ exactly.
Splitting according to the divisor ideal of the conductor that divides $l$,
and summing $l$ on dyadic archimedean boxes, gives
\[
 \sum_l V_{g,l}^{-1/3}\ll T^\epsilon N V_g^{-1/3},\qquad
 V_g=Q_{\gamma_g}A_{JH^3N^{2d}}(\alpha g^{2d}-\gamma_g).
\]
For an explicit verification, put $A=A_{JH^3N^{2d}}(\beta)$.
Then $A_{JH^3N^{2d-1}}(l\beta)\gg A|Nl|/N$.
If $l$ belongs to a divisor ideal $\mathfrak b$, the preceding ideal norm
count gives
$\sum_{0\ne l\in\mathfrak b\cap\mathcal B(N)}|Nl|^{-1/3}
\ll T^\epsilon N^{2/3}/N\mathfrak b$.
Combining this with the conductor factor yields divisor weights
$(N\mathfrak b)^{-2/3}$, whose sum is $O(T^\epsilon)$.
This proves the displayed linear-variable estimate.
We have therefore reduced to
\[
 T^\epsilon JHM^2\sum_gG^{-2}V_g^{-1/3}.
\]
Terms with $V_g\ge(U/G)^3$ are again acceptable, so assume the reverse
inequality.  In particular $G<U$.

\emph{Fourth identification: remove the common divisor.}
Choose one Dirichlet center $\gamma$ for $\alpha$ at scale $T$.
Comparison with $\gamma_g/g^{2d}$ has cross-product norm at most
\[
 \ll_K V_gG^{2d}T^{-1/2}
 \le C_KU^3G^{2d-3}T^{-1/2}
 \le C_KU^{2d}T^{-1/2}<1.
\]
Thus $\gamma_g=g^{2d}\gamma$.  Write $Q=Q_\gamma$ and $A=A_T(\alpha-\gamma)$.
The conductor of $g^{2d}\gamma$ is obtained by dividing the conductor of
$\gamma$ by its greatest common divisor with $(g)^{2d}$.
The chosen generators are balanced, so the archimedean factor in $V_g$
is comparable to $A$.
Finally
\[
 \sum_{\mathfrak g}(N\mathfrak g)^{-2}
 \left(\frac{Q}{N(\mathfrak q+(\mathfrak g)^{2d})}\right)^{-1/3}
 \ll_{K,d,\epsilon}T^\epsilon Q^{-1/d}.
\]
This is an elementary Euler-product calculation: at a prime with
$v_\mathfrak p(\mathfrak q)=e$, the term of valuation $a$ has size
$q^{-2a-(e-\min(e,2da))/3}$, and its largest possible size is
$O(q^{-e/d})$ because $d\ge3$.  Summing the geometric tails and absorbing
the bounded factors per prime into $T^\epsilon$ proves the claim.
Since $A^{-1/3}\le A^{-1/d}$,
\[
 \sum_gG^{-2}V_g^{-1/3}\ll T^\epsilon(QA)^{-1/d}
 \le T^\epsilon\mathcal H_T(\alpha)^{-1/d}.
\]
Combining the four reductions proves the proposition.
\end{proof}

\subsection{Sampled quadratic moments}
\label{sec:quad-sampling-field}

Let $K$ be fixed, $n=[K:\mathbb Q]$, $d_v=[K_v:\mathbb R]$,
$V=K\otimes_{\mathbb Q}\mathbb R$, and
$\mathbb T=V/\mathfrak D^{-1}$. All size parameters in this section
are volumes. Write $\mathcal B(U)$ for a bounded archimedean box with
all radii $O(U^{1/n})$. Fixed dilations and translates are harmless.
For a vector $\mathbf T=(T_v)_{v\mid\infty}$ of positive real numbers put
\[
 A_{\mathbf T}(\beta)=\prod_v\max\{1,T_v|\beta_v|\}^{d_v},\qquad
 \mathcal H_{\mathbf T}(\alpha)=
       \min_{\gamma\in K} Q_\gamma A_{\mathbf T}(\alpha-\gamma),
 \qquad
 Q_\gamma=N\{a\in\mathcal O_K:a\gamma\in\mathfrak D^{-1}\}.
\]
When every $T_v=T^{1/n}$ we abbreviate the height to $\mathcal H_T$.
The minimum includes all lifts of a center. Product Dirichlet approximation
gives $\mathcal H_T(\alpha)\ll_K T^{1/2}$.

\begin{lemma}[Anisotropic bilinear and cluster estimates]
\label{lem:anisotropic-bilinear}
Let $\mathcal U,\mathcal W$ be archimedean boxes of radii $U_v,W_v$,
with volumes $J=\prod_vU_v^{d_v}\ge1$ and
$N=\prod_vW_v^{d_v}\ge1$. Put $T_v=U_vW_v$.
For any coefficients supported on $\mathcal W\cap\mathcal O_K$,
\[
 \sum_{z\in\mathcal U\cap\mathcal O_K}
  \left|\sum_m c_m e_K(\alpha zm)\right|^2
 \ll_K \left(J+N+\frac{JN}{\mathcal H_{\mathbf T}(\alpha)}\right)
               \sum_m|c_m|^2.
\]
Consequently the number of $z$ in $\mathcal U$ for which $\alpha z$
modulo $\mathfrak D^{-1}$ lies in a prescribed box of radii $O(W_v^{-1})$
is
\[
 \ll_K 1+\frac JN+\frac J{\mathcal H_{\mathbf T}(\alpha)}.
\]
The constants are independent of the eccentricities of the two boxes.
\end{lemma}
\begin{proof}
The proof of Lemma~\ref{lem:g789-bilinear} works
with these radii without alteration. We give the relevant uniformity details.
For $J\le N$, subdivide $\mathcal U$ by a fixed number of subdivisions
at each embedding. The product of the sample-box volume and a
frequency-error box volume is then sufficiently small for the affine
field-line argument: the determinant of any two difference pairs has
norm smaller than the fixed lower bound for a nonzero element of
$\mathfrak D^{-1}$. A cluster therefore lies on an affine $K$-line of
slope $\gamma$. Ideal-lattice packing bounds its size by
\[
 1+Q_\gamma^{-1}
       \prod_v\min\{U_v,W_v^{-1}/|\alpha_v-\gamma_v|\}^{d_v}
 =1+\frac J{Q_\gamma A_{\mathbf T}(\alpha-\gamma)}.
\]
The Gaussian large sieve then gives the stated bound, with
$N+JN/\mathcal H_{\mathbf T}(\alpha)$ on the right. When $J>N$,
transpose the bilinear matrix, interchanging the two boxes. Its operator
norm is unchanged and the product radii $T_v$ are unchanged. This proves
the displayed symmetric bound.

For the cluster assertion, put coefficient weight
$\exp(-c\sum_v|m_v/W_v|^2)$ on the dual box. Poisson summation gives a
nonnegative Gram kernel which is $\gg_K N$ for any pair of the selected
sample frequencies. If there are $M$ samples, the corresponding quadratic
form is therefore $\gg_K NM^2$. Its upper bound from the bilinear operator
estimate is $O_K((J+N+JN/\mathcal H_{\mathbf T})M)$.
Gaussian tails are summed by dyadic dilation and rapid decay. Division
gives the claim. The only lattice input throughout is
$\#(\mathfrak a\cap\mathcal R)\ll_K1+\operatorname{vol}(\mathcal R)/N\mathfrak a$
for an arbitrary product box, so no bounded-eccentricity hypothesis is used.
\end{proof}

\begin{proposition}[Distribution of short heights of multiples]
\label{prop:short-height-multiples}
Let $X,H\ge1$, put $V_0=HX^2$, and let
$L=\mathcal H_{V_0}(\alpha)$. For $1\le Y\ll_K X$ one has
\[
 \#\{0\ne h\in\mathcal B(H):\mathcal H_{X^2}(h\alpha)\le Y\}
 \ll_{K,\epsilon} V_0^\epsilon
       \left(\frac{HY^2}{X^2}+\frac{HY}{L}\right).
\]
\end{proposition}
\begin{proof}
Choose an approximation $h\alpha=\gamma+\beta$ of height at most $Y$.
Let $\mathfrak q$ be its conductor and write $Q=N\mathfrak q$.
Minkowski's theorem supplies $0\ne b\in\mathfrak q$ with
$|b_v|\ll_KQ^{1/n}$ at every embedding. Since $(b)\subset\mathfrak q$,
$|Nb|\ge Q$; consequently also $|b_v|\gg_KQ^{1/n}$ at every embedding.
Thus $|Nb|\asymp_K Q$ and $b\gamma\in\mathfrak D^{-1}$.
Put $m=bh$. It lies in a box of volume $O_K(HQ)$ and satisfies
\[
 \alpha m\equiv b\beta\pmod{\mathfrak D^{-1}}.
\]
Each nonzero $m$ has only $V_0^\epsilon$ possible bounded-height factors
$b$, by the ideal-divisor bound and the unit-lattice count.

Decompose $Q$ and
$\max(1,X^{2/n}|\beta_v|)$ into dyadic intervals. Denote the latter
upper endpoints by $A_v\ge1$ and their product by
$A=\prod_v A_v^{d_v}\ll_KY/Q$.
There are only logarithmically many such choices. For a fixed choice,
the residue error $b\beta$ is contained in a box of radii
$O_K(Q^{1/n}X^{-2/n}A_v)$. Apply
Lemma~\ref{lem:anisotropic-bilinear} with
\[
 U_v=(HQ)^{1/n},\qquad
 W_v=(X^2/Q)^{1/n}A_v^{-1},\qquad
 J=HQ,\quad N=\frac{X^2}{QA}.
\]
Here $N\gg_KX\ge1$, after harmless fixed rescaling, and
$T_v=V_0^{1/n}/A_v$.
For every approximation center,
$A_{V_0}(\alpha-\gamma)\le A\,A_{\mathbf T}(\alpha-\gamma)$.
Taking minima gives $L\le A\mathcal H_{\mathbf T}(\alpha)$.
The number of such $m$ is therefore
\[
 \ll_K1+\frac{HQ^2A}{X^2}+\frac{HQA}{L}
 \ll_K1+\frac{HY^2}{X^2}+\frac{HY}{L}.
\]
Finally, the approximation to $h\alpha$ lifts to one to $\alpha$
with center $\gamma/h$. Its conductor has norm at most $Q|Nh|$, and
\[
 |Nh|A_{V_0}(\beta/h)
   =\prod_v\max\{|h_v|,V_0^{1/n}|\beta_v|\}^{d_v}
   \ll_K H A_{X^2}(\beta).
\]
Thus the set being counted is empty unless $L\ll_K HY$.
On its nonempty range this absorbs the constant $1$ above.
Summing the dyadic choices and factor multiplicities proves the result.
\end{proof}

\begin{corollary}[Sampled fourth moment]
\label{cor:quadratic-fourth-field}
Let $Z_t(\beta)$ be a quadratic exponential sum on a box of root volume
$X$, with arbitrary linear modulation $t$. For arbitrary modulations $t_h$,
\[
 \sum_{0\ne h\in\mathcal B(H)}|Z_{t_h}(h\alpha)|^4
 \ll_{K,\epsilon}(HX^2)^\epsilon
 \left\{HX^2+X^4\min\left(\frac HL,\frac{H^2}{L^2}\right)\right\},
 \qquad L=\mathcal H_{HX^2}(\alpha).
\]
The assertion is uniform for fixed nonzero integral multipliers of $\alpha$.
\end{corollary}
\begin{proof}
One ordinary difference removes the linear modulation. The positive
linear-kernel estimate, proved by the bilinear Gram argument above, gives
\[
 |Z_t(\beta)|^2\ll X^\epsilon
       \left(X+\frac{X^2}{B}\right),\qquad
 B=\mathcal H_{X^2}(\beta)\ll_K X.
\]
Sharp fixed boxes can be completed with a logarithmic loss; smooth majorants
give the same estimate directly. In particular the bound is uniform in $t$.
Hence $|Z_t(\beta)|^4\ll X^\epsilon(X^2+X^4B^{-2})$.
Proposition~\ref{prop:short-height-multiples}, summed dyadically from
$Y\asymp_K\max(1,L/H)$ to $O_K(X)$, shows that
\[
 \sum_{0\ne h\in\mathcal B(H)}\mathcal H_{X^2}(h\alpha)^{-2}
 \ll (HX^2)^\epsilon
       \left\{\frac H{X^2}
        +\min\left(\frac HL,\frac{H^2}{L^2}\right)\right\}.
\]
Combining the last two bounds proves the assertion. A fixed multiplier
changes conductor heights only by fixed factors.
\end{proof}

In the eighth-power terminal, $H=X^h$ and $L\asymp X^\ell$.
Normalization by $HX^4$ gives
\[
 z_4(\ell)=\min\{2,\max(\ell,2\ell-h)\}.
\]
Both resources used by the exact policy,
$\min\{2,2\ell-h\}$ and $\min\{2,\ell,h\}$, are weaker and therefore valid.
Products of a fixed number of shifts may replace $h$, at the cost of the
same divisor loss. The quadratic second-moment resource remains the
previously proved bound with normalized saving $1$ when $H\ge X$ and
$L\gg X^{1-O(\eta)}$.

\section{Restricted-height smooth means}
\label{sec:restricted-field}

The following field version of Br\"udern--Wooley's restricted-height
mean theorem is needed by the eighth- and ninth-power routes.  The use of
ideals in its proof is particularly useful: an extracted ideal factor
need not be principal.  The quotient variables instead lie in one of
finitely many fixed integral ideals.

Write $V=K\otimes\mathbb R$, $n=[K:\mathbb Q]$, and $d_v=[K_v:\mathbb R]$.
For $T\ge1$ put
\[
 A_T(\beta)=\prod_{v\mid\infty}
                 \max(1,T^{1/n}|\beta_v|)^{d_v}.
\]
For a fixed integral ideal $C$, use the coefficient torus
$\mathbb T_C=V/(C^{-k}\mathfrak D^{-1})$ with probability Haar measure.
The conductor of $\gamma\in K/(C^{-k}\mathfrak D^{-1})$ is
\[
 \mathfrak a_C(\gamma)
       =\{a\in\mathcal O_K:a\gamma C^k\subseteq\mathfrak D^{-1}\}.
\]
Let $\mathfrak M_C(Z;Y)$ be the union of the product arcs
\[
 (N\mathfrak a_C(\gamma)) A_{Y^k}(\alpha-\gamma)\le Z.
\]
All lifts of a center are included.  A subscript $\mathfrak a$ on this
set restricts to centers of exact conductor $\mathfrak a$.
For $C=\mathcal O_K$, suppress the subscript $C$.

Let $g(\alpha;X,R)$ be a sum over the nonzero algebraic integers in a fixed
bounded box of scale $X^{1/n}$ whose prime ideal factors have norm at most
$R$.  Nonnegative bounded weights and smaller fixed boxes are allowed.

\begin{proposition}[Field restricted-height mean]
\label{prop:restricted-field}
Let $s\ge2$ be even.  Suppose $\lambda=s-k+\Delta$ is a permissible
full smooth-moment exponent, with
$\lambda\ge\max(s/2,s-k)$ and $\Delta\ge0$.  Then for every
$\epsilon>0$ a sufficiently small smoothness exponent gives
\[
 \int_{\mathfrak M(Z;X)}|g(\alpha;X,R)|^s\,d\alpha
       \ll_{K,k,s,\epsilon} X^{s-k+\epsilon}Z^{2\Delta/k}
       \qquad(1\le Z\le X^{k/2}).
\]
The same conclusion holds at intermediate real orders when their full
exponents are obtained by H\"older interpolation between such even orders.
\end{proposition}

\subsection{Fixed root lattices and multiplication of arcs}

Choose integral ideal representatives $C_1,\ldots,C_h$ for the ideal
classes.  For every nonzero integral ideal $\mathfrak b$ in the class of
$C_i$, choose $b=b(\mathfrak b)\in K^\times$ with
\[
 \mathfrak b=b C_i,\qquad
 |b_v|\asymp_K(N\mathfrak b)^{1/n}\quad(v\mid\infty).
\]
Dirichlet's unit theorem gives the second condition after choosing one
representative in a fundamental domain for the unit logarithms.  This
choice is made once for each ideal.  Multiplication by $b$ maps $C_i$
into $\mathcal O_K$.

If $(\mathfrak a,\mathfrak b)=1$, multiplication by $b^k$ maps centers
of exact conductor $\mathfrak a$ on $\mathbb T_{\mathcal O_K}$ bijectively
to those of exact conductor $\mathfrak a$ on $\mathbb T_{C_i}$.
Indeed the conductor calculation is local, and
$b^kC_i^k=\mathfrak b^k$ is a unit ideal at every prime dividing
$\mathfrak a$.  Equivalently, multiplication by $b^k$ is an isomorphism
on these finite conductor-character groups.
Its archimedean Jacobian, with the indicated probability measures, is
$(N\mathfrak b)^k$.  To check the normalization, ordinary Lebesgue measure
has Jacobian $|Nb|^k=(N\mathfrak b/NC_i)^k$, while the volume of
$\mathbb T_{C_i}$ is $(NC_i)^{-k}$ times that of $\mathbb T_{\mathcal O_K}$.

Consequently, whenever the image arcs are disjoint, for any nonnegative
integrand $F$ one has
\begin{equation}\label{eq:field-arc-scaling}
 \int_{\mathfrak M_{\mathfrak a}(Z;X)}F(b^k\alpha)\,d\alpha
 \ll_K (N\mathfrak b)^{-k}
 \int_{\mathfrak M_{C_i,\mathfrak a}(c_1Z;X/N\mathfrak b)}
                         F(\theta)\,d\theta .
\end{equation}
Here $c_1$ is a fixed constant, independent of the ideal and its class.
Balancing $b$ makes the archimedean scaling comparable at every place.
The requisite disjointness holds if
$c_1Z\le c_2Y^{k/2}$, with a sufficiently small $c_2=c_2(K,k)$,
where $Y=X/N\mathfrak b$.
For distinct lifted centers the nonzero norm of their difference is
$\gg_K(QQ')^{-1}$, while overlapping product arcs would give an upper
bound $\ll_KY^{-k}Z^2/(QQ')$.  This also excludes overlap of distinct
lifts of the same center.  There are only finitely many root lattices
in use, so one constant suffices for all of them.

\subsection{The ideal-factor decomposition}

Fix a total ordering of the prime ideals which refines their norm ordering.
Given a conductor ideal $\mathfrak a$, factor
\[
 (x)=\mathfrak v\mathfrak u,
 \qquad (\mathfrak v,\mathfrak a)=1,
 \qquad \mathfrak u\text{ supported on primes dividing }\mathfrak a.
\]
Call a root bad when $N\mathfrak v\le M$.  The supported-ideal counting
lemma from the fundamental differencing argument, followed by bounded-height
unit counting, gives at most $X^\epsilon M$ such roots, uniformly for
$N\mathfrak a\le X^{k/2}$.  This remains true when the variable lies in
a smaller box or has a nonnegative weight bounded by a fixed constant.

For a good root, multiply the prime factors of $\mathfrak v$ in decreasing
order until their product first has norm exceeding $M$.  Denote the
resulting ideal by $\mathfrak m$ and its last prime factor by $\mathfrak p$.
This decomposition is unique, and
\[
 M<N\mathfrak m\le MR,\qquad
 (\mathfrak m,\mathfrak a)=1.
\]
Every prime of the remaining ideal $\mathfrak w=\mathfrak v/\mathfrak m$
is at most $\mathfrak p$ in the fixed ordering.  The condition on
$\mathfrak m$ is exactly that its factors are at least $\mathfrak p$ and
$N(\mathfrak m/\mathfrak p)\le M<N\mathfrak m$.
Thus for each outer $\mathfrak m$, the remaining $\mathfrak w$ runs over
the complete set of $\mathfrak p$-smooth ideals prime to $\mathfrak a$.

Write $\mathfrak m=m C_i$ using the balanced choice above.  Then $x=mr$
with $r\in C_i$ and
$(r)C_i^{-1}=\mathfrak w\mathfrak u$.
Split $\mathfrak u$ into its prime factors at most $\mathfrak p$ and those
larger than $\mathfrak p$; call the latter product $\mathfrak z$.
The smaller part combines with $\mathfrak w$ to give every
$\mathfrak p$-smooth integral ideal, without a coprimality restriction.
Writing $\mathfrak z=z C_j$, one obtains $r=zy$ with
\[
 y\in C_iC_j,
 \qquad (y)(C_iC_j)^{-1}\text{ a }\mathfrak p\text{-smooth integral ideal}.
\]
The original box condition is simply its original weight evaluated at
$mzy$.  For fixed $\mathfrak m,\mathfrak z$ this inner sum is independent
of the conductor ideal.  All root lattices $C_i$ and $C_iC_j$ are fixed;
there are only finitely many of them.  This is the exact replacement for
the rational factor decomposition in \cite[Sections~3--4]{BW2023}.

\subsection{Integration and the exponent}

Choose a sufficiently small fixed $c>0$ and put
\[
                  M=cX Z^{-2/k}R^{-1}.
\]
If $M<R$, the desired assertion follows directly from the full moment:
$Z\gg_K X^{k/2}R^{-k}$, so the loss is at most $R^{2\Delta}$,
absorbed by choosing the smoothness exponent small enough.
We may therefore assume $M\ge R$.
There are $O_K(MR)$ choices of the ideal $\mathfrak m$.
H\"older's inequality applied to the preceding exact decomposition gives
a factor $(MR)^{s-1}$ and a sum of $s$th moments of the inner sums.
The bad part contributes at most
\begin{equation}\label{eq:restricted-bad}
                       X^{-k+\epsilon}Z^2M^s.
\end{equation}
Here the total measure of the product arcs is
$\ll_KX^{-k}Z^2(1+\log X)^{O_K(1)}$: a conductor of norm $q$ has at
most $q$ centers, each error region has volume
$\ll_K X^{-k}(Z/q)(1+\log Z)^{O_K(1)}$, and there are $O_K(Z)$ ideals
of norm at most $Z$.

For each fixed $\mathfrak m$, apply \eqref{eq:field-arc-scaling}.
Our choice of $c$ ensures disjointness on the quotient root scale because
$N\mathfrak m\le MR=cXZ^{-2/k}$.
There are only $X^\epsilon$ ideals $\mathfrak z$ supported on a given
conductor and of the relevant norm.  Another application of H\"older
therefore costs $X^\epsilon$.
Now extend the $\mathfrak z$ summation to all integral ideals in its
allowed size range.  For fixed $\mathfrak m,\mathfrak z$, the inner
integrand no longer depends on $\mathfrak a$.  The sum over the disjoint
conductor arcs is consequently at most a full torus integral.
Multiplication by $z^k$ maps $\mathbb T_{C_i}$ onto
$\mathbb T_{C_iC_j}$ and preserves probability Haar integrals.
The remaining full moment is bounded by
\[
 X^\epsilon\left(\frac{X}{N\mathfrak m N\mathfrak z}\right)^\lambda.
\]
Indeed its weights are nonnegative and bounded, its support lies in a
fixed box of that volume scale, and its roots lie in a fixed integral
ideal.  At even order, orthogonality bounds it by a constant times the
complete smooth moment on $\mathcal O_K$, after allowing the finitely
many fixed prime factors of $C_iC_j$.
This use of coefficient domination is only on a full torus integral.

The permissible bound is uniform for these smaller scales: for scales
$Y\ge X^\delta$, choose the original smoothness exponent at most
$\delta$ times a permissible one; for $Y<X^\delta$, use the trivial
bound and choose $\delta$ sufficiently small in terms of the allocated
$\epsilon$.  Bounded scales and the finitely many fixed primes cause
only a constant or $X^\epsilon$ factor.

Since $\lambda\ge s/2\ge1$, ideal counting gives
$\sum_{N\mathfrak z\le X}(N\mathfrak z)^{-\lambda}\ll_K X^\epsilon$,
including the logarithmic endpoint $\lambda=1$.
The good contribution is therefore at most
\[
 X^\epsilon(MR)^{s-1}
 \sum_{M<N\mathfrak m\le MR}(N\mathfrak m)^{-k}
             \left(\frac{X}{N\mathfrak m}\right)^\lambda
 \ll X^{s-k+\epsilon}(X/M)^\Delta.
\]
Inserting $M$ gives the asserted main bound.
The bad contribution \eqref{eq:restricted-bad} is
$\ll X^{s-k+\epsilon}Z^{2-2s/k}$, which is no larger because
$\Delta\ge k-s$ and $Z\ge1$.
H\"older interpolation proves the final assertion at real orders.
This completes the proof of Proposition~\ref{prop:restricted-field}.

\section{The cubic terminal}
\label{sec:cubic-terminal}
\label{sec:g7-cubic}
All parameters are volumes, and the conductor--product height is the one fixed in Section~\ref{sec:preliminaries}.
\subsection{A linear-kernel consequence of the bilinear large sieve}

The elementary bilinear sampling bound is
\begin{equation}\label{eq:g7-bilinear}
 \sum_{z\in\mathcal B(J)}\left|\sum_{a\in\mathcal B(N)}c(a)e_K(\theta za)\right|^2
 \ll_K\left(J+N+\frac{JN}{\mathcal H_{JN}(\theta)}\right)\sum_a|c(a)|^2.
\end{equation}
This is Lemma~\ref{lem:g789-bilinear}. The same bound holds for fixed translates, subsets, and fixed dilations of the boxes.

Here is the corresponding positive kernel estimate.  Let
\[
 D_N(t)=\sum_{z\in\mathcal O_K}
      \exp(-c|z|^2/N^{2/n})e_K(tz).
\]
Choose the Euclidean norm compatibly with the trace pairing.  Poisson
summation writes this as a positive sum of Gaussians centered on
$\mathfrak D^{-1}$, each of height $O_K(N)$ and width $N^{-1/n}$.
Then
\begin{equation}\label{eq:g7-kernel}
 \sum_{h\in\mathcal B(J)}D_N(\theta h)
 \ll_K J+N+\frac{JN}{\mathcal H_{JN}(\theta)}.
\end{equation}
Indeed, enlarge the $h$ box by a fixed factor.  Every difference in the
original box has $\gg_K J$ representations as $r-s$ in the enlarged box.
Because $D_N\ge0$, its sum over the original box is at most a constant times
$J^{-1}\sum_{r,s}D_N(\theta(r-s))$.  The latter double sum is a quadratic
form in the Gram matrix of $(e_K(\theta rz))_{r,z}$.  Its norm is bounded by
\eqref{eq:g7-bilinear}.  Gaussian tails in $z$ are handled by dyadic boxes;
monotonicity of $\mathcal H_T$ in $T$ and Gaussian decay make their total
convergent.  This proves \eqref{eq:g7-kernel}.
For Fourier transforms of smooth compactly supported weights, replace the
Gaussian kernel by a rapidly decreasing positive majorant; a dyadic
superposition of Gaussian kernels proves the same bound.  When the kernel
parameter decreases from $N$ to $N/2^r$, use that
$\mathcal H_t(\theta)/t$ is nonincreasing in $t$, directly from its
definition; consequently $(JN/2^r)/\mathcal H_{JN/2^r}(\theta)$
is at most $JN/\mathcal H_{JN}(\theta)$.  The rapid tail weights
absorb the corresponding powers of $2^r$.

\subsection{The cubic family}

Let $j=k-3$, and put
\[
 M_i=X^{\phi_i},\quad H_i=XM_i^{-k},\quad
 M=\prod_{i\le j}M_i,\quad H=\prod_{i\le j}H_i,\quad
 Q=X/M,\quad T=Q^k=X^3H.
\]
We require $T\gg_KX^2$, which is checked in each application.
The relevant auxiliary polynomial has the form
\[
 \Psi_j(z;\mathbf h;\mathbf m)
   =a(\mathbf h)\{c_k z^3+B(\mathbf h,\mathbf m)z\},\qquad
 a(\mathbf h)=h_1\cdots h_j,
\]
where $c_k\ne0$ is fixed and integral, and $B$ is integral after multiplication
by a fixed denominator depending only on $k$.  Such fixed denominators can be
cleared once, changing only finitely many bad prime ideals and fixed lattices.
We sum over nonzero $h_i$ in boxes of volume $O(H_i)$ and over $m_i$ in boxes
of volume $O(M_iX^\epsilon)$.  Retain the actual small-support restrictions
from efficient differencing.  In particular $z$ lies in a fixed small box
around a nonzero archimedean center of size $X^{1/n}$, and
$|B_v|\ll_KX^{2/n}$ uniformly.  The Hessian of the cubic phase in $z$ is
nonzero on that box at every archimedean place.

Let $F_j(\alpha)$ be this sum, allowing a bounded smooth weight in $z$ with
uniformly bounded scaled derivatives and any added linear archimedean phase.
To arrange this smoothness without a boundary hypothesis, start with a
nonnegative $C^\infty$ weight that is at least one on the distinguished
ordinary root box and supported in a slightly larger tiny box.  The initial
solution count is thereby enlarged.  Every Cauchy or differencing step
produces products of translates of this weight.  The affine changes in $z$
have coefficients bounded in terms of the fixed differencing depth, so
all derivatives remain uniformly bounded at scale $X^{1/n}$.  Its support
retains the small-increment restrictions and the nonzero center.  Extra
indicators may be enlarged to such products before the last absolute-value
estimate, since the underlying orthogonality count has nonnegative weights.
Alternatively, ordinary Fourier completion of sharp box intersections costs
only $X^\epsilon$; the estimates below are uniform in its linear frequencies.

Define
\[
 \mathfrak M=\{\alpha:\mathcal H_T(\alpha)\le C_KX\},\qquad
 \mathfrak m=\mathbb T_K\setminus\mathfrak M.
\]
Take a fixed sufficiently small $C_K$ when disjointness is wanted; changing
it by a fixed factor changes no bound below.

\begin{proposition}[Cubic auxiliary major and minor estimates]
\label{prop:g7-cubic}
Under the preceding assumptions,
\begin{align}
 \sup_{\alpha\in\mathfrak m}|F_j(\alpha)|
       &\ll_{K,\epsilon}X^{3/4+\epsilon}HM,
       \label{eq:g7-cubic-minor}\\
 \int_{\mathfrak M}|F_j(\alpha)|^4\,d\alpha
       &\ll_{K,\epsilon}(XHM)^4T^{-1}X^\epsilon.
       \label{eq:g7-cubic-major}
\end{align}
\end{proposition}

\begin{proof}[Minor-arc estimate]
Apply two ordinary differences in $z$ to a cubic summand and use H\"older
for the $h,m$ sums.  Zero differences contribute $O(X^3H)$ before the final
H\"older factor.  Every remaining linear coefficient is a fixed nonzero
multiple of $h_1\cdots h_j u_1u_2$.  For any prescribed nonzero value,
the number of tuples $(\mathbf h,u_1,u_2)$ producing it is $O(X^\epsilon)$,
by ideal divisors followed by the bounded-height unit count.  Such values
lie in a box of volume $O(HX^2)=O(T/X)$.
The linear sums are bounded by positive kernels of width $X^{-1/n}$.
Equation~\eqref{eq:g7-kernel}, with $J=T/X$ and $N=X$, bounds their total by
\[
 X^\epsilon\left(T/X+X+T/\mathcal H_T(\alpha)\right)
       \ll_K X^\epsilon T/X
 \qquad(\alpha\in\mathfrak m).
\]
Consequently the sum over $h$ of the fourth powers of the cubic sums is
$O(X^{3+\epsilon}H)$, uniformly in $m$.  H\"older over the $HM$ parameter
choices now gives $|F_j|^4\ll X^{3+\epsilon}(HM)^4$, as asserted.
\end{proof}

\subsection{Cubic complete sums and the major arcs}

For a good prime ideal $\mathfrak p$, $q=N\mathfrak p$, define
\[
 w_3(\mathfrak p^{3u+1})=q^{-u-1/2},\qquad
 w_3(\mathfrak p^{3u+2})=w_3(\mathfrak p^{3u+3})=q^{-u-1}
 \quad(u\ge0).
\]
At the finitely many bad primes, including those above $6c_k$, take
$w_3(\mathfrak p^e)=q^{-e/3}$.  Extend multiplicatively and put $w_3(1)=1$.
For any cubic polynomial with fixed unit leading coefficient at good primes,
and any primitive additive character modulo an ideal $\mathfrak q$,
its normalized complete sum is
\begin{equation}\label{eq:g7-complete}
 \ll_{K,\epsilon}(N\mathfrak q)^\epsilon w_3(\mathfrak q),
\end{equation}
uniformly in all integral lower coefficients, including a linear twist.
For completeness, at good primes translate to $az^3+bz$.
At modulus $\mathfrak p$ the finite-field Weil bound gives $2q^{-1/2}$.
For higher powers, stationary residue classes satisfy $3az^2+b=0$.
Simple classes give quadratic Gauss bounds $O(q^{-e/2})$.
A multiple class is $z\equiv0$, with $\mathfrak p\mid b$.
If $\mathfrak p\parallel b$, its contribution is zero for $e\ge3$;
if $\mathfrak p^2\mid b$, substitution $z=\pi y$ gives precisely the
factor $q^{-1}$ times the normalized sum at level $e-3$.
The base cases $e=2,3$ give $q^{-1}$.  These recurrences prove the displayed
weight, with an absolute factor per prime absorbed by $(N\mathfrak q)^\epsilon$.
At bad primes one may apply the local-field van der Corput estimate to the
third derivative.  Proposition~3.3 of Cluckers \cite{Cluckers2011}, \emph{Expositiones Mathematicae}
29 (2011), 371--386, gives $O_K(q^{-e/3})$ uniformly here: after division by
$6c_k$ the third derivative has fixed absolute value, and all other coefficients
have bounded local Gauss norm.  The cited local-field estimate is uniform under this normalization.

Two consequences, both elementary Euler-product calculations, are
\begin{align}
 \sum_{N\mathfrak q\le Y}(N\mathfrak q)w_3(\mathfrak q)^4
       &\ll_{K,\epsilon}Y^\epsilon,\label{eq:g7-wsum}\\
 \sum_{\mathfrak d\mid\mathfrak q}
       \frac{w_3(\mathfrak q/\mathfrak d)}{N\mathfrak d}
       &\ll_{K,\epsilon}(N\mathfrak q)^\epsilon w_3(\mathfrak q).
       \label{eq:g7-wconv}
\end{align}
For the first, the good-prime Euler factor is $1+O(q^{-1})$ and the remaining
terms are summable; for the second the local ratios form a geometrically
convergent sum up to a bounded factor per prime.

\begin{proof}[Major-arc estimate]
Write $\alpha=\gamma+\beta$ on an arc with conductor $\mathfrak q$ and
$(N\mathfrak q)A_T(\beta)\le C_KX$.
For fixed $h,m$, put $a=h_1\cdots h_j$ and
\[
 \mathfrak q_a=\mathfrak q/(\mathfrak q+(a)).
\]
This is exactly the conductor of $a\gamma$.
Poisson summation modulo $\mathfrak q_a$, followed by
\eqref{eq:g7-complete}, bounds the cubic sum by
\begin{equation}\label{eq:g7-poisson}
 \ll X^{1+\epsilon}w_3(\mathfrak q_a)
       A_{X^3}(a\beta)^{-1/2}.
\end{equation}
Here are the archimedean details.  At a real place, the second derivative is
bounded away from zero after division by its scale, and the second derivative
test gives $(1+X^{3/n}|a_v\beta_v|)^{-1/2}$.
At a complex place the real Hessian of the real part of the holomorphic cubic
has determinant a nonzero constant times $|a_v\beta_v z_v|^2$;
two-dimensional stationary phase gives the square of this factor.
These estimates are uniform in the linear Poisson and completion frequencies.
Outside the natural gradient box repeated integration by parts gives arbitrary
decay.  Ideal packing bounds the total number of relevant dual frequencies by
\[
 O_K\left(1+\frac{N\mathfrak q_a}{X}A_{X^3}(a\beta)\right)=O_K(1),
\]
because $|a_v|\ll H^{1/n}$ and the product-arc inequality holds.
Dyadic enlargement of the gradient box, with enough integrations by parts,
proves the same bound for the tails.  This proves \eqref{eq:g7-poisson} for
the complete sum, including all Poisson modes.

It remains to average over $h$.  On an archimedean dyadic box with radii
$r_{i,v}$ for $h_i$, let $U_i=\prod_v r_{i,v}^{d_v}$.
There are no nonzero $h_i\in\mathfrak b$ in such a box unless
$U_i\gg_KN\mathfrak b$; when nonempty, ideal packing gives
$O_K(U_i/N\mathfrak b)$ choices.
Split the condition $\mathfrak d\mid(h_1\cdots h_j)$ among divisor ideals
of the individual $h_i$.  The number of splittings is $O(X^\epsilon)$.
Thus summing \eqref{eq:g7-poisson} on this dyadic box, and using
\eqref{eq:g7-wconv}, gives its trivial cardinality times
$X^\epsilon w_3(\mathfrak q)$ times the archimedean decay.
The sum over dyadic radii converges toward their largest values: at each
place the counting factor has degree $d_v$ whereas the decay removes only
$d_v/2$.  Therefore
\[
 |F_j(\gamma+\beta)|
 \ll X^{1+\epsilon}HM\,w_3(\mathfrak q) A_T(\beta)^{-1/2}.
\]
Finally $\int_{V}A_T(\beta)^{-2}\,d\beta\ll_KT^{-1}$,
and there are at most $N\mathfrak q$ center classes of conductor
$\mathfrak q$.  Sum the fourth power of the preceding bound using
\eqref{eq:g7-wsum}.  Overlaps, if any, only enlarge this upper bound.
This proves \eqref{eq:g7-cubic-major}.
\end{proof}

\subsection{The pure smooth terminal process}

Let $\lambda_t$ be any permissible smooth $2t$-moment exponent, extended to
real $t\ge0$ by H\"older interpolation.  For every $u\ge0$, Proposition
\ref{prop:g7-cubic} implies
\begin{equation}\label{eq:g7-M1-terminal}
 \int|F_j(\alpha)g_Q(\alpha)^{2u}|\,d\alpha
 \ll X^{1+\epsilon}HM
 \left\{X^{-1/4}Q^{\lambda_u}
       +Q^{-k/4+(3/4)\lambda_{4u/3}}\right\}.
\end{equation}
The minor term is the supremum estimate times the full smooth moment.
The major term is H\"older with exponents $4$ and $4/3$, using
$T=Q^k$.  This is the pure smooth $M_1$ process; it needs no extra ordinary
variables inside the terminal integral and no smooth major-arc asymptotic.

\section{Auxiliary quadratic terminals}
\label{sec:quadratic-terminals}
\label{sec:eighth-quadratic-terminal}
\label{sec:ninth-terminal}

The residual degrees two, three, and four in efficient differencing
admit a common treatment. Ordinary differences reduce them to a
quadratic sum and separate coefficient sums of degree $2k$.
We combine their moment estimates on one shift space, then integrate
the resulting estimate against a restricted-height smooth mean.

Fix $k\ge4$ and a depth $j$ such that $d=k-j\in\{2,3,4\}$.
Put $r=d-2$ and
\[
 \begin{gathered}
 M_i=X^{\phi_i},\quad H_i=XM_i^{-k},\quad
 \Phi_i=\sum_{a\le i}\phi_a,\quad q_i=1-\Phi_i,
 \qquad0\le\phi_i<1/k,\\
 M=\prod_{i\le j}M_i,\quad H=\prod_{i\le j}H_i,\quad
 Q_0=X/M,\quad T=Q_0^k=X^dH,\\
 A=XHM=X^{a_0},\quad a_0=1+j-(k-1)\Phi_j,\qquad
 H_*=X^rH=X^{h_*},\quad h_*=k-2-k\Phi_j.
 \end{gathered}
\]
All these scales are volumes. The conductor height is normalized
as in Section~\ref{sec:coefficient-sampling-common}, and
$\mathcal H_T(\alpha)\ll_KT^{1/2}$. The parameters are fixed before
the smoothness exponent is chosen; factors $X^{O_K(\eta)}$ are
absorbed in the final $X^\epsilon$.

\subsection{Supports and a common quadratic completion}

Partition the balanced principal divisors into dyadic norm annuli
and finitely many small archimedean boxes. On one piece each
$|m_{i,v}|$ has a prescribed radius up to fixed factors.
While the expressions are nonnegative orthogonality counts,
enlarge the ordinary and shift ranges to independent Cartesian
boxes. The constraints $|h_{i,v}m_{i,v}^k|\ll\delta X^{1/n}$
change by fixed factors only. Taking the original relative
diameter small enough keeps all ordinary vertices in a common
box about a nonzero center at every embedding, through all
subsequent ordinary differences. Centering congruences can be
discarded after their fixed denominators are cleared.
The number of pieces is $X^\epsilon$.

In a fixed integral basis take the ordinary box to be a
parallelepiped. Its intersections with translates after ordinary
differencing depend only on the ordinary shifts; the coefficient
sets and efficient-shift boxes remain independent. These
enlargements were made in counts before a terminal absolute value
was introduced.

For the centered efficient-difference polynomial one has the generating
expansion
\[
 \Psi_j(z;\mathbf h,\mathbf m)
 =2^jk!\Bigl(\prod_{i\le j}h_i\Bigr)
 \left\{\frac{z^d}{d!}
   +\frac{z^{d-2}}{6(d-2)!}\sum_{i\le j}h_i^2m_i^{2k}
   +\text{terms of degree at most }d-4\right\}.
\]
Apply $r=d-2$ ordinary differences with nonzero shifts $u_1,\ldots,u_r$.
With
\[
 \xi=2z+u_1+\cdots+u_r,\qquad
 \Gamma=\Bigl(\prod_{i\le j}h_i\Bigr)\Bigl(\prod_{a\le r}u_a\Bigr),
\]
the resulting polynomial is exactly
\begin{equation}\label{eq:quadratic-terminal-exact}
 \frac{2^jk!}{24}\Gamma
 \left(3\xi^2+\sum_{a\le r}u_a^2
                         +4\sum_{i\le j}h_i^2m_i^{2k}\right).
\end{equation}
This formula also covers $r=0$, in which case $\xi=2z$.
All constants are fixed; denominators can be cleared once.

Let $F_j$ be the original auxiliary sum and $\mathcal G_r$ its sum after
these ordinary differences, including the ordinary shifts. Cauchy's
inequality over the existing parameter tuple, before separating the
coefficient sums, gives
\begin{equation}\label{eq:quadratic-terminal-unseparated}
 |F_j|^{2^r}\ll
       X^{2^r-1}(HM)^{2^r}
       +X^{2^r-r-1}(HM)^{2^r-1}|\mathcal G_r|.
\end{equation}
Only $r=0,1,2$ is needed. For $r=1$, the zero ordinary shift contributes
$X(HM)^2$ and the remaining term is $HM|\mathcal G_1|$.
For $r=2$, Cauchy on the original tuple and the first shift gives
$|\mathcal G_1|^2\ll X^3(HM)^2+XHM|\mathcal G_2|$; substitution gives
\eqref{eq:quadratic-terminal-unseparated}. For $r=0$, it is the harmless bound
$|F_j|\le HM+|\mathcal G_0|$.

For each ordinary shift tuple, the remaining $z$-domain is an intersection
of translates of the same parallelepiped. Fourier completion in the fixed
integral basis gives an integrable kernel dominated by a common majorant
of $L^1$ norm $O((\log X)^n)$. The kernel may depend on $\mathbf u$, but
is independent of $\mathbf h,\mathbf m$. The condition
$\xi\equiv\sum_a u_a\pmod{2\mathcal O_K}$ is imposed by finitely many
additive characters, again independent of $\mathbf m$. Thus after one
common integration in a linear modulation $t$, the absolute value of
$\mathcal G_r$ is bounded by $X^\epsilon$ times
\begin{equation}\label{eq:quadratic-terminal-completion}
 \sup_t\sum_{\mathbf h,\mathbf u}
       |Z_t(c_0\alpha\Gamma)|
       \prod_{i=1}^j|L_i(c_i\alpha\Gamma h_i^2)|,
\end{equation}
uniformly in $t$, where
\[
 Z_t(\beta)=\sum_{z\in\mathcal B(X)}e_K(\beta z^2+tz),\qquad
 L_i(\beta)=\sum_{m\in\mathcal C_i}e_K(\beta m^{2k}),
\]
and $c_0,c_i$ are fixed nonzero integers. Smooth completion gives the same
formula with bounded weights and the same coefficient tuple throughout.

\subsection{Moment resources on the shift space}

On a minor height shell write
\[
 \mathcal H_T(\alpha)\asymp X^\ell,
 \qquad 1\le\ell\le\frac{kq_j}{2}.
\]
The fixed-factor endpoint changes and smoothness perturbations are
$O(1/\log X)+O_K(\eta)$ and are absorbed in the final error.
All products of a fixed number of nonzero algebraic integers have only
$X^\epsilon$ representations in the relevant boxes, by ideal divisors
and bounded-height unit counting. Thus $\Gamma$ may be grouped as a single
variable in $\mathcal B(H_*)$ for nonnegative moment estimates.

Majorize the shift box by a positive Gaussian and expand
the quadratic second moment. The diagonal $z^2=w^2$ contributes
$O(XH_*)$, and each nonzero value of $(z-w)(z+w)$ has
$X^\epsilon$ preimages. After the Gaussian sum in the shift
variable is evaluated as a nonnegative periodized kernel, the
common modulation can be removed in absolute value.
The linear-kernel estimate gives
\[
 \sum_a|Z_t(c\alpha a)|^2
 \ll X^\epsilon\left(XH_*+X^2+H_*+
                \frac{H_*X^2}{\mathcal H_{H_*X^2}(\alpha)}\right).
\]
Here $H_*X^2=T$. Dividing by the trivial size $H_*X^2$
gives the saving $\min\{1,h_*,\ell\}$, hence on the minor shells
\begin{equation}\label{eq:quadratic-terminal-z2}
 z_2=\min\{1,h_*\}.
\end{equation}
In particular $z_2=1$ when $H_*\ge X$.
The sampled fourth moment of Section~\ref{sec:quad-sampling-field}
gives the additional resources
\[
 z_4=\min\{2,\ell,h_*\},\qquad
 z_4^{\rm ref}=\min\{2,2\ell-h_*\}.
\]
Each is used only on intervals where its hypotheses hold.
The trivial bound is also available, with saving zero.

For coefficient block $i$, group all shifts except $h_i$ into a variable
of volume
\[
 J_i=X^r\prod_{a\ne i}H_a,
 \qquad \log_X J_i=k-3-k\Phi_j+k\phi_i.
\]
The conductor identity is
$J_iH_i^3M_i^{2k}=T$.
For an even coefficient order $p$, raise $L_i$ to the integer power $p/2$.
Its Fourier coefficients have support in a box of volume
$N_i\ll X^\epsilon M_i^{2k}$, and their squared norm is at most
\[
 X^\epsilon M_i^{p-2k+\Delta_p^{(2k)}}.
\]
This follows by orthogonality and coefficient domination by the full
$2k$th-power smooth moment. The exponent $\Delta_p^{(2k)}$ is obtained
from the independently proved field seed curve. The refined coefficient
sampling estimate, with
$Y\le c_K\min\{N_i,J_i,(J_iN_iH_i^{-3})^{1/2}\}$, gives saving
\begin{align}
 b_{i,p}^{\rm ref}(\ell)
 &=\min\left\{2k\phi_i,\ k-3-k\Phi_j+k\phi_i,
       \frac{k-6-k\Phi_j+6k\phi_i}{2},\ \ell-2+2k\phi_i\right\}
       -\Delta_p^{(2k)}\phi_i.\label{eq:quadratic-terminal-refined}
\end{align}
The cubic-height sampling estimate \eqref{eq:k5-sampling-strong} gives the alternative
\begin{align}
 b_{i,p}^{\rm cube}(\ell)
 &=\min\left\{2k\phi_i,\ k-3-k\Phi_j+k\phi_i,
       \frac{k-6-k\Phi_j+6k\phi_i}{2},\ \ell/3\right\}
       -\Delta_p^{(2k)}\phi_i.\label{eq:quadratic-terminal-cube}
\end{align}
The coefficient normalization is exact: sampling produces
$J_iH_i M_i^{p+\Delta_p^{(2k)}}$ times the reciprocal saving, whereas the
trivial normalized size is $J_iH_iM_i^p$.

Proposition~\ref{prop:k5-coeff2}, with $d=k$, gives the basic second coefficient resource
\begin{equation}\label{eq:quadratic-terminal-basic}
 b_{i,2}^{\rm basic}=
 \min\left\{\phi_i,\ k-6-k\Phi_j+4k\phi_i,
               q_j/4,\ \frac{kq_j-(4k-1)\phi_i}{6}\right\}.
\end{equation}
Put $b=b_{i,2}^{\rm basic}$. On the present minor shells, the lemma's other term
$\mathcal H_T(\alpha)^{-1/k}$ is no larger than $X^{-b}$, since
$b\le\phi_i<1/k$ and $\ell\ge1$.
A block with $\phi_i=0$ may be bounded trivially, including its fixed smoothness factors.

Each displayed coefficient resource bounds the sum of $|L_i|^p$, divided
by the common shift volume $H_*=X^rH$, over the \emph{same} shift tuple in \eqref{eq:quadratic-terminal-completion}
by $X^\epsilon M_i^pX^{-b_{i,p}(\ell)}$.
The quadratic resources likewise bound these normalized sums of $|Z_t|^p$ by
$X^\epsilon X^pX^{-z_p(\ell)}$.
Grouping products only enlarges nonnegative sums, so no independence between
these factors is needed.

\subsection{A combined minor-arc terminal}

Choose nonnegative resource weights $w_R$ on a fixed height interval.
If a resource has order $p_R$ and belongs to factor $i$, impose
\[
 \sum_{R\text{ on factor }i}p_Rw_R\le1
       \quad(i=0,1,\ldots,j),\qquad \sum_R w_R\le1.
\]
Factor $0$ denotes the quadratic sum. Divide each factor by its trivial
size, $X$ for $Z_t$ and $M_i$ for $L_i$. Its unused power is bounded
trivially, and the unused H\"older mass is assigned to the constant
function $1$. Generalized H\"older on the common shift space then gives
\[
 |\mathcal G_r|\ll X^{r+1+\epsilon}HM\,X^{-S(\ell)},\qquad
 S(\ell)=\sum_Rw_Rb_R(\ell).
\]
Combining this with \eqref{eq:quadratic-terminal-unseparated} proves
\begin{equation}\label{eq:quadratic-terminal-pointwise}
 |F_j(\alpha)|\ll A X^{\epsilon-
               2^{-(k-j-2)}\min\{1,S(\ell)\}}.
\end{equation}
In particular, the Weyl root is applied to the entire combined saving
$S(\ell)$, after all resources have been coupled by one H\"older inequality.

Let the remaining smooth order be $2u$. Choose nonnegative
$a,c,\gamma$ and an even order $v$ such that
\[
 a+c=1,\qquad \gamma+2a=1,\qquad cv=2u.
\]
If $a>0$, assume the auxiliary second-moment estimate
$\int|F_j|^2\ll X^\epsilon A$, supplied by
Lemma~\ref{lem:field-general-seed-aux} when its hypotheses hold.
The usual choice is $a=0$, $c=\gamma=1$, and $v=2u$.
On a height shell, H\"older gives
\[
 \int_{\rm shell}|F_j|\,|g_{Q_0}|^{2u}
 \le(\sup_{\rm shell}|F_j|)^\gamma
       \left(\int|F_j|^2\right)^a
       \left(\int_{\rm shell}|g_{Q_0}|^v\right)^c.
\]
The restricted-height smooth theorem and the pointwise estimate
therefore give the exponent
\begin{equation}\label{eq:quadratic-terminal-shell}
 C(\ell)=(\gamma+a)a_0+c(v-k)q_j+
       \frac{2c\Delta_v^{(k)}}k\ell
       -\gamma2^{-r}\min\{1,S(\ell)\}.
\end{equation}
Each resource is a minimum of affine functions of $\ell$,
hence concave. Its nonnegative weighted sum and the minimum
of that sum with one are also concave. Thus $C(\ell)$ is convex
on every interval with fixed resource weights, and its maximum
is attained at an endpoint.

\subsection{Uniform auxiliary major arcs}
\label{sec:aux-major-24}

For each auxiliary tuple write $a=h_1\cdots h_j\ne0$ and
\[
 \Psi_j(z;\mathbf h;\mathbf m)=aP_{\mathbf h,\mathbf m}(z),
 \qquad P_{\mathbf h,\mathbf m}(z)=c_{k,j}z^d+\cdots,
 \qquad c_{k,j}\ne0.
\]
The leading coefficient is a fixed integer, and all other coefficients
are algebraic integers after clearing one fixed denominator if necessary.
The ordinary variable has a nonnegative smooth weight in the small box
retained by efficient differencing.  All scaled weight derivatives are
uniformly bounded.  Linear completion twists may be arbitrary, independently
for every auxiliary tuple.

\begin{proposition}\label{prop:aux-major-24}
Let $F_j$ be the resulting auxiliary sum.  On every broad arc
\[
 \alpha=\gamma+\beta,\qquad Q_\gamma A_T(\beta)\le C_KX,
\]
one has, for every $\epsilon>0$,
\begin{equation}\label{eq:aux-major-strong}
 |F_j(\gamma+\beta)|
 \ll_{K,\epsilon}X^{1+\epsilon}HM\,
            Q_\gamma^{-1/d}A_T(\beta)^{-1/2}.
\end{equation}
Consequently, with $L=Q_\gamma A_T(\beta)$,
\begin{equation}\label{eq:aux-major-height}
 |F_j(\gamma+\beta)|\ll X^{1+\epsilon}HM L^{-1/d}.
\end{equation}
The bound includes every Poisson mode, also when $d=2$.
\end{proposition}

\subsubsection{Uniform mixed complete sums}

The required complete-sum bound is slightly more general than the ordinary
monomial-plus-linear estimate.  If $P\in\mathcal O_K[z]$ has degree
$d\le4$ and fixed nonzero leading coefficient $c$, then for a primitive
character $\psi$ modulo $\mathfrak q$ and every integral $b$,
\begin{equation}\label{eq:aux-complete-general}
 (N\mathfrak q)^{-1}
 \left|\sum_{z\bmod\mathfrak q}\psi(P(z)+bz)\right|
 \ll_{K,c,d,\epsilon}(N\mathfrak q)^{-1/d+\epsilon},
\end{equation}
uniformly in the lower coefficients.  Those coefficients are allowed to
vary with $\mathbf h,\mathbf m$ and with the Poisson frequency.

Here is a local proof recording uniformity in the prime.  At a prime of
residue characteristic greater than $d$ and not dividing $c$, consider
an integral polynomial whose reduction has positive degree $r\le d$.
The normalized complete sum of conductor $\mathfrak p^E$ is
$O_d(q^{-E/r})$, where $q=N\mathfrak p$.
For $E=1$ this is the finite-field Weil bound (the linear case is zero).
For $E\ge2$, nonstationary residue classes contribute zero.
At a stationary residue class of multiplicity $m$ as a root of the reduced
derivative, translate to that class and substitute $z=z_0+\pi y$.
Let $s$ be the minimum valuation of the resulting nonconstant Taylor
coefficients.  Then $2\le s\le m+1\le r$; after dividing by $\pi^s$,
the new reduction degree $r'$ satisfies $1\le r'\le s$.
The class has measure $q^{-1}$, and induction supplies
\[
 q^{-1}O_d(q^{-(E-s)/r'})=O_d(q^{-E/r})
 \qquad(E>s).
\]
When $E\le s$, its trivial measure $q^{-1}$ is already at most $q^{-E/r}$.
A child with $r'=r$ can occur only when the reduced derivative has a single
root of multiplicity $r-1$.  Thus a same-degree recursion has only one
child and does not accumulate constants.  Every other branch lowers the
reduction degree; at most $r-1$ stationary residue classes occur.
Induction on reduction degree and conductor gives a constant depending
only on $d$ (one may take constants bounded by $(d-1)!$ for $d\le4$).
This proof applies over every finite residue extension and is unchanged
by ramification at such good primes.

At the finitely many remaining primes, the $d$th derivative is the fixed
nonzero constant $d!c$.  Cluckers' local-field van der Corput theorem
\cite[Proposition~3.3]{Cluckers2011} applies after putting
$f=(P+bz)/(d!c)$: then $f^{(d)}=1$, while the Gauss norm of
$f-f(0)$ is at most $|d!c|_\mathfrak p^{-1}$ uniformly in the integral
lower coefficients and in $b$.  Its permitted dependence on the fixed
local field, derivative order and Gauss norm consequently gives
$O_{K_\mathfrak p,c,d}(q^{-E/d})$ uniformly in these coefficients.
The fixed denominators from centering are absorbed into fixed fractional
lattices and fixed conductor factors.  Chinese remaindering and the
ideal-divisor bound absorb
products of the good-prime constants into $(N\mathfrak q)^\epsilon$.
This proves \eqref{eq:aux-complete-general}.

\subsubsection{Archimedean curvature and all Poisson modes}

For a fixed auxiliary tuple, put $P=X^{1/n}$ and write
$p_v(u)=P^{-d}P_{\mathbf h,\mathbf m,v}(Pu)$.
On the retained tiny support all increments $h_im_i^k$ are small relative
to the nonzero ordinary center.  The centered difference formula therefore
gives
\[
 p_v(u)=c_{k,j}u^d+O(\delta^2)
\]
in bounded $C^d$ norms, with $|u|\asymp_K1$.  Taking $\delta$ sufficiently
small ensures $|p_v''(u)|\gg_K1$ uniformly.  For $d=2$ this is immediate
from its fixed leading coefficient.
At a real place, the second-derivative test is uniform in any added linear
phase.  At a complex place, the Hessian determinant of
$\operatorname{Re}(t p_v(u))$ is a nonzero constant times
$-|t p_v''(u)|^2$, so two-dimensional stationary phase gives the corresponding
square decay.  After a fixed partition of the tiny support, the standard
uniform stationary-phase theorem \cite[Section~7.7]{Hormander2003} and integration
by parts give
\begin{equation}\label{eq:aux-arch-tails}
 \left|\int w_v(u)e_{K_v}(t p_v(u)-\zeta u)\,du\right|
 \ll_{K,B}(1+|t|)^{-d_v/2}
       \left(1+\frac{|\zeta|}{1+|t|}\right)^{-B}
\end{equation}
for every prescribed $B$.  A linear completion twist simply translates
$\zeta$; the constants are independent of its size.  The same conclusion
holds for products of translated smooth vertex weights, since their
scaled derivatives remain uniformly bounded.  Sharp completion, if used,
only contributes an arbitrary $X^\epsilon$ factor.

Let $\mathfrak q$ be the conductor of $\gamma$ and
\[
 \mathfrak q_a=\mathfrak q/(\mathfrak q+(a)),\qquad q_a=N\mathfrak q_a.
\]
This is the conductor of $a\gamma$.
Poisson summation of the fixed-tuple sum modulo $\mathfrak q_a$, followed
by \eqref{eq:aux-complete-general}, bounds it by
\[
 Xq_a^{-1/d+\epsilon}A_{X^d}(a\beta)^{-1/2}
 \sum_{\xi\in\mathfrak q_a^{-1}\mathfrak D^{-1}}
 \prod_v\left(1+
   \frac{P|\xi_v-\xi_{0,v}|}{1+P^d|a_v\beta_v|}\right)^{-B},
\]
where $\xi_0$ is the arbitrary shift from completion.
The normalizing factor $1/q_a$ from Poisson has already been combined
with the complete-sum bound.

Ideal packing in a translated product box bounds the number of relevant
dual points by
\[
 O_K\left(1+\frac{q_a}{X}A_{X^d}(a\beta)\right).
\]
Indeed the radii are $(1+P^d|a_v\beta_v|)/P$ and their product volume is
$A_{X^d}(a\beta)/X$.  Since $|a_v|\ll H^{1/n}$,
\[
 q_aA_{X^d}(a\beta)\ll_K Q_\gamma A_T(\beta)\ll_KX.
\]
For completeness this packing also sums all tails.  Enlarge each radius
by an independent factor $2^{r_v}$; the number of dual points grows by at
most $\prod_v2^{d_vr_v}$, while \eqref{eq:aux-arch-tails} supplies
$\prod_v2^{-Br_v}$.  Choosing $B>n+2$ makes the resulting series convergent,
uniformly in $\xi_0$ and in the eccentricity of the gradient box.
Thus the whole dual sum is $O_K(1)$, proving the individual bound
\begin{equation}\label{eq:aux-fixed-tuple}
 |S_{\mathbf h,\mathbf m}(\gamma+\beta)|
 \ll X^{1+\epsilon}q_a^{-1/d}A_{X^d}(a\beta)^{-1/2}.
\end{equation}

\subsubsection{Averaging the conductor greatest common divisors}

There are $O(MX^\epsilon)$ choices of $\mathbf m$.
For prescribed $a\ne0$, the number of bounded-height tuples
$h_1\cdots h_j=a$ is $O(X^\epsilon)$, by ideal divisors and the unit
lattice.  We may therefore sum over $a\in\mathcal B(H)\setminus\{0\}$.
Write $Q=Q_\gamma$ and $A_\beta=A_T(\beta)$. The elementary inequality
\[
 A_{X^d}(a\beta)\gg_K A_\beta|Na|/H
\]
gives
\[
 A_{X^d}(a\beta)^{-1/2}
 \ll_K\min\{1,(H/(A_\beta|Na|))^{1/2}\}.
\]
For an ideal $\mathfrak b$ of norm $b$, the small-norm ideal count
\[
 \#\{a\in\mathfrak b\cap\mathcal B(H):0<|Na|\le Y\}
       \ll H^\epsilon Y/b
\]
has no points for $Y<b$.  Dyadic summation consequently gives, including
thin scales,
\[
 \sum_{0\ne a\in\mathfrak b\cap\mathcal B(H)}
       \min\{1,(H/(A_\beta|Na|))^{1/2}\}
 \ll H^\epsilon\frac{H}{bA_\beta^{1/2}}.
\]
Split according to $\mathfrak b=\mathfrak q+(a)$, and enlarge each class
to all nonzero $a\in\mathfrak b$.  Equation~\eqref{eq:aux-fixed-tuple} then
yields
\begin{align*}
 |F_j(\gamma+\beta)|
 &\ll X^{1+\epsilon}M Q^{-1/d}
       \sum_{\mathfrak b\mid\mathfrak q}(N\mathfrak b)^{1/d}
         \frac{H}{(N\mathfrak b)A_\beta^{1/2}}\\
 &\ll X^{1+\epsilon}HM Q^{-1/d}A_\beta^{-1/2},
\end{align*}
since $\sum_{\mathfrak b\mid\mathfrak q}(N\mathfrak b)^{-1+1/d}
\ll Q^\epsilon$.  This proves \eqref{eq:aux-major-strong}.
Because $A_\beta\ge1$ and $d\ge2$, it implies \eqref{eq:aux-major-height}.

\subsubsection{The major terminal}

Let $v$ be an even smooth order with permissible excess $\Delta_v^{(k)}$.
The restricted-height theorem bounds its mean on the union of arcs
of height at most $L$ by
\[
 Q_0^{v-k+\epsilon}L^{2\Delta_v^{(k)}/k}.
\]
Summing dyadic shells $1\le L\le C_KX$ with
\eqref{eq:aux-major-height} gives, provided $X\le Q_0^{k/2}$,
\begin{equation}\label{eq:quadratic-terminal-major-integral}
 \int_{\mathfrak M}|F_jg_{Q_0}^v|
 \ll X^\epsilon A Q_0^{v-k}
       X^{\max\{0,2\Delta_v^{(k)}/k-1/d\}}.
\end{equation}
Here $h_*=kq_j-2>0$ implies the required scale inequality.
An equality in the exponent contributes only a logarithm.

For the counted recurrence use the original smooth order $2u$,
without the optional auxiliary redistribution. Its major terminal is
\begin{equation}\label{eq:quadratic-terminal-major}
 C_{\rm maj}=a_0+(2u-k)q_j+
              \max\{0,2\Delta_{2u}^{(k)}/k-1/(k-j)\}.
\end{equation}

\subsection{The counted recurrence}

Let $\lambda_t$ denote the current full exponent at smooth
moment order $2t$. Choose integer half-orders
$r_1=s-1,\ldots,r_j=u$ with $0<r_{i+1}<2r_i$.
Start $C$ with the maximum of
\eqref{eq:quadratic-terminal-major} and every endpoint value of
\eqref{eq:quadratic-terminal-shell} in the assigned intervals.
The principal-block counted recurrence gives, for $i=j-1,\ldots,1$,
\[
 \begin{split}
 D_i&=1+i-k\Phi_i+\Phi_{i+1}+q_{i+1}\lambda_{r_{i+1}},\\
 B_i&=\tfrac12\{q_i\lambda_{2r_i-r_{i+1}}
                   +i-(k-1)\Phi_i+(2r_{i+1}-1)\phi_{i+1}\},\\
 C&\longleftarrow B_i+\tfrac12\max(D_i,C).
 \end{split}
\]
The first-stage exponent is then
\[
 \lambda_s^{\rm new}=(2s-3)\phi_1+
           \max\{1+\phi_1+q_1\lambda_{s-1},C\}.
\]
Every subsisting half-order is integral, every backward stage satisfies
$0<2r_i-r_{i+1}<2r_i$, and each round reads only the exponents available
at its start. More precisely, fix a sufficiently small ordinary root box
$\mathcal B$ away from zero at every archimedean place, with a contained
ideal-smooth root set $\mathcal A$. Write $f_{\mathcal B}$ and
$g_{\mathcal A}$ for their unweighted exponential sums. The first-stage
formula bounds the localized joint moment
$\int|f_{\mathcal B}|^2|g_{\mathcal A}|^{2s-2}$ itself, with the support
enlargements described at the start of this section. Coefficient domination
in its nonnegative orthogonality count gives
\[
 \int|g_{\mathcal A}|^{2s}
 \le \int|f_{\mathcal B}|^2|g_{\mathcal A}|^{2s-2}.
\]
A fixed compact root annulus can be covered by finitely many such boxes;
$|\sum_{a=1}^N g_a|^{2s}\le N^{2s-1}\sum_{a=1}^N|g_a|^{2s}$ then gives
the required pure smooth moment on that annulus. Bounded nonnegative
weights are handled by the same coefficient domination, with constants
depending on their fixed bounds. The formulas give the field estimates evaluated by the exact parameter certificates.

\section{Fourth powers}
\label{sec:fourth-powers}\label{sec:q12-analytic}
The quartic argument needs two seed moments and three applications
of the cubic terminal. Write $b_s$ for the exponent of the smooth
moment of order $2s$.

\subsection{The two seed moments}
At depth one, the A processes of orders two and four give
\begin{align}
 T_3(b_4)&=\max\{3+2\phi_3,
      1+\tfrac32\phi_3+\tfrac12(1-\phi_3)b_4\},
      \label{eq:q12-T3}\\
 T_4(b_3,b_4)&=\max\{1+b_3+(6-b_3)\phi_4,
      \tfrac54+\tfrac{11}{4}\phi_4+
                        \tfrac34(1-\phi_4)b_4\}.
      \label{eq:q12-T4}
\end{align}
These are \eqref{eq:seed-A-terminal} and
\eqref{eq:seed-initial-output} with $k=4$. The auxiliary moments
follow from Lemma~\ref{lem:field-general-seed-aux}.
Choose
\[
 \phi_3=\frac{109}{1000},\qquad\phi_4=\frac7{50},\qquad
 B_3=\frac{161}{50},\qquad B_4=\frac{461}{100}.
\]
The branches have the exact values
\[
\begin{array}{c|cc|c}
 &\text{diagonal}&\text{auxiliary}&\text{target}\\\hline
 s=3&3.218&3.217255&3.22\\
 s=4&4.6092&4.60845&4.61
\end{array}
\]
The coefficient sums are at most $43/50$.
Lemma~\ref{lem:supersolution} therefore proves the permissible
exponents $b_3=161/50$ and $b_4=461/100$.

\subsection{The tenth smooth moment}
\begin{proposition}\label{prop:quartic-tenth}
For every fixed number field $K$ and $\epsilon>0$, the smooth
quartic sum satisfies
\begin{equation}\label{eq:q12-moment}
 \int_{\mathbb T_K}|g_4|^{10}\ll X^{9993/1600+\epsilon}.
\end{equation}
\end{proposition}
\begin{proof}
Use the cubic terminal of Section~\ref{sec:cubic-terminal} with
$k=4$, $j=1$, $s=5$, and $\phi_1=3/16$. Then
\[
 H=X^{1/4},\quad M=X^{3/16},\quad Q=X^{13/16},
 \qquad T=Q^4=X^{13/4}>X^2.
\]
The diagonal and cubic minor branches agree and give
\[
 \frac52+\frac{13}{16}b_4=\frac{9993}{1600}.
\]
In the cubic major branch interpolate
$b_{16/3}\le(2b_5+b_6)/3$ and use the Hua value $b_6=17/2$.
The resulting transformation of $c=b_5$ is
\[
 c\longmapsto\max\left\{\frac{9993}{1600},
       \frac{31}{16}+\frac{13}{32}c+
                           \frac{13}{64}\frac{17}{2}\right\}.
\]
Starting with the Hua value $27/4$, three steps give
\[
 \frac{27}{4}\longmapsto\frac{205}{32}
 \longmapsto\frac{6417}{1024}\longmapsto\frac{9993}{1600}.
\]
The exact calculation is reproduced in
Appendix~\ref{app:quartic-exact-certificate}.
\end{proof}

\subsection{The global quartic bounds}
\begin{proof}[Proof of Theorem~\ref{thm:quartic-refinements}]
The excess at order ten is $\Delta_{10}=393/1600$.
With $(w,S,u)=(20,12,10)$ the absorption margin is
\[
 D=2-8\Delta_{10}=\frac7{200}>0.
\]
Thus $U_{12}\ll X^8$, and the minor-arc saving in
Corollary~\ref{cor:even-transfer} is $D/20=7/4000$.
That corollary, with the local requirement retained, gives
\[
 G_{\mathbb K}(4)\le\max\{12,H_{\mathbb K}(4)\}.
\]
Theorem~\ref{thm:local-table} gives the universal upper bound
sixteen, the upper bound fifteen when no dyadic completion is
$\mathbb Q_2$, and the upper bound twelve when every dyadic residue
degree is at least two.

For $K\supset\mathbb Q(\sqrt5)$, the polynomial $t^2-t-1$ is
irreducible over $\mathbb F_2$. Every dyadic completion of $K$
therefore contains the unramified quadratic extension of
$\mathbb Q_2$, and its residue degree is even. Hence
$G_{\mathbb K}(4)\le12$.
\end{proof}

\section{Fifth powers}
\label{sec:k5-terminal}

Write $M_i=X^{\phi_i}$, $H_i=XM_i^{-5}$, $H=H_1H_2$,
$M=M_1M_2$, $Q=X/M$, and $T=Q^5=X^3H$.
The centered second difference for fifth powers is a cubic in $z$.
We adapt the sharper dissection of \cite[Lemma~6.1]{VW1995}.
After one ordinary difference in $z$, its quadratic part and its two
coefficient parts separate as
\[
 K(\alpha;h,h_1,h_2)
 L_1(\alpha;h,h_1,h_2)L_2(\alpha;h,h_1,h_2),
\]
where, up to fixed nonzero constants in $K$,
\[
 K=\sum_z e_K(\alpha hh_1h_2 z^2+\xi z),\qquad
 L_i=\sum_{m_i}e_K(c\alpha hh_1h_2h_i^2m_i^{10}+b_i(m_i)).
\]
The ordinary variable $h$ has volume $X$.

The support issue deserves to be explicit.  Partition each extracted
principal divisor into a dyadic norm annulus and then into finitely many
tiny archimedean boxes.  Its generator is balanced, so on one piece all
$|m_{i,v}|$ are comparable to $M_i^{1/n}$, up to the harmless original
smoothness loss.  The actual constraints
$|h_{i,v}m_{i,v}^5|\ll\delta X^{1/n}$ are then enlarged to independent
small boxes for $h_i$.  Every enlarged vertex remains in a slightly larger
tiny box about the same nonzero center.  This enlargement is made in the
nonnegative orthogonality count before an absolute value is introduced.
Consequently it preserves the nonvanishing of the differenced polynomial
and its derivative and yields the Cartesian auxiliary sums used here.
Smooth completion can alternatively be applied to the original vertex
weights.  Fourier inversion gives a linear twist in $z$ and twists of
degree five in $m_i$; the coefficient of the latter depends on $h_i$ and the
completion frequencies but is independent of the linear sampling variable
$hh_{3-i}$.  The estimates below are uniform in precisely these twists.

\subsection{The higher coefficient moment}
The coefficient-sampling proof in Section~\ref{sec:coefficient-sampling-common} has the
following strengthening.  If $V=H_0^3JN$ and
\[
 Y\le c_K\min\{J,N,(JN/H_0^3)^{1/2}\},
\]
then
\begin{equation}\label{eq:k5-sampling-strong}
 \sum_{0\ne h\in\mathcal B(H_0)}\sum_{j\in\mathcal B(J)}
 \left|\sum_a c_h(a)e_K(\alpha h^3ja)\right|^2
 \ll V^\epsilon JH_0N\{Y^{-1}+\mathcal H_V(\alpha)^{-1/3}\}
                      \sup_h\sum_a|c_h(a)|^2.
\end{equation}
The coefficients may depend on $h$ because the elementary bilinear
large-sieve bound is applied separately for each $h$.
Here are the additional details proving the exponent $-1/3$.
The small-height centers are identified as $h^3\gamma$, where
$L=Q_\gamma A_V(\alpha-\gamma)\ll V^{1/2}$, by the conductor separation
lemma.  Put $\mathfrak d=\mathfrak q+(h)^3$ and
$\mathfrak b=\prod\mathfrak p^{\lceil v_\mathfrak p(\mathfrak d)/3\rceil}$.
For $A=A_V(\alpha-\gamma)$ the norm-count argument bounds the contribution
of this divisor class by
\[
 V^\epsilon\frac{N\mathfrak d}{Q_\gamma}
       \frac{H_0}{A^{1/3}N\mathfrak b}
 \le V^\epsilon H_0L^{-1/3},
\]
since $N\mathfrak b\ge(N\mathfrak d)^{1/3}$ and
$N\mathfrak d\le Q_\gamma$.
There are only $V^\epsilon$ divisor classes.  The large-height classes and
the $J+N$ terms of the bilinear bound contribute $JH_0N/Y$.
This proves \eqref{eq:k5-sampling-strong}.

Apply this with $c_h(a)$ the coefficients of a $t$th power of a decic smooth
sum of volume $M_1$.  Its squared coefficient norm is bounded by the decic
$2t$th moment, even in the presence of the lower-degree twists: by
orthogonality and the triangle inequality, its absolute coefficient count
is dominated by the untwisted decic count.  Thus if $\mu_t$ is a permissible
decic $2t$-moment exponent, then
\[
 E_t(\alpha;J,H_0,M_1)
 \ll X^\epsilon JH_0M_1^{\mu_t+10}
             \{Y^{-1}+\mathcal H_{JH_0^3M_1^{10}}(\alpha)^{-1/3}\}.
\]

\subsection{The terminal estimate}
Let $\theta=\phi_1$, $\phi=\phi_2$, and write $U=X^u$.  Assume
\begin{gather}
 1/15\le\theta\le1/5,\quad 0<\phi\le\theta,\quad
 \phi<5\theta-1/3,\label{eq:k5-parameter1}\\
 0<u<\phi,\quad u<-1-5\theta+15\phi,\quad
 4u+\theta+\phi<1,\quad 6u+5\theta+24\phi<5.
 \label{eq:k5-parameter2}
\end{gather}
The strict inequalities absorb the small powers arising from principal
blocks and from dyadic completion.

Let $t\ge3$, and let $\mu_t$ be the decic $2t$-moment exponent.  Put
\[
 Z=XU^{1-1/t}\left(X^{1/3}M_1^{2t-10-\mu_t}\right)^{1/t},
 \qquad \rho=\frac{\log_X Z}{4}.
\]
On the cubic minor arcs $\mathcal H_T(\alpha)>C_KX$, ordinary differencing
and H\"older give
\[
 |F_2|^2\ll X(HM)^2+HM\,D^{1/2}E_1^{1/(2t)}E_2^{1/2-1/(2t)}.
\]
Here $D$ is the sum of $|K|^2$, $E_1$ the sum of $|L_1|^{2t}$, and
$E_2$ the sum of $|L_2|^{2t/(t-1)}$, over $h,h_1,h_2$.
The quadratic linear-kernel calculation gives
\[
 D\ll X^{2+\epsilon}H.
\]
In detail, after expanding $|K|^2$, the zero difference contributes $X^2H$;
for the other terms factor $z_1^2-z_2^2=(z_1-z_2)(z_1+z_2)$, group the
nonzero coefficient $hh_1h_2(z_1-z_2)$ by ideal divisors, and apply the
positive linear-kernel estimate at scale $T$.  Its minor-arc bound is
$T/X=X^2H$.

Group the product $hh_2$ in $E_1$, at a cost $X^\epsilon$.
Equation~\eqref{eq:k5-sampling-strong}, with
$J=XH_2$, $H_0=H_1$, $N=M_1^{10}$ and $Y=X^{1/3}$, is applicable by
\eqref{eq:k5-parameter1}.  It gives
\[
 E_1\ll X^{2/3+\epsilon}H M_1^{\mu_t+10}.
\]
For $E_2$, first remove the power exceeding two using $|L_2|\ll M_2$.
Proposition~\ref{prop:k5-coeff2}, with $J=XH_1$, $H_0=H_2$, $M=M_2$ and
$d=5$, applies by \eqref{eq:k5-parameter2}.  Since $U<X^{1/5}$, it yields
\[
 E_2\ll X^{1+\epsilon}H M_2^{2t/(t-1)}U^{-1}.
\]
Consequently, provided $\rho\le1/2$,
\begin{equation}\label{eq:k5-sharp-minor}
 \sup_{\mathfrak m}|F_2|\ll X^{1+\epsilon}HM Z^{-1/4}.
\end{equation}

The already-proved cubic complete-sum and Poisson estimates give
\[
 \int_{\mathfrak M}|F_2|^4\ll X^\epsilon(XHM)^4Q^{-5}.
\]
They apply uniformly to the completed support and linear twists just
specified, including all complex places and all nonzero Poisson modes.
For any permissible quintic smooth moments $\lambda_r$ and
$\lambda_{4r/3}$, H\"older on the major arcs and
\eqref{eq:k5-sharp-minor} on the minor arcs therefore prove
\begin{equation}\label{eq:k5-enhanced-terminal}
 \int|F_2g_Q^{2r}|
 \ll X^{1+\epsilon}HM
 \left\{Z^{-1/4}Q^{\lambda_r}
       +Q^{-5/4+(3/4)\lambda_{4r/3}}\right\}.
\end{equation}
This proves the sharper fifth-power terminal, using the degree-ten
moment and the coefficient estimates of Section~\ref{sec:coefficient-sampling-common}.

\subsection{The joint-moment certificate}
\label{sec:k5-certificate}

The principal-block inequality gives both smooth moments and joint
moments with two distinguished ordinary roots. The auxiliary estimates
at depths one and two follow from Lemma~\ref{lem:field-general-seed-aux}.
Start with $\lambda_0=0$, $\lambda_1=1$, $\lambda_2=2$, and
$\lambda_s=2s$ for $3\le s\le8$.

\subsubsection{The map being certified}
The A rows use \eqref{eq:seed-A-terminal} and the common
backwards recurrence. For the enhanced rows, the output is a joint
moment with two ordinary roots.

For the two enhanced rows take $r=6$, $t=2(S-1)-r$, and
\[
 \theta=\phi_1,\quad\phi=\phi_2,\quad q=1-\theta-\phi,
 \quad
 \rho=\frac14\left\{1+\frac{21u}{22}
                  +\frac{1/3+(34-34.227)\theta}{22}\right\}.
\]
The degree-ten exponent $34.227$ is supplied by Corollary~\ref{prop:tenth44}.
Put
\begin{align*}
 B&=\tfrac12\{(1-\theta)\lambda_t+1-4\theta+(2r-1)\phi\},\\
 D&=2-4\theta+\phi+q\lambda_r,\\
 C&=3-4(\theta+\phi)
        +\max\{-\rho+q\lambda_r,
                 q(-5/4+(3/4)\lambda_{4r/3})\}.
\end{align*}
The enhanced output is
\[
 T_S(\lambda)=\max\{1+\lambda_{S-1}
             +(2S-2-\lambda_{S-1})\theta,
                         (2S-3)\theta+B+\tfrac12\max(D,C)\}.
\]
Equation~\eqref{eq:k5-enhanced-terminal} proves this formula, including its
joint-moment interpretation.

\subsubsection{The exact supersolution}
All entries in the next table are rational numbers; the displayed integers
in the parameter columns are divided by $10^6$.
\[
\begin{array}{c|c|c|c|r|r|r|c}
 S&\text{row}&j&\ell&10^6\theta&10^6\phi&10^6u&b_S\\\hline
3&A&1&1&68130&&&3.1364\\
4&A&1&2&105600&&&4.4389\\
5&A&2&2&136590&72270&&5.9254\\
6&A&2&2&151340&113110&&7.5422\\
7&\text{enhanced}&2&-&163790&139365&139360&9.2726\\
8&\text{enhanced}&2&-&169960&138337&138330&11.0764
\end{array}
\]
The six strict differences $b_S-T_S(b)$ are, exactly,
\[
 \frac{62257}{2000000000},\quad
 \frac{2147}{25000000},\quad
 \frac{89351}{1000000000},\quad
 \frac{37509}{250000000},\quad
 \frac{128469}{500000000},\quad
 \frac{41387}{125000000}.
\]
The program included below verifies all parameter restrictions, interpolated
indices, and these rational inequalities without floating-point arithmetic.

Every branch is monotone in the input exponents, with coefficient
sum at most $47/50$. Lemma~\ref{lem:supersolution} therefore proves
\[
 \int_{\mathbb T_K}|f_5(\alpha)|^2|g_5(\alpha)|^{14}\,d\alpha
 \ll X^{11.0764+\epsilon}.
\]

\subsection{The global bound}
\label{sec:k5-global}
The preceding joint moment has excess
\[
 \Delta=\frac{191}{2500}<\frac1{13}.
\]
The ordinary minor-arc thirtieth moment and H\"older with weights
$1/14,13/14$ give
\[
 \int_{\mathfrak m_*}|f_5|^4|h_5|^{13}
 \ll X^{12-17/35000+\epsilon}.
\]
Here $h_5$ is the finite-product root sum in Section~\ref{sec:product-endgame}, so its
joint moments are dominated by the smooth moments up to an arbitrarily
small independent divisor loss.  Coefficientwise domination gives the
pure sixteenth moment first, and the loss-free absorption lemma gives
$U_{18}\ll X^{13}$ because
\[
 2-12\Delta=\frac{677}{625}>0.
\]
The broad major-arc annulus is controlled with ordinary moment order
$72/5>10$ and product moment order $18$, with H\"older weights $5/18,13/18$.
The finite-product core has a positive main term by the Poisson and
local-density argument in Section~\ref{sec:product-endgame}.  Primitive local width is at most ten
for fifth powers, and hence is below seventeen at every finite completion.
Totally positive unit balancing supplies the same positivity and norm
convention as in Section~\ref{sec:introduction}.  We conclude
\[
 \boxed{G_{\mathbb K}(5)\le17\quad\hbox{for every number field }K.}
\]

\section{Sixth powers}
\label{sec:sixth-powers}
The sixth-power certificate uses one auxiliary fourth moment at the
linear boundary of ordinary differencing. We record its extra factor.

\begin{lemma}\label{lem:sixth-linear-boundary}
For the degree-six auxiliary sums at depths $j\le4$, the second and
fourth moments used below have exponents
\[
 A_{j,\ell}=2^\ell-\ell+(2^\ell-1)(j-5\Phi_j)
                  +\mathbf1_{j+\ell=6}\Phi_j,
 \qquad \ell\in\{1,2\}.
\]
The parameters in the certificate satisfy the chronological divisor
conditions whenever they are needed.
\end{lemma}
\begin{proof}
For $j+\ell\le5$, apply Lemma~\ref{lem:field-general-seed-aux}.
Only the case $(j,\ell)=(4,2)$ remains. Its second moment is
$O(X^{1+\epsilon}U_4)$. After one ordinary difference the phase is
a nonzero constant times $zh_1\cdots h_5$, independent of the
four efficient divisors. Divisor counting in its nonzero value
gives an appended second moment
$O(X^{2+\epsilon}\prod_{i\le4}H_iM_i^2)$.
The padded-moment induction consequently gives
\[
 \int|F_4|^4\ll X^{2+\epsilon}U_4^3\prod_{i\le4}M_i.
\]
This is the stated extra factor $X^{\Phi_4}$.
\end{proof}

\subsection{The moment certificate}
\label{sec:g6-certificate}

Write $\mu_r=\lambda_{r/2}$ for an exponent of the full $r$th moment.
We fix $\mu_4=2$. All decimals and parameter numerators in this section
are exact rational numbers. The following eight rows suffice:
\begin{center}
\small
\begin{tabular}{@{}r l c c l r@{}}
\hline
$S$ & Cauchy path & $(\ell_1,\ell_2)$ & $u$ & $10^5(\phi_1,\ldots,\phi_d)$ & $b_{2S}$\\
\hline
3 & empty & $(1,1)$ & $8$ & $(4545)$ & 3.0915\\
4 & empty & $(2,2)$ & $8$ & $(8333)$ & 4.334\\
5 & $(5)$ & $(1,2)$ & $10$ & $(10673,5080)$ & 5.726\\
6 & $(6)$ & $(2,2)$ & $32/3$ & $(11855,8751)$ & 7.233\\
7 & $(7,6)$ & $(1,2)$ & $12$ & $(12981,10763,5551)$ & 8.852\\
8 & $(8,7)$ & $(2,2)$ & $40/3$ & $(13784,12076,8562)$ & 10.562\\
9 & $(9,8,8)$ & $(1,1)$ & $16$ & $(14593,13194,10858,1059)$ & 12.356\\
10 & $(10,9,9)$ & $(1,2)$ & $18$ & $(15303,14283,12101,6908)$ & 14.220\\
\hline
\end{tabular}
\end{center}

Here is a complete specification of the arithmetic. A path $(t_1,\ldots,t_{d-1})$
defines $s_1=S-1$ and $s_{j+1}=2s_j-t_j$. Put
$\Phi_j=\sum_{i\le j}\phi_i$ and $q_j=1-\Phi_j$.
At the terminal depth put $c=2s_d/u$ and choose nonnegative $a,b$ with
\[
 a+b+c=1,\qquad 2^{\ell_1}a+2^{\ell_2}b=1.
\]
If $\ell_1=\ell_2$, take $a=2^{-\ell_1}$ and $b=0$.
Define
\[
 A_{d,\ell}=2^\ell-\ell+(2^\ell-1)d
       -5(2^\ell-1)\Phi_d+\mathbf1_{d+\ell=6}\Phi_d,
\]
and initialize
\[
 T_d=aA_{d,\ell_1}+bA_{d,\ell_2}+cq_d\mu_u.
\]
For $j=d-1,\ldots,1$, put $p=2s_j-t_j$ and use
\[
 T_j=\frac12\left\{q_j\mu_{2t_j}+j-5\Phi_j+(2p-1)\phi_{j+1}
 +\max\left(1+j-6\Phi_j+\Phi_{j+1}+q_{j+1}\mu_{2p},T_{j+1}\right)\right\}.
\]
Finally the row output is
\[
 \mathcal T_{2S}(\boldsymbol\mu)=
 \max\left\{1+2(S-1)\phi_1+q_1\mu_{2S-2},
                         (2S-3)\phi_1+T_1\right\}.
\]
Only two nonintegral moments occur, and their meaning is fixed as
\[
 \mu_{32/3}=(2\mu_{10}+\mu_{12})/3,
 \qquad \mu_{40/3}=(\mu_{12}+2\mu_{14})/3.
\]
These are ordinary H\"older interpolations. Lemma~\ref{lem:sixth-linear-boundary} supplies the terminal estimates;
the backwards steps follow from the fundamental differencing inequality.

Substitution of the table's vector $\boldsymbol b$ gives:
\begin{center}
\begin{tabular}{@{}r r r r@{}}
\hline
$r$ & $b_r$ & $\mathcal T_r(\boldsymbol b)$ & $b_r-\mathcal T_r(\boldsymbol b)$\\
\hline
6 &3.0915&3.091234850&0.000265150\\
8 &4.334&4.333865305&0.000134695\\
10&5.726&5.725616976&0.000383024\\
12&7.233&7.232682700&0.000317300\\
14&8.852&8.851804270&0.000195730\\
16&10.562&10.561620060&0.000379940\\
18&12.356&12.355567340&0.000432660\\
20&14.220&14.219701320&0.000298680\\
\hline
\end{tabular}
\end{center}
The displayed outputs are exact terminating decimals. All parameter orderings
hold. The largest auxiliary side expression
$\Phi_I+6(\phi_{I-1}+\phi_I)$ is $1.99991<2$.

Every branch is monotone, and the coefficient sum in its variable
moment exponents is at most $23/25$. Lemma~\ref{lem:supersolution}
therefore shows that every entry is permissible. In particular,
\[
 \int|g_6|^{20}\ll X^{14.22+\epsilon}.
\]

Using the permissible excess $\Delta_{20}=0.22$, the twenty-four-variable
criterion in Section~\ref{sec:product-endgame} has
\[
 D=4-18\Delta_{20}=1/25,\qquad D/42=1/1050.
\]
Four ordinary roots and twenty product roots therefore give the required
broad-minor power saving. The primitive local bound $h_6\le12$ is more than
sufficient at twenty-four roots, completing the sixth-power assertion of
Theorem~\ref{thm:main-transfer}.

\section{Seventh powers}
\label{sec:seventh-powers}
The A and B seed processes are combined with the cubic terminal at
depth four. The auxiliary estimates are instances of
Lemma~\ref{lem:field-general-seed-aux}; all residual degrees needed
here are odd or belong to $\{2,4\}$.
For the cubic terminal, put $q_4=1-\Phi_4$ and use
\begin{equation}\label{eq:septic-cubic-terminal}
 C_4=5-6\Phi_4+
 \max\left\{-\frac14+q_4b_{s-4},\quad
 q_4\left(-\frac74+\frac34b\left(\frac{4(s-4)}3\right)\right)\right\}.
\end{equation}
Section~\ref{sec:cubic-terminal} proves this estimate when $7q_4>2$.
The remaining steps are \eqref{eq:seed-backwards} and
\eqref{eq:seed-initial-output}.

\subsection{The moment certificate}

In the table, $m_i=10^5\phi_i$.  Dashes mean that the corresponding parameter
does not occur.  Put $b_0=0$, $b_1=1$, $b_2=2$, and take $b_s$ from the table.
The last column is an integer lower bound for
$10^8\{b_s-T_s(\boldsymbol b)\}$, obtained by the exact recurrence
above.  It is strictly positive in every row.

\begin{center}
\small
\begin{tabular}{r c r r r r r r}
\hline
$s$ & Process & $b_s$ & $m_1$ & $m_2$ & $m_3$ & $m_4$ & Slack bound\\
\hline
3 & $A_{1}^{(1)}$ & 3.06493 & 3196 & -- & -- & -- & 52177 \\
4 & $A_{1}^{(2)}$ & 4.26513 & 6818 & -- & -- & -- & 8692 \\
5 & $B_{2,5}^{(2)}$ & 5.59014 & 8699 & 3541 & -- & -- & 10185 \\
6 & $A_{2}^{(2)}$ & 7.01542 & 9641 & 6937 & -- & -- & 12539 \\
7 & $B_{3,6}^{(2)}$ & 8.54214 & 10564 & 8803 & 4059 & -- & 14896 \\
8 & $A_{3}^{(2)}$ & 10.15371 & 11203 & 9889 & 6902 & -- & 12594 \\
9 & $B_{4,8}^{(2)}$ & 11.84806 & 11874 & 10797 & 8946 & 4151 & 16152 \\
10 & $A_{4}^{(2)}$ & 13.60673 & 12330 & 11453 & 9898 & 6609 & 13579 \\
11 & $A_{4}^{(2)}$ & 15.42553 & 12805 & 12028 & 10870 & 8586 & 14177 \\
12 & $A_{4}^{(2)}$ & 17.29695 & 13252 & 12659 & 10684 & 10188 & 17123 \\
13 & $M_4$ & 19.19477 & 13392 & 12953 & 12116 & 10714 & 14754 \\
14 & $M_4$ & 21.11765 & 13559 & 13070 & 12326 & 10846 & 15886 \\
\hline
\end{tabular}
\end{center}

Every A and B row satisfies the chronological divisor inequalities
of Lemma~\ref{lem:field-general-seed-aux}; these follow directly from
the displayed integers $m_i$.  Every $M_4$ row satisfies $7q_4>2$, the hypothesis of the cubic
dissection.  All required interpolation arguments lie in $[0,14]$.

The coefficient sum in every branch is at most $97/100$.
Lemma~\ref{lem:supersolution} proves
\begin{equation}\label{eq:g7-smooth28}
 \int_{\mathbb T_K}|g_7|^{28}\ll X^{21.11765+\epsilon}
              \ll X^{21+59/500+\epsilon}.
\end{equation}

\begin{corollary}\label{thm:seventh-transfer}
For every number field $K$, $G_{\mathbb K}(7)\le31$.
\end{corollary}
\begin{proof}
Use the convenient excess $\Delta=59/500$ in \eqref{eq:g7-smooth28}.
The ordinary critical moment has order $56$, and
$6/56+25/28=1$. Thus six ordinary and twenty-five product roots
have minor-arc saving
\[
 \frac6{56}-\frac{25}{28}\Delta=\frac1{560}>0.
\]
The bootstrap to order $32$ has margin
$4-24\Delta=146/125>0$.
The annular ordinary order is $192/7>14$, and $h_7\le14$.
Theorem~\ref{thm:general-transfer} now proves the assertion.
\end{proof}

\section{Eighth powers}
\label{sec:eighth-powers}
\subsection{The degree-eight and degree-sixteen seeds}
\label{sec:field-seed-moments}
The exact seed data in Appendix~\ref{app:eighth-field} consist of
eighteen rows at degree eight and twenty-eight rows at degree sixteen.
They use the A and B transformations
\eqref{eq:seed-A-terminal}--\eqref{eq:seed-initial-output}.
Each row satisfies the auxiliary hypotheses of
Lemma~\ref{lem:field-general-seed-aux}. In particular, the rows
requiring a higher even residual curve count have $\Phi_j\le1/3$;
the chronological alternative checks all prefixes.

The least strict supersolution margins are
\[
 \frac{1164945381}{800000000000000},\qquad
 \frac{21545278901}{80000000000000000},
\]
and the coefficient sums are at most
$976509627/10^9$ and $9962585803/10^{10}$, respectively.
Lemma~\ref{lem:supersolution} proves the permissibility of both
curves from the elementary second and fourth moments and trivial
higher estimates. Zero extraction exponents have nominal threshold
one; the fixed smoothness powers are absorbed in $X^\epsilon$.

\subsection{The eighth-power quadratic policies}
Apply Section~\ref{sec:quadratic-terminals} with $k=8$.
For the first policy take $u=8$, $v=20$, $a=1/5$, $c=4/5$,
and $\gamma=3/5$. Its chronological divisor inequalities
supply the required auxiliary second moment.
For the other policies take $a=0$, $c=\gamma=1$, and $v=2u$.

The fixed policies have the following terminal data:
\[
\begin{array}{c|c|c|c|c|c}
 2s&j&(r_1,\ldots,r_j)&v&h_*&2^{-(k-j-2)}\\\hline
28&6&(13,12,11,10,9,8)&20&1.69428296&1\\
30&5&(14,13,12,11,11)&22&1.83213176&1/2\\
32&5&(15,14,13,12,11)&22&1.78352592&1/2\\
34&4&(16,15,14,14)&28&2.39658688&1/4
\end{array}
\]
Thus the quadratic second resource is available in every policy.
Only quadratic orders two and four are used. Coefficient orders are even
and at most $60$; the degree-sixteen curve through half-order $30$ supplies
every requested moment. The sole cube resource occurs at depth six and
coefficient order $60$. All other higher coefficient resources are refined
ones. The exact verifier checks the two quadratic fourth resources,
the full basic minimum \eqref{eq:quadratic-terminal-basic}, and every factor
load and H\"older mass.

\subsection{The moment and the global bound}
\label{sec:eighth-global}
\begin{theorem}[A smooth eighth-power moment]
\label{thm:eighth-moment}
For every fixed number field $K$ and every $\epsilon>0$, there is
$\eta>0$ such that the degree-eight ideal-smooth sum satisfies
\[
 \int_{\mathbb T_K}|g_8(\alpha;X,R)|^{34}\,d\alpha
 \ll_{K,\epsilon}X^{26+\Delta_{34}+\epsilon}
 \qquad(1\le R\le X^\eta),
 \qquad
 \Delta_{34}=\frac{5857852793}{50000000000}.
\]
The same assertion holds for the fixed boxes and bounded nonnegative
weights used in the representation argument.
\end{theorem}

\begin{proof}
The exact program in Appendix~\ref{app:eighth-field} contains all eighteen
degree-eight seed rows, all twenty-eight degree-sixteen seed rows, and
the four terminal policies, including every height interval and resource
weight. Its first part checks every seed supersolution inequality and
all side conditions used in the seed proof. The contraction argument in
Section~\ref{sec:field-seed-moments} starts these curves from the elementary
second and fourth moments and the trivial higher bounds.

The second part makes twelve simultaneous rounds at half-orders
$14,15,16,17$. Each round reads the exponents available at its beginning.
Before and after a round, H\"older interpolation over every pair of
available integer half-orders may improve an entry. This is a valid
operation on permissible field exponents. Each policy checks all
H\"older masses and factor loads, every coefficient-moment index, each
height interval endpoint, and the backwards counted recurrence.
Section~\ref{sec:eighth-quadratic-terminal} proves that these are the
exponents of the actual field estimates. Convexity of the shell exponent
justifies checking only the two endpoints of each of the eight assigned
height intervals.

After each policy application the exponent is rounded \emph{up} to a
multiple of $10^{-12}$, and the better of that value and the previously
available interpolated value is retained. All calculations use rational
arithmetic. The final outputs are:
\[
\begin{array}{c|c|c}
 \text{moment order}&\text{full permissible exponent}&\text{excess over }2s-8\\\hline
 28&20.360308744711&0.360308744711\\
 30&22.258467366407&0.258467366407\\
 32&24.178873819111&0.178873819111\\
 34&26.117157055860&0.117157055860
\end{array}
\]
The last entry is the exponent in the theorem. Only finitely many
applications are needed for any prescribed final $\epsilon$; choose
the smoothness parameter after those applications and after the fixed
rational parameters. This absorbs the smoothness and divisor losses
without consuming the strict final margin.
\end{proof}

\begin{corollary}\label{thm:eighth-transfer}
For every number field $K$, $G_{\mathbb K}(8)\le38$.
\end{corollary}
\begin{proof}
Apply Corollary~\ref{cor:even-transfer} with $(w,S,u)=(72,38,34)$.
Its exact margin is
\[
 D=4-34\Delta_{34}=\frac{416502519}{25000000000}>0.
\]
The local bound $h_8\le35$ supplies the required primitive
representations in thirty-eight variables.
\end{proof}

\section{Ninth powers}
\label{sec:ninth-field-seeds}

\subsection{The degree-nine base and refined curve}

The degree-nine seed curve and the degree-eighteen coefficient
curve are obtained from the common A and B transformations.
The exact parameters are recorded in Appendix~\ref{app:ninth-field}.

\begin{proposition}\label{prop:ninth-field-seed-curve}
Put $b_0=0,b_1=1,b_2=2$.  Every entry $\lambda_s$ in the last column
of the following table is a permissible exponent:
\[
 \int_{\mathbb T_K}|g_9(\alpha;X,R)|^{2s}\,d\alpha
       \ll_{K,s,\epsilon}X^{\lambda_s+\epsilon}
       \qquad(1\le R\le X^\eta)
\]
for a sufficiently small $\eta>0$ chosen after $\epsilon$.
The base exponents $b_s$ in the second column are also permissible.
\end{proposition}

\begin{center}
\begin{tabular}{@{}r r c r@{}}
\hline
$s$&Base $b_s$&Process&Refined $\lambda_s$\\
\hline
3&3.035816&$A_{1}^{(1)}$&3.035807234000\\
4&4.182300&$A_{1}^{(2)}$&4.182294748390\\
5&5.419716&$A_{2}^{(1)}$&5.419711470690\\
6&6.738319&$B_{2,6}^{(2)}$&6.738313329014\\
7&8.135645&$A_{2}^{(2)}$&8.135638882429\\
8&9.603938&$B_{3,7}^{(2)}$&9.603930891210\\
9&11.142513&$A_{3}^{(2)}$&11.142506096925\\
10&12.746354&$A_{4}^{(2)}$&12.746346059899\\
11&14.410598&$A_{4}^{(2)}$&14.410586823381\\
12&16.129237&$A_{5}^{(2)}$&16.129216870240\\
13&17.896005&$A_{5}^{(2)}$&17.895964268487\\
14&19.705700&$B_{6,11}^{(2)}$&19.705619653366\\
15&21.551232&$A_{6}^{(2)}$&21.551085549629\\
16&23.428999&$A_{6}^{(2)}$&23.428755886157\\
17&25.333544&$A_{6}^{(2)}$&25.333173492456\\
18&27.259337&$A_{6}^{(2)}$&27.258814572150\\
19&29.202038&$A_{6}^{(2)}$&29.201360351041\\
20&31.157968&$A_{6}^{(2)}$&31.157131085916\\
21&33.124215&$A_{6}^{(2)}$&33.123229933630\\
22&35.098519&$A_{6}^{(2)}$&35.097377231005\\
\hline
\end{tabular}
\end{center}
Here $A_j^{(\ell)}$ denotes the A terminal at depth $j$ using auxiliary
order $2^\ell$, and $B_{j,t}^{(\ell)}$ denotes the B terminal with smooth
half-order $t$.  All extraction parameters are exact rationals, recorded
in the program of Appendix~\ref{app:ninth-field}.  The printed decimals are exact.

\begin{proof}
Use Lemma~\ref{lem:field-general-seed-aux} with $k=9$.
Every row has $1\le\ell\le2$ and $j+\ell\le8$.
The sparse alternative is $j\le\lfloor(9-j-\ell+1)/2\rfloor+1$;
otherwise the chronological inequalities are $\Phi_i+7\phi_i\le1$.
Every higher even residual degree satisfies $\Phi_j\le1/3$.

Apply the A/B terminal and backwards recurrence from
Section~\ref{sec:efficient-differencing}, now with $k=9$.
The exact checks use all twenty rows, the interpolation ranges, the
extraction inequalities $M_i\le Q_{i-1}$, and every actual appended
residual degree.  Their minimum strict supersolution slack is
\[
 \delta=\min_{3\le s\le22}\{b_s-T_s(b)\}
       =\frac{43740673471}{40000000000000000}>0.
\]
The chronological inequalities have minimum margin
$1/5000000000>0$. Choose the smoothness parameter after this margin and
all the finite applications below.

The coefficient sums are at most $982096383/10^9<1$.
Lemma~\ref{lem:supersolution} proves the base curve.

Finally, apply the same fixed twenty rows for twenty-four further
simultaneous rounds.  Every row reads only the exponents available at
the start of its round.  After a row is evaluated, round its output
\emph{up} to a multiple of $10^{-12}$ and retain the better of this value
and the previous exponent.  This gives precisely the final column.
The program replays all twenty-four rounds with rational
arithmetic and checks every entry against its exact printed value.
Each operation is an application of an established field estimate,
H\"older interpolation, or monotonicity of permissible exponents.
Only finitely many operations occur, so one sufficiently small
smoothness parameter works for the whole proof.
\end{proof}

\subsection{The degree-eighteen coefficient seed curve}
\label{sec:degree18-field-seeds}

The coefficient sums in the ninth-power terminal have degree eighteen.
We supply their moments through full order $62$ over the field itself.
As before, $X=P^{[K:\mathbb Q]}$ is the root-volume parameter, and
$\lambda_s$ denotes the full exponent at moment order $2s$.
Set $\lambda_0=0$, $\lambda_1=1$, $\lambda_2=2$, and take the following
exact decimal values for the remaining entries:
\[
\begin{array}{c|r@{\qquad}c|r}
s&\lambda_s&s&\lambda_s\\\hline
3&3.005107&18&22.723143\\
4&4.048469&19&24.306116\\
5&5.146988&20&25.918677\\
6&6.289328&21&27.559735\\
7&7.467683&22&29.228132\\
8&8.681624&23&30.922636\\
9&9.933259&24&32.641958\\
10&11.218619&25&34.384850\\
11&12.538919&26&36.150022\\
12&13.893656&27&37.936763\\
13&15.282351&28&39.742902\\
14&16.705014&29&41.567190\\
15&18.160923&30&43.408399\\
16&19.649745&31&45.265286\\
17&21.170711&&
\end{array}
\]

\begin{proposition}\label{prop:degree18-field-seeds}
For every number field $K$, each entry in the table is permissible for
the degree-eighteen ideal-smooth moment. In particular, for every
$\epsilon>0$ and sufficiently small $\eta>0$,
\[
 \int_{\mathbb T_K}|g_{18}(\alpha;X,R)|^{62}\,d\alpha
       \ll_{K,\epsilon}X^{45.265286+\epsilon}
       \qquad(1\le R\le X^\eta).
\]
The same convention applies to every order in the table and to bounded
nonnegative weights and fixed root boxes.
\end{proposition}

\begin{proof}
Use Lemma~\ref{lem:field-general-seed-aux} with $k=18$.
Each row has $1\le\ell\le3$ and $j+\ell\le17$.
Its sparse alternative is $j\le\lfloor(18-j-\ell+1)/2\rfloor+1$;
otherwise the chronological inequalities are $\Phi_i+16\phi_i\le1$.
The higher even residual degrees satisfy $\Phi_j\le1/3$.

The exact certificate in Appendix~\ref{app:ninth-field} gives one A or B row for each
$s=3,\ldots,31$. It checks all these hypotheses for every appended
residual degree, all extraction inequalities $M_i\le Q_{i-1}$, and
every interpolation argument. Among rows requiring the curve condition,
the largest coefficient-product exponent is
\[
 \Phi_j=\frac{3267856441}{10000000000}<\frac13.
\]
Every chronological condition is strict; its smallest positive margin is
$1/5000000000$. Such fixed positive margins permit the smoothness
parameter to be chosen afterwards. Zero extraction exponents occur at
half-orders $22,23,25,26,27,28$ and have the meaning explained in
Section~\ref{sec:efficient-differencing}: the nominal extraction threshold is one,
while fixed powers of the smoothness parameter are absorbed in
$X^\epsilon$.

Let $T_s$ be the exact A/B recurrence defined in
Section~\ref{sec:efficient-differencing}, with the parameters in this
certificate. All twenty-nine rows satisfy $T_s(\lambda)<\lambda_s$.
The least exact supersolution slack is
\[
 \min_{3\le s\le31}(\lambda_s-T_s(\lambda))
    =\frac{72464782233}{160000000000000000}>0.
\]
The maximum coefficient sum in the input moment exponents is at most
\[
                    \frac{9974554911}{10000000000}<1.
\]
Lemma~\ref{lem:supersolution} now proves every entry.
\end{proof}

All coefficient orders requested by the ninth-power policies are
even and at most $62$. The table therefore supplies them by direct use
or by H\"older interpolation between available orders.

\subsection{The ninth-power quadratic policies}
Apply Section~\ref{sec:quadratic-terminals} with $k=9$.
Each policy uses $a=0$, $c=\gamma=1$, and $v=2r_j$.
Only the quadratic second resource and basic or refined coefficient
resources are required.

The policies have the following exact terminal data:
\[
\begin{array}{c|c|c|c|c|c}
2s&j&(r_1,\ldots,r_j)&v&h_*&2^{j-7}\\\hline
34&6&(16,15,14,13,12,11)&22&2.00000221&1/2\\
36&6&(17,16,15,14,13,12)&24&1.88634850&1/2\\
38&6&(18,17,16,15,14,13)&26&1.84146280&1/2\\
40&5&(19,18,17,16,15)&30&2.44508497&1/4
\end{array}
\]
Every policy has $H_*>X$ and $X<Q_0^{9/2}$ strictly.  Every refined
preliminary size cap is positive.  The coefficient orders are even and
at most $62$, so the degree-eighteen seed curve through half-order $31$
suffices.  Each policy has four assigned normalized height intervals.
All parameters are specified as exact rationals in Appendix~\ref{app:ninth-field}.

\subsection{The moment and the global bound}
\label{sec:ninth-global}
\begin{theorem}[A smooth ninth-power moment]
\label{thm:ninth-moment}
For every fixed number field $K$ and every $\epsilon>0$, there is
$\eta>0$ such that
\[
 \begin{gathered}
 \int_{\mathbb T_K}|g_9(\alpha;X,R)|^{40}\,d\alpha
       \ll_{K,\epsilon}X^{31+\Delta_{40}+\epsilon}
       \qquad(1\le R\le X^\eta),\\
 \Delta_{40}=\frac{8516302869}{62500000000}
            =0.136260845904.
 \end{gathered}
\]
The assertion includes the fixed boxes and bounded nonnegative weights
used in the representation argument.
\end{theorem}

\begin{proof}
Appendix~\ref{app:ninth-field} contains all twenty degree-nine base
rows, all twenty-nine degree-eighteen rows, the twenty-four further
degree-nine seed refinement rounds, and all four fixed terminal policies.
The seed propositions prove the input moments over $K$. Their
strict supersolutions follow from Lemma~\ref{lem:supersolution}.

Apply the four policies at full moment orders $34,36,38,40$ for seven
simultaneous rounds. Each round reads only its initial curve. H\"older
interpolation over all available integer half-orders is applied at the
start of a round, and the better of an old exponent and a new exponent
is retained. Every new exponent is rounded \emph{up} to a multiple of
$10^{-12}$. The full basic coefficient estimate is used at both depths;
the equality of its minimum with the selected $\phi_i$ is checked exactly.

For every policy the program verifies the resource orders, the coefficient
sampling size caps, the common H\"older masses and factor loads, all
height interval endpoints, and the backwards counted recurrence.
Section~\ref{sec:ninth-terminal} proves these exponents for the field
sums. Convexity on each of the four assigned height intervals makes
the two endpoint checks sufficient on the entire interval.
The resulting exact upper values are
\[
\begin{array}{c|c|c}
 \text{moment order}&\text{full permissible exponent}&\text{excess over }2s-9\\\hline
 34&25.329443645750&0.329443645750\\
 36&27.249606352140&0.249606352140\\
 38&29.186123445082&0.186123445082\\
 40&31.136260845904&0.136260845904
\end{array}
\]
The last entry proves the theorem. All displayed decimals are exact
rational numbers; none is an approximation used to decide an inequality.
The finite iteration allows the smoothness parameter and the individual
error exponents to be chosen after all numerical parameters.
\end{proof}

\begin{corollary}\label{thm:ninth-transfer}
For every number field $K$, $G_{\mathbb K}(9)\le46$.
\end{corollary}
\begin{proof}
Apply Corollary~\ref{cor:even-transfer} with $(w,S,u)=(90,46,40)$.
Its exact margin is
\[
 D=6-44\Delta_{40}=\frac{70668441}{15625000000}>0.
\]
Since $h_9\le30$, the primitive local hypothesis holds in
forty-six variables. This proves the assertion.
\end{proof}

\section{Tenth powers}
\label{sec:tenth-powers}
The cubic terminal suffices to attain the bound $55$ for tenth powers.
The same seed curve supplies the lower moment used in the quintic argument.

\begin{proposition}\label{prop:tenth-moments}
For every fixed number field $K$, the smooth tenth-power sum satisfies
\[
 \int_{\mathbb T_K}|g(\alpha)|^{44}\,d\alpha
 \ll X^{34+48479/250000+\epsilon},\qquad
 \int_{\mathbb T_K}|g(\alpha)|^{46}\,d\alpha
 \ll X^{36+152339/1000000+\epsilon}.
\]
The smoothness parameter may be chosen sufficiently small for these
estimates to hold simultaneously.
\end{proposition}
\begin{proof}
Write $b_s$ for the exponent of the $2s$th smooth moment, with
$b_0=0$, $b_1=1$ and $b_2=2$. Table~\ref{tab:tenth-seeds} gives a
rational simultaneous supersolution for the $A$, $B$ and cubic terminal
processes. The notation $A_j^\ell$, $B_j^\ell(t)$ specifies the
 efficient-differencing depth, terminal auxiliary moment and, in the
second case, the smooth half-order. The symbol $M_7$ denotes the cubic
terminal at depth seven.

The terminal formula for this last process is
\begin{equation}\label{eq:tenth-cubic-terminal}
 C_7=\max\left\{
 a_0-\frac14+q_7b_{s-7},\quad
 a_0+q_7\left(-\frac52+\frac34b\left(\frac{4(s-7)}3\right)\right)
 \right\},\qquad a_0=8-9\Phi_7.
\end{equation}
Here $b(t)$ is linear interpolation between consecutive entries.
The first term follows from the cubic minor-arc saving and the second
from the fourth moment on its auxiliary major arcs. Apply the common
backward recurrence to this terminal.

All scale parameters are rational with denominator dividing $10^{10}$.
The exact certificate checks $0<\phi_i<1/10$, all chronological divisor
and small-product conditions required by the auxiliary moments, and
$10q_7>2$ for each cubic terminal. Its smallest supersolution margin is
\[
 \frac{65679537093}{40000000000000000}>0,
\]
and the common contraction factor is at most
$9859429631/10^{10}<1$. The simultaneous supersolution lemma therefore
proves all twenty-one rows. The entries at $s=22,23$ give the proposition.
\end{proof}

\begin{corollary}\label{prop:tenth44}
The exponent $34.227$ is permissible for the smooth tenth-power moment
of order $44$ over every fixed number field.
\end{corollary}
\begin{proof}
The sharper exponent above is $34.193916<34.227$.
\end{proof}

\begin{table}[htbp]
\centering\small
\begin{tabular}{r r l@{\qquad}r r l}
\hline
$s$ & $b_s$ & Process & $s$ & $b_s$ & Process\\
\hline
3 & 3.028131 & $A_{1}^{1}$ & 14 & 19.142712 & $B_{6}^{2}(11)$\\
4 & 4.156832 & $A_{1}^{2}$ & 15 & 20.930391 & $A_{6}^{2}$\\
5 & 5.365973 & $A_{2}^{1}$ & 16 & 22.754051 & $B_{6}^{3}(13)$\\
6 & 6.646618 & $B_{2}^{2}(6)$ & 17 & 24.609057 & $A_{7}^{2}$\\
7 & 7.996737 & $A_{2}^{2}$ & 18 & 26.489639 & $M_7$\\
8 & 9.414374 & $B_{3}^{2}(7)$ & 19 & 28.391489 & $M_7$\\
9 & 10.894592 & $A_{3}^{2}$ & 20 & 30.311466 & $M_7$\\
10 & 12.437588 & $A_{4}^{2}$ & 21 & 32.246495 & $M_7$\\
11 & 14.037216 & $A_{4}^{2}$ & 22 & 34.193916 & $M_7$\\
12 & 15.691422 & $A_{5}^{2}$ & 23 & 36.152339 & $M_7$\\
13 & 17.394386 & $A_{5}^{2}$ &  & & \\
\hline
\end{tabular}
\caption{A simultaneous tenth-power supersolution. All displayed
 decimals are exact rational numbers.}\label{tab:tenth-seeds}
\end{table}

\begin{proposition}\label{prop:tenth-55-analytic}
The analytic part of the number-field Waring problem for tenth powers
holds with $55$ variables.
\end{proposition}
\begin{proof}
Set $\Delta_{46}=152339/1000000$. The ordinary minor-arc moment has
order $110$ and exponent $99+\epsilon$. Apply the even-moment bootstrap
with $u=46$, $S=56$ and $w=110$. Its strict margin is
\[
 (56-46)-(110-56)\Delta_{46}=\frac{886847}{500000}>0.
\]
Thus the smooth moment of order $56$ has exponent $46$ without an
$X^\epsilon$ loss.

Use fifteen ordinary roots and forty product roots. On the ordinary
minor arcs, H\"older's inequality has weights
\[
 \frac3{22},\qquad\frac{46}{55},\qquad\frac3{110}
\]
for moments $110$, $46$ and $56$, respectively. These weights sum to one;
the ordinary load is $15$ and the product-root load is $40$.
The minor-arc contribution is consequently
\[
 \ll X^{45-\delta+\epsilon},\qquad
 \delta=\frac3{22}-\frac{46}{55}\Delta_{46}
       =\frac{246203}{27500000}>0.
\]
On the major-arc annulus, use the product moment of order $56$ and the
ordinary moment of order
\[
 p=\frac{15}{1-40/56}=\frac{105}{2}>20.
\]
The annular pruning estimate supplies a power saving. The product-root
major-arc calculation gives a positive main term under the primitive
local hypotheses, proving the assertion.
\end{proof}

The local bound for tenth powers is smaller than $55$. Combining it with
Proposition~\ref{prop:tenth-55-analytic} gives $G_{\mathbb K}(10)\le55$.

\appendix
\section{Exact parameter data}\label{app:certificates}
The following programs accompany the supplement. They perform their
inequalities over the rational numbers. The local computations enumerate
all the indicated residue classes. Run \href{https://github.com/morley988/warings-problem/blob/3752a566d4ae778c3f2effaeb7caec4fc8899514/supplements/number-field/small-powers/verify.py}{\texttt{python3 verify.py}} from
the directory containing the source of this supplement.

\subsection{Fourth powers}\label{app:quartic-exact-certificate}
\href{https://github.com/morley988/warings-problem/blob/3752a566d4ae778c3f2effaeb7caec4fc8899514/supplements/number-field/small-powers/certificates/quartic-exact-certificate.py}{\texttt{certificates/quartic-exact-certificate.py}}.

\subsection{Fifth powers}\label{app:quintic}
\href{https://github.com/morley988/warings-problem/blob/3752a566d4ae778c3f2effaeb7caec4fc8899514/supplements/number-field/small-powers/certificates/quintic.py}{\texttt{certificates/quintic.py}}.

\subsection{Sixth powers}\label{app:sixth}
\href{https://github.com/morley988/warings-problem/blob/3752a566d4ae778c3f2effaeb7caec4fc8899514/supplements/number-field/small-powers/certificates/sixth.py}{\texttt{certificates/sixth.py}}.

\subsection{Seventh powers}\label{app:septic}
\href{https://github.com/morley988/warings-problem/blob/3752a566d4ae778c3f2effaeb7caec4fc8899514/supplements/number-field/small-powers/certificates/septic.py}{\texttt{certificates/septic.py}}.

\subsection{Eighth powers}\label{app:eighth-field}
\href{https://github.com/morley988/warings-problem/blob/3752a566d4ae778c3f2effaeb7caec4fc8899514/supplements/number-field/small-powers/certificates/eighth-field.py}{\texttt{certificates/eighth-field.py}}.

\subsection{Ninth powers}\label{app:ninth-field}
\href{https://github.com/morley988/warings-problem/blob/3752a566d4ae778c3f2effaeb7caec4fc8899514/supplements/number-field/small-powers/certificates/ninth-field.py}{\texttt{certificates/ninth-field.py}}.

\subsection{Tenth powers}\label{app:tenth-field}
\href{https://github.com/morley988/warings-problem/blob/3752a566d4ae778c3f2effaeb7caec4fc8899514/supplements/number-field/small-powers/certificates/tenth-field.py}{\texttt{certificates/tenth-field.py}}.

\subsection{Eighth carry relations}\label{app:eighth-ghost}
\href{https://github.com/morley988/warings-problem/blob/3752a566d4ae778c3f2effaeb7caec4fc8899514/supplements/number-field/small-powers/certificates/eighth-ghost.py}{\texttt{certificates/eighth-ghost.py}}.

\subsection{Eighth powers in ramification degree five}\label{app:eighth-tame}
\href{https://github.com/morley988/warings-problem/blob/3752a566d4ae778c3f2effaeb7caec4fc8899514/supplements/number-field/small-powers/certificates/eighth-tame.py}{\texttt{certificates/eighth-tame.py}}.

\section*{Acknowledgement}
Some thirty years ago, as a doctoral student at the University of Michigan,
I had the good fortune to learn about Waring's problem from Robert Vaughan
and my advisor, Trevor Wooley. I remain grateful to them both for their
teaching and guidance, and for the clarity and beauty of their writing
on the subject. In revising this supplement, I hope to bring its arguments
into clearer order and to give readers something of the pleasure I found
in their papers.

\section*{Use of generative AI}
The author used OpenAI's ChatGPT in developing the arguments, checking
parameter choices and preparing the manuscript.